\documentclass[12pt]{book}
\usepackage[spanish]{babel}
\usepackage[latin1]{inputenc}
\usepackage{courier}

\usepackage{pstricks}
\usepackage{pst-all}

\usepackage{pst-plot}

\usepackage{makeidx}
\usepackage{makeidx} 
\makeindex
\newcounter{ecu}
\begin{document}

\title{ELECTROMAGNETISMO Y GEOMETRIA}         
\author{JOSE DEL CARMEN RODRIGUEZ SANTAMARIA }        
\date{ }          
\maketitle

\bigskip
\bigskip

$$PROLOGO$$
\date{}
\bigskip

\maketitle

\
\
\

\bigskip
\normalsize
\color{black}

El presente trabajo ofrece una introducción a la física matemática moderna. Es física, hablamos de electromagnetismo,  es matemática, hablamos de geometría, y es moderna: unimos las dos cosas.

El objetivo de la física-matemática  es el de describir y entender el universo a nivel fundamental. Describirlo significa tener una teoría que prediga resultados respaldados por  el experimento. Entender nuestro universo significa ponerlo en contextualización de las infinidades de mundos posibles que tan sólo pueden explorarse matemáticamente. La contextualización puede resultar como un subproducto de la descripción de algún aspecto del universo real por medio de una estructura matemática, la cual automáticamente queda en relación con todas las otras.

La subcultura de la física-matemática  hace  énfasis tanto en el poder descriptivo (el cual incluye el  predictivo), como en la naturalidad, sencillez, elegancia, generalidad, rigor, autoanálisis,  proyección y modernidad de las teorías. La presente contribución clasifica como introducción por empezar a partir del cálculo vectorial y por tener un poquito de todo, pero no a tal grado  que intoxique a los interesados.

Nosotros  seguimos  en el presente trabajo la propuesta   de dividir la materia en dos partes: partículas y campos.

Como partículas distinguidas tenemos el electrón, el protón,  el neutrón, los neutrinos, las cuales son suficientes para explicar los átomos. En la actualidad se acepta que tanto el protón como el neutrón no son elementales sino que están compuestas de quarks.   Asumimos que las partículas   interactúan entre ellas pero no directamente   sino que lo hacen a través del campo que sirve de mediador. Las observaciones experimentales sugieren que hay cuatro tipos fundamentales de interacciones: la electromagnética (responsable de la formación de  moléculas y de que haya luz), la débil (causa de la radioactividad), la fuerte (la que hace que el núcleo atómico no estalle debido a la repulsión eléctrica entre protones) y la gravitatoria (la que mantiene la luna pegada a la tierra para que no se vaya).

Tenemos una buena descripción de las interacciones elementales, formulada por El Modelo Estándar, en el cual se pegan, casi a la fuerza, varios parches provenientes de una familia de teorías que se consideran muy prometedoras: las teorías gauge.  Nuestro objetivo es presentar un estudio completo de uno de dichos parches: la teoría gauge del electromagnetismo, la cual no se puede entender sin conocer un poco de relatividad general.

Nuestro estudio de la interacción electromagnética comenzará desde los más simple (las ecuaciones de Maxwell)  hasta lo más sofisticado (el fibrado electromagnético ). Veremos en los cinco primeros capítulos las leyes de Maxwell e introduciremos la mecánica cuántica para poder probar el teorema fundamental: La teoría  electromagnética es una teoría gauge con grupo de invariancia U(1). Algo tan abstracto contiene un aspecto sencillo: la fase de la función de onda tiene en cada punto una arbitrariedad en el ángulo que se debe tomar como  cero. Por lo tanto,  las predicciones físicas no deberán depender de tal arbitrariedad.

Al desarrollar la anterior discusión puede verse, casi de paso, de qué manera las leyes de Maxwell logran una unificación de la óptica, la electricidad, el magnetismo y la electrodinámica. La gran ilusión de la física-matemática  moderna es la de encontrar la forma de extender la covertura de dicha unificación hasta envolver todas las cuatro interacciones fundamentales. Nosotros vamos a enfrentar el terrible problema de volver dicha ilusión una realidad y lucharemos por una pre-unificación: si realmente las cuatro interacciones son casos especiales de una única interacción, entonces las diversas interacciones tienen que poderse expresar, en esencia, utilizando el mismo formalismo matemático.

Debido a que la gravitación admite   una formulación geométrica que ha probado ser muy exitosa, la llamada relatividad general, se ha investigado intensamente la posibilidad de geometrizar todas las interacciones. Eso significa dos cosas. Primera: describir todas las interacciones en términos de un punto,   o tal vez de una cuerda, y de lo que sucede en su inmediata vecindad.  Segunda: aplicar la metodología de la relatividad general al estudio de espacios extendidos que puedan reflejar estados internos de la partícula.

Ese gran proyecto tiene un punto de partida muy natural que se desarrolla en la segunda parte:  expresar las leyes de Maxwell en el lenguaje de la relatividad, que es aquella que estudia la reformulación de las  leyes de la física en diferentes marcos de referencia. Esto lo haremos en los capítulos seis y siete.

 En el capítulo ocho introduciremos los elementos de geometría básicos que permiten entender los conceptos implicados en la versión geométrica de la gravitación, capítulo nueve, para reformular, en la última  parte,  el electromagnetismo en el mismo lenguaje de la gravitación: derivada covariante y curvatura.

La versión geométrica moderna de todo lo anterior utiliza la estructura matemática conocida como fibrado principal. Hemos introducido el tema en los capítulos 10 y 11 sacando provecho  de la sencillez del electromagnetismo: el fibrado principal resultante es simplemente el producto cartesiano del espacio-tiempo   por el grupo $ U(1)$.

El presente trabajo se empeña en conectar las grandes ideas que mueven la ciencia de las partículas fundamentales con el sentido común, con el instinto que movió a los creadores de los significados contextuales. Y también representa un pequeño adelanto  en el complicado arte de convertir intuiciones difusas en ideas claras, precisas, bien formuladas y elegantes gracias al denuedo de unos 5000 años de cultura matemática. He procurado  que los prerrequisitos de física y de geometría sean lo mínimo posible, explicando casi todo lo necesario. Pero a decir verdad, convendría al lector tomar paralelamente un curso de geometría diferencial.

Apreciado lector:  siéntase libre de tomar la fuente Latex del presente documento y de utilizarla, en todo o en parte, en original o modificada, para hacer su propio material, que le quede a su gusto y medida.

\bigskip

\vspace{\baselineskip}
\begin{flushright}\noindent
Bogot\'a, Colombia \hfill {\it Jos\'e Rodr\'\i guez}\\
Junio 2008 \hfill {\it }\\
\end{flushright}

\tableofcontents

\chapter{EL CAMPO ELECTROMAGNETICO}       

La expresión 'campo electromagnético' es un orgullo de la modernidad:  denota una gran \textbf{unificación} \index{unificación} entre cuatro interacciones bien diferentes: la eléctrica, la magnética, la electrodinámica (entre partículas cargadas que se mueven), y la óptica.

\
\
\

\section{ELECTROSTATICA}    

El \textbf{campo eléctrico} \index{campo eléctrico} es, por un lado, una gran abstracción inventada para explicar lo más elemental de la conducta de la materia a nuestro alrededor. Y por el otro, es un teoría que ha podido ser actualizada hasta llegar a gozar del respaldo de la naturaleza de una manera espectacular.

\

Reinventemos el campo eléctrico. Podemos empezar con una simple observación: los gases se dejan comprimir. Gracias a eso tenemos neumáticos, la gran idea que nos permite viajar en bicicleta o en carro como si fuese en una nube. Los gases se dejan comprimir porque están compuestos de pequeñas partículas que se repelen con una fuerza finita y, seguramente, no muy grande. Hoy sabemos que dicha fuerza no es elemental, pues se explica como la resultante de fuerzas entre moléculas. Decimos que una interacción es elemental cuando descomponer el esquema  no adiciona ni sencillez ni poder explicativo. Creemos que la fuerza eléctrica entre electrones es elemental, pero no la fuerza eléctrica entre neutrones.

En nuestro intento de reinventar el campo eléctrico, denominemos  a la fuerza de repulsión entre partículas como fuerza eléctrica.  A las partículas que se repelen se denominan 'cargadas', o sea que tienen una carga eléctrica, la cual crea una fuerza de repulsión en otras cargas colocadas en su vecindad.

Cuando hablamos de fuerza eléctrica de repulsión, pensamos en dos partículas cargadas que se repelen. La idea de campo se crea así: el campo eléctrico es la fuerza experimentada por una carga de prueba de magnitud fija debido a la presencia de otras cargas que se consideran inmóviles. Como dicha fuerza es un vector, cada carga modifica el espacio que la rodea poniendo un vector en cada punto. Ese conjunto de vectores se llama campo eléctrico.

Vemos así que nuestro campo repulsivo necesita como base a la teoría atómica. Dicha teoría, inventada por Demócrito de Abdera en el s.-v., en la edad de oro de los griegos, se basa en dos elementos. Uno experimental, que uno puede partir un pedazo de palo en dos partes y cada una de ellas en otras dos y así sucesivamente. El segundo elemento es teórico: que dicha división no se puede seguir indefinidamente. La teoría atómica es interesante por un postulado filosófico que da realismo y universalidad al materialismo y que fue proclamado aún por los antiguos: todo lo que existe    consta de átomos (Aún los dioses!).

Para pasar del materialismo atómico como sistema filosófico al atomismo como reduccionismo científico hay mucho trecho. Comencemos con el problema  de explicar, basados en la teoría atómica, por qué hay tanta variedad de objetos en el universo. Eso sólo se pudo hacer con el advenimiento de la química que es el arte, el delicado arte, de obtener substancias casi puras para poder así estudiar sus propiedades. Se dice que   John Dalton es el padre de la teoría atómica por haber formulado una versión cuantitativa  de dicha teoría en 1808. El siguiente gran paso en esa dirección fue el dado por Dmitry Ivanovich Mendeleiev al proponer en 1869 y 1871 la tabla periódica: hay unas substancias básicas que se diferencian entre sí por su peso atómico (referido al peso del hidrógeno) y cuyas propiedades químicas son  función periódica, aunque algo difusa, de dicho peso.

Tenemos entonces unos átomos, los cuales existen en variadas versiones, que se repelen en el gas. De hecho, eso demuestra que la interacción repulsiva no es la gravitacional, la cual es siempre atractiva. Además, la fuerza de repulsión entre partículas debe ser muchísimo más fuerte que la gravitatoria, pues si no fuese así, la gravitatoria colapsaría los cuerpos.

Pero no todos los cuerpos se dejan comprimir fácilmente. Tal es el caso del agua. Gracias a ello, el agua puede transmitir ondas sonoras por miles de kilómetros y, si no fuese por la contaminación acústica, permitiría que las ballenas se comunicasen entre ellas a lo largo del Atlántico de un hemisferio de la tierra al otro. Será acaso que el agua no se deja comprimir por la misma razón por la cual cada vez es más difícil aumentarle la presión a un neumático? Podríamos suponer entonces que  la fuerza de repulsión tiene que ser inversamente proporcional a alguna función decreciente de la distancia.

Pero si se trata tan sólo de repulsión, entonces qué mantiene pegadas a las partículas para que existan los sólidos?  Y por qué hay que calentar a los  líquidos para evaporarlos? y por qué hay repulsión entre los átomos de agua a corta distancia, para que impidan comprimirla, y atracción a larga distancia, para que al caer la lluvia y resbalar por la ventana, dos gotitas se peguen en una sola? Y además de todo, una gota de agua  que cae sobre una plancha caliente chisporretea  encimita de ella pero sin pegársele hasta que se evapora totalmente.

Dijimos de paso que las gotas de agua se atraen a corta distancia. Debemos reconocer por tanto que además de interacción de repulsión existe interacción de atracción.

Atracción y repulsión son términos que describen una interacción entre elementos. Sucede dicha interacción entre pares de elementos, entre tripletas o entre todas las cargas al tiempo?

Pregunta crítica: Cómo es posible que existan objetos cuya estructura sea estable, una mesa, una botella, un libro,  siendo que además entre ellos o entre sus partes no parece existir atracción ni repulsión? A estos cuerpos llamémoslos neutros.

\section{ESPACIOS VECTORIALES Y DEMAS}

Cualquier oficio crea su propio lenguaje para poder proceder con claridad y eficiencia y la física-matemática no es la excepción y en realidad es mucha la terminología oficial que uno utiliza, a lo mejor sin darse cuenta. Para el neófito, eso puede ser un problema. Por eso, hemos insertado algunos comentarios sobre prerrequisitos. Cada comentario  empieza con el símbolo $<<<<<$ y se termina con $>>>>>$. Tiene un título que indica su contenido. Estos comentarios son muy comprimidos y sólo pueden servir para recordar los conceptos, a quienes ya lo saben, o  para que el neófito  tenga claridad sobre qué se requiere investigar en otros textos. El siguiente es el primer comentario.

\color{blue}

\

\

$<<<<<<<<<<<<<<<<<<<<<<<<<<<<<<<<<<<<<<<<<$

\

Algebra lineal y :   grupo, espacio vectorial, función escalar,  base, transformación lineal, función vectorial, diferencial y gradiente.

\
Operaciones sobre conjuntos: se trata de generalizar la suma, la resta, la multiplicación y la división. Nos referimos a un conjunto no vacío $U$.

Se denomina operación interna o cerrada a cada función $f: U \times U \rightarrow U$. Cuando $f$ es una operación interna, se dice de $f$ que opera sobre $U$ y además que es binaria, pues opera sobre pares de elementos. Por ejemplo: La suma es una operación interna cuando opera sobre los naturales $N$, los enteros $Z$, los reales $R$, los complejos $C$. Pero $f(a,b) = (a+b)/7$ no es operación interna sobre $Z$, pues $f(1,1) \not\in Z$. También se dice que $f$ no es cerrada.

Una operación interna $f$ puede cumplir o no cada una de las siguientes propiedades:

- Asociatividad, cuando f puede extenderse a una operación n-aria sin problema alguno. Por ejemplo $f(a,f(b,c)) = f(f(a,b),c)= f(a,b,c)$.

- Existencia de elemento neutro $e$ que no hace nada: $f(a,e) = f(e,a) = a$

- Existencia de todos los  inversos: para cada elemento $a \in U$ existe un elemento $b \in U$ tal que $f(a,b)=f(b,a)=e$. A $b$ se le denomina el inverso de $a$.

- Conmutatividad, cuando el orden no importa: $f(a,b) = f(b,a)$.

Decimos que $(G,f$) es un \index{grupo} \textbf{grupo} si f es una operación que es cerrada, asociativa, tiene elemento neutro y para cada cual existe su inverso. Si además la operación es conmutativa, al grupo se le llama conmutativo a Abeliano. Es usual notar como $(G,+) $ a un grupo conmutativo. Ejemplos de grupos: los reales con la suma $(R,+)$, los complejos con la suma $(C,+)$. El conjunto de los naturales con la suma $(N,+)$ no es un grupo pues no existe en $N$ el inverso de $-1$. El conjunto de los reales con la multiplicación $(R,\times)$ no es grupo pues $0$ no tiene inverso, pero los reales sin el cero, $(R\sp{*}, \times)$ si son grupo. Las traslaciones forman un grupo Abeliano: la operación binaria es la composición.

Debido a que la definición de grupo es tan importante, es conveniente hacer una reformulación de su definición en otra notación que se usa casi siempre exceptuando en lógica, donde se usa la anterior:

\textbf{Grupo} \index{grupo}: es un conjunto $G$ provisto de una suma generalizada. Es decir es un conjunto con una operación binaria (que opera pares de elementos) $\odot$ que cumple:

a. Es cerrada: si $a \in G, $ y $b\in G,$ entonces, $a\odot b \in G$.

b. Es asociativa, o sea que se pueden operar tres elementos sin ninguna ambigüedad:

$a \odot b \odot c = (a \odot b) \odot c= a \odot (b \odot c)$

c. Tiene un elemento neutro $e$ , que como el cero de la suma, no hace nada: $a \odot e = e \odot a= a$.

d. Cada elemento $a$ tiene su inverso $a\sp{-1}$ que restituye el elemento neutro: $a \odot a\sp{-1}= a\sp{-1} \odot  a = e$

Si además la operación  es conmutativa $(a\odot b= b \odot a)$, al grupo se le llama conmutativo o Abeliano.

Si se tiene un grupo $(G,f)$ a una función $g: A\times G \rightarrow G$ se le llama operación externa sobre $U$. Asumamos que $A$ es no vacío y $A$ puede tener estructura de grupo, pensemos en $(R,+)$ $(R\sp{*},\times)$  o en $(C\sp{*},\times)$. La operación externa puede cumplir:

- Idempotencia: $g(1,a) = a$

- Indiferencia: $g(\alpha,g(\beta,a))=g(\alpha\times\beta,a))$

- Distributividad: $g(\alpha,f(a,b))=f((g(\alpha,a),g(\alpha,b))$

A una operación externa $g: E\times G \rightarrow G $ se llama multiplicación escalar si g es idempotente, indiferente y distributiva y si $E$ es $R$ o $C$. También se puede decir que los elementos de $E$ operan sobre $G$  como escalares.

Vectores y \textbf{Espacio Vectorial} \index{Espacio Vectorial}:  'vector' abstrae el concepto de desplazamiento, fuerza, velocidad, los cuales son objetos que se pueden representar naturalmente en un plano o en el espacio tridimensional por medio de flechas que salen del origen. En este caso nos referimos a fuerzas que se aplican sobre un objeto puntual.

Resulta que todas las propiedades geométricas de los vectores dependen de su comportamiento algebraico: los vectores se pueden sumar y restar entre ellos y también alargar, o sea, multiplicar por una constante o escalar, que puede ser un número real o complejo.

Informalmente: un conjunto no vacío que tenga una suma, $+$, y una multiplicación escalar, la cual no tiene notación propia sino sobreentendida, es un espacio vectorial si

1. El conjunto con su suma forma un grupo conmutativo: la suma es cerrada (no produce basura), es asociativa (siendo binaria permite sin ambigüedades la suma de cualquier escogencia de  n elementos), tiene un cero que no hace nada y tiene el concepto de inverso aditivo (con el cual la resta se reduce a la suma del inverso aditivo) y es conmutativa.

2. La operación multiplicación por un escalar es compatible con la suma, cumple la ley distributiva, y satisface los requerimientos de sentido común: si un vector se multiplica por 1 queda igualitico y si se alarga $r$ veces y después $s$ veces, es lo mismo que se alargara $rs$ de una sóla vez.

Formalmente:

Decimos que (V,+,E) es un \textbf{espacio vectorial} sobre E si (V,+) es un grupo y si $E$ es igual a $R$ o a $C$, y si $E$ opera sobre $G$ con multiplicación escalar.

El concepto de conjunto de números lo podemos extender al de conjunto de escalares, que son  objetos que se comportan como números: los reales, los complejos, en general cuando $(V,+,E\sp{*})$ es un espacio vectorial, donde hemos notado al elemento neutro de la suma como $e$ y a $E-\{e\}$ como $E\sp{*}$. No es necesario que la multiplicación sea conmutativa.

\textbf{Función escalar}: A las funciones que toman valores en los reales $R$ o en los complejos $C$, o en general en cualquier conjunto de escalares, se les da el nombre de  funciones escalares, como el peso, la masa, la rapidez, la probabilidad, la amplitud (ver más abajo). Si el dominio de la función es $R\sp{n}$ o $C\sp{n}$ a la función se le denomina campo escalar. Si el dominio de la función es un espacio de funciones, hablamos de un funcional. Ejemplos:

$f(x) = x\sp{2}$

$g(x) = senx$

$h(z) = z\sp{2}$

$F(x,y,z)= x+y-z$

$G(l)=\int\sp{2}\sb{-8} l(x)dx$

Las funciones $f$ y $g$ son funciones escalares reales, $h$ es una función escalar de variable compleja, $F$ es un campo escalar (la variable $z$ insinúa números complejos para $h$, pero significa la tercera coordenada $'z'$ en la definición de $F$) y $G$ es un funcional.

\

\

A una función $T$ de un espacio vectorial $V$ en otro espacio vectorial $W$ que respete la suma y el alargamiento o multiplicación escalar se llama una \index{transformación lineal}  \textbf{transformación lineal}:

$T(\vec u + \vec v) = T(\vec u) + T(\vec v)$

$T(\lambda \vec v) = \lambda T(\vec v)$

donde $\lambda$ es un escalar, y $ \vec u$, $ \vec v$ son vectores.

Una función es transformación lineal si y sólo si su gráfica es una línea recta, o un plano o un hiperplano (plano en n-dimensiones) que pasa por el origen. Ejemplos:

1. $f(x) = x\sp{2}$ no es lineal porque su gráfica no es una línea recta que pase por el origen.

2. $f(x)= 3x$   es transformación lineal de $R$ en $R$ pues su gráfica es una línea recta de pendiente $3$ que pasa por el origen.

3. $F(x,y) = 3x +4y $ si es lineal porque su gráfica es el plano $z= 3x + 4y$ que  pasa por el origen, es decir, el origen (0,0,0) satisface la ecuación del plano: $0=3(0) + 4(0)$. Puede probarse además que el vector $(3,4,-1)$ es perpendicular a todo punto de dicho plano.

4. $F(x,y)=(3x +4y, 8x-7y)$ también es una transformación lineal. Su gráfica es un plano de dos dimensiones en un espacio de 4 dimensiones y  $(x,y,z,w)$ está en dicho plano si $z=3x+4y$ y $w= 8x-3y$.

Las transformaciones lineales pueden representarse en forma matricial. A $F$ le corresponde la matriz siguiente:

$$
\displaylines{
\pmatrix{3&4\cr
    8&-7\cr }\qquad}$$

Esta matriz opera sobre un vector $(x,y)$

$$
\displaylines{
\pmatrix{3&4\cr
    8&-7\cr }\qquad
\pmatrix{x\cr
    y\cr}\qquad} $$

produciendo lo que todos esperamos:

$$
\displaylines{
\pmatrix{3x+4y\cr
    8x-7y\cr }\qquad} $$

Hay una correspondencia 1-1 entre transformaciones lineales y \index{matrices} \textbf{matrices}.

\

\textbf{Base} \index{Base}: es un conjunto de vectores de un espacio vectorial, con el cual se puede reconstruir todo el espacio por medio de combinaciones lineales (alargamientos y sumas), pero de manera única.

Un espacio vectorial tiene infinidad de bases diferentes pero todas tienen el mismo número de elementos. A dicho número  se le llama dimensión. Esta definición de dimensión concuerda con la del sentido común para la línea, dimensión 1, para el plano, dimensión 2, para el espacio, dimensión 3 y para el espacio-tiempo, dimensión 4.

\

Una transformación lineal es conocida completamente si se conoce que hace ella sobre una base.

\

A las funciones que toman valores en un espacio vectorial se denominan \textbf{funciones vectoriales}.  Si el dominio es igual al codominio y es $R\sp{n}$ o $C\sp{n}$,  a la función se le llama \index{campo!vectorial} \textbf{campo vectorial}. Un campo vectorial asocia a cada punto del espacio un vector o una flechita. Ejemplos: el campo de velocidades de un fluido, el campo gravitacional, el campo eléctrico, el campo magnético. Todos ellos toman valores en $R\sp{3}$. El campo electromagnético toma valores en $R\sp{6}$, puesto que a cada punto del espacio le asocia 2 vectores de $R\sp{3}$, el eléctrico y el magnético. Si los campos son variables, hay que agregar una $R$ al codominio, el cual hace las veces de reloj. Un \index{operador} \textbf{operador} es una transformación lineal de un espacio vectorial   en  sí mismo. Muy especiales son los operadores sobre espacios de funciones.   Ejemplos:

$\vec F(\vec r) = \vec F(x\vec i + y\vec j) = \vec F([x,y]) = [x\sp{2}, y-x]$

$G(g)=\int  g(x)dx$

$F$ es un campo vectorial sobre $R\sp{2}$ en tanto que $G$ es un operador. Muchas veces no se menciona el dominio del operador, sino que se asume como aquel conjunto que es el más grande que da sentido a la definición del operador. Tratándose de $G$, el dominio es el conjunto de todas las funciones que van de $R$ en $R$ y que son integrables, es decir cuya integral de Riemann existe como número finito. Dicho conjunto es un espacio vectorial de dimensión infinita. Para ver que su dimensión es infinita piénsese en que cualquier monomio de la forma $y=x\sp{n}$, el cual   no puede expresarse como combinación lineal de los monomios $y=x\sp{m}, m\ne n$.

\

A una función escalar $f$ se le asocia la \index{diferencial} \textbf{diferencial} $df$ y el gradiente $\nabla f$ de la siguiente manera:

$\Delta f$ es el cambio que sufre la función $f$ debido a un cambio $\Delta \vec u$.

$\Delta f =  f(\vec v+\Delta \vec u)- f(\vec v)$

Obervemos que  $\vec v$ y $\Delta\vec u$    no tienen que ser en la misma dirección.

$df$ es el cambio que sufre la función $f$ debido a un cambio infinitesimal en el argumento $\vec v$, el cual se nota $d\vec u$. Puesto que $d\vec u$ denota un cambio, se asume diferente de cero. Puede entenderse a $df$ como el límite de $\Delta f $ cuando $\Delta \vec u \rightarrow \vec 0$:

$df= f(\vec v+d\vec u)-f(\vec v)$

$f(\vec v+d\vec u)= f(\vec v)+ df$

$f(\vec v+\Delta \vec u)= f(\vec v)+ df + error$

Se acostumbra a pensar que $df$ es infinitesimal y que $\Delta f$ es el cambio producido en $f$ por un cambio minúsculo en la variable de entrada. Sin embargo, a veces ni se hace la diferencia.

Por ejemplo, si $f$ denota el volumen de un cilindro con radio $r$ y altura $h$, $df$ resulta ser el volumen de pintura necesario para pintarlo con una capa de espesor $dr $ en la pared lateral  y de espesor $dh$ en las tapas.

Los \index{infinitesimales} \textbf{infinitesimales} aparecen en nuestro discurso como elementos muy subjetivos. Existen dos formalizaciones del concepto de infinitesimal y diferencial, una del análisis no estándar basada en sucesiones infinitas de números y otra en geometría diferencial, la cual es completamente operacional y basada en el concepto de que 'los árboles se conocen por su fruto', es decir, los objetos que se manejan se definen por lo que ellos hacen. Es la misma directiva de la teoría de distribuciones o funciones generalizadas.

\

Busquemos ahora la forma de calcular la diferencial para funciones de varias variables. Bien sabemos que cuando una curva plana es suave, se puede aproximar por su línea tangente con pendiente $m$. En términos oficiales:

Si tenemos una función $f$ de $R$ en $R$ y además   se cumple que existe una función lineal $m$ de $R$ en $R$ tal que:

$f(x+h) = f(x) + m(h) + \epsilon (x,h)$

tal que $\epsilon /(h\sp{2})  \rightarrow K, K\in R$, si $h\rightarrow 0$

decimos que f es diferenciable. Por supuesto que la función lineal $m$ es la línea tangente a la función pero como si el origen fuese el punto de tangencia: $m(h) = f'(x)h $. Obsérvese que no estamos exigiendo solamente que el error $\epsilon $ con que la línea tangente aproxima a la curva $f$ sea despreciable. Exigimos que dicho error sea super-despreciable con respecto a la perturbación en la variable independiente. Eso permite sumar muchos errores, por ejemplo en el teorema del cambio de variable en integración, sin temor de que una suma grande de errores peque\~{n}os se convierta en una amenaza.

Para generalizar al caso de argumentos con varias variables independientes,   en vez de la recta tangente y su pendiente, se tendrá un plano o un hiperplano y varias pendientes, una por cada eje.  Se procede así:

Sea $f$ una función de $R\sp{n}$ en $R$. Considerando $\vec v$ fijo, \index{diferencial} \textbf{la diferencial} $df$ es el cambio en la salida de $f$ debido a un cambio infinitesimal $d \vec u$ en la entrada de $f$:

$df= df(\vec v,d \vec u) = f(\vec v+ d \vec u) - f(\vec v)$

En algunos casos $df$ puede aproximarse por un plano o hiperplano tangente. Concretamente, si existe una función lineal $L$ de $R\sp{n}$ en $R$ tal que:

$f(\vec v+  \vec u) = f( \vec v) + L(\vec u) + \epsilon(\vec u)$

tal que $||\epsilon|| /|| \vec u||\sp{2}\rightarrow K, K \in R, $ si $||\vec u||\rightarrow 0 $ entonces decimos que f es diferenciable. (Hemos escrito el error de tal forma que la definición se pueda  aplicar directamente a funciones de $R\sp{n} \rightarrow R\sp{m}$). La función $L$ tendrá la expresión correspondiente a un hiperplano que pasa por el origen:

$L(\vec u) = L(u\sb{1},.. ..,u\sb{n})= m\sb{1}u\sb{1} +.. ..+  m\sb{n}u\sb{n}$

y cada $m\sb{i}$ será la pendiente del hiperplano tangente en la dirección de $x\sb{i}$

$m\sb{i}= \partial f/\partial x\sb{i}$

\

Al vector de las pendientes, o derivadas parciales, se denomina vector \index{gradiente}  \textbf{gradiente} y se nota $\nabla f$:

$\nabla f = [\partial f/\partial x\sb{1},.. ..,\partial f/\partial x\sb{n}]$

Por lo tanto

$L(\vec u) = L(u\sb{1},.. ..,u\sb{n})= m\sb{1}u\sb{1} +.. ..+  m\sb{n}u\sb{n}= \nabla f \cdot \vec u$

Por ejemplo, si $f(x,y,z) = x - y +xyz$ entonces $\nabla f = [1+yz,-1+xz,xy]$ el cual, evaluado en el punto [1,2,3] da el vector [7,2,2].

El vector gradiente de una función escalar $f$ tiene una dirección y una norma. La dirección apunta hacia el lugar en el cual la función $f$ crece lo más rápido posible. Eso se argumenta así:

$\Delta f= f(\vec v+\Delta \vec u) - f(v) = L(\Delta \vec u) + \epsilon(\Delta \vec u) $

$= \nabla f\cdot \Delta \vec u + \epsilon(\Delta \vec u)= ||\nabla f|| \hspace{0.2 cm} ||\Delta \vec u|| cos\theta + \epsilon(\Delta \vec u)$

donde $\theta $ es el ángulo entre el vector gradiente y el vector $\Delta \vec u$. Para maximizar el valor del cambio de la función con respecto a la dirección tenemos que dejar la norma de $\Delta \vec u$ fija pero que tienda a cero y variar sólo su dirección, es decir el ángulo $\theta$. Eso implica maximizar el valor del coseno, que es uno, y se produce cuando el ángulo $\theta $ es cero, o sea que el vector gradiente debe ser paralelo al vector $\Delta \vec u$. Por eso decimos que el gradiente indica la dirección de máximo crecimiento o máxima derivada direccional.

Si se trata de subir por un monte, los atletas más bravos tomarán la dirección del gradiente (maximización del cambio debido a la dirección) pero si el monte es más inclinado, el esfuerzo será mayor (efecto de la norma del vector gradiente). Similarmente, el calor fluye de donde haya más temperatura hacia donde haya menos. Por eso el calor no fluye en la dirección del gradiente sino en dirección contraria.

Cuando una función $F$ tiene derivadas de orden $k$ y son continuas, decimos que  $f$ es de clase $C\sp{k}$.

\

$>>>>>>>>>>>>>>>>>>>>>>>>>>>>>>>>>>>>>>>>>$

\color{black}

\section{UN PRIMER MODELO}

En una búsqueda de sencillez podemos estudiar el alcance de la  teoría que consta de los siguientes postulados:

1. Las partículas elementales generadoras de campo, denominadas cargas, se dividen en dos clases. Denotémoslas + y -, positivas y negativas.

2.Cargas de igual signo se repelen, cargas de signo contrario se atraen.

3. La fuerza es radial y su magnitud es inversamente proporcional a alguna función decreciente de la distancia.

\

Consideremos ahora el siguiente experimento: tomemos una peinilla de plástico, la frotamos contra el cabello o contra la ropa, y después la acercamos a unos trocitos de papel. Vemos que la peinilla atrae a los papelitos. Cuál es nuestra explicación? Que la peinilla era neutra como un todo, es decir, que las cargas positivas se balanceaban con las negativas, pero por la fricción algunas cargas pasaron de la peinilla al cabello y el balance de cargas se rompió. Quedó pues con un tipo de cargas en exceso. El papel también tiene cargas, pues por ejemplo, al fabricarlo ocurre fricción con la maquinaria que lo produce. Se da la casualidad de que los excesos de carga entre el papel y la peinilla son de signo contrario y es por eso que se atraen.

De paso, conozcamos un problema técnico: las hojas de papel se pegan en las fotocopiadoras y en las impresoras. Se pegan por la deficiencia de carga y por la naturaleza polar del agua: aunque a gran escala las moléculas de agua son neutrales, las cargas que la componen tienen una magnitud y una separación que a corta distancia se crea un campo no nulo. Por lo tanto, uno de los problemas que los fabricantes de papel tienen que resolver es el siguiente: si no se carga eléctricamente al papel, éste no se repele y se pegará en la fotocopiadora, pero si se carga mucho al papel, pues se repele, pero puede ocurrir una chispa que queme todo. O bien, puede atraer demasiada humedad y bajar la calidad del papel. A los fabricantes de pegantes les toca resolver el problema contrario: aumentar la capacidad de atracción,  lo cual es más fácil para pegantes líquidos que para cintas adhesivas que tienen que ser hidrofóbicas.

\

Estas ideas de \index{electricidad} \textbf{electricidad} nacida de frotación fueron las que le permitieron a Willian Gilbert en 1600 la introducción del término '\textbf{elektron}' \index{elektron} que en griego significa ámbar para relacionar los fenómenos que nosotros denominamos de electrostática. El ámbar y sus propiedades eléctricas que aparecían por frotación ya fueron descritas por los griegos. Sin embargo, la asociación  de electrón  con una partícula elemental portadora de electricidad fue un invento de Hendrik Antoon Lorentz en 1892. Pero el electrón como realidad experimental oficialmente vió la luz como resultado de investigaciones con tubos de alto vacío. Eso sucedió al puro final del siglo XIX, gracias a sir Joseph Thompson y su equipo. La medición de su carga fue hecha por Robert Millikan en 1909 midiendo la desviación que sufre una gota de aceite al ir cayendo en la presencia de un campo eléctrico.

\

Bien, aunque nos parece que tenemos un cierto éxito con nuestra teoría atómica, debemos aclarar que nuestra teoría cuenta con un grave desperfecto: si todas las cargas fuesen de un mismo signo, se repelerían hacia el infinito. Y si se balancearan en cuanto a signos, entonces, qué impide que colapsen y que todo sea irremediablemente neutral? Así que nuestra teoría no puede explicar los cuerpos neutrales que no se dejan comprimir. Tal vez podamos hacerlo con cargas en movimiento, pero eso deberá ser analizado más tarde, cuando sepamos la dependencia exacta de la fuerza con la distancia. Mientras tanto, tenemos que añadir un nuevo postulado:

\

4. Existen cargas neutrales.

\

Este postulado es muy vergonzoso por la siguiente razón: estamos tratando de dar cuerpo a la teoría atómica. Pero ser neutro significa que ni atrae ni repele. Por lo tanto, cómo es posible que exista un cuerpo neutro que no se deje comprimir y sea formado de átomos? Como si fuera poco, los cuerpos neutros además de no tener carga, deben tener propiedades ad hoc que nos permitan dar la siguiente explicación a los diferentes niveles de atracción o repulsión que se observan en la naturaleza:

Consideremos un dipolo, es decir, dos cargas de signo contrario que se fijan a una cierta distancia. Por supuesto que  las cargas no se pueden fijar en el vacío, pero podrán fijarse  sobre un substrato neutro y que impida el movimiento de las cargas. Decimos que  el substrato es mal conductor. De qué forma un cuerpo neutro puede sujetar a una carga eléctrica sin atraerla ni repelerla es algo que debemos no preguntarnos. Estudiemos ahora el campo eléctrico creado por un dipolo. Gráficamente podemos ver que el dipolo presenta un campo fuerte en su cercanía pero débil en su lejanía.

Para fijar ideas, imaginemos que nuestro dipolo tiene una carga +1 a la izquierda y -1 a la derecha. El vector $\vec F$ representa la fuerza de repulsión de una carga de prueba +1 colocada en un punto P del plano. El vector $\vec G$ representa la fuerza de atracción de la misma carga de prueba en el mismo punto por la carga -1. Qué es $\vec F+\vec G$? Es la diagonal  del paralelogramo generado por $\vec F$ y $\vec G$ y que parte del mismo vértice, notada $\vec d$. (Entre tanto que $\vec F- \vec G$ es la otra diagonal, notada $\vec D$).

Para distancias grandes comparadas con la distancia entre las cargas del dipolo tenemos que $\vec G$ es prácticamente $-\vec F$, por lo tanto tenemos las aproximaciones $\vec d=\vec F+ \vec G=\vec F+(- \vec F)=\vec 0$ entre tanto que $\vec D=\vec F-\vec G = \vec F-(-\vec F) = 2\vec F$. Inferimos que a lo lejos, el campo generado por un dipolo, la resultante $\vec F+\vec G=\vec d$, es de bastante menor  magnitud que el campo generado por una carga libre y única.

Un \textbf{dipolo} \index{dipolo} produce un campo muy fuerte en su cercanía. De hecho en el centro del dipolo, el campo resultante es el doble del campo generado por una sola carga, pues $ \vec F = \vec G$ y entonces $\vec d=\vec F+ \vec G=\vec F+\vec F=2\vec F$ . Imaginemos ahora el problema de transportar iones a través de la membrana celular. Ella está rodeada de iones a lado y lado y sin embargo a toda hora hay transporte. No en vano, las células transportan los iones por pares para minimizar la energía gastada. Puede ser que un ión salga y otro entre o bien que dos iones de signo contrario entren o salgan.

\

Hemos podido demostrar que con cargas eléctricas de una misma magnitud podemos producir objetos que a la distancia se ven como neutros, pero que a mediana distancia producen campos de una magnitud menor que la de campos generados por una sola carga. Campos de variada magnitud se pueden construir regulando la distancia entre las cargas del dipolo o considerando varios dipolos juntos, o desbalanceando el dipolo. En principio, entonces, podemos explicar los variados rangos de atracción y de repulsión asumiendo que hay cargas fundamentales todas con el mismo valor. Por supuesto, siempre y cuando se puedan formar estructuras neutras.

\section{LA CUANTIZACION DE LA CARGA}

Nuestro estudio del dipolo  nos ha librado de introducir el postulado adicional de que la magnitud de las cargas se mueve en un variadísimo rango. Hemos removido dicho postulado demostrando que con juegos espaciales de partículas con la misma carga se pueden producir campos electrostáticos de variada magnitud. Por lo tanto, inspirados en el postulado de que la naturaleza debe ser sencilla, declaramos que las cargas vienen con signo positivo y negativo y que todas las cargas fundamentales vienen con  la misma magnitud.

Si Usted quisiera ponerle un poco de experimentación a nuestro estudio, explore las posibilidades que ofrece el caldo de gallina caliente. Resulta que en la superficie del caldo aparecen nadando unos botoncitos de grasa. Algunos botoncitos pueden fusionarse con otros, pero no todos: hay botoncitos que se repelen, porque están cargados con cargas de igual signo. Acaso encontrará Usted la forma de asegurarse que la fuerza de repulsión está cuantizada?

\

Aunque hemos producido una sencillez conceptual muy poderosa, no deja de ser preocupante el pensar que tal vez la naturaleza no sea tan simplista como nosotros. En efecto, aunque la cuantización dicotómica de la carga es una verdad experimental irreprochable, se encuentra un desafío en las modernas teorías gauge: hay elementos, llamados quarks, que no pueden existir en estado aislado y cuya carga sería de 1/3 o 2/3 la carga del electrón.

\section{LA LEY DE COULOMB}

Nos parece muy agradable reconsiderar  un experimento llevado a cabo por Charles Augustin de Coulomb en 1786 para determinar la forma exacta de la dependencia  de la fuerza con la distancia. Para verificar y extender un estudio de  Joseph Priestley ejecutado en 1766 que anunciaba que la dependencia debería de ser cuadrática, Coulomb utilizó una balanza de torsión.

Ese es un instrumento muy sencillo: se trata de un hilo que pende del techo, al cual le colgamos una varilla mala conductora en posición horizontal y que forma un dipolo. Ponemos una carga de prueba en la cercanía y observamos el comportamiento de la varilla horizontal: según la ley de los signos, el hilo girará hasta que la atracción o repulsión que viene desde afuera  quede balanceada por la fuerza elástica del hilo.

Lo que para un experimentalista suene como una dependencia 'cuadrática', para un matemático o físico idealista, eso no puede ser más exacto que aproximadamente cuadrática. Pero el matemático no puede quedarse ahí: debe encontrar la manera de dilucidar el marasma que el mismo creó. Para lograrlo, necesitamos de los avances de la tecnología.

\

\

\color{blue}

$<<<<<<<<<<<<<<<<<<<<<<<<<<<<<<<<<<<<<<<<<$

\

Maquinaria pesada: \index{divergencia} \textbf{divergencia} de un campo, \index{teorema de Gauss}  \textbf{teorema de Gauss}, significado de la divergencia.

\
Al hablar de componentes en el espacio nos referimos con $\vec i$, $\vec j$, $\vec k$ las 3 direcciones fundamentales.

La \textbf{divergencia} de un campo con componentes $(f,g,h)= f\vec i + g \vec j + h \vec k $ es

$$\nabla  \cdot (f,g,h) =\partial (f)/\partial x  + \partial (g)/\partial y + \partial (h)/\partial z
 \addtocounter{ecu}{1}   \hspace{5cm} (\theecu )      $$

Por ejemplo, si el campo es:

$\vec F = (x,y-z,cosx)$

entonces su divergencia es:

$\nabla  \cdot (x,y-z,cosx) =\partial (x)/\partial x  + \partial (y-z)/\partial y + \partial (cosx)/\partial z $

$=1+1=2$, el cual es un escalar, un número.

\

En el \textbf{Teorema de Gauss} tenemos que considerar un volumen $V$, circundado por su frontera  $\partial V$ que es una superficie suave y orientable (a la cual se le puede marcar sin ambigüedad cual es el exterior y cual el interior),  y con vector  unitario normal hacia afuera  $\vec n$, un campo vectorial  $\vec E $ cuya divergencia es $\nabla  \cdot \vec E $. En ese caso el teorema nos dice que:

$$ \int\sb{V} \nabla  \cdot \vec E
=\int\sb{\partial V} \vec E
=\int\sb{\partial V}\vec E \cdot \vec n dS
\addtocounter{ecu}{1}   \hspace{7cm} (\theecu ) $$

No hay que dejarse asustar de esta expresión tan misteriosa. Si tenemos en cuenta que la divergencia es un escalar, digamos que se trata de una densidad, podemos decir que

$ \int\sb{V} \nabla  \cdot \vec E $

corresponde a la integración de dicho escalar o densidad por toda el volumen $ V$. En suma, se trata de un cierto tipo de carga total esparcida por el volumen $V$.   Por otra parte, la expresión:

$\int\sb{\partial V}\vec E \cdot \vec n dS $

puede interpretarse imaginando que $\vec E$ representa el campo de velocidades de un fluido que atraviesa a $V$. En ese caso, la integral de dicho campo vectorial sobre la frontera de $V$ es el flujo total que sale de $V$. Entonces el teorema de Gauss es un simple juego de contabilidad: el flujo total que sale de $V$ es el resultado de lo que le pase al campo en el interior de $V$. Busquemos la forma de ser un poco más explícitos.

\

Como  hemos dicho, interpretamos el campo como la velocidad de un fluido, entonces el flujo lo definimos como la cantidad de fluido que atraviesa la superficie por unidad de tiempo. El flujo se calcula por una integral de superficie. Al dividir el flujo que atraviesa una superficie cerrada entre el volumen encerrado, tenemos flujo que sale de la superficie por unidad de volumen. Al tomar el límite cuando el área de la superficie tiende a cero, tenemos la divergencia del campo. Para ver esto en  más detalle, apliquemos el Teorema de Gauss a una esfera para después tomar el límite del radio hacia cero:

Sea $S$ una esfera centrada en $P$, considerada con su volúmen interior y sea $\partial S$ su frontera. Tomando como  $\vec n$  el vector normal exterior unitario a la superficie, el cual es paralelo al radiovector de la esfera, el teorema de Gauss no dice:

$ \int\sb{S} \nabla  \cdot \vec E
=\int\sb{\partial S} \vec E
=\int\sb{\partial S}\vec E \cdot \vec n dS  $

Pero si hacemos que el radio tienda a cero, entonces considerando que la divergencia es una función continua podemos decir que sobre una esfera microscópica la divergencia es constante e igual a su valor en el centro o sea $\nabla  \cdot \vec E (P)$ y por lo tanto sale de la integral:

$ \int\sb{S} \nabla  \cdot \vec E
=\int\sb{S} \nabla  \cdot \vec E dV= \nabla  \cdot \vec E (P) \int\sb{S} dV$

Tomando el otro lado del teorema, tenemos:

$\nabla  \cdot \vec E (P) \int\sb{S} dV= \int\sb{\partial S} \vec E $

Despejando:

$\nabla  \cdot \vec E (P) =
(\int\sb{\partial S} \vec E )/\int\sb{S} dV$

\

En resumen: la \index{divergencia} \textbf{divergencia} es la tasa (microscópica) de flujo que sale de una superficie por unidad de volumen.

\

Podemos ahora parafrasear el Teorema de Gauss muy sencillamente: si se tiene una superficie cerrada la cual es atravesada por un fluido, lo que salga de su superficie es simplemente el resultado de lo que surja o se expanda en su interior menos lo que se consuma o comprima. En particular, en todo lugar de un campo eléctrico en el cual haya una carga, tenemos una fuente si la carga es positiva y un sumidero si es negativa. En ambos casos la divergencia no puede ser nula.

Además de fuentes y sumideros, el fluido se puede comprimir o expandir haciendo que su divergencia tampoco sea cero. Por ejemplo, un campo que denota un fluido que marcha en la dirección $Y$ y que se va comprimiendo es el siguiente:

$\vec F (x,y,z)= (0,1/y,0)$

su divergencia es:

$\nabla \cdot \vec F = -1/y\sp{2}$

la cual no es nula en ningún lado.

\

$>>>>>>>>>>>>>>>>>>>>>>>>>>>>>>>>>>>>>>>>>$

\

\

\color{black}

El siguiente teorema nos ayudará mucho para pasar de una dependencia \index{aproximadamente cuadrática}  \textbf{aproximadamente cuadrática} a cuadrática.

\

\

\textit{Teorema \theecu: Si un campo radial tiene una dependencia aproximadamente cuadrática inversa entonces: la divergencia del campo es cero fuera del origen si y sólo si tal dependencia es exactamente cuadrática inversa.}

Quizás sirva aclarar que tener divergencia cero no implica ser radial ni tener una dependencia cuadrática inversa. Por un ejemplo, un campo constante en el espacio tiene divergencia cero y sin embargo no tiene dependencia cuadrática.

\

Demostración:

Una  dependencia cuadrática pura tendría la forma

$\vec E(\vec r)=k \|\vec r\|\sp{-2} \vec r\sb{u}   $

donde $\vec r\sb{u}$ es el vector unitario en la dirección radial. Dicho vector es $\vec r\sb{u}= \vec r/\|\vec r|| $, por lo tanto, la dependencia cuadrática pura toma la forma:

 $\vec E(\vec r)=k \|\vec r\|\sp{-2} \vec r/\|\vec r\|   =k \|\vec r\|\sp{-3} \vec r\  $

Así que una expresión  aproximadamente cuadrática para un campo radial , para $\vec r\not=\vec 0$ sería:

$\vec E(\vec r)=k \|\vec r\|\sp{-3-p} \vec r   $
$=k \|\vec r\|\sp{-3-p} (x,y,z)$

Como estamos estudiando la desviación con respecto a un cuadrado exacto, permitimos una ligera perturbación: ese es el $p$ (que en principio puede ser cualquier valor). Al repartir queda como:

$$\vec E(\vec r)=
k \|\vec r\|\sp{-3-p} x \vec i
+ k \|\vec r\|\sp{-3-p} y \vec j +
k \|\vec r\|\sp{-3-p} z \vec k  \hspace{2cm}(\theecu)  $$

Calculemos la divergencia de este  campo:

$\vec E(\vec r)= k \|\vec r\|\sp{-3-p} x \vec i $
$+ k \|\vec r\|\sp{-3-p} y \vec j $

$+ k \|\vec r\|\sp{-3-p} z \vec k$

Puesto que $\|\vec r\|= (x \sp{2} + y \sp{2} + z \sp{2}) \sp{1/2}$, la divergencia es entonces:

$\nabla \ \cdot \vec E $
$=\partial /\partial x (x\|\vec r\|\sp{-3-p}) $
$+ \partial /\partial y (y\|\vec r\|\sp{-3-p})$
$+ \partial /\partial z (z\|\vec r\|\sp{-3-p})$

$=\|\vec r\|\sp{-3-p} $
$+ x(-3-p)\|\vec r\|\sp{-4-p}(2x)/(2\|\vec r\|)$

$+\|\vec r\|\sp{-3-p} $
$+ y(-3-p)\|\vec r\|\sp{-4-p}(2y)/(2\|\vec r\|)$

$+\|\vec r\|\sp{-3-p} $
$+ z(-3-p)\|\vec r\|\sp{-4-p}(2z)/(2\|\vec r\|)$

$=3\|\vec r\|\sp{-3-p} $
$+ \|\vec r\|\sp{2}(-3-p)\|\vec r\|\sp{-4-p}/\|\vec r\|$

Podemos ir concluyendo:

$$\nabla \ \cdot \vec E =3\|\vec r\|\sp{-3-p}
+ (-3-p)\|\vec r\|\sp{-3-p}
\addtocounter{ecu}{1}   \hspace{4cm} (\theecu ) $$

$$\nabla \ \cdot \vec E=-p\|\vec r\|\sp{-3-p} \addtocounter{ecu}{1}   \hspace{7cm} (\theecu ) $$

Vemos entonces por el resultado (\theecu ) que la divergencia discrimina entre  una dependencia aproximadamente cuadrática y cuadrática exacta: en un lugar no ocupado con la carga la divergencia de un campo aproximadamente cuadrático inverso es cero si y sólo si la dependencia es cuadrática inversa exacta, es decir si $p=0$.

\

Tenemos que inventarnos la forma de \index{medir exactamente} \textbf{medir exactamente} la divergencia del campo eléctrico. Un diseño inmortal, debido a Franklin, el mismo del pararrayos, y quien se basó en lo que se sabía del campo gravitatorio, es el siguiente:

\

Al cargar una esfera hueca conductora, puede medirse el campo eléctrico en cualquier punto de su interior y el resultado es cero. Eso sólo es posible, y lo demostraremos, cuando la fuerza siendo aproximadamente cuadrática tiene una dependencia cuadrática exacta, o sea cuando su divergencia es cero.

En concreto, el experimento es así: se tienen dos esferas concéntricas, conductoras y conectadas por un alambre conductor. La esfera exterior está compuesta de dos hemisferios. Se carga la esfera exterior. Se desconecta la esfera interior. La esfera exterior se separa en sus dos mitades, y se investiga la esfera interior para ver si tiene carga o no.

\

El resultado experimental es negativo: la esfera interior no se cargó porque la esfera exterior no produjo ningún campo neto en su interior.

\

Veamos cómo podemos demostrar que eso es posible  en el único caso en el cual la dependencia aproximadamente cuadrática  sea cuadrática exacta. Consideremos pues las dos esferas concéntricas, la exterior y la interior.

El campo que existe en un punto determinado de la esfera interior es la suma de todos los campos creados por todas las cargas de la esfera exterior. Con el teorema de Gauss  podemos probar que el campo en cualquier punto de la esfera interior debe ser cero para el caso cuadrático y no cero para cualquier otro caso aproximadamente cuadrático. En efecto: sea $S$ la esfera interior, considerada con su volumen interior y sea $\partial S$ su frontera, la esfera propiamente dicha a la cual nos estamos refiriendo por defecto. Tomando como  $\vec n$  el vector normal exterior unitario a la superficie, el cual es paralelo al radiovector de la esfera, el teorema de Gauss nos dice:

$$ \int\sb{S} \nabla  \cdot \vec E
=\int\sb{\partial S} \vec E
=\int\sb{\partial S}\vec E \cdot \vec n dS
\addtocounter{ecu}{1}   \hspace{7cm} (\theecu ) $$

Ahora bien, la simetría en la geometría del problema nos indica que cada punto en la superficie de la esfera interior es indistinguible de cualquier otro. Por tanto, el vector $\vec E$ no debe preferir ninguna dirección ni magnitud sobre dicha superficie. Deducimos que  debe ser un múltiplo constante del radiovector de la esfera interior y entonces el ángulo entre los dos debe ser $0$. Tenemos:

$\int\sb{\partial S}\vec E \cdot \vec n dS $
$=\int\sb{\partial S}\|\vec E\| \|\vec n\| cos(0) dS =$
$\int\sb{\partial S}\|\vec E\| dS $
$= \|\vec E\| \int \sb{\partial S} dS $
$= \|\vec E\|(4\pi)\|\vec r\|\sp{2}$

En conclusión:

$ \int\sb{S} \nabla \ \cdot \vec E $
$= \|\vec E\|(4\pi)\|\vec r\|\sp{2}$

Por tanto,  puesto que el campo medido experimentalmente es cero, se concluye que

$ \int\sb{S} \nabla \ \cdot \vec E =0$

Ahora bien, puesto que se trata de despejar la divergencia del campo, hay que encontrar la forma de sacarla de la integral: eso puede hacerse cuando el radio es casi cero. Para un radio así la divergencia es constante (pues es continua) y sale de la integral:

$ \int\sb{S} \nabla \ \cdot \vec E = \nabla \ \cdot \vec E \int\sb{S}dS =0$.

\

De donde deducimos que

$\nabla \ \cdot \vec E = 0$

que era lo que tanto habíamos buscado, demostrando de esa forma que es justo decir que la intensidad del campo eléctrico creado por cargas tiene una dependencia cuadrática inversa exacta.

\

El mérito del test recién estudiado es que es muy sensible pues según nuestra igualdad (5) es suficiente disminuir el radio de las esferas para hacer que la divergencia crezca, si es que es no nula.

Los resultados experimentales, para 1971, dieron para $p$ un valor menor a $10\sp{-15}$. Por lo tanto, es costumbre asumir sin aclaración ninguna que el campo decae cuadráticamente con la distancia con una dependencia cuadrática exacta (lo cual ya no volverá a enfatizarse).

\

Nosotros hemos definido el campo eléctrico como una abstracción para explicar la repulsión de los átomos. Hemos caracterizado el campo eléctrico como un campo que decae cuadráticamente. Por supuesto que tal decaimiento se estableció experimentalmente para distancias del orden de un metro. Es interesante saber si la ley se cumple en el resto del cosmos.

El campo magnético de Júpiter, estudiado desde el punto de vista cuántico en relación con el campo eléctrico, ha permitido decir que por esos lados también se cumple la misma ley. En general, no hay mucho interés entre los  astrónomos en declararle la guerra a la ley cuadrática.

Y qué pasa a la \index{escala microscópica} \textbf{escala microscópica}?

Los experimentos con choques de partículas cargadas muy aceleradas demuestran que las leyes de atracción y repulsión cuadrática se rompen a escalas del núcleo atómico. Eso ha permitido descubrir y estudiar la naturaleza de otras fuerzas a las cuales se les da el nombre de nucleares. Acerca de dichas leyes, las teorías gauge se han mostrado muy elocuentes.

En cuanto a la forma cuadrática inversa de la gravitación, podemos decir lo siguiente: los experimentos que tenemos hasta hoy para la ley cuadrática inversa de la caída de la fuerza gravitatoria se ha verificado sólo hasta el décimo de milímetro. No se ha podido menos por la incapacidad de silenciar las fluctuaciones eléctricas que son muchísimo más fuertes. Estos experimentos han recibido renovado interés pues la teoría de cuerdas postula que el espacio es de más dimensiones que 3 y que las constantes y fuerzas fundamentales dependen del número de dimensiones del espacio y de su topología, es decir, si las dimensiones son cerradas, como en un anillo, o si son planas como en una mesa.

\bigskip

La ley cuadrática nos permite ahora matizar mejor nuestro problema de explicar la existencia de los cuerpos neutros. Como cargas de signo contrario se atraen, la ley cuadrática predice que entre menor sea la distancia, la fuerza de atracción es mayor y tiende a infinito cuando la distancia tiende a cero. Por lo tanto, todas las cargas deberían neutralizarse mutuamente. Resulta entonces un verdadero misterio entender por qué existen cuerpos neutros que se pueden cargar por fricción. Probablemente existan muchas maneras como uno podría tratar de resolver este enigma, pero una muy sencilla y atractiva es proponer que las cargas se estabilizan entre ellas por un efecto de campo global.

Lamentablemente nuestra ilusión del campo global estabilizante carece de fundamento: es imposible construir un campo eléctrico fijo (se denomina campo electrostático) que pueda mantener en equilibrio estable a una carga eléctrica en el espacio vacío. Para entender la razón de dicha proposición necesitamos  algunos resultados matemáticos.

\

\

\color{blue}

$<<<<<<<<<<<<<<<<<<<<<<<<<<<<<<<<<<<<<<<<<$

\

La función \index{potencial} \textbf{potencial}.

\

Vimos anteriormente que a una función escalar diferenciable se le puede asociar un campo vectorial, el campo gradiente, que en cada punto produce un vector que indica la dirección en la cual la derivada direccional es máxima. Para fijar ideas, pensemos en una función escalar $f: R\sp{2} \rightarrow R$. Tal función la podemos representar como una monta\~{n}a y podemos visualizar cómo cae la lluvia y su forma de deslizarse cuesta abajo. Cuando la  monta\~{n}a es muy decente, la lluvia rueda en la dirección contraria al gradiente. El gradiente es un campo en el piso, no es un campo que vaya encima de la monta\~{n}a.

\

El problema del potencial es el inverso: dado un campo vectorial en el piso, hay que construir una monta\~{n}a tal que su gradiente reproduzca el campo vectorial dado. Debe ser una monta\~{n}a tal que cuando caiga la lluvia, ruede en dirección contraria al campo dado, y con la misma magnitud. Si tal monta\~{n}a existe, a la función que la genera se llama la función potencial.

Por ejemplo, si el campo es

$\vec F = (3-yz,4-xz,-xz) $

entonces su función potencial es:

$\phi= 3x+4y-xyz + k $

La función potencial no puede averiguarse exactamente: siempre hay una constante indeterminada,  o más bien, hay una constante que uno puede fijar a gusto personal. Por eso, si en alguna teoría con contenido físico aparece un potencial, los resultados observables no deben depender de dicha constante. Se tiene una teoría gauge con grupo de invariancia $R$.

\bigskip

No siempre existe la función potencial y, por tanto, no siempre puede hablarse de grupo de arbitrariedades. Eso nos hace pensar que no hay ninguna obligación por parte de la naturaleza en que sus interacciones sean representadas por una teoría gauge. Y por tanto, todas nuestras pretensiones han de demostrarse.

 \bigskip

 Como ejemplo de un campo que no admite un potencial, consideremos  $\vec F= (3y,4x,0)$ no tiene potencial. En efecto, si lo tuviera entonces sería $\phi$ y se cumpliría que:

$\vec F = \nabla \phi = (\partial \phi/\partial x,\partial \phi/\partial y,\partial \phi/\partial z) = (3y,4x,0) $

Por tanto:

$\partial \phi/\partial x= 3y$

$\partial \phi/\partial y= 4x$

De la primera ecuación se tiene después de integrar

$\phi = 3yx + G(y,z)$

Por consiguiente

$ \partial \phi/\partial y = 3x + \partial G(y,z)/\partial y$

Pero también debe ser igual a:

$ \partial \phi/\partial y=4x$

Por lo que resulta que 3=4. Contradicción.

\

$>>>>>>>>>>>>>>>>>>>>>>>>>>>>>>>>>>>>>>>>>$

\

\

\color{black}

Intriga: Un caso ejemplo concreto de una fuerza que no viene de un potencial es el de la fricción. Este tipo de fuerza es disipativo, que por demás produce calor, el cual puede ser indeseable o puede ser puesto a hacer algo útil, como en el caso que sigue. Existe un método para soldar metales que es como sigue: las piezas  se liman bien por las junturas a ser soldadas. Se ponen a girar en dirección contraria y se hace que friccionen. Se espera a que se eleve la temperatura lo suficiente, digamos hasta cuando los metales brillen al rojo, y entonces se para de girar, pero lo más rápido posible y luego se deja el conjunto quietecito. Al enfriarse el material, ya está soldado.

Ahora bien, la fricción viene de la interacción entre manojos de átomos y moléculas. Como éstos no son fundamentales, no hay ningún problema en suponer que la fricción no es una interacción fundamental. De donde sale dicha interacción? Nuestra única opción es suponer que la fricción puede explicarse por la interacción electromagnética entre sus constituyentes básicos. Pero si eso fuese así, cómo podemos explicar que de una interacción fundamental conservativa, como lo es el electromagnetismo entre partículas fundamentales, salga algo no conservativo, que disipa energía, como es la fricción?

Y como si fuera poco, sin fricción sería imposible caminar. Bueno, volvamos ahora a nuestra discusión.

\

Como la fuerza eléctrica tiene la misma forma que la gravitatoria, pues también goza de la propiedad de tener una función potencial (la razón matemática consiste en tener 'rotacional cero', expresión cuyo significado se explicará más luego): existe una función escalar $\phi$ tal que su gradiente es igual al campo eléctrico:

$-\nabla \cdot \phi = E$

El signo menos se debe a que el campo eléctrico es de repulsión, entre cargas iguales, y apunta hacia donde el potencial decrece, mientras que el gradiente apunta hacia donde el potencial crece. Con gravitación sucede como sigue: el gradiente señala la dirección por la cual suben los grandes atletas, pero el campo da la dirección por donde resbala el agua, la cual es movida por el campo gravitatorio, el cual apunta en dirección contraria al radiovector.

\

Ya hemos visto  que en el espacio vacío la divergencia del campo eléctrico es $0$. En donde hay carga, revisando un cálculo anterior (4), la divergencia de un campo radial con dependencia cuadrática inversa da $ \infty- \infty$ que es indeterminado. Aprovechamos esa indeterminación para darle un valor que nos convenga. Teniendo en cuenta que la divergencia es el flujo por unidad de volumen, consideremos  una esfera de radio $r$, por lo que  el campo eléctrico en cualquier parte de la superficie de la esfera tiene norma $Q/r \sp{2}$ y es paralelo al vector normal. Por consiguiente, el flujo es $(Q/r \sp{2}) \times 4\pi r \sp{2}= 4 \pi Q$ el cual es independiente del radio de la esfera. Usamos esa independencia para extender la validez de  nuestra discusión al caso $r=0$.  Podemos decir:

$$\nabla \cdot E = 4\pi Q
\addtocounter{ecu}{1}   \hspace{4cm} (\theecu )      $$

Si la carga es cero, entonces la divergencia es cero. Por eso la ecuación anterior es perfectamente válida tanto para el espacio vacío como para  cuando hay una carga en un punto determinado. Recalcamos que esta ecuación es equivalente a la ley de Coulomb de decaimiento cuadrático. Por eso se considera que (\theecu) es una ecuación fundamental y la recordaremos como \textbf{la tercera ley de Maxwell}. Formalicemos:

\

\

\color{red}

\addtocounter{ecu}{1}
\textit{Teorema \theecu: Si un campo eléctrico tiene una dependencia inversamente proporcional al cuadrado de la distancia (ley de Coulomb) entonces también cumple la tercera ley de Maxwell: }

$$\nabla \cdot E = 4\pi Q
\addtocounter{ecu}{1}   \hspace{4cm} (\theecu )      $$

\color{black}

\

\

Como el campo es el gradiente de la función potencial $\phi$, resulta que

$$\nabla \cdot E = \nabla \cdot \nabla \phi
=\nabla \sp{2} \phi = 4\pi Q
\addtocounter{ecu}{1}   \hspace{4cm} (\theecu )      $$

A esta ecuación (\theecu ) se le llama la ecuación de Poisson.

En una región que no encierra ninguna carga, la divergencia es cero, en términos de potencial se reescribe:

$$\nabla \sp{2} \phi = 0
\addtocounter{ecu}{1}   \hspace{4cm} (\theecu )      $$

La ecuación resultante  (\theecu ) se denomina la ecuación de Laplace.

\

Permitámonos explicitar nuestra exigencia básica gauge de realismo acerca de nuestros formalismos: las leyes matemáticas deben contener en todos los sistemas de coordenadas la misma ley física. Nuestra ecuación anterior puede leerse de dos maneras: como la expresión de una ley física o como su particularización en un sistema de coordenadas específico.

Veamos: la expresión $\nabla \cdot E $ representa la divergencia de un campo vectorial. Es una cantidad física que no se refiere a ningún sistema de coordenadas en particular. Pero si introducimos coordenadas específicas ya hemos pasado a la particularización de la ley al sistema de coordenados tomado. Por ejemplo, si las coordenadas son cartesianas, tenemos que

$\nabla \ \cdot E $
$=\nabla \ \cdot (E \sb{x},E \sb{y},E \sb{z} ) $
$=\partial E \sb{x}/ \partial x  $
$+ \partial E \sb{y}/ \partial y  $
$+ \partial E \sb{z}/ \partial z  $

La ecuación de Poisson representa una ley física si entendemos $\nabla \sp{2} \phi$ como la divergencia del campo gradiente. Tanto la divergencia de un campo vectorial como el gradiente de una función escalar son entidades que existen independientemente de todo sistema de coordenadas. Pero tan pronto pensemos en derivadas parciales con respecto a $x$ y demás, ya estamos en un sistema de coordenadas y nuestros resultados  serán válidos para dicho sistema. A una formulación que no dependa de ningún sistema de coordenadas se la llama formulación intrínseca. Si depende de un sistema de coordenadas se llama formulación en coordenadas cilíndricas o lo que corresponda.

En el espacio vacío, la función potencial cumple la \index{ecuación de Laplace} \textbf{ecuación de Laplace}:

$\nabla \sp{2} \phi = 0$

\

\

\color{blue}

$<<<<<<<<<<<<<<<<<<<<<<<<<<<<<<<<<<<<<<<<<$

\

\textbf{Funciones armónicas} \index{Funciones armónicas}:

\

A una función $\phi$ que cumpla la ecuación de Laplace $\nabla \sp{2} \phi = 0$ se le llama \textbf{función armónica}. Dichas funciones, que son escalares reales, tienen las propiedades siguientes:

\

1. Una función armónica no puede tener extremos, máximos o mínimos, en una región abierta libre de carga.

\

2. Si una función es armónica en una región libre de carga con frontera, entonces tiene sus extremos sobre la frontera.

\

3. Si la función $\phi (x,y,z)$ satisface la ecuación de Laplace, entonces el valor promedio de $\phi $ sobre la superficie de cualquier esfera es igual al valor de $\phi $ en el centro de la esfera. Es lo mismo que decir que un cuerpo esférico atrae a otro puntual como si toda la masa del cuerpo esférico estuviese concentrada en su centro.

\

Demostremos la primera propiedad por contradicción: supongamos que una función armónica $\phi$ tenga un extremo en P, que vamos a suponer que sea un mínimo. Puesto que se trata de un extremo local, las primeras derivadas parciales deben ser cero y cada una de sus segundas derivadas parciales, en cualquier sistema de coordenadas, debe ser positivo. La ecuación de Laplace contiene del lado izquierdo  la suma de las segundas derivadas puras, $\partial \sp{2}\phi / \partial x \sp{2} $,$\partial \sp{2}\phi / \partial y \sp{2} $, $\partial \sp{2}\phi / \partial z \sp{2} $ y como cada una es positiva, su suma también lo es, dejando de cumplir la ecuación de Laplace.

Contradicción.

\

Prueba de la segunda propiedad. Un dominio es compacto cuando todo recubrimiento abierto tiene un subrecubrimiento finito. Eso es equivalente a decir que toda sucesión tiene un subsucesión convergente a un punto dentro del conjunto. También es lo mismo que decir que el conjunto es cerrado, con su frontera, y acotado, que se puede amurallar por una bola de radio finito. Si la región de definición de la función armónica $\phi $ tiene frontera, entonces todo el conjunto es compacto. Como $\phi $ es continua y su dominio es compacto, entonces  $\phi $ tiene que tomar sus valores extremos, existirán puntos P y Q que sean mínimo y máximo. Como ellos no pueden estar sobre el interior del dominio de $\phi $, entonces deben estar sobre el borde.

\

Aclaremos eso de que una función continua definida sobre un compacto toma sus valores extremos. Primero, pensemos que para la función $f(x) = 1/x$ no existe ningún punto tal que en ese punto la función tome el valor máximo. En este caso, la función es continua sobre el abierto (0,1) pero no es continua sobre el cerrado y compacto [0,1]. Consideremos ahora la función $f(x) =x$ definida sobre (0,1). La función  $f$ es creciente y el candidato para el máximo es 1, pero no existe en el dominio ningún $P$ tal que $f(P)= 1$. En este caso la función es continua sobre un abierto, el cual no es compacto. Si tomamos como dominio [0,1], todo se arregla, pues $f(1)=1$, y se tiene que existe un punto $P$ en el dominio de $f$ tal que $f(P) \ge f(x)$ para cualquier $x$ del dominio.

Qué significa que la función continua tome sus extremos, digamos el máximo, $Max$? Eso significa que  existe un punto $P$ en $K$, el dominio de la función, tal que $ Max = \phi(P) \ge \phi(x)$ para cualquier $x \in K$.

Entonces veamos la idea de la demostración: supongamos que nuestra función toma K y lo convierte en un subconjunto de los reales, acotado, limitado por $Max$. En los reales, podemos armar una sucesión que converja a $Max$ pero cuyos elementos pertenezcan a la imagen de $K$. La imagen inversa de dicha sucesión es una sucesión en $K$.  Como $K$ es compacto, tiene la propiedad de que toda sucesión tiene una subsucesión convergente a un punto dentro de K, por lo tanto, la imagen inversa de dicha sucesión contiene una sucesión contenida en K, la cual tiene una subsuseción convergente. Pongamos que sea $P$ su límite, el cual está en $K$. Puesto que $\phi$ es función sobre $K$ entonces $ \phi(P)$ debe existir. Por construcción, $\phi(P)$ debe ser más grande que todos lo demás valores de $\phi$ .

\

$>>>>>>>>>>>>>>>>>>>>>>>>>>>>>>>>>>>>>>>>>$

\

\

\

\color{black}

Ahora podemos probar de inmediato de qué forma  es imposible \index{estabilizar un sistema} \textbf{estabilizar un sistema} de cargas por medio de un campo eléctrico fijo.  Supongamos que tenemos un campo electrostático en el cual una carga positiva pueda estar en equilibrio estable en algún punto, P. Puesto que la carga está en equilibrio estable, una ligera desviación es restaurada por el campo. Pero para lograr la estabilidad de la carga, el potencial debe ser menor en P que en una vecindad infinitesimal que rodee a P, en particular sobre la superficie de una  cierta esfera. Pero por las propiedades de las funciones armónicas, eso no es posible pues la estabilidad de la carga en P indica que el potencial tendría un mínimo en P.

\

En conclusión, ninguna carga puede ser estabilizada por un campo electrostático, pues se violaría la propiedad de las funciones armónicas.

\

Por consiguiente, desde el punto de vista de la teoría atómica, implementada como electrostática, es muy misteriosa la existencia de cuerpos estables, sean neutros o no. Esto ya era claro hacia 1920, produciendo un terreno fértil a la búsqueda de nuevas visiones, resultado de lo cual la mecánica cuántica fue aceptada ipso facto por muchos, hacia 1928, exceptuando pensadores de hilo fino, como Einstein. Naturalmente que tenemos que pasar más luego a considerar cargas en movimiento a ver ellas que dan de sí. 

\

\

\color{blue}

\

$<<<<<<<<<<<<<<<<<<<<<<<<<<<<<<<<<<<<<<<<<$

\textbf{Blindaje antigravitatorio} \index{Blindaje antigravitatorio}:

\

Tanto la fuerza gravitacional como la eléctrica decaen cuadráticamente. Ambas admiten la noción de potencial. Sin duda que ambas obedecen la ecuación de Laplace. Será acaso que dichas fuerzas son similares en todo? Si eso es verdad,  entonces podremos construir blindaje antigravitacional, semejante al blindaje que se hace para aislar un aparato  del ruido eléctrico proveniente del exterior.

\

El \textbf{blindaje eléctrico} \index{blindaje eléctrico}  se explica como sigue:

Sea $S$ la superficie que sirve de blindaje, la cual es una superficie cerrada, compacta. Ella puede ser de papel aluminio, buena conductora y  rica en cargas. Ser buena conductora significa que toda ella está en corto circuito, es decir al mismo potencial, a potencial fijo y constante, si está en equilibrio. La ponemos en medio de un campo eléctrico cuya función potencial obedece la ecuación de Laplace. Como condición de frontera para dicha ecuación tenemos el potencial sobre $S$. Pero, como ya dijimos, dicho potencial debe ser constante en toda la superficie $S$ pues de lo contrario ésta no sería buena conductora. Tenemos $\phi = \phi \sb{o}$ sobre $S$.

Ahora bien, la ecuación de Laplace tiene una solución única. Por lo tanto, en el interior de $S$,  $\phi (x,y,z)= \phi \sb{o}$ satisface la condición de frontera y también la ecuación de Laplace. Se deduce que ella es la solución buscada. Pero si $\phi$ es constante, su gradiente es cero y por consiguiente el campo eléctrico $\vec E=-\nabla \phi$ es cero. No hay campo: hemos demostrado  que la existencia de blindaje eléctrico es posible. Qué nos prohíbe pensar en blindaje gravitacional?

Hemos propuesto una demostración alternativa de que en el interior de una esfera electrizada no hay campo eléctrico, siempre y cuando se satisfaga la ecuación de Laplace, la cual se cumple para campos cuya divergencia es cero, lo cual es cierto para campos aproximadamente cuadráticos sólo cuando la dependencia cuadrática es exacta.

Dijimos de paso que la ecuación de Laplace tiene una solución única. Eso no es obvio y su demostración es como sigue: la ecuación de Laplace es lineal. Sean $\phi$ y $\psi $ dos soluciones. La resta $\theta = \phi - \psi $ también es solución de dicha ecuación pero satisface la condición de borde $\theta =0$ sobre S.

Ahora bien, la única función que satisface tal problema es la función nula. En efecto, una función constante no nula no satisface las condiciones de borde, por tanto, decir que no es nula es lo mismo que decir que no es constante. Como la solución está definida sobre un compacto, entonces tomaría su valor extremo, digamos en P. Dicho punto P no puede estar sobre S porque $\theta $ es constante sobre S. Pero como $\theta $ es armónica,  P no puede estar en el interior de la región bordeada por S. Es decir que P no puede existir en ningún lado, o sea que P no existe. Contradicción.

\

$>>>>>>>>>>>>>>>>>>>>>>>>>>>>>>>>>>>>>>>>>$

\

\

\color{black}

Finalmente, tengamos en cuenta que la electricidad es una conquista conceptual  que unifica muy variados aspectos de la naturaleza. Veamos otro más.

La energía eléctrica puede almacenarse en un tanque al igual que el agua. Esto se demostró gracias a que en la universidad holandesa de Leyden,  Pieter van Musschenbroek inventó en 1745 el primer condensador europeo, la botella de Leyden, o sea el primer tanque para almacenar electricidad. Es probable que el primer condensador en el mundo fuese el Arca de la Alianza de los Israelitas, que constaba de capas de oro separadas por madera resinosa (Exodo 25:10) que formaban un dieléctrico muy poderoso, permitiendo descargas que podían matar a una persona (2 Samuel 6:6-7), quizás unos 20.000 voltios.

Descargas tan elevadas bien pueden inducirlo a uno a pensar que los rayos y relámpagos son manifestaciones de naturaleza eléctrica. En efecto, las corrientes de aire, unas para un lado y otras para el otro, representan excelentes ocasiones  para la fricción y la consecuente separación de cargas. Bellísimo el experimento de la cometa en 1752 por Benjamín Franklin, el cual demuestra que los rayos son descargas eléctricas: al amenazar un atormenta eléctrica se eleva hacia las nubes una cometa sujeta a un hilo buen conductor, el cual termina en dos pedacitos de papel.

Se observa que los pedacitos de papel se repelen cada vez que se descarga un rayo: los rayos son descargas eléctricas que cargan el hilo y los pedacitos de papel, causando que éstos se repelan.

Explicar los \textbf{rayos} \index{rayos} como un efecto de separación de cargas por fricción es algo que suena bien. Pero la satisfacción plena sólo se logra si en el laboratorio se pueden producir descargas tan fuertes como los rayos parecen que lo sean. Tal construcción la logró Van der Graff quien inventó un generador, basado en fricción, que lleva su nombre y que puede producir potenciales de millones de voltios, cuyas descargas se utilizan para bombardear núcleos atómicos y estudiar el producto resultante.

\section{ELECTRODINAMICA}

Estamos buscando implementar la teoría atómica  con cargas eléctricas. Ya hemos visto que con campos estáticos nos ha sido imposible explicar algo tan elemental como cuerpos neutros o estructuras estables: los átomos son misterios para la \index{electrostática} \textbf{electrostática}. En la presente sección investigaremos cómo se comportan las cargas en movimiento, cuyo estudio se llama \index{electrodinámica}  \textbf{electrodinámica}. En particular, nos encantará estudiar el modelo planetario del átomo, para ver si se puede utilizar como piedra fundamental de estructura neutra y estable.

\

Recordemos que el campo eléctrico es en mucho semejante al gravitacional. Puesto que hay ríos que corren desde los montes, quemando su energía potencial, deberá entonces existir corrientes eléctricas de un potencial eléctrico mayor a uno menor. Cómo construiremos el equivalente eléctrico a los montes en gravitación? Pues ya hemos visto los condensadores, como el Arca de la Alianza del 1500 antes de Cristo, las pilas de Volta, de 1800, y los generadores de Van der Graff de 1931.  Teniendo un generador, podremos producir  ríos de corriente eléctrica para estudiarla.

Entramos entonces a la elaboración de la electrodinámica, el estudio de las cargas en movimiento. Un hecho impresionante es que la electrodinámica absorbe a las corrientes eléctricas, al magnetismo, a la óptica y por supuesto a la electrostática.

Las corrientes eléctricas nacieron históricamente  a raíz del descubrimiento hecho por el fisiólogo Luigo Galvani en 1790 de que un escalpelo electrificado causa contracciones en una pierna de rana. Fue debido a eso que Volta dilucidó que dos metales diferentes puestos en contacto  podían servir de generadores para producir corriente (creaban migración de cargas y por lo tanto, diferencia de potencial).  El estudio de la electrofisiología también inspiró la invención del cable, construido por Werner von Siemens quien en 1847 fue cofundador de la firma Siemens und Halske para la construcción de equipo telegráfico, la que después se volvió simplemente la Siemens.

\section{MAGNETISMO}

El \index{magnetismo} \textbf{magnetismo} nace de la experiencia corriente con imanes, los cuales tienen polos por los cuales bien se atraen o bien se rechazan.

El primer experimento con magnetismo consiste en partir un magneto con la intención de separar los polos. El resultado es asombroso: cuando un magneto se parte en dos, resulta otro magneto, el cual también tiene dos polos. Nadie ha encontrado hasta el día de hoy la forma de separar los polos de un imán. Sin embargo, la creencia en la existencia de monopolos, es decir de polos separados que den origen al magnetismo, es muy atractiva porque eso predeciría inmediatamente la cuantización de la carga eléctrica. 

En términos de cálculo vectorial,  el resultado negativo en la producción de monopolos lo expresamos imaginando las líneas del campo magnético formando figuras cerradas que denotan circulación. Por supuesto que al cortar un imán y separar sus pedazos, las líneas se reorganizan en curvas cerradas que de todas formas se cierran por los polos correspondientes. En otra terminología, las líneas sólo circulan pero no tienen nacimiento ni muerte.

\

\

\color{red}
En nuestra terminología oficial,  el \index{campo magnético} \textbf{campo magnético} tiene divergencia cero. Denotaremos al campo magnético por el símbolo $\vec B$ y tenemos la ley que dice que el campo magnético sólo circula, pues no tiene fuentes ni sumideros, ni compresiones, ni expansiones:

$$\nabla \ \cdot \vec B =0
\addtocounter{ecu}{1}   \hspace{4cm} (\theecu )      $$

Esta ley (\theecu ) se considera fundamental y figura como \textbf{la cuarta ley de Maxwell}.

\color{black}

\

\

El segundo experimento importante consiste en demostrar que la electricidad estática y el magnetismo son cosas muy diferentes. Eso se demuestra muy fácilmente cargando una peinilla de plástico, probando que atrae papel y después mirando a ver si el imán atrae o no al papel o a la peinilla. El experimento es un fracaso: el magnetismo y la electrostática son en principio diferentes.

\

En 1819 Hans Christian Oersted descubrió que una corriente eléctrica creaba un campo  de naturaleza vectorial, pues atraía la aguja de una brújula si el alambre por donde fluía la corriente era paralelo a la aguja y la corriente iba en un cierto sentido. Porque si se cambiaba el sentido de la corriente la aguja se deflectaba hacia el otro lado. Naturalmente que dicho campo debe llamársele magnético.

El  experimento decisivo para dilucidar qué relación existe entre las cargas y la \index{corriente eléctrica} \textbf{corriente eléctrica} consiste en demostrar que cargas colocadas sobre un aislante en forma de disco y puestas en movimiento pasivo, cuando el disco gira, generan exactamente el mismo tipo de fenómenos que una corriente eléctrica. Ese experimento fue llevado a cabo por Henry Rowland en 1878 y demostró que una carga en movimiento genera un campo magnético, es decir, que deflectaba una aguja imanada colocada apropiadamente. Se considera que su montaje experimental es fenomenal porque tenía que medir un campo magnético que era 5 órdenes de magnitud menor que el de la tierra.

Por otro lado, ya desde 1830, Faraday había descubierto la inducción electromagnética. Hacia 1840 James Prescott Joule y Hermann Ludwig Ferdinand von Helmholtz probaron que la electricidad es una forma de energía que además  se conserva. Que dicha energía se disipa ya había sido probado por Georg Simon Ohm en 1827 con su famosa ley de Ohm: la energía consumida por una resistencia es proporcional a la corriente que la atraviesa. En 1855 James Clerck Maxwell propuso sus ecuaciones. Estas ecuaciones predecían la existencia de ondas, algo que pudo ser verificado en 1888 por Heinrich Rudolf Hertz. En 1896 Guglielmo Marconi produjo el primer telégrafo inalámbrico. Las ondas eran consideradas perturbaciones de un medio llamado éter. El experimento de Michelson y Morley demostró que el tal éter era una fantasía y en 1905 apareció la teoría de la relatividad, que sin duda puede considerarse como una hija histórica del electromagnetismo.

\

La teoría matemática del electromagnetismo empieza con la formulación de las leyes de Maxwell. Para digerirlas con alegría nos conviene mejorar nuestro entendimiento de un  concepto muy importante, el rotacional, y de un teorema que lo ayuda a aclarar.

\

\color{blue}

$<<<<<<<<<<<<<<<<<<<<<<<<<<<<<<<<<<<<<<<<<$

\

Importante: el \index{rotacional} \textbf{rotacional} y el \index{teorema de Stokes} \textbf{teorema de Stokes}.

\
La corriente eléctrica que nos llega por cable es para todos el ejemplo más evidente de lo útil que pueden ser las cargas en movimiento. El cableado por el que transita la corriente forma una red cuyo bloque elemental es un circuito, es decir que se obliga a la corriente a recorrer un camino cerrado. Pues bien, una corriente eléctrica es un río de electrones en movimiento. Ohm determinó que la corriente crea calor, es decir, desgasta su energía. Por consiguiente, la corriente debe ser mantenida por alguna fuente de energía. Con el advenimiento de la superconducción, pudimos hacer circular una corriente por una camino cerrado sin gastar energía de mantenimiento. Pero sucede que en muchos casos de la vida real, a uno le gusta despreciar la energía gastada por la resistencia del cable y en esa caso uno puede aplicar   el teorema de Stokes. Veamos pues su enunciado y su significado.

Consideremos un objeto que flota en el agua de un río. A medida que el agua se lo lleva, el objeto va rotando porque el agua, en contacto íntimo con toda la superficie del objeto, produce un efecto rotacional. Podemos predecir exactamente cuando un río hará que los objetos semisumergidos roten al avanzar. Eso se logra gracias a un dispositivo que se llama el rotacional.

\

Dado un campo vectorial $\vec F$ definimos el \textbf{rotacional} de $\vec F= (P,Q,R) $, si existe,  como

$\nabla \times \vec F = (\partial /\partial x, \partial /\partial x,\partial /\partial x,) \times (P,Q,R)$

$=(\partial R /\partial y -\partial Q/\partial z )\vec i
- (\partial R /\partial x -\partial P/\partial z )\vec j
+ (\partial Q /\partial x -\partial P/\partial y )\vec k$

El rotacional de un campo vectorial es un vector. Para entender su significado, calculemos el rotacional de un campo vectorial que forme un remolino. El más sencillo es $F=(-y,x,0)$ y su rotacional es:

$\nabla \times \vec F $

$=(\partial 0 /\partial y -\partial x/\partial z )\vec i
- (\partial 0 /\partial x -\partial (-y)/\partial z )\vec j
+ (\partial x /\partial x -\partial (-y)/\partial y )\vec k$

$=(0 -0 )\vec i
- (0 -0 )\vec j
+ (1 +1 )\vec k$

$= 2\vec k$

Vemos que este remolino que remolinea horizontalmente en dirección contraria a las manecillas del reloj tiene un rotacional no nulo que apunta hacia arriba. Por lo tanto, la mismísima naturaleza de remolino  ha sido captada por el rotacional.

Debemos aclarar que un río puede hacer girar un objeto flotante aunque no haya un remolino a la vista. Basta con que haya un gradiente de velocidades, o sea  que en un lado del objeto el agua tenga una velocidad menor que en el lado opuesto, y aunque ambas velocidades sean paralelas, el objeto tenderá a rotar a medida que avanza. Por qué? Porque si ponemos un sistema de coordenadas que avance a la velocidad promedio, que es a la cual avanza el objeto, se verá que si en un lado la  velocidad relativa tira para un lado entonces en el otro lado, la velocidad relativa tira para el lado contrario, creando así un par de fuerzas en remolino. Calculemos el rotacional de un campo de ese tipo.

Sea $\vec F = (0,x,0)$. Para $x>0$ este campo denota un río que avanza paralelo al eje 'Y' pero su velocidad aumenta a medida que nos alejamos de dicho eje. Su rotacional es:

$\nabla \times \vec F =$

$=(\partial 0 /\partial y -\partial x/\partial z )\vec i
- (\partial 0 /\partial x -\partial 0/\partial z )\vec j
+ (\partial x /\partial x -\partial 0/\partial y )\vec k$

$=(0 -0)\vec i
- (0 -0 )\vec j
+ (1 )\vec k$

$=\vec k$

lo cual dice que habrá una rotación cuyo eje será el eje 'Z'.

\

Stokes divulgó un teorema que lleva su nombre ( y que el no inventó) y que tiene un significado muy gráfico: cuando un objeto semisumergido en un río rota, todo sucede como  si hubiese una cuerda que hiciera rotar al objeto, y aún más, podemos poner la cuerda precisamente sobre la frontera que separa la parte sumergida de la emergente. Más oficialmente:

\

Sea una superficie $S$ que sea simplemente conexa, (que no tenga huecos, y por consiguiente cada curva cerrada se pueda  encoger hasta desaparecer). $S$ tiene un borde $\gamma$ que es una curva que no tiene demasiadas puntas. Sea un campo vectorial $\vec F$ que cubre un dominio abierto en $R\sp{3}$ que encierra a $S$ y que se le pueda definir el rotacional, el cual ha de ser continuo. Entonces se cumple que:

$\oint \sb{\gamma} \vec E = \int\sb{S} \nabla \times \vec E$

\

Todos esos garabatos tienen un significado muy directo. El primer término, $\oint \sb{\gamma} \vec E $ se llama una integral cerrada, y se llama así pues se trata de un camino cerrado $\gamma $. El campo $\vec E$ tiene unidades de fuerza y toda la integral representa  trabajo, o sea  energía, que se gasta para recorrer el camino cerrado $\gamma $ una vuelta completa en contra del campo $\vec E $.

El segundo término, $\int\sb{S} \nabla \times \vec E$, toma el rotacional y lo integra sobre toda la superficie semisumergida. El rotacional es un dispositivo microscópico  e indica si localmente, o sea microscópicamente, el campo $\vec E $ tiene capacidad de hacer rotar a un objeto con el cual entre en interacción. El campo puede tener un rotacional en un lado que sea contrario al que haya en otro de tal forma que es difícil predecir qué pasará en definitiva con el objeto semisumergido, si rotará o no.
La integral es una sumatoria y en este caso, la integral del rotacional sobre toda la superficie de interacción nos dice que pasará en definitiva.

\

Por lo tanto el teorema de Stokes nos dice con toda sencillez:  los efectos microscópicos del campo sobre la superficie de un objeto semisumergido en el campo se suman unos con otros haciendo que el objeto tenga un patrón definido de rotación. Pero hay que ver que el eje de rotación no nos lo dice el teorema de Stokes. Dicho teorema tan sólo dice que el campo, como portador de energía está en capacidad de entregarle energía al objeto semisumergido o de ofrecer resistencia en contra de un rotación inducida. La tal energía se puede calcular muy fácilmente recorriendo el borde del objeto que separa la parte sumergida de la emergida, que es el camino $\gamma $, y mirando que tanto trabajo hay que hacer para recorrer el camino en contra del campo.

\

El teorema de Stokes nos permite refinar nuestro entendimiento del rotacional. Supongamos que tenemos un campo que tiene un remolino, digamos el campo ya descrito como $F=(-y,x,0)$. Imaginemos que el campo puede interactuar con los electrones de un alambre circular de cobre, poniendo los electrones a girar alrededor del alambre,  creando una corriente eléctrica. Nos preguntamos cómo debemos orientar el alambre para obtener la máxima corriente posible.

Para maximizar la corriente lo que tenemos que hacer es maximizar el recurso energético que hace mover los electrones. La energía está dada por la integral cerrada $\oint \sb{\gamma} \vec E $ y tenemos que encontrar la forma de orientar a $\gamma $ para que dicha energía sea máxima.

Este problema tan complicado se resuelve rápidamente gracias al teorema de Stokes, si tenemos en cuenta que el alambre encierra una superficie que es un círculo $S$:

$\oint \sb{\gamma} \vec E = \int\sb{S} \nabla \times \vec E$

como el rotacional es $2\vec k$ nos queda:

$\oint \sb{\gamma} \vec E = \int\sb{S} 2 \vec k $

Tenemos que maximizar esta integral cuadrando la orientación del alambre que hace las veces de antena. Habíamos dicho que la integral de un campo sobre una superficie se puede entender como el flujo que atraviesa la superficie. En este caso el flujo es el debido al rotacional que es $2 \vec k $ y que apunta hacia arriba.

Razonamos entonces diciendo: para maximizar el flujo que atraviesa el círculo encerrado por el alambre tenemos que orientar al alambre para que el círculo abarque la máxima cantidad de campo. Es decir, debemos orientar el alambre perpendicularmente al eje 'Z', que está en la misma dirección del rotacional.

\

Nos atrevemos entonces a generalizar el resultado particular anterior y decimos informalmente: el rotacional de un campo vectorial es un vector cuya norma indica si hay un efecto local de rotación y cuya dirección indica la forma en la cual hay que orientar una antena para extraer del campo la máxima energía posible.

Recordando que las antenas de radio necesitan ser orientadas para que funcionen bien, debemos esperar que en las leyes del electromagnetismo aparezcan rotacionales de cuando en cuando.

\

A los términos de la forma

$\partial (\partial \phi/\partial z ) /\partial y
= \partial \sp{2}\phi/\partial z \partial y $

se le llaman derivadas mixtas o cruzadas.

\

Probemos ahora un resultado útil: el rotacional del gradiente de una función escalar  es siempre cero con sólo garantizar que sus segundas derivadas cruzadas sean continuas y por lo tanto iguales

Todas esa condiciones de continuidad siempre se asumen por defecto y a veces no se mencionan pues se cumplen sin mayor problema. Por ejemplo, si $f(x,y) = x^3y^2$ entonces

$\partial f/\partial x = 3x^2 y^2$

$\partial^2 f/\partial y \partial x  = 6x^2 y$

Por otro lado, obtenemos el mismo resultado derivando en orden inverso, lo cual ilustra que las derivadas cruzadas son iguales:

$\partial f/\partial y = 2 x^3 y$

$\partial^2 f/\partial x \partial y = 6x^2 y$

\bigskip

Ahora,  sea $\phi$ la función escalar cuyo gradiente es:

$\nabla \phi = (\partial \phi/\partial x,\partial \phi/\partial y,\partial \phi/\partial z )$

Si un campo vectorial es $\vec F= (P,Q,R)$ entonces su rotacional es:

$\nabla \times \vec F = $

$=(\partial R /\partial y -\partial Q/\partial z )\vec i
- (\partial R /\partial x -\partial P/\partial z )\vec j
+ (\partial Q /\partial x -\partial P/\partial y )\vec k$

Por lo tanto,

$\nabla \times \nabla \phi =
(\partial (\partial \phi/\partial z ) /\partial y -\partial (\partial \phi/\partial y)/\partial z )\vec i
- (\partial (\partial \phi/\partial z ) /\partial x -\partial (\partial \phi/\partial x)/\partial z )\vec j
+ (\partial ()\partial \phi/\partial y /\partial x -\partial (\partial \phi/\partial x)/\partial y )\vec k$

Vemos que el cálculo del rotacional del gradiente en cada coordenada se encuentran las derivadas cruzadas las cuales se han supuesto iguales y se aniquilan, dando como rotacional el vector cero: $\vec 0$. Hemos demostrado que, con las condiciones dadas,

$\nabla \times \nabla \phi = 0$

lo cual se lee: el rotacional de un gradiente de una función escalar siempre es cero.

\
Pero carambas, eso pareciera implicar que cuando una hoja de un arbol cae en la tierra, ella debería caer sin rotar. Pero eso no es cierto. La razón es que la tierra gira y para ese caso no es válido nuestro teorema.

\

$>>>>>>>>>>>>>>>>>>>>>>>>>>>>>>>>>>>>>>>>>$

\

\

\color{black}

Amable lector, permítame proponerle ahora la siguiente intriga para que logre entender mejor el teorema de Stokes:

El campo eléctrico $\vec E $ admite un campo escalar $\phi $, el potencial electrostático, del cual es el gradiente, $\vec E = \nabla \phi$. Puesto que el rotacional de todo gradiente es cero $\nabla \times \nabla  \phi =0 $ (que equivale a decir que las segundas derivadas cruzadas del potencial son continuas e iguales), por el teorema de Stokes ningún campo electrostático puede crear una corriente circular. En efecto, el campo eléctrico es una función vectorial definida en cualquier parte donde no haya carga. Las cargas son puntuales y son finitas en número. Por tanto, el conjunto que queda de quitar de $R\sp{3}$ los puntos ocupados por las cargas es simplemente conexo. Eso quiere decir que cada curva cerrada puede encogerse hasta desaparecer totalmente. Sea $\gamma$ una curva cerrada y $S$ la superficie creada por un proceso de encogimiento hasta desaparecer. Tenemos:

$\oint \sb{\gamma} \vec E = \int\sb{S} \nabla \times \vec E$
$ =  \int\sb{S} \nabla \times \nabla \phi =\int\sb{S} 0= 0$

Por tanto, el campo eléctrico no puede hacer trabajo circulatorio y por consiguiente no podrá mover ninguna corriente circulante. Eso implica que las pilas eléctricas no podrían generar corriente en un  circuito  eléctrico. Lo cual es una gran falsedad. Cómo se soluciona este misterio?

\

El teorema de Stokes forma con el teorema de Gauss una herramienta para el cálculo vectorial supremamente poderosa con   ilimitadas aplicaciones. Sus demostraciones pueden encontrarse en cualquier libro de cálculo vectorial o de electromagnetismo. Hay además una identidad del cálculo vectorial que vamos a necesitar y que se relaciona mucho con estos dos teoremas: es la fórmula de integración por partes en varias dimensiones. Veámosla, pues la necesitaremos para un resultado importante del cálculo variacional con respecto al campo electromagnético. 

\

\color{blue}

$<<<<<<<<<<<<<<<<<<<<<<<<<<<<<<<<<<<<<<<<<$ 
 
\

\textbf{La fórmula de integración por partes}

\

  Una versión liberal  de  la fórmula de integración por partes en una variable dice así:
 
 $\int^b_a udv = uv|^b_a -\int^b_a vdu$
 
 fórmula que se deriva de la regla del producto: si $u,v$ son integrables y con  derivada en un intervalo abierto que contenga el intervalo (a,b), entonces

$d(uv) = (du)v+udv$

$udv = d(uv) - (du)v$

Integrando obtenemos:

$\int^b_a udv =  \int^b_a d(uv) -\int^b_a vdu$

y aplicando el primer teorema fundamental del cálculo, se obtiene:

 $\int^b_a udv = uv|^b_a -\int^b_a vdu$
 
 que es una fórmula que intercambia los roles de $u$ y $v$ en una integración en la cual hay multiplicaciones de varios términos.
 
 \
 
Necesitaremos saber cómo se extiende esta fórmula a funciones y dominios con varias variables. Debemos reemplazar el dominio de integración y derivación de un intervalo por un conjunto de dimensión $n$ -pero no cualquier conjunto, sino uno que no ponga problema ni para derivar parcialmente, ni para integrar tanto en volumen como en superficie. Será suficiente si el conjunto es acotado, sin huecos y si está bordeado por  una superficie suave, que tenga vector normal bien definido en todo lado. (Se conocen generalizaciones sobre condiciones mucho más débiles.)     
 
Supongamos entonces que $\Omega$ es un conjunto abierto de $\Re^n$ acotado por una frontera que es suave (su parametrización es diferenciable). Si $ u$ y $v$ son dos funciones escalares tal que ellas y sus derivadas parciales son integrables sobre $\Omega$  y sobre su frontera $\partial \Omega$, entonces la fórmula de \index{integración por partes} \textbf{integración por partes} dice así:

\

$\int_\Omega \frac{\partial u}{\partial x_i}v dx = $
$\int_{\partial \Omega} uv \vec \nu_i d\sigma -$
$\int_\Omega u \frac{\partial v}{\partial x_i}dx$

\

donde $dx$ es el elemento de integración en $\Re^n$, $d\sigma $ es el elemento de integración sobre $\partial \Omega$ que es la superficie que bordea a $\Omega$,  $\vec \nu$ es el vector normal hacia afuera de la superficie  y  $\vec \nu_i$ es la i-ésima componente de $\vec \nu$. 

\

 Reemplazando $v$ por $\vec v_i$ para un campo vectorial $\vec v$ se obtiene:
 
 $\int_\Omega \frac{\partial u}{\partial x_i}\vec v_i dx = $
$\int_{\partial \Omega} u\vec v_i \vec \nu_i d\sigma -$
$\int_\Omega u \frac{\partial \vec v_i}{\partial x_i}dx$
 
y sumando sobre $i$ desde 1 hasta $n$, obetenemos:

\

$\int_\Omega \nabla u \cdot \vec v  dx = $
$\int_{\partial \Omega} u   \vec v \cdot \vec \nu d\sigma -$
$\int_\Omega u  \nabla \cdot \vec vdx$

\

donde $\nabla$ representa el operador gradiente.

\

Si ponemos $u=1$ obtenemos:

 $\int_\Omega \nabla 1 \cdot \vec v  dx = $
$\int_{\partial \Omega} 1   \vec v \cdot \vec \nu d\sigma -$
$\int_\Omega 1  \nabla \cdot \vec vdx$
 
es decir

 $0 = $
$\int_{\partial \Omega}    \vec v \cdot \vec \nu d\sigma -$
$\int_\Omega   \nabla \cdot \vec vdx$

que también se puede escribir como

$\int_{\partial \Omega}    \vec v \cdot \vec \nu d\sigma = $
$\int_\Omega   \nabla \cdot \vec vdx$

o bien

\

$\int_{\partial \Omega}    \vec v   = $
$\int_\Omega   \nabla \cdot \vec v $

\

que es el \index{teorema!de la divergencia o de Gauss} \textbf{teorema de la divergencia o de Gauss}.

\

Si partimos de

 $\int_\Omega \nabla u \cdot \vec v  dx = $
$\int_{\partial \Omega} u   \vec v \cdot \vec \nu d\sigma -$
$\int_\Omega u  \nabla \cdot \vec vdx$

y reemplazamos $\vec v$ por un gradiente $\nabla v$  teniendo en cuenta que $\Delta$ denota el operador  Laplaciano de segundas derivadas,  obtenemos:

\

$\int_\Omega \nabla u \cdot \nabla v  dx = $
$\int_{\partial \Omega} u \nabla v \cdot \vec \nu d\sigma -$
$\int_\Omega u  \Delta   vdx$

\

que se denomina la \index{primera identidad de Green} \textbf{primera identidad de Green}.

\

Como vemos, la integral por partes es algo muy poderoso, así que revisemos los pasos operacionales de su demostración.

Para demostrar que 

$\int_\Omega \frac{\partial u}{\partial x_i}v dx = $
$\int_{\partial \Omega} uv \vec \nu_i d\sigma -$
$\int_\Omega u \frac{\partial v}{\partial x_i}dx$

es suficiente demostrar que

$\int_\Omega \frac{\partial M}{\partial x_i} dx = $
$\int_{\partial \Omega} M \vec \nu_i d\sigma $

pues poniendo $M = uv$ se tiene 

$\int_\Omega \frac{\partial (uv)}{\partial x_i} dx = $
$\int_{\partial \Omega} uv \vec \nu_i d\sigma $

y elaborando la regla de derivación de un producto obetenemos:

$\int_\Omega \frac{\partial u}{\partial x_i}v dx   $
$+\int_\Omega u\frac{\partial v}{\partial x_i} dx =$
$\int_{\partial \Omega} uv \vec \nu_i d\sigma $

por lo que 

$\int_\Omega \frac{\partial u}{\partial x_i}v dx = $
$\int_{\partial \Omega} uv \vec \nu_i d\sigma -$
$\int_\Omega u \frac{\partial v}{\partial x_i}dx$

que es la integración por partes.

\

Demostremos entonces que 

$\int_\Omega \frac{\partial M}{\partial x_i} dx = $
$\int_{\partial \Omega} M \vec \nu_i d\sigma $

Para fijar ideas, pensemos que estamos en el espacio tridimensional y que aplicamos la fórmula tomando la derivada parcial con respecto a la dirección $z$.

Para ello, lo primero que hacemos es   dividor el sólido $\Omega$ en subsólidos de la forma

$Q = \{(x,y,z)| g(x,y) \leq z \leq h(x,y), (x,y) \in D \}$

 es decir, cada subsólido tiene un piso $g$ y un techo $h$ que lo definen totalmente, donde $g, h$ son funciones. Por tanto, el sólido $Q$ no tiene huecos, ni vacíos en la dirección $z$. Tanto el piso como el techo de $Q$ proyectan una sombra $D$ sobre el plano XY. Nuestra demostración será para sólidos de éste tipo, pero por un procedimiento de pegado nuestro resultado puede extenderse a sólidos muy complicados que se puedan particionar en subsólidos del tipo considerado. 
 
 La frontera de $Q$ tiene: techo, piso y pared lateral. Pero esta partición se hace de tal forma que la pared lateral sea precisamente vertical y por tanto su normal queda horizontal, por lo que su componente $z$ es nula: $\nu_z = 0$ sobre la pared lateral. Eso quiere decir que:
 
 $\int_{\partial \Omega} M \vec \nu_i d\sigma = $
  $\int_{Techo} M \vec \nu_z d\sigma + $
   $\int_{Piso} M \vec \nu_z d\sigma + $
    $\int_{Pared \ lateral} M \vec \nu_z d\sigma  $
    
     $\int_{\partial \Omega} M \vec \nu_i d\sigma = $
      $\int_{Techo} M \vec \nu_z d\sigma + $
   $\int_{Piso} M \vec \nu_z d\sigma $
   
\

Por consiguiente, nuestra tarea se reduce a demostrar que

$\int_Q \frac{\partial M}{\partial z} dV = $
$\int_{Techo} M \vec \nu_z d\sigma + $
   $\int_{Piso} M \vec \nu_z d\sigma $

donde $dV$ denota el elemento de volumen sobre $\Re^3$. Eso es lo mismo que demostrar que

\

$\int_{Techo} M \vec \nu_z d\sigma + $
   $\int_{Piso} M \vec \nu_z d\sigma = $
   $\int_Q \frac{\partial M}{\partial z} dV$

\

Elaboremos la integral de superficie sobre el techo. Esta superficie tiene una parametrización natural

$Techo = \{(x,y,z)| z= h(x,y), (x,y) \in D \}$

$Techo = \{(x,y,z)| z- h(x,y)=0, (x,y) \in D \}$

lo cual dice que el techo es una superficie de nivel de  $l(x,y,z) = z- h(x,y)$ y por tanto su gradiente es perpendicular al techo y al normalizarlo nos queda el vector normal unitario hacia el exterior:

$\vec \nu = \frac{\nabla l}{||\nabla l||} = \frac{-h_x \vec i - h_y\vec j + 1\vec k}{\sqrt{(h_x)^2 + (h_y)^2 + 1}} $

donde $h_x$, $h_y$ son las derivadas parciales en las direcciones correspondientes.

De acá tenemos que la componente $z$ de dicha normal es

$\vec \nu_z =\frac{1}{\sqrt{(h_x)^2 + (h_y)^2 + 1}} $

Transformemos la integral de superficie sobre el techo en una integral sobre la sombra $D$ que está en el plano $XY$. Tenemos que cambiar el elemento de integración sobre el techo en un elemento de integración sobre la sombra. Eso se logra teniendo en cuenta que  el área $dxdy$ de un rectángulo sobre el plano $XY$ y el  área correspondiente al romboide en el techo, $||d \sigma||$, se relacionan por 

$dxdy = cos \theta ||d \sigma||$

donde $\theta $ es el  (máximo) ángulo del elemento de área sobre el techo con el piso.  Ahora bien, el ángulo (máximo) $\theta$ entre la horizontal y el elemento de área del techo es el mismo que el ángulo entre la normal al techo y el eje $Z$. Como estamos sobre una normal unitaria, su componente $z$ es lo mismo que el coseno   $\theta$, pues las componentes de un vector unitario corresponden a sus cosenos directores. Nos queda:

$dxdy = \frac{||d \sigma||}{\sqrt{(h_x)^2 + (h_y)^2 + 1}}  $

$  ||d \sigma ||  =  \sqrt{(h_x)^2 + (h_y)^2 + 1}dxdy$:

\

Tenemos entonces:

$\int_{Techo} M \vec \nu_z d\sigma = $
$\int_{D} \frac{M}{\sqrt{(h_x)^2 + (h_y)^2 + 1}} \sqrt{(h_x)^2 + (h_y)^2 + 1} dxdy $

$\int_{Techo} M \vec \nu_z d\sigma = $
$\int_D M(x,y,h(x,y)) dxdy$

\

Procedemos de igual forma con el piso, con  la salvedad de que la normal unitaria sobre el piso debe apuntar hacia abajo. Obtenemos:

$\int_{Piso} M \vec \nu_z d\sigma = $
$-\int_{D} \frac{M}{\sqrt{(g_x)^2 + (g_y)^2 + 1}} \sqrt{(g_x)^2 + (g_y)^2 + 1} dxdy  $

$\int_{Piso} M \vec \nu_z d\sigma = $
$-\int_D M(x,y,g(x,y)) dxdy$

\

Sumemos ahora las dos igualdades obtenidas:

$\int_{Techo} M \vec \nu_z d\sigma + \int_{Piso} M \vec \nu_z d\sigma   $

$=\int_D M(x,y,h(x,y)) dxdy -\int_D M(x,y,g(x,y)) dxdy$

$=\int_D [M(x,y,h(x,y))  - M(x,y,g(x,y))] dxdy$

$=\int_D [ \int^{h(x,y)}_{g(x,y)} \frac{\partial M}{\partial z}dz)] dxdy$

$ = \int_Q \frac{\partial M}{\partial z}dxdydz $

$= \int_Q \frac{\partial M}{\partial z} dV$

que era lo que teníamos que demostrar.
 
 \
 
$>>>>>>>>>>>>>>>>>>>>>>>>>>>>>>>>>$

\

\color{black}

\section{LAS LEYES DE MAXWELL Y SU SIGNIFICADO}

 Las \textbf{leyes de Maxwell} \index{leyes de Maxwell}  rigen el comportamiento del \index{campo electromagnético} \textbf{campo electromagnético}. Las estudiaremos restringiéndonos al caso del vacío. Utilizamos la letra $\vec E$ para el campo eléctrico, la $\vec B$  para el magnético, la $c$ para la velocidad de la luz, $t$ para el tiempo, $Q$ para la carga eléctrica, $\vec J$ para la densidad de corriente: cantidad de carga por unidad de área que atraviesa una cierta superficie.

\

\

\color{red}
Ya sabemos que un imán y una carga eléctrica no interactúan entre ellos si están quietos. Pero la cosa es distinta cuando las cargas o los imanes se mueven. La primera ley nos dice que un imán en movimiento,  el cual  crea un campo magnético variable, simula un campo eléctrico, es decir, un imán que se mueve interactúa con una carga eléctrica en su vecindad:

$$\nabla  \times \vec E =-(1/c) \partial \vec B / \partial t
\addtocounter{ecu}{1}   \hspace{4cm} (\theecu )      $$

El término $1/c$ puede entenderse como el factor de proporcionalidad que fija la escala y también como una premonición de que del electromagnetismo sale  la teoría de la relatividad.

La siguiente intriga es muy delicada y la resolveremos enseguida: habíamos dicho que el rotacional de un campo eléctrico es cero, pues todo campo eléctrico $\vec E$ puede expresarse como el gradiente de un potencial $E=\nabla \phi$ y el rotacional de un gradiente es siempre cero: $\nabla  \times \vec E= \nabla  \times \nabla  \phi= 0 $. Pero ahora decimos que un imán en movimiento simula un campo eléctrico, el cual aparece escondido detrás de un rotacional. Como quien dice: un imán en movimiento no hace nada, y como si fuera poco a eso le llamamos la primera ley de Maxwell.

Tanto embrollo se resuelve con un poco de claridad: a un campo se le llama eléctrico cuando interactúa con la materia a través de la carga eléctrica que ésta pueda tener. En ese caso, el campo eléctrico ejercerá una fuerza igual a la carga multiplicada por el campo y por lo tanto la carga se acelerará. En esta definición no se exige ni se deduce que el campo pueda expresarse como el gradiente de alguna función potencial. Cuando el campo es producido por cargas eléctricas fijas, entonces sí se puede expresar como un gradiente, pero eso es un bono adicional.

La primera ley de Maxwell nos dice entonces que un imán en movimiento produce un campo eléctrico que no puede simularse  por cargas eléctricas y que no puede expresarse como el gradiente de ninguna función escalar.

\

La segunda ley nos dice que las cargas en movimiento generan campos magnéticos. Más claramente: una carga eléctrica que se mueve puede interactuar con un imán.  Puede hacerse referencia directamente a campos eléctricos variables, in abstracto,  o a las cargas que se mueven en una corriente eléctrica con densidad de corriente $\vec J$ :

$$\nabla  \times \vec B =-(1/c) \partial \vec E / \partial t
+ (4\pi /c)\vec J
\addtocounter{ecu}{1}   \hspace{4cm} (\theecu )      $$

Inferimos que un campo magnético, en general, no puede expresarse como el gradiente de una función escalar.

\

La tercera ley es toda una maravilla:

$$\nabla  \cdot \vec E = 4 \pi Q
\addtocounter{ecu}{1}   \hspace{4cm} (\theecu )      $$

Pudimos demostrar que el campo eléctrico creado por cargas eléctricas fijas cumple con esta ecuación: los únicos puntos donde la divergencia no es cero es donde hay cargas. Por eso podemos decir 'campo creado por cargas'. Lo increíble es que el campo eléctrico creado por imanes en movimiento también tiene la misma virtud: en todo lugar su divergencia es cero. Eso se basa en la identidad que dice: la divergencia de todo rotacional es cero. Sería bonito demostrar eso.

\

La cuarta ley nos dice que las líneas de campo magnético son circulantes sin tener nacimiento ni muerte, tanto las creadas por imanes como las creadas por cargas en movimiento:

$$\nabla  \cdot \vec B =0
\addtocounter{ecu}{1}   \hspace{4cm} (\theecu )      $$

\color{blue}

\

\

$<<<<<<<<<<<<<<<<<<<<<<<<<<<<<<<<<<<<<<<<<$

\

Punto seguido: sobre la consistencia de las leyes de Maxwell

\

Las dos ecuaciones

$x+y=4$

$3x+3y=9$

son inconsistentes.

\

De una sola inconsistencia entre un millón de verdades consistentes se puede probar cualquier cosa. Probemos a partir del anterior sistema que 1=0. En efecto:

Simplificando la segunda ecuación por 3 nos da:

$x+y=3$

que conjuntamente con la primera ecuación nos permite concluir que

$3=4$

Restando $3$ de ambos lados de la ecuación nos queda $0=1$ y como la igualdad es una relación simétrica nos queda lo deseado.

\

Serán consistentes las leyes de Maxwell?

\

Sinceramente, yo ni había pensado en tanto complique. Con todo, me parece razonable la fe ciega que todos tenemos en su consistencia. Tal fe se apoya en los siguientes hechos, seguramente discutibles:

\

1. Si un sistema axiomático puede realizarse en un modelo, es decir, describe un universo de juguete, entonces es consistente.

2. Las leyes de Maxwell han sido copiadas de la naturaleza sin sufrir tergiversaciones.

3. Por lo tanto, podemos decir que las leyes de Maxwell tienen un modelo, que además no es de juguete sino real.

4. El tal modelo (nuestro universo) es muy antiguo y por consiguiente ha tenido tiempo suficiente para poner a prueba su consistencia interna.

\

Nuestro argumento demuestra no sólo la consistencia de las leyes de Maxwell sino además su seriedad.

\

Con todo, hagámosle un pequeño \index{test} \textbf{test} a la primera ley, el cual siendo inofensivo, nos da seguridad.

La primera ley reza:

$\nabla  \times \vec E =-(1/c) \partial \vec B / \partial t$

Nuestro problema de \index{autoconsistencia} \textbf{autoconsistencia} es el siguiente: estamos seguros que al mover un imán se produce un campo vectorial, pero cómo asegurar que dicho campo es el rotacional  de otro al cual llamamos eléctrico?

Para ayudarnos, nos basamos en el siguiente teorema: para que un campo vectorial continuo se pueda expresar como el rotacional de otro campo, se requiere que su divergencia sea cero.

\

Lo que tenemos que hacer es calcularle la divergencia a la expresión

$ \partial \vec B / \partial t$

que da

$\nabla \cdot (\partial \vec B / \partial t)$
$= (\partial /\partial x,\partial /\partial y,\partial /\partial z) \cdot (\partial B\sb{x} / \partial t,\partial B\sb{y} / \partial t,\partial B\sb{z} / \partial t,)$

$=\partial (\partial B\sb{x} / \partial t)/\partial x+ \partial (\partial B\sb{y} / \partial t)/\partial y+ \partial (\partial B\sb{z} / \partial t)/\partial z $

$=\partial (\partial B\sb{x} / \partial x)/\partial t + \partial (\partial B\sb{y} / \partial y)/\partial t + \partial (\partial B\sb{z} / \partial z)/\partial t $

$=\partial /\partial t (\partial B\sb{x} / \partial x+ \partial B\sb{y} / \partial y + \partial B\sb{z} / \partial z)$

$=\partial /\partial t (\nabla \cdot \vec B)$

$=0$

donde hemos utilizado la existencia y continuidad de todas las derivadas, lo mismo que la cuarta ley de Maxwell, de que el campo magnético no tiene fuentes ni cosas parecidas. Hemos demostrado entonces que las leyes de Maxwell han pasado este primer test de autoconsistencia.

\

Vamos bien.

\

$>>>>>>>>>>>>>>>>>>>>>>>>>>>>>>>>>>>>>>>>>$

\

\

\color{black}

De acuerdo a la tercera ley, las cargas eléctricas producen campo eléctrico y cuando se mueven, según la segunda ley, producen campo magnético. Por tanto, las leyes de Maxwell reducen el magnetismo a la electrodinámica. Estas leyes son fenomenológicas y resumen los trabajos de Faraday sobre inducción. Algo muy diferente es deducir teóricamente que una corriente eléctrica produce efectivamente un campo magnético. Eso se puede hacer rigurosamente con ayuda de la relatividad especial. También, las leyes de Maxwell nos permiten predecir la existencia de ondas en el vacío, estudiar nuestro modelo planetario de átomo, formular la teoría de la relatividad y motivar la necesidad de una mecánica más general que la de Newton. Veamos de qué forma las leyes de Maxwell predicen la existencia de \index{ondas}  \textbf{ondas}.

\section{ONDAS ELECTROMAGNETICAS}

En el vacío y sin cargas ni corrientes, las leyes de Maxwell toman la siguiente forma:

$$\nabla  \times \vec E =-(1/c) \partial \vec B / \partial t
\addtocounter{ecu}{1}   \hspace{4cm} (\theecu )      $$

$$\nabla  \times \vec B =-(1/c) \partial \vec E / \partial t \addtocounter{ecu}{1}   \hspace{4cm} (\theecu )      $$

$$\nabla  \cdot \vec E = 0
\addtocounter{ecu}{1}   \hspace{4cm} (\theecu )      $$

$$\nabla  \cdot \vec B =0
\addtocounter{ecu}{1}   \hspace{4cm} (\theecu )      $$

Para demostrar que estas ecuaciones predicen la existencia de ondas, podemos transformar las ecuaciones de Maxwell en ecuaciones de onda, lo cual haremos un poco más adelante,  o podemos demostrar que las ecuaciones de Maxwell  se satisfacen por una función de onda apropiada. Tomando la última opción, construyamos una función que represente una onda viajera. La onda más simple se representa, digamos, por un seno que se alza verticalmente y que es paralela al eje 'Y':

$\vec E=\vec E\sb z =  \vec k sen(y)$

Esta es una onda que vale $0$ cuando $y=0$. Si le añadimos una fase, la onda se desplaza:

$\vec E=\vec E\sb {z} = \vec k sen(y + \phi)$

Si la fase es positiva, la onda se desplaza hacia la izquierda, es decir se retrasa. Si la fase es negativa, la onda se desplaza hacia la derecha, es decir, se adelanta. Tomemos entonces una onda de amplitud $E \sb{0}$ y una fase que dependa linealmente del tiempo, para lograr que la onda se desplace uniformemente hacia la derecha, es decir, que sea viajera con velocidad $v$:

$\vec E=\vec E\sb {z} = \vec k E\sb{0}sen(y - vt)$

Podemos entender a  $v$ como la rapidez de desplazamiento de la onda y a $E\sb{0}$ como una constante, la amplitud máxima del campo eléctrico. En cuanto al campo magnético $\vec B$, también le damos la forma de onda, pero que se levante en la dirección del eje 'X' con amplitud $B\sb{0}$ y al igual que $\vec E$, que se desplace en la dirección del eje 'Y':

$\vec B=\vec B\sb {x} = \vec i B\sb{0}sen(y - vt)$

Por construcción, el campo eléctrico es perpendicular al campo magnético. Verifiquemos ahora que nuestras ondas satisfacen todas las cuatro ecuaciones de Maxwell:

La divergencia de $\vec E$ es cero pues tiene sólo una componente en la dirección del eje 'Z' y ésta no depende de z. Similarmente, la divergencia de $\vec B $ es cero, pues tiene una sola componente en la dirección del eje 'X' y ésta no depende de x. Un poco más complicado es verificar las otras dos leyes. Tomemos la primera (15):

$\nabla  \times \vec E =-(1/c) \partial \vec B / \partial t$

Calculemos cada lado. El lado izquierdo es:

$\nabla  \times \vec k E\sb{0}sen(y - vt)$
$=\vec i E\sb{0}cos(y - vt) $

El lado derecho es:

$-(1/c) \partial / \partial t (\vec i B\sb{0}sen(y - vt))  $
$= -(1/c)(-v\vec i B\sb{0}cos(y - vt))$

Para que los dos resultados sean iguales, es suficiente suponer que

$E\sb{0}=(vB\sb{0})/c$

\

La segunda ley (16) dice:

$\nabla  \times \vec B =-(1/c) \partial \vec E / \partial t $

Reemplazando el valor del campo magnético:

$\nabla  \times \vec i B\sb{0}sen(y - vt) $
$=-\vec k B\sb{0}cos(y - vt)$

En tanto que:

$-(1/c) \partial \vec E / \partial t $
$=(1/c)(-v\vec k E\sb{0}cos(y - vt))$

Igualando las dos expresiones anteriores debe cumplirse que:

$B\sb{0}=vE\sb{0}/c$

Recordando que

$E\sb{0}=(vB\sb{0})/c$

Tenemos que

$B\sb{0}=v\sp{2}B\sb{0}/c\sp{2}$

Por tanto:

$v = \pm c$ y $B\sb{0}= \pm E\sb{0}$

\

Hemos probado  que la onda desplaza su cresta a la velocidad c, si viaja en sentido positivo y a la velocidad -c si viaja en sentido negativo. Esta velocidad debe medirse experimentalmente y se encuentra que concuerda con la velocidad de la luz, la cual ya era conocida por métodos astronómicos, pero fue medida en el laboratorio por Armand Hippolyte Louis Fizeau. Por otra parte, con ayuda del telégrafo sin hilos Abraham y Dufour et Ferrier midieron la velocidad de una onda electromagnética entre París y Washington. Para todos los efectos, las velocidades resultaron ser las mismas.

Por construcción, nuestra onda electromagnética es transversal a la dirección de propagación. Al igual que Maxwell, deducimos que la luz, que tiene propiedades de onda transversal, también es una onda electromagnética. Esta conclusión pudo parecer apresurada y ha sido estudiada al detalle.

En un experimento por Nichols y Tear en 1922, se utilizó la propiedad siguiente: la ondas luminosas infrarrojas se pueden producir con una onda tan larga como hasta de 314 micras, mientras que las ondas Hertzianas se pueden producir con una onda tan corta como de 30 micras.  Ellos estudiaron el diapasón intermedio y no pudieron discernir entre la luz y la onda electromagnética. Hemos unificado la óptica y el electromagnetismo.

\section{LA FUERZA DE LORENTZ}

Nosotros hemos trabajado muy campantes con el campo magnético y lo hemos mezclado con el campo eléctrico. Además, con Maxwell hemos propuesto las reglas para comerciar con dichos campos. Sin embargo, es necesario hacer un reclamo conceptual: Es esa mezcla de la que hablamos algo natural o es simplemente una entelequia matemática? La legitimidad de todos nuestros conceptos nos la concede la denominada fuerza de Lorentz. En general, fuerza, $\vec F$, se define como la tasa de cambio en el momento lineal, masa por velocidad, lo cual se puede observar y medir con relativa sencillez.

Imaginemos ahora que hay una partícula cargada que se mueve con velocidad $\vec v$ en un punto en donde hay un campo eléctrico $\vec E$ y otro magnético  $\vec B$. El campo eléctrico lo podemos medir fácilmente mediante la relación $\vec F= Q\vec E$, la cual es una definición de campo eléctrico, pero queda el interrogante de integrar el campo magnético a nuestros esquemas.

Definimos, esto también es una definición, el campo magnético $\vec B$ por la siguiente ecuación:

$$\vec F= Q\vec E + (Q/c)\vec v \times \vec B
\addtocounter{ecu}{1}   \hspace{4cm} (\theecu )      $$

Esta es la ecuación de la \index{fuerza de Lorentz}  \textbf{fuerza de Lorentz}, la cual nos permite tener una esperanza en cuanto a la resolución del misterio de la existencia de cuerpos neutros. En efecto, estudiemos qué le pasa a una partícula que pasa con velocidad $\vec v$ paralela al plano 'XY' a una distancia $d$ del eje 'Z' cuando hay un campo magnético uniforme y en la dirección 'Z'. Pues ella va experimentar una fuerza $F$ que es perpendicular tanto a $\vec v$ como a $\vec B$. Dicha fuerza resulta ser paralela al plano 'XY' y como es perpendicular a la velocidad, pues hace el mismo papel que una fuerza centrípeta de la gravitación. Por tanto, el campo magnético pone a girar la partícula.

\

Es la primera vez que podemos imaginar una estructura estable formada por una carga y estabilizada por un imán. A estas construcciones se les da el nombre de trampas. Este resultado es importante para el estudio experimental de las propiedades fundamentales de la materia, por ejemplo permitiendo la construcción del ciclotrón o de trampas para iones. No será esa una pista que nos lleve a entender cómo es posible imaginar un átomo estable?

\section{EL MODELO PLANETARIO DE ATOMO}

Pasemos ahora a la consideración del \index{modelo planetario} \textbf{modelo planetario} de átomo, a ver si podemos descubrir alguna manera cómoda y sencilla de garantizar la existencia de estructuras neutras (al menos si se miran desde lejos) y estables. Imaginemos, por ejemplo, una partícula positiva girando alrededor de una carga negativa bastante grande. La fuerza es atractiva pero la partícula no se cae al centro negativo por tener una apropiada energía cinética.

Un modelo planetario nos permite muchas cosas hermosas, digamos explicar la tabla de Mendeleiev o imaginar que el campo magnético de los imanes es el resultado del movimiento de cargas planetarias. Lamentablemente fue descubierto que tal modelo es ilusorio: el movimiento planetario implica que hay una fuerza centrípeta y por lo tanto una aceleración. Ahora bien, las partículas cargadas irradian energía cuando son aceleradas. Por lo tanto, la energía cinética del planeta disminuirá y terminará cayendo contra su núcleo que la atrae. Seguimos sin entender cómo existen partículas neutras. Nos veremos a caso enfrentados a la vergüenza de tener que suponer que hay partículas neutras elementales?

Recalquemos que la existencia de cuerpos astronómicos eléctricamente neutros es obligatoria, para verlo, consideremos el sistema tierra-luna. Si hubiese una descompensación eléctrica en dicho sistema, la fuerza de atracción no sería la gravitatoria y eso afectaría la dinámica muy rápida y notoriamente. Pero eso no se ve: la luna es neutra. Nuestro desafío sigue firme.

Para ver con más detalle de qué forma una partícula acelerada irradia energía estudiemos la caricatura en la cual ésta avanza a velocidad constante hasta el punto P y es frenada hasta el reposo en un tiempo $\tau$ y permanece en el punto Q. Como la velocidad de la luz es finita, en el tiempo $\tau$ dado y a una distancia mayor que la que la luz recorre en ese tiempo, hay un campo, el irradiado por la partícula desde P. Ahora,  cuando la partícula está en Q irradia de una forma que sus radios cortan a aquellos que salían desde P. Como la partícula describe una trayectoria continua, entonces el campo total debe variar continuamente y sus líneas describen   un cierto zig-zag. Ese zig-zag es un pulso, que por el teorema de Fourier puede descomponerse en ondas, senos  y cosenos.

Ahora bien, una onda electromagnética conlleva energía, como lo saben muy bien las culebras que salen a los claros del bosque a tomar el sol para calentarse y acelerar su metabolismo. La magnitud de la energía de una onda puede calcularse directamente de las leyes de Maxwell notando que $\vec E \times \vec B$ tiene unidades de potencia.

Como la energía ni se crea ni se destruye, nos preguntamos: De dónde procede dicha energía electromagnética producida por la partícula que se frena? Pues de la energía cinética que la partícula perdió al ser frenada.

Pasar de nuestra caricatura de una partícula que se frena rápidamente al estudio del modelo planetario sólo suaviza las cosas, pero el resultado final es el mismo: el planeta caería irremediablemente a su sol. Sigue siendo para nosotros un misterio entender de manera sencilla la posibilidad de reducir la existencia de cuerpos neutros a la de cuerpos cargados que se neutralizan mutuamente.

\

Nuestro modelo planetario que parecía tan prometedor no es más que un castillito de arena.  Sin embargo, la curiosa verdad es que dicho modelo no pudo ser desechado por la ciencia porque volvió a ser encontrado en la mecánica cuántica, si bien de una manera bastante transformada. Para abrirnos paso hasta ella y luego hasta la formulación gauge del electromagnetismo, nos conviene llenar algunos prerequisitos. Por ahora, formalicemos lo que hemos dicho de pasada, que una onda es portadora de energía. Veamos el cálculo de la energía de un campo electromagnético.

\section{EL TEOREMA DE POYNTING}

Mantener irradiando a una emisora de radio o de TV demanda muchos kilowatios de potencia, los cuales son producidos, digamos, por una caída de agua que  quema la energía cinética del agua adquirida por un intercambio de energía potencial, la cual a su  vez fue adquirida gracias a la luz del sol que siendo una onda electromagnética, evaporó el agua de los mares. Queda demostrado que la energía electromagnética es tan \index{energía} \textbf{energía} como la cinética o la potencial.

Busquemos unidades de energía entre nuestros elementos. La carga $Q$ se mide en culombios, y uno de ellos contiene $6,25 \times 10\sp{18}$ electrones. La  corriente en un alambre se mide en amperios. Un amperio es un culombio por segundo. La diferencia en el potencial eléctrico se mide en voltios mediante la siguiente definición: un julio de energía se requiere para mover una carga de un culombio a través de una diferencia de potencial de un voltio. Dicho de otra manera: un voltio es el potencial eléctrico ganado por una carga de un culombio al aplicarle un julio en contra del campo. El campo eléctrico se mide en voltios por metro. En las ecuaciones de Maxwell la corriente $\vec J$ ya no va por un alambre sino que puede tratarse de unos rayos cósmicos que vienen hacia la tierra. Por lo tanto, la corriente puede medirse en amperios por metro cuadrado.

Por consiguiente, la cantidad $\vec E \cdot \vec J$ tiene dimensiones de (voltio/metro) x culombio/segundo/metro cuadrado= [julio/(metro x culombio)] x culombio/metro cuadrado = (julio/seg)/metro cúbico. Los julios son energía, mientras que los julios por segundo son potencia o sea la capacidad de entregar energía a gran rapidez. En resumen,  $\vec E\cdot \vec J$ tiene unidades de potencia por unidad de volumen o de producción de energía por unidad de tiempo y de volumen. Su valor lo podemos calcular con ayuda de las ecuaciones de Maxwell reformuladas no sólo para el vacío sino para cualquier medio isotrópico. La primera ley:

$\nabla  \times \vec E =-\mu \partial \vec B / \partial t$

La constante $\mu$ aparece para indicar, entre otras cosas, que la velocidad de la luz no es la misma en todos los medios. La segunda ley dice:

$\nabla  \times \vec B $
$=-\epsilon \partial \vec E / \partial t +  \kappa \vec J$

La constante $\epsilon$ indica que el medio puede aumentar o disminuir la intensidad del campo magnético debido, digamos, a átomos que funcionan como imancitos. Las otras dos ecuaciones son las mismas. Para calcular $\vec E \cdot \vec J$ multiplicamos, en producto escalar, la segunda ecuación por $\vec E$

$\vec E \cdot \nabla  \times \vec B $
$=\vec E \cdot (-\epsilon \partial \vec E / \partial t) +  \vec E \cdot \kappa \vec J$

y despejamos:

$\kappa \vec E \cdot \vec J$
$= \vec E \cdot \nabla \times \vec B $
$- \epsilon \vec E \cdot \partial \vec E / \partial t$

Teniendo en cuenta que:

$\nabla\cdot (\vec E\times \vec B) $
$= \vec B \cdot \nabla \times \vec E $
$- \vec E \cdot\nabla \times \vec B$

despejamos

$ \vec E \cdot\nabla \times \vec B $
$=    \vec B \cdot \nabla \times \vec E $
$- \nabla\cdot (\vec E\times \vec B)$

y substituimos:

$\kappa \vec E \cdot \vec J$
$= \vec B \cdot \nabla \times \vec E$
$ - \nabla\cdot (\vec E\times \vec B) $
$- \epsilon \vec E \cdot \partial \vec E/\partial t$

Como podemos substituir el rotacional del campo eléctrico según la primera ley de Maxwell nos queda:

$\kappa \vec E \cdot \vec J$
$= - \mu \vec B \cdot \partial \vec B/\partial t $
$- \nabla\cdot (\vec E\times \vec B) $
$- \epsilon \vec E \cdot \partial \vec E/\partial t$

Ahora bien, tanto para $\vec B$ como para $\vec E$ tenemos que:

$\vec B \cdot \partial \vec B/\partial t $
$= (1/2)(\vec B \cdot \partial \vec B/\partial t $
$+ \partial \vec B/\partial t \cdot \vec B) $
$= (1/2) \partial (\vec B \cdot \vec B /\partial t) $
$= (1/2) \partial \|\vec B\|\sp{2} /\partial t$

por tanto tenemos que

$ \kappa \vec E \cdot \vec J$
$= - (\mu /2) \partial \|\vec B\|\sp{2} /\partial t $
$- \nabla\cdot (\vec E\times \vec B) $
$- (\epsilon /2) \partial \|\vec E\|\sp{2} /\partial t$

Como $\vec E \cdot \vec J$ tiene unidades de potencia por unidad de volumen, cuando se trate  de una región cerrada, debemos integrar sobre todo el volumen para saber la cantidad total de potencia que se genera o se consume. Tenemos por tanto:

$-\int \sb {V} \kappa \vec E \cdot \vec J d \nu $
$= \partial /\partial t [\int \sb{V} (1/2) (\mu \|\vec B\|\sp{2}$
$ + \epsilon   \|\vec E\|\sp{2})] $
$+ \int \sb{V} \nabla\cdot (\vec E\times \vec B)$

Aplicando el teorema de Gauss tenemos:

$\int \sb{V} \nabla\cdot (\vec E\times \vec B) $
$= \int \sb{\partial V} \vec E \times \vec B$
entonces:

$-\int \sb {V} \kappa \vec E \cdot \vec J d \nu $
$= \partial /\partial t [\int \sb{V} (1/2) (\mu \|\vec B\|\sp{2}$
$ + \epsilon   \|\vec E\|\sp{2})] $
$+ \int \sb{\partial V} \vec E \times \vec B $

En general, la corriente $\vec J$ tiene dos componentes: la que es mantenida activamente $\vec J \sb{a}$ y da origen a campos electromagnéticos y la que es pasiva $\vec J \sb{p}$ puesto que es inducida por el campo en algún medio susceptible,  $\vec J \sb{p} = \sigma \vec E$. Por tanto, normalizando $\kappa=1$ pero sin olvidar las   unidades para que todo quede en términos de potencia, tenemos:

$-\int \sb {V} \vec E \cdot \vec J d \nu $
$= -\int \sb {V}  \vec E \cdot \vec J \sb{a} d \nu $
$-\int \sb {V}  \vec E \cdot \vec J \sb{p} d \nu $
$= -\int \sb {V}  \vec E \cdot \vec J \sb{a} d \nu $
$-\int \sb {V}  \vec E \cdot \sigma \vec E d \nu $
$= -\int \sb {V}  \vec E \cdot \vec J \sb{a} d \nu $
$-\int \sb {V}   \sigma \|\vec E\|\sp{2} d \nu $

$ = \partial /\partial t [\int \sb{V} (1/2) (\mu \|\vec B\|\sp{2}$
$ + \epsilon   \|\vec E\|\sp{2})] $
$+ \int \sb{\partial V} \vec E \times \vec B$

Este es el \index{teorema de Poynting} \textbf{teorema de Poynting} y se interpreta término por término como sigue:

\

\addtocounter{ecu}{1}
\textit{Teorema \theecu } Teorema de Poynting

\textit{El flujo neto de potencia, que siendo un pérdida va con signo menos,   $-\int \sb {V} \vec E \cdot \vec J d \nu $  es el generado por las corrientes activas $ -\int \sb {V}  \vec E \cdot \vec J \sb{a} d \nu $ mas el consumido por las pasivas $-\int \sb {V}   \sigma \|\vec E\|\sp{2} d \nu $, es igual al flujo debido a los campos eléctrico y magnético $  \partial /\partial t [\int \sb{V} (1/2) (\mu \|\vec B\|\sp{2}$
$ + \epsilon   \|\vec E\|\sp{2})] $ mas el flujo de potencia debido a la propagación del campo $\int \sb{\partial V} \vec E \times \vec B$ }.

\section{EL POTENCIAL VECTOR}

Anteriormente dijimos que la representación del campo electromagnético requiere de un espacio de 6 dimensiones, 3 para el eléctrico y otras 3 para el magnético. Es posible que eso le hubiese hecho pensar a alguien que, digamos,  dado un campo eléctrico uno podía  cuadrar el campo magnético que a bien quisiese. Pues eso no es cierto. Precisamente, las leyes de Maxwell nos dicen que el campo electromagnético es un organismo tal que si uno recuadra un término, debe esperar modificaciones en los demás. Por lo tanto, estamos en el deber de preguntarnos: qué se puede decir del subconjunto de valores permitidos en ese espacio de 6 dimensiones para un experimento dado?

Como estamos pensando en un experimento arbitrario, no podemos decir algo tan certero como: ese espacio es una esfera. Pero hay algo que es cierto para todos los casos: el subconjunto de valores permitidos formado por la representación de un campo electromagnético cualquiera forma una variedad de dimensión 4. Lo de ser variedad sale de ser la solución a una sistema de ecuaciones diferenciales con solución única. Lo de tener dimensión 4 es algo que veremos con mucho detalle: demostraremos que existe una entidad con 4 grados de libertad, 4 variables a los que se les puede dar valores arbitrarios,  que satisface  las ecuaciones de Maxwell.  Para un experimento real, con condiciones iniciales (al iniciar el experimento) y de frontera (que se mantienen mientras corre el experimento), los valores permitidos formarán una trayectoria en dicho espacio de dimensión cuatro. A dicha entidad se le llama el potencial vector.

\

Para demostrar nuestra aseveración necesitamos una maquinaria poderosa que iremos construyendo poco a poco y que nos llevará más allá de la relatividad. Por ahora, hablaremos del potencial vector clásico, el cual sale del siguiente teorema: Si $\vec F$ es un campo vectorial cuya divergencia es cero, entonces existe otro campo $\vec G$ tal que $\vec F$ es su rotacional. La demostración del teorema se hace por construcción. Lo aplicamos al campo magnético, pues una de las leyes de Maxwell nos dice precisamente que el campo magnético no tiene fuentes, o sea que su divergencia es cero.

Existe entonces al menos un campo, llamado el potencial vector, $\vec A$,  tal que $\vec B = \nabla \times \vec A$. Por supuesto que no olvidaremos que el potencial vector no es un escalar sino que es un vector del cual sale el campo magnético.

\chapter{LA LIBERTAD GAUGE}

Habíamos dicho que el campo magnético no tiene por qué ser elemental, sino que puede ser considerado como derivado de otra entidad más elemental, el potencial vector: existe  un campo, llamado el potencial vector, $\vec A$,  tal que $\vec B = \nabla \times \vec A$.

En realidad hay un número infinito de campos $\vec A$ con la propiedad de que $\vec B = \nabla \times \vec A$. En efecto, si existe un $\vec A$ tal que $\vec B = \nabla \times \vec A$, entonces $\vec A + \nabla f$ también sirve. Eso se debe a que $\nabla \times \nabla f=0$.  Mejor dicho:

$$\nabla \times ( \vec A + \nabla f)= \nabla \times  \vec A + \nabla \times  \nabla f = \nabla \times  \vec A =\vec B  $$

Este es un detalle insignificante que sorpresivamente se ha tornado en una piedra angular de la física moderna. Por eso le dedicaremos la atención conveniente.

\section{EL GAUGE MAGNETICO}

Al igual que el potencial gravitatorio o eléctrico que están definidos módulo una constante, el potencial vector está definido módulo el gradiente de una función escalar.

Por construcción, el gradiente adicional no tiene ninguna influencia física, puesto que lo único que importa es la fuerza de Lorentz y ella depende tan sólo de $\vec E$ y de $\vec B$, que no perciben ninguna diferencia al cambiar de gradiente. Cuando se cambia de gradiente, muchos autores dicen que se cambia de gauge o que hacemos una transformación gauge. Nosotros preferimos decir: Hemos descubierto que las leyes de Maxwell permiten una arbitrariedad matemática y   que existe un principio gauge de sentido común: las arbitrariedades matemáticas no pueden tener ninguna implicación observable. Hemos dicho de paso que el potencial vector no puede observarse por ninguna razón.

\

Enfatizamos:

\

\

\color{red}
\emph{Las leyes clásicas del electromagnetismo  son invariantes ante las transformaciones \index{gauge} \textbf{gauge}, o sea aquellas que fijan a capricho la arbitrariedad nacida de que el vector potencial magnético no puede definirse exactamente sino tan sólo módulo el gradiente de una función escalar.}

\

\

\color{black}

Es costumbre reescribir este tipo de ideas ayudándose del concepto de grupo. Para ello, busquemos nuestro grupo en el resultado anterior. Nuestro grupo es el conjunto de funciones   escalares que tienen gradiente.   
\

\

\color{black}

 El conjunto de funciones   escalares que tienen gradiente forma un grupo con la suma, pues la suma de tales funciones es conmutativa,  asociativa, tiene elemento neutro la función cero, cada elemento tiene su inverso aditivo. Demostremos que la suma sea cerrada. Para ello, tomamos dos funciones escalares que tengan gradiente, y tenemos que demostrar que su suma también tiene gradiente. Pero eso se reduce a decir que $\nabla (f + g) = \nabla f + \nabla g$.

\bigskip

Podríamos entonces redactar nuestro resultado anterior como sigue:

Las leyes clásicas del electromagnetismo  son una teoría \index{gauge} \textbf{gauge}  cuyo grupo de invariancia tiene como elementos a las funciones escalares que tienen gradiente.

Pero, cuáles son las funciones que tienen gradiente? Pues aquellas que son diferenciables. Por lo tanto, decimos

\

\color{red}
\emph{Las leyes clásicas del electromagnetismo  son una teoría \index{gauge} \textbf{gauge}  cuyo grupo de invariancia tiene como elementos a las funciones escalares $\phi : \Re^3 \rightarrow \Re$ diferenciables.}

\

\color{black}

Lamentablemente, hay muy poco poder creativo en este enunciado. Es decir, nos queda difícil imaginar a partir de lo dicho cómo ha de ser una interacción fundamental. Nos conviene estudiar el tema más a fondo. En lo que sigue, entenderemos muy claramente algo mucho más luminoso y es que:

La \index{arbitrariedad} \textbf{arbitrariedad} gauge del electromagnetismo clásico,  en permitir adicionar el gradiente de una función escalar al potencial vector sin producir cambios observables, tiene su contraparte en la mecánica cuántica:   un cambio de gauge clásico está ligado en mecánica cuántica a un cambio de fase o exponente imaginario. Dichas fases definen un grupo, el famoso $\textbf U(1)$, el cual tampoco crea efecto físico en la formulación cuántica del electromagnetismo. Por eso diremos que la interacción entre el campo electromagnético y las partículas cargadas se describe por una teoría gauge $\textbf U(1)$.

\

El trabajo que sigue nos demandará mucho esfuerzo y bien vale la pena preguntarse si hay alguna justificación para seguir hasta demostrar algo tan simple como lo ya expresado. Pues bien, lo que nos anima es que estamos trabajando para entender las leyes fundamentales de la materia y estamos elaborando exhaustivamente un ejemplo de un principio que tiene que ser universal: las leyes de la física, es decir, las ecuaciones que se refieren a resultados experimentales, no pueden depender de ninguna arbitrariedad que resulte de los formalismos matemáticos utilizados para describirlas. Por consiguiente, podemos decir ahora con más entendimiento:

\

\

\color{red}
Si cambiar de arbitrariedad es cambiar de gauge, entonces es absolutamente necesario que  una teoría física  sea invariante ante toda transformación gauge para que tenga algún futuro científico.

\

\

\color{black}

En relación con el campo electromagnético encontraremos un grupo, $\textbf U(1)$. La visión gauge nos demanda que para cada clase de fenómenos fundamentales encontremos unas leyes que permitirán un juego de arbitrariedades, las cuales se describirán por un grupo adecuado \emph{G}.

\

Nos quedan dos  tormentos. Primero:  acaso existirá alguna teoría que no sea gauge? Segundo: En qué es más importante la libertad gauge que la libertad de un guitarrista  de improvisar una variación en su canción predilecta? Mejor dicho, para qué tanto alboroto con las arbitrariedades?

Lo que sucede es que uno puede formular la pregunta al revés: dado un \index{grupo de arbitrariedades} \textbf{grupo de arbitrariedades}, un grupo   $G$, un grupo que de antemano se ha escogido sabia y prudentemente, que teorías interesantes hay por ahí que tengan como grupo gauge a $G$? Con esa reformulación, las arbitrariedades gauge adquieren poder creativo, cuyo impacto ha sido demostrado sobre el estudio de las interacciones débil y fuerte. Nosotros estamos labrando los prerrequisitos donde todo eso se ha edificado.

Vale la pena anotar que, en retrospectiva, el poder creativo de la idea gauge ya fue exhibido en la historia de la ciencia. Se trata de la invención de la segunda ley de Newton: Fuerza = masa por aceleración. Entendamos eso mejor:

En primer término, la idea de fuerza es un concepto matemático y no físico. En física, la idea de fuerza  representa lo mismo que lo que el dinero representa en las transacciones humanas. El dinero no es necesario pero ha servido de mucho. Así mismo, el concepto de fuerza y de muchas otras definiciones más, tales como el campo electromagnético, han resultado supremamente útiles. Concretamente, la fuerza no es más que una entelequia matemática que permite igualar diversos objetos físicos que tengan las mismas dimensiones que masa por aceleración, por ejemplo, el estiramiento de un resorte y el producto masa por la aceleración en caída libre que dé el peso de un cuerpo.

En segundo término, la segunda ley de Newton, para cuerpos con masa constante,  puede escribirse como

$\vec F = m \vec a = m d(\vec v)/dt $

Escrita de esa manera, queda claro que hay una arbitrariedad en la definición de fuerza: uno puede adicionarle a la velocidad considerada una velocidad inicial arbitraria y el resultado físico no cambia. El primero que reportó haberse dado cuenta de eso fue Galileo: el iba en un barco y se dio cuenta que en algunas ocasiones el pensó que estaba quieto pero al revisar su entorno se dio cuneta que se estaba moviendo.  Eso significaba que el organismo de Galileo no podía percibir las velocidades absolutas. Le quedaba por saber si ese era problema de su organismo o de la naturaleza en general.  Para dilucidarlo, se puso a experimentar con una jarra y un pocillo: si uno llenaba la jarra con el agua que caía del pocillo, daba lo mismo si uno se movía con respecto al mar o si estaba quieto. Por lo tanto, la naturaleza como un todo era insensible a una velocidad absoluta.

Adicionalmente, Galileo también logró hacer que el tiempo pasara en cámara lenta para poder estudiar mejor algunos fenómenos relacionados con la fuerza: se inventó el plano inclinado.

Galileo se sentía tan impresionado por su método de trabajo, pensar y experimentar para volver a repensar, que declaró una idea que mi subconsciente recuerda así: de ahora en adelante, cualquiera que tenga ojos y un tris de cerebro, podrá encontrar el camino a la verdad. Al método de Galileo se le llama el método científico y hoy en día ha sido tan sobredarrollado que se aplica de una manera un poco distinta a la de Galileo, por ejemplo así: lo que hemos encontrado, lo encontramos con el método científico, entonces hemos encontrado la verdad y si no es la verdad absoluta, entonces es algo físicamente indistinguible de ella, y lo demás, por supuesto, no tiene importancia.

Como conclusión, Galileo infirió,  usando nuestra terminología, que las leyes de la mecánica son una teoría gauge con grupo $R^3$, el grupo de las velocidades. Este grupo puede extenderse a $R^6$, pues tampoco podemos definir una posición privilegiada en el espacio. Newton propuso además que la ley fundamental fuese dada por una ecuación diferencial de segundo orden, en la cual no aparecieran ni las posiciones ni las velocidades. Entonces,  Newton formuló su segunda ley. Esa ley da una ecuación diferencia lineal. Me da la impresión de que dicha linealidad puede inferirse considerando un barquito que nada en una piscina que está en otro barco más grande, es decir, considerando movimiento con referencia a sistemas inerciales encadenados y aduciendo que todas las descripciones tienen que ser equivalentes.

\section{GAUGES USUALES}

Tener o no tomar en cuenta la libertad de calibración o gauge tiene efectos cuantitativos a la hora de cuantizar tanto en teoría de campos como en teoría de cuerdas. Por lo tanto, cuando hay libertad gauge es necesario tomarla en cuenta. Pero, cómo hacerlo? En algunos casos, es suficiente  implementar la libertad gauge de la manera que a uno le guste y se procede. Veamos varias maneras usuales de fijar la libertad gauge: el gauge temporal, el gauge axial, el de Coulomb o de radiación, el de Lorentz. Queda el de Feynman, pero ese lo veremos más luego. En todos los casos se  empieza con una descripción cualquiera del campo: un potencial eléctrico $\phi_o$ y un potencial vector $\vec A$. Con dicha descripción, uno fabrica otra $(\phi_k, \vec A_{\Lambda})$ añadiendo un escalar al potencial y un gradiente al vector potencial; $\phi_k=\phi_o+k, \vec A_{\Lambda}=\vec A +\nabla \Lambda$

\begin{enumerate}
\item Gauge temporal: en este gauge se define el potencial $\phi$ como $0$ y se toma $\Lambda = \int_{-\infty}^t\phi_o(x,y,z,t')dt'$.

\item Gauge Axial: en este gauge se  orienta el mundo de tal manera que el eje Z queda perpendicular al campo electromagnético: $A^z=0$. Este gauge este generado por $\Lambda = \int_{-\infty}^z\phi_o(x,y,z',t)dz'$.

\item Gauge de Coulomb o de radiación: El gauge temporal no determina totalmente el gauge de trabajo. Puede adicionársele la condición: $\nabla \cdot \vec A=0$. Todos estos gauges tienen el inconveniente de no ser covariantes, es decir, de requerir diccionarios de traducción complicados entre dos investigadores cuyos laboratorios se muevan en marcos inerciales y que quieran intercambiar y comparar resultados. El siguiente gauge no tiene ese problema.

\item Gauge de Lorentz. Este gauge es muy importante por dos razones: es covalente y sale naturalmente al querer reexpresar las ecuaciones de Maxwell como ecuaciones de onda.

 \end{enumerate}

Transformemos, pues, las ecuaciones de Maxwell en \index{ecuaciones de onda}  \textbf{ecuaciones de onda} y veamos qué forma toma el gauge de Lorentz. Las ecuaciones de Maxwell con cargas y corrientes son:

\

Ley de Gauss $$\nabla\cdot \vec E = \rho $$

No a los monopolos magnéticos $$\nabla \cdot \vec H = 0$$

Ley de Faraday $$\nabla \times \vec E + \partial \vec H/\partial t=0$$

Ley de Amp\`ere $$\nabla \times \vec H -\partial \vec E/\partial t =\vec J$$

Como ya sabemos, $\nabla \cdot H = 0$ implica que $\vec H=\nabla \times \vec A$. Substituyendo esto en la ecuación de Faraday y teniendo en cuenta que las derivadas parciales conmutan, tenemos:

$0 = \nabla \times \vec E + \partial /\partial t(\nabla \times \vec A) =$
$\nabla \times \vec E + \nabla \times \partial \vec A/\partial t = \nabla \times (\vec E + \partial \vec A/\partial t)$

lo cual dice que hay un campo cuyo rotacional es cero. Por tanto, dicho campo es conservativo: existe un potencial $\phi$ tal que

$$\vec E +\partial \vec A/\partial t = -\nabla \phi $$

o lo que es lo mismo

$$\vec E =-\partial \vec A/\partial t  -\nabla \phi $$

Por otro lado, la ley de Amp\`ere puede escribirse como     $\nabla \times (\nabla \times \vec A) -\partial \vec E/\partial t =\vec J$

que gracias a la identidad $\nabla \times (\nabla \times \vec A) = \nabla (\nabla \cdot \vec A)-\nabla^2 \vec A$  se transforma en

    $$\vec J= \nabla (\nabla \cdot \vec A)-\nabla^2 \vec A -\partial \vec E/\partial t $$

recordando que $\vec E =-\partial \vec A/\partial t  -\nabla \phi $ esta ley queda como

$$\vec J= \nabla (\nabla \cdot \vec A)-\nabla^2 \vec A -\partial /\partial t(-\partial \vec A/\partial t  -\nabla \phi ) $$

Reordenando

$$\nabla^2 \vec A-\partial^2  /\partial t^2( \vec A)= -\vec J+ \nabla (\nabla \cdot \vec A)   +\partial ( \nabla \phi ) /\partial t$$

intercambiando derivadas parciales en el último término

$$\nabla^2 \vec A-\partial^2  /\partial t^2( \vec A)= -\vec J+ \nabla (\nabla \cdot \vec A)   -\nabla\partial   \phi  /\partial t = -\vec J+ \nabla [(\nabla \cdot \vec A)   +\partial   \phi  /\partial t] $$

Esta es una ecuación complicada. El truco para resolverla consiste en hacer que esta ecuación se simplifique a la fuerza, aquí se hace uso de la libertad gauge,   haciendo  que

$$\nabla \cdot \vec A   +\partial   \phi  /\partial t =0$$

o su equivalente

$$\nabla \cdot \vec A   =-\partial   \phi  /\partial t $$

A esta selección se le llama el gauge de Lorentz. Por supuesto, es necesario verificar que esta restricción viene naturalmente de restar del potencial vector el gradiente de una función escalar. Aceptado que eso es posible, la ecuación de Gauss también puede transformarse en una ecuación de onda:

$$\rho = \nabla \cdot \vec E = \nabla \cdot (-\nabla \phi -\partial \vec A/\partial t) = -\nabla^2 \phi  -\nabla \cdot \partial \vec A/\partial t = - \nabla^2 \phi -\partial/\partial t(\nabla \cdot \vec A)$$

Utilizando el gauge de Lorentz para la divergencia de $\vec A$ tenemos:

$$\rho = - \nabla^2 \phi -\partial/\partial t(-\partial \phi/\partial t)= - \nabla^2 \phi +\partial^2 \phi/\partial t^2$$

si se reordena para que quede como ecuación de onda, su sencillez salta a la vista:

$$\partial^2 \phi/\partial t^2- \partial^2 \phi/\partial t^2 = \rho$$

 Uno de los méritos de este  resultado es queda demostrado que hay un gauge  gracias al cual las ecuaciones de Maxwell se convierten en ecuaciones de onda.

Con esto terminamos nuestra exposición del campo electromagnético clásico y pasamos enseguida al estudio de la acción, la cual nos permitirá hacer la formulación cuántica del electromagnetismo. La bibliografía para este capítulo está al final del siguiente, pues hay una sola para ambos capítulos.

\section{HISTORIA DE LAS TEORIAS GAUGE }

La palabra del idioma inglés \index{gauge}  '\textbf{gauge}' significa 1. Medida estándar 2. Distancia entre los rieles de una carrilera 3. Instrumento para medir el nivel o cantidad de algo, 4. Hecho o circunstancia que puede usarse para caracterizar a alguien ( Oxford Advanced Learner's Dictionary). El primer significado nos diría  que 'gauge' es una unidad que se utiliza para calibrar todos los instrumentos. Las teorías gauge o de calibración se denominan así porque, como se verá más abajo,    nacieron como un resultado histórico de una preocupación relacionada con las unidades de medida.
\

\

\

Veamos un detalle que ilumina la historia de la forma como nació el nombre de las teorías gauge.

\

Los países poderosos tienen sus propios sistemas de medidas y terminan imponiéndolos a los otros. En la actualidad tenemos dos sistemas usuales: el inglés (pie - libra -segundo)  y el francés  o cegesimal (centímetro -gramo-segundo). Pues bien, a finales del 1999 hubo un desastre aeroespacial debido a que el conjunto total fue armado por partes, una en un laboratorio y otra en otro, pero dichos laboratorios manejaban unidades de medida distintos: uno el sistema cegesimal y otro el sistema inglés. Pero no notaron eso sino hasta que el sistema total colapsó, todo por falta de una tabla de conversión del programa de control global.

Imaginémonos ahora, para ver como nacieron las teorías gauge, que el mundo está lleno de países infinitesimales y que cada uno tiene sus propias unidades de medida. Lo que nos interesa es, digamos, describir la presión atmosférica a nivel de tierra. De tanto querer comunicarse entre países y no poder, las unidades de medida tendrán que evolucionar. Al llegar a un equilibrio, un observador neutral podría describir la situación como que hay una unidad de medida que varía continua y suavemente de un lugar a otro.

Si en $\vec x$ la unidad de medida se define como $1$, entonces en $\vec x + \vec dx$ la unidad de medida, convertida a las unidades que se usan en $\vec x$, sería $1 +s(\vec x,\vec dx)$. Teniendo como mira al campo electromagnético, H. Weyl  se permitió reescribir, hacia 1920, $s(\vec x,\vec dx)= \vec S(\vec x)\cdot \vec dx$, un producto punto entre el diferencial de posición $\vec dx$ y la función vectorial $\vec S(\vec x)$ que no tenía que involucrar  ningún gradiente. Ahora solicitemos a cada país $\vec x$ que nos comunique  la presión atmosférica que se registra en su localidad $\vec x$. Lo que cada país reporte lo llamaremos $f(\vec x)$. Naturalmente que $ f(\vec x)$ vendrá en la medida nacional de cada lugar. Denotemos por $g(\vec x)$ la conversión de $f(\vec x)$ a las unidades que se usan en el país donde nosotros vivimos.

Calculemos   $g(\vec x+ \vec dx)$, que es la conversión de $f(\vec x+\vec dx)$ a las unidades que se usan en $\vec x$:

\begin{center}
$g(\vec x + \vec dx) = f(\vec x+\vec dx) (1 + \vec S(\vec x)\cdot \vec dx) = (f(\vec x) + \nabla f \cdot \vec dx) (1 + \vec S(\vec x) \cdot \vec dx)
= f(\vec x) + \nabla f \cdot \vec dx + f(\vec x)\vec S(\vec x)\cdot \vec dx + (\nabla f \cdot \vec dx)(\vec S(\vec x) \cdot \vec dx) $
\end{center}

que desechando los términos de segundo orden se aproxima por

\centerline {$ g(\vec x+\vec dx) = f(\vec x) + (\nabla + \vec S(\vec x))f(\vec x)\cdot \vec dx.$ }

Esa sería nuestra teoría gauge más primitiva. Lo es tanto que es una simple tabla de conversión infinitesimal de unidades de medida o gauge.

\

Como veremos, la expresión $\vec S + \nabla $, que apareció directamente como tabla de conversión de unidades de medida o gauge, se encuentra en cálculos que involucran al  campo electromagnético y Weyl procuró, en vano, ver en la función vectorial $\vec S$ el equivalente al potencial vector electromagnético, al cual conoceremos a fondo. Con la mecánica cuántica, Fock y London se dieron cuenta de que la expresión $\nabla  + \vec S$ reaparecía, pero esta vez a la función $\vec S$ se  la comparaba con la cantidad $i$ veces  el potencial vector . Desde entonces las ideas se han aclarado mucho mejor y ahora podemos enunciar el principio gauge, como lo hicimos en el resumen del presente trabajo y notar que es algo sencillo y perfectamente natural.

Debido a que la teoría que hoy en día describe la interacción electromagnética es el prototipo de teoría gauge, nos sentimos obligados a hacer una revisión global y entendible de dicha teoría, reformulando sus principios básicos al derecho y al revés. Por supuesto que centraremos la atención en las leyes de Maxwell, el significado de sus conceptos y su aplicación. Notaremos su poder predictivo, el cual incluso  da lugar a la relatividad  especial. Destacaremos un experimento que demostró que el electromagnetismo dentro del marco de la mecánica clásica presenta falsas predicciones e invitaremos a la mecánica cuántica a solucionar la inconsistencia. En este último formalismo es donde veremos a la teoría de la interacción electromagnética como una teoría gauge con grupo de simetría \textbf{U(1)}, el cual es un grupo de Lie.

El detalle que hemos narrado se enmarca dentro de una gran historia. Ver, por favor, 

http://arxiv.org/abs/hep-ph/9810524

\section{LA VISION GAUGE}        

Quisiéramos enfatizar que las teorías gauge tienen su soporte en el más elemental sentido común. Veamos por qué.

La \textbf{física} \index{física} tiene como objetivo describir la naturaleza de forma sencilla, elegante e irreducible. La \textbf{física matemática} \index{física matemática} además tiene que fabricar una plataforma desde la cual sea normal hacerse preguntas como: Por qué la naturaleza es como es? La construcción de tales plataformas es algo muy tedioso pero dado el enorme trabajo de las generaciones pasadas es mucho lo que podemos adelantar invirtiendo un esfuerzo razonable.

La pregunta básica de la física - matemática es la siguiente: cuáles son los elementos fundamentales que forman la materia y cómo son sus interacciones? Se ha buscado intensamente la respuesta a esta pregunta y hoy por hoy muchos creen que la respuesta ha de buscarse dentro de las teorías gauge. Desde el punto de vista clásico, los ingredientes más destacados de estas teorías son: cambios de coordenadas, fibrados principales, conexiones y curvaturas. Aunque nos llevará tiempo entender todo eso, hay un elemento de las teorías gauge que es fácil de asimilar y que sin embargo es muy importante:

Partimos de la consideración de  que una teoría física debe por sobre todas las cosas describir un \index{universo objetivo} \textbf{universo objetivo} y cognoscible, es decir, que se pueda observar por diferentes investigadores y que entre ellos haya mutuo acuerdo. A esto le vamos a llamar el principio del realismo. Para implementarlo, las teorías deben tener en cuenta que toda medición implica sistemas de coordenadas arbitrarios, los cuales son impuestos por el investigador. Por consiguiente, la teoría física debe tener estabilidad ante cambios de coordenadas. Dichas coordenadas incluyen el espacio-tiempo y también otras coordenadas que se requieran para representar los estados internos de los sistemas que se estudien.

Más modernamente, se han creado métodos en geometría diferencial que permiten una formulación intrínseca de la física, pues toma verdades que son objetivamente ciertas y cuya formulación  no hace referencia a ningún sistema de coordenadas en particular.

El principio del realismo, estabilidad ante todo cambio de coordenadas, debe ser observado por toda teoría. Oficialmente, a este principio se le llama  covarianza generalizada. Algunos piensan que este principio  es inútil, pues uno puede describir un objeto desde un sistema de coordenadas y después utilizar las reglas reconocidas para cambiar de sistema de referencia, y así evitar toda contradicción. Con todo, nosotros veremos este principio a toda hora y en todo lugar. Que baste un ejemplo:

\

 Consideremos  el movimiento de una piedra que hemos lanzado: cómo podremos referirnos a su posición y velocidad si no es por medio de un  sistema de referencia que de por sí es arbitrario? Como hemos dicho, tenemos entonces que descubrir leyes matemáticas que describan la misma ley física en  todos los sistemas de referencia. O tal vez, formulaciones intrínsecas que no hagan referencias a ningún sistema coordenado en particular. También podemos descubrir que  la energía total de la piedra se conserva. Dicha energía representa en nuestra discusión el espacio interno, con dos componentes: energía cinética y energía potencial. A pesar de haber ya definido un sistema de coordenadas,  aún se tiene la libertad de fijar una altura arbitraria como el nivel cero de  la energía potencial, adjudicando, por ejemplo, dicho nivel a la superficie de la tierra.

Fijar el nivel cero de la energía potencial es una arbitrariedad introducida por la maquinaria matemática, que equivale a sumar una constante arbitraria a la energía potencial. El principio del realismo exige que dicha arbitrariedad no modifique en modo alguno las leyes físicas que describen el movimiento de los objetos. Eso se cumple automáticamente debido a que las leyes de la física se reducen en el presente caso a la segunda ley de Newton en un campo de fuerzas conservativo,

$$\vec F = ma = -\nabla \phi$$

por lo tanto la fuerza $\vec F $ sólo depende del gradiente de la energía potencial, $\phi$, que implica derivadas y que aniquila toda constante $k$:

$$ \nabla (\phi + k) =\nabla \phi $$

Notemos que el conjunto de constantes arbitrarias que podemos sumar a la energía potencial forman un conjunto, el de los números reales. Este conjunto es un ejemplo de lo que se denomina \index{grupo} \textbf{grupo}, concepto que representa un conjunto cuyos elementos se pueden operar, con propiedades semejantes a las de la suma ordinaria de los números reales. Puede parecer  algo tonto, pero vamos a decir que nuestra descripción de la dinámica de una partícula admite a los reales como grupo gauge. Esto insinúa que una teoría gauge contiene un grupo de arbitrariedades, \index{grupo gauge} \textbf{el grupo gauge}, que uno puede fijar a gusto ante las cuales las predicciones físicas no cambian. Lo que no es nada obvio, y lo veremos claramente para el electromagnetismo, es que el grupo de arbitrariedades de una interacción formulada al estilo gauge es el corazón mismo de la teoría.

\

El estudio del electromagnetismo es ciertamente más complicado que el de una piedra que vuela hasta caer, pero  el principio del realismo también estará ahí. Nos será fácil demostrar que la formulación clásica, dada por las leyes de Maxwell, tiene un grupo gauge, pero hacer lo equivalente para la versión cuántica será toda una obra de arte.

Al dejar la interacción electromagnética para pasar a interacciones más complicadas (que no veremos) tendremos que reconocer que el filtro del realismo es demasiado burdo y no da ideas para proponer la forma de dichas interacciones, es decir, no produce métodos específicos para hacer predicciones. Es en ese momento   cuando comenzamos a apreciar el poder de los otros conceptos de las \index{teoría gauge} \textbf{teorías gauge}: fibrados principales, conexiones y curvaturas, los cuales podemos considerar como un legado de la relatividad general. Todo esto lo veremos en detalle pero enfocados hacia la interacción electromagnética.

\section{CONCLUSION}

Para que una teoría matemática tenga futuro como teoría científica, es decir, como teoría que represente algún aspecto del universo físico, se requiere que las predicciones de la teoría sean estables ante toda arbitrariedad nacida sea de la observación, de la descripción o del tratamiento matemático. En particular, debe cumplir con el principio de covarianza generalizada, la estabilidad ante todo cambio de coordenadas, y con la invariancia gauge, la estabilidad de las predicciones ante todo grupo de arbitrariedades introducido por la formalización matemática, el grupo gauge. Las modernas teorías gauge de la física fundamental cumplen todo esto pero además involucran elementos serios de la geometría:  fibrados principales, conexiones y curvaturas, lo cual veremos a lo largo del curso.

\chapter{LA ACCION }

Teniendo en claro que la \index{mecánica clásica} \textbf{mecánica clásica} no explica todo, debemos preguntarnos por qué funciona. La manera más fácil de responder a ese interrogante es presentar un formalismo del cual pueda deducirse cómodamente las leyes de la mecánica clásica y que nos lleve directamente al corazón de la mecánica cuántica. Aunque tenemos varias propuestas, la del principio de la mínima acción  de Hamilton es sin duda alguna una de las más sencillas, versátiles y poderosas.

\

Pensemos en una piedra que vuela libremente después de que la tiramos hacia la guayaba. Ella describe una trayectoria. Preguntar por qué la mecánica clásica funciona es lo mismo que preguntarse por qué la piedra sigue el camino que siguió. Vamos a desarrollar un antropomorfismo para captar la forma como  Hamilton nos responde dicho interrogante: la piedra sigue el camino que siguió porque dicho camino era el que a la naturaleza le salía el más barato posible. Todos los caminos son posibles, pero cada camino tiene un peaje y el que la naturaleza escoge es aquel cuyo peaje total sea mínimo.

En realidad, no se requiere que sea mínimo sino simplemente extremo local  sobre el espacio de todos los caminos. Quizá hubiésemos deseado que la naturaleza funcionase como nosotros, que prefiriera aquella carretera, que en igualdad de condiciones, tuviese  el peaje mínimo. Sin embargo, no somos tan exigentes: un cuerpo está en equilibrio bien sea sobre un valle  o sobre una cresta. Por tanto, no hay razón para exigir del peaje que tome un mínimo. Por eso, sólo necesitamos un extremo. Claro que si el peaje total sólo tiene un extremo, entonces uno puede arreglar las cosas para que ese extremo sea un mínimo.

Un peaje es dictaminado por un alcalde o sus delegados. Si nosotros no hemos cuestionado la posibilidad de cambiar de alcalde, es porque ese alcalde es irremovible: está dado por la estructura del espacio-tiempo y por las interacciones fundamentales. A la estructura del espacio-tiempo se le acostumbra a identificar con la gravitación, pues equivale a una curvatura. El estudio general de sistemas de coordenadas sobre espacios curvos se denomina geometría diferencial y su propósito es formular leyes intrínsecas, o sea que tengan la misma forma en todos los sistemas de coordenadas pero que se puedan bajar directamente a cada sistema de coordenadas en particular.

\

Volviendo al principio de la mínima acción, el precio del peaje por unidad de camino se llama el Lagrangiano $L$, que es una función escalar (en general es real, pero para ciertos efectos puede ser complejo). Dado un camino, el peaje total será la integral  del Lagrangiano a lo largo de todo el camino. Al peaje total a lo largo del camino $\gamma$ se le llama acción, $S (\gamma)$. El camino elegido por la mecánica clásica produce un extremo de la acción. Y ese es el principio de Hamilton.

\

Cuál es la diferencia de la mecánica clásica con la mecánica cuántica? La mecánica clásica produce un camino, uno solo, aquel que da un valor extremo a la acción, y los demás los desprecia. En contraste, la mecánica cuántica toma en cuenta todos los caminos y un enjambre de   caminos parecidos da  un aporte particular al peaje total.

Por qué la mecánica clásica es tan exitosa siendo que no es fundamental? Porque en los casos donde acierta es porque de todos los caminos el que más importa es el tomado por la mecánica clásica y los aportes de todos los demás caminos se cancelan entre ellos.

Ahora viene la implementación de estas ideas intuitivas.

Comenzamos enfrentando el siguiente contraste: se sabe que las leyes de Newton son obedecidas por la naturaleza, al menos dentro de ciertas escalas. Dichas leyes están representadas por ecuaciones diferenciales. Puesto que para calcular una derivada tan sólo importa el comportamiento infinitesimal de las variables en juego, nosotros decimos que las leyes de la mecánica son locales. Es una felicidad que así sean, o de lo contrario nos resultaría imposible hacer un experimento pues no podríamos controlar todo el universo al mismo tiempo. Pero por otro lado, el principio de Hamilton toma toda una trayectoria como un todo indivisible y, como si fuese poco, analiza todas las trayectorias al tiempo y elige aquella que produce un peaje total extremo. Evidentemente, el criterio de Hamilton no es local.

Hay alguna contradicción? Debemos deducir del principio de Hamilton una ecuación diferencial ordinaria para rescatar la localidad. Hagámoslo. Pongamos el \index{Lagrangiano}  \textbf{Lagrangiano} o peaje por unidad de longitud de camino como:

$$L = T-U
\addtocounter{ecu}{1}   \hspace{4cm} (\theecu )      $$

donde $T$ es la energía cinética y $U$ la potencial. Para un camino $\gamma$ ya parametrizado con respecto al tiempo, desde tiempo inicial $t \sb{0}$ a tiempo final $t \sb{1}$, definimos la acción, o peaje total, $S$ como

$$S(\gamma) = \int \sb{t \sb{0}}\sp{t\sb{1}} L dt
\addtocounter{ecu}{1}   \hspace{4cm} (\theecu )      $$

Entendamos esta terminología y otros conceptos necesarios para seguir.

\

\

\color{blue}

$<<<<<<<<<<<<<<<<<<<<<<<<<<<<<<<<<<<<<<<<<$

\

Caminos que cuestan: El peaje como una \index{integral de línea} \textbf{integral de línea}, \index{funcionales} \textbf{funcionales}, extremos de funcionales.

\

Seguramente tenemos una imagen clara de lo que representa una integral común y corriente al estilo de

$\int senxdx $

Tratemos de ver dicha integral con otros ojos, con los de un chofer que tiene que andar una carretera.

Nuestro chofer va por una carretera que es una línea recta perfecta: la recta real. Y su carro lleva un taxímetro, que mide las unidades de carretera recorridas. En los taxis usuales los taxímetros le suman una cierta cantidad fija al costo total por cada unidad de camino recorrido. Para hacer que un taxímetro calcule la integral de la función $f(x)$ se necesita hacerle un cambio: cada vez que el taxímetro registre una nueva unidad de camino se suma al peaje acumulado la cantidad dada por $f(x)$. Por ejemplo, si llevamos 2 unidades de camino y acabamos de registrar la tercera, le sumamos al peaje acumulado $f(3)$. Si $f(x)=x\sp{2}$ entonces hay que sumar 9.

Lo que el taxímetro ha hecho se denomina técnicamente como suma de Riemann y nos da una aproximación a la integral. Para calcular la integral exactamente hay que hacer el  trabajo tedioso de considerar unidades de camino cada vez más pequeñas y calcular el   límite.  Si el límite existe, se  le llama la integral de Riemann de la función $f(x)$. Podemos garantizar que la integral de Riemann existe para funciones continuas o que sean discontinuas excepto en un número finito de puntos.

\

En la \index{integral de Riemann} \textbf{integral de Riemann} el punto de vista es el determinado por la carretera. Por otro lado, el taxímetro tiene que ser de aguja, esto es, tiene que poder medir en números reales, a precisión infinita.

\

Sucede que las carreteras no son líneas rectas perfectas.  Las carreteras son curvas, suben, bajan, dan la vuelta, regresan. Para nosotros no hay ningún problema: tomamos la carretera en línea recta y la enroscamos sobre la carretera que nos indiquen.

Es muy útil imaginarse que la carretera en línea recta es el tiempo. Cuando vamos enrollando el tiempo sobre la carretera, lo que hacemos es simplemente recorrerla: a veces vamos más rápido, a veces más despacio.

Al proceso de enroscamiento se le denomina \index{parametrización} \textbf{parametrización}:

$\gamma : R\rightarrow R\sp{3}$

Hay tantas parametrizaciones como carros y aún más: cada uno puede recorrer la carretera variando la velocidad muy caprichosamente.

Y cada carro lleva su taxímetro, el cual ha sido programado para calcular una integral dada, o un peaje, o el costo de una carrera o expreso, según ya hemos explicado: si se acaba de recorrer la unidad $x$, por esa unidad se agrega $f(x)$ al peaje ya acumulado. Después se refinan las unidades y se averigua el l\'imite. Si tal l\'imite existe, se denomina  integral de línea.

Los peajes pueden cobrarse de mil maneras posibles. Puede cobrarse una tarifa de una naranja por metro recorrido. Puede cobrarse más si es de noche. Puede cobrarse menos si la velocidad se mantiene entre 40km/h y 50 km/h. Puede duplicarse si la velocidad sobrepasa los 90 km/h.  O puede cobrarse según se haya gastado energía para hacer el expreso: eso fue lo que hicimos antes cuando hablábamos de trabajo para recorrer un camino. En fin, no hay límite a la creatividad para definir peajes.

\

Cuál es el peaje de la naturaleza desde el punto de vista de la aproximación clásica?

Al igual que un punto puede ser el mínimo de muchas curvas, de igual forma, un camino puede ser un mínimo de muchos peajes totales. Pero se ha encontrado que hay un peaje que resulta de maravilla: por unidad de camino se llama el Lagrangiano L, la energía cinética menos la energía potencial. Como la energía cinética involucra a la velocidad, el peaje de la naturaleza castiga fuertemente cualquier velocidad que ya de por sí es un exceso.

Al peaje total  se le llama \index{acción} \textbf{acción}, $S$ y es la integral de línea del peaje por unidad de camino:

$S(\gamma) = \int \sb{t \sb{0}}\sp{t\sb{1}} L dt$

Los límites indican que hemos tomado la parametrización natural, en la cual la carretera es el resultado de enroscar el tiempo, y por tanto el peaje se cobra desde el tiempo que inició la carrera hasta el tiempo en el cual se llegó al destino.

Si la parametrización es $\gamma (t) = \vec F(t)$, donde $\vec F$ es un vector de coordenadas, en general tenemos que $L=L(\vec F(t),\dot {\vec F}(t),t)$, lo cual quiere decir que el peaje se cobra de acuerdo al camino recorrido dependiendo del lugar por donde vaya y también de la velocidad a la cual transite, e incluso de la hora. Como quien dice, se cobra por todo.

\

Tenemos ahora que encontrar condiciones para los extremos de la acción. Para atacar un problema tan complicado, revisemos rápidamente el caso conocido, de hallar un extremo a funciones de R en R, cuando la función es diferenciable.

Para que una función de $R$ en $R$ sea diferenciable se requieren dos elementos: que uno pueda aproximar la función por la línea tangente, que define una transformación lineal, $df=mdx$,  y que el error sea super-despreciable con respecto a $dx$.

\

Definimos que una función $f(x) : R \rightarrow R$, que sea diferenciable, tiene un extremo en $x$ si la línea tangente tiene pendiente cero. Eso implica que en $x$ la diferencial cumple $df= f(x+dx)-f(x)= (lim_{ h \rightarrow 0})[(f+h)-f(x)]= mdx = 0$.

\

Similarmente, para estudiar  los extremos de la acción   tenemos que variar el camino un poquito y ver cómo cambia dicha acción: tendremos un extremo cuando la acción no cambia al deformar el camino un poquito. Formalmente:

\

\emph{Definición: La acción o funcional $S$ es diferenciable en $\gamma$ si   $S(\gamma +h)-S(\gamma)$ puede aproximarse por el efecto de un funcional lineal $ F(h) $  y si el error de aproximación es super-despreciable. El término $h$ denota una perturbación, es decir, un camino que adicionado al camino $\gamma$ lo perturba muy poquito, en particular, no le cambia los extremos. A $F$ se le llama la diferencial de $S$.}

\

Aclaramos que

\

1. Los caminos se pueden sumar y multiplicar por una constante pues denotan funciones de $R$ en el espacio vectorial $R\sp{3n}$. El $3$ es por las coordenadas espaciales y el $n$ porque el sistema en estudio puede tener n partículas. En general el $3n$ debe sobrentenderse y sólo describiremos una partícula que se mueve en un espacio de una dimensión, pero como el efecto es lineal sobre cada coordenada, se podrá extender inmediatamente al caso 3n.

\

2. Un funcional de caminos $S$ es lineal si es compatible con la suma y con la multiplicación escalar:

$S(\alpha + \beta) = S(\alpha) + S(\beta)$

$S(c\alpha)=cS(\alpha)$

donde $\alpha, \beta$, son caminos y $c$ es un número real.

Ejemplo: el funcional $S(\alpha) = \int\sb{0}\sp{1} g\sp{2}(\alpha(t))dt$, donde $g$ es una función que toma valores reales, no es lineal, en tanto que $S(\alpha) = \int\sb{0}\sp{1} 3g(\alpha(t))dt $ si lo es siempre y cuando $g$ sea lineal.

\

3. Decir que el error es super-despreciable implica un mundo de cosas. Primero hay que reconocer que los caminos, como funciones de $R$ en un espacio vectorial, en este caso $R\sp{n}$, forman un conjunto que acepta el concepto de norma o longitud que permite hablar de 'grande' y de  'peque\~{n}o', al igual que los números reales, o los vectores en $R\sp{n}$ que tienen su norma definida por el valor absoluto o por la distancia de la punta del vector hasta el origen.

Mejor dicho, los caminos forman un espacio normado (con concepto de longitud, $\|x\|$) que además lo hace métrico (con concepto de distancia, $d(x,y)= \|x-y\|$:

Si $\gamma$ es un camino $\gamma:R \rightarrow R\sp{n}$ entonces definimos la norma del camino $\| \gamma\|$ como

$\| \gamma\| = sup \| \gamma(t) \|$

En esa definición, $\gamma(t)$  es un vector y $\|\gamma(t) \|$ es su norma que es la misma longitud. La expresión $sup \| \gamma(t) \|$ dice que debemos tomar como norma del camino a su máxima desviación efectiva a partir del origen.

Por lo tanto, decir que $h$ es un camino que es una perturbación es lo mismo que decir que su norma es casi cero y que $h$ se mueve muy cerca del origen. Si el codominio de $h$ fuese $R$ entonces la gráfica de $h$ sería una curva pegada al eje 'X'.

Por otra parte, si $\gamma $ es un camino, un enroscamiento del tiempo, entonces por ser una función que toma valores en un espacio vectorial donde está definida la suma, la resta y la división por un escalar, pues también se le puede definir su derivada $\dot \gamma (t)$ que representa al vector velocidad del móvil que describe la parametrización dada por $\gamma $:

$\dot\gamma (t)= lim_{ h\rightarrow 0} (\gamma(t+h)-\gamma(t))/h$ donde $t \in R$.

\

Ahora  podemos ser exactos y formalizar el concepto de 'error superdespreciable' en el caso de un funcional $S$. Si tenemos que:

 $S(\gamma +h)-S(\gamma) = F(h)+ R(\gamma,h)$

y si $F$ es lineal en $h$, decimos que el error $R(\gamma,h)$ es super-despreciable si para $\|h\|<\epsilon$ y $\|\dot h\|<\epsilon$ se tiene que existe una constante $K$ positiva tal que $\mid R \mid   <K\epsilon\sp{2} $, cuando $\epsilon \rightarrow 0 $.

\

$>>>>>>>>>>>>>>>>>>>>>>>>>>>>>>>>>>>>>>>>>$

\

\

\color{black}

Para determinar los extremos de un funcional $S$ tenemos entonces que proceder así: se expresa el funcional en un camino cualquiera. Se perturba ligeramente el camino y se observa el cambio. Se extrae la parte lineal y se iguala a cero. Comencemos:

$S(\gamma + h) -S(\gamma )  $
$=\int \sb{t \sb{0}}\sp{t\sb{1}} [L(x+h,\dot{x}$
$ + \dot{h} ,t)dt - \int \sb{t \sb{0}}\sp{t\sb{1}}L(x,\dot{x},t)] dt$

Usando la linealidad de la integral nos da:

$S(\gamma + h) -S(\gamma )  $
$=\int \sb{t \sb{0}}\sp{t\sb{1}} [L(x+h,\dot{x}$
$ + \dot{h} ,t)-L(x,\dot{x},t)] dt$

Recordando la definición de diferencial de una función escalar, eso es simplemente

$=\int \sb{t \sb{0}}\sp{t\sb{1}} (dL)dt$

Pero $dL = dL(d\vec u)$ y en este caso $\vec u$ tiene dos grupos de componentes, los que corresponden  a la posición $\vec x$ y los que corresponden a su velocidad $\dot {\vec x}$. Como sólo anotamos una coordenada: $dL = dL(x,\dot x) = m\sb{x}dx + m\sb{\dot x} d\dot x$, donde $m\sb{x}$ y $m\sb{\dot x}$ son las pendientes en la dirección de $x$ y de $d\dot x$, o sea

$m\sb{x}= \partial L/\partial x $

$m\sb{\dot x}=\partial L/\partial \dot{x} $

denotando $dx$ como $h$ y $d\dot x$ como $\dot h$ nos queda:

$dL = (\partial L/\partial x )h $
$+ (\partial L/\partial \dot{x} )\dot {h}$

Por tanto

$S(\gamma + h) -S(\gamma )  $
$=\int \sb{t \sb{0}}\sp{t\sb{1}}[(\partial L/\partial x )h $
$+ (\partial L/\partial \dot{x} )\dot {h} ]dt  $

que captura el efecto lineal de la perturbación. Para hallar el error hay que expandir hasta todos los órdenes, pero no lo hacemos pues al error   lo consideramos super-despreciable.

Usando la integración por partes tenemos:

$\int \sb{t \sb{0}}\sp{t\sb{1}}$
$(\partial L/\partial \dot{x} )\dot {h}dt $
$=(h \partial L/\partial \dot{x})$
$ - \int \sb{t \sb{0}}\sp{t\sb{1}}h d/dt$
$(\partial L/\partial \dot{x} ) dt$

El término libre $(h \partial L/\partial \dot{x})$ se evalúa entre $t \sb{0}$  y $t\sb{1}$, pero como $h$ es una perturbación que deja los extremos iguales, entonces $h(t\sb{1})=h(t\sb{1})=0$. Eso implica que el término libre no da nada. Simplificando:

$S(\gamma + h) -S(\gamma )  $
$= \int \sb{t \sb{0}}\sp{t\sb{1}}[(\partial L/\partial x )h $
$-h d/dt(\partial L/\partial \dot{x} )]dt$

$= \int \sb{t \sb{0}}\sp{t\sb{1}}h[(\partial L/\partial x ) $
$- d/dt(\partial L/\partial \dot{x} )]dt$

Para que $\gamma $ sea un extremo, esta diferencia debe ser cero y lo debe ser para todo $h$, de lo cual se deduce que el integrando libre de $h$ debe ser cero:

$[(\partial L/\partial x ) - d/dt(\partial L/\partial \dot{x} )]=0$

A esta ecuación se le llama la ecuación de Euler-Lagrange.
\

\

\color{blue}

$<<<<<<<<<<<<<<<<<<<<<<<<<<<<<<<<<<<<<<<<<$

\

No tan de  prisa: en realidad no es tan obvio que el integrando libre de $h$ deba ser cero. Por eso formalizamos un poco mejor:

Si una función continua f(t), definida en el domino [a,b] satisface $\int \sp{b} \sb{a} f(t)h(t)dt = 0$ para toda función continua $h(t)$ tal que $h(a) = h(b) = 0$, entonces f(t) = 0 en todo el intervalo [a,b].

Prueba por contradicción: supongamos que $f(t)$ no sea idénticamente nula en el intervalo dado. Eso significa que existe algún punto $c$ tal que $f(c) \not= 0$. Puede ser negativa o positiva. Sin pérdida de generalidad podemos suponer que es positiva en $c$. Pero como $f$ es continua, entonces es positiva en una vecindad $U$ de $c$. Aún más: existe un intervalo $V$ contenido en $U$ tal que sobre $V$ la función $f$ es mayor que una cierta constante, $K$, tal vez muy cercana a cero.

Como el enunciado es cierto para toda $h$, escogemos una $h$ que nos guste, una que no anule el producto $fh$. Es suficiente con que sea continua, siempre positiva, valga uno en la vecindad $U$ y que se anule en los extremos del intervalo.

Sea B la longitud de V. Teniendo en cuenta que $fh>0$ y que sobre $V$ , $h$ vale 1, tenemos:

$\int \sp{b} \sb{a} f(t)h(t)dt > \int \sb{V} f(t)dt >  \int \sb{V} Kdt> BK>0$

\

Hemos llegado a una contradicción, pues partimos de que dicha integral debía ser cero: $f$ tiene que ser idénticamente nula.

\

$>>>>>>>>>>>>>>>>>>>>>>>>>>>>>>>>>>>>>>>>>$

\

\

\color{black}

Hemos demostrado entonces que el camino que produce un extremo en la acción necesita cumplir la \index{ecuación de Euler-Lagrange } \textbf{ecuación de Euler-Lagrange}:

$$d/dt(\partial L/\partial \dot{x} )-\partial L/\partial x  =0 \addtocounter{ecu}{1}   \hspace{4cm} (\theecu )      $$

En una dimensión, la ecuación se lee directamente, pero en dimensión $3n$, tenemos $3n$ ecuaciones, todas del mismo tipo. Estas  ecuaciones son diferenciales, o sea que hemos recuperado la naturaleza local de la mecánica clásica. Podemos además  recuperar la segunda ley de Newton. Veamos:

Para todos los casos, tomamos el Lagrangiano como la energía cinética $T$ menos la energía potencial $U$. Consideremos  una partícula que se mueve bajo el potencial, $V(x)$. El Lagrangiano propuesto es $L=(m/2)(\dot x\sp{2}\sb{1} + \dot x\sp{2}\sb{2} +\dot x\sp{2}\sb{3}) - V(x)$. Tenemos tres ecuaciones de Euler-Lagrange:

$d/dt(\partial L/\partial \dot{x\sb{i}} )-(\partial L/\partial x\sb{i} ) $
$= (m/2) d/dt(\dot x\sp{2}\sb{i}) -(-\partial V(x)/\partial x\sb{i} )=0 $

$= m \ddot x \sb{i}+\partial V(x)/\partial x\sb{i}= 0 $

Hemos deducido que

$Fuerza\sb{i} = m \ddot x \sb{i}= -\partial V/ \partial x \sb{i}$, como debía esperarse. Evidentemente, la forma del Lagrangiano se ha escogido  para que esto pueda hacerse de manera natural a partir de consideraciones energéticas.

\section{INVENTANDO LAGRANGIANOS}

El formalismo lagrangiano se ha llegado a convertir en una importante puerta de entrada a la gran física de interacciones fundamentales. \textbf{Inventar lagrangianos} \index{inventar lagrangianos} es un arte complicado del cual no se sabe mucho. La idea del grupo gauge se considera interesante al respecto. Por ahora,  vamos a hacer un ejercicio de inducción para ilustrar una forma como se inventan lagrangianos: teniendo la ley de Newton, vamos a deducir la forma del Lagrangiano y del principio variacional del cual dicha ley podría deducirse.

Consideremos un sistema mecánico clásico descrito por $n$ variables y sus derivadas $q_1,.. ..,q_n, \dot{q_1},.. ..,\dot{q_1}$. La ley de Newton dice

$m_i\ddot q_i = -\partial V/\partial q_i$

Multiplicando por $dq_i $ y sumando sobre todas las coordenadas

$\sum m_i\ddot q_i  dq_i  = -\sum (\partial V/\partial q_i) dq_i = -dV $

$\sum m_i\ddot q_i  \delta q_i  + dV  = 0$

Ahora integramos sobre el tiempo entre $r$ y $s$ pero con la condición de que en dichos tiempos no se haya ninguna variación : $dq_i(r)=dq_i(s)=0$. Tenemos:

$\int^s_r (\sum m_i\ddot q_i  \delta q_i  + \delta V)dt  = 0$

Introducimos $\delta$ en vez de $d$ para diferenciar entre las variables independientes de las dependientes.

Ahora integramos por partes el primer término y el segundo lo dejamos de lado por un momento. Tomamos como $u=\delta q_i$ y como  $dv=\ddot q_i$. Por tanto, $v= \dot q_i$, pero necesitamos saber qué es $du$. Tengamos en cuenta que

$\delta(dq/dt)=d(q+\delta q)/dt - dq/dt = d/dt(q+\delta q-q)=d(\delta q)/dt$

por tanto  a $du/dt$ lo tomamos como $d(\delta q)/dt = \delta(dq/dt)$. Entonces

$\int^s_r (\sum m_i\ddot q_i  \delta q_i) dt = \sum m_i \delta q_i \dot q_i |^s_r -  \int^s_r (\sum m_i\dot q_i  \delta \dot{q}_i) dt $

utilizando las condiciones de frontera se aniquila el primer término y sólo queda el segundo, el cual al ser reemplazado en

$\int^s_r (\sum m_i\ddot q_i  \delta q_i  + \delta V)dt  = 0$

nos queda en

$-\int^s_r (\sum m_i\dot q_i  \delta \dot{q}_i  + \delta V)dt  = 0$

multiplicando por -1 y maquillando tenemos

$\int^s_r \delta(\sum m_i  (1/2) \dot{q}^2_i -   V)dt  = 0$

$\delta \int^s_r (\sum m_i  (1/2) \dot{q}^2_i -   V)dt  = 0$

$\delta \int^s_r L dt  = 0$

donde $L=\sum m_i  (1/2) \dot{q}^2_i -   V$

Hemos demostrado que si la ley de Newton es válida para un sistema, entonces dicho sistema funciona como si ejecutara un valor crítico para una integral que hemos denominado la acción. El argumento es reversible siempre y cuando se tenga en cuenta que el $\delta q_i$ pueda ser cualquiera con las condiciones de frontera adecuadas.

\section{LEYES DE CONSERVACION}

Si uno sabe que el \textbf{Lagrangiano es invariante} \index{Lagrangiano invariante } ante el efecto de un cierto grupo de transformaciones, uno puede deducir una ley de conservación. Veamos algunos ejemplos, el primero muy intuitivo.

Imaginemos que vamos por un territorio el cual es muy uniforme. Podemos inferir que como nos fue en el recorrido pasado así nos irá en el futuro, por ejemplo, si antes recorrimos 30 km en un día, en los siguientes también recorreremos 30 km por día. Decir que el territorio es muy uniforme es asegurar que los descriptores del ambiente son invariantes ante traslaciones, las cuales forman un grupo. Por lo tanto, el recorrido diario se conserva.

Si el Lagrangiano es invariante ante un translación cualquiera a lo largo de un determinado eje, entonces $\partial L/\partial x\sb{i} = -\partial V(x)/\partial x\sb{i}= 0$ y por tanto, $m \ddot x \sb{i}=0$ lo que implica que $m \dot x \sb{i} =k$, o sea que el momento lineal en esa dirección se conserva. También se puede probar que si el Lagrangiano es invariante en el tiempo, entonces la energía total debe conservarse. Veamos esto al detalle.

\section{CONSERVACION DE LA ENERGIA}

Averigüemos bajo qué condición en mecánica clásica se conserva la energía. Tomemos un \textbf{sistema mecánico conservativo} \index{sistema mecánico conservativo} en el cual la fuerza motriz sea independiente del tiempo.

La fuerza, $\vec F$, se define como $\vec F = d \vec p/dt $ donde $\vec p$ es el momento lineal $\vec p= m \vec v$. Esta definición de fuerza también es válida para cuerpos de masa variable, como los cohetes. Pero si consideramos que la masa no cambia, entonces $d \vec p/dt = d (m\vec v)/dt= md \vec v/dt=m \vec a$. Se dice que el sistema es mecánico si $\vec F= \vec F(\vec q,\vec p,t)$, donde $\vec q$ es la posición. La función $\vec F$ debe ser continua como función de  sus variables. La fuerza puede depender del momento lineal, y de la velocidad, como sucede con la fuerza de Lorentz de una partícula cargada en un campo magnético. Decimos que el campo es conservativo si la fuerza viene de un potencial $V$ : $\vec F = d \vec p/dt = - \nabla V$.

Probemos ahora que existe una propiedad que se conserva, es decir, una observable, una función que describe alguna propiedad del sistema, que no depende del tiempo y hallemos su forma. Trabajaremos en una dimensión, pero la generalización a mas dimensiones es inmediata.

Sea $\phi = \phi(q,p)$ la función que queremos descubrir y que no depende del tiempo, es decir, que es constante en el tiempo, que se conserva. Por tanto:

$d\phi /dt = 0 $
$= (\partial \phi/\partial q) (d q/d t) $
$+  (\partial \phi/\partial p) (d p/d t)$

Como
$ dq/dt=v=p/m$ y además $d p/d t$
$= dp/dt= - \nabla V $, entonces:

$d\phi /dt = 0 = (\partial \phi/\partial q) (p/m) $
$+  (\partial \phi/\partial p) (-\nabla V)$

Que en una dimensión implica:

$(\partial \phi/\partial q) (p/m) $
$=  (\partial \phi/\partial p) (dV/dq)$

Separando variables, $q$ a un lado y $p$ al otro:

$(\partial \phi/\partial q) /(dV/dq)^{-1}$
$= (\partial \phi/\partial p) /(p/m)^{-1}$

Y nos damos cuenta que hay una igualdad entre dos funciones, una función sólo de $q$ y la otra función sólo depende de $p$. Tal igualdad puede ser posible sólo si las funciones son iguales a una constante, que podemos tomar como uno, fijando una escala de medida:

$(\partial \phi/\partial q) /(dV/dq)$
$= (\partial \phi/\partial p) /(p/m) =k = 1$

(Cuando tenemos muchas coordenadas, en el paso anterior hay que dejar las variables que corresponden a un sola coordenada a un lado de la ecuación y las demás juntas al otro lado y se aplica lo mismo, coordenada por coordenada).

De la primera ecuación se deduce que:

$\partial \phi/\partial q = dV/dq \Rightarrow  \phi $
$= V + g(p) \Rightarrow \partial \phi/\partial p= dg/dp $

Entre tanto, de la segunda ecuación tenemos:

$\partial \phi/\partial p= p/m = dg/dp \Rightarrow g$
$= p\sp{2}/2m + k = mv\sp{2}/2+ k $

Por consiguiente $\phi= V + g(p) = V + mv\sp{2}/2 + k $, donde $k$ is un número real.

Hemos descubierto la conservación de la energía total. Mientras que los componentes cinético o potencial pueden cambiar, su suma total tiene que permanecer invariante. Qué fue lo que necesitamos para demostrarlo? Tuvimos necesidad de un campo conservativo que no cambiase con el tiempo. Hay una forma de decir lo mismo que es muy fructífera:

\

Si las leyes que rigen un sistema son invariante en el tiempo, entonces existe una cantidad que se conserva, la cual denominamos energía. Esa forma de hablar la debemos a \textbf{Emmy Noether} \index{Emmy Noether} quien descubrió en 1918 un teorema que en forma general dice: \emph{si el Lagrangiano de un sistema es invariante ante un grupo de transformaciones, eso implica que una cierta observable se conserva}. En este caso,  la invariancia respecto al tiempo implica que existe una cierta observable, que llamamos la energía, que se conserva en el tiempo. Similarmente, anteriormente vimos que si un sistema es invariante a la rotación, entonces conservaba el momento angular.

Al definir la energía total como $\phi= V + g(p) = V + mv\sp{2}/2 + k $ constatamos que nace una arbitrariedad, dada por la constante $k$. Dicha arbitrariedad no puede hacerse parte de la energía cinética $mv\sp{2}/2 $, puesto que la velocidad es una observable objetiva que no admite más que variaciones en la escala de medida. Pero dicha arbitrariedad queda bien en la energía potencial, que no puede observarse. Cumple ese arreglo el test gauge o de calibración? Por supuesto que sí: los resultados observables se determinan por la segunda ley de Newton, fuerza = masa x aceleración, y la fuerza está ligada al potencial mediante una derivada y por lo tanto, las constantes no importan: $\vec F = - \nabla V= - \nabla (V+k)$.

Pudimos probar entonces que una arbitrariedad en los modelos matemáticos para nada interfiere con el contenido físico de dichos modelos. El modelo mecánico conservativo que hemos establecido cumple el principio de calibración. Hay una arbitrariedad, sumar una constante real arbitraria. Dichas constantes forman un grupo, el de los números reales, $R$. Aunque no se acostumbra, podemos decir que la teoría de los sistemas mecánicos conservativos es una teoría gauge con grupo gauge o de invariancia gauge igual a $R$.

\section{CONSERVACION DEL MOMENTO ANGULAR}

\

La \textbf{formulación} \index{formulación intrínseca} de la mecánica clásica basada en la acción es \textbf{intrínseca}, es decir, las ecuaciones de Euler-Lagrange no cambian de forma aunque cambiemos de coordenadas. Esto es algo espectacular. Eso se debe a que en su desarrollo se usaron propiedades algebraicas más no hubo referencia a ningún eje, ni recto ni curvo. Utilicemos esta propiedad  para calcular, por ejemplo, la dinámica de una partícula sometida a rotación uniforme. En ese caso, el Lagrangiano se escribe en coordenadas polares y las ecuaciones Euler-Lagrange son las mismas, de lo cual se deduce la ley de la conservación del momento angular, siempre y cuando la energía potencial sólo dependa de la distancia al origen más no del ángulo polar. Veamos:

Sea $\vec e\sb{r} = cos \theta \vec i + sen \theta \vec j$ el vector unitario en la dirección de $\vec r$ y  $\vec e \sb{\theta} =\dot {\vec e}\sb{r} =-sen \theta \vec i + cos \theta \vec j$ el vector unitario en la dirección angular. Es decir que mientras que $\vec e\sb{r}$ marca la dirección radial, $\vec e \sb{\theta} $ nos indica la dirección perpendicular al radio, tangencial a toda circunferencia.

Como el vector posición se descompone en $\vec{r} = r \vec e\sb{r}$, entonces,

$\dot {\vec r} = d/dt(r \vec e\sb{r}) =  \dot{r} \vec e\sb{r} $
$+ r \dot {\vec e}\sb{r} =$
$ = \dot{r} \vec e\sb{r} + r  \vec e\sb{\theta} \dot \theta $.

Teniendo en cuenta que  $\vec e\sb{r}$ y $\vec e \sb{\theta}$ son ortogonales y de norma uno, la energía cinética toma el valor

$T=(1/2)m (\dot {\vec r} \cdot \dot {\vec r)} $
$=(1/2)m (\dot{r} \vec e\sb{r} + r  \vec e\sb{\theta} \dot \theta) \cdot (\dot{r} \vec e\sb{r} + r  \vec e\sb{\theta} \dot \theta) $
$=(1/2)m(\dot{r} \sp{2} + r \sp{2}  \dot \theta \sp{2}).$

 Por lo tanto, el Lagrangiano queda

 $L=T-U=(1/2)m(\dot{r} \sp{2} $
$+ r \sp{2}  \dot \theta \sp{2})-U $

donde la energía potencial depende sólo de la posición. Observemos que el Lagrangiano es insensible al valor del ángulo polar. Tenemos una teoría con una arbitrariedad en la definición del ángulo de referencia. En cambio el valor absoluto del radio si es tremendamente importante. Probemos que la insensibilidad a las rotaciones fijas  implica que hay  una cantidad que se conserva. Para ello, usemos las ecuaciones de Euler-Lagrange:

$d/dt(\partial L/\partial \dot{x}\sb{i} )$
$-(\partial L/\partial x \sb{i}) =0$

La variable que nos interesa es la velocidad angular, $\dot \theta $, por lo tanto especificamos la ecuación correspondiente a dicha variable:

$d/dt(\partial ((1/2)m(\dot{r} \sp{2} $
$+ r \sp{2}  \dot \theta \sp{2})-U  ) /\partial \dot{\theta} )$
$-(\partial ((1/2)m(\dot{r} \sp{2} $
$+ r \sp{2}  \dot \theta \sp{2})-U  ) /\partial \theta) =0$

$ md/dt(r\sp{2}\dot \theta) + \partial U/\partial \theta = 0$

Si además la energía potencial no depende de $\theta$ tenemos que

$ d/dt(m r\sp{2}\dot \theta) =0$

La última ecuación nos dice que hay una cantidad conservada, a la cual se le llama momento angular. Gracias a esta ley, las balletistas pueden acelerar su giro si comienzan a girar con los brazos extendidos y después recogen los brazos. El mismo efecto se usa en la regulación de la velocidad angular de cuerpos en rotación: se le ponen unos brazos. Si se acelera mucho la rotación, los brazos se levantan absorbiendo energía y manteniendo la rotación uniforme o incluso se podrían acondicionar para que regulen la fuente de energía. La dinámica completa se calcula adicionando  las ecuaciones Euler-Lagrange correspondientes a $r$ y a su derivada.

\section{EL TEOREMA DE EMMY NOETHER }

Vimos que  si las leyes que gobiernan la dinámica de una partícula eran  invariantes ante rotaciones, entonces eso implicaba la conservación del momento angular. En efecto, la energía cinética es en ese caso proporcional a la norma de un vector y por tanto es invariante a rotaciones y la energía potencial depende tan solo de la distancia al centro y no del ángulo polar, y por consiguiente también es invariante ante rotaciones. Eso nos permitió descubrir  una cantidad conservada, el momento angular. Veamos ahora que eso no es una coincidencia: eso tiene que ser así por necesidad, lo cual fue probado por Emmy Noether alrededor de 1918:  cada vez que el Lagrangiano es invariante ante un grupo de transformaciones, debe existir una cantidad conservada.

Emmy Noether fue discípula de Felix Klein quien hizo de su vida una cruzada para enseñarle al mundo la importancia de la teoría de grupos tanto en matemáticas como en física. Mientras que haya hombres se recordará a Emmy Noether por haber formulado un teorema que convierte algo intuitivamente obvio en una clara realidad matemática capaz de iluminar toda la física. Hay varias versiones del teorema de Noether: en mecánica clásica, en mecánica cuántica, en teoría de campos y en general con ecuaciones de Lie. Nosotros veremos la demostración para un caso particular de la mecánica clásica y para la mecánica cuántica.

\bigskip

\addtocounter{ecu}{1}

\textit{Definición \theecu}
\addtocounter{ecu}{1} La física de un sistema es invariante ante la acción de un grupo si la  derivada del lagrangiano con respecto a cada uno de los grados de libertad del grupo es cero. Por ejemplo, si el grupo es el de las traslaciones en el plano, un grado de libertad es el de las traslaciones a lo largo del eje X.

\bigskip

\textit{Teorema de Noether \theecu} \emph{Si el lagrangiano de un sistema es invariante ante la acción de un grupo, entonces el sistema tiene al menos una cantidad que se conserva.}

Demostración:  Consideremos la acción del grupo a lo largo de uno de sus grados de libertad, al cual lo cuantificamos como $s$ denotando, por ejemplo, la cantidad de kms que se traslada el sistema observado, y al lagrangiano lo notamos $L$.

 Partimos de un  lagrangiano con forma $L=L( q, \dot q) $, que quiere decir que el peaje de la naturaleza depende del punto y de la velocidad en el punto. Como la acción del grupo no perturba la física tenemos que $ \partial L/\partial s = 0$. Pero eso también puede usando la regla de la cadena asumiendo que $L=L( q(s), \dot q(s)) $.  Denotando la deriva con respecto a $t$ con un punto y la deriva con respecto a $s$ con una prima, la regla de la cadena da:

$ \partial L/\partial s = 0 = (\partial L/\partial q)(q')+ (\partial L/\partial \dot q)\partial \dot{q'}$

pero como la física del sistema debe ser un extremal de la acción, entonces deben cumplirse las ecuaciones de Euler-Lagrange:

$$d/dt(\partial L/\partial \dot{q} )-\partial L/\partial q  =0 \addtocounter{ecu}{1}        $$

o sea

$$d/dt(\partial L/\partial \dot{q} )=\partial L/\partial q  \addtocounter{ecu}{1}        $$

Si reemplazamos este resultado en la ecuación que dio la invariancia del lagrangiano ante la acción del grupo, tenemos:

$ 0 = (\partial L/\partial q)(q')+ (\partial L/\partial \dot q)\partial \dot{q'}= $
$d/dt(\partial L/\partial \dot{q} )(q')+ (\partial L/\partial \dot q)\partial \dot{q'}$

y notando que nos queda la deriva de un producto, reescribimos esta ecuación como:

$d/dt((\partial L/\partial \dot{q}) q')=0$

Integrando nos queda:

$(\partial L/\partial \dot{q}) q'=k $

que da a partir del Lagrangiano la cantidad conservada.

Comparemos esta respuesta general con la ya conocida para el caso en el cual es lagrangiano es invariante ante las rotaciones. Sabemos que el momento angular se conserva, o sea $ d/dt(m r\sp{2}\dot \theta) =0$. Veamos que tenemos aquí. Nuestro lagrangiano es

 $L=T-U=(1/2)m(\dot{r} \sp{2} $
$+ r \sp{2}  \dot \theta \sp{2})-U $

y la cantidad que el teorema de Noether predice que ha de conservarse es

$(\partial L/\partial \dot{q}) q'=k $

en la cual tomamos $q = \theta$ y el grupo es el de las rotaciones, con grado de libertad, $\theta$, por lo que $q' = d\theta/d\theta = 1$. Por otro lado $(\partial L/\partial \dot{q})$ se convierte en

$\partial L/\partial \dot{\theta}= m r^2 \dot \theta$

por lo que la cantidad conservada es

$(\partial L/\partial \dot{q}) q'=(\partial L/\partial \dot{\theta}) d\theta/d\theta= m r^2 \dot \theta= k $

\

Nuestra demostración es particular en el sentido en que hemos utilizado a $R^n$ como el espacio base. En el caso general, hay que tomar una variedad cualesquiera para representar sistemas con restricciones arbitrarias (holónomas). Mas luego veremos el teorema de Noether en mecánica cuántica. También hemos asumido que el grupo de invariancia admite el concepto de derivación, o sea que estamos tratando con un grupo de Lie. Por eso nuestra demostración no es válida para reflexiones del espacio o para inversiones del tiempo. ¿Se podrá extender el teorema a dichos casos?

\chapter{LA MECANICA CUANTICA}

De dónde salió la \index{mecánica cuántica} \textbf{mecánica cuántica}? De la incapacidad de la mecánica clásica para explicar experimentos como el  llamado \index{efecto fotoeléctrico} \textbf{efecto fotoeléctrico}, el cual es el siguiente: dentro de un tubo al vacío se coloca una placa metálica y al frente una rejilla positiva. Si por cualquier razón, digamos por calentamiento o por estimulación con luz, algunos electrones  abandonan la placa, éstos son acelerados hacia la rejilla positiva que está al frente. Algunos quedarán presos al chocarse contra ella, otros pasarán de largo y podrán ser registrados en su camino.

Alumbramos la placa metálica fría y observamos que hay un efecto umbral que depende de la frecuencia de la onda de luz con que se irradia:

Antes de una cierta frecuencia no se observa evaporación de electrones, pero después sí. Este efecto umbral es imposible imaginar desde el punto de vista clásico, pues como hemos visto en el teorema de Poynting, toda onda conlleva energía, y por tanto, si con poca intensidad no se alcanza a hacer algo, debería ser suficiente aumentar la intensidad adecuadamente para lograrlo. Tal efecto no se observa (a intensidades moderadas pero que deberían ser suficientes).

\

Esto ya es motivo de sobra para formular una modificación a la teoría clásica  de la interacción entre el campo electromagnético y la materia. La teoría cuántica ha demostrado ser una modificación al fin de cuentas sencilla pero  muy poderosa. Veamos cómo se introduce su formalismo matemático.

\section{LA INTEGRAL DE FEYNMAN}

El camino parametrizado que describe la evolución de un sistema en mecánica clásica es aquel que hace extrema la acción, peaje total o tarifa del taxi lagrangiano. Sólo ese camino importa. En mecánica cuántica las cosas suceden de otro modo: todos los caminos importan.

De dónde nace esta idea tan extra\~{n}a? Nace de la óptica: si uno tiene un fuente de luz que ilumina una pantalla e interpone una tabla con dos rendijas, la luz iluminará la pantalla  después de haber pasado por ambas rendijas al tiempo.

Uno podría decir: la luz pasa por ambas rendijas al tiempo, pues un poco de luz pasa por una rendija y otro poco, distinto del primero, pasa por la otra. Tratando de entender lo que significa eso de 'un poco de luz', uno podría pensar que se trata de una porción de onda lumínica e inventar un mundo entero al respecto. La sencillez científica depurada por la historia nos muestra otra alternativa:

En primer término, la luz dejó de ser considerada una simple onda: Einstein explicó los elementos cuantitativos del efecto fotoeléctrico asumiendo que la luz viene en gránulos, llamados fotones, cuantos de energía, y que cada fotón se choca con un electrón haciéndolo que brinque.

Alegamos: los fotones, siendo gránulos minúsculos de energía, pasan unos por una rendija y otros por la otra.

Para dilucidar eso, se baja la intensidad de luz  hasta cuando la luz conste de gránulos muy aislados a ver qué pasa. Por supuesto que eso requiere métodos sofisticados para medir la intensidad de luz.

Un método para registrar bajas intensidades de luz utiliza a los fotomultiplicadores, unos tubos inmensos que operaban según el siguiente propósito: si la luz logra hacer que un fotón se evapore de una placa metálica, entonces al ser acelerado por un campo eléctrico adquirirá energía cinética. Se interpone en su camino otra placa metálica. Cuando el electrón evaporado se choca contra la placa, sale un chispero de electrones que caen bajo la acción del campo eléctrico acelerante. Ese juego se repite varias veces y lo que antes era un sólo electrón se convierte en una avalancha. De esa forma tenemos una manera de medir intensidades de luz supremamente bajas.

Gracias a los fototransistores hoy en día ya se tienen diseños que involucran paneles que en todo se parecen a la retina de un ojo: una superficie tapizada de elementos fotosensibles e independientes.

Se acumulan los registros de luz que recibe cada elemento y se gráfica  la densidad total en cada uno de tres casos: cuando hay una sola de las dos rendijas abiertas y cuando están ambas.

Decimos que la luz es de baja intensidad cuando el panel de elementos fotosensibles es estimulado de cuando en cuando y cada vez en uno solo de sus elementos fotosensibles.

Qué pasa con los cuantos de luz? Acaso pasan ellos por una sola rendija?

\

Si los cuantos de luz pasaran sólo por una rendija, entonces debemos esperar que la intensidad registrada con las dos rendijas abiertas fuese simplemente la suma de las intensidades registradas al dejar las rendijas abiertas por separado.

\

El resultado experimental es contundente: lo que se observa a mediana intensidad también se observa a baja intensidad: la intensidad de luz que llega a un punto con las dos rendijas abiertas no es la suma de las intensidades de luz que deja pasar cada rendija por separado. Decimos que hay interferencia.

Esto es válido a intensidad mediana, registrable con una placa de fotografía y también es válido a muy baja intensidad, registrable con retinas de estado sólido.

Debemos concluir que los cuantos de luz, o fotones, interfieren consigo mismo, para lo cual se exige que cada fotón pase por ambas rendijas al tiempo cuando las dos están abiertas. Hablamos entonces de la dualidad onda-partícula de la luz: avanza como onda pero se registra como partícula.

\

Por una sugerencia que Louis Victor de Broglie propuso en 1924, este experimento de óptica también se ha llevado a cabo con electrones y el resultado es el mismo: los electrones se comportan de igual forma que la luz: se crea interferencia. Esto lo lograron Clinton Joseph Davisson, Lester Harbert Germer y George Paget Thomson: hicieron pasar electrones por entre las rendijas interatómicas de la red de un cristal y obtuvieron un espectro en todo semejante al que se observa si uno hace pasar rayos X por el mismo cristal, o  luz visible sobre rendijas creadas por carpintería.

Me parece muy pedagógico el cine mudo presentado por la enciclopedia Encarta bajo el título mecánica cuántica (quantum mechanics) que ilustra estos aspectos. Es digno de verse.

\

Ahora viene la genialidad pedagógica de \index{Feynman} \textbf{Feynman}: Qué es un medio transparente? Es un medio totalmente opaco al cual primero se le abrieron dos rendijas, lo cual creó dos caminos, ambos importantes. Y después vinieron otras dos rendijas, que crearon otros dos caminos que son tan importantes como los dos primeros. Siguiendo este proceso vemos que un medio transparente es un medio opaco totalmente invadido de rendijas. Y por consiguiente un medio transparente es un mundo de caminos todos importantes.

Cómo se cuantifica el aporte de los caminos? Una idea algo aventurada de Dirac, publicada en 1933, fue la clave para que Feynman inventara su punto de vista,  y lo trabajara hasta hacerlo cuantitativo y en armonía con lo que ya se había descubierto, produciendo su versión de la mecánica cuántica, publicada en 1948, la cual entramos a explicar.

\color{blue}

\

\

$<<<<<<<<<<<<<<<<<<<<<<<<<<<<<<<<<<<<<<<<<$

\

Los complejos:

\

\

Los \index{números complejos}  \textbf{números complejos} $C$: en la matemática hay grandes obras de ingeniería. El conjunto de los números complejos es una de ellas. Esos números fueron diseñados para que resistieran el embate de la siguiente pregunta: cuál es aquel número cuyo cuadrado es menos uno?

Los números reales no fueron capaces de resistir tanta presión: todo cuadrado es positivo pero $-1$ es negativo.

Hay que diseñar entonces un conjunto cuyos elementos sean  números: que se puedan sumar, restar, multiplicar, dividir, sacar raíz cuarta y logaritmo y hallar el seno y el coseno y todo eso. Pero que además permita a los números negativos tener raíces.

Resulta que tal construcción existe y es simplemente el plano armado de una estructura aritmética. Cada punto del plano se direcciona con dos coordenadas $(x,y)= x + iy$ , que representan la distancia al eje vertical 'Y' desde el origen y al eje horizontal 'X', respectivamente pero con su signo. El símbolo $x$ nos dice que el número complejo está a $x$ unidades en la dirección horizontal, mientras que $iy$ declara que hay que avanzar en la dirección vertical $y$ unidades.

La suma se define coordenada por coordenada. Los complejos con la suma forman un grupo. Como los reales se identifican con el eje 'X' o sea los números de la forma $(x,0)$, y también forman un grupo con su suma que corresponde a la misma suma de todos los días, decimos que los reales son un subgrupo de los complejos.

La multiplicación se ve nítida en coordenadas polares. Dichas coordenadas definen un punto del plano por dos valores: la longitud del radio vector que parte del origen,   el complejo (0,0),y llega  al punto dado, que se llama el módulo o norma, y el ángulo polar que dicho radiovector forma con la parte positiva del eje 'X' en dirección contraria a las manecillas del reloj. Debemos advertir que el ángulo polar no se puede definir de manera única: dado un ángulo polar, si se le suma cualquier número de vueltas, o sea $2n\pi$, con $n \in Z$, tenemos el mismo número complejo. Como veremos, la elección del semieje positivo $X\sp{+}$ como el origen del ángulo polar es la causa de la mayor parte del material del presente trabajo.

Existe una notación de los números complejos que dan inmediatamente su estructura multiplicativa: si $z \in C$, es decir, si $z$ es un punto del plano complejo, entonces $z=re\sp{i\theta}$ donde $r$ es el modulo de $z$ y $\theta$ es el ángulo polar.

Si tenemos dos complejos

$z= re\sp{i\theta}$

 $w=se\sp{i\phi}$

entonces la multiplicación simplemente multiplica los módulos y suma los ángulos:

$zw=re\sp{i\theta} se\sp{i\phi}=rse\sp{i(\theta+\phi)}$

Se dice que un número es unitario si tiene módulo uno:

$z= e\sp{i\theta}$

Multipliquemos un número cualquiera $w=se\sp{i\phi}$ por uno unitario $z= e\sp{i\theta}$:

$zw=e\sp{i\theta} se\sp{i\phi}=se\sp{i(\theta+\phi)}$

lo que pasa es que el número sufre una rotación. Por supuesto que el radio no varía. Como veremos eso es clave.

El módulo de un complejo se lee directamente en su expresión polar pero también puede recuperarse mediante el siguiente truco:

Definimos el conjugado de $z=re\sp{i\theta}$ como $z\sp{*}=re\sp{-i\theta}$, es decir que tiene el mismo módulo pero el ángulo está medido  al revés, en dirección negativa, en la misma dirección que las manecillas del reloj. Ahora bien:

$zz\sp{*}=re\sp{i\theta}re\sp{-i\theta}=r\sp{2}e\sp{i(\theta-\theta)} = r\sp{2}e\sp{i(0)} $

el resultado que nos da es un número complejo que tiene norma $r\sp{2}$ y ángulo polar cero. Cuál es ese número? como tiene ángulo polar cero está sobre la parte positiva del eje 'X' por lo que es un número real y como tiene módulo $r\sp{2}$, pues es en definitiva el número real $r\sp{2}$. Por consiguiente, el módulo cuadrado, largo cuadrado o norma cuadrado de un complejo es $zz\sp{*}$. Lo recordaremos.

Los complejos con la multiplicación no forman un grupo pues el cero no tiene inverso multiplicativo. Pero si quitamos el cero este problema se arregla.

\

$>>>>>>>>>>>>>>>>>>>>>>>>>>>>>>>>>>>>>>>>>$

\

\

\color{black}

\section{LA AMPLITUD}

La forma como la luz se difracta se explica asumiendo que ella es una onda. Como la materia también posee la posibilidad de tener interferencia, entonces la materia también tiene características de onda. Por otro lado, los números complejos tienen ondas en su infraestructura según lo demostró Euler:

$ z= r e \sp{i\theta} = rcos \theta + i rsen \theta$

Por consiguiente, empecemos asociando un número complejo a cada proceso físico. (Lo podemos hacer como un ensayo a ver qué pasa, o podemos hacerlo para poder deducir el principio de incertidumbre, porque de otra forma no se puede). A dicho número complejo lo llamaremos amplitud. Si un proceso puede ocurrir por dos canales alternativos, por el primero con amplitud $z$ y por el segundo con amplitud $w$, entonces, y éste es el postulado fundamental de la mecánica cuántica, al permitir los dos canales se tendrá como amplitud la suma de las dos primeras, $z + w$.

\

Hay un problema relacionado con el tema de este trabajo: un número complejo tiene un módulo y una fase o argumento. Sucede que la fase puede definirse tan sólo si ya se ha definido de antemano un sistema de coordenadas: la fase del número complejo $z$ se define como el ángulo entre $z$ y el eje generado por el número complejo $1$. Acaso el $1$ es más natural que $1+4i$ o que $7-8i$? Si cambiamos de generador,  cambiamos de sistema, cambiamos de fase. Por lo tanto, la fase, y por consiguiente la amplitud, no puede tener ningún significado físico, pues consideramos que la materia existe objetivamente y que su comportamiento puede describirse independientemente de los sistemas de coordenadas que utilicemos para direccionarla. La fase está descrita por un número complejo de norma uno, y esos números forman un grupo de Lie: el $\textbf{U(1)}$.

\

\color{blue}
$<<<<<<<<<<<<<<<<<<<<<<<<<<<<<<<<<<<<<<<<<$

\

\

Glosario:   Espacio topológico, variedad, grupo de Lie, simetría, grupo de invariancia, grupo de invariancia gauge.

 \textbf{Espacio topológico:} \index{espacio topológico} es un conjunto en el cual se puede hablar de abiertos y cerrados. Generaliza los conceptos de intervalo abierto y cerrado de la línea real, o de bolas abiertas o cerradas del plano: uniones arbitrarias de abiertos da un abierto, e intersecciones finitas de abiertos da abierto. Un espacio topológico es un conjunto al cual se le adjudica una familia de subconjuntos la cual es cerrada a uniones arbitrarias y a intersecciones finitas.

\textbf{Variedad: } \index{variedad} Es la abstracción de la tierra  desde el punto de vita de los cartógrafos. Cada región se pinta fielmente sobre un mapa o carta, y para cartografiar toda la tierra hay que hacer muchos mapas solapándolos unos con otros hasta que quede bien cubierta pero cuidando de que haya compatibilidad de unos mapas con otros.

Se acostumbra a de decir: una \textbf{variedad} es un conjunto localmente isomorfo a un abierto de $R\sp{n}$. Localmente significa que se puede decir sólo en un pedacito o carta. Isomorfo significa que en cierto modo no hay diferencia entre lo uno y lo otro (entre la superficie de la tierra y el mapa local). Para definir qué es un abierto, primero definimos qué es una esfera abierta: Decimos que una esfera es abierta cuando no tiene el cascarón, borde o frontera. Se dice que un conjunto es abierto cuando cada punto puede rodearse con una esfera  abierta totalmente contenida en el conjunto (si estamos en $R\sp{2}$ por supuesto que hablaremos de círculo. Y en $R$ hablaremos de intervalo, una esfera de dimensión 1). Al estudio de los abiertos y de sus propiedades se llama topología. Los físicos matemáticos prefieren definir la topología como el estudio de los huecos, como los de una dona.  Se dice que la variedad es diferenciable cuando es suave. Con el concepto de variedad se inaugura la geometría diferencial. 

Ejemplos: La recta real $R$  es una variedad pues es un esfera abierta e  infinita. En general cualquier abierto de $R\sp{n}$ es una variedad trivial. La esfera del espacio tridimensional puede recubrirse con unos seis pedazos de periódico, por lo tanto también es una variedad.

Formalmente tenemos:

Al subconjunto de $R\sp{n}$ denotado como $B\sb{r}(\vec a)=\{\vec x:\|\vec x-\vec a\|<r \}$ se le denomina la bola abierta de radio $r$ centrada en $\vec a$. Se dice que un conjunto $U$ de $R\sp{n}$ es abierto si cada punto $\vec x$ de $U$ posee alguna  bola abierta centrada en $\vec x$ totalmente contenida en $U$. Se comprende, entonces, por qué a los conjunto abiertos se les dibuja sin frontera, pues no la tienen. Decimos que un conjunto es cerrado cuando su complemento es abierto.

Los conjuntos abiertos tienen las siguientes propiedades:

0. Tanto el conjunto vacío como todo  $R\sp{n}$ son abiertos.

1. Una intersección finita de abiertos es un conjunto abierto. Si consideramos $I\sb{n} = [-1/n,1/n]$ entonces tenemos que la intersección de todos esos conjuntos, que son cerrados, da el punto $0$, el cual no es abierto.

2. Cualquier unión arbitraria de conjuntos abiertos da un conjunto abierto.

Dado un conjunto cualquiera $X$ , no vacío, para el cual se define una familia de subconjuntos $\tau $de $X$ que cumpla las tres propiedades anteriores se denomina una topología. Por ejemplo, si $X=\{a,b,c\}$ y definimos $\tau$ como $\tau = \{\O,X\}$ entonces $\tau$ es una topología sobre $X$ y decimos que $(X,\tau)$ es un espacio topológico. Naturalmente que un conjunto dado puede tener muchas topologías incompatibles entre ellas (un conjunto es abierto en una topología pero no en la otra). A los elementos que son miembros de la topología se les llama abiertos y si contienen a un punto dado $x$ también se les llama vecindades de $x$.

Decimos que una función $f:X\rightarrow Y$ es continua en $x \in U$ cuando $(X,\tau\sb{X})$, $(Y,\tau\sb{Y})$ son  espacios topológicos y cuando para toda vecindad $V$ de $ f(x)$ puede encontrarse una vecindad $U$ de $x$ tal que $f(U)\subset V$. La continuidad puede imaginarse como un juego de control, de estabilidad: hay dos jugadores, el operario y el auditor. El operario está en $X$ y el auditor en $Y$. El auditor propone una vecindad $V$ de $f(a)$, la que desee. En respuesta, el operario debe fabricar una vecindad de operación $U$ tal que al moverse dentro de dicha vecindad no desestabiliza el sistema afuera de la restricción dada por el auditor. Naturalmente que el auditor querrá poner restricciones muy severas y si el operario siempre puede ajustarse es porque el sistema responde ''continuamente'' a los cambios del operario.

Camino a definir lo que es una variedad, volvamos a recordar que la tierra es cartografiable, es decir, que puede dibujarse sobre un pedazo de papel los accidentes que cada región pueda tener y que toda la tierra requiere muchos mapas compatibles entre ellos. Empezamos con  un espacio topológico  $(X,\tau)$, entonces tomamos una familia de abiertos de $X$ que recubre a todo $X$ y que define un conjunto de  regiones cartografiables, definimos por cada región $U\sb{i}$ el protocolo $\phi\sb{i}$ que define su mapa correspondiente $\phi\sb{i}(U\sb{i} )$ . El conjunto de mapas  $(U\sb{i},\phi\sb{i})$ forma un atlas. Decir que existe un protocolo de mapeo es lo mismo que exigir que se cumpla que para cada  región $U\sb{i}$ existe funciones continuas $ \phi\sb{i}, \phi\sb{i}\sp{-1}$ tales que:

$\phi\sb{i}: U\sb{i}\rightarrow \phi\sb{i}(U\sb{i}) \subset R\sp{n}$ (esto significa que hay una forma de dibujar a la región $U\sb{i}$

$\phi\sb{i}\sp{-1}: \phi\sb{i}(U\sb{i})\rightarrow U\sb{i} \subset X $ (esto implica que el mapa resultante es una fiel copia de la región).

Es natural que exijamos que en los lugares en los cuales hay solapamiento de mapas exista plena compatibilidad entre las informaciones que cada uno aporta. Por eso exigimos que para ser variedad se requiera que:

la función de transición dada por $\phi\sb{i} \circ \phi\sb{j}\sp{-1} : \phi\sb{j}(U\sb{i} \cap U\sb{j}) \rightarrow \phi \sb{i} (U\sb{i} \cap U\sb{j})$ es un \textbf{difeomorfismo} \index{difeomorfismo}  de clase $C\sp{k}$ (tiene inversa y tanto ella como su inversa  es  de clase $C\sp{k}$).

Advertencia Gauge: el principio gauge que nosotros estamos desarrollando dice que existe una realidad estable ante todo elemento subjetivo. Desde siempre se ha considerado que el objeto que uno desea cartografiar es real y que los mapas o cartas que se usan para su descripción son subjetivos. Por consiguiente falta exigir que dos cartógrafos puedan ponerse de acuerdo sobre lo que ellos hablan guiados por sus mapas. Lo corriente es que cada cartógrafo produzca mapas muy particulares quizás bien distintos de los producidos por todos los demás. Uno de ellos puede utilizar coordenadas rectangulares y otro polares, etc. Lo que se necesita es que ellos puedan encontrar un lenguaje de traducción entre sus mapas. Eso se garantiza con la siguiente condición que es lo menos que se le puede pedir a dos  mapas para que sean compatibles:

Dos atlas son equivalentes ssi entre ellos hay plena compatibilidad, es decir si al hacer un solo atlas  usando todos los mapas al tiempo hay una forma de ponerse de acuerdo. Concretamente: Dos atlas de clase $C\sp{k}$, $(U\sb{i},\phi\sb{i})$ $(V\sb{j},\psi\sb{j})$ son equivalentes ssi para cualquier par de protocolos de mapeo $\phi\sb{i}$ y $ \psi\sb{j}$ se tiene que  $\phi\sb{i} \circ \psi\sb{j}\sp{-1} $ es un difeomorfismo de $\psi\sb{j}(U\sb{i}\cap V\sb{j}) $ sobre $\phi\sb{i}(U\sb{i} \cap V\sb{j})$.

Como cada atlas puede refinarse todo lo que se quiera adicionando nuevos atlas, por tanto definimos un atlas maximal $M$ generado por un atlas dado $(U\sb{i},\phi\sb{i})$ como aquel que contiene a todos los refinamientos de dicho atlas. Surge una pregunta delicada: Todos los atlas dan origen al mismo atlas maximal? Por ahora contentémonos con poner de presente que nuestro principio gauge exige que las propiedades objetivas sean independientes tanto de los mapas, como de sus refinamientos como de las distintas cosmovisiones que puedan dar origen a atlas maximales incompatibles.

Ahora si podemos definir: decimos que $X$ es una variedad cuando es un espacio topológico y está provisto de un atlas maximal que la describe.

Las variedades se hacen necesarias para simplificar el estudio de sistemas con restricciones, digamos un tren que debe moverse sobre una carrilera, o un cuerpo compuesto cuyas partes deben permanecer unidas aunque con posición relativa variable.

Veamos algunas razones de por qué el círculo unitario $S\sp{1}= \{(x,y) \in R\sp{2}: x\sp{2} + y\sp{2}\}$ es una variedad de clase $C\sp{\infty}$ es decir que es de clase $C\sp{k}$ para todo $k$. En efecto: el círculo es un espacio topológico donde los abiertos de su topología son los abiertos del plano intersecados con el $S\sp{1}$. Un atlas consta de cuatro mapas generados por un protocolo de cartografía que se denomina proyección estereográfica.

Simplemente se toma cada uno de los polos, norte y sur, oriente y occidente y desde cada polo se traza una línea recta que atraviesa tanto al círculo como al eje de coordenadas que le es transversal en donde produce la imagen del círculo sobre $R$. Concretamente, tomamos como regiones a cartografiar la parte superior del círculo, $U\sb{1}$, en tanto que la parte inferior del  círculo   es $U\sb{2}$. De otra parte, hay que tomar también la parte derecha y la parte izquierda.

La proyección estereográfica sobre $U\sb{1}$ es $\phi\sb{1}(x,y) = x/(1-y)$ y la que va sobre $U\sb{2}$ es $\phi\sb{2}(x,y) = x/(1+y)$, lo cual se verifica por triángulos semejantes. De igual forma se construye las otras dos funciones de mapeo. De la geometría del problema se deduce que cada función es invertible. Después hay que verificar que las funciones de transición son continuas e invertibles y además que sobre las intersecciones todo es infinitamente derivable.

\

\textbf{Grupo de Lie} \index{Grupo de Lie} : es un conjunto que tiene tanto la  estructura de grupo como la de variedad  diferenciable  con completa compatibilidad entre las dos. La compatibilidad implica que la operación de grupo es continua y muy suave, lo mismo la asignación del inverso. El adjetivo 'Lie' honra a Sophus Lie, un matemático que tuvo mucho que enseñarnos referente al uso de la teoría de grupos aplicada a la solución de las ecuaciones diferenciales. Ejemplos:

a. $R\sp{n}$: la estructura de grupo está provista por la suma natural. La de variedad por las esferas abiertas. Decimos que la suma es continua porque si $\vec u + \vec v = \vec w$, y si se fija un límite a las deformaciones de $\vec w$ entonces se puede deformar levemente tanto a $\vec u$ como a $\vec v$ y si se hace la suma, el nuevo resultado no sobrepasará los límites preestablecidos. Similarmente, si el inverso de  $\vec u$ es  $-\vec u$ y si se fija un límite a las deformaciones de $-\vec u$ entonces se puede deformar levemente a $\vec u$ de tal manera que el inverso aditivo de dicha deformación  no sobrepase los límites preestablecidos.

b. El conjunto de las traslaciones en el plano: tienen estructura de grupo con operación binaria la composición o sea la aplicación de traslaciones en serie. Para todos los efectos, el grupo de traslaciones es isomorfo al mismo plano, pues cada vector del plano genera un translación y viceversa. De ahí se deduce que la traslaciones tienen estructura de variedad y de grupo de Lie puesto que el plano es $R\sp{2}$, posee ambas estructuras.

c. Las rotaciones  en el plano forman un conjunto que es un grupo de Lie. La operación es la composición. La estructura de variedad  no es igual  a la de $S\sp{1}$ pero se le parece, pues el círculo unitario es isomorfo al subconjunto de rotaciones de ángulo menor o igual a $2\pi$. Para ver la diferencia, consideremos dos definición: decimos que un subconjunto de un espacio topológico es compacto cuando todo recubrimiento abierto tiene un recubrimiento finito.

En $R\sp{n}$ un conjunto es compacto ssi es cerrado y acotado ssi toda sucesión contiene una subsucesión convergente. De esto se desprende que el círculo unitario es compacto, pues toda sucesión en dicho círculo es una sucesión en la bola cerrada unitaria, la cual es cerrada y acotada. Por lo tanto debe haber una subsucesión convergente a un elemento de dicha bola cerrada. Pero el límite de dicha subsucesión debe estar en el círculo. En cambio $R$ no es compacto pues se puede recubrir con intervalos abiertos de largo $2$ que se solapan en una unidad, pero dicho recubrimiento deja al descubierto al menos un punto si se le quita cualquiera de sus intervalos. Por lo tanto no tiene ningún subrecubrimiento finito y no es compacto. Similarmente, la sucesión $\{1,2,3,4,.....\}$ no tiene ninguna subsucesión convergente.

Como el conjunto de rotaciones en el plano es isomorfo a $R$, entonces no es compacto y por eso dicho conjunto no tiene una estructura de variedad igual a la de $S\sp{1}$.

d. \textbf{U(1)} \index{\textbf{U(1)}} ( se lee u-uno): a este grupo de Lie lo vamos a encontrar en la presente introducción a toda hora. Así que procuremos entenderlo lo mejor que se pueda:

Hemos visto que $R$ es un grupo: es la recta real,  con la suma natural y la estructura de variedad determinada por los intervalos abiertos. Ahora bien, imaginemos a $R$  como un lazo  infinito y enrosquémoslo alrededor del círculo unitario, conservando tanto la distancia entre números como las propiedades de la suma. Es evidente que el enroscamiento no le quita ni le pone nada a la estructura de grupo de Lie que los reales tenían antes del enroscamiento. A esta nueva estructura se le llama \textbf{U(1)} no compacto, naturalmente.

Para producir la estructura \textbf{U(1)} como grupo de Lie compacto se necesita simplemente darle a $S\sp{1}$ una estructura algebraica, la heredada de $\textbf{U(1)}$ no compacto con la salvedad de que  $2\pi$ se iguala a cero. \textbf{U(1)} compacto puede tener versiones con cualquier número de vueltas y sigue siendo compacto. Los expertos consideran que la naturaleza ve una diferencia entre las múltiples versiones de \textbf{U(1)}.

Los elementos de  \textbf{U(1)} tienen varias representaciones. Una muy usual es la exponencial:

$z \in  \textbf{U(1)}$ ssi existe un número real $\theta$ tal que $z= e\sp{i\theta}$. Por supuesto que $\theta $ es el ángulo del radio vector que determina a $z$. Para pasar de la notación exponencial a la cartesiana se utiliza  una expresión debida a Euler:

$e\sp{i\theta}= cos\theta + isen\theta$.

\

\textbf{Simetría} \index{Simetría}: decimos que un objeto es simétrico cuando al sufrir una reflexión queda igual que antes, es decir, invariante. Una reflexión dada que se aplica dos veces da la identidad y entre las dos forman un grupo, el grupo generado por la reflexión. Por lo tanto la definición oficial de simétrico es: un cuerpo es simétrico cuando es invariante ante el grupo generado por una reflexión dada.

En física, y al parecer en arquitectura, la palabra simetría se ha liberado de referirse a una reflexión y se aplica a cualquier grupo que tenga sentido. Ejemplo:

Consideremos al plano $R\sp{2}$ originalmente blanco, pero lo pintamos de negro a 1 galón de pintura por $m\sp{2}$. No hay en un plano así punto de referencia alguno. Si corremos a cualquier lado se ve igualitico. O si lo rotamos. Pero si lo  expandimos, entonces el negro se tornará en gris. Se dice: el grupo de simetrías del plano pintado de negro incluye entre sus generadores tanto a las traslaciones como a las rotaciones como a las reflexiones, pero no a las magnificaciones o estiramientos. En general, el grupo de simetría de un objeto, así sea una teoría, está generado por todas las operaciones que dejan al objeto invariante desde un punto de vista prefijado.

\

El \textbf{Grupo de invariancia} \index{Grupo de invariancia} gauge es un grupo que representa una arbitrariedad matemática en la formulación de una ley física. Tomemos el caso de la energía, la cual ni se crea ni se destruye, sólo se transforma. La forma como lo hace implica que la energía potencial tiene una arbitrariedad en el nivel cero. Dicha arbitrariedad está representada por un número real. Dichos números forman un grupo que además es de Lie. Por tanto, podemos resumir: la teoría clásica de la transformación de la energía  es una teoría con grupo de invariancia gauge $R$, el conjunto de los números reales.

Explicaremos detalladamente de qué forma la teoría cuántica de la interacción electromagnética es una teoría gauge con grupo de invariancia \textbf{U(1)}, lo cual  significa que la teoría  predice que no habrá cambios físicos u observables si una cierta cantidad se multiplica por cualquier número de \textbf{U(1)}.

\

$>>>>>>>>>>>>>>>>>>>>>>>>>>>>>>>>>>>>>>>>>$

\

\

\color{black}

Hemos introducido la amplitud para descubrir que la fase no contiene física alguna. En realidad, la física si está contenida en la amplitud pero de forma oculta:

Recordemos  el efecto fotoeléctrico: si uno ilumina un metal con luz para estudiar si salen o no electrones, sucede un efecto escalón: hay una frecuencia de la luz tal que si se ilumina con luz de mayor frecuencia, entonces salen electrones, pero si se iluminan con luz de menor frecuencia entonces no salen electrones.

En general, los experimentos miden no un electrón sino infinidades de ellos y después se compara la intensidad relativa. Es lo que sucede en la pantalla de un televisor: de donde caen muchos electrones sale mucha luz y de donde caen pocos electrones sale poca luz. Un modelo adecuado para describir tal situación se logra introduciendo probabilidades.

Deseamos saber si la formulación probabilística de nuestra mecánica cuántica en ciernes   es una astucia matemática o un reflejo de una propiedad fundamental de la materia. Eso se averigua haciendo experimentos con luz a muy baja intensidad y utilizando fotoreceptores de muy alta sensibilidad tales que puedan captar un solo fotón. Se descubre que los fotoreceptores se encienden de forma aleatoria y concluimos que las leyes fundamentales de la luz son probabilísticas.

\

Cómo relacionamos nuestra amplitud con las probabilidades? Tenemos tres leyes muy sencillas:

\

a) Si un proceso se puede desarrollar por varios canales independientes y si por el canal $j$ la amplitud es $w\sb{j}$ (lo cual se averigua cerrando los otros canales), entonces al abrir todos los canales se tiene la amplitud $\sum w\sb{j}$ (esta suma también puede significar una integral).

\

b)Si la amplitud de que un proceso se desarrolle por determinado canal o conjunto de canales es $w$, la probabilidad de que el proceso ocurra  por dichos canales es el cuadrado del  módulo de la amplitud, $\|w\| \sp{2}$.

\

c) El concepto de amplitud debe ser compatible con el de probabilidad, es decir, cuando se tiene un proceso que consta de dos subprocesos en tándem, el primero con amplitud $w$ y el segundo con amplitud $z$, entonces el proceso total tiene amplitud $ws$.

\

 Estas recetas son de aplicación universal y en lo que sigue las desarrollaremos para partículas libres o en potenciales de fuerzas.

\section{LA MISMISIMA INTEGRAL}

La mecánica clásica se puede describir por el principio de la mínima acción de Hamilton: para ir de un estado a otro la naturaleza escoge una trayectoria que es extremal de la acción, la cual es la integral del Lagrangiano a lo largo de la trayectoria. Frecuentemente, dicho extremo también es un mínimo, pero no necesariamente, por ejemplo: en la órbita de un planeta uno toma dos puntos cualesquiera y por tanto la órbita queda dividida en dos tramos, uno produce un mínimo de la acción y el otro un máximo.

La mecánica cuántica es una extensión muy simple de la clásica si se mira desde el punto de vista de la acción: para ir de un estado a otro, ya no se escoge un camino, sino que cada camino posible se torna en un obligado canal de evolución. A cada canal se le asocia una amplitud y la amplitud total es la suma o integral sobre todos los caminos posibles. Por supuesto que hay diferencia de camino a camino y unos aportan una amplitud parcial y otros otras. Lo que resulta de todo eso se denomina como integral de Feynman.

\

Concretamente, consideremos una partícula que se mueve en una dimensión sujeta a un potencial $V(x)$. Su Lagrangiano o peaje por unidad de camino es la energía cinética menos la potencial:

 $L= (1/2)m v\sp{2} -V(x)$

La acción o peaje total para un camino $\gamma$ ya parametrizado con respecto al tiempo, desde tiempo inicial $t \sb{0}$ a tiempo final $t \sb{1}$ es :

$S(t \sb{0},t\sb{1}, \gamma) = \int \sb{t \sb{0}}\sp{t\sb{1}} L dt$
$=\int \sb{a}\sp{b} L dt = S(a,b, \gamma) $

Ahora no tomamos un extremo sino que cada camino $\gamma$ que lleva de $a$ a $b$ aporta una amplitud parcial $w$ dada por:

$w= e \sp{(i/\hbar)S(a,b,\gamma)}$

Todos los caminos que llevan de $a$  a $b$ producen una amplitud total $ K(b,a) $ dada por la integral de todas las amplitudes sobre todos los caminos:

$$ K(b,a)
= \int \sb{a} \sp{b} e\sp{(i/\hbar)S(a,b,\gamma)}D\gamma (t) \addtocounter{ecu}{1}   \hspace{4cm} (\theecu )      $$

Todo eso es un símbolo conocido como la \index{integral de Feynman} \textbf{integral de Feynman}, al cual se le ha podido dar una formulación totalmente rigurosa sólo en algunos casos. Nosotros tomaremos una aproximación heurística y en particular evitaremos hablar sobre la medida en el espacio de caminos, excepto que debe cumplir condiciones mínimas como por ejemplo que permita que la probabilidad total siempre sume uno.

\

\

\color{blue}

$<<<<<<<<<<<<<<<<<<<<<<<<<<<<<<<<<<<<<<<<<$

\

Taxímetros generalizados: medida e integración

\

Cuando analizábamos el concepto de integral de línea encontramos que había que programar un taxímetro para que adicionara un cierta cantidad por cada unidad de espacio recorrido. Para algunos tipos de integración el peaje total depende de dos factores independientes: una debido al carro, digamos el número de ejes, y otra debida a la carretera, que esté o no bien asfaltada. Se llama teoría de la medida al estudio general de estos tópicos, pero en particular al que tiene que ver con la carretera.

Podemos generalizar nuestra idea de integral para incluir mucho más que costo de expresos. En general, cualquier cantidad que pueda taxarse puede reformularse en términos de una integral. Nada más consideremos el punto de vista de los que tienen que asfaltar la carretera. Ellos no pensarían en una dimensión, como el taxista. Ellos pensarían en dos dimensiones, pues su trabajo se relaciona con áreas. O inclusive en tres dimensiones pues es más barato, más duradero y menos traumático, asfaltar una carretera con una capa de asfalto bien gruesa que asfaltarla dos veces con una capa de la mitad de espesor. Tendríamos entonces integrales en dos y tres dimensiones.

La forma más sencilla de taxar por el uso de la carretera es imponer una tarifa proporcional al camino recorrido. En dos dimensiones, lo correspondiente sería cobrar proporcionalmente al área, y en tres al volumen.

Una compañía encargada de la seguridad de un edificio cobraría, por ejemplo, por el tamaño del edificio, una integral de volumen, y por el tiempo que dure el servicio, es decir, su tarifa estaría dada por una integral en 4 dimensiones.

La energía total consumida por un ser humano involucra una integral de dimensión cuatro. En efecto: hay que alimentar el volumen del cuerpo de la persona, por el tiempo que ella viva.

Consideremos la energía lumínica consumida por un ecosistema de algas en un lago: tiene que ver con el volumen del lago, una integral triple, con el tiempo, otra dimensión más, y con el espectro de absorción: unas algas absorben en rojo y otras en azul. Cinco variables que producen una integral en cinco dimensiones.

\

Nosotros estamos trabajando con caminos. No precisamente con carreteras, sino con formas de viajar: escoger una carretera y recorrerla a la velocidad deseada. Debido a que en cada lugar y estado de viaje en donde uno se encuentre se puede tomar decisiones, eso quiere decir que cada lugar es una dimensión. El espacio de caminos es  infinito-dimensional. Por supuesto que definir la medida debe ser algo desafiante. Más luego veremos cómo enfrentamos ese problema.

\

$>>>>>>>>>>>>>>>>>>>>>>>>>>>>>>>>>>>>>>>>>$

\

\

\color{black}

\section{INVARIANCIA GAUGE}

Para nosotros, la palabra 'gauge' o 'calibración' está asociada a un principio de sentido común: las leyes de la física no pueden depender de las arbitrariedades que nazcan bien sea para describir el universo o en los formalismos matemáticos que describen la evolución de los sistemas. Hemos formulado la mecánica cuántica. Puede ella satisfacer el filtro gauge?

Comenzamos preguntándonos si es que acaso hemos dado lugar a alguna arbitrariedad dentro de nuestra formulación de la amplitud total $K(a,b)$. Pues, si. Si la hay: se trata de que en su puro corazón hay una exponencial imaginaria, o sea una fase. No es posible definir una fase o argumento de un número complejo a menos que uno tenga una sistema de coordenadas en el plano. Tal sistema es totalmente arbitrario y cada quien tiene derecho a escoger el sistema que le guste y por tanto su fase también cambiaría: cambiará la física? Debemos probar que ante un cambio de fase, la física no sufre ningún cambio. En efecto, si en un sistema tenemos que la amplitud $K(a,b)$ es:

$ K(b,a)$
$= \int \sb{a} \sp{b} e\sp{(i/\hbar)S(a,b,\gamma)}D\gamma (t) $

Entonces en otro sistema de coordenadas del plano complejo podríamos tener:

$ K \sb{n} (b,a)= $
$\int \sb{a} \sp{b} e\sp{(i/\hbar)[S(a,b,\gamma) + \phi ]}$
$D\gamma (t) $

$=\int \sb{a} \sp{b} e\sp{(i/\hbar)S(a,b,\gamma)}$
$e\sp{(i/\hbar)\phi }D\gamma (t)  $

$=e\sp{(i/\hbar)\phi } \int \sb{a} \sp{b} e\sp{(i/\hbar)S(a,b,\gamma)}D\gamma (t)  $

$=e\sp{(i/\hbar)\phi }K(b,a)$

Pero $K$  no es observable, lo que es observable es la probabilidad de transición de $a$ a $b$ que ella genera, la cual es $\|K\| \sp{2}$ y es claro que tal magnitud no depende de un corrimiento global de fase, pues multiplicar por un número complejo, $e\sp{(i/\hbar)\phi }$, lo único que causa es una rotación de ángulo $\phi$ y eso no causa estiramiento o encogimiento de ningún complejo. Por lo tanto, la probabilidad, que es la norma cuadrado, tampoco cambiará.

Vamos bien.

Podemos observar que dimos lugar al nacimiento de una arbitrariedad en la definición de la fase y que la física sin embargo no varía. El cambiar de sistemas de coordenados en el plano complejo implica que la amplitud se multiplica por una exponencial imaginaria

$$K \sb{n} (b,a)=e\sp{(i/\hbar)\phi }K(b,a)
\addtocounter{ecu}{1}   \hspace{4cm} (\theecu )      $$

Por lo tanto, decir que la mecánica cuántica pasa el filtro gauge es lo mismo que decir que las probabilidades de transición, cuantificadas por la norma cuadrado de la amplitud $K$, no varían ante la multiplicación de $K$ por una exponencial imaginaría. Esas exponenciales imaginarias forman un grupo con la operación multiplicación, que además es de Lie. La operación de grupo corresponde a la multiplicación de complejos que se interpreta como  hacer primero un cambio de fase y después otro. Ese grupo  es el $\textbf{U}(1)$, con la $\textbf{U}$ de unitario, que quiere decir que las probabilidades no cambian y la suma total sobre todos los canales da un 1.

\section{RECUPERANDO LA MECANICA CLASICA}

Le damos la bienvenida a la  mecánica cuántica porque explica cosas que la clásica no puede, a saber, la interferencia entre partículas. Pero queda entonces por explicar el misterio de por qué funciona la mecánica clásica sabiendo que no es una teoría fundamental. La razón es la siguiente: cada camino aporta conforme a la acción que éste define. Como son tantos caminos, cada uno de ellos aporta algo insignificante, o sea que lo importante no es lo que aporte cada uno sino la forma como entre ellos se ayuden  o interfieran. Ahora bien, si el camino es un extremo de la acción, una deformación leve, un camino ligeramente diferente,  producirá una acción que a primer orden no puede diferenciarse del camino extremal. Por tanto, todas las pequeñas perturbaciones del camino extremal produce perturbaciones constructivas, que refuerzan la acción del camino extremal: las fases $e \sp{(i/\hbar)S(a,b,\gamma)}$ interferirán constructivamente.

Pero una variación de una camino que no es un extremo produce una acción que si puede diferenciarse a primer orden de la acción del camino original. En consecuencia, las fases variaran alocadamente con sólo variar el camino un poquito. Eso se debe a que la acción se mide en unidades de $\hbar$, la cual es extremadamente peque\~{n}a: en términos ergonómicos tiene dimensiones del orden de $10 \sp{-27} (gramos \times cm \sp{2} /segundo)$. Por lo tanto, cuando el camino no produce un extremo de la acción, su aporte es aniquilado por el de los caminos cercanos, con los cuales  interferirá destructivamente y  no aportarán nada a la amplitud total $K(a,b)$.

En resumen, el camino que domina la situación es aquel que da un extremo de la acción. Hemos recuperado la mecánica clásica y hemos probado que la mecánica clásica es tan cuántica como lo más cuántico que exista. Sin embargo,   se dice que uno tiene un efecto cuántico  cuando uno no puede explicar los resultados sin apelar a muchos caminos totalmente diferentes. Por algún tiempo se pensó  que los sistemas macroscópicos no podían tener propiedades de coherencia, donde las fases se sumen a gran escala, pero ahora se sabe que eso no es cierto: un laser es el ejemplo perfecto: un mundo de fotones todos con la misma fase y recorriendo un camino que puede ser de aquí a la luna. Otro ejemplo es la superconducción, en la cual no se cumple la ley de Ohm y por consiguiente una corriente puede recorrer caminos circulares sin gastar energía y sin violar el teorema de Stokes. Más luego veremos los detalles de un caso específico: el efecto Aharonov-Bohm.

\section{LA ECUACION DE SCHR\"{O}EDINGER}

Para entender lo que significa una \index{integral de Feynman} \textbf{integral de Feynman} es necesario calcular alguna  y nosotros vamos a hacer un cálculo que tiene como objetivo básico responder a la siguiente inquietud: si todos los caminos de evolución importan, entonces debe ser que el futuro depende del pasado en toda su extensión. Pues bien, vamos no sólo a mostrar que en mecánica cuántica el futuro depende única y exclusivamente del presente, si es que éste existe, sino además calcularemos la dependencia que el futuro inmediato tiene del presente. Para hacer eso, deduciremos que la integral de Feynman implica una ecuación diferencial que se denomina la \index{ecuación de Schr\"{o}edinger} \textbf{ecuación de Schr\"{o}edinger}.

Que sirva de aclaración decir que una ecuación diferencial del tipo $\dot y=f(t,y)$ denota una ley o algoritmo cuya evolución sólo depende del presente porque si se sabe qué es $y$ en el tiempo $t$, entonces sabremos qué es $y$ en el tiempo
$t + \bigtriangleup t$, a saber:

$y(t + \bigtriangleup t) = y(t) + \dot y \bigtriangleup t $
$= y(t) + f(t,y) \bigtriangleup t $

Por eso buscamos una ecuación diferencial, para demostrar que el futuro sólo depende del presente.

\section{FUNCION DE ONDA}

Suponemos que el presente existe, como entidad separada del pasado, o del futuro, y se asume en mecánica cuántica que puede describirse de manera exhaustiva por lo que se  denomina función de onda $\psi $ y se construye como sigue.

En primer término, la función de onda toma valores complejos para que pueda representar una amplitud. Denota la amplitud de que en el tiempo $t$ el sistema esté en un estado específico, por ejemplo, en una vecindad del espacio tridimensional. De acuerdo con lo dicho, tal función se calcula por

$$\psi (x,t) =
\int \sb{y} \sp{x} e\sp{(i/\hbar)S(y,x,\gamma)}D\gamma (t)
\addtocounter{ecu}{1}   \hspace{4cm} (\theecu )      $$

y tenemos que sumar el aporte de todos los caminos habidos y por haber que desde algún lugar $y$, que el sistema ocupó en el pasado, conducen hasta $x$, ocupado en el presente.

Debido a que el presente existe, todos los caminos que del pasado lleguen hasta el futuro tienen que pasar por el presente, es decir, tienen que hacer escala en algún punto $y$ del espacio de estados ocupado en el presente.

Por tanto, si conocemos el presente, $\psi(x,t)$, entonces $\psi(x,t+\epsilon)$ se calcula, esquemáticamente, por:

\

$\psi(x,t+\epsilon)=$ amplitud de estar en $x$ en el tiempo $t+\epsilon $

 $= \sum $(amplitud de estar en algún lugar $y$ en el tiempo $t=\psi(y,t))$ x (amplitud de transición de $y$ a $x=e\sp{(i/\hbar)S(y,x,\gamma)})$.

donde la suma se toma sobre todos los puntos $y$ y todos los caminos que desde $y$ llevan hasta $x$. Eso se acostumbra a notar como:

$$\psi(x,t+\epsilon)=
\int \sb{y} \sp{x} e\sp{(i/\hbar)S(y,x,\gamma)} \psi(y,t)D\gamma (t)
\addtocounter{ecu}{1}   \hspace{4cm} (\theecu )      $$

Ahora mostraremos, en un caso especial, que podemos calcular la forma como el futuro evoluciona a partir del presente.

\section{PARTICULA EN UNA DIMENSION}

Tomaremos como ejemplo de estudio a una partícula que se mueve en una dimensión bajo la acción de un campo de fuerzas conservativo. En ese caso, el Lagrangiano es  $L= (1/2)m v\sp{2} - V(x,t)$, por tanto la acción correspondiente para un camino $\gamma$  parametrizado con respecto al tiempo, desde tiempo inicial $t \sb{0}$ a tiempo final $t \sb{0}$ es :

$S(\gamma) = \int \sb{t \sb{0}}\sp{t\sb{1}} L dt$

de tal forma que la amplitud total para ir de $a$ a $b$ es:

$ K(b,a)= \int \sb{a} \sp{b} e\sp{(i/\hbar)S(a,b,\gamma)}D\gamma (t) $

$=\int \sb{a} \sp{b} $
$e\sp{(i/\hbar)\int \sb{t \sb{0}}\sp{t\sb{1}} Ldt}D\gamma (t) $

$= \int \sb{a} \sp{b} $
$e\sp{(i/\hbar)\int \sb{t \sb{0}}\sp{t\sb{1}}[(1/2)m v\sp{2} - V(x,t)]dt} D\gamma (t)$

Por consiguiente:

$\psi(x,t+\epsilon)= \int \sb{y} \sp{x} $
$e\sp{(i/\hbar)\int \sb{t}\sp{t+\epsilon}[(1/2)m v\sp{2} - V(v,t)]dt}$ $\psi(y,t)D\gamma (t) $

Hay que tener en cuenta que para $\epsilon $ pequeño,

$\int \sb{t}\sp{t+\epsilon}  f(t)dt = base \times altura $
$= \epsilon f(t)$

Aplicando eso a la integral del Lagrangiano, nos queda:

$\psi(x,t+\epsilon)=\int \sb{y} \sp{x} $
$e\sp{(i \epsilon /\hbar)[(1/2)m v\sp{2} - V(v,t)]}$
$ \psi(y,t)D\gamma (t) $

Tenemos que hacer la integral teniendo en cuenta todos los caminos que llevan de algún lugar $y$ a $x$. Argumentemos que tan sólo hay que tener en cuenta los caminos directos. En efecto, un camino no directo en una dimensión implica retrocesos y por tanto, para llegar en el tiempo requerido tendría que tener aceleraciones y velocidades mayores, por eso una pequeña deformación en el camino cambiaría la acción notablemente, es decir la fase cambiaría alocadamente. Por tanto, caminos vecinos que no sean directos aportan amplitudes que se aniquilan mutuamente.

Por consiguiente, teniendo en cuenta los caminos directos y nada más, podemos tomar como velocidad $(x-y)/ \epsilon$ y como lugar para evaluar el potencial, tomamos el punto intermedio entre $y$ y $x$, o sea, $(x+y)/2$. De igual forma, el potencial se evalúa en el tiempo intermedio, que es $t +\epsilon /2$. En consecuencia:

$\psi(x,t+\epsilon)=$
$\int \sb{-\infty} \sp{\infty} (1/A)$
$e\sp{(i \epsilon /\hbar)[(1/2)m [(x-y)/ \epsilon] \sp{2}- V((x+y)/2,t +\epsilon /2)]} \psi(y,t)dy $

Hemos puesto el factor $1/A$ porque no conocemos cómo se define la medida en el espacio de caminos. Lo único que requerimos es que se pueda normalizar adecuadamente, en particular que la suma de todas las probabilidades de un 1. Separando los términos de la exponencial obtenemos:

$\psi(x,t+\epsilon)= $
$\int \sb{-\infty} \sp{\infty} (1/A)$
$e\sp{(i  /\hbar)  [(1/2)m [(x-y)\sp{2}/ \epsilon]]}  $
$e\sp{(i \epsilon /\hbar)  [- V((x+y)/2,t +\epsilon /2)]}\psi(y,t)dy $

Podemos esperar que las amplitudes más importantes sean las provenientes de puntos cercanos a $x$ debido a que para puntos lejanos la velocidad tiene que aumentar, y como figura al cuadrado, la fase oscilará rápidamente de tal forma que el aporte de cualquier camino se anulará con cualquiera de sus deformaciones. Por tanto, hacemos el cambio de variable $y=x+\eta$, o bien $y-x=\eta$,  con el objetivo de hacer un desarrollo en serie como función de $\eta$. Tenemos que $x+y = y+x = y-x+2x= \eta +2x= 2x+\eta$. Como  $ y=x +\eta$,  y además $x$ es constante,  entonces $dy = d \eta$:

$\psi(x,t+\epsilon)= \int \sb{-\infty} \sp{\infty} (1/A)$
$e\sp{(i  /\hbar)  (1/2)m [(-\eta) \sp{2}/ \epsilon]}  $
$e\sp{(i \epsilon /\hbar)  [- V(x+\eta/2,t +\epsilon /2)]}\psi(y,t)dy $

Teniendo en cuenta que $\epsilon $ buscará su límite a cero y que $V$ es continua, aproximamos $ t +\epsilon /2$ por $t$. Hacemos una aproximación similar con $\eta$, por lo que podemos reemplazar $V(x+\eta/2,t +\epsilon /2)$ por  $V(x,t)$:

$$\psi(x,t+\epsilon)= \int \sb{-\infty} \sp{\infty} (1/A)
e\sp{i  m \eta \sp{2}/2 \hbar \epsilon}
   e\sp{-(i \epsilon /\hbar)  V(x,t)}\psi(x + \eta,t)d\eta   \hspace{1cm} (\theecu )      $$

Expandimos en serie hasta primer orden en el lado izquierdo. Si mantenemos $x$ constante:

$\psi(x,t+\epsilon)= \psi(x,t) + \epsilon \partial \psi/\partial t$

Teniendo en mente el lado derecho de la ecuación en elaboración, la expansión de la exponencial, a primer orden, es:

$e\sp{-x} = 1-x$

Aplicando esto:

$ e\sp{-(i \epsilon /\hbar)  V(x,t)}= 1- (i \epsilon /\hbar)  V(x,t)$

 Recordemos el desarrollo en serie de Taylor hasta orden 2 para después aplicarlo a la amplitud:

$f(x+\eta)=f(x) + \eta df/dx + (1/2) \eta \sp{2} d\sp{2}f/dx\sp{2}$

Y si mantenemos $t$ fijo en la amplitud, nos queda como función de una sola variable y la podemos expandir en la otra:

$\psi(x + \eta,t)= \psi(x,t) + \eta \partial \psi/\partial x $
$+ (1/2)\eta \sp{2} \partial \sp{2} \psi /\partial x\sp{2}$

Hemos desarrollado hasta segunda potencia de $\eta$ porque en la acción $\eta $ es proporcional a la velocidad y ésta va al cuadrado. Substituyendo todo en (33) nos queda en conclusión:

$\psi(x,t) + \epsilon \partial \psi/\partial t = $

$=\int \sb{-\infty} \sp{\infty} (1/A)$
$e\sp{i  m \eta \sp{2}/(2 \hbar \epsilon)}      $
$[1- (i \epsilon /\hbar)  V(x,t)] $
$[ \psi(x,t) + \eta \partial \psi/\partial x $
$+ (1/2)\eta \sp{2} \partial \sp{2} \psi /\partial x\sp{2}]d\eta$

De acá salen 6 integrales. Hay que tener en cuenta que la variable de integración es $\eta$ y todo lo que no dependa de ella es constante que sale de la integral. La tarea se reduce entonces a evaluar varias integrales, la primera de las cuales es:

$\int \sb{-\infty} \sp{\infty} (1/A)e$
$\sp{i  m \eta \sp{2}/(2 \hbar \epsilon)} d\eta$

Esta integral es la de la campana de Gauss que se resuelve por el truco de hallar su cuadrado con dos variables mudas independientes, reescribir como una doble integral y cambiar a coordenadas polares. Para todos los efectos, el número $i$ hace las veces de constante y se trata como tal. Dicha integral nos da:

$= (1/A)(2 \pi i \hbar \epsilon / m)\sp{1/2}$

Eso nos permite comparar el coeficiente de $\psi$ en cada lado, pues la ecuación es verdadera para todo incremento de la variable $x$, y como conclusión sacamos que

$\psi(x,t) =\int \sb{-\infty} \sp{\infty} (1/A)$
$e\sp{i  m \eta \sp{2}/(2 \hbar \epsilon)}\psi(x,t)d\eta $

$=\psi(x,t) \int \sb{-\infty} \sp{\infty} (1/A)$
$e\sp{i  m \eta \sp{2}/(2 \hbar \epsilon)}d\eta $

Por consiguiente, la integral indicada, la de la campana de Gauss, debe ser uno, pero ya sabíamos su valor en términos de $A$:$ (1/A)(2 \pi i \hbar \epsilon / m)\sp{1/2}$

$\int \sb{-\infty} \sp{\infty} (1/A)$
$e\sp{i  m \eta \sp{2}/(2 \hbar \epsilon)}d\eta =1$$= (1/A)(2 \pi i \hbar \epsilon / m)\sp{1/2}$.

Despejando A nos da:

$A= (2 \pi i \hbar \epsilon / m)\sp{1/2}$

Entre tanto,

$\int \sb{-\infty} \sp{\infty} (1/A)$
$e\sp{i  m \eta \sp{2}/(2 \hbar \epsilon)} \eta d\eta= $

$lim_{ L\rightarrow \infty} \int \sb{-L} \sp{L} (1/A)$
$e\sp{i  m \eta \sp{2}/(2 \hbar \epsilon)} \eta d\eta= 0$

porque se trata de una función impar de $R$ en $C$. Los límites de integración se toman simétricos porque de otra forma la integral no converge. Eso es equivalente a tomar las ecuaciones en el sentido de las distribuciones que usan funciones test con soporte compacto y simétrico. Por otra parte, después de aplicar integración por partes y reencontrar la densidad de la campana de Gauss tenemos:

$\int \sb{-\infty} \sp{\infty} (1/A)$
$e\sp{i  m \eta \sp{2}/(2 \hbar \epsilon)} \eta \sp{2} d\eta$
$= i \hbar \epsilon /m $

Reemplazando en las 6 integrales nos queda:

$\psi + \epsilon \partial \psi/\partial t $
$= \psi - (i \epsilon /\hbar) V \psi $
$- (\hbar \epsilon/2im) \partial \sp{2} \psi /\partial x\sp{2}$

Simplificando, primero $\psi$ y después $\epsilon$ tenemos:

$$\partial \psi/\partial t = - (i  /\hbar) V \psi
- (\hbar /2im) \partial \sp{2} \psi /\partial x\sp{2}
\addtocounter{ecu}{1}   \hspace{4cm} (\theecu )$$

Antes de esta ecuación figuraba, en forma invisible, un número infinito de términos con altas potencias de $\epsilon$, que al tomar el límite desaparecen. Esta ecuación es un algoritmo con el cual la naturaleza construye el futuro inmediato a partir del presente. En efecto:

$(\psi(x,t+\epsilon)-\psi(x,t))/\epsilon \approx  \partial \psi/\partial t = - (i  /\hbar) V \psi
- (\hbar /2im) \partial \sp{2} \psi /\partial x\sp{2}$

de donde podemos despejar el valor de $\psi(x,t+\epsilon)$. De hecho, tenemos ahí una manera para mimificar la naturaleza y predecir la evolución de un sistema cuántico.

\

A pesar de todo el enredo de cosas que entran en la integral de Feynman, de la ecuación de Schr\"{o}edinger parece deducirse que lo que prevalece es un efecto local. Eso se debe a que el futuro en una vecindad se determina por un efecto de las varias derivadas, las cuales toman sólo en consideración lo que pasa en la inmediata vecindad de la partícula. Pero atención: eso es válido para la función de onda y sucede que ella no es observable y no puede serlo pues no se puede determinar sino módulo una fase inobservable. Lo que es observable es la probabilidad de transición de un estado a otro, o cambios de energía, y ellos se calculan por productos internos, que son integrales y por lo tanto la mecánica cuántica predice que el futuro se construye a partir de lo que pasa en todo el mundo.

Por supuesto que cuando un problema de  mecánica cuántica puede reducirse a otro de tipo clásico, los efectos globales se destruyen unos a otros en favor de los efectos locales. Pero lo terrible es que  la mecánica cuántica permite fuertes violaciones a la idea de que el futuro del universo es algo que se computa por efectos locales.

La predicción de violar el principio de localidad ha sido discutido profundamente y al tema se le conoce como las interacciones EPR ( de Einstein, Podolski y Rosen). Los estudios experimentales apoyan las predicciones cuánticas: el universo funciona como un todo, como un organismo y no parte por parte.

El \index{holismo cuántico} \textbf{holismo cuántico} plantea interrogantes como el siguiente: si el mundo entero interfiere con mis decisiones, cómo podré yo declararme responsable de lo que hago?

Me parece que hay pocas personas que toman tan en serio este tipo de preguntas. Una razón puede ser que una respuesta rigurosa y a prueba de toda objeción sobrepasaría todo poder de cálculo y terminaría siendo objeto de caprichos subjetivos.

\bigskip

La ecuación  ($\theecu$) se escribe más usualmente como :

\addtocounter{ecu}{1}

$$-(\hbar/i)\partial \psi/\partial t =
- (\hbar \sp{2} /2m) \partial \sp{2} \psi /\partial x\sp{2} + V \psi \hspace{4cm} (\theecu ) $$

Esta es la famosísima ecuación de Schr\"{o}edinger, que se convirtió en la piedra angular de la mecánica cuántica a partir de su descubrimiento por Edwin Schr\"{o}edinger en 1925, quien la descubrió por otro camino.  Es usual sintetizar dicha ecuación en la forma:

$$-(\hbar/i)\partial \psi/\partial t = H \psi
\addtocounter{ecu}{1}   \hspace{4cm} (\theecu )      $$

Donde $H$ es el operador llamado Hamiltoniano y definido por:

$$H = - (\hbar \sp{2} /2m) \partial \sp{2} /\partial x\sp{2} + V
\addtocounter{ecu}{1}   \hspace{4cm} (\theecu )      $$

\

La  ecuación de Schr\"{o}edinger que hemos descrito está desarrollada para una dimensión. Su generalización a tres dimensiones se  puede hacer de dos maneras. Partiendo del Lagrangiano en tres dimensiones y siguiendo el mismo camino que hemos recorrido. Ese es un problema tedioso con 18 integrales.  Nosotros vamos a hacer una inferencia: teniendo en cuenta que no hay ninguna dirección privilegiada en la naturaleza, entonces, lo que es válido en la dirección $X$, también lo es en la $Y$ o en la $Z$. Por otra parte, al tener en cuenta las tres direcciones, el desarrollo en serie de Taylor, a primer orden, le concede igual valor a cada una de forma independiente. Por tanto, la adivinanza (se dice Ansatz)que proponemos, reza de la siguiente forma:

$-(\hbar/i)\partial \psi/\partial t = $
$- (\hbar \sp{2} /2m) [\partial \sp{2} \psi /\partial x\sp{2}$
$+ \partial \sp{2} \psi /\partial y\sp{2} $
$+ \partial \sp{2} \psi /\partial z\sp{2}] + V \psi $

Identificando el operador Laplaciano, podemos reescribirla como:

$-(\hbar/i)\partial \psi/\partial t = - (\hbar \sp{2} /2m) \nabla \sp{2} \psi  + V \psi $

Por tanto, el \index{operador Hamiltoniano} \textbf{operador Hamiltoniano} tiene la forma general:

$$H= - (\hbar \sp{2} /2im) \nabla \sp{2}   + V
\addtocounter{ecu}{1}   \hspace{4cm} (\theecu )      $$

\

Podría pensarse que $H$ es la traducción del Lagrangiano, pero eso no es cierto. Siguiéndole la pista al Lagrangiano a través de la acción uno se da cuenta que la energía cinética se convierte en una frecuencia angular y que las segundas derivadas parciales tienen un origen muy diferente. De paso, entre más grande sea la energía cinética, mayor es la frecuencia y tenemos por tanto una asociación partícula-onda.

Si el operador Hamiltoniano $H$ no es la versión cuántica del Lagrangiano, qué significa entonces $H$? Respondamos lentamente a esa pregunta, estableciendo primero cual es su dominio de definición.

\section{EL ESPACIO PROYECTIVO }

Revisamos muy rápidamente algunos conceptos que se usan usualmente y que necesitamos en esta sección.

\

\color{blue}

$<<<<<<<<<<<<<<<<<<<<<<<<<<<<<<<<<<<<<<<<<$

\

Ideogramas: producto interno, espacio normado completo o de Banach, espacio de Hilbert, integral de Lebesgue, espacios de Lebesgue.

\

\textbf{Producto interno} \index{producto interno} o escalar: se denomina producto a cualquier operación binaria que distribuya la suma. El producto escalar es un producto que da como resultado un escalar. Por ejemplo, si tomamos $R\sp{2}$ tenemos un producto interno definido entre $\vec u = (u\sb{1}, u\sb{2})$ y $\vec v = (v\sb{1}, v\sb{2})$ por

$<\vec u,\vec v> = u\sb{1}v\sb{1}+ u\sb{2} v\sb{2}$

Este producto cumple:

a. $<\vec u,\vec v>=<\vec v,\vec u>$

b.$<\lambda \vec u,\vec v>= \lambda  <\vec u,\vec v>$

c. $<\vec u,\vec v + \vec w>=<\vec v,\vec u> + <\vec v,\vec w>$

d. $<\vec u,\vec u> >0$ si  $\vec u \not= \vec 0$

\

Gracias al producto interno se pudo definir una norma o largo de cada vector como la raíz cuadrada del producto escalar del vector con él mismo:

\

$ \|\vec u\|= (<\vec u,\vec u>)\sp{1/2} $

\

Un \textbf{espacio vectorial se llama completo} \index{espacio completo} cuando toda sucesión   de vectores (un conjunto de vectores numerados) que sea de Cauchy (que amenace con convergir) converge en realidad a un vector en el espacio.

\

Al número de direcciones independientes que hay en un espacio vectorial se llama dimensión: por eso $R\sp{3}$ tiene dimensión 3. Todos los espacios vectoriales de dimensión finita son completos.

\

Hay espacios vectoriales de dimensión infinita?

Si y abundantemente: los espacios de funciones son ejemplos muy apropiados. Tomemos las funciones de variable real que van de $R$ en $R$. Dichas funciones se pueden sumar, como por ejemplo, $senx +cosx$, y se pueden alargar, como $5senx$. También se puede definir un producto escalar con todas las propiedades deseadas:

$<f,g>= \int\sb{R} f(x)g(x)dx$

y por consiguiente podemos definir la norma cuadrado de una función como

$\|f\|\sp{2} = <f,g>$

\

Espacio vectorial normado completo: si ya sabemos que podemos definir una suma y un alargamiento en un conjunto dado de funciones escalares, y además una norma, y si el espacio es completo, lo llamamos espacio vectorial normado completo o de Banach.

\

\textbf{Espacio de Hilbert} \index{Espacio de Hilbert}: es un espacio vectorial con producto interno o escalar, que define una norma que lo hace completo o de Banach.

Hemos visto cómo se define un producto interno de funciones escalares reales. Si las funciones escalares son complejas, es decir, toman valores complejos, el producto interno se define ayudándose del conjugado como:

$<f,g>= \int\sb{R} f(x)g(x)\sp{*}dx$

y cumple propiedades parecidas al producto interno de funciones escalares reales. Pero hay diferencias: en especial recalcamos las siguientes dos propiedades válidas para $\lambda$ complejo. Tenemos 		que los escalares complejos salen igual desde la primera entrada del producto interno, pero salen conjugados de la segunda:

$<\lambda f,g>= \int\sb{R} \lambda f(x)g(x)\sp{*}dx $
$ =  \lambda< f,g>$

$< f,\lambda g>= \int\sb{R} f(x)(\lambda g(x))\sp{*}dx$
$ =  \lambda^*< f,g>$

\

La razón por la cual involucramos al conjugado es para poder definir la norma (que cumpla las propiedades de largo de un vector) como:

$||f||^2 = <f,f>$

tal como se hace en espacios con  producto interno real.

\

Integral  de Lebesgue:

\

Vimos que para definir la integral de Riemann se programaba el taxímetro secuencialmente, a medida que uno avanzaba por la carretera y que el taxímetro tenía en todo momento precisión infinita.

En la integral de Lebesgue se procede de otra forma: el taxímetro que se usa es ralo, con pocos valores permitidos, pero de densidad ajustable hacia la plenitud. Y cómo se hace la integral de Lebesgue de $f(x)$? Para funciones que no tengan demasiados altibajos ni discontinuidades, funciona el siguiente proceso:

Se ajusta una baja precisión en el taxímetro y se recorre el camino: se va comparando a $f(x)$ con los valores permitidos del taxímetro y se aproxima $f(x)$ al valor más cercano posible permitido por el taxímetro. Eso define tramitos en los cuales se tiene una función constante cuyo aporte a la integral es base x altura.

Pero atención, la base tiene una manera de ser medida que hace caso omiso de fallas como la siguiente: si al intervalo [0,1] le quitamos todos los números racionales, de todas formas su medida es uno.

Se refina la precisión del taxímetro y se repite el proceso hasta que se demuestre que se llega a un límite al cual se  denomina integral de Lebesgue (se lee lebeg).

Cuando la integral de Riemann existe, se puede garantizar que coincide con la de Lebesgue. Pero hay funciones que no tienen integral de Riemann y en cambio si tienen integral de Lebesgue, como es el caso de la función que sobre [0,1] está definida como 1 sobre los irracionales y 0 sobre los racionales. Su integral de Lebesgue vale 1, puesto que los racionales por ser un conjunto que se puede numerar no generan área alguna y no le quita nada a la integral.

\

\textbf{Espacios de Lebesgue} \index{Espacios de Lebesgue}: uno está acostumbrado a errores tipográficos menores y a identificar una gráfica dada, digamos, con $senx$ aunque falten puntos aquí y allá. Los espacios de Lebesgue son espacios de funciones a los cuales se les hace la misma concesión.

Habiendo definido la integral de Lebesgue se procede a notar que se puede definir un producto interior entre dos funciones como la integral del producto de las dos funciones. Se define la norma de una función como la raíz cuadrada del cuadrado de la función.

Resulta que hay pares de funciones que siendo diferentes, su diferencia no se nota al ser medida por la norma. Por ejemplo, una función continua y otra que sea igual a ésta en todas partes excepto en un punto difieren por la integral sobre un punto que puede interpretarse como el área de la función dada sobre un punto, o sea cero. Es decir que desde el punto de vista de la norma no difieren en nada de valor o sea que son equivalentes.

\

Los espacios de Lebesgue son los que resultan de no hacer diferencia entre pares de funciones equivalentes.
El más famoso de eso espacios es $L\sb{2}$ que es con el que hemos venido y seguiremos trabajando. En general, una función $f$ pertenece a $L\sb{p}, p \in R $ si y sólo si $f(x)\sp{p}$ es lebesgue-integrable, es decir si su integral como límite existe y es finito. En ese caso, la norma es la raíz p-ésima de la integral de $f(x)\sp{p}$:

$\|f\| \sb{p} = ( \int \sb{R} f(x) \sp{p} dx )\sp{1/p}$

Lo que hemos hablado en $R$ se generaliza de manera natural a funciones escalares sobre $R\sp{n}$.

\

Como nosotros nos referimos, casi siempre, a funciones que toman valores complejos y trabajamos en $L\sb{2}$ no escribimos el subíndice 2 de la norma, la cual se define como:

$\|\psi\|  = ( \int  \psi \psi \sp{*}  dV )\sp{1/2}$

En mecánica cuántica no se usa la norma, sino la norma cuadrado.

Las materias encargadas de formalizar estos conceptos se llaman  análisis, teoría de la medida y análisis funcional.

\

$>>>>>>>>>>>>>>>>>>>>>>>>>>>>>>>>>>>>>>>>>$

\color{black}

\

Evidenciemos ahora que hay una \index{arbitrariedad} \textbf{arbitrariedad} matemática en la ecuación de Schr\"{o}edinger.

En principio el operador $H$ se define sobre toda función con valores complejos que tenga segundas derivadas. Como las funciones que toman valores complejos y tienen segundas derivadas se pueden sumar y multiplicar por escalares complejos sin perder dichas propiedades, el conjunto de tales funciones forma un espacio vectorial. Definiendo un producto interior también podemos definir una norma y nos restringimos a funciones cuya norma sea finita. Al ver que dicho espacio es completo tenemos un espacio de Hilbert que de hecho resulta separable, gracias a lo cual podemos utilizar el computador para hallar soluciones aproximadas a los problemas que resulten.  A ese espacio lo notamos $L\sb{2}$.

Procedemos a definir el producto interno. Si notamos el conjugado de $\psi$ como $\psi \sp{*}$ podemos definir un producto interno entre dos funciones como:

$$<\psi \sb{1}, \psi \sb{2}> = \int \psi \sb{1} \psi \sb{2}\sp{*} dV
\addtocounter{ecu}{1}   \hspace{4cm} (\theecu )      $$

donde la integral se sobrentiende que corre sobre todo el universo.

Entonces, la expresión $|\psi |\sp{2}= \psi  \psi\sp{*}$ corresponde a una densidad de probabilidad, la cual integrada sobre todo el universo debe dar uno. Pero esa es precisamente la norma en el espacio $L\sb{2}$ de la función de onda:

$$\|\psi \| \sp{2}= <\psi, \psi > = \int \psi  \psi\sp{*} dV
= \int |\psi |\sp{2} dV = 1
\addtocounter{ecu}{1}   \hspace{4cm} (\theecu )      $$

Por lo tanto, las funciones de onda son elementos que pertenecen a la bola unitaria de un espacio de Hilbert, que también es métrico.

Si observamos la ecuación de Schr\"{o}edinger, ella es lineal. O sea que permite ser multiplicada por cualquier escalar complejo. Eso significa que uno puede tomar arbitrariamente el escalar que a uno le guste, y tal escalar jamás aparece en los cálculos. Entonces, pues no puede tener significado físico. Por tanto, una función y un múltiplo de ella son equivalentes desde el punto de vista físico. Eso implica que la física se lee no sobre funciones sino sobre clases de equivalencia, módulo una constante. El espacio resultante se llama \index{espacio proyectivo} \textbf{espacio proyectivo}. Explícitamente, la relación de equivalencia $\Re$ es:

$ \psi \sb{1} \Re  \psi \sb{2}  $ si y sólo si $ \psi \sb{1}$
$= \lambda \psi \sb{2}, \lambda \in C \addtocounter{ecu}{1}   \hspace{3cm} (\theecu) $

Dicha relación tiene una metáfora muy simple: imaginemos el plano como espacio vectorial. Una forma de decir que lo único que importa de un vector  es su dirección y no su norma es decir que un vector está relacionado con otro si se obtiene del primero por una alargamiento o un acortamiento, como en el caso que acabamos de ver.   Eso quiere decir que ahora tenemos un espacio conformado por las líneas que pasan por el origen, las cuales pueden ser indexadas por  los  puntos del hemisferio superior del círculo unidad e igualando los dos puntos que quedan sobre el eje $X$, pues ellos definen la misma línea. Ese es el equivalente del espacio proyectivo cuántico. 

\

\section{FEYNMAN Y SU GENIALIDAD}

La integral de Feynman que hemos visto, que  construimos por sentido común  a partir de la óptica,  números complejos y el lagrangiano, apareció en nuestro mundo como la herencia de un solo hombre: Richard Philip Feynman.  Algunos lo consideran el máximo genio de la física del siglo XX.  Hagamos una pequeña disquisición sobre su genialidad y veamos si podemos aprender algo sobre cómo ser un genio. Todas las referencias no anotadas vienen de la biografía de Feynman escrita por Mehra (1994).

\bigskip

\textbf{El Challenger}

\bigskip

El brillante porvenir de la carrera cosmonáutica del Challenger  recibió un duro golpe  el 28 de Enero de   1986, día  en el  cual murieron 7 astronautas.  Pero  a las dos semanas, el  11 de febrero,    en una reunión transmitida en directo   por televisión para millones de personas, apareció el premio Nobel de Física Richard Feynman explicando la causa del accidente:  el Nobel  hizo una demostración con agua casi helada y un pedazo de caucho para mostrar a todos, aún a los niños, que un caucho elástico a temperatura ambiente pierde su elasticidad si se  mete en agua a punto de congelarse. Es decir, la elasticidad del caucho depende de la temperatura y se pierde cuando ésta es muy baja.  Al haber perdido elasticidad, el caucho se rompió al experimentar presión, apareció una grieta por donde salieron llamas, las cuales alcanzaron una parte clave del equipo de control lo mismo que un tanque de combustible, el cual se separó del resto. Después se destruyó todo por  causas aerodinámicas (ver 'Challenger', Wikipedia, 2008).

La demostración de Feynman fue un momento de gloria no sólo para él sino para toda la ciencia.  La razón es que un equipo conformado por expertos ingenieros no pudo ni prever con certeza ni encontrar la causa del accidente y eso a pesar de tenerse indicios sobre la importancia crítica de los tales cauchos de sellamiento. Y  en cambio un teórico de la ciencia si pudo encontrar la causa en menos de una semana.  Después de eso, el programa   del Challenger reanudó  una carrera que llega hasta nuestros días y que sigue adelante aún a pesar  de desastres mortales que muestran la complejidad del oficio.

\bigskip

\textbf{¿Es la ciencia o es Feynman?}

\bigskip

Quisiera uno entender  qué fue  lo que pasó.  Es decir,  quisiera uno responder a la siguiente pregunta: ¿qué le faltaba  a la ingeniería de alto nivel que sí tenía  la ciencia pura para que Feynman hubiese notado algo que estaba a la vista de todos durante tanto tiempo?  

La  respuesta es un tanto desesperante:  la ciencia pura no tiene nada que ver en el asunto.  Sabemos eso con certeza por lo siguiente: había un   enigma en  la ciencia, a la vista de todos  durante 7 años,  el cual había sido propuesta por Dirac en 1933  en relación con la exponencial  imaginaria del Lagrangiano, y nadie lo había descifrado. Pues Feynman lo descifró en 1941 a su 23 años tan pronto como se enteró de que existía.  Esto le abrió la puerta a la Integral y a los Diagramas que llevan su nombre y que lo  hicieron famoso con una fama que durará mientras que haya hombres.  

Así que Feynman tenía algo que lo hacía grande tanto  entre los grandes ingenieros como entre los grandes  genios de la física teórica.  ¿Qué era?   ¡Quién sabe!  No lo sabemos, pues si lo supiésemos estaríamos produciendo genios en serie.  Con todo, hay algunos factores que son de indiscutible significado.

Para empezar, notemos  que si bien las ideas revolucionarias le nacen en un instante,  trabajarlas le lleva tiempo y mucho tiempo de un esfuerzo continuado.  Eso lo relata muy bien en su charla del Nobel: su duro trabajo le permitió formar variadas y diversas visiones sobre el universo en que se movían sus preguntas.  Quizá convenga resaltar que el duro trabajo incluía la perfección del detalle como la generación y fortalecimiento de   una gran contextualización, la cual permite, entre otras cosas, utilizar diversas herramientas para atacar un mismo problema y ayudarse de las diversas visiones para detectar y corregir rápidamente los errores causados por el mal manejo de una de ellas.

El duro trabajo produce fruto. Ese es un axioma al que  Feynman le hace mucha propaganda a través de su   slogan preferido: lo que un bobo puede hacer, otro bobo podrá hacerlo mucho mejor ('what one fool can do, another can do better').  Feynman enseña eso con autoridad, pues fue la manera como él se educó a sí mismo.

Miremos ahora   un detalle de su temprana juventud que nos da otra luz para dilucidar lo que ese gran hombre tenía:  en unas vacaciones, Feynman consiguió un trabajo como químico de superficies. Uno de los problemas de su tiempo era dar un baño de plata a los plásticos.    Por ejemplo, hay un plástico muy barato a base de acetato de celulosa y que no era posible platear con los métodos tradicionales. Feynman resolvió el problema  modificando un poco los procedimientos que funcionaban para otros plásticos. El procedimiento que llegó a funcionar era así:  primero se pone el plástico durante unos momentos en hidróxido de sodio, después se pone por muchas horas en cloruro de estaño y después se si se platea, usando un proceso de reducción a partir de nitrato de plata. Superbarato y efectivo.  

Lo alegre de esta narración es que se presta a una generalización  automatizada: todas las grandes  ideas resultan de un proceso de pequeñas modificaciones  a ideas ya existentes unidas a una superposición de   unas cuantas palabras o de otras   pequeñas ideas. Esto  puede hacerse con un computador y un filtro semántico. 

Pues bien, conocí a un profesor alemán que  programó un computador   para que combinara palabras claves de los escritos de Hegel y produjera frases. Su esperanza era producir grandes ideas al estilo Hegeliano.  El hacía de filtro semántico leyendo las frases y seleccionaba las que le parecían interesantes. De tanto en tanto  yo le preguntaba al profesor cómo le iba en su tarea y siempre me decía: no salió nada especial.  Hasta que él mismo se cansó de que nada saliera y abandonó su proyecto. Podemos aprender que poner ideas unas sobre otras no lleva a ningún lugar y sin embargo, eso es lo único que hace todo genio. ¿Acaso haya alguna manera inteligente de hacerlo? 

La vida de Feynamn nos propone dos  respuestas. La primera viene de su hermana que le dijo: lo que debes hacer es entender en tus propios términos,   es decir, contextualizar dentro de tu propio conocimiento.  Para entender mejor qué significa eso,  ayudémonos de la  siguiente narración de la vida del niño  Feynman, en relación a su profesión de reparador de radios cuando él era apenas un adolescente:  

'What`s the matter with the radio set? Asked Feynman. 'Oh', the man said, 'It makes noises when you first turn it on, but the noise the stops and everything becomes all right.' Richard  turned the radio on; it gave out  a boold-curdling noise, a terrifying racket, and one could see why the guy wanted to have it fixed. After a while, the noise quited down, and the radio would play all right. So Richard turned it on, listened to it, turned it off, and started to walk back and forth, thinking. 'What are you doing? Can you fix it?' the man asked. 'I'am thinking,' Richard answered. He tried to figure out how there could be noise that disappears? Something is changing with time, something is heating up before something else. He guessed that the amplifiers were heating up before the information came in from the grids, from the early circuits, thus it was picking up some kind of noise. So, it was `robably due to heating up of the tubes; if he could reverse the earlies tubes, maybe it would be all right. It would then heat up the other way around. He changed the tubes around and put them back, and switched the radio on. It was just as quiet as you please!

La segunda idea de cómo generar ideas brillantes viene de su esposa: cuando Fyenman fue comisionado para investigar el desastre del Challenger, el le dijo a su esposa Gweneth  que él no quería ir:  ' Look, anybody can do it. They can get somebody else'. Su esposa   dijo. 'No. I f you don't do it, there will be twelve people, all going in a group, going around from place to place together. But if you join the commission, there will be eleven people  -all in a group, going around from place to place- while the twelfth one, you, runs around all over the place, checking all kinds of unusual things. There  probably won't be anything, but if there is, you will find it.  There isn't anyone else who can do that like you can.'

Resaltamos: 'Además de ti, nadie andará por ahí haciendo preguntas alocadas.. ..' es decir, hay que tener una contextualización lo suficientemente amplia para producir ideas en contracorriente con las usuales pero que sonando a locura puedan también ser muy eficaces para resolver los problemas propuestos. 

\bigskip

\textbf{El Sultan} 

\bigskip

Esto de querer saber qué tiene Feynman ha hecho que se repita la historia del Sultán, el cual le pidió al rey Ricardo que le mandara su espada pensando que eso lo haría buen espadachín (por favor, no malinterprete al Sultán:  desde su punto de vista, el rey Ricardo no   parecía hacer cosas sobrenaturales, así que la causa de su excelente desempeño debía venir de su espada. Tengamos ahora en cuenta que una espada difiere de otra como mínimo en el tensor de inercia que entra a formar parte de la conducta de la espada al rotar con respecto a la muñeca.  Añadiendo todas las otras complejidades, llegamos a la conclusión de que la mejor manera de juzgar una espada es peleando con ella. Y esa fue la intención del Sultán).  En nuestro caso, lo que todo el mundo pensó, incluyendo al propio Feynman, es que habría un gran adelanto en la enseñanza de la física si se ponía su forma de pensar al alcance de todos. Hicieron entonces un gran esfuerzo editorial y pusieron por escrito las clases de Feynman con todo lo que él decía.

Ese material fue usado   en muchas partes. Pero no pasó mucho tiempo antes de ser sustituido por otros. Y, ¿cuál pudo ser la razón?  El material está  hermosamente contextualizado. Por ejemplo, en esa colección  uno puede aprender que los rayos no solamente caen sino que rebotan de nuevo hacia el cielo.  Pero la colección no es fuerte en las cosas que lo hacen a uno fuerte para los exámenes o para la investigación.     En cambio, los  textos que han prevalecido contextualizan las cosas con un único objetivo: fortalecer la conexión entre los conceptos fundamentales  y las aplicaciones inmediatas,  adicionando algunas pocas  salidas al mundo de la gran contextualización.   

Con todo, la gran variedad de estilos pedagógicos de tantos textos indica claramente que enseñar física es un problema demasiado grande y difícil. Tal y como lo dijo Feynman.  Por tanto, nada mejor que propiciar muchas y variadas soluciones. Lo que esto quiere decir es que el lector también puede hacer su propio esfuerzo editorial.  Para  el peor de los casos,   recordemos  que   un texto es como un caballo y que no hay nada mejor   caballo que aquel que uno mismo ha amaestrado. Por eso, los textos deben venir con derecho y poder para que el lector pueda modificarlo a su gusto. El presente texto es así: amable lector, tome la fuente Latex y utilícela para hacer su propio texto tal y como le venga en gana. Aún más: sería una felicidad saber que mi trabajo le pueda representar alguna  ganancia económica. 

\bigskip

\textbf{Nuestra lección}

\bigskip

De todo esto, aprendamos lo siguiente:
 
\begin{itemize}
	\item	Hay un tipo de inteligencia que es buena y muy buena para lo que sea, sea ciencia o ingeniería o aplicaciones prácticas.

	\item	Ese tipo de inteligencia  tiene cuatro pares de ojos, dos sobre los árboles y dos sobre el bosque.  Dos sobre los árboles, para mirar el detalle. Y otros dos para mirar el bosque, la contextualización, el gran contexto que todo lo encierra y todo lo interrelaciona de manera natural, sencilla y eficaz. Esa segunda manera de mirar produce el ordenamiento adecuado para poner las ideas unas sobre otras para dar explicaciones coherentes, efectivas y verdaderas. 
	\item	Este tipo de inteligencia disfruta trabajando y el trabajo duro   es el fundamento de la diversión. Además de lograr manejar las pequeñas cosas como lo hace un maestro, el resultado natural del trabajo duro  y diversificado  es una amplia  contextualización. 
	\item	El trabajo duro produce fruto y es propiedad de todo el mundo.  En mi propia experiencia docente, a mi me ha tocado experimentar eso por rutina, al ver cómo una estudiante de primer semestre de Lenguas  o de filosofía o de historia o de abogacía o de ciencia política  tiene que aprender a toda prisa a utilizar la estadística   a un nivel que hace 30 años sólo lo podían hacer los profesores de posgrado.  
	\item	Por muy buena que sea una inteligencia, no es buena para todo. Metiéndome en lo que no me importa,  me atrevo a decir que el segundo matrimonio de Feynman con  Marie Lou fue un fracaso anunciado que se hubiese podido evitar conservando una hermosa amistad.  Y como si fuese poco, un flirteo  con Gweneth  en una playa vacía se tradujo en un magnífico tercer matrimonio que le duró toda la vida y que le dió dos maravillosos hijos. ¿Es eso científico?
	\item	No todo en una gran inteligencia es misterioso. Buena parte se reduce  a  entender, lo cual significa simplemente poder relacionar los conceptos entre ellos y  con las cosas simples de este mundo y con las propias ideas y conceptos personales. 
	\item	A la inteligencia hay que ponerla a prueba,  hay que desafiarla dándole tareas difíciles  y poniéndola a que defienda sus resultados en público, ante la comunidad, ante el mundo, ante la historia de la humanidad.

\end{itemize}
 
\textbf{ Referencias}
\begin{itemize}
	\item The Feynman Lectures on Physics (with Leigton and Sands). 3 volumes 1964, 1966. Hay una edición bilingue, inglés -español (Feynman también hablaba español).
	
		\item Feynman R, (1965) The development of the space-time view of quantum electrodynamics.  Charla con ocasion del Premio Nobel.  Se consigue por Internet.  
		
			\item Mehra J., (1994) The beat of a different drum, The life and science of Richard Feynman.  Clarendon Press-Oxford. 
\end{itemize}

Sigamos ahora en la búsqueda de una interpretación del Hamiltoniano.

\section{EL GRUPO DE GALILEO}

La ecuación de Schr\"{o}edinger  fue deducida de una forma artificiosa a partir de la óptica y del formalismo lagrangiano. Todo eso tan complicado, tan circunstancial, se parece a la aventura de un ciego en terreno desconocido. Pero pasado algún tiempo, aún los ciegos son capaces de formarse una imagen geométrica de su entorno. Después de eso, todo aparece iluminado por una luz interna que lo guía a uno a donde sea.

Pues si eso pasa con los ciegos, no menos debemos esperar de la comunidad físico-matemática. Lo que vamos a presentar ahora es una imagen muy sofisticada que tiene un  solo objetivo: la sencillez. En efecto: vamos a demostrar que la mecánica cuántica tiene su origen en la teoría de grupos. Es decir, usando la teoría de grupos y los números complejos, redescubriremos la ecuación de Schöedinger. 

Hay que pagar un precio y es conocerse con estructuras matemáticas cuyo estudio puede ser tedioso pero que a la larga se vuelven rutinarias. Comencemos oficializando la terminología  sobre la que hemos venido trabajando.

\

\

\color{blue}

$<<<<<<<<<<<<<<<<<<<<<<<<<<<<<<<<<<<<<<<<<$

\

Términos y más términos:   isometría, operador unitario, operador adjunto, operador autoadjunto, valores y vectores propios, ortogonalidad. El operador derivada no es autoadjunto pero el operador segunda derivada sí. El Hamiltoniano es autoadjunto.

\textbf{Isometría}\index{Isometría}: iso quiere decir igual, y metro significa medida. Una isometría es una transformación que conserva las medidas, como una rotación o una translación.

Una isometría debe ser inyectiva, 1-1, conservar la identidad, pero no tiene porque ser sobre: consideremos la semirecta real positiva: un corrimiento de 3 unidades a la derecha representa una isometría pero no es sobreyectiva. Deducimos que una isometría no tiene porque ser invertible. Formalmente:

Definición: un operador lineal $T$ de un espacio de Hilbert $L$ en si mismo se denomina \textbf{isometría} cuando conserva la norma. Como la norma se define en términos de producto interior, una isometría también conserva el producto interior:

 $ <T(\psi),T(\phi)> = <\psi,\phi>$ para $\phi, \psi \in L$.

\

\textbf{Operador unitario} \index{Operador unitario}:

Decimos que un operador lineal es unitario cuando es una isometría sobreyectiva.

\

\textbf{Operador adjunto} \index{Operador adjunto}: consideremos   dos operadores $T$ y $T\sp{\dagger} $, tal que $T$ va de un espacio de Hilbert $L$ en otro $M$: $T: L \rightarrow M$, en tanto que $T\sp{\dagger}: M \rightarrow L $.   $T$ y $T\sp{\dagger} $ son contravariantes, o sea que van en contravía.

Si se cumple que

$<T(\phi),\psi>=<\phi,T\sp{\dagger} (\psi)>$

 para todo $\phi \in L$ y $\psi \in M$ decimos que $T\sp{\dagger} $ es el adjunto de $T$. La existencia de $T$ está garantizada cuando $T$ es acotado, o sea que transforma  la bola unidad  en un conjunto acotado, que cabe en otra bola de radio finito. Toda isometría lineal es acotada y por lo tanto tiene adjunto.

\

Hablando de un espacio de Hilbert $L$, decimos que un operador lineal $T$ de $L$ en sí mismo es autoadjunto si cumple con la igualdad

$<T\phi,\psi>=<\phi,T\psi>$

El conjunto de operadores autoadjuntos no es vacío: la identidad es autoadjunto. Estas definiciones generalizan otras de álgebra lineal: adjunto generaliza a matriz transpuesta, o sea la matriz resultante de cambiar las columnas por filas y al revés. Operador auto-adjunto generaliza a matriz simétrica, la cual es aquella que es su propia transpuesta. Ejemplo:

$$
\displaylines{
\pmatrix{1&4&5\cr
    4&2&6\cr
    5&6&3\cr}\qquad} $$

Esta matriz es simétrica pues es su propia transpuesta. Y cumple además con la igualdad:

$<T\vec u,\vec v>=<\vec u,T\vec v>$

para elementos $\vec u, \vec v$ de $R\sp{3}$.

\

SI $M$ es una matriz cuadrada cuyas entradas sean números complejos, la matriz adjunta de $M$, notada $M^\dagger$ se obtiene transponiendo a $M$ y conjugando a cada entrada. Se dice: $M^\dagger$ es la transpuesta conjugada. En este caso, $M$ define un operador lineal de $C^n$ en sí mismo y el producto interior de dos vectores de $C^n$ se define como sigue:

$\vec Z = (z_1,....z_n)$

$\vec W = (w_1,....w_n)$

$<Z,W> = z_1w^*_1 +.. .. + z_n w^*_n$ donde la $*$ significa conjugación.

\

Ahora volvemos a nuestro tema.

\

 Sea $V$ un espacio de Hilbert (vectorial con producto interno y completo en la norma inducida). Notamos como $L\sb{2}(V)$ al espacio de las funciones cuadrado integrable que de $V$ van sobre $C$.  Definimos el operador $A$ sobre $L\sb{2}(V)$   por:

$A(\phi) = \alpha \phi$, donde $\alpha$  es una función cualquiera del espacio $V$ en los complejos.

Veamos si $A$ es autoadjunto o no, es decir si se cumple:

$<A\phi,\psi>=<\phi,A\psi>$ donde $\phi$ y $\psi$ son funciones de $V$ sobre los complejos.

Reemplazando:

$<A\phi,\psi>=\int \alpha(x)\phi(x) (\psi(x)\sp{*})dx$

$<\phi,A\psi>= \int \phi(x) (\alpha(x)\psi(x))\sp{*}dx$
$=\int \phi(x) \alpha(x)\sp{*} \psi(x)\sp{*}dx$

de lo cual se concluye que $A$ es autoadjunto cuando $\alpha$ toma valores reales. Si $\alpha$ toma valores complejos, este operador no es autoadjunto. En particular, la multiplicación por un complejo define un operador que no es autoadjunto, sino anti-autoadjunto, es decir que pone un signo menos de demás. Para que haya orden en la sala se requiere que todas las integrales de las que se hable existan. Eso se puede garantizar cuando $\alpha $ es continua y además tiene soporte compacto (es nula fuera de algún intervalo cerrado y acotado).

\

Si se cumple la ecuación:

$A(\phi) = \lambda \phi$

para un operador lineal sobre un espacio de Hilbert, decimos que el escalar $\lambda$ es un valor propio y que $\phi \not= 0$ es un vector propio de $A$. No es necesario que todo operador siempre tenga valores y vectores propios. Por ejemplo, el operador $A$, de multiplicar por una función escalar, no los tiene a menos que la tal función sea una constante, y ese sería su único valor propio. En ese caso, cualquier función es un vector propio asociado a dicho valor propio.

\

Los \textbf{operadores autoadjuntos} \index{operador autoadjunto} tienen tres propiedades espectaculares:

1. Los valores propios son reales.

2. A valores propios reales corresponden vectores propios ortogonales.

3. Puesto que sobre vectores propios el operador es tonto, $A(\phi) = \lambda \phi$, los vectores propios son preciosos como elementos de una base para estudiar el efecto del operador.

Demostración de la primera propiedad.

Supongamos que para el operador lineal autoadjunto $A$ se tiene que

$A(\phi) = \lambda \phi$ , entonces,

$<A\phi,\phi>=<\lambda \phi, \phi> = \lambda <\phi, \phi> $
$=\lambda\| \phi\|\sp{2}$

pero como $A$ es autoadjunto,

$<A\phi,\phi> = <\phi,A\phi> = <\phi, \lambda \phi>$

ahora bien, estamos bajo una integral y el segundo factor entra conjugado, por eso:

$<\phi, \lambda \phi>=   \lambda \sp{*} <\phi, \phi>= \lambda \sp{*}\| \phi\|\sp{2}$

\

Uniendo todo tenemos:

$\lambda \| \phi\|\sp{2}= \lambda \sp{*}\| \phi\|\sp{2}$

Pero como el vector $\phi$ es no nulo, su norma no es cero y se puede dividir por ella, demostrando que $\lambda$ debe ser necesariamente un real puesto que es igual a su conjugado.

Para demostrar la segunda propiedad se toman dos valores propios distintos, $\lambda , \mu$, que ya sabemos que son reales:

$ A(\phi)= \lambda \phi$

$A(\psi) = \mu \psi$

y hacemos el producto escalar:

$<A(\phi),\psi> = <\lambda \phi,\psi>= \lambda <\phi,\psi>$

$<A(\phi),\psi> = <\phi,A(\psi)> = <\phi,\mu\psi>= \mu\sp{*} <\phi,\psi>= \mu <\phi,\psi>$

porque al pasar al segundo renglón usamos la propiedad de ser autoadjunto. Por transitividad:

$\lambda <\phi,\psi>=\mu <\phi,\psi>$

$(\lambda - \mu )<\phi,\psi>=0$

como los valores propios son diferentes, su resta no es cero y se puede tachar de la igualdad, de lo cual se deduce que

$<\phi,\psi>=0$

en un caso como ese, decimos que los dos vectores son ortogonales, al igual que se hace en dimensión 2 o 3.

\

Expliquemos ahora el siguiente slogan: el operador derivada no es autoadjunto, pero el operador segunda derivada, sí. Como hablamos de derivadas, estamos en un espacio de funciones y como nuestro tema es la mecánica cuántica, nuestras funciones van de $\Re$ (con generalización inmediata a  $\Re^3$) en los complejos. Tanto el  conjunto $C$ de los números complejos como el de los reales $\Re$ tienen suma, multiplicación y división sin problemas. También son completos en el sentido que toda \textbf{sucesión de Cauchi} converge a un elemento de $C$. Y eso es todo lo que se requiere para definir derivadas e  integrales. 

\

Se dice que una \index{sucesión!de Cauchi} \textbf{sucesión es de Cauchi} si para todo $\epsilon >0$ existe un $N$ tal que si $n, m > N $ entonces $|a\sb{n} -a\sb{m}| <\epsilon$. En un \index{espacio!completo} \textbf{espacio completo}, toda serie de Cauchi converge a un límite dentro del espacio, sea $L= lim a\sb{n}=a\sb{\infty}$.

\

Nuestro producto interior entre funciones $f$ y $g$ lo definimos como sigue:

$<f,g>  = \int_\Re f(t)g^*(t)dt$

La función $f$ podría ser $f(t) = e^{3it}$ en tanto que la $g$ podría ser $g(t) = e^{8it}$. En ese caso,

$<f,g>  = \int_\Re e^{3it}(e^{8it})^*(t)dt $

$=\int_\Re e^{3it}(e^{-8it})dt = \int_\Re e^{-5it} dt  = (1/(-5i))e^{-5it} |^\infty_{-\infty} = (1/(-5i)) (0-\infty ) = \infty $

decimos que la integral diverge. 

Los problemas de convergencia los resolvemos así: nuestro objetivo es modelar partículas por medio de funciones que van del espacio en los complejos, las funciones de onda. Como nuestras partículas están localizadas, suponemos que las funciones de onda tienen a cero en el infinito y lo mismo sus derivadas. Con eso podemos aplicar la integración por partes de una forma espectacular para poder explicar lo que nuestro slogan significa. 

Sobre nuestro espacio de funciones consideremos el operador derivada y verifiquemos que no es autoadjunto:

$<d\psi/dt, \phi > = \int_R (d\psi(t)/dt ) \phi^*(t)dt $
$ = \int^\infty_{-\infty} ( d\psi(t)/dt)  \phi^*(t)dt$

Ahora aplicamos integración por partes:

$  \int^\infty_{-\infty} ( d\psi(t)/dt)  \phi^*(t)dt = $
$ \psi(t)\phi^*(t)|^\infty_{-\infty} - \int^\infty_{-\infty}   \psi(t)   (d\phi^*(t)/dt)dt$

Aplicando las condiciones en el infinito, el primer término se anula y nos queda:

$<d\psi/dt, \phi > = - \int^\infty_{-\infty}   \psi(t)   (d\phi^*(t)/dt) dt = -<\psi, d\phi /dt> $

Vemos que el operador derivada es anti-autoadjunto. 

\

Podemos decir también  que el operador segunda derivada es adjunto pues la primera derivada pone un signo menos que es aniquilado por la segunda y al final nos queda signo positivo.  

\

En física y en química pueden aparecer operadores muy complicados. La técnica general para probar que un operador es autoadjunto se basa en resolver el problema    por pedazos:

\

Si dos operadores $T,S$ son autoadjuntos entonces su suma $T+S$ también lo es:

Por hipótesis:

$<T(u),v>= <u,T(v)>$

$<S(u),v>= <u,S(v)>$

Sumando miembro a miembro:

$<T(u),v> + <S(u),v> = <u,T(v)> + <u,S(v)>$

Como el producto interior es un producto, factorizamos:

$<Tu + Su,v>  = <u,Tv + Sv>$

$<(T + S)(u),v>  = <u,(T + S)(v)>$

Lo cual termina la demostración.

Similarmente se demuestra que un escalar real por un operador autoadjunto también es autoadjunto. Pero para escalares complejos recordemos que el producto interior de dos funciones de onda es la integral de la primera por la segunda conjugada. Por tanto, hay que tener presente que los escalares complejos pueden entrar y salir sin problema de la primera entrada del producto interior, pero entran y salen conjugados de la segunda.  Por tanto, multiplicar por un imaginario no da autoadjunto: da anti-autoadjunto. La demostración es así: 

$<iu,v>   = i<u,-iv> = <u,-iv> =  -<u,iv> $

Tomemos ahora  dos operadores   autoadjuntos y estudiemos su composición. Supongamos que

 $<T(u),v>= <u,T(v)>$

$<S(u),v>= <u,S(v)>$

entonces

 $<ST(u),v>= <T(u),S(v)> = <u,TS(v)$
 
 lo cual dice que la compuesta de dos autoadjuntos es autoadjunta sólo cuando los operadores conmutan. De lo contrario, no.

Podemos componer ahora dos operadores anti-autoadjuntos. El resultado es autoadjunto, pues cada uno de ellos aporta un menos que en conjunto se multiplican dando un mas. En particular, podemos decir: la derivada es anti-autoadjunta, la multiplicación  por $i$ es anti-autoadjunta. Por consiguiente 

\

el operador $ikd/dt$ donde $k$ es un real es autoadjunto.

\

  La importancia de ser autoadjunto radica en que de acuerdo a la tradición a partir de la interpretación de Copenhague, las mediciones experimentales de una observable corresponden a los  valores propios del operador autoadjunto que la representa. Pero operadores no autoadjuntos también se encuentran por todos lados en mecánica cuántica, el problema es que sus valores propios son imaginarios o complejos y así uno no los puede medir directamente sino que uno tiene que conformarse con el estudio de efectos colaterales. 

\

Podemos cerrar esta sección demostrando que  el Hamiltoniano de una partícula que se mueve en una dimensión bajo un potencial conservativo es un operador autoadjunto. Veamos:

\

\bigskip

El Hamiltoniano es la suma de dos operadores, el operador que representa a la energía potencial y el que representa a la energía cinética:

$H= V - (\hbar \sp{2} /2m) \partial \sp{2} /\partial x\sp{2}$

El operador $V=V(x)$ es un operador que lo único que hace es multiplicar cada función por la función escalar real $V(x)$. Es autoadjunto:

$<V\psi,\phi>=\int V\psi \phi\sp{*} dx = \int \psi V\phi\sp{*} dx  $

Pero como $V$ es real, el es igual a su propio conjugado:

$\int \psi V(x)\phi\sp{*} dx= \int \psi (V(x)\phi)\sp{*} dx= <\psi,V(x)\phi>$

\

También sabemos que  la segunda deriada es   autoadjunta. Si lo es, entonces un múltiplo real también lo es. Después decimos:   el operador Hamiltoniano es autoadjunto pues es la suma de dos autoadjuntos. Y eso termina la demostración de que el Hamiltoniano de una partícula que se mueve en una dimensión bajo un campo conservativo es autoadjunto y por consiguiente sus valores propios  deberán ser reales y por consiguiente medibles.

$>>>>>>>>>>>>>>>>>>>>>>>>>>>>>>>>>>>>>>>>>$

\color{black}

 Después de haber armado toda una maquinaria, retomemos nuestro propósito: deducir la mecánica cuántica a partir de la teoría de grupos. 
 
 Empecemos enfatizando que ahora los estados de la materia se representan por funciones de onda. Así que preguntémonos: ¿Cuáles son las condiciones básicas de un proceso dinámico sobre las funciones de onda?

Teniendo en mente el espacio de las funciones de onda, definimos un operador de evolución $E$ como una asignación que a un estado inicial y a un tiempo t, cualesquiera, les corresponde un estado, el estado en el tiempo t:

$$(\psi\sb{0},t) \rightarrow E\sb{t}(\psi \sb{0})=\psi \sb{t}
\addtocounter{ecu}{1}   \hspace{4cm} (\theecu )      $$

Al escribir la expresión $E\sb{t}$ lo que tenemos es que para cada $t$ existe un operador que transforma la condición inicial $\psi \sb{0}$ en una condición final $\psi \sb{t}$, por eso decimos que una evolución está definida por una familia uniparamétrica de operadores que ha de cumplir ciertas condiciones.

A la evolución de los sistemas físicos se le exige las siguientes propiedades:

$$E\sb{0}= I = identidad
\addtocounter{ecu}{1}   \hspace{4cm} (\theecu )      $$

$$E\sb{t} E\sb{s} = E\sb{t+s}
\addtocounter{ecu}{1}   \hspace{4cm} (\theecu )      $$

La primera condición dice que una evolución hasta el tiempo cero es lo mismo que nada, y la segunda dice que si dejamos evolucionar un sistema hasta el tiempo $s$ y después hasta el tiempo $t$, o al revés, es lo mismo que dejarlo evolucionar hasta el tiempo $t+s$. Nada más natural.

Decimos que la dinámica de un sistema es reversible si para cada $t$ existe $E\sb{-t}$. En ese caso:

$$E\sb{0}= I=E\sb{t-t}=
 E\sb{t}E\sb{-t} \Rightarrow  E\sb{-t}=E\sb{t}\sp{-1}
\addtocounter{ecu}{1}   \hspace{4cm} (\theecu )      $$

En palabras oficiales, una dinámica reversible de evolución está dada por una \index{representación de un grupo} \textbf{representación} de $R$, interpretado como tiempo, en el espacio de operadores lineales de un espacio de Hilbert dado. Formalicemos el concepto de representación.

\

\color{blue}

$<<<<<<<<<<<<<<<<<<<<<<<<<<<<<<<<$

\

\textbf{Representación de un grupo} $G$ \index{representación de un grupo}  sobre un espacio vectorial: cuando uno tiene un grupo $G$ y un espacio vectorial $V$,   decimos que tenemos una representación de $G$ sobre $V$ si a cada elemento de $G$ le corresponde un operador, el cual es una transformación lineal de $V$ en sí mismo. Por ejemplo, tenemos la recta real $R$, que forma un grupo con la suma. Tenemos ahora el plano, visto como espacio vectorial. El grupo de traslaciones en el plano es $T$. En un sentido muy natural, cada número real $r$ genera una translación, exactamente, la translación que a cada elemento lo corre $r$ unidades hacia la derecha. En general, en vez de decir que tenemos la representación de un grupo sobre un espacio vectorial, podemos decir la representación de un grupo sobre otro. Formalmente:

Sea los grupos $G,H$, decimos que $f$ es un \index{homomorfismo} \textbf{homomorfismo} de $G$ en $H$ si  $f: G\rightarrow H$ y si $f$ preserva la operación de grupo:

\

f(a+b)= f(a) + f(b)

\

Muy de anotar es que la suma en $G$ no tiene que ser la suma en $H$, ni cosa parecida: en el lado izquierdo pueden aparecer los reales con la suma y en el derecho puede tratarse de matrices con su multiplicación.  

Decimos que $f$ es una \index{representación} \textbf{representación} de un grupo $G$ sobre  un espacio vectorial $V$, si existe un homomorfismo inyectivo de $G$ sobre el espacio de los operadores   de $V$. Es decir, $f$ produce una copia al carbón de $G$ sobre los operadores de $V$. 

Ejemplo: Sea el grupo $R$ con la suma y sea $V$ el plano como espacio vectorial. Sea $M$ una matriz cualquiera $2\times 2$, la cual representa un operador sobre $V$. La función $f(r) = rM$ es una representación de $R$ en $V$ puesto que

$f(r+s) = (r+s)M = rM + sM = f(r) + f(s)$ 

\

$>>>>>>>>>>>>>>>>>>>>>>>>>>>>>>>>>>>>>>>>>$

\

\color{black}

Ahora vamos a desarrollar la siguiente propuesta: la mecánica cuántica no es más que una representación del  los números  reales, vistos como grupo con su suma,  sobre el espacio de Hilbert de las funciones de onda.

\

Veamos eso con detalle.

\

Nuestro grupo está conformado por los reales, el cual representa el tiempo. Y, ¿por qué es tan natural que una representación del tiempo sea todo lo necesario para descubrir la mecánica cuántica?

Sin desear resolver el enigma del tiempo, podemos pensar que la naturaleza tiene a nivel fundamental una dinámica de la cual el ser humano puede abstraer el concepto de tiempo como 'esa entidad que fluye linealmente'.  

\

Veamos las implicaciones de lo que es una representación sobre las funciones de onda. 

La norma cuadrado de una función de onda debe representar la integral sobre todo el universo de una densidad de probabilidad. Por lo tanto, el operador de evolución de un sistema cuántico no puede crear ni destruir la probabilidad total, la cual debe ser uno. Como la probabilidad total asociada a una función de onda $\psi$ es la norma cuadrado en $L\sb{2}$, entonces las funciones de onda son elementos de su bola unitaria. El operador de evolución debe moverse por la bola unitaria, o mejor dicho, no modifica la norma, por eso decimos que es una isometría:

$$<\psi\sb{1},\psi\sb{2}> = <E\sb{t} \psi\sb{1},E\sb{t} \psi\sb{2}>
\addtocounter{ecu}{1}   \hspace{3cm} (\theecu ) $$

Invocando el operador adjunto:

$$<E\sb{t} \psi\sb{1},E\sb{t} \psi\sb{2}> =
< \psi\sb{1},E\sb{t}\sp{\dagger} E\sb{t} \psi\sb{2}> \addtocounter{ecu}{1}   \hspace{3cm} (\theecu )      $$

O sea

$$<\psi\sb{1},\psi\sb{2}> =  < \psi\sb{1},E\sb{t}\sp{\dagger} E\sb{t} \psi\sb{2}> $$

Como eso debe ser cierto para toda función de onda, podemos poner $\psi\sb{1}  = \psi\sb{2}$ y nos queda:

$$1= <\psi\sb{2},\psi\sb{2}> =  < \psi\sb{2},E\sb{t}\sp{\dagger} E\sb{t} \psi\sb{2}> $$

entonces tiene que cumplirse que $E\sb{t}\sp{\dagger} E\sb{t}=I$, y como $E\sb{t}$ tiene inverso, se cumple que $E\sb{t}\sp{\dagger} = E\sb{t}\sp{-1}$. Es decir que estamos hablando de un operador de evolución  unitario.

Lo que nosotros hemos llamado una dinámica de evolución es realmente una familia uniparamétrica de operadores con propiedades de grupo cuya operación es la composición: es cerrado, tiene una identidad y es conmutativo y asociativo, tiene inverso.  La representación del tiempo la hemos logrado mediante un grupo abeliano: a cada tiempo dado le corresponde un  operador unitario que tiene inverso.

La unitariedad de la evolución cuántica implica que la existencia del ser ni se crea ni se destruye. Puede elegirse entre pensar que eso siempre ha sido así, o que todas las leyes del ser y de la física nacieron con el Big-Bang, o que se fueron creando a medida que el Big-Bang se fraguaba y desenvolvía, o  que el Big-Bang es tan sólo una creencia emocionante que pronto pasará a la historia. En todos los casos, el misterio de la existencia es el más grande de todos los misterios al cual no se le puede añadir ni quitar nada.

\

Todas las condiciones enunciadas y satisfechas por la dinámica de evolución nos invitan a aplicar el Teorema de Stone:  así como un número complejo $z$ de módulo uno puede expresarse como  $z=e\sp{i\theta}$, un operador unitario  admite una expresión prácticamente equivalente. Tratándose de operadores, una expresión así sería algo maravilloso. Para ver que también puede ser algo natural, tomemos unos minutos para una observación.

\

\

\color{red}

La ecuación diferencial de la evolución:

\

Para una función $f$  de $R$ en $R$ definimos la derivada como

$f'(x) = lim_{ h\rightarrow 0} (f(x+h) -f(x))/h$

Notemos que para definir la derivada de $f$ se requiere  que en el dominio se pueda sumar, que en el codominio se pueda restar y multiplicar por un número real, y que la noción de límite exista.

Por lo tanto, si tenemos una función de $R$ en algún  espacio vectorial completo, estamos hechos: también podemos definir la derivada. Nuestra dinámica tiene todas esas buenas cualidades, pues los operadores sobre un espacio de Hilbert se pueden restar y multiplicar por un escalar y el mismo espacio es completo, o sea que los límites que debieran existir existen en realidad. Calculemos entonces la derivada de la evolución cuántica, formulemos una ecuación diferencial y resolvámosla. De antemano podemos inferir qué nos ocurrirá: la evolución cumple la propiedad : la evolución hasta $a$ y después hasta $b$ debe ser lo mismo que evolución hasta $a+b$:

$E(a+b) = E(a)E(b)$

Esas son exactamente las mismas propiedades de la exponencial. Bueno, pues eso es lo que buscamos. Veamos entonces:

$E'(0) = lim_{h\rightarrow 0} (E(0+h) -E(0))/h$

      $=  lim_{h\rightarrow 0} (E(h) -E(0))/h$

      $=  lim_{h\rightarrow 0} (E(h) -I)/h$

donde $I$ es la identidad, pues la evolución hasta cero deja las cosas como están.

\

En general:

$E'(t) =  lim_{h\rightarrow 0}(E(t+h) -E(t))/h$

      $=  lim_{h\rightarrow 0} (E(t) E(h) -E(t))/h$

      $=  lim_{h\rightarrow 0}E(t) (E(h) -I)/h$

      $=  E(t)  lim_{h\rightarrow 0}(E(h) - I )/h$

      $ = E(t)E'(0)$

 \

Por consiguiente, la evolución es derivable en toda momento ssi es derivable en cero. Lo que hemos hecho es válido en ámbitos muy generales, donde además se pueda factorizar: cuando eso puede hacerse se tiene un álgebra. Sea el álgebra de los números reales o de los operadores lineales. Si se tratara de reales, la derivada en cero sería un número, $\lambda$, que se puede poner a la derecha o a la izquierda, pues la multiplicación por reales es conmutativa:

$ E'(t)= E(t)E'(0)=E(t)\lambda $

de aquí sale inmediatamente que $E(x)$ debe ser la función exponencial y

$E(t)=e\sp{\lambda t}$

Cuando estamos en el álgebra de operadores, la derivada en cero es un operador, digamos $A$, al cual se le da el nombre de \index{generador infinitesimal} \textbf{generador infinitesimal de la evolución}:

$ E'(t)= E(t)E'(0)=E(t)A $

de lo cual inferimos que

$ E(t)= e\sp{At} = e\sp{tA} $

como $t$ es un número real se puede poner a cualquier lado.
Eso estaría muy bien si supiésemos que significa la exponencial de un operador. Resulta que todo lo válido para reales respecto a la exponencial se puede extender, con algunas precauciones por la no conmutatividad, al álgebra de operadores lineales:

$e\sp{tA}= \sum \sb{k} (tA)\sp{k}/k! = \sum \sb{k} t\sp{k}A\sp{k}/k!$

$A\sp{k}$ es el operador que resulta de iterar o componer a $A$ $k $ veces sobre sí mismo. En ese caso el producto es conmutativo puesto que $A$ conmuta consigo mismo.  

\

\

\

Lo que hemos estado haciendo fue justificado rigurosamente por Stone quien nos legó el siguiente teorema: 

\color{red}

\

\addtocounter{ecu}{1}

\textit{Teorema de Stone \theecu}

\

 Todo operador de evolución, que sale de una representación  unitaria continua de los reales dentro del algebra de operadores lineales de un espacio de Hilbert, tiene un generador infinitesimal que es de la forma $-iH/\hbar $ donde $H$ es un operador autoadjunto. Mejor dicho:

 $$E\sb{t}= e\sp{-itH/\hbar}
\addtocounter{ecu}{1}   \hspace{4cm} (\theecu )      $$

es decir:

$$\psi(x,t) = E\sb{t} \psi(x,0)= e\sp{-itH/\hbar} \psi(x,0)
\addtocounter{ecu}{1}   \hspace{4cm} (\theecu )      $$

 donde $H$ es autoadjunto: $<H\psi,\phi>=<\psi,H\phi>$.

Por supuesto,

$$e\sp{-itH/\hbar}= \sum \sb{k} (-itH/\hbar)\sp{k}/k!
\addtocounter{ecu}{1}   \hspace{4cm} (\theecu )      $$

\

Podemos ver la necesidad de que $H$ sea autoadjunto para que el operador de evolución sea unitario.  Si el operador de evolución es unitario, entonces

$1 = E^\dagger-{t} E\sb{t} = (e\sp{-itH/\hbar})^\dagger e\sp{-itH/\hbar}$ 
$= e\sp{itH/\hbar} e\sp{-itH/\hbar} $
$ = e\sp{it(H^\dagger - H)/\hbar}$

Por tanto $H^\dagger - H = 0$, lo cual implica que $H^\dagger = H$, es decir, que $H$ es autoadjunto.

\

\color{black}

En este  teorema de Stone la exponencial que resulta es imaginaria. ¿De donde resulta el número $i$? En realidad, eso no implica nada especial pues uno siempre tiene que   todo operador $T$ se puede reescribir como $T = -i\sp{2}T = -i(iT)$. Lo que sí es importante es tener presente que el operador $H$ del teorema es autoadjunto.  

El signo menos da la impresión de convergencia, pero no hay que presumir pues se trata de un exponencial imaginaria, la cual causa rotaciones: en general, debido a que estamos trabajando con espacios de dimensión infinita, el análisis funcional, la materia encargada de estudiar estos temas, es algo delicado y lleno de sorpresas. Por ejemplo, por ninguna razón es obvio que la representación del tiempo tenga que ser continua y por eso hay que ponerla como una condición explícita (de hecho, no ha de faltar alguien que niegue  que el tiempo fluye y que prefiera pensar que da saltos, los cuales pueden ser hacia adelante o hacia atrás y sin correlación alguna entre un lugar y otro).

Hay por lo menos dos maneras importantes y no siempre equivalentes de definir la continuidad cuando se trata de operadores. El primero comienza caracterizando el tamaño de cada operador por la máxima distorsión que causa sobre la bola unitaria de su dominio. Esa forma de medir es una norma que permite definir una distancia, como la norma de la resta. Cuando se usa esta norma para definir la continuidad como una propiedad de no causar roturas, hablamos de continuidad en el sentido fuerte o de la norma.

La otra forma de definir continuidad se denomina débil e involucra el espacio dual, o sea el conjunto de funciones escalares lineales definidas sobre  el espacio de Hilbert en cuestión. La convergencia fuerte,  en el sentido de la norma, implica la convergencia débil, pero no al revés.

 Nosotros ignoramos estas sofisticaciones en nuestra discusión por lo que a veces se notarán vacíos en nuestros argumentos. Lo mejor sería remitirse a libros que también tengan análisis funcional.

\

Podemos deducir la  ecuación de Schr\"{o}dinger a partir de nuestra evolución. Tenemos:

$E\sb{t} \psi\sb{0}= \psi\sb{t} $
$=\psi(x,t)= e\sp{-itH/\hbar} \psi\sb{0}$

Derivando a ambos lados con respecto al tiempo queda:

$\partial \psi(x,t)/\partial t = (-iH/\hbar) e\sp{-itA} \psi\sb{0}$
$= (-i/\hbar)H\psi(x,t)$

Por tanto tenemos  que:

$$\partial \psi(x,t)/\partial t  = (-i/\hbar)H\psi(x,t)
\addtocounter{ecu}{1}   \hspace{4cm} (\theecu )      $$

que también se puede escribir como la ecuación de Schr\"{o}edinger:

$-(\hbar/i)\partial \psi/\partial t = H \psi $

\

En conclusión: la dinámica cuántica  que sale al representar los reales sobre el espacio de Hilbert de las funciones de onda está dada por

$$E\sb{t} \psi\sb{0}= \psi\sb{t}
=\psi(x,t)= e\sp{-itH/\hbar} \psi\sb{0}
\addtocounter{ecu}{1}   \hspace{4cm} (\theecu )      $$

donde $H$ es un operador autoadjunto. De acá podemos deducir la forma de la ecuación de Schröedinger:

$-(\hbar/i)\partial \psi/\partial t = H \psi $

\

De cierta manera, hemos resuelto la ecuación de Schr\"{o}edinger y lo hemos hecho de tal forma que nos es factible programar un computador para hallar una solución aproximada. Anteriormente ya habíamos encontrado otro método, pero en mi opinión, el que acabamos de hallar es mejor pues para un tiempo fijo se puede estimar el error de la inexactitud usando el criterio de Cauchi para convergencia de series. Un detalle que ayuda a entender eso es el siguiente:

Se dice que una serie es de Cauchi si para todo $\epsilon >0$ existe un $N$ tal que si $n, m > N $ entonces $|a\sb{n} -a\sb{m}| <\epsilon$. En un espacio completo, toda serie de Cauchi converge a un límite, sea $L= lim a\sb{n}=a\sb{\infty}$. Entonces tenemos: $|a\sb{n} -a\sb{m}| = |a\sb{n} -a\sb{\infty} + a\sb{\infty} - a\sb{m}| = |a\sb{n} +a\sb{\infty}| + | a\sb{\infty} - a\sb{m}| <2\epsilon$. Por consiguiente $|a\sb{n} +a\sb{\infty}|  <2\epsilon - | a\sb{\infty} - a\sb{m}| < 2\epsilon. $

\

Ahora debemos hallar la expresión para el operador autoadjunto $H$.  Y debemos hacerlo a partir de la teoría de grupos. Hagámoslo.

Demostremos primero que  $H$ es energía pura. Para ello, podemos  verificar sus unidades físicas. Desarrollando la exponencial a primer orden, tenemos:

$\psi(x,\triangle t)= e\sp{-i \triangle t H/\hbar} \psi\sb{0}$
$= (1-i\triangle tH/\hbar)\psi\sb{0}$

Para igualar unidades de medición, $\triangle t H $ y $ \hbar $ deben tener las mismas unidades, para que el quebrado  $\triangle tH/\hbar$ no tenga ninguna unidad y se pueda restar del número 1. Eso implica necesariamente que $H$ viene en unidades de energía. En efecto, energía = trabajo =  fuerza x espacio = masa x aceleración x espacio = masa x (espacio/tiempo cuadrado) x espacio = masa x espacio cuadrado / tiempo cuadrado. Por lo tanto, $\hbar$ debe venir en unidades de masa x espacio cuadrado /tiempo, lo cual es cierto.

De lo dicho se deduce que el operador $H$ debe ser el operador Hamiltoniano cuántico. Probemos que a partir de la teoría de grupos podemos deducir la forma del Hamiltoniano para una  partícula en un potencial.

\

En retrospectiva podemos decir:

Como un adolescente que logra subirse a un carro ya en marcha, así fuimos nosotros con la mecánica cuántica: primero fue inventada por    Schr\"{o}edinger y sus contemporáneos y después fue reformulada por Feynman, quien definió nuestro punto de partida.  

Todo eso ha sido muy complejo desde variados puntos de vista. Pero en la sección pasada vimos que la ecuación de Schr\"{o}edinger puede deducirse a partir casi de la nada, postulando únicamente que 

\begin{itemize}
\item Los estados de un sistema físico se representan por funciones a valores complejos cuyo cuadrado tenga norma uno. Tales funciones forman un espacio de Hilbert, el espacio de las funciones de onda. 
\item Las leyes naturales de evolución de un sistema físico cumplen con una representación continua de los números reales con la suma sobre los operadores del espacio de Hilbert de las funciones de onda con la composición. 
\end{itemize}

\bigskip

Para enfrentar nuestra tarea de especificar la forma del Hamiltoniano a partir de la teoría de grupos, investiguemos la siguiente pregunta:  ¿Qué pasará si estudiamos la representación de los reales interpretados no como tiempo sino como el espacio de las traslaciones espaciales?  

Lo que hicimos antes fue estudiar el efecto de las traslaciones en el tiempo y dedujimos que

$\psi(x,\triangle t)= e\sp{-i \triangle t H/\hbar} \psi\sb{0}$
$= (1-i\triangle tH/\hbar)\psi\sb{0}$

Ahora bien, $\triangle t H $ y $ \hbar $ deben tener las mismas unidades, para que el quebrado no tenga ninguna unidad. Eso implica necesariamente que $H$ viene en unidades de energía. 

Si repetimos el mismo procedimiento, pero interpretamos ahora los reales no como tiempo sino como espacio, llegamos al siguiente resultado: 

$\psi(x + \triangle x, t)= e\sp{-i \triangle x H/\hbar} \psi\sb{0}$
$= (1-i\triangle x H/\hbar)\psi\sb{0}$

Ahora tenemos que $H$ es un operador autoadjunto cuyas dimensiones multiplicadas por espacio se igualan a las dimensiones de $\hbar$. 

Como $\hbar$ viene en unidades de masa x espacio cuadrado /tiempo, entonces el nuevo operador $H$ debe venir en unidades de masa x espacio   /tiempo, es decir, masa por velocidad, es decir, momento, que se nota $P$. Tenemos entonces que a primer orden:

$\psi(x + \triangle x, t)= e\sp{-i \triangle x P/\hbar} \psi(x,t)$
$= (1-i\triangle x P/\hbar)\psi(x,t)$

Hemos ligado a las traslaciones espaciales con el momento. Suena excelente. Busquemos ahora la forma de $P$:

 $\psi(x + \triangle x, t) - \psi(x, t) $
$= (1-i\triangle x P/\hbar)\psi(x,t) - \psi(x, t)  $

$\psi(x + \triangle x, t) - \psi(x, t)= -(i/\hbar)\triangle x P \psi(x,t)$

Dividiendo por $\triangle x$ y tomando el límite obtenemos:

$\partial \psi(x,t)/\partial x =  -(i/\hbar)  P  \psi(x,t)$

 $i\partial \psi(x,t)/\partial x =   (1/\hbar) P  \psi(x,t)$
 
  $i\hbar \partial \psi(x,t) /\partial x =    P \psi(x,t)$

y como esto debe ser cierto para toda $\psi(x,t)$ podemos decir que 

$P = i\hbar \partial  /\partial x$

\bigskip

Notemos que este operador $P$ es autoadjunto pues es la compuesta de dos anti-autoadjuntos, la derivada y la multiplicación por un imaginario. 

\

Sin embargo, tenemos un problema, y es que en todos los libros aparece la representación de $P$, no como la tenemos sino con signo menos:

 $P = -i\hbar \partial  /\partial x$

y esta nueva representación del momento no es autoadjunta sino anti-autoadjunta. ¿Qué podrá estar pasando? Puesto que no he encuentro error aritmético en mis procedimientos, mi propuesta es interpretativa: 

Estamos buscando una representación de los reales entendidos como traslaciones en el espacio. Esas traslaciones generan velocidad, generan momento, tal como nos dio. Pero hay un detalle más, a saber:

Tomemos una función de onda $\psi(x,t)$: si queremos trasladarla $\triangle x$ hacia la derecha, hacia el lado positivo del eje $X$, tenemos que someterla a un traslación como sigue:

$\triangle x \rightarrow \psi(x,t)\rightarrow  \psi(x-\triangle x,t)$

De tal manera que debemos releer el teorema de Stone de la siguiente manera:

$\psi(x-\triangle x,t) =   e\sp{-i \triangle x P/\hbar} \psi(x,t)$
$= (1-i\triangle x P/\hbar)\psi(x,t)$

Despejemos $P$:

$\psi(x,t)-\psi(x-\triangle x,t) =   \psi(x,t)-
  (1-i\triangle x P/\hbar)\psi(x,t)$

Dividiendo por $\triangle x$ y tomando el límite hacia cero, nos queda:

$\partial \psi(x,t) /\partial x = i\triangle x P/\hbar \psi(x,t)$

de donde obtenemos $P$ multiplicando por $i$ y $\hbar$ en ambos lados:

$i\hbar \partial \psi(x,t) /\partial x =   - P \psi(x,t)$

y como esto debe ser cierto para toda $\psi(x,t)$ podemos decir que 

$P = -i\hbar \partial  /\partial x$

Adoptaremos esta interpretación, mediada por una anti-representación del espacio, como la oficial. Pero es de notar que uno puede tomar cualquiera de las dos interpretaciones y lo que sigue no podrá distinguirlo. 

\

Conociendo la representación del momento, masa por velocidad,  podemos inmediatamente formular la representación de la energía cinética $K$:

$K = \frac{1}{2} m v^2 = \frac{1}{2m}  p^2 $

Ahora bien, si el momento se representa con una derivada, ¿cómo ha de representarse el momento cuadrado? Ha de representarse por un  operador autoadjunto que tenga dimensiones de momento cuadrado. La elección natural es el operador al cuadrado, es decir, él aplicado sobre él mismo. Proponemos entonces que la energía cinética se represente por 

 $K =   \frac{1}{2m} (-i\hbar \partial  /\partial x)^2 $
 
 $ K = -\frac{\hbar^2}{2m}   \partial^2  /\partial x^2 $ 

Si hay energía potencial constante, $V$, podemos proponer su representación por el operador autoadjunto

$\psi(x,t) \rightarrow V\psi(x,t)$

Al final nos queda que la energía cinética  puede representarse por el operador Hamiltoniano que ya conocíamos:

$K + V  \rightarrow -\frac{\hbar^2}{2m}   \partial^2  /\partial x^2  + V$

\bigskip

Hemos demostrado que el operador Hamiltoniano tiene su base en la representación del grupo uniparamétrico de los reales con una doble lectura: como el grupo de las traslaciones en el tiempo y como el grupo de las traslaciones en el espacio. 

Seguramente no tendremos problemas en generalizar nuestra ecuación a tres dimensiones. Pero hay que notar que la generalización será nuestra responsabilidad y no la del teorema de Stone, pues ese teorema nos permite tomar grupos uniparamétricos, mientras que el las traslaciones en el espacio forman un grupo triparamétrico. Lo que hay que hacer es generalizar el teorema de Stone a grupos cualesquiera y llegamos a la necesidad de estudiar una teoría general de la representación de grupos.

El significado de los grupos en mecánica clásica nos lo da el teorema de Noether: si la acción es invariante ante un grupo, ahí hay una entidad que se conserva: para las traslaciones espaciales, tenemos el momento.  Para las rotaciones, tenemos la conservación del momento angular.

Todo eso se resume diciendo, informalmente,  que la mecánica clásica es generada por el \index{grupo de Galileo} \textbf{grupo de Galileo}, el cual consta, por definición, de: traslación en el tiempo, traslaciones en el espacio y rotaciones en el espacio. Este grupo tiene 7 parámetros.

El   programa que hemos iniciado nos conduce a predecir que: la mecánica cuántica resulta de la representación en un espacio de Hilbert complejo del grupo de Galileo. Para completar dicho programa nos faltaría estudiar cómo se representa el grupo de las rotaciones. Sucede que las rotaciones en el espacio se representan naturalmente por matrices. Po ende, una representación natural de las rotaciones en un espacio de Hilbert complejo también ha de ser por medio de matrices. Eso quiere decir que los estados de nuestras partículas han de ser representados por vectores con varias coordenadas. Este estudio fue llevado a cabo por Pauli, quien nos legó la noción de spin que es el equivalente cuántico del momento angular. 

Como nuestro interés está en el campo electromagnético, no nos conviene profundizar en todos los temas al mismo tiempo. Por eso, nos contentaremos con el estudio de las implicaciones de nuestra mecánica cuántica que describe partículas por medio de funciones de onda con una coordenada. Y como veremos, con eso ya tenemos bastante.

\section{PODER PREDICTIVO}

La belleza matemática  de un formalismo físico no lo convierte en verdad. La única entidad que valida un formalismo es la naturaleza misma. Uno de los grandes logros de la mecánica cuántica fue demostrar que podía predecir las observaciones espectroscópicas de átomos excitados que irradiaban luz visible, dando una clara sustentación al modelo atómico planetario,  estudiada exhaustivamente tanto para el hidrógeno como para el Helio. Fue la espectroscopía la que puso el visto bueno a la propuesta de interpretar los valores propios del Hamiltoniano (en general, debe tomarse el espectro del operador dado) con los valores de la energía de excitación de un átomo. Más concretamente:

Eso de que $H$ es autoadjunto es excelente pues sus valores propios toman valores reales y entonces pueden interpretarse directamente como valores experimentales de la energía. De hecho, la mecánica cuántica conquistó al mundo cuando quedó claro la verificación de las predicciones cuánticas para el átomo de hidrógeno al ser contrastadas con las diferencias energéticas de los niveles de energía observadas en espectroscopía.

En realidad, los niveles de energía son inaccesibles al experimento, pero en cambio un cambio en los niveles produce la emisión de un fotón, cuya energía puede medirse muy exactamente: la energía da el color del fotón y ese puede separarse con un simple prisma. Por tanto, el contacto experimental entre nuestra matemática y la espectroscopía se hace calculando los valores propios de la energía, restándolos entre ellos y calculando la longitud de onda asociada al dicha diferencia.

Pero además, un espectro de emisión  indica la intensidad con que se irradia a una frecuencia determinada. Por lo tanto, las predicciones cuánticas deberían decir no sólo las frecuencias sino además su intensidad.

Como ya hemos dicho, las predicciones cuánticas son probabilísticas y un espectro debe leerse como una densidad de probabilidad. Por consiguiente preguntamos: ¿con qué probabilidad se producirá determinada emisión? Esa respuesta fue predicha por von Nuemann antes de que se supiera que la mecánica cuántica iba a explicar el espectro de radiación atómica.

Cabe advertir que la formulación correcta de los valores observables no se hace sobre el conjunto de valores propios sino sobre otro más general denominado el espectro. Nosotros nos hacemos los de la vista gorda pues cuando el espectro es discreto, de valores aislados, entonces coincide con el conjunto de valores propios. Nuestro caso es el de la luz y todo indica que viene cuantizada, o sea que estamos en un caso en el cual el espectro es discreto y por lo tanto con los valores propios nos basta.

\

Tenemos entonces el siguiente postulado:

Con qué probabilidad ocurrirá la medición de un valor propio $\lambda$? Con la determinada por la función de onda. Si $(\lambda, \psi)$ es un par propio del Hamiltoniano, entonces si el sistema está en el estado $\phi$, la probabilidad $p$ de que al medir la energía se encuentre el valor $\lambda$ es:

$$p= |<\psi,\phi>|\sp{2}
\addtocounter{ecu}{1}   \hspace{4cm} (\theecu )      $$

\

Para nosotros, es importante que el modelo atómico planetario haya sido validado por la mecánica cuántica: es la primera vez que podemos aspirar a entender la existencia de cuerpos neutros y estables. Podremos tener paz.. ..por algún tiempo.

\section{UN MISTERIO}

La vida sin misterio no es nada. Aunque la mecánica cuántica de partículas  ha logrado explicarnos lo que desde la mecánica clásica parecía inexplicable, como la existencia de átomos que forman las unidades básicas de cuerpos neutros, es necesario saber que el misterio se sigue. En efecto:

\

Teorema: Las reflexiones y la identidad forman un grupo de orden 2, $GR$. El  \textbf{Hamiltoniano molecular} \index{Hamiltoniano molecular} es invariante ante $GR$.

\

Demostración: El Hamiltoniano describe  a la energía y como tal tiene dos partes. La primera es  debida a la energía potencial de interacción entre las partículas, la cual tiene una dependencia única y exclusiva de la distancia. Dicha distancia es un elemento objetivo y en nada se ve afectado por los adjetivos subjetivos, derecha, izquierda, que son los que definen a $GR$.

La otra parte, la correspondiente a la energía cinética,  tiene un segunda derivada. La segunda derivada es la derivada de la primera derivada. La primera derivada, en una dimensión, indica la pendiente. Por lo tanto, su segunda derivada indica la forma como cambia la pendiente, o sea la forma como se curva la gráfica. Por eso a la segunda derivada se le llama curvatura. Ese es otro elemento objetivo que tampoco discierne entre derecha e izquierda. De ahí resulta la invariancia solicitada.

Como el Hamiltoniano es invariante ante $GR$, entonces las funciones de onda que son solución  a la ecuación de Schr\"{o}edinger tienen que ser invariantes ante cambios mentales de derecha-izquierda. Por consiguiente, todos los sistemas moleculares deberían ser simétricos.

\

Pues bien, la química nos ha enseñado que todas las biomoléculas son asimétricas. La gran mayoría  tienen un tipo de asimetría que se llama levógira. La más popular entre las dextrógiras es la D-dextrosa, aunque hay muchas que en el momento necesario son fabricadas de sus versiones levógiras.

Vemos entonces que la química ha falsificado a la mecánica cuántica que considera que una molécula está compuesta de sus átomos que la forman y de nada más. No se sabe qué es lo que pasa. Sin embargo, se cree que es suficiente modificar el Hamiltoniano incluyendo a la radiación electromagnética, la cual tiene infinitos grados de libertad. Debido a eso, un sistema descrito por el nuevo Hamiltoniano puede sufrir una \index{rotura  de la simetría}  \textbf{rotura espontánea de la simetría} (dejar de ser invariante ante la totalidad del grupo, en este caso $GR$ ) y producir estructuras asimétricas. Con todo, no se ha podido demostrar ni que eso es falso ni que es verdadero.

No es irrazonable pensar que con el campo electromagnético incluido dentro del Hamiltoniano molecular se pueda explicar la chiralidad, como se llama a la asimetría geométrica:  el campo electromagnético viene cuantificado y sus partículas tiene spin, o modo de girar, si se habla a la usanza clásica, lo cual ya se\~{n}ala una dirección como predilecta y el espacio comienza a tener una orientación objetiva. Claro que el spin orienta al espacio de la partícula que lo posee direccionando 'adelante vs atrás', que también se parece a   'izquierda vs derecha' que es lo que aparece en la chiralidad.

\section{EL PRINCIPIO DE INCERTIDUMBRE}

Nuestras elucubraciones teóricas necesitan ser relacionadas con el experimento para poder convertise en ciencia. Necesitamos por ello reglas de interpretación que nos liguen operadores con valores experimentales. El contexto podría ser como sigue:

Los valores propios de un operador autoadjunto son reales, por lo tanto se pueden conectar directamente con mediciones experimentales. Los espacios propios asociados a valores propios diferentes son ortogonales y definen una base de todo el espacio. Si los valores propios son repetidos, de toda formas uno puede formar una base ortogonal para todo el espacio. Gracias a ellos podemos hacer aproximaciones numéricas y tratarlas en el computador.

La mecánica cuántica no puede predecir cuál de todos los valores propios o de sus diferencias será medida en el próximo experimento. Tan sólo se pueden predecir las probabilidades con que aparecerán los valores propios en mediciones experimentales independientes. 

Todo parece cuadrar bien si tomamos la siguiente regla de interpretación de lo que es una medición:

Al comenzar el experimento, el sistema observado está descrito por la función de onda $\phi$, que puede ser cualquiera. Se mide algo. Ese algo está representado en mecánica cuántica por una observable, un operador autoadjunto, $B$, con pares propios $(\lambda_i, \psi_i)$. El proceso de medición 'arroja' al sistema desde $\phi$ hasta alguno de los vectores propios de $B$, digamos $\psi_i$ y el resultado de la medición es el valor propio  $\lambda_i$. Ahora bien, no hay forma de predecir a cuál vector propio será arrojado el sistema. Lo único que puede predecirse es la probabilidad $p_i$ con la cual eso pasará:

Si los pares propios ( en el espacio proyectivo) de un operador autoadjunto $\textbf{B} $ son $(\lambda \sb{i}, \psi\sb{i})$ entonces la probabilidad $p\sb{i}$ de que al medir $\textbf{B} $ en un sistema que fue preparado sobre la función de onda $\phi$ se obtenga el valor  $\lambda \sb{i}$ es

$p\sb{i}= |<\psi\sb{i},\phi>| \sp{2}$

En particular, si el sistema fue preparado en $\psi_i$ y después de ser medido se establizó en $\psi_i$, entonces se mide $\lambda_i$ con probabilidad $|<\psi\sb{i},\psi_i>| \sp{2}$ = 1, puesto que todas las funciones de onda tienen norma uno, pues representan una partícula. Eso implica que si un sistema está en un valor propio, el operador de medición no lo cambia de estado. 

Ahora bien, si un sistema fue preparado en $\phi$, la suma de todas las  probabilidades debe sumar uno sobre todos los $\lambda_i$. ¿Es eso cierto? Sí. Veamos por qué:

Tenemos que demostrar que 

$1 = \sum_i p\sb{i} = \sum_i |<\psi\sb{i},\phi>| \sp{2} =  \sum_i <\psi\sb{i},\phi> <\psi\sb{i},\phi>^*  $ 

Para ello, expresemos la función de onda $\phi$ en la base propia:

$\phi = \sum c_i \psi_i$

multiplicando en cada lado por la derecha por $\psi_j$ obtenemos:

$<\phi,\psi_j > = \sum c_i <\psi_i,\psi_j>$

por la ortogonalidad se aniquilan todos los términos de la suma excepto el j-ésimo. Obtenemos:

$<\phi,\psi_j > =   c_j$

Por consiguiente, la descomposición de $\phi$ en la base se lee:

$\phi = \sum <\phi,\psi_i > \psi_i$

La norma de $\phi$ es uno y su cuadrado también:

$1= ||\phi||^2 = ||\sum <\phi,\psi_i > \psi_i||^2$

$ = \ \ <\sum <\phi,\psi_i > \psi_i,\sum <\phi,\psi_i > \psi_i>$

Por ortogonalidad sólo sobreviven los índices repetidos:

$1 = \sum_i <\phi,\psi_i >  <\phi,\psi_i >^* <\psi_i.\psi_i>$

pero como la norma cuadrado de cada función de onda es uno, se simplifica:

$1 = \sum_i <\phi,\psi_i >  <\phi,\psi_i >^*  $

que era lo que queríamos demostrar.

\ 

Hemos asumido que los vectores propios de un operador autoadjunto tiene un número enumerable de vectores propios que forman una base y por eso usamos una suma al expresar un vector cualquiera en dicha base. Decimos que el operador tiene un espectro discreto. Pero en otras ocasiones   resulta un espectro continuo y se usa una integral en vez de la suma. El estudio riguroso de esta temática se basa sobre el teorema espectral, que puede ser encontrado en los libros de análisis funcional. 

\

Definimos el \index{valor esperado   de un operador } \textbf{valor esperado $E(A)$ de un operador $A$ dado una función de onda $\psi(x)$} como el valor medido por $A$ (por un instrumento que mida $A$) si el sistema se prepara en  $\psi(x)$ y si se registra en  $\psi(x)$ después de la medición:

$E(A(\psi(x))) = <A\psi(x), \psi(x)> = \int( A\psi(x)) \psi^*(x)dx$

El valor esperado o esperanza cumple las propiedades:

a) $E(A+B) = E(A) +E(B)$, $A,B$ operadores.

b) $E(kA) = kE(A)$, $A$ operador, $k$ una constante. 

c) $E(c) = c$, para el operador que multiplica por $c$,  una constante cualquiera.

d) La esperanza evaluada sobre un vector propio da el valor propio correspondiente.

e) El valor esperado de un operador autoadjunto sobre un vector propio es un número real.

f)El valor esperado de un operador anti-autoadjunto sobre un vector propio es un imaginario puro. Decimos que $A$ es anti-autoadjunto si $A^\dagger = -A$.

\

La demostración de la cuarta propiedad (d) es como sigue: consideremos que $\psi(x)$ es vector propio de $A$: $A\psi(x) = \lambda \psi(x)$ con $\lambda $ complejo. Tenemos: 

$E(A(\psi(x))) = \int (A\psi(x)) \psi^*(x)dx = \int \lambda \psi(x) \psi(x)^* dx = \lambda\int  \psi(x)\psi^*(x) dx $

$= \lambda$

\

Demostración de la quinta propiedad (e): un operador autoadjunto tiene todos sus valores propios reales y, por la propiedad anterior, la esperanza del operador sobre un vector propio es el valor propio correspondiente, que es real. 

\

Para demostrar la última propiedad, comenzamos averiguando cómo son los valores propios de un operador anti-autoadjunto:

Si $A^\dagger = -A$ y si $A\psi = \lambda \psi$ entonces:

$<A\psi, \psi > = <\lambda\psi, \psi > = <  \psi, A^\dagger \psi >  = < \psi, -A\psi > = -< \psi,\lambda \psi > = -\lambda^* < \psi,\psi > $

Obtenemos que

$<\lambda\psi, \psi > = -\lambda^* < \psi,\psi > $

es decir:

$\lambda<\psi, \psi > = -\lambda^* < \psi,\psi > $

de donde 

$\lambda = -\lambda^*$

lo cual dice que $\lambda$ es un imaginario puro, como en $i = -i^* = -(-i) = i$.

Ahora bien,  si calculamos la esperanza de $A$ sobre un vector propio  obtenemos:

$E(A) = \int (A\psi(x)) \psi^*(x) dx = \int (\lambda\psi(x)) \psi^*(x) dx =\lambda \int (\psi(x)) \psi^*(x) dx  = \lambda $

lo cual da imaginario puro o un múltiplo de $i$: la esperanza de un operador anti-autoadjunto relativa a un vector propio es un imaginario puro. 

\

Como el proceso de medición es caótico, dando a veces una cosa y a veces otra, usamos una medida de la dispersión de los resultados como sigue. Definimos ahora la \index{varianza  de un operador} \textbf{varianza $V_A$ de un operador $A$} dada una función de onda $\psi(x)$:

$V_A = E[(A-E(A))^2] $

La varianza cumple con las propiedades:

a) $V_A = E(A^2) - E^2(A)$

b) La varianza evaluada sobre un vector propio es cero.

\

La demostración de la primera propiedad es la siguiente:

$V_A = E[(A-E(A))^2]  = E( A^2-2AE(A) + E^2(A))   $

$= 
 E( A^2)-2E(A)E(A) + E(E^2(A)) =  E( A^2)-2E^2(A)  +  E^2(A)$
 
 $  = E(A^2) - E^2(A)$

Para demostrar la segunda propiedad, consideremos que $\psi(x)$ es vector propio de $A$: $A\psi(x) = \lambda \psi(x)$. Tenemos: $A^2\psi(x)  = \lambda^2 \psi(x)$.

$V_A = E(A^2) - E^2(A) =  \lambda^2-  \lambda^2 = 0$.

\

Esta última propiedad se interpreta así: cuando un sistema está en un estado que es valor propio del operador autoadjunto que se está midiendo, con toda seguridad y sin lugar a dudas, la medición del operador dará el valor propio correspondiente. Esta propiedad es la que permite hacer ciencia, es decir, relacionar matemáticas con el resultado de un experimento. 

\

Le queda a uno la impresión de que la teoría con respecto a la medición de una observable está prácticamente entendida. 

 \
 
 ¿Y qué pasa cuando uno desea hacer mediciones de dos observables, una después de la otra?
 
 \
 
 Esta inocente pregunta  ha llevado a la ciencia a develar una de las más sencillas pero contundentes diferencias entre la mecánica clásica y la mecánica cuántica, la cual se expresa por  el \index{principio de incertidumbre }\textbf{principio de incertidumbre}, el cual tiene la siguiente motivación:

Se enciende un bombillo de una luz extremadamente débil (también se puede con una débil fuente de electrones), del cual salen fotones muy de vez en cuando. Incluso uno puede hacer eso con un laser para garantizar que la longitud de onda es siempre la misma. Se hace pasar dicha luz por un rotico circular (sobre   una tabla opaca): este rotico funciona como un detector de posición de los fotones  incidentes. Pero al pasar, la probabilidad de que el fotón colisione con el rotico es mayor cuando menor sea el radio del rotico. La interferencia entre el rotico y el fotón no le causa nada al rotico pero si puede hacer que el fotón cambie, bien su dirección o bien su energía o ambos. En realidad, lo mas probable es que cambie la dirección, se deflecte, pero sin cambiar su energía, es decir, cambia el momento, el cual es un vector con magnitud, dirección y sentido. 

Uno podría pensar que la deflexión es determinista, es decir siempre la misma. Esa es la intuición al estilo de la física clásica. Pero eso no es así: al poner fotoreceptores detrás del rotico a veces se enciende uno de ellos y a veces otro y de manera azarosa. Es decir, hay una probabilidad de por medio, motivo que incita a modelar a un sistema así por la mecánica cuántica. Y con todo, se requiere demostrar que todo lo que pueda pasar se resume simplemente diciendo que disminuir la incertidumbre de posición aumenta la incertidumbre en el momento y al revés. 

\

El modelamiento en mecánica cuántica de un experimento de este tipo es como sigue: a una medición del espacio se le asigna el operador $X$  definido por

$X\psi(x,t) = x \psi(x,t)$.

Observemos que $x$ es real, pues da una posición. Es por eso que el operador $X$ es autoadjunto:

$<X\psi, \phi> = \int_R (X\psi(x)) \phi(x)^*dx =  \int_R x \psi(x) \phi(x)^*dx = \int_R  \psi(x) (x\phi(x))^*dx = \int_R \psi(x) (X\phi(x))^*dx = <\psi, X \phi> $

Para entender qué hace el operador $X$, consideremos una partícula fuertemente concentrada alrededor de un punto $x_o$.  Debido a que todo se concentra cerca de $x_o$, podemos tomar la variable $x$, de la definición de $X$ como constante e igual a $x_o$:

$X\psi(x,t) = x_0 \psi(x,t)$.

Así que leemos: el operador $X$ da la posición de una partícula puntual. 

Podemos considerar ahora el operador momento $P$ e imaginar partículas que viajan con un momento bien definido y repetir lo mismo que hicimos con el operador $X$.

El siguiente paso es imaginar  que podemos hacer un trabajo simultáneo con $X$ y $P$ de la manera descrita. Pues eso es imposible y ese hecho se considera una de las diferencias fundamentales entre la mecánica clásica y la cuántica. La razón consiste en que el operador $X$ y el operador momento $P$ no conmutan.  Veamos todo al detalle.

\

Se postula una relación entre la incertidumbre de posición $\triangle x$ y la de momento  $\triangle p$ tal que:

$\triangle x \triangle p  = k$

Demostrar que esta relación tiene un lugar natural dentro de la mecánica cuántica fue uno de los méritos de Werner Heisenberg, gloria que le fue reconocida por medio del Premio Nobel de 1932. Este memorable físico nacido en 1901 en Munich  no emigró fuera de Alemania cuando se veía venir el programa bélico nazi, sino que desarrolló para Hitler   el estudio de la energía atómica, causa por la cual todos estamos resentidos con él. Existe sin embargo la posibilidad de que sea verdad lo que él alega: que se quedó para asegurarse de que el programa atómico no fuera a producir nada útil para el estado nazi. En contraparte, sus compañeros que emigraron a América desarrollaron el programa de los Alamos que culminó con las bombas sobre el Japón y la subsecuente rendición del Emperador. Este gesto del emperador Hiro-Hito es inmarcesible: prefirió la vida de su pueblo y de sus soldados antes que morir por honor. Y él mismo admitió ser despojado de su  divinidad para someterse a una constitución.

\

Estos hechos muestran claramente que la física y la ciencia en general no son políticamente neutrales y que \textit{el científico está en la obligación de decidir a quién sirve en un momento de crísis y que su elección se hace con mucho tiempo de antelación  en tiempo de paz}.

 \
 
 Volvamos ahora a la elaboración del principio de incertidumbre. Nuestro propósito es estudiar el efecto y resultado de dos mediciones en tándem hechas sobre un sistema. La estructura matemática que capta nuestra preocupación es el conmutador.
 
 \

Dados dos operadores $A,B$, el \index{conmutador} \textbf{conmutador $[A,B]$} es el operador:

$[A,B] = AB -BA $ 

Si los operadores conmutan, el conmutador es cero. Por eso decimos que el conmutador $[A,B]$ es el operador que mide le deficiencia de conmutación.

Calculemos el conmutador del operador posición $X$ con el operador momento $P$:

$[X,P]\psi(x) = [X, -ih\partial/\partial x]\psi(x)$
 
$ = X (-ih\partial\psi(x)/\partial x ) - (-ih\partial/\partial x ( X \psi(x))) $

$= x (-ih\partial\psi(x)/\partial x ) - (-ih\partial/\partial x ( x \psi(x)))$

$ = -ihx\partial\psi(x)/\partial x  + ih  \psi(x) + ihx\partial\psi(x)/\partial x $ 
 
$ =  ih \psi(x) $

Como esto es cierto para toda $\psi(x)$, obtenemos que el operador momento y el operador posición no conmutan:

$[X,P] = ih$

\

El sentido físico del conmutador nos lo fue explicado por el \index{principio de incertidumbre} \textbf{principio de incertidumbre} propuesto por Heinsenberg, cuya expresión general entramos a elaborar.

Vamos a elaborar ahora una versión moderna del principio de incertidumbre, válida para cualquier par de observables, es decir, de operadores autoadjuntos, por lo cual los valores propios de cada operador en separado son reales y factibles de ser medidos en un experimento.  El problema es que cuando dos operadores van en secuencia uno no necesariamente obtiene un operador autoadjunto. Eso se deduce de la siguiente identidad:

Supongamos que tenemos dos operadores autoadjuntos $A,B$. Entonces se cumple que

$AB = \frac{1}{2}[A,B] + \frac{1}{2}\{A,B\}$

donde $ [A,B] = AB-BA$ es el conmutador, mientras que $\{A,B\} = AB+BA$ es el anticonmutador. Lo dicho es verdad puesto que

$\frac{1}{2}[A,B] + \frac{1}{2}\{A,B\} = \frac{1}{2}(AB-BA) + \frac{1}{2}(AB+AB) = AB$

Tenemos ahora que  el conmutador $[A,B]$ de dos operadores autoadjuntos $A,B$ es anti-autoadjunto (por lo tanto, su esperanza  es puramente imaginaria):

$ ([A,B])^\dagger  = (AB-BA)^\dagger = (AB)^\dagger-(BA)^\dagger = B^\dagger A^\dagger -A^\dagger B^\dagger $

$= BA  - AB = - (AB - BA) $

$= - [B,A]$

 Por otro lado, el anticonmutador si es autoadjunto ( su esperanza es real):

$ (\{A,B\})^\dagger  = (AB+BA)^\dagger = (AB)^\dagger+(BA)^\dagger = B^\dagger A^\dagger +A^\dagger B^\dagger $

$= BA  + AB =  AB + BA $

$=  \{B,A\}$

Vemos claramente ahora que la compuesta de dos operadores autoadjuntos es autoadjunto sólo cuando los operadores conmutan. Por esa razón, el estudio de un experimento en el cual se mide primero el observable $B$ y después el $A$ se modela por el operador $AB$ que no es necesariamente autoadjunto y cuyo comportamiento hay que estimar por medios aproximados. Con todo, es conveniente tener en claro que nuestros experimentos miden cantidades que quedan bien modeladas por números reales.  Por tanto, a un  operador cualquiera $T$ hay que meterlo dentro de una estructura que produzca números reales. Sin excepciones, se le mete dentro de un producto interior de la forma

$<T\psi, \phi>$

y en particular a nuestro $AB$ se le involucra dentro de la esperanza:

$<AB\psi, \psi> =<B\psi, A\psi> $

Expresión válida para operadores autoadjuntos.

\

Lo primero que necesitamos para nuestro estudio de $AB$ es  la \index{desigualdad de Cauchi-Schwarz} \textbf{desigualdad de Schwarz}, válida en general para cualquier espacio vectorial con escalares complejos  y   con   producto interno complejo:

$||a||^2 ||b||^2 \geq |a \cdot b|^2 = |<a , b>|^2$

 cuya demostración es como sigue:
 
 Para cualquier par de vectores $a$, $b$ y para cualquier número complejo $\lambda$ tenemos:
 
 $||(a -\lambda b)||^2 \geq 0$
 
 $||(a -\lambda b)||^2= <a -\lambda b,a -\lambda b> $
 
 $= <a,a>  - <a, \lambda b> - <  \lambda b,a > + < \lambda b ,   \lambda b> $
 
 $ =<a,a> -  \lambda^* <a,   b> - \lambda <  b,a > + \lambda \lambda^*<  b ,    b> $ 
 
 $=<a,a> -  \lambda^* <a,   b> - \lambda <  a,b >^* + \lambda \lambda^*<  b ,    b> $

Ahora ponemos

$\lambda = \frac{ <a,b>}{<b,b>}  $

y nos queda:

$||(a -\lambda b)||^2 = <a,a> -  (\frac{ <a,b>}{<b,b>})^* <a,   b> $
$- (\frac{ <a,b>}{<b,b>}) <  a,b >^* $

$+ (\frac{ <a,b>}{<b,b>})(\frac{ <a,b>}{<b,b>})^*<  b ,    b> $

\
Como $<b,b> = ||b||^2$ es real e igual a su conjugado, podemos simplificar:

$||(a -\lambda b)||^2 = <a,a> -  \frac{<a,b>^*<a,   b>}{<b,b>}     $
$- \frac{<  a,b >^* <a,b> } { <b,b> } $
$+ \frac{(<a,b> <a,b>^*}{ <b,b>}  $

\

$||(a -\lambda b)||^2 = <a,a> -  \frac{<a,b>^*<a,   b> }{ <b,b>} \geq 0  $

\

de lo cual obtenemos:

$<a,a><b,b>  \ \ \geq  \ \  <a,b>^*<a,   b> $

lo cual es equivalente a la desigualdad buscada:

$||a||^2 ||b||^2 \geq |<a , b>|^2$

\

Demostremos ahora la siguiente desigualdad:

$ ||A\psi||^2 ||B\psi||^2 \ \  \geq \  \ \frac{1}{4} | < (AB-BA) \psi ,  \psi>|^2 = \frac{1}{4}|<[A,B]\psi, \psi>|^2 $ donde $[A,B]= AB-BA$.

\

Veamos. Sean $A,B$ operadores autoadjuntos. Sobre la desigualdad de Cauchi-Schwarz:

$||a||^2 ||b||^2 \ \ \geq \ \ |<a , b>|^2$

hacemos las sustituciones:

$a \rightarrow   B\psi$

$b \rightarrow   A \psi$

y nos queda:

$||B\psi||^2 ||A\psi||^2  = ||A\psi||^2 ||B\psi||^2 \ \ \geq \ \ |< B \psi ,  A \psi>|^2$

Ahora bien:

$|< B \psi ,  A \psi>|^2 \geq  |Im (< B \psi ,  A \psi>|^2)|$

$ = |2\frac{1}{2}Im (< B \psi ,  A \psi>|^2)|$
$ =\frac{1}{4} |2Im (< B \psi ,  A \psi>|^2)|$

\

Tengamos en cuenta que si $z = a+bi$ entonces $z-z^* = 2bi = 2Im z$, por lo que

$ \frac{1}{4} |2Im (< B \psi ,  A \psi>|^2)$
$ =\frac{1}{4} | < B \psi ,  A \psi>- < B \psi ,  A \psi>^*|^2 $

$=\frac{1}{4} | < B \psi ,  A\psi>- < A \psi ,  B \psi>|^2  $

$=\frac{1}{4} | < AB \psi ,  \psi>- < BA \psi ,   \psi>|^2  $

$=\frac{1}{4} | < (AB-BA) \psi ,  \psi>|^2  $

$ = \frac{1}{4}|<[A,B]\psi, \psi>|^2 $
\

Ahora unimos la cadena de desigualdades para obtener

$||A\psi||^2 ||B\psi||^2  \ \  \geq \ \  \frac{1}{4} \ |<[A,B]\psi, \psi>|^2 $

tal como lo habíamos deseado.

\

Esa desigualdad es válida siempre que los operadores $A$ y $B$ sean autoadjuntos, digamos: $<A\psi, \phi > = <\psi, A\phi>$. Por consiguiente, podemos usar la misma desigualdad para   otros operadores autoadjuntos, cualesquiera que sean. Lo hacemos siguiendo la inteligente notación de Wikipedia: 

$A \rightarrow \triangle_\psi A = A-E_\psi(A)$, autoadjunto.

$B \rightarrow \triangle_\psi B = B-E_\psi(B)$ autoadjunto.

donde hemos enfatizado que la esperanza de un operador se define con respecto a un estado dado, a una función de onda $\psi$

La desigualdad
 
 $||A\psi||^2 ||B\psi||^2  \  \geq  \ \frac{1}{4}|<[A,B]\psi, \psi>|^2 $

se lee ahora como

$||\triangle_\psi A \psi||^2 ||\triangle_\psi B \psi||^2 \ \geq \ \frac{1}{4}|<[\triangle_\psi A ,\triangle_\psi B]\psi, \psi>|^2   $

\

Como $\triangle_\psi A$ es autoadjunto, obtenemos:

$||\triangle_\psi A||^2 = <\triangle_\psi A, \triangle_\psi A> $ 

$ = <(\triangle_\psi A)^2, I>  = \int (\triangle_\psi A)^2 \psi(x) \psi^*(x) dx  = E_\psi[( \triangle_\psi A)^2] $ 

similarmente

$||\triangle_\psi B||^2 = E_\psi[( \triangle_\psi B)^2] $ 

\

Tenemos entonces que

$E_\psi[( \triangle_\psi A)^2] E_\psi[( \triangle_\psi B)^2] \ \geq \ \frac{1}{4}|<[\triangle_\psi A ,\triangle_\psi B]\psi, \psi>|^2   $
 
 \
 
Por otro lado:

$[\triangle_\psi A,\triangle_\psi B] = [A - E_\psi(A),B - E_\psi(B)] $

$= (AB - AE_\psi(B)  - E_\psi(A)B + E_\psi(A)E_\psi(B)) $

$ - (BA - BE_\psi(A)  - E_\psi(B)A + E_\psi(B)E_\psi(A)) = AB - BA $

$ = [A,B]$

Reemplazando $[\triangle_\psi A,\triangle_\psi B]$ por $[A,B]$ en

$E_\psi[( \triangle_\psi A)^2] E_\psi[( \triangle_\psi B)^2] \ \geq \ \frac{1}{4}|<[\triangle_\psi A ,\triangle_\psi B]\psi, \psi>|^2   $

  nos queda:

$E_\psi[( \triangle_\psi A)^2] E_\psi[( \triangle_\psi B)^2] \ \geq \ \frac{1}{4}|<[  A , B]\psi, \psi>|^2   $

\

Podemos decir que con respecto a cualquier función de onda $\psi$:

\

$V_A V_B \geq \frac{1}{4}|E([  A,  B])|^2$

\

Podemos verbalizar esta desigualdad enunciando el 
\textbf{Principio de Incertidumbre} como sigue:

 Consideremos un experimento en el cual uno prepara un sistema en una función de onda dada y el objetivo es medir dos observables, $A,B$, una después de la otra. Una cierta cantidad de veces se mide: primero una observable   y a continuación la otra. Eso está representado por la expresión 
 
 $<BA\psi, \psi> = <A\psi, B\psi>$
 
 que se encuentra al principio de la demostración.
 
 De las mediciones se cuantifica la varianza de cada una de las observables, $V_A, V_B$. Se tiene:

 $V_A V_B \geq \frac{1}{4}|E([  A,  B])|^2$

que se lee: el producto de las varianzas es supra-proporcional a la esperanza del conmutador. Esto implica que es imposible disminuir al mismo tiempo las varianzas de las mediciones de las dos observables  tanto como se desee, pues medir uno de los operadores con precisión implicará que las mediciones del otro sean muy dispersas, de tal forma que el producto de las varianzas se mueve, difusamente, encima de una hipérbole.

\

\section{INVARIANCIA GAUGE DEL HAMILTONIANO}

Resulta que un átomo no sólo tiene niveles energéticos, asociados a los valores propios del Hamiltoniano. También tiene otras observables, digamos el momento angular. Dichas observables también evolucionan y el culpable de todo es el Hamiltoniano que es el operador que causa la evolución de todo el sistema.

\

Por todo esto, debemos preguntarnos la compatibilidad de los formalismos con nuestro principio de invariancia gauge:

\
1. Qué efecto cuántico observable tiene la arbitrariedad del nivel cero en la energía potencial?

2. Está la física asociada a los valores propios libre de toda arbitrariedad matemática?

3. Al evolucionar un sistema cuántico, las demás observables son inmunes a las arbitrariedades matemáticas admitidas por el Hamiltoniano?

\

Vimos anteriormente que las arbitrariedades asociadas al espacio proyectivo se solucionan de una vez y para siempre alegando que toda ecuación de la mecánica cuántica es lineal y que por tanto los escalares entran, salen y desaparecen sin ningún problema. 

Ahora revisemos los problemas relacionados al cambio de fase ocasionado por la arbitrariedad nacida de  definir la energía potencial módulo una constante. Primero para el Hamiltoniano y después para las demás observables. En resumen, lo que tenemos que estudiar es en qué se diferencia la física de dos evoluciones distintas, una debida a $H$ y otra debida a $H+k$:

En el primer caso,  $\psi(x,t)= e\sp{-itH/\hbar} \psi\sb{0}$

En el segundo,  $\psi(x,t)= e\sp{-it(H+k)/\hbar}  \psi\sb{0}=$
$ e\sp{-itk/\hbar}  e\sp{-itH/\hbar}  \psi\sb{0}$

En la última igualdad utilizamos la conmutatividad de los operadores exponenciales cuando los exponentes conmutan. Todo lo que obtuvimos fue un cambio de fase. Como ya sabemos, eso no cambia las probabilidades, pero queda el problema de los valores propios del Hamiltoniano. Dichos valores están conectados a la física a través de la siguiente regla:

La física no depende de los valores absolutos de la energía sino de las diferencias de energía. Esto es cierto en mecánica clásica y tiene que ser cierto en mecánica cuántica pues la primera es un caso asintótico (en el espacio de caminos) de la segunda. De hecho, lo que se mide en espectroscopía es la energía que arrastra un fotón y esa depende de los saltos de energía entre los diferentes niveles de los átomos y estos están relacionados con diferencias entre valores propios del Hamiltoniano.

\

Cambian las diferencias entre valores propios debidos a un cambio de fase?

Probemos que un cambio global de fase no cambia las diferencias entre valores propios:

Si $H\psi = \lambda \psi$ entonces, multiplicando por el escalar $e\sp{-itk/\hbar}$ a ambos lados tenemos:
$e\sp{-itk/\hbar} H\psi = e\sp{-itk/\hbar} \lambda \psi $.  Conmutando:

$ H e\sp{-itk/\hbar}\psi =  \lambda e\sp{-itk/\hbar}\psi $

Lo cual quiere decir que si $\lambda $ es un valor propio del Hamiltoniano original, entonces un corrimiento da fase en el Hamiltoniano para todos los efectos equivale a un corrimiento de fase de las funciones de onda. El corrimiento es global y biyectivo. Por lo tanto, la diferencia entre valores propios se conserva. Concluimos que la física del Hamiltoniano es inmune a las arbitrariedades matemáticas de nuestros formalismos.

Y qué de las demás observables? Todas las observables están ligadas al Hamiltoniano, pues es éste operador el causante de la evolución del sistema. Por lo tanto, cualquier culpa no puede ser de este operador. Veamos esto al detalle.

\section{ INVARIANCIA GAUGE DE LAS OBSERVABLES}

Para poder hacer mecánica cuántica se necesita un espacio de Hilbert y en ese caso a los operadores autoadjuntos se les puede atribuir propiedades físicas, las cuales son libres de todo gravamen dictaminado por las arbitrariedades de origen matemático. Veamos:

\

¿Qué pasa cuando sometemos dicha expresión a evolución y permitimos una arbitrariedad en el Hamiltoniano $H \rightarrow H +k $?

\

La forma cómo evoluciona el sistema bajo $H$ es:

$\phi(x,t)= e\sp{-itH/\hbar} \phi\sb{0}$

Por tanto, las probabilidades en el tiempo $t$ tendrán el valor:

$p\sb{i,t}= |<\psi\sb{i},e\sp{-itH/\hbar} \phi > \sp{2}|$

Si el sistema evoluciona bajo $H+k$ entonces:

$\psi(x,t)= e\sp{-it(H+k)/\hbar}  \psi\sb{0}$
$= e\sp{-itk/\hbar}  e\sp{-itH/\hbar}  \psi\sb{0}$

Y las probabilidades correspondientes son:

$$p\sb{i,k,t}= |<\psi\sb{i},e\sp{-itk/\hbar} e\sp{-itH/\hbar} \phi >| \sp{2}
\addtocounter{ecu}{1}   \hspace{4cm} (\theecu )      $$

El escalar $e\sp{-itk/\hbar}$ puede salir del producto interno y como tanto él como su conjugado tienen norma 1, al tomar el cuadrado nada cambia. Como consecuencia, uno puede tomar la constante que se le antoje pero la física permanece igual:

$p\sb{i,k,t}= p\sb{i,l,t}$

\

Hemos demostrado, parcialmente, que la física de las observables es inmune a las arbitrariedades matemáticas nacidas de nuestros formalismos. Nuestro objetivo es demostrar que cumplimos con el principio gauge, pero hasta ahora sólo hemos podido librar al Hamiltoniano de toda acusación. Nos quedan las arbitrariedades propias de cada observable. Esto lo haremos muy explícitamente para el campo electromagnético. A propósito, cuál será el grupo de invariancia? Pues como se trata de multiplicar por exponenciales imaginarias de números reales, estamos trabajando con \textbf{U(1)}.

En nuestros cálculos hemos estudiado la arbitrariedad de definir la energía potencial aumentando o disminuyendo una constante cualquiera. La constante la hemos tomado igual para todo el espacio y para todos los tiempos. Será eso demasiada restricción? Cómo hacemos para permitir que investigadores diferentes tomen constantes diferentes en diferentes lugares o momentos? Cuál podrá ser la influencia de tales libertades?

\section{INVARIANCIA GLOBAL vs INVARIANCIA LOCAL }

Todo lo que hemos estudiado parece resumirse muy simplemente: cualquier arbitrariedad que origine un cambio global de fase no afecta la física. O bien, la mecánica cuántica de cualquier sistema cuya evolución está definida por   un operador autoadjunto es una teoría gauge con grupo de invariancia global \textbf{U(1)}.

\

Me place narrar ahora uno de los experimentos más sublimes de toda la física, en la cual veremos que cambiar de manera no uniforme las fases de un sistema físico en un proceso dinámico si tiene resultados observables.

En mecánica cuántica cada posible canal de evolución, cada camino, en general cuenta poco. Pues porque son muchos y cada uno tiene vecinos que valen aproximadamente lo mismo. Lo que importa es la manera como dichos caminos trabajan en conjunto, sea en forma coherente, sumando insignificancias o en interferencia destructiva, aniquilándose mutuamente.

Hemos visto que un cambio global de fase no produce ningún cambio en los resultados observables debido a los patrones de interferencia. Ahora podemos advertir que un cambio local de fase, en este lugar de una manera y allá de otra, si puede causar efectos observables.

Un cambio global, homogéneo, de fase era originado por el cambio $H \rightarrow H+k$. Podemos lograr un cambio no homogéneo cambiando $k$ a medida que cambia el lugar. Una grandiosa coincidencia de ingenio y de casualidades de órdenes de magnitud permitió en 1975 a R. Collela, A. Overhauser, y S.A. Werber correr un  experimento basado en la siguiente idea:

Consideremos el experimento, que se puede hacer en una mesa, de hacer pasar neutrones por los bordes de un rectángulo, comenzando en una esquina y terminando en el vértice opuesto por la diagonal. El montaje es una variante del experimento de interferencia de luz que pasa por entre dos rendijas. En este caso, los neutrones pueden pasar sea por 'abajo' o sea por 'arriba'. Resulta que el potencial gravitatorio depende de la altura y hace las veces de constante $k$ en nuestro razonamiento anterior.

Para variar no homogéneamente $k$, el potencial gravitatorio, lo único que se requiere es girar el rectángulo de tal forma que el lado de abajo y el de arriba queden a diferentes alturas.

La parte de neutrón que viaje por arriba experimentará, en promedio, un potencial mayor, una $k$ mayor, que la parte del neutrón que viaje por abajo. Por lo tanto, la parte de arriba oscilará más velozmente que  la de abajo. Al juntarse las partes en la otra esquina interferirán ora constructivamente ora destructivamente. El grado y tipo de interferencia depende de la diferencia de alturas, la cual se puede graduar cambiando la pendiente del rectángulo.

Los investigadores antes mencionados estudiaron la variación de la interferencia como función del ángulo de giro, o sea la pendiente del rectángulo. La coincidencia entre los resultados experimentales y los teóricos es asombrosa. Cálculos muy semejantes haremos un poco más tarde en relación con el campo magnético.

Es por eso, que distinguimos fuertemente entre un cambio global de fase y otro local, punto por punto. Mientras que el primero es imperceptible, el segundo puede crear resultados asombrosos, por no decir mágicos.

\chapter{ELECTROMAGNETISMO COMO TEORIA GAUGE}

Hemos visto que el electromagnetismo clásico, la teoría que estudia la interacción entre el campo electromagnético y la materia, admite una arbitrariedad: el potencial vector sólo puede definirse módulo un gradiente de una función escalar. Por construcción se cumple el principio gauge pues cambiar de  gradiente no tiene ninguna consecuencia física. El objetivo de la presente sección es formular el Hamiltoniano correspondiente, probar que sigue admitiendo la misma arbitrariedad del campo clásico, y verificar que  cumple  el criterio gauge. Nuestro punto de partida es el Lagrangiano para pasar por una transformación al Hamiltoniano.

\section{EL LAGRANGIANO ELECTROMAGNETICO}

¿Qué quiere decir que el Lagrangiano electromagnético sea tal o cual? Eso quiere decir que al ponerlo en el funcional de la acción, y al estudiar las variaciones de dicho funcional debidas a perturbaciones infinitesimales, y al suponer que la variación a primer orden en el funcional sea cero, debemos obtener las ecuaciones de Euler Lagrange, las cuales deben describir, como mínimo, lo que hace una partícula cargada   ante la acción del campo.

En una investigación con el Lagrangiano, podríamos partir de un Lagrangiano para deducir ecuaciones de campo junto con su interacción con la materia. Pero a uno siempre le queda la espina de saber cómo se inventan los Lagrangianos.  Debido a que ambas posiciones son intrigantes, las ilustraremos a ambas.    Comencemos a partir de las ecuaciones del campo electromagnético y lleguemos al Lagrangiano oficial.

\subsection{De la Fuerza de Lorentz al Lagrangiano}

Para reinventar el Lagrangiano electromagnético, nuestra estrategia será llegar a las ecuaciones de Euler Lagrange de las cuales uno puede leer el Lagrangiano, pues dichas ecuaciones son

\bigskip

$\frac{\partial L}{\partial x_i} - \frac{d}{dt}(\frac{\partial L}{\partial \dot{x_i}})=0$

\bigskip

o bien

\bigskip

$\frac{\partial L}{\partial x_i} = \frac{d}{dt}(\frac{\partial L}{\partial \dot{x_i}}) $

\

 Empecemos. El campo eléctrico $\vec E$ y el potencial eléctrico $\phi$ están relacionados por

$$\vec E = -\nabla \phi
\addtocounter{ecu}{1}   \hspace{4cm} (\theecu )      $$

Similarmente, el campo magnético $\vec B $ y el potencial vector $\vec A$ se relacionan por

$$\vec B = \nabla \times \vec A
\addtocounter{ecu}{1}   \hspace{4cm} (\theecu )      $$

La fuerza de Lorentz es

$$\vec F = e(\vec E + (1/c) \vec v \times \vec B)
     $$

$$\vec F = e(\vec E + (1/c) \vec v \times (\nabla \times \vec A))
     $$

Si $\vec A = (A_1,A_2,A_3)$, entonces

\bigskip

$\nabla \times \vec A = \pmatrix{ \frac{\partial A_3}{\partial y}-\frac{\partial A_2}{\partial z} \cr
\frac{\partial A_1}{\partial z}-\frac{\partial A_3}{\partial x} \cr
\frac{\partial A_2}{\partial x}-\frac{\partial A_1}{\partial y}}$

\bigskip

Y por tanto

\bigskip

$$\vec v \times \vec B =\pmatrix{\dot{y}[\frac{\partial A_2}{\partial x}-\frac{\partial A_1}{\partial y}] -\dot{z} [\frac{\partial A_1}{\partial z}-\frac{\partial A_3}{\partial x}]  \cr
\dot{z}[\frac{\partial A_3}{\partial y}-\frac{\partial A_2}{\partial z}]-\dot{x}[\frac{\partial A_2}{\partial x}-\frac{\partial A_1}{\partial y}]   \cr
 \dot{x}[\frac{\partial A_1}{\partial z}-\frac{\partial A_3}{\partial x}]-\dot{y}[\frac{\partial A_3}{\partial y}-\frac{\partial A_2}{\partial z}]    } $$

     \bigskip

Eso implica que la fuerza de Lorentz

     $$m \pmatrix{\ddot{x} \cr \ddot{y} \cr \ddot{x}} =  \vec F = e(-\nabla \phi + (1/c) \vec v \times (\nabla \times \vec A))
     $$

     es lo mismo que

$$m \pmatrix{\ddot{x} \cr \ddot{y} \cr \ddot{x}} =   e\{-\pmatrix{\frac{\partial\phi}{\partial x} \cr
\frac{\partial\phi}{\partial y} \cr
\frac{\partial\phi}{\partial z}  } + (1/c) \pmatrix{\dot{y}[\frac{\partial A_2}{\partial x}-\frac{\partial A_1}{\partial y}] -\dot{z} [\frac{\partial A_1}{\partial z}-\frac{\partial A_3}{\partial x}]  \cr
\dot{z}[\frac{\partial A_3}{\partial y}-\frac{\partial A_2}{\partial z}]-\dot{x}[\frac{\partial A_2}{\partial x}-\frac{\partial A_1}{\partial y}]   \cr
 \dot{x}[\frac{\partial A_1}{\partial z}-\frac{\partial A_3}{\partial x}]-\dot{y}[\frac{\partial A_3}{\partial y}-\frac{\partial A_2}{\partial z}]    } \}
     $$

Miremos qué tenemos en la primera coordenada y reescribámosla para hallar la primera ecuación de Euler Lagrange:

$m   \ddot{x} = e \{ -\frac{\partial\phi}{\partial x} + (1/c)[\dot{y}[\frac{\partial A_2}{\partial x}-\frac{\partial A_1}{\partial y}] -\dot{z} [\frac{\partial A_1}{\partial z}-\frac{\partial A_3}{\partial x}]]        \}$

     $m   \ddot{x} = -e   \frac{\partial\phi}{\partial x} + (e/c)[\dot{y}[\frac{\partial A_2}{\partial x}-\frac{\partial A_1}{\partial y}] -\dot{z} [\frac{\partial A_1}{\partial z}-\frac{\partial A_3}{\partial x}]]        $

  $m   \ddot{x} = -e   \frac{\partial\phi}{\partial x} + (e/c)[\dot{y} \frac{\partial A_2}{\partial x}-\dot{y} \frac{\partial A_1}{\partial y} - \dot{z} \frac{\partial A_1}{\partial z} +\dot{z}\frac{\partial A_3}{\partial x}]        $

$m   \ddot{x} = -e   \frac{\partial\phi}{\partial x} + (e/c)[\dot{y}\frac{\partial A_2}{\partial x}+\dot{z}\frac{\partial A_3}{\partial x} -(\frac{\partial A_1}{\partial y}\dot{y} +\frac{\partial A_1}{\partial z}\dot{z}) ]        $

\bigskip

Nos damos cuenta que hay una simetría frustrada que se puede remediar como sigue

$m   \ddot{x} = -e   \frac{\partial\phi}{\partial x} + (e/c)[\dot{x}\frac{\partial A_1}{\partial x}+\dot{y}\frac{\partial A_2}{\partial x}+\dot{z}\frac{\partial A_3}{\partial x} -(\dot{x}\frac{\partial A_1}{\partial x}+\frac{\partial A_1}{\partial y}\dot{y} +\frac{\partial A_1}{\partial z}\dot{z}) ]        $

\bigskip

Utilizando la regla de la cadena obtenemos:

$m   \ddot{x} = -e   \frac{\partial\phi}{\partial x} + (e/c)[\dot{x}\frac{\partial A_1}{\partial x}+\dot{y}\frac{\partial A_2}{\partial x}+\dot{z}\frac{\partial A_3}{\partial x} -\frac{dA_1}{dt} ]        $

$m   \ddot{x} +(e/c) \frac{dA_1}{dt}= -e   \frac{\partial\phi}{\partial x} + (e/c)[\dot{x}\frac{\partial A_1}{\partial x}+\dot{y}\frac{\partial A_2}{\partial x}+\dot{z}\frac{\partial A_3}{\partial x}  ]        $

$ \frac{d}{dt}( m\dot{x} + (e/c)A_1)=    \frac{\partial }{\partial x}(-e\phi + (e/c)[    \dot{x} A_1 + \dot{y}A_2 + \dot{z}A_3] ) $

\bigskip

Esta ecuación para la primera coordenada tiene una generalización inmediata para todas las demás:

$ \frac{d}{dt}( m\dot{x_i} + (e/c)A_i)=    \frac{\partial }{\partial x_i}(-e\phi + (e/c)[    \dot{x} A_1 + \dot{y}A_2 + \dot{z}A_3] ) $

Como las ecuaciones de Euler Lagrange son

$\frac{\partial L}{\partial x_i} = \frac{d}{dt}(\frac{\partial L}{\partial \dot{x_i}}) $

o lo que es lo mismo:

$ \frac{d}{dt}(\frac{\partial L}{\partial \dot{x_i}}) = \frac{\partial L}{\partial x_i}$

tenemos que encontrar $L$ tal que sea simétrico en todas las coordenadas y que cumpla:

$\frac{\partial L}{\partial \dot{x_i}} = m\dot{x_i} + (e/c)A_i$

$L= -e\phi + (e/c)[    \dot{x} A_1 + \dot{y}A_2 + \dot{z}A_3] + \phi(\dot{x_i})$

\

De la primera ecuación tenemos:

$L_o = (1/2)m\dot{x}^2  +  .. .. +  (e/c)  (\dot{x} A_1+ .. ..) + \phi(x,y,z)$

Comparando con la segunda ecuación, llegamos a:

$$L_o= (1/2)m (\dot{x}^2 + \dot{y}^2+\dot{x}^2) - e \phi + (e/c)  (\dot{x} A_1+ \dot{y} A_2+\dot{z} A_3) $$

$$L_o= (1/2)m (\dot{x}^2 + \dot{y}^2+\dot{x}^2) - e \phi + (e/c)  (\dot{x}, \dot{y},\dot{z})   \cdot (A_1,A_2,A_3)$$

$$L_o= (1/2)m \|\vec v \|\sp{2} - e \phi + (e/c) \vec v \cdot \vec A
        $$

el cual es \index{Lagrangiano!electromagnético} \textbf{Lagrangiano electromagnético} oficial y al cual nos referiremos por defecto.

\

Recalquemos que la acción correspondiente es una integral de camino determinada por la trayectoria de la partícula en estudio, tal que las condiciones iniciales y finales permanecen constantes:

\

$A = \int [(1/2)m \|\vec v \|\sp{2} - e \phi + (e/c) \vec v \cdot \vec A] dt $
 
\
  
Hemos deducido esta acción partiendo de la fuerza de Lorentz para tener al menos una idea de cómo fabricar Lagrangianos. Destaquemos que nos sirvió de guía el haber tratado de conservar la simetría entre las variables. La abstracción de una simetría se hace por medio de un grupo: si queremos simetría entre variables, el grupo es el grupo de permutaciones entre ellas y lo que se requiere es que el Lagrangiano no varíe ante la acción del grupo. En general, lo que se hace es  fijar un grupo y construir un Lagrangiano que sea invariante ante su acción. 

Pero si uno ya tiene un Lagrangiano, entonces uno hace  el estudio variacional para deducir las leyes que rigen el sistema. Ilustraremos eso en  la siguiente subsección.

\bigskip

\subsection{La acción de los Cosserat}

Proponer una acción o un  Lagrangiano corresponde a postular la conducta fundamental de un sistema. Por ello, proponer Lagrangianos es un arte delicado y que lo puede hacer a uno famoso, si es que uno da con algo sencillo, hermoso y poderoso. Además del Lagrangiano oficial para el campo electromagnético  que acabamos de re-inventar, existen otras propuestas. La que vamos a ver enseguida es debida a los dos hermanos Cosserat hacia 1909, quienes eran ingenieros civiles franceses con un acendrado interés en lo práctico y por ende en lo teórico de muy alto nivel (Pommaret, 1987). 

\

Los Cosserat parten de la acción siguiente:

\

$A = \int (\frac{\epsilon_o}{2} \vec E \cdot \vec E - \frac{1}{2\mu_o} \vec B \cdot \vec B + \rho \vec v \cdot \vec A - \rho V + m \frac{\vec v \cdot \vec v }{2})dxdydzdt $

\

donde $A$ es la acción, $\vec A$ es el potencial vector, $\vec E$ es el campo eléctrico, $\vec B$ es el campo magnético,  $V$ es el potencial eléctrico, $\vec v $ es la velocidad de una partícula cargada que puede ser el electrón, cuya masa es $m$ y su carga es $\rho$ y que está en la posición $\vec r = (x,y,z)$. Hay otras constantes por ahí cuyo  papel principal es hacer que todo quede en las mismas unidades de energía. 

La idea de esta formulación es representar el campo electromagnético por los dos primeros términos, $\frac{\epsilon_o}{2} \vec E \cdot \vec E - \frac{1}{2\mu_o} \vec B \cdot \vec B $, la interacción entre la partícula cargada y el campo por  $\rho \vec v \cdot \vec A - \rho V$, y el término de la energía cinética de la partícula por $ \rho \frac{\vec v \cdot \vec v }{2}$. 

Obsérvese la forma simétrica de la acción con respecto al  trato del espacio y del tiempo (para 1909 ya se conocía la relatividad especial). 

No se especifican las condiciones de frontera pues, como lo hemos aprendido con el tratamiento de la acción para una partícula libre, las condiciones de frontera se aniquilan haciendo que las pequeñas variaciones que se tienen en cuenta no las cambien. De esa forma al integrar por partes, ellas no cuentan para nada. Claro que la situación no es tan simple  cuando se hace esto en un espacio de 4 dimensiones: equivale a decir que el universo es finito y sin frontera (quizá insertado en un espacio de dimensión 5 y cerrándose sobre sí mismo como un círculo) o que es infinito y que todo decae rápidamente a cero. Ambas salidas son traumáticas. Sin embargo, fue demostrado por Kaluza-Klein que  es buena idea pensar en un universo de 5 dimensiones.

\

La idea es ahora calcular el efecto a primer orden sobre la acción de una pequeña perturbación en las variables fundamentales, igualar dicho efecto a cero y ver qué pasa. Pero, ¿cuáles son las variables fundamentales? Como sabemos, el campo magnético no es fundamental sino que viene del potencial vector, mediante la identidad:

$\vec B = \nabla \times  \vec A$

De igual forma:

$\vec E = - \nabla V - \frac{\partial \vec  A}{\partial t} $

Reemplazando estas identidades en la acción

$A = \int (\frac{\epsilon_o}{2} \vec E \cdot \vec E - \frac{1}{2\mu_o} \vec B \cdot \vec B + \rho \vec v \cdot \vec A - \rho V + m \frac{\vec v \cdot \vec v }{2})dxdydzdt $

nos queda:

\

$A = \int (\frac{\epsilon_o}{2} (- \nabla V - \frac{\partial \vec A}{\partial t})^2 - \frac{1}{2\mu_o}  (\nabla \times  \vec A)^2 + \rho \vec v \cdot \vec A - \rho V + m \frac{(\vec v)^2  }{2})dxdydzdt $

\

donde hemos utilizado una notación taquigráfica para los cuadrados.

Vemos que las variables fundamentales son $V, \vec A, \vec v, \vec r $. Pero esto se debe a que nosotros estamos modelando una partícula, porque si estuviésemos modelando un líquido o un gas, como para el estudio de los efectos electromagnéticos sobre  el citoplasma celular o sobre la evolución de la galaxia (cada estrella se toma como una molécula) o sobre las mareas solares o sobre  un motor de plasma, entonces también deberíamos permitir que la carga y que la masa pudiesen variar.

Ahora, nuestro propósito es calcular 

\

$\delta A = A(V + \delta V, \vec A +\delta \vec A, \vec v + \delta \vec v, \vec r ) - A(V, \vec A, \vec v, \vec r) $.

\

para después tomar los términos a primer orden, igualar a cero y ver que pasa. 

\

Nuestra  acción

$A = \int (\frac{\epsilon_o}{2} (- \nabla V - \frac{\partial \vec A}{\partial t})^2 - \frac{1}{2\mu_o}  (\nabla \times  \vec A)^2 + \rho \vec v \cdot \vec A - \rho V + m \frac{(\vec v)^2  }{2})dxdydzdt $

puede partirse en términos asociados al campo y otros asociados a la partícula o a su interacción con el campo:

$A = \int (\frac{\epsilon_o}{2} (- \nabla V - \frac{\partial \vec A}{\partial t})^2 - \frac{1}{2\mu_o}  (\nabla \times  \vec A)^2 + L_o)dxdydzdt $

donde

$L_o= (1/2)m \|\vec v \|\sp{2} - \rho V + \rho \vec v \cdot \vec A $

Tenemos la fuerte tentación de decir que, por los resultados del ejercicio anterior, de este Lagrangiano $L_o$ sale la fuerza de Lorentz. Seguramente sí, pero hay un problema que es como un dolor de muela: la acción del ejercicio anterior era una integral de línea mientras que la acción que estamos considerando es una integral sobre todo el espacio-tiempo de 4 dimensiones. 

Nos conviene dividir el estudio de todo ese complejo  problema en partes adecuadas, así sea por razones de tipografía. Nuestro primer trabajo   es estudiar el problema variacional sobre  la acción asociada al campo en sí:

\

$A_c = \int (\frac{\epsilon_o}{2} (- \nabla V - \frac{\partial \vec A}{\partial t})^2 - \frac{1}{2\mu_o}  (\nabla \times  \vec A)^2  )dxdydzdt $

\

en la cual los campos varían levemente como entidades que tienen existencia en sí mismos.

Para calcular

$\delta A_c = A_c(V + \delta V, \vec A +\delta \vec A) - A_c(V, \vec A)) $

primero calculamos $A_c(V + \delta V, \vec A +\delta \vec A)$:

\

$A_c(V + \delta V, \vec A +\delta \vec A)=$

$=\int [\frac{\epsilon_o}{2} (- \nabla (V+ \delta V) - \frac{\partial (\vec A+\delta \vec A)}{\partial t})^2 - \frac{1}{2\mu_o}  (\nabla \times  (\vec A +\delta \vec A) )^2  ]dxdydzdt $

\

\

$A_c(V + \delta V, \vec A +\delta \vec A)=$

$=\int [\frac{\epsilon_o}{2} (- \nabla V- \frac{\partial \vec A}{\partial t}+(- \nabla \delta V - \frac{\partial (\delta \vec A)}{\partial t}))^2 - \frac{1}{2\mu_o}  (\nabla \times  (\vec A +\delta \vec A) )^2  ]dxdydzdt $

\

Recordemos que los cuadrados realmente denotan producto punto de una entidad vectorial por ella misma. Expandimos:

\

$A_c(V + \delta V, \vec A +\delta \vec A)$

$=\int [\frac{\epsilon_o}{2} ((- \nabla V- \frac{\partial \vec A}{\partial t})^2+(- \nabla \delta V - \frac{\partial (\delta \vec A)}{\partial t})^2 + 2(- \nabla V- \frac{\partial \vec A}{\partial t})\cdot (- \nabla \delta V - \frac{\partial (\delta \vec A)}{\partial t})) $

$ \hspace{8cm}- \frac{1}{2\mu_o}  (\nabla \times  (\vec A +\delta \vec A) )^2  ]dxdydzdt $

\

$=\int [\frac{\epsilon_o}{2} ((- \nabla V- \frac{\partial \vec A}{\partial t})^2+(- \nabla \delta V - \frac{\partial (\delta \vec A)}{\partial t})^2 + 2(- \nabla V- \frac{\partial \vec A}{\partial t})\cdot (- \nabla \delta V - \frac{\partial (\delta \vec A)}{\partial t})) $

$ \hspace{3.5cm}- \frac{1}{2\mu_o}  [(\nabla \times  \vec A)^2 +(\nabla \times \delta \vec A) ^2 + 2(\nabla \times  \vec A)\cdot (\nabla \times \delta \vec A)]  ]dxdydzdt $

\

$=\int [\frac{\epsilon_o}{2} ((- \nabla V- \frac{\partial \vec A}{\partial t})^2+(- \nabla \delta V - \frac{\partial (\delta \vec A)}{\partial t})^2 + 2(- \nabla V- \frac{\partial \vec A}{\partial t})\cdot (- \nabla \delta V - \frac{\partial (\delta \vec A)}{\partial t})) $

$ \hspace{3.5cm}- \frac{1}{2\mu_o}  [(\nabla \times  \vec A)^2 +(\nabla \times \delta \vec A) ^2 + 2(\nabla \times  \vec A)\cdot (\nabla \times \delta \vec A)]  ]dxdydzdt $

\

Como sólo estamos interesados en los términos a primer orden, simplificamos:

$A_c(V + \delta V, \vec A +\delta \vec A)=$

$=\int [\frac{\epsilon_o}{2} ((- \nabla V- \frac{\partial \vec A}{\partial t})^2  + 2(- \nabla V- \frac{\partial \vec A}{\partial t})\cdot (- \nabla \delta V - \frac{\partial (\delta \vec A)}{\partial t})) $

$ \hspace{3.5cm}- \frac{1}{2\mu_o}  [(\nabla \times  \vec A)^2   + 2(\nabla \times  \vec A)\cdot (\nabla \times \delta \vec A)]  ]dxdydzdt $

\

Detrás de nuestra inocente simplificación que conserva sólo los términos a primer orden hay una delicada suposición: una pequeña perturbación en una función, sea vectorial o sea escalar, puede tener derivadas tan grandes como se desee. Por lo tanto, cuando despreciamos sus cuadrados realmente estamos diciendo que las pequeñas perturbaciones son pequeñas tanto en los cambios de la función como en el cambio de sus derivadas. Y así todo queda arreglado salvo que uno se pregunta:  ¿quién hace esos arreglos en la naturaleza? Uno puede ignorar esa pregunta e incluir por defecto una  especie de panteísmo,  o bien, lo que uno puede decir es que nuestras leyes no son tan fundamentales como parecen y que un cambio en las derivadas de los campos exige más energía en el espacio y/o en el tiempo. Es decir, estamos haciendo un estudio de baja energía, de bajos gradientes de energía  y de baja potencia. Por ende, nuestros resultados pueden ser falsos para altos niveles de energía o de sus gradientes como en una supernova. 

\

Resumimos: como

$A_c = \int (\frac{\epsilon_o}{2} (- \nabla V - \frac{\partial \vec A}{\partial t})^2 - \frac{1}{2\mu_o}  (\nabla \times  \vec A)^2  )dxdydzdt $

entonces

\

$\delta A_c =\int [\frac{\epsilon_o}{2} (  2(- \nabla V- \frac{\partial \vec A}{\partial t})\cdot (- \nabla \delta V - \frac{\partial (\delta \vec A)}{\partial t})) $

$ \hspace{3.5cm}- \frac{1}{2\mu_o}  [  2(\nabla \times  \vec A)\cdot (\nabla \times \delta \vec A)]  ]dxdydzdt $

\

Utilizando la nomenclatura

$\vec D = \epsilon_o \vec E = \epsilon_o(- \nabla V- \frac{\partial \vec A}{\partial t})$

$\vec H = \frac{1}{\mu_o} \vec B$

podemos reexpresar $\delta A_c$ como:

$\delta A_c =\int [\vec D \cdot (- \nabla \delta V - \frac{\partial (\delta \vec A)}{\partial t})- \vec H \cdot (\nabla \times \delta \vec A) ]dxdydzdt $

$\hspace{0.75cm} =\int [-\vec D \cdot   \nabla \delta V - \vec D \cdot \frac{\partial (\delta \vec A)}{\partial t}- \vec H \cdot (\nabla \times \delta \vec A) ]dxdydzdt $

\

Ahora recordemos una fórmula que sale de la de integración por partes:

$\int_\Omega \nabla u \cdot \vec v  dx = $
$\int_{\delta \Omega} u   \vec v \cdot \nu d\sigma -$
$\int_\Omega u  \nabla \cdot \vec vdx$

Si en esta fórmula olvidamos las condiciones de frontera, quizá porque no ella no exista, nos queda:

$\int_\Omega \nabla u \cdot \vec v  dx =  - \int_\Omega u  \nabla \cdot \vec vdx$

Apliquemos esta versión sobre las coordenadas espaciales de:

$\int [-\vec D \cdot   \nabla \delta V]dxdydzdt = \int [-\nabla \delta V \cdot \vec D   ]dxdydzdt$

Nos queda:

\

$\int [-\vec D \cdot   \nabla \delta V]dxdydzdt = \int [\delta V \nabla\cdot \vec D   ]dxdydzdt$

\

Por otra parte, si en la siguiente fórmula de integración por partes

$\int_\Omega \frac{\partial u}{\partial x_i}v dx = $
$\int_{\delta \Omega} uv \nu_i d\sigma -$
$\int_\omega u \frac{\partial v}{\partial x_i}dx$

olvidamos la componente de la frontera, nos da:

$\int_\Omega \frac{\partial u}{\partial x_i}v dx =  -$
$\int_\omega u \frac{\partial v}{\partial x_i}dx$

Aplicando esto sobre

$\int ( - \vec D \cdot \frac{\partial (\delta \vec A)}{\partial t})dxdydzdt$

obtenemos:

\

$\int ( - \vec D \cdot\frac{\partial (\delta \vec A)}{\partial t})dxdydzdt$ 
$= \int  \frac{\partial  \vec D}{\partial t}   \cdot \delta \vec A dxdydzdt$

\

Probemos ahora que

$\int \vec H \cdot \nabla \times G dxdydz =$
$ \int (\nabla \times \vec H) \cdot \vec G dxdydz$

En efecto:

$\int \vec H \cdot \nabla \times G dxdydz =$

$=\int \vec H \cdot [\frac{\partial G_3}{\partial y } - \frac{\partial G_2 }{\partial z} ,
\frac{\partial G_1}{\partial z} - \frac{\partial G_3}{\partial x },
\frac{\partial G_2}{\partial x } - \frac{\partial G_1}{\partial y} ] dxdydz$

\

$=\int [H_1(\frac{\partial G_3}{\partial y } - \frac{\partial G_2 }{\partial z} ) 
+ H_2 (\frac{\partial G_1}{\partial z} - \frac{\partial G_3}{\partial x })
+ H_3(\frac{\partial G_2}{\partial x } - \frac{\partial G_1}{\partial y} )] dxdydz$

\

$=\int [H_1\frac{\partial G_3}{\partial y } - H_1\frac{\partial G_2 }{\partial z} 
+ H_2 \frac{\partial G_1}{\partial z} - H_2\frac{\partial G_3}{\partial x }
+ H_3 \frac{\partial G_2}{\partial x } - H_3 \frac{\partial G_1}{\partial y} ] dxdydz$

\

Ahora aplicamos nuestra versión de la fórmula de integración por partes y nos da

\

$=\int [-G_3\frac{\partial H_1}{\partial y } + G_2\frac{\partial H_1 }{\partial z} 
- G_1 \frac{\partial H_2}{\partial z} + G_3\frac{\partial H_2}{\partial x }
- G_2 \frac{\partial H_3}{\partial x } + G_1 \frac{\partial H_3}{\partial y} ] dxdydz$

\

Ordenando adecuadamente obtenemos

$=  \int [   G_1 \frac{\partial H_3}{\partial y} -G_1 \frac{\partial H_2}{\partial z}+ G_2\frac{\partial H_1 }{\partial z}  -
G_2 \frac{\partial H_3}{\partial x }  
 + G_3\frac{\partial H_2}{\partial x }- G_3\frac{\partial H_1}{\partial y }
 ] dxdydz$

$=  \int [G_1 (\frac{\partial H_3}{\partial y}-\frac{\partial H_2}{\partial z} )
 +G_2 (   \frac{\partial H_1 }{\partial z}-\frac{\partial H_3}{\partial x } )
+ G_3(  \frac{\partial H_2}{\partial x }-\frac{\partial H_1}{\partial y })
 ] dxdydz$

$= \int \vec G \cdot ( \nabla \times H) dxdydz $

\

$ =  \int (\nabla \times \vec H) \cdot \vec G dxdydz$

\

Hemos demostrado entonces que:

$\int \vec H \cdot \nabla \times G dxdydz  =  \int (\nabla \times \vec H) \cdot \vec G dxdydz$

\

Este resultado aplicado sobre las coordenadas espaciales de

$\int \vec H \cdot (\nabla \times \delta \vec A)  dxdydzdt $

nos la convierte en

$\int \vec H \cdot (\nabla \times \delta \vec A)  dxdydzdt =  $
$\int \nabla \times \vec H \cdot   \delta \vec A  dxdydzdt$ 

\

 Reuniendo todas las aplicaciones de la metodología de integración por partes, partimos de   
 
 \
 
 $\delta A_c =\int [\vec D \cdot (- \nabla \delta V - \frac{\partial (\delta \vec A)}{\partial t})- \vec H \cdot (\nabla \times \delta \vec A) ]dxdydzdt $

$\hspace{0.75cm} =\int [-\vec D \cdot   \nabla \delta V - \vec D \cdot \frac{\partial (\delta \vec A)}{\partial t}- \vec H \cdot (\nabla \times \delta \vec A) ]dxdydzdt $

\

y obtenemos

\

$\delta A_c = \int [\delta V \nabla\cdot \vec D  + \frac{\partial  \vec D}{\partial t}   \cdot \delta \vec A  -\nabla \times \vec H \cdot   \delta \vec A ]dxdydzdt $

$\hspace{0.75cm} = \int [ \nabla\cdot \vec D \delta V  + (\frac{\partial  \vec D}{\partial t}     -\nabla \times \vec H) \cdot   \delta \vec A ]dxdydzdt $

\

Igualando esto a cero y asumiendo que las variaciones de las variables fundamentales $V$ y $\vec A$ son independientes, concluimos que

1) $\nabla\cdot \vec D = 0$: la divergencia del campo eléctrico en el vacío es cero y que

2) $\frac{\partial  \vec D}{\partial t}     =\nabla \times \vec H$: la variación temporal del campo eléctrico simula un campo magnético.

\

En realidad, la acción propuesta por los Cosserat no fue diseñada para describir la conducta del campo electromagnético en el vacío sino que también incluye la posibilidad de que existan cargas y de que hay interacción entre el campo y las cargas. Resultamos en el vacío debido a que estamos estudiando el problema por pedazos. Tratemos de volver a la acción original:

\

$A = \int (\frac{\epsilon_o}{2} (- \nabla V - \frac{\partial \vec A}{\partial t})^2 - \frac{1}{2\mu_o}  (\nabla \times  \vec A)^2 + \rho \vec v \cdot \vec A - \rho V + m \frac{(\vec v)^2  }{2})dxdydzdt $

\

La primera parte de esta acción ya la estudiamos, nos queda por estudiar la acción asociada a la partícula y a su interacción con el campo:

\

$A_p =   \int (    \rho \vec v \cdot \vec A - \rho V + m \frac{(\vec v)^2  }{2})dxdydzdt $

\

Calculemos 

$\delta A_p = A_p(\vec A + \delta \vec A, V+ \delta V, \vec v + \delta \vec v,  )$
$ -  A_p(\vec A, V,    \vec v)$.

\

Si hacemos bien nuestro cálculo, deberíamos llegar a un resultado que incluya a la fuerza de Lorentz. Y como ya sabemos a dónde llegar, cuadramos las cosas para lograrlo fácilmente. Esto convierte nuestra empresa en algo semejante a un arte, haciendo cosas precisamente para que todo funcione. 

Pongamos pues atención a la siguiente observación: cuando decimos que la acción es una integral sobre el espacio-tiempo, estamos diciendo que si algo pasa, entonces no pasa por fuera del espacio-tiempo. En particular, cuando estudiamos la acción $A_c$ asociada al campo en solitario, consideramos $\delta V$ y $\delta \vec A$. Pero eso sucede en el espacio-tiempo. Como esos campos lo llenan todo, por derecho propio representan una variación en el espacio, pero cuando decimos que pueden variar como entidades que existen en sí mismas, estamos diciendo que pueden variar en el tiempo. Por consiguiente, siempre estamos contabilizando las variaciones en el espacio y en el tiempo.
En los siguientes cálculos, eso se verá separando las variaciones debidas al espacio de las que suceden en el tiempo, para lo cual se utilizarán los sub-índices $E$ y $t$:

\

Ahora bien, si

$A_p =   \int (    \rho \vec v \cdot \vec A - \rho V + m \frac{(\vec v)^2  }{2})dxdydzdt $

entonces

$A_p( \vec A + \delta \vec A, V+ \delta V,   \vec v + \delta \vec v )$

$ = \int (\rho  (\vec v + \delta \vec v) \cdot (\vec A + \delta_t \vec A+ \delta_E \vec A) -  \rho   (V + \delta_t V + \delta_E V) + m\frac{(\vec v + \delta \vec v)^2  }{2})dxdydzdt $

$ = \int [  \rho \vec v \cdot \vec A +   \rho \delta \vec v \cdot \vec A  + \rho \vec v \cdot  \delta_t \vec A + \rho \vec v \cdot  \delta_E \vec A$

 $ \hspace{1.3cm} -\rho   V  -\rho  \delta_t V-\rho  \delta_E V $
 $ + m\frac{ \vec v  ^2}  {2} +    m  \vec v  \cdot  \delta \vec v      ]dxdydzdt $

donde sólo hemos tenido en cuenta las perturbaciones a primer orden.

\

Por tanto

 $\delta A_p = \int [      \rho \delta \vec v \cdot \vec A  + \rho \vec v \cdot  \delta_t \vec A + \rho \vec v \cdot  \delta_E \vec A$

 $ \hspace{1.3cm}   -\rho  \delta_t V-\rho  \delta_E V $
 $ + m  \vec v  \cdot  \delta \vec v      ]dxdydzdt $

 \

 Sumando los aportes a la variación total tanto de la acción del campo como de la acción asociada a la partícula y a su interacción con el campo, obtenemos:
 
 \
 
 $\delta A = \int [ \nabla\cdot \vec D \delta V  + (\frac{\partial  \vec D}{\partial t}     -\nabla \times \vec H) \cdot   \delta \vec A $
 
 $\hspace{1.5cm} + \rho \delta \vec v \cdot \vec A  + \rho \vec v \cdot  \delta_t \vec A + \rho \vec v \cdot  \delta_E \vec A $

 $ \hspace{1.5cm}  -\rho  \delta_t V-\rho  \delta_E V $
 $ +     m  \vec v  \cdot  \delta \vec v  ]dxdydzdt $ 
 
 \

 Ahora decimos que la variación en el tiempo de los campos realmente corresponde a su variación intrínseca:
 
  $\delta A = \int [ \nabla\cdot \vec D \delta V  + (\frac{\partial  \vec D}{\partial t}     -\nabla \times \vec H) \cdot   \delta \vec A $
 
 $\hspace{1.5cm} + \rho \delta \vec v \cdot \vec A  + \rho \vec v \cdot  \delta \vec A + \rho \vec v \cdot  \delta_E \vec A $

 $ \hspace{1.5cm}  -\rho  \delta V-\rho  \delta_E V $
 $ +     m  \vec v  \cdot  \delta \vec v  ]dxdydzdt $

 \
  
   $\delta A = \int [ \nabla\cdot \vec D \delta V  + (\frac{\partial  \vec D}{\partial t}     -\nabla \times \vec H) \cdot   \delta \vec A $
 
 $\hspace{1.5cm}-\rho  \delta V  +  \rho \vec v \cdot  \delta \vec A$

 $ \hspace{1.5cm} + \rho \delta \vec v \cdot \vec A   + \rho \vec v \cdot  \delta_E \vec A -\rho  \delta_E V $
 $ +     m  \vec v  \cdot  \delta \vec v  ]dxdydzdt $ 
 
 \
 
 que también puede escribirse como:
 
 \

   $\delta A = \int [ (\nabla\cdot \vec D  -\rho) \delta V  + (\frac{\partial  \vec D}{\partial t}     -\nabla \times \vec H +\rho \vec v ) \cdot   \delta \vec A $
 
 $\hspace{1.5cm} + \rho \delta \vec v \cdot \vec A     + \rho \vec v \cdot  \delta_E \vec A -\rho  \delta_E V $
 $ +     m  \vec v  \cdot  \delta \vec v  ]dxdydzdt $

 \

 Igualando a cero las componentes en $dV$ y en $dA$   obtenemos:
 
 $\nabla\cdot \vec D  =\rho $: la carga es la única fuente de campo eléctrico.
 
 $\frac{\partial  \vec D}{\partial t} +\rho \vec v    = \nabla \times \vec H $: una carga en movimiento o una corriente eléctrica, lo mismo que un campo eléctrico variable, crean un campo magnético.

 \
 
 Es muy interesante notar que en la 'deducción' de dos de las leyes de Maxwell que acabamos de hacer  hemos mezclado el campo electromagnético, el cual está definido sobre todo el espacio-tiempo, con la velocidad de una partícula, la cual no es un campo definido sobre todo el espacio-tiempo sino sólo sobre su trayectoria. 
 
 \
 
 Es   vergonzoso que utilicemos algo sin sentido  para construir una maquinaria que nos funcione.  Es por eso que hablamos de arte. ¿Qué es lo que permite que eso pueda suceder?  Aparte de que de una mentira siempre puede deducirse una verdad, ¿será posible que aquí tengamos una premonición de la mecánica cuántica en la cual a la partícula se le asocia un campo (complejo) definido sobre todo el espacio-tiempo?
 En cualquier caso, vemos que no es necesario llegar hasta la teoría de campos para ver cosas sin sentido pero que funcionan. 
 
 \

El mismo misterio de mezclar campos y trayectorias lo encontramos en la variación asociada a la partícula, la cual nos falta por estudiar:
 
 \
 
    $\delta A_p = \int [  \rho \delta \vec v \cdot \vec A     + \rho \vec v \cdot  \delta_E \vec A -\rho  \delta_E V $
 $ +     m  \vec v  \cdot  \delta \vec v   ]dxdydzdt  $
  
  \
    
donde integramos sobre un espacio de 4-dimensiones algo qué sólo está bien definido en un subconjunto, quizá un subvariedad, de dimensión uno. Nuestra salida  es utilizar el teorema de Fubini para reescribir todo como una integral con respecto al tiempo sobre una trayectoria  $\vec r(t)$ cuya velocidad sea $\vec v(t)$. Las otras integrales  son ignoradas. 

Pero en ese caso, estamos precisamente estudiando la variación de Lagrangiano al que llegamos en el ejercicio anterior: si nos regresamos  por el mismo camino desde el Lagrangiano debemos llegar a la expresión de la fuerza de Lorentz. Podemos creerlo pero es más instructivo deducirlo. Así que deduzcámoslo. Supondremos que  las variaciones son independientes coordenada por coordenada mientras que se conservan invariantes las condiciones de borde. Por ello podremos utilizar nuestra fórmula ultracorta para integración por partes, en la cual no hay efectos de borde:

\

$\delta A_p = \int [  \rho \delta \vec v \cdot \vec A     + \rho \vec v \cdot  \delta_E \vec A -\rho  \delta_E V $
 $ +     m  \vec v  \cdot  \delta \vec v   ] dt  $
  
  $\delta A_p = \int [-\rho  \delta_E V    + \rho \vec v \cdot  \delta_E \vec A +  \rho \delta \vec v \cdot \vec A   $
 $ +     m  \vec v  \cdot  \delta \vec v   ] dt  $
 
\

  $ = \int [   - \Sigma\rho\frac{\partial V}{\partial x_i}  dx_i  + \Sigma \rho  \vec v \cdot \frac{\partial  \vec A    }{\partial x_i}  dx_i    
 +  \rho \delta \vec \frac{d\vec r}{dt} \cdot \vec A            +     m \vec v   \cdot  \delta \frac{d\vec r}{dt}       ] dt  $

  $ = \int [  -  \Sigma\rho\frac{\partial V}{\partial x_i}  dx_i  + \Sigma \rho  \vec v \cdot \frac{\partial  \vec A    }{\partial x_i}  dx_i    
 +  \rho \vec \frac{d\delta \vec r}{dt} \cdot \vec A   +     m \vec v   \cdot  \frac{  d \delta \vec r}{dt}       ] dt  $

  $ = \int [ -   \Sigma\rho\frac{\partial V}{\partial x_i}  dx_i  + \Sigma \rho  \vec v \cdot \frac{\partial  \vec A    }{\partial x_i}  dx_i    
 +  \Sigma \rho  \frac{d\delta  x_i}{dt}    A_i   +     \Sigma  m  v_i    \frac{  d \delta   r_i}{dt}       ] dt  $

  $ = \int [  -  \Sigma\rho\frac{\partial V}{\partial x_i}  dx_i  + \Sigma \rho  \vec v \cdot \frac{\partial  \vec A    }{\partial x_i}  dx_i   -
   \Sigma \rho \frac{d      A_i }{dt}     \delta  x_i        -      \Sigma  \frac{d ( m   v_i) }{dt}      \delta  x_i      ] dt  $

\

Igualando a cero el integrando coordenada por coordenada obtenemos el siguiente sistema de ecuaciones:

\

 $    -\rho\frac{\partial V}{\partial x_i}  dx_i  +   \rho  \vec v \cdot \frac{\partial  \vec A    }{\partial x_i}  dx_i   -
     \rho \frac{d      A_i }{dt}     \delta  x_i        -        \frac{d ( m   v_i) }{dt}      \delta  x_i   =0  $    
   
  $- \rho \frac{\partial V}{\partial x_i}   +   \rho  \vec v \cdot \frac{\partial  \vec A    }{\partial x_i}     -
     \rho \frac{d      A_i }{dt}               =         \frac{d ( m   v_i) }{dt}    $
                    
                    $        \frac{d ( m   v_i) }{dt}  = -\rho  \frac{\partial V}{\partial x_i}    +   \rho  \vec v \cdot \frac{\partial  \vec A    }{\partial x_i}     -
     \rho \frac{d      A_i }{dt}              $

       $     (\vec F)_i  = -   \rho \frac{\partial V}{\partial x_i}  +   \rho  \vec v \cdot  \frac{\partial   \vec A  }{\partial x_i}        
                 -\rho\frac{d      A_i  }{dt}       $
                 
    \

 Utilicemos la regla de la   cadena:

$\frac{d      A_i }{dt} =\frac{\partial A_i}{\partial x}\dot{x}+\frac{\partial A_i}{\partial y}\dot{y} +\frac{\partial A_i}{\partial z}\dot{z}=\dot{x}\frac{\partial A_i}{\partial x}+\dot{y}\frac{\partial A_i}{\partial y} +\dot{z}\frac{\partial A_i}{\partial z}$

\

$ (\vec F)_i  = -\rho   \frac{\partial V}{\partial x_i} + \rho[\dot{x}\frac{\partial A_1}{\partial x_i}+\dot{y}\frac{\partial A_2}{\partial x_i}+\dot{z}\frac{\partial A_3}{\partial x_i} -(\dot{x}\frac{\partial A_i}{\partial x }+\dot{y}\frac{\partial A_i}{\partial y} +\dot{z}\frac{\partial A_i}{\partial z}) ]        $

\

Para la coordenada $x$, esta ecuación es:

$ (\vec F)_1  = -\rho   \frac{\partial V}{\partial x } + \rho[\dot{x}\frac{\partial A_1}{\partial x }+\dot{y}\frac{\partial A_2}{\partial x }+\dot{z}\frac{\partial A_3}{\partial x } -(\dot{x}\frac{\partial A_1}{\partial x }+\dot{y}\frac{\partial A_1}{\partial y} +\dot{z}\frac{\partial A_1}{\partial z}) ]        $

$ (\vec F)_1  = -\rho   \frac{\partial V}{\partial x } + \rho[\dot{y}\frac{\partial A_2}{\partial x}+\dot{z}\frac{\partial A_3}{\partial x} -(\dot{y} \frac{\partial A_1}{\partial y}+\dot{z}\frac{\partial A_1}{\partial z}) ]        $

\bigskip

Reunamos todas estas ecuaciones para coordenadas en una sola ecuación para  vectores: 

\

$$\vec F =   \rho\{-\pmatrix{\frac{\partial V}{\partial x} \cr
\frac{\partial V}{\partial y} \cr
\frac{\partial V}{\partial z}  } +   \pmatrix{\dot{y}[\frac{\partial A_2}{\partial x}-\frac{\partial A_1}{\partial y}] -\dot{z} [\frac{\partial A_1}{\partial z}-\frac{\partial A_3}{\partial x}]  \cr
\dot{z}[\frac{\partial A_3}{\partial y}-\frac{\partial A_2}{\partial z}]-\dot{x}[\frac{\partial A_2}{\partial x}-\frac{\partial A_1}{\partial y}]   \cr
 \dot{x}[\frac{\partial A_1}{\partial z}-\frac{\partial A_3}{\partial x}]-\dot{y}[\frac{\partial A_3}{\partial y}-\frac{\partial A_2}{\partial z}]    } \}
     $$

     $$  \vec F = \rho(-\nabla  V +   \vec v \times (\nabla \times \vec A))
     $$

 Teniendo en cuenta que
  
   $\vec E = -\nabla  V $

y que  

$\vec B = \nabla \times \vec A $

Obtenemos la Fuerza de Lorentz:

$$\vec F = (\vec E +   \vec v \times \vec B)
     $$

 Al final de cuentas tenemos que el Lagrangiano de los Cosserat contiene tanto dos de las ecuaciones del campo electromagnético como la fuerza de Lorentz. Nos demoraremos un buen trecho antes de proponer una forma de obtener las otras dos leyes que nos faltan. 
 
Quedamos intrigados por haber trabajado un programa exitoso donde nos inventamos las reglas a medida que las necesitábamos. Es posible que eso tenga una explicación, es decir pueda  que exista un modelo en el cual nuestro arte se cambie por una deducción rigurosa. Pero es posible que tal modelo no exista: gracias a  los teoremas de incompletitud de Gödel  hemos logrado saber que la naturaleza, siendo un ente complejo y mucho más que los números naturales, no puede ser reducida por completo a un modelo consistente, sino que hay que aparece un compromiso: o uno hace un modelo que lo explique todo pero con  contradicciones o bien, uno tiene modelos incompletos y autoconsistentes a los cuales uno tiene que adicionar reglas ad hoc que uno debe pegar a fuerza bruta para que le funcionen las cosas.

\addtocounter{ecu}{1}

\section{EL HAMILTONIANO}

Tenemos que hallar el  \index{Hamiltoniano} \textbf{Hamiltoniano} que corresponde al Lagrangiano propuesto. Podríamos ir por el formalismo de la acción y después por la integral de Feynman sobre todos los caminos. Aunque eso es posible, es muy tedioso. Contamos, sin embargo, con otras alternativas para llegar al Hamiltoniano. La que vamos a describir ahora acorta mucho el trabajo.

Puesto que el Hamiltoniano es el operador que describe la observable energía, entonces debemos hallar la correspondiente expresión clásica, o sea la energía total. Lo deseable es hallar una expresión que sea intrínseca, es decir, que tenga la misma expresión en todos los sistemas de coordenadas. Eso se hace como sigue:

Si

$$L=L(q\sb{1},..,q\sb{n},\dot q\sb{1},..,\dot q\sb{n})
\addtocounter{ecu}{1}   \hspace{4cm} (\theecu )      $$

Entonces se cambia de variable,  las velocidades $\dot q\sb{i}$ se reexpresan utilizando  los momentos

$$p\sb{i}=\partial L/\partial \dot q\sb{i}
\addtocounter{ecu}{1}   \hspace{4cm} (\theecu )      $$

y se define $H$ por

$$H = H( q\sb{1},..,q\sb{n},p\sb{1},..,p\sb{n}) = \sum p\sb{i}\dot q\sb{i}-L
\addtocounter{ecu}{1}   \hspace{4cm} (\theecu )      $$

Por ejemplo, si estamos en el caso de una partícula que se mueve en un campo conservativo,

$L=T-V = (1/2)m \|\vec v \|\sp{2} - V $
$=(1/2m) ( p\sb{x}\sp{2} + p\sb{y}\sp{2} +p\sb{z}\sp{2})-V $

Entonces $H$ se calcula como

$H= \sum p\sb{i}\dot x\sb{i}-L =$
$ \sum p\sb{i}\dot x\sb{i}-(1/2)m \|\vec v \|\sp{2} + V =$
$ \sum mv\sb{i}v\sb{i}-(1/2)m \|\vec v \|\sp{2} + V = $
$m\|\vec v \|\sp{2}-(1/2)m \|\vec v \|\sp{2} + V= $
$(1/2)m \|\vec v \|\sp{2} + V$

$=(1/2m) ( p\sb{x}\sp{2} + p\sb{y}\sp{2} +p\sb{z}\sp{2})+V$

Vemos que se recupera la energía total.

\

Veamos otra forma de cuantizar, es decir de hallar la expresión del operador Hamiltoniano. Es una receta que no hace más que automatizar los delicados procesos que ya hemos visto.

Teniendo el Hamiltoniano clásico, la regla de cuantización, es decir, la forma de crear el Hamiltoniano cuántico tiene tres puntos:

\

1) Se escribe el Hamiltoniano clásico en coordenadas cartesianas.

2) Al potencial clásico $V$ le corresponde el operador cuántico $V$, el cual lo único que hace es tomar una función de onda y multiplicarla por la función escalar $V$.

3) Al momento en una dirección dada, le corresponde el operador derivada parcial en la misma dirección multiplicado por $\hbar / i$.

\

Pongamos por caso que  se trate de una partícula que se mueve bajo la acción de un campo de fuerzas conservativo, entonces la energía total clásica es

$H = (1/2)m \|\vec v \|\sp{2} + V = (1/2m) \|\vec p\|\sp{2} + V $
$=(1/2m) ( p\sb{x}\sp{2} + p\sb{y}\sp{2} +p\sb{z}\sp{2})+V $

Eso significa que el Hamiltoniano cuántico es:

$H=(1/2m) [( (\hbar/i)\partial /\partial x) \sp{2} $
$+ ( (\hbar/i)\partial /\partial y) \sp{2}+ $
$( (\hbar/i)\partial /\partial z) \sp{2}]+V $

$=(1/2m) (\hbar/i)\sp{2}[( \partial /\partial x) \sp{2} $
$+ ( \partial /\partial y) \sp{2}+ (\partial /\partial z) \sp{2}]+V $

$=-(\hbar\sp{2}/2m)\nabla \sp{2} +V$

Que es nuestro viejo amigo. Ahora repetiremos lo mismo para el campo electromagnético:

$L= (1/2)m \|\vec v \|\sp{2} - e \phi + (e/c) \vec v \cdot \vec A $

$=(1/2)m ( v\sb{x}\sp{2} + v\sb{y}\sp{2} +v\sb{z}\sp{2})-$
$ e \phi + (e/c) \vec v \cdot \vec A $

$p\sb{i} = \partial L/\partial \dot x\sb{i}= m v\sb{i} + (e/c) A\sb{i}$

 Reescribiendo esta igualdad en forma vectorial obtenemos:

$\vec p = m\vec v +(e/c)\vec A \Rightarrow \vec v = (1/m)(\vec p - (e/c)\vec A)$

Como $H$ se define por

$H = \sum p\sb{i}\dot q\sb{i}-L = \vec p \cdot \vec v -L$

$=(m\vec v +(e/c)\vec A)\cdot \vec v - ((1/2)m \|\vec v \|\sp{2} $
$- e \phi + (e/c) \vec v \cdot \vec A) $

$=m\vec v \cdot \vec v+(e/c)\vec A\cdot \vec v $
$- (1/2)m \vec v \cdot \vec v + e \phi - (e/c) \vec v \cdot \vec A$

$= (1/2)m \vec v \cdot \vec v + e \phi $

Pero, puesto que $\vec v = (1/m)(\vec p - (e/c)\vec A)$, entonces

$H= (1/2)m (1/m)(\vec p - $
$(e/c)\vec A) \cdot (1/m)(\vec p - (e/c)\vec A) + e \phi $

$$H= (1/2m)(\vec p -(e/c)\vec A) \cdot (\vec p -(e/c)\vec A) + e \phi
\addtocounter{ecu}{1}   \hspace{4cm} (\theecu )      $$

Conociendo  $H$, el operador de evolución queda como de costumbre:

$\phi(x,t)= e\sp{-itH/\hbar} \phi\sb{0}$

Ahora estamos listos para estudiar el efecto de la arbitrariedad en la definición de $\vec A$.

Probaremos el siguiente teorema:

\

\

\color{red}
\textit{El electromagnetismo es una teoría gauge con \index{grupo de invariancia local} \textbf{grupo de invariancia local} } \textbf{U(1)}.
\

\color{black}

\bigskip

En palabras más entendibles, lo que esa jerga significa es lo siguiente:

\

\color{red}
\textit{El Hamiltoniano cuántico asociado a la interacción electromagnética admite una arbitrariedad en la definición del potencial vector, pues se le puede sumar el gradiente de una función escalar cualquiera. La física de la interacción electromagnética no cambia aunque cambiemos de arbitrariedad, pues eso equivale a un cambio de fase que siendo local, tiene que ser universal. Cambiar de fase equivale a una multiplicación de la función de onda por un elemento de} \textbf{U(1)}, \textit{es decir por una exponencial imaginaria. Al cambiar de punto se puede también cambiar de fase (cambio local), pero en el mismo punto todas las funciones de onda deben tener el mismo cambio de fase (es un universal).}

\

\

\color{black}
Este teorema exige una demostración rigurosa porque, como ya sabemos, un cambio local de fase en un proceso puede crear resultados observables. Acá estamos hablando de cambios locales de fase y sin embargo estamos prometiendo que eso no importa.

No probaremos el teorema de forma directa sino que utilizaremos una maquinaria muy poderosa y muy útil en varios contextos. Cada término utilizado en el enunciado del teorema quedará nítidamente definido a lo largo de la demostración y sólo ahí.

\section{EL CONMUTADOR}

Sean $A,B,C$ operadores cualesquiera definidos en un espacio de Hilbert dado, $L$. Definimos el \textbf{conmutador} \index{conmutador} entre $A$ y $B$, notado $[A,B]$, como:

$$[A,B]=AB-BA
\addtocounter{ecu}{1}   \hspace{4cm} (\theecu )      $$

El espacio de Hilbert podría ser $R\sp{n}$ y los operadores transformaciones lineales, o matrices. El conmutador mide qué tanto no conmutan  dos operadores dados. Las propiedades básicas del conmutador están consignadas en el siguiente teorema, válido para operadores lineales:

\

\addtocounter{ecu}{1}
\textit{Teorema \theecu. El conmutador de operadores lineales cumple:}

1. $[A,A]=0$

2. $[A,B] = - [B,A]$

3. $[A,c]=0, c \in C$, en general, $[A,f(x)]=0$, donde $f(x)$ es una función escalar.

4. $[A,B+C] = [A,B] + [A,C]$

5. $[A,BC] =[A,B]C + B[A,C]$

6. $[A,B^2] =[A,B]B + B[A,B]$

7. $[A,[B,C]] +[B,[C,A]] +[C,[A,B]] =0 $, la identidad de Jacobi.

\

Demostración:

1. $[A,A]= AA-AA=0$

2. $[A,B] = AB-BA= -(BA-AB)= - [B,A]$

3. Sea $\vec v \in L$, y $c$ un escalar, entonces
$[A,c](\vec v) = A(c\vec v) - cA(\vec v)=  cA(\vec v)- cA(\vec v)=0$, de igual forma, $[A,f(x)]=0$, donde $f(x)$ es una función escalar.

4. $[A,B+C] = A(B+C) - (B+C)A= AB + AC - BA - CA $
$= AB - BA + AC  - CA =[A,B] + [A,C]$

5. $[A,BC] =ABC - BCA= ABC-BAC + BAC - BCA $
$= (AB-BA)C + B(AC-CA)= [A,B]C + B[A,C]$

6. $[A,B^2] =[A, BB] =  [A,B]B + B[A,B]$

7. $[A,[B,C]] +[B,[C,A]] +[C,[A,B]] =[A,BC-CB] +[B,CA-AC] +[C,AB-BA]$

$=[A,BC]-[A,CB] +[B,CA]-[B,AC] +[C,AB]-[C,BA]$

$=ABC-BCA-ACB+CBA+BCA-CAB - BAC+ACB + CAB-ABC-CBA+BAC $

$=ABC-ABC-BCA+BCA-ACB+ACB +CBA-CBA-CAB + CAB- BAC+BAC =0 $

\

Ahora vamos a demostrar una serie de propiedades a tono con las funciones de onda, es decir, tomaremos como espacio de Hilbert, $L = L\sb{2}$, el conjunto de las funciones que van de $R\sp{3n}$ en los complejos $C$ y que sean  integrables en el sentido de la norma inducida por el producto escalar. Tales funciones pueden ser multiplicadas por funciones escalares, de tal forma que tiene estructura de módulo. En particular, pueden ser multiplicadas por las coordenadas de una función vectorial. Por lo tanto, dichas multiplicaciones definen operadores.

\

Definimos el 'operador  coordenada sub-i', $\textbf{x} \sb{i}$ por

$$ \textbf{x}\sb{i}(\psi) =x\sb{i}\psi
\addtocounter{ecu}{1}   \hspace{4cm} (\theecu )      $$

El operador 'coordenada momento sub -j' $\textbf{p}\sb{j}$, por

$$\textbf{p}\sb{j}(\psi) = (\hbar/i)\partial \psi /\partial x\sb{j}
\addtocounter{ecu}{1}   \hspace{4cm} (\theecu )      $$

Los operadores 'posición' y 'momento' se definen como

$$\textbf{x}= \textbf{x}\sb{1}\vec i+ \textbf{x}\sb{2}\vec j + \textbf{x}\sb{3}\vec{k}
\addtocounter{ecu}{1}   \hspace{4cm} (\theecu )      $$

y por

$$\textbf{p} = \textbf{p}\sb{1}\vec i+ \textbf{p}\sb{2}\vec j + \textbf{p}\sb{3}\vec{k}
\addtocounter{ecu}{1}   \hspace{4cm} (\theecu )      $$

Podemos entonces formular y probar el siguiente

\

\addtocounter{ecu}{1}
\textit{Teorema \theecu.  Relaciones de conmutación de los operadores posición y momento:}

1. $[\textbf{x}\sb{i},\textbf{x}\sb{j} ]=0$

2. $[\textbf{p}\sb{i},\textbf{p}\sb{j} ]=0$

3. $[\textbf{x}\sb{i},\textbf{p}\sb{j} ]=i\hbar \delta \sb{ij}$

4.$[f(\textbf{x}),g(\textbf{x})]=0$

5. $[f(\textbf{p}),g(\textbf{p})]=0$

\

Probemos la tercera identidad, teniendo cuidado de no confundir el número complejo $i$ con la coordenada número $i$:

$[\textbf{x}\sb{i},\textbf{p}\sb{j}](\psi)=$
$\textbf{x}\sb{i}\textbf{p}\sb{j} (\psi)-$
$\textbf{p}\sb{j} \textbf{x}\sb{i}(\psi)=$
$ x\sb{i}(\hbar /i) \partial \psi/\partial x\sb{j}$
$- (\hbar /i) \partial (x\sb{i}\psi)/\partial x\sb{j}$

$= x\sb{i}(\hbar /i) \partial \psi/\partial x\sb{j}- $
$(\hbar /i) (\partial x\sb{i}/\partial x\sb{j})\psi- $
$(\hbar /i) x\sb{i}\partial \psi /\partial x\sb{j}$

$=-(\hbar /i) (\partial x\sb{i}/\partial x\sb{j})\psi$

$=i\hbar\delta \sb{ij}\psi$

\

Ahora viene un teorema un poco más elaborado, en el cual $G$ y $F$  son funciones escalares del operador posición y momento, respectivamente:

\

\addtocounter{ecu}{1}
\textit{Teorema \theecu}

1.$[\textbf{x}\sb{i},F(\textbf{p})]=i\hbar \partial F/\partial p\sb{i}$

2.$[\textbf{x},F(\textbf{p})]=[\sum\textbf{x}\sb{i}\vec e\sb{i}, F(\textbf{p})]  =\sum [\textbf{x}\sb{i}, F(\textbf{p})]\vec e\sb{i}$

$=\sum i\hbar \partial F(\textbf{p})/\partial p\sb{i}\vec e\sb{i}$
$= i\hbar \nabla \sb{(p\sb{i})}F(\textbf{p})$

3. $[\textbf{p}\sb{i},G(\textbf{x})]$
$=i\hbar \partial G/\partial x\sb{i}$

4.$[\textbf{p},G(\textbf{x})]=$
$[\sum\textbf{p}\sb{i}\vec e\sb{i}, G (\textbf{x})] $
$ =\sum [\textbf{p}\sb{i}, G(\textbf{x})]\vec e\sb{i}$

$=\sum i\hbar \partial G(\textbf{p})/\partial x\sb{i}\vec e\sb{i}$
$= i\hbar \nabla G(\textbf{x})$

5.$[\textbf{p}\sp{2},G(\textbf{x})]= $
$-i\hbar (\nabla G)\textbf{p} -i\hbar \textbf{p}(\nabla G)$

\

Demostremos las propiedades 3 y 5, válidas para el caso en el cual $G$ pueda desarrollarse en serie.

Demostración de la 3. Probemos el teorema para un monomio:

$[\textbf{p}\sb{i},\textbf{x}\sb{j}\sp{k}]=?$

Para $k=0$: $[\textbf{p}\sb{i},\textbf{x}\sb{j}\sp{0}]=[\textbf{p}\sb{i},1]=0$

Para $k=1$:$[\textbf{p}\sb{i},\textbf{x}\sb{j}]=-i\hbar \delta \sb{ij}$
$= -i\hbar\partial x\sb{j}/\partial x\sb{i} $

Para $k=2$:

$[\textbf{p}\sb{i},\textbf{x}\sb{j}\sp{2}]=$
$[\textbf{p}\sb{i},\textbf{x}\sb{j}\textbf{x}\sb{i}]$

$=[\textbf{p}\sb{i},\textbf{x}\sb{j}]\textbf{x}\sb{j} +$ $\textbf{x}\sb{j}[\textbf{p}\sb{i},\textbf{x}\sb{j}]$

$=-i\hbar \delta\sb{ij}\textbf{x}\sb{j} +$
$ \textbf{x}\sb{j}(-i\hbar\delta\sb{ij}) $
$= -2i\hbar \delta\sb{ij}\textbf{x}\sb{j}$

$=-i\hbar \partial \textbf{x}\sb{j}\sp{2}/\partial x\sb{i}$

Supongamos que sea cierto para k:

$[\textbf{p}\sb{i},\textbf{x}\sb{j}\sp{k}]=$
$-i\hbar \partial \textbf{x}\sb{j}\sp{k}/\partial x\sb{i} =$
$-i\hbar k \delta\sb{ij}\textbf{x}\sb{j}\sp{k-1}$

Demostremos que también lo es para k+1:

$[\textbf{p}\sb{i},\textbf{x}\sb{j}\sp{k+1}]$
$=[\textbf{p}\sb{i},\textbf{x}\sb{j}\textbf{x}\sb{j}\sp{k}]$
$=[\textbf{p}\sb{i},\textbf{x}\sb{j}]\textbf{x}\sb{j}\sp{k}$
$+ \textbf{x}\sb{j}[\textbf{p}\sb{i},\textbf{x}\sb{j}\sp{k}]$
$=-i\hbar \delta\sb{ij}\textbf{x}\sb{j}\sp{k}$
$ + \textbf{x}\sb{j}(-i\hbar k \delta\sb{ij}\textbf{x}\sb{j}\sp{k-1})$

$=-i\hbar (k+1) \delta\sb{ij}\textbf{x}\sb{j}\sp{k}$

$=-i\hbar \partial \textbf{x}\sb{j}\sp{k+1}/\partial x\sb{i}$

Necesitamos un detalle más, y tomemos $l \not= i$:

$[\textbf{p}\sb{i},\textbf{x}\sb{j}\sp{k}\textbf{x}\sb{l}\sp{m}]$
$=[\textbf{p}\sb{i},\textbf{x}\sb{j}\sp{k}]\textbf{x}\sb{l}\sp{m}$ $+\textbf{x}\sb{j}\sp{k}[\textbf{p}\sb{i},\textbf{x}\sb{l}\sp{m}]$

$=[\textbf{p}\sb{i},\textbf{x}\sb{j}\sp{k}]\textbf{x}\sb{l}\sp{m}+0$

$=(-i\hbar \partial \textbf{x}\sb{j}\sp{k}/\partial x\sb{i})\textbf{x}\sb{l}\sp{m} $

Eso quiere decir que el operador momento en la dirección i deriva solo con respecto a i y las demás variables las ve como constantes. Podemos ver inmediatamente la estructura general del punto 5 del teorema:

\

Si $G$ es desarrollable en serie (se puede aproximar por polinomios de tal forma que la precisión puede aumentarse tanto como se desea), tenemos entonces:

$G(x,y,z)=a\sb{0} + a\sb{1}x +a\sb{2}y +a\sb{3}z $
$+a\sb{4}xy +a\sb{5}xz +a\sb{6}yz +a\sb{7}x\sp{2} $
$+a\sb{8}y\sp{2} +a\sb{9}z\sp{2} +..$

Por tanto, su derivada, digamos con respecto a y, es:

$\partial G/\partial y = a\sb{2} +a\sb{4}x +a\sb{6}z +2a\sb{8}y  +..$

Por otro lado:

$[p\sb{y},G(x,y,z)]=[p\sb{y},a\sb{0} + a\sb{1}x $
$+a\sb{2}y +a\sb{3}z +a\sb{4}xy +a\sb{5}xz $
$+a\sb{6}yz +a\sb{7}x\sp{2}+.. ] $

Aplicando la distributividad del conmutador, calculemos cada sumando resultante. Hay que recordar que  dos operadores son iguales si hacen lo mismo sobre una función inespecífica $\psi$. Por ejemplo, la proposición que dice

$[p\sb{y},a\sb{2}y]=-i\hbar a\sb{2}$

se demuestra operando sobre $\psi$, como sigue:

\bigskip

$[p\sb{y},a\sb{2}y]\psi = (h/i) \partial (a_2y\psi)/\partial y - a_2 y (h/i)\partial \psi/\partial y $

$=(h/i) a_2\psi + (h/i) a_2y\partial \psi/\partial y -a_2 y (h/i)\partial \psi/\partial y  = (h/i) a_2 y=   -i\hbar a\sb{2}y$

\bigskip

Los cálculos dan:

$[p\sb{y},a\sb{0}]= 0$

$[p\sb{y},a\sb{1}x]=0 $

$[p\sb{y},a\sb{2}y]=-i\hbar a\sb{2}$

$[p\sb{y},a\sb{3}z ]=0$

$[p\sb{y},a\sb{4}xy]=-i\hbar a\sb{4} x$

$[p\sb{y},a\sb{5}xz]=0$

$[p\sb{y},a\sb{6}yz]=-i\hbar a\sb{6} z$

$[p\sb{y},a\sb{7}x\sp{2}]=0$

$[p\sb{y},a\sb{8}y\sp{2}]= -2i\hbar a\sb{8} y$

$[p\sb{y},a\sb{9}z\sp{2}]=0$
...

Sumando todas estas igualdades a lado y lado nos da:

$[p\sb{y},G(x,y,z)]=-i\hbar a\sb{2}-i\hbar a\sb{4} x$
$-i\hbar a\sb{6} z-2i\hbar a\sb{8} y...$

Por otra parte:

$\partial G/\partial y = a\sb{2} +a\sb{4}x +a\sb{6}z +2a\sb{8}y  +..$

Inferimos  entonces que:

$[p\sb{y},G(x,y,z)]=-i\hbar \partial G/\partial y$

\

Qué le falta a nuestra inferencia para que se convierta en una verdadera demostración? Escribir la  fórmula general del polinomio de Taylor y sobre ella aplicar lo que hemos aprendido. De pronto vendría bien una explicación de por qué y cuando  la derivada de la serie de una función es la serie de la derivada de la función.

\

Demostración de la propiedad 5:

5.$[\textbf{p}\sp{2},G(\textbf{x})]= -[G(\textbf{x}),\textbf{p}\sp{2}]$
$=-[G(\textbf{x}),\textbf{p}\textbf{p}] $
$=-[G(\textbf{x}),\textbf{p}]\textbf{p}-\textbf{p}[G(\textbf{x}),\textbf{p}] ] $

$=[\textbf{p},G(\textbf{x})]\textbf{p}$
$+\textbf{p}[\textbf{p},G(\textbf{x})] $

$=-i\hbar (\nabla G)\textbf{p} -i\hbar \textbf{p}(\nabla G)$

Procuremos en los siguientes teoremas no confundir el número $e$ con la carga del electrón, denotada por la misma letra:

\

\addtocounter{ecu}{1}
\textit{Teorema \theecu}
$\nabla e\sp{ie\Lambda/\hbar c}= ie (\nabla \Lambda)/\hbar c)$ $e\sp{ie\Lambda/\hbar c}$

\

Esto se debe a que el gradiente produce un vector tal que en cada coordenada va una derivada parcial, la cual opera sobre una exponencial. La regla de la cadena da: la derivada del exponente por la misma exponencial.

\

\addtocounter{ecu}{1}
\textit{Teorema \theecu}

Sea $\Lambda = \Lambda(\textbf{x})$, entonces:

$e\sp{-ie\Lambda/\hbar c}\textbf{p}e\sp{ie\Lambda/\hbar c}$
$= \textbf{p} + e\nabla \Lambda/c$

\

Demostración:

$e\sp{-ie\Lambda/\hbar c}\textbf{p}e\sp{ie\Lambda/\hbar c}$
$=e\sp{-ie\Lambda/\hbar c}(\textbf{p}e\sp{ie\Lambda/\hbar c}$
$-e\sp{ie\Lambda/\hbar c}\textbf{p} $
$+ e\sp{ie\Lambda/\hbar c}\textbf{p})$

$=e\sp{-ie\Lambda/\hbar c}([\textbf{p},e\sp{ie\Lambda/\hbar c}]$
$ + e\sp{ie\Lambda/\hbar c}\textbf{p})$

$=e\sp{-ie\Lambda/\hbar c}(-i\hbar \nabla e\sp{ie\Lambda/\hbar c})$
$ + e\sp{-ie\Lambda/\hbar c}e\sp{ie\Lambda/\hbar c}\textbf{p}$

$=e\sp{-ie\Lambda/\hbar c}(-i\hbar \nabla e\sp{ie\Lambda/\hbar c})$
$ + \textbf{p}$

$=e\sp{-ie\Lambda/\hbar c}(-i\hbar (ie (\nabla \Lambda)/\hbar c)$ $e\sp{ie\Lambda/\hbar c})$
$ + \textbf{p}$

$=e\sp{-ie\Lambda/\hbar c}(e/ c) (\nabla \Lambda) $
$e\sp{ie\Lambda/\hbar c}$
$ + \textbf{p}$

$=(e/ c) e\sp{-ie\Lambda/\hbar c} $
$e\sp{ie\Lambda/\hbar c}(\nabla \Lambda)$
$ + \textbf{p}$

$=(e/ c) (\nabla \Lambda)$
$ + \textbf{p}$

$= \textbf{p} + (e/c) \nabla \Lambda$

\

Ahora procedemos  con $\textbf{p}\sp{2}$:

\

\addtocounter{ecu}{1}
\textit{Teorema \theecu}

Sea $\Lambda = \Lambda(\textbf{x})$, entonces:

$e\sp{-ie\Lambda/\hbar c}\textbf{p}\sp{2} e\sp{ie\Lambda/\hbar c}$
$= (\textbf{p} + (e/c)\nabla \Lambda)\sp{2}$

\

Demostración:

$e\sp{-ie\Lambda/\hbar c}\textbf{p}\sp{2} e\sp{ie\Lambda/\hbar c}$
$=e\sp{-ie\Lambda/\hbar c}(\textbf{p}\sp{2} e\sp{ie\Lambda/\hbar c}$
$-e\sp{ie\Lambda/\hbar c}\textbf{p}\sp{2} $
$ + e\sp{ie\Lambda/\hbar c}\textbf{p}\sp{2} )$

$=e\sp{-ie\Lambda/\hbar c}([\textbf{p}\sp{2} ,e\sp{ie\Lambda/\hbar c}]$
$ + e\sp{ie\Lambda/\hbar c}\textbf{p}\sp{2} )$

$=e\sp{-ie\Lambda/\hbar c}[\textbf{p}\sp{2} ,e\sp{ie\Lambda/\hbar c}]$
$ + e\sp{-ie\Lambda/\hbar c}e\sp{ie\Lambda/\hbar c}\textbf{p}\sp{2} $

$=e\sp{-ie\Lambda/\hbar c}[\textbf{p}\sp{2} ,e\sp{ie\Lambda/\hbar c}]$
$ + \textbf{p}\sp{2} $

(Ahora aplicamos $[A,B^2] =[A,B]B + B[A,B]$ y la propiedad 4 )

$=e\sp{-ie\Lambda/\hbar c}\{-i\hbar (\nabla e\sp{ie\Lambda/\hbar c})\textbf{p} $
$+ \textbf{p}(-i\hbar (\nabla e\sp{ie\Lambda/\hbar c})    )\}$
$ + \textbf{p}\sp{2} $ (por la propiedad 5)

$=e\sp{-ie\Lambda/\hbar c}\{-i\hbar (ie\nabla \Lambda/\hbar c) e\sp{ie\Lambda/\hbar c}\textbf{p} $
$+ \textbf{p}(-i\hbar (ie\nabla \Lambda /\hbar c) $
$e\sp{ie\Lambda/\hbar c})    \}$
$ + \textbf{p}\sp{2} $

$=e\sp{-ie\Lambda/\hbar c}\{ e\sp{ie\Lambda/\hbar c}(e\nabla \Lambda/ c) \textbf{p} $
$+ \textbf{p}( e\sp{ie\Lambda/ \hbar c}(e\nabla \Lambda / c))    \}$
$ + \textbf{p}\sp{2} $

$= (e/ c)\nabla \Lambda \textbf{p} $
$+ (e/c)e\sp{-ie\Lambda/\hbar c}\textbf{p} e\sp{ie\Lambda/\hbar c}\nabla \Lambda     $
$ + \textbf{p}\sp{2} $

$=(e/ c)\nabla \Lambda \textbf{p} $
$+ (e/c)((e/c)\nabla \Lambda + \textbf{p})\nabla \Lambda    $
$ + \textbf{p}\sp{2} $

$=(e/ c)\nabla \Lambda \textbf{p} $
$+ (e/c)\sp{2}(\nabla \Lambda)\sp{2} $
$+ (e/c)\textbf{p}\nabla \Lambda    $
$ + \textbf{p}\sp{2} $

$=\textbf{p}\sp{2} + (e/c)\textbf{p}\nabla \Lambda $
$+(e/ c)\nabla \Lambda \textbf{p} $
$+ (e/c)\sp{2}(\nabla \Lambda)\sp{2}      $

$=(\textbf{p} +(e/c)\nabla \Lambda )\sp{2}$

\

Analicemos una variante más con una translación local:

\

\addtocounter{ecu}{1}
\textit{Teorema \theecu}

Sea $\textbf{G} = \textbf{G\textbf(x)}$ una función vectorial  de la posición. Entonces se cumple que:

$e\sp{-ie\Lambda/\hbar c}(\textbf{p}$
$- \textbf{G} )\sp{2} e\sp{ie\Lambda/\hbar c}$
$= (\textbf{p} + (e/c)\nabla \Lambda- \textbf{G})\sp{2}$

\

Demostración:

$e\sp{-ie\Lambda/\hbar c}(\textbf{p}$
$- \textbf{G} )\sp{2} e\sp{ie\Lambda/\hbar c}$

$=e\sp{-ie\Lambda/\hbar c}(\textbf{p}\sp{2}$
$- \textbf{p}\textbf{G} -\textbf{G}\textbf{p}$
$+ \textbf{G}\sp{2}) e\sp{ie\Lambda/\hbar c}$

$=e\sp{-ie\Lambda/\hbar c} \textbf{p}\sp{2}e\sp{ie\Lambda/\hbar c}$
$- e\sp{-ie\Lambda/\hbar c}\textbf{p}\textbf{G} e\sp{ie\Lambda/\hbar c}$
$-e\sp{-ie\Lambda/\hbar c}\textbf{G}\textbf{p}e\sp{ie\Lambda/\hbar c}$

$+ e\sp{-ie\Lambda/\hbar c}\textbf{G}\sp{2}e\sp{ie\Lambda/\hbar c} $

$=e\sp{-ie\Lambda/\hbar c} \textbf{p}\sp{2}e\sp{ie\Lambda/\hbar c}$
$- e\sp{-ie\Lambda/\hbar c}\textbf{p} e\sp{ie\Lambda/\hbar c}\textbf{G}$
$-\textbf{G} e\sp{-ie\Lambda/\hbar c}\textbf{p}e\sp{ie\Lambda/\hbar c}$

$+ e\sp{-ie\Lambda/\hbar c}\textbf{G}\sp{2}e\sp{ie\Lambda/\hbar c} $

$=(\textbf{p} +(e/c)\nabla \Lambda )\sp{2}$
$-(\textbf{p} +(e/c)\nabla \Lambda )\textbf{G}$
$-\textbf{G}(\textbf{p} +(e/c)\nabla \Lambda )$
$+\textbf{G}\sp{2}$

$= (\textbf{p} + (e/c)\nabla \Lambda- \textbf{G})\sp{2}$

\

Aplicando el anterior teorema (\theecu) cuando
$\textbf{G}= (e/c)(\textbf{A}+\nabla \Lambda)$ tenemos:

\

\addtocounter{ecu}{1}
\textit{Teorema \theecu}

Definiendo

$\textbf{H}=(\textbf{p} - (e/c)\textbf{A})\sp{2}  $

y

$\textbf{H}\sb{\Lambda}=(\textbf{p}  $
$- (e/c)(\textbf{A}+ \nabla \Lambda))\sp{2}$

Se cumple que:

$e\sp{-ie\Lambda /\hbar c} \textbf{H}\sb{\Lambda}e\sp{ie\Lambda/\hbar c}$
$=\textbf{H}$

\

Demostración:

$e\sp{-ie\Lambda/\hbar c} \textbf{H}\sb{\Lambda}e\sp{ie\Lambda/\hbar c}$

$= e\sp{-ie\Lambda/\hbar c}(\textbf{p}$
$- (e/c)(\textbf{A}+\nabla \Lambda) )\sp{2} e\sp{ie\Lambda/\hbar c}$
$= (\textbf{p} + (e/c)\nabla \Lambda$
$- (e/c)(\textbf{A}+\nabla \Lambda))\sp{2}$
$=(\textbf{p} - (e/c)\textbf{A})\sp{2}  $
$=\textbf{H}$

\

Toda la maquinaria que hemos desarrollado hasta ahora nos permite demostrar muy cómodamente que el electromagnetismo es una teoría gauge con grupo de invariancia local $\textbf{U(1)}$. En efecto:

Consideremos el Hamiltoniano

$\textbf{H}=(\textbf{p} - (e/c)\textbf{A})\sp{2}  $

Su ecuación de Schr\"{o}edinger correspondiente es:

$ \textbf{H}\psi =(\textbf{p}  - (e/c)\textbf{A})\sp{2}\psi $
$= (-\hbar/i)\partial \psi/\partial t $

Puesto que $\textbf{A}$ admite una arbitrariedad, pues está definida módulo un gradiente cualesquiera de una función sólo de posición, debemos comparar la física producida por $\textbf{H}$ con la producida por:

$\textbf{H}\sb{\Lambda}=(\textbf{p}  $
$- (e/c)(\textbf{A}+ \nabla \Lambda))\sp{2}$

Su ecuación de Schr\"{o}edinger es:

$ \textbf{H}\sb{\Lambda}\psi\sb{\Lambda} =(\textbf{p} - $
$(e/c)(\textbf{A}+ \nabla \Lambda))\sp{2}\psi \sb{\Lambda}$
$= (-\hbar/i)\partial \psi\sb{\Lambda}/\partial t $

\

Tomemos la primera ecuación y  utilicemos los teoremas anteriores pero al revés:

$ \textbf{H}\psi =(\textbf{p}  - (e/c)\textbf{A})\sp{2}\psi $
$= (-\hbar/i)\partial \psi/\partial t $

$e\sp{-ie\Lambda/\hbar c}(\textbf{p}- $
$(e/c)(\textbf{A}+\nabla \Lambda) )\sp{2} e\sp{ie\Lambda/\hbar c}\psi $
$= (-\hbar/i)\partial \psi/\partial t $

O sea:

$e\sp{-ie\Lambda/\hbar c} \textbf{H}\sb{\Lambda} $
$e\sp{i\Lambda/\hbar c}\psi = (-\hbar/i)\partial \psi/\partial t $

Si multiplicamos por $e\sp{i\Lambda/\hbar c}$ en ambos lados:

$\textbf{H}\sb{\Lambda} e\sp{ie\Lambda/\hbar c}\psi $
$= (-\hbar/i)e\sp{ie\Lambda/\hbar c} \partial \psi/\partial t $

Finalmente obtenemos:

$\textbf{H}\sb{\Lambda} (e\sp{ie\Lambda/\hbar c}\psi) $
$= (-\hbar/i) \partial (e\sp{ie\Lambda/\hbar c}\psi)/\partial t $

\

Podemos entonces resumir:

\

\addtocounter{ecu}{1}
\textit{Teorema \theecu}

Teorema:

$\psi $ es solución de

$ \textbf{H}\psi =(\textbf{p}  - (e/c)\textbf{A})\sp{2}\psi $
$= (-\hbar/i)\partial \psi/\partial t $

ssi

$\psi\sb{\Lambda}=e\sp{ie\Lambda/\hbar c}\psi  $ es solución de

$ \textbf{H}\sb{\Lambda}\psi\sb{\Lambda} =$
$(\textbf{p} - (e/c)(\textbf{A}+ $
$\nabla \Lambda))\sp{2}\psi \sb{\Lambda}= $
$(-\hbar/i)\partial \psi\sb{\Lambda})/\partial t $

\

Verbalizamos este teorema diciendo que la arbitrariedad en la definición del vector potencial $\textbf{A}$ módulo el gradiente de una función de posición $\Lambda (\textbf{x})$ es equivalente a un cambio local de fase en la función de onda. La palabra 'local' significa que en cada punto del espacio se puede tomar una fase diferente, puesto que la función $\Lambda$ depende del lugar en donde se esté.

\

Anteriormente vimos que un cambio global de fase no cambia la física en absoluto. Cuál será el efecto de un cambio local de fase? Tenemos que investigar los valores propios del Hamiltoniano y además las probabilidades correspondientes a cada medición.

\

\addtocounter{ecu}{1}
\textit{Teorema \theecu. Un cambio local de fase es una isometría, es decir, no crea ni destruye probabilidades.}

\

Demostración :

Un cambio local de fase es equivalente a multiplicar por una exponencial imaginaria. Esto viene del hecho de que todo número complejo tiene un logaritmo (que no es una función, por lo que hay que tomar una rama apropiada) y por tanto, toda función de onda puede reescribirse como la exponencial de una cierta función compleja.

Cambio de fase global: $\psi \rightarrow e^{\phi_o}\psi$

Cambio de fase local: $\psi \rightarrow e^{\phi(x)}\psi = e\sp{ie\Lambda (x)/\hbar c}\psi$

Una probabilidad importante se da por el producto interno entre dos funciones  de onda $ \psi\sb{1}$, $\psi\sb{2}$. Antes del cambio de fase, la probabilidad es simplemente $<\psi\sb{1},\psi\sb{2}>$. Con el cambio de fase local, el mismo para ambas funciones, las funciones se cambian en $e\sp{ie\Lambda (x)/\hbar c}\psi_1$ y $e\sp{ie\Lambda(x) /\hbar c}\psi\sb{2}$ respectivamente: Por lo tanto, el producto interno se convierte en:

$<e\sp{ie\Lambda (x)/\hbar c}\psi_1, e\sp{ie\Lambda(x) /\hbar c}\psi\sb{2}>=$

$=\int e\sp{ie\Lambda (x)/\hbar c}\psi_1 ( e\sp{ie\Lambda(x) /\hbar c}\psi\sb{2})\sp{*}$

$=\int e\sp{ie\Lambda (x)/\hbar c}\psi_1  e\sp{-ie\Lambda(x) /\hbar c}\psi\sb{2}\sp{*}$

$=\int \psi\sb{1} \psi\sb{2}\sp{*}$

$=<\psi\sb{1},\psi\sb{2}>$

\

Por lo tanto, un cambio local de fase no crea ni destruye probabilidades. Revisemos qué pasa con los valores propios:

\

\addtocounter{ecu}{1}
\textit{Teorema \theecu}

Teorema:

 $(\psi, \lambda )$ es un par propio de $\textbf{H}$ ssi

$(e\sp{ie\Lambda(x)/\hbar c}\psi, \lambda)$ es un par propio de $\textbf{H}\sb{\Lambda} \textbf{H}\sb{\Lambda} e\sp{ie\Lambda/\hbar c}$

\

Demostración: Nos basamos en un resultado anterior:

Si

$\textbf{H}=(\textbf{p} - (e/c)\textbf{A})\sp{2}  $

y

$\textbf{H}\sb{\Lambda}=(\textbf{p}  $
$- (e/c)(\textbf{A}+ \nabla \Lambda))\sp{2}$

Se cumple que:

$e\sp{-ie\Lambda /\hbar c} \textbf{H}\sb{\Lambda}e\sp{ie\Lambda/\hbar c}$
$=\textbf{H}$

\

Utilizando este resultado tenemos:

 $(\psi, \lambda )$ es un par propio de $\textbf{H}$ ssi

$\textbf{H}\psi = \lambda \psi$ ssi

$e\sp{-ie\Lambda/\hbar c}\textbf{H}\sb{\Lambda} $
$e\sp{ie\Lambda/\hbar c}\psi = \lambda \psi$

ssi

$\textbf{H}\sb{\Lambda} e\sp{ie\Lambda/\hbar c}\psi $
$= \lambda e\sp{ie\Lambda/\hbar c}\psi$

ssi

$\textbf{H}\sb{\Lambda} (e\sp{ie\Lambda/\hbar c}\psi) $
$= \lambda (e\sp{ie\Lambda/\hbar c}\psi)$

Lo cual dice que si    $(\psi, \lambda )$ es un par propio de $\textbf{H}$, entonces  $(\lambda , e\sp{ie\Lambda/\hbar c}\psi)$ es un par propio de otro operador Hamiltoniano que conlleva la misma física, pues todo lo que hace es involucrar un gradiente en la definición del potencial vector. En este caso, la física está dada por los saltos de energía, dados por los valores propios, y por las probabilidades de transición, dadas por los productos interiores.

Por lo tanto, los valores propios se conservan aunque haya un cambio local de fase. Pero como el Hamiltoniano es el operador de evolución, un cambio local de fase podría modificar las mediciones de otra observables. Será ese el caso?

\

\addtocounter{ecu}{1}
\textit{Teorema \theecu. Un \textbf{cambio local de fase} \index{cambio local de fase} no cambia las probabilidades de ninguna medición. }

\

Demostración:

Sea $\textbf{B}$ un operador autoadjunto que represente una observable determinada, con pares propios $(b\sb{i}, \phi\sb{i})$.   Entonces, la probabilidad $p\sb{i}$ de que  se obtenga como medición de $\textbf{B}$ el valor $b\sb{i}$ cuando el sistema está en el estado $\psi$ es:

$p\sb{i}=|<\phi\sb{i},\psi >|\sp{2}$

Supongamos entonces que el sistema se prepara en la condición inicial $\phi\sb{0}$. Si lo sometemos a evolución según $\textbf{H}$ y hacemos una medición de $\textbf{B}$ en el tiempo $t$, la probabilidad  $p\sb{t,i}$ de que  se obtenga como medición de $\textbf{B}$ el valor $b\sb{i}$ es

$p\sb{t,i}=|<\phi\sb{i},e\sp{-it\textbf{H}/\hbar} \psi\sb{0} >|\sp{2}$

Recordando que $e\sp{-it\textbf{H}/\hbar} \psi\sb{0} =\psi(\textbf{x},t)$,

$p\sb{t,i}=|<\phi\sb{i},e\sp{-it\textbf{H}/\hbar} \psi\sb{0} >|\sp{2}$
$=|<\phi\sb{i},\psi(\textbf{x},t)>|\sp{2}$

Ahora bien, si el sistema evoluciona bajo la acción de $\textbf{H}\sb{\Lambda}$, entonces en el tiempo $t$ el sistema se encontrará en el estado  $(e\sp{ie\Lambda/\hbar c}\psi(x,t))$. Observemos ahora por qué no nos funcionan las cosas:

Si el sistema en el tiempo $t$ está en el estado $\phi\sb{i}$, las probabilidades asociadas a las mediciones serían:

$p\sb{t,i,\Lambda}= $
$|<\phi\sb{i},e\sp{ie\Lambda/\hbar c} \psi(\textbf{x},t)>|\sp{2}$

$=|\int \phi\sb{i} e\sp{ie\Lambda/\hbar c} \psi(\textbf{x},t)dV |\sp{2}$

\

No podemos sacar la fase, puesto que ella ya no es una constante, sino que ahora varía, puesto que el cambio de fase es local o sea que $\Lambda$ es función de la posición $\vec{x}$. Para que las cosas funciones claramente, revisemos el siguiente detalle: 

\

Hacer un cambio de fase, equivale rotar los ejes del plano complejo. Si hacemos eso, cómo podremos modificar una función de onda y no otra en el mismo sitio? Podemos modificar la fase en sitios diferentes, pero no podemos modificar en el mismo sitio la fase de funciones de onda diferentes, puesto que las funciones de onda no son observables. Si lo fuesen, uno podría aspirar a separar una función de onda de otra y modificarlas a voluntad. Pero no lo son y lo que le pase a una en un  lugar, necesariamente ha de pasarle a otra cualquiera en ese mismo lugar.

Por lo tanto, no podemos escribir

$p\sb{t,i,\Lambda}= $
$|<\phi\sb{i},e\sp{ie\Lambda/\hbar c} \psi(\textbf{x},t)>|\sp{2}$

Debemos transformar a  $\phi$ en $e\sp{ie\Lambda/\hbar c}\phi$ (esta precaución no fue tomada anteriormente cuando estudiábamos el cambio global de fase: por favor, corrija la demostración adecuadamente) y con esa aclaración las probabilidades asociadas a las mediciones de $\textbf{B}$ son:

$p\sb{t,i,\Lambda}= |<e\sp{ie\Lambda/\hbar c}$ $\phi\sb{i},e\sp{ie\Lambda/\hbar c} \psi(\textbf{x},t)>|\sp{2}$

Pero como un cambio local de fase que afecte a todas las funciones por igual es una isometría, entonces podemos cerrar la discusión:

$p\sb{t,i,\Lambda}= |<\phi\sb{i},\psi(\textbf{x},t)>|\sp{2}$
$=p\sb{t,i}$

Hemos demostrado  que un cambio de fase que sea local  no afecta la física para nada siempre y cuando no discrimine entre funciones de onda en el mismo punto.

\

\

\color{red}

Nuestro resumen es:

  En mecánica cuántica, siempre podemos agregar al potencial vector el gradiente de una función escalar en el Hamiltoniano que representa el campo electromagnético. Eso es equivalente a  un  \index{cambio local de fase} \textbf{cambio local de fase}, el cual debe ser universal para todas las funciones de onda.

  En otras palabras, la \index{invariancia gauge} \textbf{invariancia gauge} del electromagnetismo clásico, de permitir una arbitrariedad en el potencial vector módulo el gradiente de una función escalar, se transformó en mecánica cuántica en la invariancia gauge de permitir una arbitrariedad en la definición de la fase módulo una función escalar que puede variar de punto a punto.

Hacer un cambio local de fase es lo mismo que multiplicar localmente por un elemento de \textbf{U(1)}. Por tanto, el resumen de todo es:

\

La interacción electromagnética es una teoría Gauge con grupo de invariancia local \textbf{U(1)}.

\

\color{black}

\section{LA LUZ Y  LA GUERRA}

Un registro antiguo cuenta así: "Y se corrompió la tierra delante de Dios, y estaba la tierra llena de violencia". Podemos ver que eso sigue siendo cierto hasta el día de hoy. La primera bomba atómica fue construida por un grupo de genios que podían expresarse con tal claridad, precisión y sabiduría que no necesitaban repetirse las cosas dos veces para que éstas quedaran en claro. Uno de ellos era \index{Feynman} \textbf{ Feynman}, en ese tiempo, apenas un jovencito tamborilero. 

Por esta época, las bombas de destrucción masiva no están de moda. Ahora lo que se quiere es poder para hacer operaciones estratégicas de alta efectividad. Como punta de lanza se trabaja en el perfeccionamiento del programa de posicionamiento global, gracias al cual se espera poder hacer operaciones quirúrgicas de alta precisión en cualquier parte del planeta.

Surge entonces la necesidad de poder medir con exactitud. La medida fundamental es la del tiempo. Si uno puede medir el tiempo, uno puede medir la distancia utilizando la velocidad de la luz como referencia.

Para medir el tiempo se cuenta los periodos de oscilación de un rayo de luz. La luz sale de átomos excitados que regresan a un punto de potencial más bajo. La frecuencia de dicha luz, vista desde un sistema de coordenadas que viaja con el átomo, es siempre la misma. Pero los átomos se mueven y son como sirenas que pasan: cuando se acercan agudizan el sonido y cuando se alejan lo agravan. Se denomina efecto Doppler a este fenómeno. Debido a que los átomos reales siempre se mueven, incluso casi en el cero absoluto, por causa del efecto Doppler no se puede aún medir el tiempo con precisión absoluta.

Se sigue buscando: se ha logrado ya hacer que la luz coherente (laser) frene partículas casi al cero absoluto. Pero la precisión aún no es suficiente para que los gerentes de la guerra estén tranquilos.

Se sigue trabajando y es un gran orgullo para cualquiera poder  mostrar algo nuevo al respecto.

\section{EL TEOREMA DE NOETHER}

Vimos que existe una versión en mecánica clásica del teorema de Noether. Podemos inferir  que en mecánica cuántica debe ser lo mismo. La razón   es que la mecánica clásica es un caso límite de la cuántica y que además hay un método automático para pasar de la mecánica clásica a la cuántica. En esta sección formalizaremos esa idea.

Lo primero que debemos hacer es preguntarnos cómo podremos codificar en mecánica cuántica la proposición de que la física de un sistema sea invariante ante la acción de un grupo.

En mecánica cuántica la física está ligada a la observación, la cual es predicha, de un parte, por los valores propios de los operadores que representan las observables y, por otra,  por las probabilidades de transición. Los valores propios dan las mediciones posibles que la observable puede arrojar y las probabilidades de transición permiten predecir la dinámica observacional del sistema, pues la dinámica en sí está descrita por el Hamiltoniano que es el operador de la energía que es el operador de evolución.

Consideremos entonces un sistema con Hamiltoniano $H$ que no depende del tiempo. La probabilidad de que un sistema que se encuentre en el estado $\psi $ sea registrado en el estado $\phi$ está dada por

$\|<\phi,\psi>\|^2$

Consideremos la acción de un grupo. Para fijar ideas,  digamos que se trata del grupo de las  rotaciones, y consideremos que la acción del grupo no afecta la física: poco nos importa que el sistema se estudie en el centro del salón o al final del pasillo, o si lo orientamos hacia la ventana o hacia la estrella polar (para experimentos eléctricos esto no es perfectamente cierto: por las paredes y el piso van cables que transportan cargas en movimiento que crean campos magnéticos cuyo efecto no ha de ser despreciable en todas las direcciones).

El sistema observado antes de la rotación está descrito por la función de onda $\psi$, y después de la rotación por $\psi '$. Debe haber una manera biunívoca de pasar de una descripción a la otra, es decir, debe ser posible especificar un operador $U$ que transforme una descripción en la otra:

$ \psi \rightarrow \psi '$ = $U\psi$

Hemos escogido la letra $U$ para designar dicho operador teniendo en cuenta que ha de ser un operador unitario, es decir tiene que conservar las probabilidades de transición.

$\|<\phi,\psi>\|^2 = \|<\phi ',\psi '>\|^2$

$=\|<U\phi ',U\psi '>\|^2= \|<\phi ',U\sp{\dag} U\psi '>\|^2$

Por otro lado, estos operadores están directamente asociados a las rotaciones y, en cierto sentido, el conjunto de operadores es una copia del conjunto de rotaciones, de tal manera que las relaciones algebraicas se conservan: se dice que el grupo de operadores representa al grupo de las rotaciones.

Para poder decir que la física no cambia con una rotación, debemos exigir que las leyes de la naturaleza se vean exactamente de la misma forma esté o no el sistema rotado. Tenemos sola una ley, la ley de evolución la cual está dada por el Hamiltoniano:

$i\partial \psi/\partial t = H\psi$

$i\partial\psi '/\partial t = H\psi '$

además, las probabilidades de transición deben ser las mismas:

$<\phi ', H\psi'> = <\phi , H\psi> $

Pero teniendo en cuenta que hay una relación entre las dos funciones de onda, la cual está mediada por un operador unitario $U$:

$<\phi ', H\psi'> = <U \phi ,  HU\psi>$ , usando el adjunto de $U$ nos queda:

 $<U \phi ,  HU\psi> = <\phi , U \sp{\dag} H U\psi>$

por transitividad:

 $<\phi , H\psi>= <\phi ', H\psi'> = <U \phi ,  HU\psi> = <\phi , U \sp{\dag} H U\psi>$

y como eso es cierto para cualquier sistema físico, se deduce que

$H=U \sp{\dag} H U$

multiplicando por $U$ a ambos lados y teniendo en cuenta que

$U\sp{\dag}U = UU\sp{\dag} = I$ se tiene que

$UH=UU \sp{\dag} H U= HU$

lo que implica que:

$HU-UH=0= [U,H]$

\bigskip

Hemos demostrado entonces el siguiente

 \addtocounter{ecu}{1}

\textit{Teorema (\theecu) Decir que la física es invariante ante la acción de un grupo es lo mismo que decir, en mecánica cuántica, que la representación del grupo en el espacio de operadores unitarios forma un grupo con la misma estructura que el grupo dado y que cada uno de sus elementos conmuta con el Hamiltoniano.} Eso es equivalente a:

\bigskip

\addtocounter{ecu}{1}

\textit{Teorema (\theecu) El valor esperado de U es una constante de la dinámica.}

$i(d/dt)<\psi, U\psi(t)>= <\psi, UH-HU\psi>=0$.

Demostración:

Ignorando la constante de Planck

$i(d/dt)<\psi, U\psi(t)>= i<(\partial /\partial t)\psi, U\psi(t)> + i <\psi, (\partial /\partial ) U\psi(t)>$

=$<i(\partial /\partial t)\psi, U\psi(t)> +  <\psi, -i(\partial /\partial t )U\psi(t)>$

=$<H\psi, U\psi(t)> -  <\psi, HU\psi(t)>$

=$<\psi, HU\psi(t)> -  <\psi, HU\psi(t)>$

$ =<\psi, (HU-HU)\psi>=0$.

donde hemos utilizado el hecho de que el Hamiltoniano es autoadjunto: $H= H \sp{\dag} $. Como el Hamiltoniano conmuta con $U$, la derivada es cero, el valor esperado es constante, por lo que es una constante del movimiento.

\bigskip

Es eso el teorema de Noether?

No, el teorema de Noether es algo mucho más fino:

Para fijar ideas, sigamos pensando en las rotaciones. Ellas se pueden sumar y restar, es decir, componer hacia la derecha o hacia la izquierda de forma asociativa pero no conmutativa. Además, la rotación nula recobra la identidad, pues no hace nada y cada rotación tiene su inversa. Pero  resulta que toda rotación es una sucesión de rotaciones infinitesimales.

Para hablar de rotaciones en mecánica cuántica, necesitamos representar al grupo en el grupo de los operadores unitarios. Como estamos trabajando con una representación,  la rotación nula se representa con la identidad. Por otro lado, una rotación infinitesimal ha de representarse por un operador ligeramente diferente de la identidad. Consideremos un eje fijo, y sobre el pongamos una rotación infinitesimal, el operador correspondiente es de la forma:

$U= I-i\epsilon J$

donde $\epsilon $ es un número real que uno hace tender a cero para producir una rotación infinitesimal, en tanto que $J$ es el operador que representa la 'dirección' de la rotación izquierda determinada por el eje escogido. El número $i$ se eligió para poder asegurar que $J$ sea autoadjunto o sea que represente una observable:

\bigskip

\addtocounter{ecu}{1}

\textit{Teorema (\theecu) Si $U$ es unitario y $U= I-i\epsilon J$ entonces J es autoadjunto}

Demostración:

Si $U$ es unitario, entonces $I= U\sp{\dag}U $. Tenemos:

$ U\sp{\dag} = (I-i\epsilon J)\sp{\dag} =  (I+i\epsilon J\sp{\dag})$

por lo que

$I= U\sp{\dag}U = (I+i\epsilon J\sp{\dag})(I-i\epsilon J)$

 $= I +i\epsilon J\sp{\dag}   -i\epsilon J + \epsilon^2 J^2$

como $\epsilon $ es infinitesimal, su cuadrado es despreciable al tomar el límite:

$I = I +i\epsilon (J \sp{\dag} - J)$

simplificando:

$0= i\epsilon (J \sp{\dag} -J)$

es decir

$J \sp{\dag} = J$

Jerga: a $J$ se le llama un \index{generador del grupo} \textbf{generador del grupo}. El conjunto de generadores forma un espacio vectorial (ejercicio) que se llama el \textbf{álgebra de Lie} \index{álgebra de Lie} asociada al grupo de Lie dado. Obsérvese que, en particular, en vez de escribir $U= I-i\epsilon J$ hubiésemos podido definir $U= I+i\epsilon J$.

\bigskip

Ahora, hagamos  una generalización de algo que ya hemos visto:

\addtocounter{ecu}{1}

 \index{conservación de la carga}  \textbf{Teorema sobre la conservación de la carga} (\theecu).\textit{ Si $H$ es autoadjunto, es decir, si $H= H^{\dag}$, o lo que es lo mismo, si $<H\psi, \psi>=<\psi, H\psi>$, entonces la dinámica definida por $H$ conserva la norma de cada función de onda, lo cual implica que también se conserva la probabilidad total y por tanto la carga.}

\bigskip

Demostración:

$i\partial  \psi /\partial t = H\psi$ es la ecuación que da la dinámica. Multiplicando por $-i$:

$\partial  \psi /\partial t = -iH\psi$

Por otra parte y utilizando primas para denotar derivada temporal:

$(d/dt) || \psi ||^2 = $

$(d/dt)<\psi,\psi>= <\psi ',\psi>+<\psi,\psi'>= <-iH\psi,\psi>+<\psi,-iH\psi>$

$=-i<H\psi,\psi>+i<\psi,H\psi>$

(pero como $H$ es autoadjunto, $<H\psi,\psi>= <\psi,H\psi>$)

$=-i<\psi,H\psi>+i<\psi,H\psi>=0$

Esto da la conservación de la norma cuadrado.  Como la norma cuadrado se interpreta como la probabilidad total, dicha probabilidad se conserva. Multiplicando en todas partes por la carga del electrón o del protón, deducimos que  un operador autoadjunto no crea ni destruye partículas cargadas. Como las partículas masivas neutras las consideramos no elementales, inferimos que un operador autoadjunto no  crea ni destruye partículas en general.

Notemos ahora que el operador de evolución cambia funciones de onda en funciones de onda, eso genera una dinámica descrita por el grupo de operadores unitarios dados por

$U= e^{-iHt}$

desarrollando en serie, para $t$ muy pequeño,

$U=I-itH$

Notemos que U es unitario, $U^{\dag}U=UU^{\dag}=I$, y H es autoadjunto.

\bigskip

\addtocounter{ecu}{1}

\textbf{Teorema de Noether } \index{Teorema de Noether } (\theecu) \textit{Sea $G$ un grupo ante el cual la física es invariante. EL álgebra de Lie del grupo de representaciones de $G$ está formado por operadores autoadjuntos y que por consiguiente representan observables, cuyos valores propios son constantes en el tiempo, o sea son constantes de movimiento.}

\bigskip

Demostración:

Sea $U$ un elemento del grupo de invariancia infinitesimalmente diferente de la identidad: $U= I - i\epsilon  J$,y sea  $\lambda $ un valor propio de $J$, es decir $J\psi = \lambda \psi$. Nuestra tarea es demostrar que la derivada temporal de dicho valor propio es cero.

Recordando que el valor esperado de $U$ es constante en el tiempo nos queda:

$ i(d/dt)<\psi, U\psi> = 0$

pero $U= I - i\epsilon  J$, por lo cual:

$ i(d/dt)<\psi, (I - i\epsilon J)\psi> = 0$

$ i(d/dt)<\psi, \psi - i\epsilon \lambda \psi> =0 $

$ i(d/dt)<\psi, (1-i\epsilon \lambda) \psi> =0 $

$ i(d/dt)[(1+i\epsilon \lambda)<\psi,  \psi>] =0 $

puesto que $\lambda$ es real, pues $J$ es autoadjunto.

Como la norma cuadrado se conserva y es no nula

$ i<\psi,  \psi>(d/dt)[(1+i\epsilon \lambda)] =0 $

$ (d/dt)[(1+i\epsilon \lambda)] =0 $

$ i\epsilon(d \lambda/dt] =0 $

por lo cual la derivada de $\lambda $ es nula: los valores propios de los operadores del álgebra de Lie del grupo de invariancia se conservan.

\section{LA CONSERVACION LOCAL DE LA CARGA}

Nosotros ya vimos que toda evolución dada por un grupo de evolución de operadores unitarios, necesariamente conserva la carga eléctrica. Pero sucede que tal conservación es global, es decir, si en un momento dado la carga total del universo es cero, para un universo globalmente neutro, podemos garantizar que por la eternidad la situación seguirá siendo la misma.

Todos los resultados experimentales conocidos hasta la fecha predicen que la carga se conserva no sólo globalmente sino también localmente. Es decir, en cada rincón del universo donde haya carga nula, seguirá habiendo carga nula, si dicho rincón se mantiene aislado del resto del mundo. Eso no implica que el número de cargas positivas siga invariante ni tampoco que el número de cargas negativas siga invariante. Significa que si se crean 3 partículas positivas, entonces necesariamente se crearán 3 partículas negativas.

Para generalizar la conservación de la carga a cualquier rincón aislado del mundo se requiere   que los cambios globales de fase también se puedan generalizar a cambios locales. Veamos por qué dichos cambios locales implican la conservación local de la carga.

De acuerdo al teorema de Noether, si la física de un sistema es invariante ante un cambio global de fase, debe haber alguna cantidad que se conserva. Dicha cantidad es evidente: como la evolución cuántica es una isometría, la norma en $L\sb{2}$ de cada función de onda se conserva. Dicha norma es uno. Por lo tanto, ese uno siempre se conserva, de lo cual deducimos que la mecánica cuántica ni crea ni destruye partículas. El uno de la norma se multiplica por la carga del electrón y se concluye que en mecánica cuántica ni se crea ni se destruye carga eléctrica. Y lo que hemos dicho lo hemos argumentado para un cambio global de fase.

La norma resulta de evaluar una integral sobre todo el volumen de integración. Hemos asumido que se trata de todo el espacio $R\sp{3}$, lo cual tiene sentido para una onda que se extiende por todo el universo. Pero se considera que las partículas elementales, como creadoras de campo, son objetos puntuales. Por lo tanto, no es necesario tomar todo el universo como dominio de integración. Es suficiente tomar una esfera de radio infinitesimal centrada en la partícula y entonces el término 'global' se convierte en 'local'.

Quedamos entonces ante la necesidad de admitir que la conservación de la carga debe ser una ley que se cumple tanto global como localmente. Por eso es que hemos dado la pelea para garantizar que los cambios de fase también puedan ser locales, y que se pueda variar la fase de lugar a lugar. De esa forma varios experimentadores en diferentes lugares podrán comparar sus resultados de laboratorio sin necesidad de guardar protocolos complicadísimos de orientación en espacio curvo del uno con respecto al otro.

\section{EL EFECTO DE AHARANOV-BOHM}

Se acostumbra atribuir a  \index{Aharonov y Bohm} \textbf{Aharonov y Bohm} la gloria de imaginar en 1959 un experimento  que podría evidenciar un efecto electromagnético no clásico (diez a\~{n}os antes W. Ehrenberg y R. E. Siday ya habían propuesto la misma idea, pero como al turco del Principito, nadie les paró bolas). El experimento se ha llevado a cabo con éxito por otros investigadores, Chambers fue el primero en 1960,  A.Tonomura también lo repitió en 1982. Hacia 1998 se hizo uno muy refinado usando técnicas de estado sólido.

 El arreglo experimental incluye un hilito magnetizado, delgadito pero muy largo, que cuelga verticalmente. Dicho hilo tiene un campo magnético en su interior que tiene la misma dirección que el hilo, pero su magnetismo alcanza a influenciar el espacio que le es exterior tan sólo en los polos del hilo, que se les deja por allá bien lejos del experimento. Se hace pasar horizontalmente una corriente de electrones por el espacio que circunda el hilo desde un punto de origen 'a' hasta un punto de registro 'b', al otro lado del hilo. El hilo está blindado con una substancia aislante y los electrones no tocan ni el hilo  ni el campo magnético que hay en su interior.

\

La pregunta que hay que responder es la siguiente: De qué forma depende la corriente de electrones del flujo del campo magnético  dentro del hilo?

\

De acuerdo a la mecánica clásica, el grado de magnetización del hilo para nada debe influir la corriente de electrones, pues la fuerza de Lorentz depende del campo magnético y este para nada afecta a ningún electrón, pues el campo está dentro del hilo y los electrones viajan por fuera del hilo.

Estudiemos ahora la predicción cuántica.

Primero tendremos que calcular la amplitud de ir de 'a' a 'b', $K(a,b)$. Esta amplitud se puede dividir en dos grandes contribuciones, la debida a los caminos que pasan por un lado del hilo y las que pasan por el otro lado. Llamémoslas por arriba y por abajo, respectivamente:

$K(a,b)= K(a,b)\sb{arriba} + K(a,b)\sb{abajo}$

Nosotros vamos tras la física, es decir, tras efectos observables.  Por eso, calculamos la probabilidad $P$ de que una partícula llegue a su objetivo, el punto 'b':

$P=|K(a,b)|\sp{2} = |K(a,b)\sb{arriba} + K(a,b)\sb{abajo}|\sp{2}$

$=(K(a,b)\sb{arriba} + K(a,b)\sb{abajo})(K(a,b)\sb{arriba} + K(a,b)\sb{abajo})\sp{*}$

$=(K(a,b)\sb{arriba} + K(a,b)\sb{abajo})(K(a,b)\sb{arriba}\sp{*} + K(a,b)\sb{abajo}\sp{*})$

$=K(a,b)\sb{arriba}K(a,b)\sb{arriba}\sp{*} +K(a,b)\sb{arriba}\sp{*}K(a,b)\sb{abajo} $

$+ K(a,b)\sb{arriba}K(a,b)\sb{abajo}\sp{*}+K(a,b)\sb{abajo}K(a,b)\sb{abajo}\sp{*}$

$=|K(a,b)\sb{arriba}|\sp{2}+ |K(a,b)\sb{abajo}|\sp{2} $

$+K(a,b)\sb{arriba}\sp{*}K(a,b)\sb{abajo} +  K(a,b)\sb{arriba}  K(a,b)\sb{abajo}\sp{*}$

Un número complejo siempre se puede escribir como su módulo por la exponencial de su argumento:

$K(a,b)\sb{arriba}= a e\sp{i\phi}= a cos\phi + i asen\phi$ por lo que

$|K(a,b)\sb{arriba}|\sp{2} = a^2$

$K(a,b)\sb{abajo}= b e\sp{i\psi}= b cos\psi + i bsen\psi$ por lo que

$|K(a,b)\sb{abajo}|\sp{2} = b^2$

De tal forma que

$K(a,b)\sb{arriba}K(a,b)\sb{abajo}\sp{*}$
$= ab e\sp{i\phi} e\sp{-i\psi} $

$= ab e\sp{i(\phi -\psi)} = ab cos(\phi -\psi) - i absen(\phi -\psi)$

En tanto que

$K(a,b)\sb{arriba}\sp{*}K(a,b)\sb{abajo}$
$= ab e\sp{-i\phi} e\sp{i\psi} $
$=
 ab e\sp{i(-\phi +\psi)} $

$= ab cos(-\phi +\psi) - i absen(-\phi +\psi) = ab cos(\phi -\psi) + i absen(\phi -\psi) $

porque hemos utilizado la paridad del coseno  y la imparidad del seno. Sumando las dos expresiones buscadas nos queda:

$K(a,b)\sb{arriba}\sp{*}K(a,b)\sb{abajo} +  K(a,b)\sb{arriba}  K(a,b)\sb{abajo}\sp{*}$

$= 2ab cos(\phi -\psi)$ .

Por consiguiente

$P=a\sp{2} + b\sp{2} + 2ab cos(\phi -\psi)$

Tenemos una predicción para elaborar: la probabilidad de llegar a 'b' tiene un componente sinusoidal que depende de la diferencia de fase de las contribuciones por arriba y por abajo. Averigüemos dicha diferencia. Para eso nos sirve el formalismo Lagrangiano:

El Lagrangiano de la interacción electromagnética está dado por:

$L= (1/2)m \|\vec v \|\sp{2} - e \phi + (e/c) \vec v \cdot \vec A $

Como en nuestro experimento no hay campo eléctrico, el Lagrangiano se reduce a :

$L= (1/2)m \|\vec v \|\sp{2} + (e/c) \vec v \cdot \vec A $
$=L\sb{0} + (e/c) \vec v \cdot \vec A $

donde $L\sb{0}$ representa el Lagrangiano de una partícula libre, cuya acción correspondiente notamos como $S\sb{0}$.

El Lagrangiano representa el peaje por unidad de camino. La acción es el peaje total a lo largo de un recorrido. Si tenemos un camino $\gamma$ ya parametrizado, la acción es:

$S(a,b,\gamma)= \int\sb{\gamma} L$
$= \int\sb{\gamma}(L\sb{0} + (e/c) \vec v \cdot \vec A) $
$=S\sb{0} + \int\sb{\gamma}(e/c) \vec v \cdot \vec A $

Dicha acción clásica debida al camino $\gamma$ presenta un aporte cuántico a la amplitud total 'proporcional' a:

$e\sp{(i/\hbar) S(a,b,\gamma)}= $
$e\sp{(i/\hbar)[ S\sb{0}
+ \int\sb{\gamma}(e/c) \vec v \cdot \vec A]}$

De tal forma que la amplitud total $K(a,b)$ resulta de la integración de los aportes parciales de todos y cada uno de los caminos posibles:

$K(a,b)=\int \sb{a}\sp{b} $
$e\sp{(i/\hbar)[ S\sb{0} + \int\sb{\gamma}(e/c) \vec v \cdot \vec A]}$
$ D\gamma $

\

Dividiendo esta amplitud en sus dos grandes contribuciones obtenemos:

$K(a,b)= K(a,b)\sb{arriba} + K(a,b)\sb{abajo}$

$=\int \sb{a}\sp{b} e\sp{(i/\hbar) S\sb{0} } $
$e\sp{(i/\hbar)\int\sb{\gamma}(e/c) \vec v \cdot \vec A} $
$D\gamma \sb{arriba}$
$+  \int \sb{a}\sp{b} e\sp{(i/\hbar) S\sb{0} } $
$e\sp{(i/\hbar)\int\sb{\gamma}(e/c) \vec v \cdot \vec A} $
$D\gamma \sb{abajo}$

\

Ahora debemos cuadrar el experimento para que el efecto de la energía cinética sea conservativo, es decir que no exista fricción de ningún tipo. Una condición de alto vacío podría servir, pero no es indispensable: lo único que se necesita es que los electrones no interactúen con el medio. Por eso, el experimento puede diseñarse en estado sólido. En ese caso, para cada camino, el efecto de la acción sólo depende de los puntos extremos y esos son fijos para todos los caminos.

Por lo tanto  el aporte de la energía cinética  produce una fase que podemos denominar como  $E(\gamma)=  S\sb{0}(\gamma)/\hbar$. Esta es una función del camino y de su parametrización, de su velocidad: entre más alejado esté el camino, la velocidad ha de ser mayor y la fase oscilará más rápidamente.

$K(a,b)= $
$ \int \sb{a}\sp{b} $
$e\sp{iE(\gamma)}e\sp{(i/\hbar)\int\sb{\gamma}(e/c)\vec v \cdot \vec A}$
$ D\gamma \sb{arriba}$
$+   \int \sb{a}\sp{b}$ $e\sp{iE(\gamma)}e\sp{(i/\hbar)\int\sb{\gamma}(e/c) \vec v \cdot \vec A} $
$D\gamma \sb{abajo}$

\

Ahora hay que considerar el efecto del vector potencial sobre todos los caminos habidos y por haber que llevan de 'a' a 'b', pero que no pasan por el hilo.

\

Uno podría pensar que por fuera del hilo no  hay vector potencial pues no hay campo magnético. Bien, eso  no es cierto, y la apreciación correcta nos la da el teorema de Stokes:

$\oint\sb{curva} A = \int\sb{superficie} \nabla \times A $
$= \int\sb{superficie} B $

A la última integral se le llama flujo magnético. Por consiguiente, cada curva cerrada que encierre al hilo evidencia un vector potencial no nulo, no importa cuán distante del hilo se encuentre. Sin embargo, si la curva cerrada no encierra al hilo, entonces, el flujo magnético encerrado por la curva es cero y se concluye que  cada camino cerrado que no encierre al hilo no aporta nada a la integral total. Eso es equivalente a decir que la acción debida al potencial vector es la misma para dos caminos cualesquiera que tengan la misma relación con respecto al hilo. En efecto:

Sea $\alpha$ y $\beta$ caminos del mismo lado del hilo y que llevan ambos de 'a' a 'b'. Sea $-\beta$ la parametrización reversa de $\beta$. Tenemos que $\alpha$ concatenado con $-\beta$ produce un camino cerrado y la integral de línea da cero:

$\oint\sb{\alpha -\beta}  =0 =\int\sb{\alpha} +\int\sb{-\beta} = \int\sb{\alpha} -\int\sb{\beta}  $

Por consiguiente las dos integrales son iguales, el camino no importa y lo único importante son los puntos inicial y final de las trayectorias. Además:

$\int\sb{\gamma} \vec v \cdot \vec A $

$=\int\sb{t\sb{0}}\sp{t\sb{f}}  \vec v \cdot \vec A dt$

$=\int\sb{t\sb{0}}\sp{t\sb{f}} \vec A \cdot (d\vec s/dt) dt$

$=\int\sb{a}\sp{b} \vec A \cdot d\vec s$

Esta igualdad es válida dentro de la misma clase de caminos. Por clase se entiende clase de homotopía a la que pertenece el camino.

\

\

\color{blue}

$<<<<<<<<<<<<<<<<<<<<<<<<<<<<<<<<<<<<<<<<<$

\

Se hace camino al andar: \textbf{clases de Homotopía} \index{clases de Homotopía} alrededor de un hilo.

\

Considerando todos los caminos habidos y por haber que van desde un lugar a otro podemos definir, siguiendo el instinto, una relación de equivalencia como sigue: si un camino va por abajo del hilo y otro también, entonces son equivalentes, porque uno de ellos puede deformarse en el otro y tal deformación no presenta rompimientos, discontinuidades.

Un camino que va por encima no puede ser equivalente con otro que vaya por abajo porque al tratar de deformar el uno en el otro habría que romper un camino para luego volverlo a pegar pero al otro lado.

Similarmente, si un camino se enrosca 8 veces alrededor del hilo, éste no puede ser equivalente a otro que se enrosque 9 veces. Y así por el estilo. El conjunto de todos los caminos que son equivalentes a uno dado forma una clase de equivalencia. La idea es reemplazar todos los caminos equivalentes por uno cualquiera de ellos, cuando eso tenga sentido.

Por cada número natural $n$ tenemos dos caminos que representan a su clase, uno que se enrosca $n$ veces por encima y otro que se enrosca $n$ veces por debajo. El uno va con las manecillas del reloj, el otro en contra.

Si los dos puntos, inicial y final, coinciden, los caminos pueden componerse, ponerse unos detrás de otros para ser recorridos en línea, para dar caminos más largos. Por otro lado, cada camino tiene su inverso, que es el camino recorrido al revés. Hay un camino, en el espacio de las clases de homotopía que hace las veces de cero: es el camino que no se enrosca ni una sola vez y es contráctil hasta desaparecer.

En ese caso, tenemos que los caminos con la operación componer forman un grupo.

\

$>>>>>>>>>>>>>>>>>>>>>>>>>>>>>>>>>>>>>>>>>$

\

\

\color{black}

Como de costumbre vamos a filtrar los caminos que no aporten nada a la amplitud total. Por ejemplo, los caminos que se enrosquen un millón de veces alrededor del hilo no aportan nada, pues una ligera modificación de dichos caminos causaría que la fase oscilara alocadamente por el efecto de la desmesurada energía cinética. Por esa razón y referente  a la energía cinética tenemos que entre más derecho vaya un camino, más notable es su aporte. Por eso nos restringimos a los caminos que van directo de 'a' a 'b'. O sea a la clase de homotopía más baja.

\

Con las aclaraciones anteriores, la amplitud total se reescribe como:

$K(a,b)= $
$ \int \sb{a}\sp{b} e\sp{iE(\gamma)}$
$e\sp{(ie/\hbar c)\int\sb{\gamma} \vec v \cdot \vec A} D\gamma \sb{arriba}$
$+   \int \sb{a}\sp{b} e\sp{iE(\gamma)}$
$e\sp{(ie/\hbar c)\int\sb{\gamma} \vec v \cdot \vec A} $
$D\gamma \sb{abajo}$

$ =\int \sb{a}\sp{b} e\sp{iE(\gamma)}$
$e\sp{(ie/\hbar c)\int\sb{a}\sp{b} \vec A \cdot d\vec s}$
$ D\gamma \sb{arriba}$
$+   \int \sb{a}\sp{b} e\sp{iE(\gamma)}$
$e\sp{(ie/\hbar c)\int\sb{a}\sp{b} \vec A \cdot d\vec s}$
$ D\gamma \sb{abajo}$

Nos restringimos a la \index{clase de homotopía } \textbf{clase de homotopía} de los caminos directos, y en ese caso, todos los caminos por arriba son de una misma clase, la cual es diferente de la clase que forman  todos los caminos por abajo. Eso nos permite tomar un camino fijo, cualquiera y, como vimos,  el aporte del potencial vector es una constante en cada uno de las dos clases, sea por arriba, sea por abajo. Como es constante, el aporte del potencial vector sale de la integral:

$K(a,b)$
$=e\sp{(ie/\hbar c)\int\sb{a}\sp{b} \vec A \cdot d\vec s (arriba)}$
$\int \sb{a}\sp{b} e\sp{iE(\gamma)} D\gamma \sb{arriba}$

$+  e\sp{(ie/\hbar c)\int\sb{a}\sp{b} \vec A \cdot d\vec s (abajo)}$
$  \int \sb{a}\sp{b} e\sp{iE(\gamma)} D\gamma \sb{abajo}$

Cuando se trata de navegar en campos magnéticos, ellos no son neutrales para nada,   de tal forma que no es lo mismo ir de un punto a otro por debajo que por encima del hilo. Pero si la partícula es libre, entonces si da lo mismo ir por un lado que por el camino reflejado en el otro lado. Por consiguiente, el término asociado a la partícula libre vale lo mismo en ambos casos y lo podemos factorizar:

$K(a,b)$
$=[e\sp{(ie/\hbar c)\int\sb{a}\sp{b} \vec A \cdot d\vec s (arriba)}$
$+  e\sp{(ie/\hbar c)\int\sb{a}\sp{b} \vec A \cdot d\vec s (abajo)} ] \int \sb{a}\sp{b} e\sp{iE(\gamma)} D\gamma \sb{abajo}$

Esta expresión tiene la estructura siguiente:

$K(a,b)  = (w + z )y$

Por tanto, la probabilidad $P$ de llegar a 'b' es:

$P=|y|\sp{2}(|w|\sp{2} + |z|\sp{2} +  wz\sp{*}+w\sp{*}z) $

Los números $w$ y $z$ tienen norma 1, pero $wz\sp{*} +w\sp{*}z= cos (\phi - \psi)$  donde:

$\phi = (e/\hbar c)\int\sb{a}\sp{b} \vec A \cdot d\vec s (arriba)$

$\psi = (e/\hbar c)\int\sb{a}\sp{b} \vec A \cdot d\vec s (abajo)$

Por lo tanto

$P = [\int \sb{a}\sp{b} e\sp{iE(\gamma)} $
$D\gamma \sb{abajo}]\sp{2}(1+1 +cos (\phi - \psi)) $

Evaluando

$\phi - \psi = $
$(e/\hbar c)\int\sb{a}\sp{b} \vec A \cdot d\vec s (arriba)$
$-(e/\hbar c)\int\sb{a}\sp{b} \vec A \cdot d\vec s (abajo)$

$=(e/\hbar c)[\int\sb{a}\sp{b} \vec A \cdot d\vec s (arriba)-$
$\int\sb{a}\sp{b} \vec A \cdot d\vec s (abajo)]$

$=(e/\hbar c)[\int\sb{a}\sp{b} \vec A \cdot d\vec s (arriba)+ $
$\int\sb{b}\sp{a} \vec A \cdot d\vec s (abajo)]$

$=(e/\hbar c)\oint\ \vec A \cdot d\vec s $

$=(e/\hbar c)\int\ \int \nabla \times \vec A  $
$=(e/\hbar c)\int\ \int \vec B $

$= (e/\hbar c)\Phi \sb{B}$.

Hemos puesto $\Phi \sb{B}$ como el flujo magnético:

$$\Phi \sb{B}= \int\ \int \vec B
 \addtocounter{ecu}{1}   \hspace{4cm} (\theecu )      $$

Por lo tanto

$$P = [\int \sb{a}\sp{b} e\sp{iE(\gamma)} D\gamma \sb{abajo}]\sp{2}(2 +cos ((e/\hbar c)\Phi \sb{B}))
 \addtocounter{ecu}{1}   \hspace{4cm} (\theecu )      $$

Tenemos ya en nuestras manos una predicción sin necesidad de calcular esa integral de Feynman que nos falta. La predicción que podemos hacer es:

El conteo de partículas que llegan a 'b' tiene un componente sinusoidal  cuyo periodo puede predecirse exactamente. Se cumple un número entero de periodos cuando

$|e/\hbar c)\Phi \sb{B}|=2n\pi$

Es decir cuando:

 $$\Phi \sb{B}= (2n\pi )(\hbar c/e) = 2n\pi \hbar c/|e|= 4.135 n \times 10\sp{-7} Gauss-cm\sp{2}
 \addtocounter{ecu}{1}   \hspace{1cm} (\theecu )      $$

Predicción verificada por varios autores: la interacción electromagnética exige necesariamente una descripción cuántica y es una \index{teoría gauge}  \textbf{teoría gauge}. Es decir, el potencial vector que fue introducido prácticamente como un artificio matemático, es el que hace andar toda nuestra maquinaria, y el que permite hacer las predicciones que figuran en el montaje divulgado por Bohm -Aharonov, y el que admite una arbitrariedad.

\bigskip

Fin de la historia.

\section{CONCLUSION}

En el electromagnetismo clásico la libertad gauge estaba descrita por la posibilidad de sumar el gradiente de una función escalar diferenciable al potencial vector.  En tanto, en mecánica cuántica   la libertad  gauge consiste en poder correr arbitrariamente  la fase de todas las funciones de onda en el mismo punto y las predicciones no cambian. Pudimos además dilucidar la relación entre las libertades gauge de las 2 formulaciones: el proceso de cuantizar, empezado desde el lagrangiano,  hace que si uno tiene el grupo de invariancia clásico, el de las funciones diferenciables, uno necesariamente tendrá el grupo de las fases, el $\textbf{U(1)}$, en mecánica cuántica. También estudiamos un experimento divulgado por Bohm-Aharonov, en el cual las predicciones se basan directamente sobre el potencial vector, el cual bien puede considerarse innecesario a nivel clásico. Pero en este experimento, y en toda la teoría cuántica, dicho potencial es fundamental. Y es dicho potencial el causante de la libertad gauge tanto a nivel clásico como cuántico. Por eso resumimos diciendo: la teoría electromagnética es una teoría gauge $\textbf{U(1)}$. Con todo, nos sentimos un poco humillados al ver que todo eso parece más bien artificioso, es decir, no percibimos el poder creativo de la visión gauge.

\

\section{BIBLIOGRAFIA}

-Arnold V.I., Mathematical Methods of Classical Mechanics. Springer, 1978.

- A. Z. Capri (1985), Nonrelativistic Quantum Mechanics, Benjamin, Cummings, Menlo Park, California.

-Encarta 1995, Microsoft. 

- R.P. Feynman, (1962), The Theory of Fundamental Processes. Benjamin Inc, N.Y.

- R.P. Feynman, (1962), Quantum Mechanics and Path Integrals, McGraw-Hill, N.Y.

- Goudet G., \'{E}lectricit\'{e}, huiti\'{e}me \'{e}dition, Masson et Cie, 1967.

-Gohberg I., S. Goldberg, (1981), Basic Operator Theory. Birkhauser.

-Chan Hong-Mo, Tsou Sheung Tsun , (1993), Some Elementary Gauge Theory Concepts, World Scientific, Singapore.

- Lee D., Electromagnetic Principles of Integrated Optics, Wilet and Sons, 1986.

-Ohanian H., Classical Electrodynamics, 1988, Allyn and Bacon Inc.

- Purcell E.M., Electricity and Magnetism, Berkeley Physics Course, volume 2,  McGraw-Hill,1985.

-S. Paycha, (1997), 'Une petite introduction aux FIBRES DETERMINANTS Pour grands débutants', Apuntes Matemáticos No 38, Universidad de los Andes, Bogotá, Colombia.

- Pommaret J.F. (1987) \textit{Lie Pseudogroups and Mechanics}. Gordon and Breach. NY, London, Paris.

- J.J.Sakurai, (1985) Modern Quantum Mechanics, Benjamin, Cummings, Menlo Park, California.

-Yosida K., (1978), Functional Analysis, fifth edition, Springer-Verlag.

\chapter{RELATIVIDAD}

\Large

\centerline{RESUMEN}

\bigskip

Se demuestra que las leyes de Maxwell y la mecánica clásica son incompatibles. Como solución a esta contradicción se introduce la teoría de la relatividad. No se asume la invariancia de la velocidad de la luz sino que, a partir de la teoría de grupos, se deduce que debe existir una velocidad invariante. Se define el intervalo y se demuestra su invariancia. Se deduce la forma de las transformaciones de Lorentz, las cuales conservan el intervalo. Se enfatiza la noción de cuadrivector.
\
\
\

\color{red}

\bigskip
\bigskip
\bigskip
\normalsize
\color{black}
\section{INTRODUCCION}

En la primera sección estudiamos la interacción electromagnética desde el punto de vista no relativista, tanto la formulación hecha por Maxwell, dentro de la mecánica clásica, como su versión cuántica.

Ahora reinventaremos la relatividad especial para poder formular el electromagnetismo en lenguaje relativista. Es una obligación hacerlo, pues, como lo dedujo Einstein a sus 16 años, la mecánica clásica es incompatible con el electromagnetismo.

\section{EINSTEIN 16}

De acuerdo con las leyes de Maxwell, en el vacío pueden existir ondas electromagnéticas viajando a la velocidad de la luz, las cuales tienen sus componentes eléctricas y magnéticas mutuamente perpendiculares y perpendiculares también a la dirección de propagación. Pues bien,  cuando \index{Einstein} Albert \textbf{Einstein} tenía 16 a\~{n}os se formuló  el siguiente experimento mental (se acostumbra a decir gedankenexperimente, del alemán 'gedanken' que significa pensado): Viajemos con la onda a la velocidad de la luz. El resultado es que se percibirá una onda que no es viajera sino que es estacionaria. Pero sucede que al reemplazar una onda estacionaria en  las ecuaciones de Maxwell se encuentra una contradicción. Veamos: Las dos primeras ecuaciones de Maxwell, determinadas en el laboratorio, es decir en un marco de referencia atado a la tierra, son:

\addtocounter{ecu}{1}

$ \nabla  \times \vec E =-(1/c) \partial \vec B / \partial t   $

$\nabla  \times \vec B =-(1/c) \partial \vec E / \partial t $

las cuales tienen como solución:

$\vec E=\vec E\sb {z} = \vec k C sen(y - ct)$

$\vec B=\vec B\sb {x} = \vec i C sen(y - ct)$

Si viajamos con la onda a la velocidad de la luz, $c$, definiendo un marco de referencia ligado a la onda, lo que se ve es una onda estacionaria, independiente del tiempo, que resulta ser:

$\vec E=\vec E\sb {z} = \vec k C sen(y)$

$\vec B=\vec B\sb {x} = \vec i C sen(y )$

como la \index{onda  estacionaria} \textbf{onda es estacionaria}, la derivada parcial respecto al tiempo tanto de su componente eléctrico como magnético es cero. Pero por otro lado, el rotacional del campo eléctrico de dicha onda es $\vec i C cos(y)$, mientras que el de la componente magnética es $-\vec j C cos(y)$. Claramente, una onda estacionaria no es solución de las ecuaciones de Maxwell.

Lo que esto quiere decir es que el campo magnético no es una entidad en sí misma, sino que es un resultado del protocolo de observación. Pareciera  que el marco de referencia determinado por la tierra es increíblemente privilegiado, a tal punto de que puede crear las leyes de Maxwell. El adolescente Albert Einstein decidió que tal cosa no podría ser y  renunció a creer en la existencia de marcos privilegiados y en cambio se acogió a la fe ciega y pujante de  que las leyes de Maxwell eran las mismas en cualquier sistema de referencia inercial, pues representaban propiedades importantes de un objeto real, cuyas propiedades no dependían del \index{protocolo de observación} \textbf{protocolo de observación}.

Quedaba entonces un cuestionamiento: si se encuentra una inconsistencia al aceptar que el cambio de marco de referencia afecta la percepción de las leyes de Maxwell tal como lo predice la mecánica clásica, entonces las trasformaciones que rigen los cambios de marco de referencia no son como   acostumbramos a imaginarnos, y como lo hemos hecho al calcular la onda estacionaria, sino de alguna manera que queda por dilucidar.

La forma como uno tiende a imaginar la ley de transformación entre distintos marcos es simple: si una partícula se mueve con respecto a un marco, digamos un tren en movimiento, con velocidad hacia la derecha $u$ y si ese marco se mueve con respecto a otro, la carrilera, con velocidad hacia la derecha $v$, entonces la partícula se moverá con relación al segundo marco, la carrilera, con velocidad $u+v$.  Naturalmente que nosotros sabemos que eso es cierto a bajas velocidades, pero la contradicción encontrada nos induce a pronosticar que no es esa la forma de sumar velocidades cuando alguna de ellas es muy alta.

Vamos a estudiar este problema partiendo de su esencia: la \index{inconsistencia} \textbf{inconsistencia} entre las leyes de Maxwell y las reglas de transformación entre diferentes marcos de referencia aceptadas por el sentido común y la mecánica clásica. 

\section{LA COMPOSICION DE VELOCIDADES}

Podemos  comenzar a explorar la  siguiente intriga: ¿Cuál es la forma más general de una transformación de marcos de referencia inerciales que son compatibles con las leyes de Maxwell, es decir, que permiten que si se tiene una solución a las leyes de Maxwell en un marco, al ser transformado a otro marco también se tenga una solución a dicho sistema? Como esta pregunta es tan complicada, resuelta tan sólo hacia 1935, vamos a atacarla con la siguiente idea: si no hay marcos de observación privilegiados, tampoco hay campos privilegiados. Eso significa que el problema que hemos encontrado estudiando el campo electromagnético también podría ser encontrado estudiando cualquier otro campo. Por consiguiente, las reglas de transformación no deben involucrar ningún campo sino tan sólo las coordenadas. Nuestro punto de partida es pues:

\addtocounter{ecu}{1}

Postulado (\theecu): no hay marcos de referencia privilegiados.

\addtocounter{ecu}{1}

Corolario (\theecu): El conjunto de transformaciones que suman velocidades y que son compatibles con las leyes de Maxwell forman un \index{grupo} \textbf{grupo} conmutativo.

Significado: Un grupo es un conjunto provisto de una operación, la cual es cerrada, conmutativa, asociativa, tiene un elemento neutro y cada elemento tiene un inverso que revierte el elemento neutro. En nuestro caso, el grupo es un conjunto de transformaciones cuya operación es la composición.

Que la composición sea cerrada significa esto:   si dos transformaciones cambian una solución en otra solución, entonces las dos transformaciones seguidas   una de otra darán una nueva transformación que cambie soluciones en soluciones.  Que sea conmutativa, significa que da lo mismo hacer una trasformación y después otra o al revés. Que sea asociativa significa que el resultado de tres transformaciones es el mismo no importa si se hacen una por una o primero una y después el resultado de las otras dos, o primero las dos primeras y después la tercera. Que tenga elemento neutro implica que hay una transformación que no hace nada: esa es la transformación que no cambia de marco de referencia, es la identidad. Que tenga un inverso significa que cada transformación tiene su inversa que deshace lo que hizo la primera: si una transformación cambia del marco A al B, entonces la transformación inversa será la que cambie de B a A.

Ahora bien, si el conjunto de transformaciones no fuese un grupo, entonces habría la posibilidad de inventar recetas para crear leyes específicas verdaderas en unos marcos pero falsas en otros.

\addtocounter{ecu}{1}

Corolario (\theecu) No cualquier conjunto de operaciones forma un grupo. Por lo tanto, ser un grupo ya es una gran restricción. Podemos preguntarnos: Cuál es el grupo más general de transformaciones de coordenadas entre marcos inerciales? Esta pregunta por ser tan general es muy complicada. Necesitamos un principio simplificador que restringa las alternativas.

\addtocounter{ecu}{1}

Principio de la \index{relatividad} \textbf{relatividad} (\theecu): la forma de componer velocidades en distintos marcos de referencia inerciales no depende de absolutamente nada excepto de las velocidades en juego.

La forma antigua de componer velocidades, se llama Galileana y era sumándolas. Ella define un grupo que además cumple con el principio de relatividad. El problema es que tal forma de sumar no es la que obedece la naturaleza, pues fue la que utilizamos para calcular la onda estacionaria. Nuestra investigación debe conducirnos a una alternativa que sí sea fiel a la naturaleza.

\subsection{Forma general}

Hallemos la forma general de la composición de velocidades entre \index{marcos de referencia} \textbf{marcos de referencia} sabiendo que el conjunto de transformaciones forma un grupo que obedece el principio de relatividad.

Para fijar ideas, y siguiendo el estilo impuesto por Einstein, consideremos que un hay un objeto que se mueve con velocidad uniforme $u$ con respecto a un marco de referencia, digamos un tren. Este se mueve con velocidad $v$ con respecto a la tierra. Todos los movimientos son horizontales, paralelos a la carrilera que es recta. El principio de relatividad dice: la velocidad con que el móvil se mueve registrada desde la tierra es $w= f(u,v)$ y no depende de nada más. Especifiquemos algunas propiedades :

El elemento neutro es cero:

$f(y,0) =f(0,y) = y $

y en particular $f(0,0) = 0$

$f$ depende de dos entradas o variables. La derivada con respecto a la segunda variable es notada como $f\sb 2 (x,y) = \partial f(x,y)/\partial y$

La asociatividad dice:

$f(f(x,y),z) = f(x,f(y,z))$

y derivando con respecto al segunda variable y aplicando la regla de la cadena:

$f\sb 2 (f(x,y),z) = f\sb 2 (x,f(y,z))f\sb 2 (y,z)$

poniendo $z= 0$ da:

$f\sb 2 (f(x,y),0) = f\sb 2 (x,f(y,0))f\sb 2 (y,0) = f\sb 2 (x,y)f\sb 2 (y,0) $

Ahora guardemos $x$ fija, como un parámetro. Entonces $f\sb 2 (y,0)$ es una función de una sola variable, $y$. Eso permite ver esta ecuación como una ecuación diferencial ordinaria en la variable $y$:

$f\sb 2 (f,0) = (df/dy) f\sb 2 (y,0) $

Separando variables:

$ dy/f\sb 2(y,0) =df / f\sb 2(f,0)$

Ambos lados tienen la forma $ dz / f\sb 2(z,0)$, por lo que podemos integrar a ambos lados, y llamando $h(z) = \int dz /f\sb 2(z,0)$ vemos que estas integrales pueden diferir tan sólo por una constante:

$h(f) = h(y) + c(x)$

Teniendo en cuenta que $f(x,y) = f(y,x)$ podemos repetir el proceso y encontrar:

$h(f) = h(x) + c(y)$

de lo cual se concluye que $h(f) = h(f(x,y)) = h(x) + h(y)$

y asumiendo que $h$ es invertible, queda

\addtocounter{ecu}{1}
$$f(x,y) = h\sp{-1}(h(x) + h(y)) \eqno{(\theecu)} $$

En particular, y recordando que $f(0,0) = 0$, tenemos:

$h(f(0,0)) = h(0) = h(0) + h(0)$ de lo cual se concluye que $h(0) = 0$.

\bigskip

Entonces para determinar $f$, la forma de componer velocidades, es suficiente determinar $h$. Lo curioso es que para saber $h$ de largo a largo es suficiente saber $f$ en la vecindad de $y=0$. Eso se debe a que el siguiente problema de valor inicial resuelve $h$ de manera única, para lo cual recordamos que

$h(z) = \int dz /f\sb 2(z,0)$

por lo que

$h'(z) = 1/(f\sb 2(z,0) ) =  1/(\partial f(z,y)/\partial y )$ con $y=0$.

Podemos por tanto formular explícitamente el problema como sigue:

\bigskip

\addtocounter{ecu}{1}

Problema de valor inicial (\theecu):

Resolver la ecuación diferencial $h'(z) = 1/(\partial f(z,y)/\partial y ) $ con $y=0$.

$h(0) = 0$

Este problema sólo depende de $h$ puesto que $f$ está en términos de $h$.

\bigskip

La existencia y unicidad de la solución viene de teoremas generales para funciones diferenciales con derivada continua. La unicidad demostrada garantiza que los resultados tendrán valor científico, pues una predicción única puede contrastarse experimentalmente aunque se tenga poca resolución, pues basta repetir muchas veces el experimento para que al promediar errores quede un resultado experimental dentro de márgenes de error tan despreciables como se desee. Por otra parte, la misma unicidad nos garantiza que para encontrar la forma única y correcta se puede proceder por cualquier método. El que vamos a utilizar nosotros es el llamado Gedankenexperimente, el cual no es más que el análisis de una situación donde se puedan estudiar determinados factores, en nuestro caso, el principio de relatividad aplicado al grupo de transformaciones de marcos de referencia.

\subsection{Forma específica}

Una de las predicciones mas aterradoras de la relatividad es la \index{relatividad del tiempo}  \textbf{relatividad del tiempo}. Entender esto puede ser mucho más fácil si uno se da cuenta que un reloj es simplemente un movimiento periódico, un prototipo del cual entramos a analizar.

Tomemos entonces una plataforma plana delimitada por paredes, y sobre ella pongamos una bola a rebotar de una pared a otra. Ese será nuestro reloj. En esencia, lo que haremos será demostrar que si la bola viaja a nuestro encuentro, nos parecerá que esta viajará más rápido de lo esperado, y uno tendrá que concluir que el tiempo corre más rápido.

La bola no sufre ningún tipo de rozamiento ni de disipación y por consiguiente su movimiento de ir y venir rebotando es periódico.  La plataforma  está sobre un tren y el tren sobre la carrilera, de tal manera que el movimiento de la pelota sea paralelo a la carrilera, la cual es recta. Todas las velocidades son constantes.

La velocidad del tren respecto a la carrilera es $v$.
La velocidad de la pelota con respecto a la plataforma es $u$.
Todas las demás mediciones son respecto a la carrilera.

Para que la bola sirva de reloj y sus rebotes contra las paredes sean los tic-tacs, éstos deben ser frecuentes con respecto a la escala de tiempo dado por la velocidad del tren, por eso elegimos $u>>v$. Los tic-tacs no los produce la bola sino un sistema de rastreo con un laser (de otra forma la bola disiparía energía).

En primer término, visto desde la plataforma, el tiempo que se demora la bola para hacer tic es el mismo que se demora para hacer tac. Si las transformaciones de Galileo fuesen las correctas, lo mismo podríamos decir desde la carrilera. Pero eso no puede ser cierto, por la incompatibilidad de dichas transformaciones con las leyes de Maxwell. Sin embargo, podemos probar el siguiente teorema:

\bigskip

 \addtocounter{ecu}{1}

\textit{Teorema (\theecu). Decimos que el reloj produce tic cuando se choca con la pared delantera y que produce tac cuando se choca contra la pared trasera, de donde salió. Sea $t\sb 1 (u,v)$ el tiempo que el reloj se demora para hacer tic  después de haber hecho tac cuando el tiempo es medido desde la carrilera, es decir, el tiempo que gasta hacia adelante. Y sea $t\sb 2 (u,v)$ el tiempo que el reloj se demora para hacer tac (después de haber hecho tic, cuando es medido desde la carrilera), es decir, el tiempo hacia atrás. Entonces, $t\sb 1 (u,v)- t\sb 2 (u,v)$ es independiente de $u$, la velocidad de la bola-reloj.}

\bigskip

Demostración. Si guardamos $v$ fijo, la función $g(u)= t\sb 1 (u,v)- t\sb 1 (u,v)$ debe ser continua en $u$, puesto que hemos asumido diferenciabilidad de todas las funciones. Para especificar una función continua es suficiente especificarla en un conjunto de puntos suficientemente denso como para poder rellenar el resto  a mano alzada. Un tal conjunto es el de los números racionales, o sea, los quebrados. Un quebrado es de la forma $m/n$ donde ambos son enteros. Aunque esto no aparezca explícitamente en la siguiente discusión, está implícito.

Decir que $g(u)$ es constante, es lo mismo que decir, entonces, que $g$ tiene el mismo valor para dos velocidades genéricas distintas.

Consideremos dos bolas con velocidades $u$ y $u'$ referidas a la plataforma  que empiezan al mismo tiempo desde la parte trasera de la plataforma. Ellas van rebotando, rebotando, la una más adelantada que la otra hasta. Hemos elegido estas velocidades de tal forma que ambas llegarán   al tiempo al frente de la plataforma después de $m+1/2$ ciclos la primera y $n+1/2$ ciclos la segunda y que en otro tanto, ambas llegarán juntas a su punto de partida. Obsérvese que son ciclos y no tiempos, es decir, que es falso que en la mitad del tiempo estén en la mitad del recorrido. Cuando ambas bolas llegan al tiempo a alguna pared, el sistema de rastreo produce un    supertic o un supertac, los cuales pregonarán para todos los marcos de referencia posibles que las dos bolas llegaron al mismo tiempo al mismo lugar.

En particular, desde la carrilera, en donde hay un reloj de alta resolución,  podemos decir que las dos bolas, la que marcamos sin primas y la marcada con primas,  suenan al tiempo en esos determinados momentos. Como la    bola   toma un tiempo $t\sb 1(u,v)$ medido desde tierra para hacer todo el recorrido hasta el frente, mientras que toma  el tiempo $t\sb 2(u,v)$ para devolverse,   el tiempo transcurrido entre la salida y el supertic de la primera bola es

$supertic = (m+1)t\sb 1(u,v) + m t\sb 2(u,v)$

pues ha habido $m+1$ viajes hacia adelante y $m$ hacia atrás. Por otro lado, el tiempo transcurrido entre el primer supertic de la primera bola y su primer supertac se relacionan por :

$supertac = m t\sb 1(u,v) + (m+1)t\sb 2(u,v)$ .

pues ha habido m viajes hace adelante y $m+1$ hacia atrás.
Por lo que

$supertic -supertac =t\sb 1(u,v)-t\sb 2(u,v)$

Similarmente

$supertic' -supertac' =t\sb 1(u',v)-t\sb 2(u',v)$

y por la transitividad que se deriva de que las pelotas dan lugar a un único supertic y a un único supertac:

$t\sb 1(u,v)-t\sb 2(u,v)= t\sb 1(u',v)-t\sb 2(u',v)$

lo cual termina la demostración.

\bigskip

Revisemos ahora la implicaciones de nuestro análisis en cuanto a la forma de la función de composición de velocidades $f$.

\bigskip
 \addtocounter{ecu}{1}

\textit{Teorema (\theecu): Al componer dos velocidades, $u,v$ la velocidad resultante $w= f(u,v)$ está determinada por $h$ de los teoremas anteriores, y $h$ obedece la ecuación diferencial $h'(u) = 1/(1-Ku\sp 2)$ con valor inicial $h(0)=0$, el cual tiene solución única.}

\bigskip

Demostración: sea $L$ la longitud de la plataforma medida desde tierra.
Consideremos una bola que se mueve respecto al tren con velocidad $u$ desde la cola hacia el frente. El tren se mueve respecto a la tierra con velocidad $v$. La velocidad de la bola respecto a tierra es $f(u,v)$ y toma un tiempo $t\sb 1(u,v)$ medido desde tierra para hacer todo el recorrido hasta el frente. Por tanto el espacio recorrido es $f(u,v)t\sb 1(u,v)$. Pero eso también es la longitud de la plataforma mas el espacio recorrido por el tren:

$f(u,v)t\sb 1(u,v)= L + vt\sb 1(u,v)$.

Cuando la pelota regresa hacia atrás, lo hace con velocidad vista desde tierra $f(u,-v)$ tomando un tiempo $t\sb 2(u,v)$ y recorriendo una distancia

$f(u,-v)t\sb 2(u,v)= L - vt\sb 2(u,v)$.

Tenemos definido  un sistema de ecuaciones:

$f(u,v)t\sb 1(u,v)- vt\sb 1(u,v)= L  $.

$f(u,-v)t\sb 2(u,v)+ vt\sb 2(u,v)= L  $.

que se puede reorganizar

$(f(u,v)- v)t\sb 1(u,v)= L  $.

$(f(u,-v)+ v)t\sb 2(u,v)= L  $.

lo que da:

$t\sb 1(u,v)= L/(f(u,v)- v)  $.

$t\sb 2(u,v)= L/(f(u,-v)+ v)  $.

Por tanto, restando:

$t\sb 1(u,v)-t\sb 2(u,v) = L/(f(u,v)- v) - L/(f(u,-v)+ v) $.

Por el teorema anterior, no hay dependencia de $u$ en ninguna parte de la ecuación anterior. Llamemos

$g(v) = 1/(f(u,v)- v) - 1/(f(u,-v)+ v)$

sacando común denominador:

$g(v) = \frac{\rm f(u,-v)+ v - (f(u,v)- v) }{\rm (f(u,v)- v)(f(u,-v)+ v}$
$=\frac{\rm f(u,-v)- f(u,v)+2 v }{\rm (f(u,v)- v)(f(u,-v)+ v)}$
$=-\frac{\rm f(u,v)- f(u,-v)-2 v }{\rm (f(u,v)- v)(f(u,-v)+ v)}$

$=-\frac{\rm f(u,v)-f(u,0) + f(u,0)- f(u,-v)-2 v }{\rm (f(u,v)- v)(f(u,-v)+ v)}$

Dividiendo en todos lados por $2v$, y tomando el límite de $v$ hacia cero, nos damos cuenta que en el numerador del lado derecho quedan dos expresiones correspondientes cada una a $(1/2) (\partial f/\partial v) $, las cuales se suman:

$\frac{\rm g(v)}{\rm 2v} =-\frac{\rm  [\partial f(u,v)/\partial v]\sb{v=0}-1}{\rm (f(u,v)- v)(f(u,-v)+ v)}$
$=\frac{\rm  [1- \partial f(u,v)/\partial v]\sb{v=0}}{\rm (f(u,v)- v)(f(u,-v)+ v)}$

Pero tomando  $ \lim_{v\to 0}$  en el denominador:

$(f(u,v)- v)(f(u,-v)+ v)\rightarrow (f(u,0)- 0)(f(u,-0)+ 0)  \rightarrow (u)(u) = u\sp 2$

Por lo tanto, el límite existe y le llamamos $K$:

$lim \sb{v\rightarrow 0} \frac{\rm g(v)}{\rm 2v}= K  =\frac{\rm  [1- \partial f(u,v)/\partial v]\sb{v=0}}{\rm u\sp 2}$

Despejando $[\partial f(u,v)/\partial v]\sb{v=0} = 1-K u\sp 2  $

pero como
$h'(z) = 1/(\partial f(z,y)/\partial y ) $ evaluada en $y=0$.

queda que

$h'(u) = 1/(1-K u\sp 2)$

\bigskip

Podemos ahora resolver el problema de valor inicial con $h(0)=0$ y obtenemos:

\bigskip

\addtocounter{ecu}{1}

\textit{Teorema de composición de velocidades} (\theecu):

$$ w= f(u,v) = (u+v)/(1+Kuv)$$

Demostración. Al integrar la ecuación diferencial para $h$ hay dos casos, el uno cuando el denominador puede partirse en fracciones simples, $K>0$ y el otro cuando el denominador es irreducible, $K<0$. La integración para $K>0$, por fracciones parciales da:

$h(u)= \frac{\rm 1}{\rm 2\sqrt{ K}} ln \frac{\rm 1 + \sqrt{ K} u}{\rm 1-\sqrt{ K}u} $

debemos reemplazar  esta expresión en

$f(x,y) = h\sp{-1}(h(x) + h(y))$

por lo que necesitamos saber $h\sp{-1}$. Sea

$r= \frac{\rm 1}{\rm 2\sqrt{ K}} ln \frac{\rm 1 + \sqrt{ K} u}{\rm 1-\sqrt{ K}u} $

$2\sqrt{ K} r =  ln \frac{\rm 1 + \sqrt{ K} u}{\rm 1-\sqrt{ K}u} $

Tomando exponencial

$\frac{\rm 1 + \sqrt{ K} u}{\rm 1-\sqrt{ K}u} =exp(2\sqrt{ K} r) $

$1 + \sqrt{ K} u = exp(2\sqrt{ K} r) (1-\sqrt{ K}u)$

$1 + \sqrt{ K} u =exp(2\sqrt{ K} r) - exp(2\sqrt{ K} r) \sqrt{ K}u $

$h^{-1}(r) = u = \frac{exp(2\sqrt{ K} r) -1}{\sqrt{ K}( 1 + exp(2\sqrt{ K} r))}$

pero como

$f(x,y) = h\sp{-1}(h(x) + h(y))$

entonces

$f(u,v) = h\sp{-1}(h(u) + h(v))$

 $ = h\sp{-1}(\frac{\rm 1}{\rm 2\sqrt{ K}} ln \frac{\rm 1 + \sqrt{ K} u}{\rm 1-\sqrt{ K}u} +
  \frac{\rm 1}{\rm 2\sqrt{ K}} ln \frac{\rm 1 + \sqrt{ K} v}{\rm 1-\sqrt{ K}v})$

componiendo con  $h^{-1}$, eso da:

$$f(u,v) = - \frac{1- \frac{\rm 1 + \sqrt{ K} u}{\rm 1-\sqrt{ K}u}\frac{\rm 1 + \sqrt{ K} v}{\rm 1-\sqrt{ K}v}}{\sqrt{k} ( 1 + \frac{\rm 1 + \sqrt{ K} u}{\rm 1-\sqrt{ K}u}\frac{\rm 1 + \sqrt{ K} v}{\rm 1-\sqrt{ K}v} )}$$

$$= - \frac{(1-\sqrt{ K}u) (1-\sqrt{ K}v) - (1+\sqrt{ K}u)(1+\sqrt{ K}v)}{\sqrt{K} ( (1-\sqrt{ K}u) (1-\sqrt{ K}v) + (1+\sqrt{ K}u)(1+\sqrt{ K}v))}$$

$$= - \frac{1-\sqrt{ K}v -\sqrt{ K}u + Kuv - 1-\sqrt{ K}u - \sqrt{ K}v - Kuv}{ \sqrt{K}( 1-\sqrt{ K}v -\sqrt{ K}u +Kuv + 1 + \sqrt{ K}u+ \sqrt{ K}v +Kuv})$$

$$  = -\frac{-2\sqrt{K}v -2\sqrt{K}u }{\sqrt{K}(2+2Kuv)}$$

que al simplificar se convierte en

$$w= f(u,v) = (u+v)/(1+Kuv)$$

La generalidad de esta expresión se comprueba por ejemplo notando que si $K=0$ se tiene la suma de Galileo. Ahora entraremos a argumentar que la \index{velocidad invariante}  velocidad de la luz es invariante, es decir que se ve igual desde todos los marcos de referencia inerciales. Veamos primero de qué manera demostraremos que existe una velocidad invariante.

\bigskip

\addtocounter{ecu}{1}

\textit{Teorema de la velocidad invariante(\theecu): Existe una velocidad invariante y su valor es $K\sp{-1/2}$, la cual es la misma en todos los sistemas de referencia. Nótese que nosotros deducimos que debe haber una velocidad invariante, al contrario de lo usual, en lo cual se asume que la velocidad de la luz es invariante: la invariancia de la velocidad de la luz se vuelve un corolario de la oscuridad.}

Este teorema es debido a Mermin (1984), quien  describió la sensación que le inspiraban sus resultados diciendo: 'working in the darkness may be illuminating'.

\bigskip

Demostración: Asumiendo que  $K>0$ entonces, substituyendo $u = K\sp{-1/2}$ en

$w= f(u,v) = (u+v)/(1+Kuv)$

se obtiene:

$w= f(u,v) = (K\sp{-1/2}+v)/(1+K K\sp{-1/2}v)$
$= K\sp{-1/2}(1+ K\sp{1/2}v)/(1+ K\sp{1/2}v)$

$w=K\sp{-1/2}$

lo que dice que la velocidad invariante es la misma aunque se mida desde dos sistemas inerciales diferentes.

Claro que nosotros no hemos demostrado la unicidad de la velocidad invariante. Para eso despejemos $u$ de

$u= (u+v)/(1+Kuv)$

$u(1+Kuv)=u+v$

$u + Ku\sp 2 v = u + v$

$Ku\sp 2 v = v$

$Ku\sp 2  = 1$

$u = K\sp{-1/2}$

lo que demuestra que hay una única velocidad invariante.

El experimento de Michelson-Morley, buscando el efecto del éter, tuvo la gloria de ser la primera indicación seria de que la velocidad de la luz era la misma en todos los sistemas. Así que ese es el valor de $K$. Como es costumbre, la velocidad de la luz se nota $c$, y entonces, de la última ecuación tenemos:

$c =  K\sp{-1/2} = 1/\sqrt{K}$

es decir $K = 1/(c^2)$, por lo que la ley de composición de velocidades queda al final como:

\addtocounter{ecu}{1}

$$w = (u+v)/(1+c\sp{- 2}uv)\eqno{(\theecu)}$$

vale la pena notar que las unidades cuadran: $w$ viene en unidades de velocidad, pues el denominador no tiene unidades.

\section{LAS TRANSFORMACIONES DE LORENTZ}

La teoría fundamental de la relatividad nos dice que la relación de composición de velocidades es la fundamental. Sin embargo,  las transformaciones entre coordenadas han sido tradicionalmente las que se han considerado como importantes y son ellas las que por antonomasia se denominan \index{transformaciones de Lorentz} \textbf{transformaciones de Lorentz}. Para nosotros hallarlas es por demás sencillo:

Todo lo que tenemos que hacer es reescribir la expresión para la composición de velocidades como un quebrado de la forma espacio sobre tiempo, de donde podemos leer en el numerador la forma para transformar coordenadas espaciales y en el denominador la forma para transformar la coordenada temporal. Hagámoslo. Interpretamos la fórmula

$$w = (u+v)/(1+c\sp{- 2}uv)\eqno{(\theecu)}$$

de la siguiente forma: $w$ es la velocidad de un móvil referida a la tierra, que usualmente se nota $v_x$, en tanto que $v'_x$ es la velocidad del móvil con respecto al tren, el cual se mueve a velocidad $v$. De igual modo, lo que tenga primas es con respecto al tren, lo que no es con respecto a la tierra. Nos da:

$$v_x = \frac{v'_x+v}{1+c\sp{- 2}v'_x v}   $$

Definimos

\

$\beta = v/c$

$\gamma = 1/(1-(v/c)\sp 2 )\sp{1/2}= 1/(1-(\beta)\sp 2 )\sp{1/2}$

\

Ahora podemos elaborar la expresión original:

$v_x = \frac{v'_x+(v/c)c}{1+c\sp{- 1} v'_x (v/c)}  $
$=  \frac{v'_x+\beta c}{1+c\sp{- 1} v'_x \beta } $
$=  \frac{\gamma v'_x+\gamma \beta c}{\gamma +\gamma c\sp{- 1} v'_x \beta }$
$ = \frac{\gamma (dx'/dt')+\gamma \beta c}{\gamma +\gamma c\sp{- 1} (dx'/dt') \beta }$

\

Sacando denominador común tanto arriba como abajo y simplificando obtenemos:

$$v_x = \frac{dx}{dt} = \frac{\gamma dx'+\gamma \beta c dt'}{\gamma dt'+\gamma c\sp{- 1}  \beta dx' }$$

lo cual presenta unidades de espacio en el numerador y de tiempo en el denominador, pues ni $\gamma$ ni $\beta$ tienen unidades. Por tanto, la forma como las coordenadas se transforman es:

\

$dx = \gamma dx'+\gamma \beta c dt' $

$dt = \gamma dt'+\gamma c\sp{- 1}  \beta dx' $

\

Ahora bien, estas diferenciales vienen de una expresión que, aparte de constantes no medibles, tienen la forma:

\

$x = \gamma x' + \gamma \beta c t'$

$t = \gamma t' + \gamma \beta x'/c$

\

que son las archiconocidas \index{transformaciones de Lorentz} \textbf{transformaciones de Lorentz}. Cuando estas ecuaciones se reescriben en forma matricial, a la matriz resultante se nota $\Lambda$.

Las transformaciones  de Lorentz forman un    un grupo con la operación de composición, el cual se denomina Grupo de Lorentz, el cual es un subgrupo del grupo  de transformaciones que es compatible con las ecuaciones de Maxwell.

\section{EL TIEMPO PROPIO Y EL INTERVALO }

El \index{tiempo propio} \textbf{tiempo propio} es el tiempo medido por un reloj que viaja con la partícula. Se supone que el reloj no se desajusta para nada con el viaje de la partícula. La idea de un reloj atómico luce atractiva y se ha utilizado en experimentos reales: se monta un reloj atómico en un avión el cual le da la vuelta al mundo. Otro reloj queda en tierra. Al terminar el tour del primero, se comparan las lecturas de los dos relojes y se contrastan con las predicciones de la teoría, la cual dice que el tiempo no es una entidad absoluta sino que depende de la velocidad relativa con respecto al sistema de observación: No ha habido reclamos.

\

Es importante notar que el tiempo propio es un invariante relativista pues podemos imaginar que el reloj que viaja con la partícula es digital y que se comunica con los demás sistemas por medio de un lenguaje hablado. Eso quiere decir que, después de un cierto momento,  ningún sistema tiene incertidumbre acerca del tiempo propio de la partícula y por lo tanto vale lo mismo en todos los sistemas. Por construcción, es un invariante de la teoría.

\

Las transformaciones de Lorentz mezclan el espacio y el tiempo en un todo que funciona como una unidad: el espacio-tiempo. Podemos, sin embargo, rastrear qué le pasa al tiempo separado del espacio. Para ello, en la segunda ecuación consideramos que $x'=0$, lo cual implica que

$t = \gamma t' $

Ahora bien, si $t'$ indica el tiempo propio, entonces vemos que lo que para una partícula que lleva su reloj puede parecerle un tiempo dado, a un observador en tierra puede parecerle un tiempo muy largo, cuánto más largo cuanto más rápido viaje la partícula. Esto ha animado a los experimentalistas que estudian reacciones entre partículas elementales, que pueden suceder muy rápidamente, a realizar experimentos con partículas a muy alta velocidad para poder verlas en cámara lenta.

\bigskip

Veamos ahora una de las consecuencias más sencillas y más poderosas de la invariancia de la velocidad de la luz: la conservación de una cantidad que se llama el intervalo y que define una estructura muy semejante a una métrica y que técnicamente se denomina pseudo-métrica, pues lo único que le falta para ser métrica es ser no negativa. El intervalo sale de un cambio de unidades en  el tiempo propio.

\addtocounter{ecu}{1}

\bigskip

\textit{Definición (\theecu): El \index{intervalo}  \textbf{intervalo} $s$ se define como: $s^2 = c^2\tau^2$ donde $\tau$  es  el tiempo propio. Observemos que el intervalo tiene unidades de espacio mientras que el tiempo las tiene de tiempo.}

\bigskip

\addtocounter{ecu}{1}

\textit{Teorema (\theecu): El intervalo cumple la propiedad:}

 $s^2 = c^2\Delta \tau^2 =  \Delta (x^o)^2 -\Delta  (x^1)^2 -\Delta  (x^2)^2-\Delta  (x^3)^2$

\textit{donde $x^o = c\tau $ y las demás coordenadas son las espaciales $x$, $y$, $z$. El $\Delta$ significa que a lo largo de la trayectoria de la partícula, uno puede definir el tiempo propio cero a cualquier momento, y a partir de ahí comenzar a contar el tiempo. Eso define también una posición inicial y otra final, entre las cuales se calcula el $\Delta$.}

\

Demostración:

\

Recordando que $\beta = v/c$ y que $\gamma^2 = \frac{1}{1-\beta^2}$ podemos escribir:

$s^2 = c^2\tau^2 = \frac{1}{1- \beta^2} (c^2\tau^2) (1- \beta^2 ) $

$= \gamma^2 c^2 \tau^2(1-\beta^2)$

$ = \gamma^2 c^2 \tau^2- \gamma^2 \beta^2c^2 \tau^2$

Como además, el tiempo $t$ medido desde el sistema de observación y el tiempo propio se relacionan por $t = \gamma \tau$, obtenemos:

$c^2\tau^2 = c^2t^2- \gamma^2 (v^2/c^2)c^2 \tau^2$

 $= c^2t^2 -  \gamma^2 v^2 \tau^2 = c^2t^2 - v^2  \gamma^2 \tau^2 = $

 $= c^2t^2 - v^2  t^2$

Como $x = vt$ nos queda:

$c^2\tau^2 = c^2t^2 - x^2$

Para hacer más explícito el hecho de que los sistemas de referencia son básicos, podemos reescribir la ecuación anterior como

$c^2\Delta \tau^2 = c^2\Delta t^2 -\Delta  x^2$

Involucrando las demás coordenadas podemos decir que, en general,

$c^2\Delta \tau^2 = c^2\Delta t^2 -\Delta  x^2 -\Delta  y^2-\Delta  z^2$

Además, el término $ c\Delta t $ tiene dimensiones de espacio, por lo que todo se reescribe como una ecuación que define el intervalo:

 $s^2 = c^2\Delta \tau^2 =  \Delta (x^o)^2 -\Delta  (x^1)^2 -\Delta  (x^2)^2-\Delta  (x^3)^2$

\

Podemos resumir todo como sigue:

\bigskip

\addtocounter{ecu}{1}

Teorema (\theecu):
\addtocounter{ecu}{1} Todos los enunciados siguientes son equivalentes:

1) El conjunto de transformaciones de un marco inercial a otro es un grupo (el grupo de Lorentz), cuya constante característica (velocidad de la luz) es finita.

2) Todos los observadores inerciales miden la misma velocidad de la luz.

3) El intervalo es un invariante relativista.

4) El tiempo propio es un invariante relativista.

\bigskip

La demostración de que la invariancia de la velocidad de la luz implica las transformaciones de Lorentz fue el resultado que prendió todo y que fue debido a Einstein en 1905.

Demostremos que si el intervalo es un invariante, entonces la velocidad de la luz también. En efecto:

Para el observador A: $(ds_A)^2 = c^2 (dt_A)^2 - (dx_A)^2$. Pero si es la luz   la que se mueve, $dx_A = cdt_A$ entonces $(ds_A)^2 = 0$.

Como el intervalo es invariante,  para el observador B tenemos

$(ds_B)^2 = c^2 (dt_B)^2 - (dx_B)^2= (ds_A)^2 =  0$ .

Es decir,  $c^2 (dt_B)^2 - (dt_Bdv_B)^2 = 0$.

de donde deducimos que $v = c$.

\section{CUADRIVECTORES}

Queremos destacar que las coordenadas ya no se transforman unas independientes de otras, sino que forman una unidad llamada  \index{cuadrivector} \textbf{cuadrivector}, el cual  se transforma como un todo. Hay otras tetrapletas que también se transforman como un todo, aunque su ley puede diferir del cuadrivector de coordenadas.

A continuación referimos algunos cuadrivectores importantes y que aparecen de sorpresa en cualquier parte.

a) El espacio-tiempo $[ct,x\sp{1},x\sp{2},x\sp{3}]= [x\sp{\mu}]$

b) La cuadri-velocidad $[cdt/d\tau,x\sp{1},dx\sp{2}/d\tau,dx\sp{3}/d\tau]$ done $\tau $ es el tiempo propio.

c) La cuadri-aceleración

d) El vector energía-momento

$[p\sp{\mu}]=mu\sp{\mu}= [E/c,p\sb{x},p\sb{y},p\sb{y}]$

e) El 4-vector corriente $[\rho,j\sb x , j\sb y, j \sb z ]$

Antes: carga=$\rho=qn$, vector corriente $=\vec j=qn\vec v$

Ahora: $j\sp \mu= qn \sb o u \sp\mu$ donde $n= n \sb o / \sqrt{1-\gamma \sp 2}$

Esta es la definición del 4-vector corriente. Su primer componente es

$j\sp o =   qnc=$ carga x c.

Las otras componentes son :

$j\sp{\mu} =  qnv $

y todo el cuadrivector se transforma como:

$\rho '=\frac{\rm \rho -Vj\sb x/c\sp{2}}{\rm \sqrt{1-\gamma\sp{2}}}$

$j'\sb x =\frac{\rm j \sb x -V\rho}{\rm \sqrt{1-\gamma\sp2 }}$

f) El cuadrivector potencial $[A\sp{\mu}]$ que será estudiado extensamente en el próximo capítulo.

\section{RESUMEN}

Después de haber constatado que las leyes de Maxwell no son compatibles con las transformaciones de Galileo sobre composición de velocidades en marcos de referencia, nos pusimos a investigar la existencia de otras formas de componer dichas velocidades que  fuesen compatibles con dichas leyes. Hallamos la  forma general de componer velocidades y se predijo la invariancia de una única velocidad específica, la cual fue determinada experimentalmente como la velocidad de la luz. Dedujimos la forma de sumar velocidades de donde elaboramos las transformaciones de Lorentz entre coordenadas. Definiendo el tiempo propio como el tiempo atado a la partícula, pudimos demostrar que al cambiar de unidades se obtiene el intervalo, un invariante relativista de primordial importancia. Enfatizamos  la noción de cuadri-vector, cuyo prototipo es $\vec u = [ct, x,y,z]$ y cuya  ley de transformación es $\vec  u'=\Lambda \vec u $ donde $\Lambda $ es una matriz de Lorentz.

\section{REFERENCIAS}

\enumerate

\item Mermin N. David, (1984), \textit{Relativity without light}, Am.J.Phys, 52(2), February 1984.

\item  Brédov M., V. Rumiántsev, I. Toptoguin, (1986) \textit{Electrodinámica clásica}, Mir, Mocú.

\item Ohanian H., (1988), \textit{Classical electrodynamics}, Allyn and Bacon, Inc. Boston.

\item Taylor, E.F. and J.A. Wheeler, 1966, \textit{Spacetime physics}, Freeman, San francisco.

\chapter{ FORMAS DIFERENCIALES }

\Large

\centerline{RESUMEN}
\bigskip
\bigskip
\bigskip
El cálculo diferencial en varias variables, sobre $\Re \sp n$, es reformulado y extendido en términos de formas diferenciales, tal como lo propuso Grassmann. Estas formas  permiten una generalización natural del cálculo a espacios curvos y variedades.   Son reformuladas en este lenguaje  las leyes de Maxwell y la libertad gauge.
\
\
\

\bigskip
\bigskip
\bigskip
\normalsize
\color{black}
\section{INTRODUCCION}

La maquinaria que vamos a presentar en esta sección fue concebida por \index{Grassmann}  \textbf{Grassmann} en la primera mitad el siglo XIX y perfeccionada por Henri Cartan  hacia principios del siglo XX. Nuestra introducción se vale de la diferencia natural que en un experimento elemental existe entre estímulo  al sistema estudiado y  la maquinaria interior de dicho sistema, el cual causa la reacción que pueda observarse.

Las observaciones y los estímulos y en algunos casos las respuestas del sistema se representan por plps (\textbf{paralelepípedos}), \index{plp(paralelepípedo)} tal como sale de la experiencia cotidiana:

Ej 1. Se mide una longitud, contando los pasos necesarios para recorrerla. Un $\vec {paso}$ es un vector, o un plp de dimensión uno: es un 1-plp.

Ej 2. Medir una área implica superponer un cuadrado de 1m de lado, el cual es un 2-plp.

Ej 3. Medir el agua contenida en un tanque consiste en contar las veces que se llena un cubo de 1dm de lado. Tal cubo es un 3-plp.

Ej 4. Aplicamos una fuerza sobre una pared y observamos que la pared se deforma de manera diferente en cada dirección. Por tanto, ante un vector, la fuerza, obtenemos una tripleta de números, o sea,  un 1-plp.

Nosotros describimos un plp por un método que elimina la redundancia: un rombo se describe por dos aristas o vectores que salen de un mismo vértice, pues las otras dos salen por paralelismo. Un cubo lo describimos por 3 vectores que salen de un mismo vértice.   También hay que tener en cuenta que, por ejemplo, el paralelepípedo puede ser de 2 dimensiones pero puede estar en cualquier parte y dirección de  $\Re\sp 3$ o de $\Re^n.$

En general, a un paralelepípedo se le designa por el grupo de aristas que parten de un vértice de tal manera que las demás aristas se puedan reconstruir por paralelismo. Formalmente,

\addtocounter{ecu}{1}

\textit{Def (\theecu): Un n-paralelepípedo (n-plp) en $\Re \sp m, $  es un elemento de la forma $(P,v_1,.. ..,v_n), $ donde $P \in \Re\sp m$ es el punto donde se pone el vértice de referencia del plp cuyas  aristas generadoras  son los $v\sb j \in \Re \sp m $. En general, el punto P quedará sin especificar y sin referir.}

\bigskip

Por derecho legítimo, esta definición también se extiende al espacio-tiempo: la posición de n-partículas está dada por $(v_1,.. ..,v_n), $ donde $v\sb j \in \Re \sp 4 = Minkowski = M$, que lo tomamos como si fuese simplemente $\Re \sp 4$.

Por otra parte, la estructura general de un experimento elemental está dada por
la siguiente cadena de causa-efecto: el experimentador hace algo, la naturaleza reacciona y produce un resultado. En un esquema, esto se representa así:

$$\pmatrix{Resultado\cr
    del\cr
    experimento\cr}=
\pmatrix{Reacci\acute on\cr
    de \:la\cr
    naturaleza\cr}
\pmatrix{Est\acute {\i} mulo\cr
    del\cr
    experimentador\cr}
$$

Nuestro esquema coincide exactamente con el del álgebra lineal: nuestro trabajo es fabricar una plataforma matemática que nos permita convertir dicha coincidencia en una realidad natural. Para tal fin, necesitamos que el efecto de un experimento real se pueda calcular como la integración de efectos infinitesimales, y que dichos efectos infinitesimales se describan por vectores y matrices o de sus correspondientes generalizaciones.

Además debemos incorporar la estructura del espacio-tiempo: los eventos siempre tendrán sus 4 coordenadas, una temporal y tres espaciales, de tal manera que las transformaciones entre diversos sistemas de referencia inerciales se efectúen por transformaciones de Lorentz.

Comencemos con un ejemplo de interés: al mover una masa un desplazamiento dado infinitesimal en contra de una fuerza, gastamos una cierta energía.

$$\pmatrix{Resultado:\cr
    energ\acute {\i}a\cr
    gastada\cr}=
\pmatrix{Reacci\acute on:\cr
    fuerza  \cr
    en \: contra\cr}
\pmatrix{Est\acute {\i} mulo:\cr
    desplazamiento \: de \cr
    una \:masa\cr}
$$

Pero atención: nosotros ya sabemos que la energía es apenas una coordenada de un cuadrivector, el de energía-momento. Por lo tanto, sería mejor no imaginarse que la ley fundamental que nos permita calcular la energía gastada venga como una ecuación por separado y especialmente destinada a nuestra pregunta puntual. Con esta nueva perspectiva, y teniendo cuidado de equiparar unidades, vemos que la estructura infinitesimal del experimento queda

$$\pmatrix{Vector\cr
    energ\acute {\i}a\cr
    momento\cr}=
\pmatrix{Descriptor \cr
     del \cr
    campo\cr}
\pmatrix{\cr
    \vec dx/dt  \cr
    \cr}
$$

EL $\vec {dx}$ es claramente un 1-plp, al igual que la cuadrivelocidad, notada abusivamente $\vec dx/dt$, lo mismo que el cuadrivector energía-momento. Por lo que inferimos que  lo que hemos llamado descriptor de campo, lo que produce la fuerza, debe estar representado por lo que haga las veces de una matriz. Maticemos esta apreciación mucho mejor:

El estímulo del experimentador se representa por un n-plp, al igual que las observaciones que éste registra. Por otra parte, la naturaleza tomará el plp que representa el estímulo y lo transformará en otro. De dicho plp resultante, el experimentador tiene que obtener una medición, producida en números reales. Pues bien, la manera más natural de obtener números a partir de un plp es por medio del determinante o elemento de volumen.

El \index{determinante} \textbf{determinante} es  un operador multilineal  alternado. Es decir: es lineal en cada arista del n-plp estímulo, y es antisimétrico en cada par de vectores o aristas del  plp estímulo: al intercambiar dos aristas, el determinante cambia de signo. Y además  cumple con una normalización o escala: el volumen del cubo unitario es uno (en todas partes asumimos la base natural con su orientación acostumbrada).

Con todo lo anterior en mente, para  caracterizar  la reacción de la naturaleza y el proceso de experimentación y de medición de respuestas, comenzamos con la siguiente

\addtocounter{ecu}{1}

\textit{Def (\theecu): una \index{forma} 1-\textbf{forma}  $\omega $ es un operador lineal que a cada 1-plp o vector $\xi$ le hace corresponder un número real $\omega (\xi)$.} 

\bigskip

Lineal significa que

$\omega (\alpha \xi_1 + \beta \xi_2) = \alpha \omega (\xi_1 )+ \beta \omega (\xi_2)$

A partir de 1-formas podemos generar k-formas que operan sobre k-plps como sigue:

\addtocounter{ecu}{1}

\textit{Def (\theecu): El \index{producto cu\~{n}a o wedge} \textbf{producto cu\~{n}a o wedge} $\wedge $ de k 1-formas $\omega  \sb 1, .. .. ,\omega \sb k $ es una k-forma monómica $\omega \sp k= \omega_1  \wedge .. .. \wedge \omega \sb k $ que opera sobre un k-plp $ \Xi \sp k= (\xi_1 ,.. .., \xi \sb k)$, donde cada $\xi \sb j$ es un 1-plp, de la siguiente manera:}

$$\omega \sp k\ (\Xi \sp k) = (\omega \sb 1 \wedge .. .. \wedge \omega \sb k)((\xi_1 ,.. .., \xi \sb k))=Det
\pmatrix{\omega \sb 1  (\xi_1 )&.. ..&\omega \sb k (\xi_1 )\cr
         \vdots &.. ..&\vdots\cr
         \omega \sb 1   (\xi \sb k)&.. ..&\omega \sb k (\xi \sb k)\cr}
\in \Re
$$

Los monomios resultantes operan como los determinantes:  lineal en cada arista, y alternada, cambia de signo al intercambiar dos aristas. Un producto wedge siempre se aniquila cuando hay elementos repetidos. Eso se argumenta, en el caso más sencillo, así: 

como las formas son alternadas, $\omega \sb i \wedge \omega \sb j = - \omega \sb j \wedge \omega \sb i$, pero si $i=j$, queda una forma que es igual a su contraria aditiva: es la forma nula. Por eso, las repeticiones aniquilan las formas.

Los monomios pueden sumarse y multiplicarse por un escalar con toda naturalidad dando formas polinómicas. Por ejemplo: $3\omega \sb 1 \wedge \omega \sb 2 + 6\omega \sb 1 \wedge \omega \sb 2 = 9\omega \sb 1 \wedge \omega \sb 2   $, dió un monomio. Por otra parte $3\omega \sb 1 \wedge \omega \sb 2 + 6\omega \sb 1 \wedge \omega \sb 3 $  es un polinomio que tiene que evaluarse averiguando el efecto de cada monomio sobre un 2-plp cualquiera. Se deduce que los polinomios de igual rango pueden sumarse y multiplicarse por un escalar: tienen estructura de espacio vectorial, con base, dimensión y todo.

Los monomios también pueden multiplicarse naturalmente, donde el producto es el acostumbrado wedge, con lo cual, el producto es cerrado, producto de formas da una forma, no es conmutativo, pero es asociativo:

$(\omega \sb 1 \wedge \omega \sb 2) \wedge \omega \sb 3 = \omega  \sb 1 \wedge (\omega \sb 2 \wedge \omega \sb 3)$

En particular, $\omega \sb 1 \wedge \omega \sb 1 \wedge \omega \sb 3 = - \omega \sb 1 \wedge \omega \sb 1 \wedge \omega \sb 3 =0$
donde hemos usado la propiedad de alternación del producto.

En $\Re \sp n $ hay una n-forma distinguida $\omega \sp n $ que mide el volumen de un n-plp de tal manera que el volumen del n-plp cubo-unitario $(e \sb 1 ,.. .., e \sb n)$   es uno. Tal forma es el producto de n 1-formas elementales $\omega \sb 1 , .. .. ,\omega \sb k $ y debe cumplir:

$$\omega \sp n\ (\Xi \sp n) = (\omega \sb 1 \wedge .. .. \wedge \omega \sb n)((e \sb 1,.. .., e \sb n))=Det
\pmatrix{\omega \sb 1 (e \sb 1 )&.. ..&\omega \sb n (e )\cr
         \vdots &.. ..&\vdots\cr
         \omega  \sb 1 (e \sb n)&.. ..&\omega \sb n (e \sb n)\cr}
= Det I = 1
$$

Se deduce que todo se cumple si $\omega \sb j (e \sb k)= \delta \sb{jk} = 1 $ si $i=j$ y $0$ si no. Se dice que  $\omega \sb j = (e \sb j) \sp * $,  y a $(e \sb j) \sp *$ se le llama el dual de $e\sb j$.

\section{FORMAS DIFERENCIALES}

Las formas operan sobre plps macroscópicos que residen en espacios planos. Necesitamos extender su formalismo al mundo microscópico, infinitesimal. Pero hay que hacerlo de tal manera que su generalización a espacios curvos o variedades sea inmediata y natural, pues es requisito de toda variedad que sea plana cuando se le mira infinitesimalmente. La idea, entonces, es que una k-forma diferencial es una k-forma que opera sobre k-plps infinitesimales que son tangentes a la variedad o espacio curvo. Históricamente no resultó tan fácil lograrlo y en primera instancia es algo extraño. Pero eso se convierte en algo sencillo si uno nota que en $\Re \sp n$  hay un isomorfismo entre vectores y operadores de derivada direccional.

Dichas \index{derivadas direccionales} \textbf{derivadas direccionales} operan por medio del gradiente así:

 $D\sb{\vec u} (f)=\vec u \cdot \nabla f  $.

Por ejemplo, si $\vec u = [-2,4]$, entonces,

$D\sb{[-2,4]} (f)=  [-2,4]\cdot\nabla f =  -2\partial f/\partial x + 4\partial f/ \partial y $

de tal manera que el vector $\vec{u} = [-2,4]$ define un operador de derivada direccional dado por

$ D\sb{\vec u} (\cdot ) = -2\partial/\partial x + 4\partial / \partial y $

el cual al operar sobre una función escalar $f$ da:

$ D\sb{\vec u} (f) = -2\partial f/\partial x + 4\partial f/ \partial y$

\bigskip

\addtocounter{ecu}{1}

\textit{Teorema(\theecu). En $\Re \sp n$
 existe un \index{isomorfismo}  \textbf{isomorfismo} natural entre vectores y operadores de derivada direccional dado por:}

$\vec u \leftrightarrow D\sb{\vec u}(\cdot )$

\bigskip

Demostración:

  Al vector $(a_1, .. .., a_n)$ le hacemos corresponder biunívocamente el operador de derivada direccional dado por
$a_1\partial /\partial x_1 + .. .. + a_n\partial / \partial x_n$, el cual opera sobre una función escalar $f$ produciendo
$a_1\partial f/\partial x_1 + .. .. + a_n\partial f/ \partial x_n$.

 Las evaluaciones de las derivadas pueden hacerse en el origen, pero el isomorfismo sigue válido para cualquier punto. Es decir, para cada punto de $\Re^n$ tenemos un isomorfismo. Si no se habla de un punto en particular, se sobreentiende que se está operando sobre un punto dado pero inespecífico.

\bigskip

No fue sino hacia 1930 que se entendió la importancia de este isomorfismo para generalizar el cálculo a espacios curvos o variedades (estamos definiendo informalmente una variedad como un espacio curvo pero suave, como una circunferencia o una esfera). Toda la idea se reduce a lo siguiente: un espacio curvo pero suave es infinitesimalmente plano. Es decir,  un espacio curvo es localmente lo mismo que su espacio tangente. Ahora bien, lo que estamos haciendo para $\Re^n$ se extiende a los espacios curvos de manera natural. A los espacios curvos los llamados variedades cuando reúnen las condiciones suficientes para hacer cálculo. La definición formal la dejamos para después. Definamos por ahora el espacio tangente a $Re^n$ en el punto $P$. Intuitivamente, el espacio tangente a $\Re^n$ es el mismo $\Re^n$:

\bigskip

\addtocounter{ecu}{1}

\textit{Teorema(\theecu). El \index{espacio tangente} \textbf{espacio tangente} a  $\Re \sp n$ en el punto $P$ es el conjunto de todos los vectores $v = v_1 \partial /\partial x_1 + .. .. + v_n\partial / \partial x_n$ tal que a cada función escalar $f$ le hace corresponder el número real dado por $v_1\partial f/\partial x_1 + .. .. + v_n\partial f/ \partial x_n$, donde las evaluaciones se hacen en el punto $P$.}

\bigskip

Resulta muy claro que de acuerdo a esta nueva definición, el espacio tangente a un plano en un punto dado es el mismo plano y que el espacio tangente a $\Re^n$ en un punto dado es él mismo.

Lo que estamos haciendo para el plano es válido para variedades, aunque no lo haremos explícito. En particular, el espacio tangente a una circunferencia será un espacio que para todos los efectos es una línea, la línea tangente a la circunferencia en el punto dado.

Podemos  redefinir ahora el prototipo de las formas diferenciales que es la diferencial $d$. Hagámoslo paso a paso:

En cálculo vectorial sobre $\Re^n$, la diferencial evaluada en un punto $P$ mide el cambio infinitesimal de una función escalar $f$ debido a una variación infinitesimal de $P$:

$df= \nabla f \cdot \vec{dr} = (\partial f/\partial x\sb1) dx\sb 1 +.. .. + (\partial f/\partial x\sb n) dx\sb n $

donde las derivadas se evalúan en el punto $P$.

\

Releemos esa definición en términos de plps: $d$ opera sobre una función que parte de $\Re \sp n$, toma valores reales $f:\Re \sp n\rightarrow \Re$,  y da una forma diferencial $df$ que opera sobre un 1-plp infinitesimal $\vec{dr}= [dx \sb 1,.. ..,dx\sb n]$ produciendo un número real infinitesimal:

$df(\vec{dr})= \nabla f \cdot \vec{dr} = (\partial f/\partial x\sb1) dx\sb 1 +.. .. + (\partial f/\partial x\sb n) dx\sb n $

Observemos que $df$, la diferencial de $f$, se descompuso en una combinación lineal de las  $dx\sb i$ que son las diferenciales de las coordenadas $x\sb i$. Por lo tanto, podremos reinterpretar $df$ siempre y cuando reinterpretemos las $dx\sb i$.

\addtocounter{ecu}{1}

\textit{Definición (\theecu): Definimos la \index{forma diferencial} \textbf{forma diferencial} $dx\sb i $ como el operador que al vector $\vec j = \partial /\partial x\sb j$ le hace corresponder}

$dx\sb i (\partial /\partial x \sb j) = \delta \sp i \sb j$

\bigskip

El vector $\vec j = \partial /\partial x\sb j$ es, por supuesto, nuestra representación del vector que tiene un uno en la coordenada $j$ y cero en las demás, $(0,.. .., 1,.. ..,0)$, el cual define la dirección del eje número $j$. Es decir, $\vec j = \partial /\partial x\sb j$ es el operador de derivada direccional en la dirección del eje $j$.

A partir de la anterior definición, todo se concatena por \index{multilinealidad alternada} \textbf{multilinealidad alternada}. Por ejemplo, tomemos  $\Re \sp 2$ y sobre ella consideremos el efecto de la 1-forma diferencial $\omega = 2x d x+ 2ydy= d(r\sp 2=x\sp 2 + y \sp 2)$ sobre el vector $\vec u = [-1,-1]$. Como la forma diferencial no es constante, hay que fijar un punto de operación, sea $P=(2,2)$. En primer término, al vector $\vec u = [-1,-1]$ le corresponde el operador,  $\vec u = -\partial /\partial x-\partial /\partial y$. Ahora operemos:

$\omega (\vec u)= (2x d x+ 2ydy)(-\partial /\partial x-\partial /\partial y)$
$=2x d x (-\partial /\partial x-\partial /\partial y)+ 2ydy(-\partial /\partial x-\partial /\partial y)  $

=$2x d x (-\partial /\partial x)+ 2x d x (-\partial /\partial y)+ 2ydy(-\partial /\partial x)+ 2ydy(-\partial /\partial y)$

=$2x d x (-\partial /\partial x)+ 2ydy(-\partial /\partial y)$

$=-2x - 2y$

Ahora evaluamos en el punto $P=(2,2)$, y $\omega (\vec u)=2(2)(-1)+2(2)(-1)= -4 -4= -8 $

Obsérvese que el resultado final se lee $d(r\sp 2)(\vec u) = -2x - 2y $. Este resultado también puede calcularse simplemente como $(\vec u)(r\sp 2)$. En efecto, puesto que al vector $\vec u$ le corresponde el operador

$\vec u = -\partial /\partial x-\partial /\partial y$

y teniendo en cuenta que $r\sp 2 = x\sp 2 + y \sp 2$ resulta

$(\vec u)(r\sp 2)=  -\partial(r\sp 2) /\partial x-\partial (r\sp 2)/\partial y= -\partial(x\sp 2 + y \sp 2) /\partial x-\partial (x\sp 2 + y \sp 2)/\partial y=-2x - 2y$

El resultado anterior es un ejemplo de un caso general de \index{dualidad} \textbf{dualidad}: una forma diferencial determinada por una función que toma valores reales, $f$, que se aplica sobre un vector, $\vec u$, produce el mismo resultado que el vector, considerado como operador de derivación direccional, aplicado  sobre la función:

$df(\vec u) = \vec u (f)$

lo cual se demuestra generalizando el siguiente razonamiento. Si notamos $\vec u = u\sb 1 \partial/\partial x + u\sb 2 \partial / \partial y$, donde $u_1$, $u_2$ son números fijos, y teniendo en cuenta que $dx(\partial /\partial x) = 1,dx(\partial /\partial y) = 0,dy(\partial /\partial x) = 0,dy(\partial /\partial y) = 1 $ resulta que

$df(\vec u)= [(\partial f/ \partial x)dx  +  (\partial f/ \partial y)dy] ((u\sb 1 \partial/\partial x + u\sb 2 \partial / \partial y) )$

 $=u\sb 1(\partial f/ \partial x) + u\sb 2   (\partial f/ \partial x)$

Por otro lado, $\vec u(f)$ es:

$\vec u(f) = [u\sb 1(\partial / \partial x) + u\sb 2(\partial / \partial x)](f)$

que coincide con la expresión en la cual se opera la forma diferencial $df$ sobre el vector $\vec u$. Concluimos:

\addtocounter{ecu}{1}

\textit{Teorema (\theecu) : Sobre $\Re\sp n$ y para una función escalar $f$ se tiene que}

$$df(\vec u) = \vec u (f)$$

Esa  expresión de dualidad es la que permitirá definir a $df$ sobre variedades.

\bigskip

En general, toda la estructura algebraica de las formas y de los plps  en $\Re \sp n$
se heredan a las formas diferenciales y a los vectores del espacio tangente. En particular, al intercambiar dos elementos en un producto wedge de formas diferenciales, se cambia de signo, lo cual implica que cuando hay dos elementos repetidos, se produce aniquilación. Los cálculos pueden hacerse de forma automática, como lo indica el ejemplo siguiente. Sea la 2-forma diferencial sobre $\Re \sp 2$

$\omega = x  dx  \wedge  dy -y dy \wedge  dx  $ la cual debe  actuar sobre un 2-plp compuesto por:

$\xi = [-1,0]= -\partial/\partial x$

$\eta = [0,1]= \partial/\partial y$

y que opera en el punto $P=[1,2]$

Qué es : $\omega(\xi, \eta)$ calculada en $P$?

$\omega  = x  dx  \wedge  dy + y dx  \wedge  dy $

la cual evaluada en $P=[1,2]$ da $\omega= 3 dx  \wedge  dy$. Por tanto:

$$\omega \ (\xi , \eta) =3 Det
\pmatrix{dx  (\xi) &dy (\xi)\cr
         dx  (\eta)&dy(\eta)\cr}
= 3Det\pmatrix{-1 &0\cr
         0&1\cr} = -3 $$

pues, por ejemplo, $ dx  (\xi) = dx (-\partial/\partial x) = -1$.

\

Similarmente, si

$\omega = (x  + y)dx  \wedge  dy  $

$\xi = [1,0]= \partial/\partial x$

$\eta = [1,1]= \partial/\partial x+ \partial/\partial y$

 $P=[1,1]$

tenemos que en P:

$\omega = 2 dx  \wedge  dy $

$dx (\xi)= 1, dy(\xi)= 0,dx (\eta)= 1,dy(\eta)= 1, $

$$2 dx  \wedge  dy\ (\xi , \eta) =2 Det
\pmatrix{dx  (\xi) &dy (\xi)\cr
         dx  (\eta)&dy (\eta)\cr}
= 2 Det\pmatrix{1 &0\cr
         1&1\cr} = 2 $$

No olvidemos que estamos generalizando a espacios curvos las estructuras naturales que permiten el cálculo en $\Re^n$. Por eso, tratemos de recobrar la antigua diferencial a partir  de las nuevas definiciones. Para fijar ideas, pensemos en el plano. Habíamos visto que, en general, para $u = (u_1,u_2)$:

$df(\vec u)= [(\partial f/ \partial x)dx  +  (\partial f/ \partial y)dy] ((u\sb 1 \partial/\partial x + u\sb 2 \partial / \partial y) )$

 $=u\sb 1(\partial f/ \partial x) + u\sb 2   (\partial f/ \partial x)$

Ahora, como $u$ podemos tomar un vector infinitesimal: $u = (\Delta u_1, \Delta u_2)$, es decir, $u = (\Delta u_1 \partial f/ \partial x, \Delta u_2 \partial f/ \partial y)$. En ese caso, la nueva definición da

$df(\vec u)  =\Delta u\sb 1(\partial f/ \partial x) + \Delta u\sb 2   (\partial f/ \partial x)$

lo cual corresponde exactamente a la antigua versión de $df$.

\section{DERIVACION}

Una  forma diferencial puede derivarse de tal manera que el resultado sea una forma diferencial. A la derivada de una forma diferencial se le llama oficialmente \index{derivada exterior}  \textbf{derivada exterior} (debido a que hay varios tipos de derivación), pero aquí sólo se llamará derivada, pues no será posible confundirse. Los plps no pueden derivarse. Como tal, cualquier cosa que tenga derivación tiene que encajar dentro de una forma diferencial. Por eso, las leyes de la naturaleza que antes venían en ecuaciones en derivadas parciales, ahora podrán venir en formas diferenciales.

\

Como hemos ligado la diferencial a derivadas parciales, al derivar una diferencial se incrementa el orden de derivación. La contabilidad de tales cambios se lleva automáticamente dictaminando el grado de una forma diferencial.

\

\textit{Definimos \index{forma} \textbf{forma} mediante el siguiente algoritmo:}

1. Una función $f$ que toma valores reales es una \textbf{0-forma}. Su diferencial es la 1-forma:

\addtocounter{ecu}{1}

$df = (\partial f/\partial x\sb1) dx\sb 1 +.. .. + (\partial f/\partial x\sb n) dx\sb n $

2. Una \textbf{1-forma} opera sobre un vector u operador de  derivada direccional. Un vector es un 1-plp, representa un estímulo dado a la naturaleza  y una 1-forma representa la forma de operar de la naturaleza, que aplicada sobre el estímulo recibido produce un resultado, en este caso un número real. 

3.La derivada de una 1-forma diferencial es una 2-forma diferencial. La derivación es muy natural e involucra el producto wedge:

Con miras a un tratamiento relativista del campo electromagnético, consideremos la 1-forma:

$A=A(x)\sb o dt +A(x\sb 1) dx\sb 1 +A(x)\sb 2 dx\sb 2 +A(x)\sb 3 dx\sb 3 $

donde cada $A\sb i$ es una función común y corriente, una 0-forma, que tiene como diferencial una 1-forma. La diferencial de $A$ da una 2-forma y se calcula así:

$dA=dA\sb o \wedge dx\sb o +dA\sb 1 \wedge dx\sb 1 +dA\sb 2 \wedge dx\sb 2 +dA\sb 3 \wedge dx\sb 3$

Quizá sea conveniente aclarar que una expresión del tipo   $dA\sb 1 \wedge dx\sb 1$ se interpreta como $ d(A\sb 1) \wedge dx\sb 1$, lo cual es permitido, pues la diferencial de una función $A\sb 1$  es una 1-forma, la cual se puede multiplicar exteriormente por la 1-forma $dx$ para que de una 2-forma, como sigue:

Teniendo en cuenta la expansión de cada diferencial, tenemos que:

$dA=[(\partial A\sb o /\partial x\sb o) dx\sb o +.. .. + (\partial A \sb o /\partial x\sb 3) dx\sb 3] \wedge dx\sb o + [ (\partial A\sb 1 /\partial x\sb o) dx\sb o +.. .. + (\partial A\sb 1  /\partial x\sb 3) dx\sb 3] \wedge dx\sb 1 +.. ..$

$ = (\partial A\sb o /\partial x\sb o) dx\sb o\wedge dx\sb o + (\partial A\sb o /\partial x\sb 1) dx\sb 1\wedge dx\sb o + (\partial A\sb o /\partial x\sb 2) dx\sb 2\wedge dx\sb o + (\partial A\sb o /\partial x\sb 3) dx\sb 3\wedge dx\sb o + (\partial A \sb 1/\partial x\sb o) dx\sb o\wedge dx_1 + (\partial A \sb 1/\partial x_1) dx_1\wedge dx_1 + (\partial A \sb 1/\partial x\sb 2) dx\sb 2\wedge dx_1 + (\partial A \sb 1/\partial x\sb 3) dx\sb 3\wedge dx_1 +.. ..  $

teniendo en cuenta que, cuando en un producto exterior hay repeticiones, resulta una aniquilación, y reescribiendo los ceros correspondientes, esta expresión tiene una representación matricial inmediata:

$$\bordermatrix{&dx\sb o &dx\sb 1&dx\sb 2&dx\sb 3\cr
dx\sb o& 0&\partial A \sb 1/\partial x\sb o
&  \partial A\sb 2 /\partial x\sb o  &
  \partial A\sb 3 /\partial x\sb o\cr
dx\sb 1& \partial A\sb o /\partial x\sb 1&0& \partial A\sb 2 /\partial x\sb 1& \partial A\sb 3 /\partial x\sb 1 \cr
dx\sb 2 &\partial A\sb o /\partial x\sb 2&\partial A\sb 1 /\partial x\sb 2&0&   \partial A\sb 3 /\partial x\sb 2 \cr
    dx\sb 3 & \partial A\sb o /\partial x\sb 3&\partial A\sb 1 /\partial x\sb 3&\partial A\sb 2 /\partial x\sb 3&0\cr}
$$

Podemos reordenar, recordando la anticonmutatividad del producto exterior:

$ dA= 0 - (\partial A\sb o /\partial x\sb 1) dx\sb o\wedge dx\sb 1 - (\partial A\sb o /\partial x\sb 2) dx\sb o\wedge dx\sb 2 - (\partial A\sb o /\partial x\sb 3) dx\sb o\wedge dx\sb 3 + (\partial A\sb 1 /\partial x\sb o) dx\sb o\wedge dx\sb 1 + 0 - (\partial A \sb 1/\partial x\sb 2) dx\sb 1\wedge dx\sb 2 - (\partial A \sb 1/\partial x\sb 3) dx\sb 1\wedge dx\sb 3 +.. ..  $

Reorganizando tenemos:

$ dA= [ \partial A\sb 1 /\partial x\sb o- \partial A\sb o /\partial x\sb 1] dx\sb o\wedge dx\sb 1+ $
$ [ \partial A\sb 2 /\partial x\sb o- \partial A\sb o /\partial x\sb 2] dx\sb o\wedge dx\sb 2+ $
$ [ \partial A\sb 3 /\partial x\sb o- \partial A\sb o /\partial x\sb 3] dx\sb o\wedge dx\sb 3+ $
$ [ \partial A\sb 2 /\partial x\sb 1- \partial A\sb 1 /\partial x\sb 2] dx\sb 1\wedge dx\sb 2+ $
$ [ \partial A\sb 3 /\partial x\sb 1- \partial A \sb 1/\partial x\sb 3] dx\sb 1\wedge dx\sb 3+ $
$ [ \partial A\sb 3 /\partial x\sb 2- \partial A\sb 2 /\partial x\sb 3] dx\sb 2\wedge dx\sb 3 $

Esta expresión nos permite representar  la diferencial de una 1-forma como un arreglo que no es una matriz, sino que es una 2-forma que opera sobre 2-plps:

$$\bordermatrix{&dx\sb o &dx\sb 1&dx\sb 2&dx\sb 3\cr
dx\sb o& 0&\partial A\sb 1 /\partial x\sb o- \partial A\sb o /\partial x\sb 1
&  \partial A\sb 2 /\partial x\sb o- \partial A\sb o /\partial x\sb 2 &
  \partial A\sb 3 /\partial x\sb o- \partial A\sb o /\partial x\sb 3\cr
dx\sb 1& &0& \partial A\sb 2 /\partial x\sb 1- \partial A\sb 1 /\partial x\sb 2& \partial A\sb 3 /\partial x\sb 1- \partial A\sb 1 /\partial x\sb 3\cr
dx\sb 2 &&&0&   \partial A\sb 3 /\partial x\sb 2- \partial A\sb 2 /\partial x\sb 3\cr
    dx\sb 3 & &&&0\cr}
$$

Todo esto se escribe  sucintamente como sigue:

\addtocounter{ecu}{1}

\textit{Teorema (\theecu). La diferencial de una 1-forma diferencial $A=\sum \sb i A\sb i dx\sb i$ es la 2-forma}

$dA = \sum \sb{i<j} (\partial A\sb j/\partial x \sb i -\partial A\sb i/\partial x \sb j )dx\sb i \wedge dx\sb j$

\bigskip

Observemos que esta 2-forma opera sobre 2-plps y que si intercambiamos el orden de sus aristas, el valor resultante cambia de signo. Tenemos, pues, que la diferencial de una 1-forma es automáticamente antisimétrica.

Observación: tanto los plps como las formas diferenciales tienen reglas precisas de transformación ante un cambio de base. Cuando decimos que la diferencial de una forma es otra forma, estamos diciendo, que la diferencial es compatible con las reglas de transformación ante cambios de base. Haber logrado eso, es el mérito de Grassmann.

\bigskip

Al ir derivando, vamos aumentando el grado de la forma diferencial resultante. Es usual anotar el grado de una forma diferencial como un superíndice. Por ejemplo, $\alpha^5$ denota una forma diferencial de grado 5, al estilo $dx_1\wedge dx_2 \wedge dx_3 \wedge dx_4 \wedge dx_5$ o también $dx_2\wedge dx_4 \wedge dx_7 \wedge dx_{11} \wedge dx_{16}$. No siempre una forma diferencial viene de una derivación, pero cuando eso pasa se llama exacta. Cuando la diferencial de una forma da cero, se llama cerrada.

\bigskip

\addtocounter{ecu}{1}

\textit{Teorema (\theecu): La derivación extiende la regla del producto de las derivadas ordinarias para 0-formas o funciones: si $\alpha \sp p$ es una p-forma y $\beta \sp q$ es una q-forma, entonces
$d(\alpha \sp p \wedge \beta \sp q) = d\alpha \sp p \wedge \beta\sp q + (-1)\sp p \alpha\sp p \wedge d\beta\sp q$}

\bigskip

Demostración para 1-formas:  Como hemos visto, para 0-formas o funciones  $d$ coincide con la diferenciación ordinaria. El producto de dos 1-formas es una 2-forma y su derivada es una 3-forma. Si

$\alpha=\sum \sb i a\sb i dx\sb i$ es una 1-forma, tenemos

$d\alpha = \sum \sb{i<j} (\partial a\sb j/\partial x \sb i -\partial a\sb i/\partial x \sb j )dx\sb i \wedge dx\sb j$

y si

$\beta=\sum \sb m b\sb m dx\sb m$ es la otra 1-forma

$d\beta = \sum \sb{m<n} (\partial b\sb n/\partial x \sb m -\partial b\sb m/\partial x \sb n )dx\sb m \wedge dx\sb n$

Por lo tanto, de acuerdo con la regla que tenemos que probar, la derivación del producto sería:

$d(\alpha  \wedge \beta ) = d\alpha  \wedge \beta + (-1) \alpha \wedge d\beta$ que al reemplazar daría:

$d(\alpha  \wedge \beta ) = (\sum \sb{i<j} (\partial a\sb j/\partial x \sb i -\partial a\sb i/\partial x \sb j )dx\sb i \wedge dx\sb j)  \wedge (\sum \sb m b\sb m dx\sb m) + (-1) (\sum \sb i a\sb i dx\sb i) \wedge (\sum \sb{m<n} (\partial b\sb n/\partial x \sb m -\partial b\sb m/\partial x \sb n )dx\sb m \wedge dx\sb n)$

Por otro lado:

$\alpha  \wedge \beta = (\sum \sb i a\sb i dx\sb i)\wedge (\sum \sb m b\sb m dx\sb m)$
$=\sum \sb {im} a\sb i b\sb m dx\sb i\wedge dx\sb m$

Por tanto, y aplicando la regla para 0-formas o funciones ordinarias:

$d(\alpha  \wedge \beta )= d(\sum \sb {im} a\sb i b\sb m dx\sb i\wedge dx\sb m)$

$= \sum \sb {i,m} d(a\sb i b\sb m) dx\sb i\wedge dx\sb m$

$=\sum \sb {i,m} [d(a\sb i) b\sb m +a\sb i (db \sb m)] dx\sb i\wedge dx\sb m$

$=\sum \sb {i,m} [(\sum \sb k (\partial a\sb i /\partial x\sb k)dx\sb k) b\sb m +a\sb i (\sum \sb k  (\partial b\sb m/\partial x\sb k)dx\sb k)]\wedge dx\sb i\wedge dx\sb m$

$=\sum \sb {i,m,k} [(\partial a\sb i /\partial x\sb k)dx\sb k b\sb m \wedge dx\sb i\wedge dx\sb m+a\sb i (\partial b\sb m/\partial x\sb kdx\sb k)\wedge dx\sb i\wedge dx\sb m]$

$=\sum \sb {i,m,k} [(\partial a\sb i /\partial x\sb k)dx\sb k  \wedge dx\sb i\wedge b\sb mdx\sb m+a\sb i (\partial b\sb m/\partial x\sb kdx\sb k)\wedge dx\sb i\wedge dx\sb m]$

$=\sum \sb {i,m,k} [(\partial a\sb i /\partial x\sb k)dx\sb k  \wedge dx\sb i\wedge b\sb mdx\sb m+(-1)a\sb i dx\sb i\wedge (\partial b\sb m/\partial x\sb kdx\sb k)\wedge dx\sb m]$

$=\sum \sb {i,k} (\partial a\sb i /\partial x\sb k)dx\sb k  \wedge dx\sb i\wedge (\sum \sb m b\sb mdx\sb m)+(-1)(\sum \sb i a\sb i dx\sb i)\wedge (\sum \sb{m,k}(\partial b\sb m/\partial x\sb k)dx\sb k)\wedge dx\sb m$

Ahora aniquilamos todos los productos repetidos y ordenamos en forma creciente los subíndices:

$=\sum \sb {k<i} (\partial a\sb i /\partial x\sb k-\partial a\sb k /\partial x\sb i)dx\sb k  \wedge dx\sb i\wedge (\sum \sb m b\sb mdx\sb m)+(-1)(\sum \sb i a\sb i dx\sb i)\wedge (\sum \sb{k<m}(\partial b\sb m/\partial x\sb k-\partial b\sb k/\partial x\sb m)dx\sb k)\wedge dx\sb m$

lo cual coincide, aparte de subíndices mudos, con la expresión predicha por el teorema.

\bigskip

Tenemos ahora un gran teorema que calcula $d(dA)$ y que permite la introducción de la cohomología.

\addtocounter{ecu}{1}

\textbf{Teorema de la \index{cohomología} cohomología} (\theecu): $d\sp 2 = 0$.

\bigskip

La demostración se basa en utilizar la igualdad de las derivadas mixtas y en alegar duplicación de términos con signos contrarios en todos lados. En nuestro caso podemos ver cómo funciona la mutua aniquilación de términos al ver la diferencial de dos de los 4 términos de $dA$.

Si  $A=\sum \sb i A\sb i dx\sb i$ su diferencial es

$dA = \sum \sb{i<j} (\partial A\sb j/\partial x \sb i -\partial A\sb i/\partial x \sb j )dx\sb i \wedge dx\sb j$

Calculemos la diferencial del primer término de $dA$ :

$ d[ \partial A_1 /\partial x\sb o- \partial A\sb o /\partial x_1] dx\sb o\wedge dx_1=$

$\partial/\partial x\sb o( \partial A_1 /\partial x\sb o- \partial A\sb o /\partial x_1) dx\sb o\wedge dx\sb o\wedge dx_1$
$+\partial/\partial x_1( \partial A_1 /\partial x\sb o- \partial A\sb o /\partial x_1) dx_1\wedge dx\sb o\wedge dx_1$
$+\partial/\partial x\sb 2( \partial A_1 /\partial x\sb o- \partial A\sb o /\partial x_1) dx\sb 2\wedge dx\sb o\wedge dx_1$
$+\partial/\partial x\sb 3( \partial A_1 /\partial x\sb o- \partial A\sb o /\partial x_1) dx\sb 3\wedge dx\sb o\wedge dx_1$

$=0 + 0 + [\frac{\rm \partial \sp 2 A_1}{\rm \partial x \sb 2 \partial x \sb o}
-\frac{\rm \partial \sp 2 A\sb o}{\rm \partial x \sb 2 \partial x_1 }]dx\sb o\wedge dx_1 \wedge dx\sb 2$
$+[\frac{\rm \partial \sp 2 A_1}{\rm \partial x \sb 3 \partial x \sb o}
-\frac{\rm \partial \sp 2 A\sb o}{\rm \partial x \sb 3 \partial x_1 }]dx\sb o\wedge dx_1\wedge dx\sb 3$

\bigskip
Calculemos ahora la diferencial del segundo término de $dA$:

$ d[ \partial A\sb 2 /\partial x_0- \partial A_0 /\partial x\sb 2] dx_0\wedge dx\sb 2=$
$\partial/\partial x\sb o( \partial A\sb 2 /\partial x_0- \partial A_0 /\partial x\sb 2) dx\sb o\wedge dx_0 \wedge dx\sb 2$
$+\partial/\partial x_1( \partial A\sb 2 /\partial x_0- \partial A_0 /\partial x\sb 2) dx_1\wedge dx_0 \wedge dx\sb 2$
$+\partial/\partial x\sb 2( \partial A\sb 2 /\partial x_0- \partial A_0 /\partial x\sb 2) dx\sb 2\wedge dx_0\wedge dx\sb 2$
$+\partial/\partial x\sb 3( \partial A\sb 2 /\partial x_0- \partial A_0 /\partial x_2) dx\sb 3\wedge dx_0\wedge dx\sb 2$

$=0 +  [-\frac{\rm \partial \sp 2 A\sb 2}{\rm \partial x \sb 1 \partial x_0 }
+\frac{\rm \partial \sp 2 A_0}{\rm \partial x \sb 1 \partial x \sb 2}]dx\sb o\wedge dx_1\wedge dx\sb 2$
$+[\frac{\rm \partial \sp 2 A\sb 2}{\rm \partial x \sb 3 \partial x_0 }
-\frac{\rm \partial \sp 2 A_0}{\rm \partial x \sb 3 \partial x \sb 2}]dx_0 \wedge dx\sb 2\wedge dx\sb 3$

Podemos ver que el término $\frac{\rm \partial \sp 2 A_0}{\rm \partial x \sb 1 \partial x \sb 2}dx\sb o\wedge dx_1 \wedge dx\sb 2$ aparece con signos diferentes en las dos diferenciales anotadas. Al calcular los demás términos, se obtiene la aniquilación total.

\section{PRODUCTO INTERIOR}

Por medio de la derivación podemos subir el grado de las formas diferenciales. Veamos ahora un mecanismo para bajarlo. Se trata simplemente de tomar una p-forma que funciona sobre p-plps y adueñarse de la primera arista poniendo un vector fijo para siempre:

\addtocounter{ecu}{1}

\textit{Def(\theecu): Si
$\vec v$ es un vector y $\alpha $ es una p-forma, definimos el \index{producto interior} producto interior de la forma y el vector como la (p-1)-forma $i\sb{\vec v} \alpha $ que funciona así: $i\sb{\vec v} \alpha $ es cero si $\alpha$ es una  función o  0-forma. Da $i\sb{\vec v} \alpha =\alpha(\vec v)$ si $\alpha$ es una 1-forma. En los demás casos, para una p-forma $\alpha \sp p$ con $p>1$,    $ i\sb{\vec v} \alpha^p $ es una forma de grado $p-1$ que opera como}

 $(i\sb{\vec v} \alpha^p )\sp{p-1} (\vec w \sb 2,.. ..,\vec w\sb p) = \alpha \sp p (\vec v,\vec w\sb 2,.. ..,\vec w\sb p)$.

\bigskip

Al producto interno así definido se le llama antiderivación por bajar el grado de la forma y por cumplir el siguiente teorema:

\addtocounter{ecu}{1}

\textit{Teorema} (\theecu). $i\sb{\vec v}(\alpha \sp p \wedge \beta \sp q)=[i\sb{\vec v}\alpha \sp p] \wedge \beta \sp q)+ (-1)\sp p \alpha  \sp p \wedge [i\sb{\vec v} \beta \sp q]$

\bigskip

Es importante tener en cuenta  que el producto interno es un operador lineal en  sus dos entradas, tanto  la que se escribe como subíndice como la que se escribe como operando.

\section{REESCRITURA}

Veamos cómo se reescriben en el lenguaje de las formas diferenciales algunos resultados conocidos del cálculo en tres variables.

Repitamos una regla pasada que hay que tener en claro: las formas son las únicas entidades que se derivan, antiderivan e integran. No es una exageración decir que las formas se hicieron para integrarlas y de hecho, la diferenciación de formas no se puede entender a menos que se enuncie su relación con integración, la cual se ha definido para que cumpla reglas muy naturales. En efecto:

Cuando uno desea una integral de línea, como hallar la masa de un alambre o el trabajo hecho al recorrer una cierta trayectoria con un peso a cuestas, lo que se integra es una 1-forma. Similarmente, si se trata de una integral de superficie, lo que se integra es una 2-forma, la cual opera sobre 2-plps que son los que generan las superficies. Y, por supuesto, al querer hallar volúmenes, masas encerradas por una superficie, lo que se integra es una 3-forma.

Por tanto, si uno tiene que reescribir una expresión vectorial en el presente lenguaje, uno primero debe entender en qué contexto va a entrar. Ejemplo: si $\vec F$ es una campo de fuerzas tangencial a una curva, y si uno desea meterla en una integral de línea, para calcular un trabajo, entonces $\vec F$ entra como una 1-forma:

\addtocounter{ecu}{1}

$$\vec F = (F\sb 1,F\sb 2,F\sb 3 ) \leftrightarrow \alpha = F\sb 1 dx\sp 1 +F\sb 2 dx\sp 2 +F\sb 3 dx\sp 3 \eqno{(\theecu)} $$

Pero si el campo de fuerzas es perpendicular a una superficie, y se trata de una integral de superficie, entonces a dicho campo le corresponde una 2-forma, puesto que las superficies se expanden con 2-plps. La 2-forma  se construye teniendo en cuenta que para un punto dado, el vector $\vec F$ sirve para completar un 3-plp al cual se le calcula su volumen. Lo que tenemos en mente es que al campo se le calcula el  flujo sobre una superficie infinitesimal, el cual resulta del volumen expandido por un 2-plp sobre la superficie y el vector definido por el campo.

De eso se infiere que la 2-forma asociada a $ \vec F$ se construye así:

Primero, la forma de volumen es $vol\sp 3= du_1\wedge du_2 \wedge du_3 $. Segundo: a un vector, digamos $ \vec v =(v\sb 1, v\sb 2, v\sb 3)$, se le asocia la 2-forma $\beta \sp 2 = i\sb{\vec v}vol$. Para saber cuánto es eso, usamos el hecho de que el producto interno sobre un producto wedge opera como una derivación. En el presente caso tenemos un triple producto, por lo que en el caso especial  $\vec v= \partial \sb j = \partial/\partial x\sb j$ y notando la forma $ vol \sp 3$ como $\alpha $, entonces

 $i\sb{\vec v} \alpha = i\sb{(\partial \sb j)}(du\sp 1)du\sp 2\wedge du\sp 3-du\sp 1\wedge i\sb{(\partial \sb j)}(du\sp 2)\wedge du\sp 3 + du\sp 1\wedge du\sp 2\wedge i\sb{(\partial \sb j)}(du\sp 3) $

$= du\sp 1(\partial \sb j)du\sp 2\wedge du\sp 3-du\sp 1\wedge du\sp 2(\partial \sb j)\wedge du\sp 3 + du\sp 1\wedge du\sp 2\wedge du\sp 3(\partial \sb j) $

$= \delta \sp 1\sb j du\sp 2\wedge du\sp 3-\delta \sp 2\sb j du\sp 1\wedge  du\sp 3 + \delta \sp 3\sb j du\sp 1\wedge du\sp 2 $

Si $j=1$ esto se convierte en:

$i\sb{\vec v} \alpha =\delta \sp 1\sb 1 du\sp 2\wedge du\sp 3-\delta \sp 2\sb 1 du\sp 1\wedge  du\sp 3 + \delta \sp 3\sb 1 du\sp 1\wedge du\sp 2 = du\sp 2\wedge du\sp 3$

Por tanto, si $\vec v = \sum v\sp i \partial \sb i$, y la 3-forma sigue siendo la del volumen, entonces,

$i\sb{\vec v} \alpha = v\sp 1 i\sb{\partial/\partial u\sp 1} (\alpha) + v\sp 2 i\sb{\partial/\partial u\sp 2} (\alpha) + v\sp 3 i\sb{\partial/\partial u\sp 3} (\alpha)$

$= v \sp 1 du\sp 2\wedge du\sp 3-v\sp 2 du\sp 1\wedge  du\sp 3 + v\sp 3 du\sp 1\wedge du\sp 2 $

$= v \sp 1 du\sp 2\wedge du\sp 3+v\sp 2 du\sp 3\wedge  du\sp 1 + v\sp 3 du\sp 1\wedge du\sp 2  $

Por tanto, al vector $\vec B =(B\sb 1, B\sb 2, B\sb 3)$, pensemos en el campo magnético, le corresponde la 2-forma:

\addtocounter{ecu}{1}

$$\vec B =(B\sb 1, B\sb 2, B\sb 3) \leftrightarrow \beta \sp 2 = B\sb 1 dy \wedge dz +B\sb 2 dz \wedge dx + B\sb 3 dx \wedge dy \eqno{(\theecu)} $$

\addtocounter{ecu}{1}

Guardemos en mente que:

A un vector se le puede representar como una 1-forma o como una 2-forma y que eso depende del contexto. Como uno puede cuadrar el contexto a su antojo, eso es una arbitrariedad. Ahora bien, el propósito de la presente contribución es el de explorar las consecuencias matemáticas del principio del realismo: existe una realidad independiente de la mente, de los métodos de observación, de las formas de descripción, del tratamiento matemático para estudiarlas.  Pues bien, cuando tenemos un vector en mecánica clásica, al pasarlo a las formas sobre el espacio de Minkowski tiene dos opciones: una 1-forma y una 2-forma. Es arbitrario decir que una representación es mejor que la otra. Por lo tanto, si en alguna parte aparece la 1-forma, en otra parte deberá aparecer la 2-forma acompañante. Ninguna de las dos representaciones puede aparecer sin la otra en un conjunto de leyes.

\bigskip

Veamos ahora qué le corresponde a un  \index{producto cruz } \textbf{producto cruz}:

Sea $\vec A =(a\sb 1, a\sb 2, a\sb 3)$ con la 1-forma asociada $\alpha \sp 1= a\sb 1 dx\sb 1 +  a\sb 2 dx\sp 2 + a\sb 3 dx\sp 3 $ y $ \vec B =(b\sb 1, b\sb 2, b\sb 3)$ con 2-forma asociada $\beta\sp 2= b \sb 1 du\sp 2\wedge du\sp 3+b\sb 2 du\sp 3\wedge  du\sp 1 + b\sb 3 du\sp 1\wedge du\sp 2  $

 Entonces,

\addtocounter{ecu}{1}

$$\vec A \times \vec B \leftrightarrow -i\sb{\vec A}\beta\sp 2 \eqno{(\theecu)}$$

Verifiquémoslo. Por un lado:

$$\vec A\times \vec B = det
\pmatrix{i&j&k\cr
         a\sb 1&a \sb 2&a \sb 3\cr
         b\sb 1&b\sb 2&b\sb 3\cr}
=(a\sb 2 b\sb 3- a\sb 3 b \sb2, -a\sb 1 b\sb 3 + a\sb 3 b \sb 1 , a\sb 1 b\sb 2 - a\sb 2 b\sb 1)
$$

pero por otra parte:

$-i\sb{\vec A}\beta\sp 2 $

$= -i\sb{\vec A}(b \sb 1 du\sp 2\wedge du\sp 3+b\sb 2 du\sp 3\wedge  du\sp 1 + b\sb 3 du\sp 1\wedge du\sp 2  )$

$=-[i\sb{\vec A}(b \sb 1 du\sp 2\wedge du\sp 3)+i\sb{\vec A}(b\sb 2 du\sp 3\wedge  du\sp 1) + i\sb{\vec A}(b\sb 3 du\sp 1\wedge du\sp 2)  ]$

$=-[b \sb 1 i\sb{\vec A}(du\sp 2)\wedge du\sp 3-b \sb 1 du\sp 2\wedge i\sb{\vec A}(du\sp 3)+b\sb 2 i\sb{\vec A}(du\sp 3)\wedge  du\sp 1-b\sb 2 du\sp 3\wedge  i\sb{\vec A}(du\sp 1)+b\sb 3 i\sb{\vec A}(du\sp 1)\wedge du\sp 2-b\sb 3 du\sp 1\wedge i\sb{\vec A}(du\sp 2)]$

$=-[b \sb 1 du\sp 2 (\vec A)du\sp 3+b \sb 1 du\sp 2 du\sp 3(\vec A)+b\sb 2 du\sp 3(\vec A) du\sp 1 +b\sb 2 du\sp 3 du\sp 1(\vec A) + b\sb 3 du\sp 1(\vec A) du\sp 2+ b\sb 3 du\sp 1 du\sp 2(\vec A) ]$

$=-[b \sb 1 a\sb 2 du\sp 3-b \sb 1 a \sb 3 du\sp 2 +b\sb 2 a\sb 3 du\sp 1 -b\sb 2 a\sb 1 du\sp 3  + b\sb 3 a\sb 1 du\sp 2- b\sb 3 a\sb 2 du\sp 1  ]$

$=-[b\sb 2 a\sb 3 du\sp 1-b\sb 3 a\sb 2 du\sp 1+ b\sb 3 a\sb 1 du\sp 2-b \sb 1 a \sb 3 du\sp 2+ b \sb 1 a\sb 2 du\sp 3  -b\sb 2 a\sb 1 du\sp 3    ]$

$=[(b\sb 3 a\sb 2-b\sb 2 a\sb 3 )du\sp 1+(b \sb 1 a \sb 3 - b\sb 3 a\sb 1 )du\sp 2+ (b\sb 2 a\sb 1-b \sb 1 a\sb 2) du\sp 3    ]$

a esta 1-forma le corresponde claramente el vector que representa $\vec A \times \vec B$

Otras asociaciones son:

\addtocounter{ecu}{1}

$$\nabla f \leftrightarrow df \eqno{(\theecu)}$$

lo cual dice que el vector gradiente se representa por la  diferencial que es una 1-forma.

\addtocounter{ecu}{1}

Por otra parte, si representamos un campo vectorial $\vec A$ por una 1-forma, $\alpha\sp 1$, entonces el \index{rotacional} \textbf{rotacional} del campo será utilizado en integrales de superficie y por lo tanto debe ser una 2-forma. Debe estar por tanto ligado a la derivación, pues esa operación es la que sube el grado de las formas:

 $$Rot \vec A = \nabla   \times A \leftrightarrow  d\alpha \sp 1 \eqno{(\theecu)}$$

\addtocounter{ecu}{1}

Similarmente, la \index{divergencia} \textbf{divergencia} de un campo se usa en integrales de volumen, por lo tanto, debe ser una 3-forma, que ha de ser el resultado de derivar la dos forma que representa al campo vectorial dado:

$$Div \vec B=\nabla \cdot \vec B \leftrightarrow d\beta \sp 2 \eqno{(\theecu)}$$

Del cálculo conocíamos una ley que decía: la divergencia del rotacional de un campo es cero. En nuestro lenguaje, eso dice simplemente $d^2A = 0$, pues al campo se le representa por una   1-forma, cuyo rotacional se calcula por medio de la diferencial, al cual se le puede sacar la divergencia, que es simplemente la segunda derivación.

\addtocounter{ecu}{1}

Finalmente, el producto punto de dos vectores $\vec B \cdot \vec A$ corresponde al producto interno de la  2-forma $\alpha \sp 1 $ que representa a $\vec A$, con el vector $\vec B$:

$$\vec B \cdot \vec A \leftrightarrow i\sb{\vec B}\alpha \sp 1 \eqno{(\theecu)}$$

Reescribamos ahora el teorema de Stokes de cálculo vectorial: cuando todo está bien definido se tiene:

$\int_S \nabla \times F = \int_{\partial S}  F$

donde $S$ es una superficie de $\Re^3$ suave y orientable cuya frontera $\partial S$ es suave. Ahora la podemos escribir como

$\int_S d F^1 = \int_{\partial S}  F^1$

donde el campo vectorial $F$ se modeló como una forma diferencial de orden uno, que se puede integrar al lado derecho sobre un camino cerrado, o bien,  cuya diferencial da una forma de grado 2, en el lado izquierdo, la cual opera sobre 2-plps de una superficie para dar una integración sobre una área.

\bigskip

Escribamos ahora el \index{teorema de Gauss} \textbf{teorema de Gauss}, el cual dice:

$\int_V \nabla \cdot F = \int_{\partial V}  F$

donde $V$ es un volumen de $\Re^3$ cuya frontera $\partial V$ es una superficie suave. Ahora la podemos escribir como

$\int_V d F^2 = \int_{\partial V}  F^2$

donde el campo vectorial $F$ se modeló como una forma diferencial de grado 2, que se puede integrar en el lado derecho sobre una superficie, o bien, cuya diferencial da una forma de grado 3, en el lado izquierdo, la cual opera sobre 3-plps de un volumen para dar una integral de volumen.

\bigskip

Observemos que en la reescritura adoptada, ambos teoremas tienen exactamente la misma forma. Sobre  una variedad $M$ de dimensión $n$ y cuya frontera $\partial M$ es una variedad de dimensión $n-1$,  se tiene:

$\int_M d F^{n-1} = \int_{\partial M}  F^{n-1}$

donde integramos una forma diferencial  $F$   de grado n-1, que se puede integrar en el lado derecho sobre una variedad de dimensión n-1, o bien, cuya diferencial da una forma de grado n, en el lado izquierdo, la cual opera sobre n-plps que aproximadamente recubren  la variedad $M$.

\bigskip

Esta reescritura es un teorema de la geometría diferencial, válido en cualquier dimensión, y que se conoce como el teorema de Stokes.

\section{REESCRITURA DE LAS LEYES DE MAXWELL}

Como una aplicación de las correspondencias recién establecidas, \index{reescritura} reescribamos las leyes de Maxwell  en el lenguaje de las formas diferenciales.

Lo primero es decidir qué le corresponde tanto al campo eléctrico como al magnético. Para eso, miremos la fuerza de Lorentz, la cual define estos campos:

$$ \vec F = q(\vec E +  \vec v \times \vec B)$$

Esta fuerza representa la reacción de la naturaleza ante una carga eléctrica. Por lo tanto, tiene que ser una forma. Empecemos asignándole una 1-forma $f\sp 1$, a ver qué pasa. Eso implica que el campo eléctrico también sea una 1-forma, $\epsilon \sp 1$:

\addtocounter{ecu}{1}

$$\vec E \leftrightarrow\epsilon \sp 1 = E\sb 1 dx + E\sb 2 dy + E\sb 3 dz \eqno{(\theecu)}$$

 Por otro lado, para poder medir el campo magnético, el experimentador debe acelerar su partícula cargada y darle una velocidad $\vec v$, de donde inferimos que a la velocidad se le puede representar por un vector. Tenemos también un producto cruz que debe estar en las mismas unidades que el campo eléctrico, por lo que debe estar representado por  una 1-forma. Lo podemos lograr si hacemos un producto interno  entre el vector velocidad, un vector, y  el campo magnético, el cual debe entonces ser una 2-forma, $\beta\sp 2 $

\addtocounter{ecu}{1}

$$\vec B \leftrightarrow \beta \sp 2 = B\sb 1 dy \wedge dz +B\sb 2 dz \wedge dx + B\sb 3 dx \wedge dy \eqno{(\theecu)} $$

Esta 2-forma al operar sobre un 2-plp encontrará la primera arista siempre ocupada por $\vec v$. Por consiguiente, la ley de Lorentz se escribe, con las asociaciones naturales, como:

\addtocounter{ecu}{1}

$$f\sp 1= q(\epsilon\sp 1 + i\sb{\vec v}\beta \sp 2) \eqno{(\theecu)}$$

Ahora reescribamos las leyes de Maxwell.

La primera:

$\nabla  \times \vec E =-(1/c) \partial \vec B / \partial t $
se reescribe, poniendo la velocidad de la luz como uno, simplemente como

\addtocounter{ecu}{1}

$$d\epsilon\sp 1 = -\partial \beta \sp 2/\partial t \eqno{(\theecu)}$$

La cuarta ley $\nabla  \cdot \vec B =0 $ se reescribe como

\addtocounter{ecu}{1}

$$d\beta \sp 2=0 \eqno{(\theecu)}$$

Pero ahora tenemos un problema, la tercera ley es:
$\nabla  \times \vec B =-(1/c) \partial \vec E / \partial t $ y en esta ecuación, el operador $\nabla$ no puede substituirse por $d$ porque $d$ sobre una 2-forma mide una divergencia y nosotros lo que queremos es un rotacional. No hay más remedio que representar al campo magnético de dos maneras: como una 2-forma para la primera y cuarta ley y como una 1-forma para la segunda y tercera ley.

Sin embargo,  en el lado izquierdo de la tercera ley tendríamos la derivada de una 1-forma, la cual sería una 2-forma, mientras que en el lado derecho tendríamos la deriva parcial de una 1-forma que seguiría siendo una 1-forma. La nueva escogencia nos obliga a representar al campo eléctrico por una 2-forma, lo cual ya sabemos hacer.

Podemos hacer todo esto sin que haya lugar a confusión. Adoptamos, en primer término, las convenciones aceptadas para la primera y cuarta ley y nos las arreglamos para la segunda y la tercera. Tomamos el campo eléctrico como una 1-forma y al magnético como una 2-forma. Definimos:

\addtocounter{ecu}{1}

$$*\epsilon \sp 2 = (*\epsilon) \sp 2 = i\sb{\vec E}vol\sp 3 = E\sp 1dx\sp 2\wedge dx\sp 3 + E\sp 2dx\sp 3\wedge dx\sp 1 + E\sp 3dx\sp 1\wedge dx\sp 2 \eqno{(\theecu)}$$

\addtocounter{ecu}{1}

$$*\beta \sp 1=(*\beta) \sp 1= B\sb 1 dx\sp 1 + B\sb 2 dx\sp 2 + B\sb 3 dx\sp 3 \eqno{(\theecu)}$$

La * se escribe antes de la forma en cuestión, pues realmente es un operador que a asocia formas a formas: a una forma de grado 1, le hace corresponder otra de grado 2, y a una de grado 2, le hace corresponder otra de grado 1. Con esa adaptación, podemos reescribir la tercera ley como

\addtocounter{ecu}{1}

$$d*\beta \sp 1 = - \partial *\epsilon\sp 2/\partial t \eqno{(\theecu)}$$

y la cuarta ley

$\nabla  \cdot \vec E = 4 \pi  Q$

la podemos escribir como

\addtocounter{ecu}{1}

$$d*\epsilon \sp 2 = 4 \pi \sigma \sp 3 \eqno{(\theecu)}$$

donde  la densidad de carga, Q,  está representada por una 3-forma, $\sigma \sp 3 $, pues la densidad de carga se integra sobre un volumen 3-dimensional.

\addtocounter{ecu}{1}

Resumen (\theecu): Las leyes de Maxwell en el cálculo de formas diferenciales se escriben como

$d \epsilon\sp 1 = -\partial \beta \sp 2/\partial t$

$d  \beta \sp 2=0$

$d  *\beta \sp 1 = 4 \pi j\sp 2 + \partial *\epsilon\sp 2/\partial t $

$d  *\epsilon \sp 2 = 4 \pi \sigma \sp 3$

La simetría de estas ecuaciones en el vacío, cuando no hay cargas ni corrientes,  exige investigar si existirá  una reescritura más avanzada en la cual se borren las diferencias entre campo magnético y eléctrico y en particular si existen cargas magnéticas, los llamados monopolos. Parte de esa investigación la realizamos sobre el espacio-tiempo.

\section{LAS LEYES DE MAXWELL EN EL ESPACIO TIEMPO}

Recordamos que en relatividad, todas las observables básicas vienen empaquetadas en cuadrivectores. El conjunto de cuadrivectores se denomina espacio-tiempo.

\

El   espacio-tiempo está provisto de un \index{producto interno} \textbf{producto interno}, tal que si  los cuadrivectores son $(u\sb o, u\sb 1, u\sb 2, u\sb 3)$ y $(v\sb o, v\sb 1, v\sb 2, v\sb 3) $  entonces su producto interno es

 $<(u\sb o, u\sb 1, u\sb 2, u\sb 3),(v\sb o, v\sb 1, v\sb 2, v\sb 3) > = -c\sp 2u\sb o v\sb o +  u\sb 1 v\sb 1+ u\sb 2 v\sb 2+ u\sb 3 v\sb 3$

 A este producto interno se le llama también la \index{métrica o pseudo-métrica de Minkowski}  \textbf{métrica o pseudo-métrica de Minkowski}, pues permite generar una norma: la norma cuadrado de un cuadrivector se define como el producto del vector por el mismo: $\|P\|\sp 2 = <P,P>$.

 Para hacer énfasis en la norma, también se puede decir. espacio de Minkowski. Es algebraicamente equivalente a decir que se tiene una métrica común y corriente pero con el tiempo   imaginario, pero en realidad se prefiere decir que el tiempo es real y que la métrica tiene \index{signatura}  \textbf{signatura} $(-+++)$.

Para considerar la reescritura  de las leyes de Maxwell en el espacio-tiempo necesitamos extender el formalismo tridimensional, considerado en el capítulo anterior, a uno tetradimensional para que incluya el tiempo.  Para empezar debemos lograr identificar la fuerza de Lorentz. Definiendo fuerza como masa por aceleración estamos suponiendo que la masa es constante y perdemos la posibilidad de analizar una bomba que al desinflarse va dando tumbos de una lado para el otro. Notando que la definición anterior es equivalente a fuerza igual a masa por la derivada de la velocidad, definimos fuerza como la derivada de la  masa por la velocidad, o sea, la derivada del momento, $P$. Eso es válido tanto en mecánica clásica como en relativista. Pero en relatividad, la derivada es respecto al tiempo propio $\tau$, o sea al tiempo medido por un reloj que viaja con la partícula test. Claramente:

Al tiempo propio lo notamos $\tau$, al tiempo medido por el observador como $t$. Los dos se relacionan mediante el factor de Lorentz, $\gamma$:

$dt\sp 2 = (1- (v/c)\sp 2)^{-1}d\tau\sp 2 = \gamma \sp 2 d\tau \sp 2 $

Otra definición alterna del factor de Lorentz es $\gamma = dt/d\tau = (1-v\sp 2/c\sp 2)\sp{-1/2}$

A la posición en 3:d la notamos $\vec x$ y al cuadrivector posición, relativa al observador, $x=(t,\vec x)$.
La velocidad en 3:d la notamos $\vec v$,  la flechita siempre marca los vectores en 3:d,
y al cuadrivector velocidad como $u = dx/d\tau =(dt/d\tau, d\vec x/d\tau) = \gamma (1,\vec v)$. Indistintamente, usaremos la equivalencia entre vector y operador de derivación direccional. Por tanto, al cuadrivector velocidad también lo notaremos: $u=\gamma \partial/\partial \tau + \gamma\vec v $.

El momento es la masa en reposo por el cuadrivector velocidad:

 $P= m\sb o u = m\sb o \gamma (1,\vec v)= m(1, \vec v) = (m, \vec p)$

La fuerza es entonces: $f=dP/d\tau = d(m,\vec p)/d\tau=(dm/d\tau, d\vec p/d\tau) =(dm/d\tau,d\vec p/dt (dt/d\tau))=  (dm/d\tau, \gamma d\vec p/dt) $, o sea,

$f= (f\sp o, \gamma \vec f)$

donde $\vec f$ es la fuerza clásica en $\Re\sp 3$ y $f\sp o$ lleva la contabilidad del cambio de masa: nuestros conceptos se aplican a cohetes que cambian de masa a cada instante.

\addtocounter{ecu}{1}

\textit{Teorema (\theecu): La norma del cuadrivector momento}  $P$ es $-m\sb o \sp 2 c \sp 2$.

\

Veamos:

El producto interno  entre dos cuadrivectores cualesquiera en el espacio -tiempo es $<(u\sb o, u\sb 1, u\sb 2, u\sb 3),(v\sb o, v\sb 1, v\sb 2, v\sb 3) > = -c\sp 2u\sb o v\sb o +  u\sb 1 v\sb 1+ u\sb 2 v\sb 2+ u\sb 3 v\sb 3$ y la norma cuadrado es  el producto del vector por el mismo:

 $\|P\|\sp 2 = <P,P>$.

De tal manera que la norma cuadrado de un cuadrivector $(t,x,y,z)$ es $\|(t,x,y,z)\| = -c\sp 2 t\sp 2 + x\sp 2 +y \sp 2 + z \sp 2$ , donde aparece la norma cuadrado en 3:d. Por lo tanto, si el momento $P= (m\sb o \gamma , m\vec v) =(m\sb o \gamma , m\sb o \gamma \vec v) $, su  norma cuadrado es:

$\|P\|\sp 2 = <P,P> =  -m\sb o \sp 2\gamma \sp 2 c \sp 2  + m\sb o \sp 2 \gamma \sp 2  v\sp 2 = -m\sb o \sp 2\gamma \sp 2 c \sp 2  (1- v^2/c\sp 2) = -m\sb o \sp 2\gamma \sp 2 c \sp 2 \gamma \sp{-2} = -m\sb o \sp 2 c \sp 2 = $ constante.

\addtocounter{ecu}{1}

\textit{Teorema (\theecu): El cuadrivector velocidad tiene las componentes}

 $$f=\gamma(c\sp{-2}\vec f\cdot \vec v, \vec f) \eqno{(\theecu)}$$

Demostración: Como la norma cuadrado del cuadrivector momento se conserva, al derivar con respecto al tiempo propio debe dar  cero:

$d/d\tau(<P,P>)=0 = <dP/d\tau, P> + <P,dP/d\tau>$

$ = 2<dP/d\tau, P> =2<f, P> =2<f,m\sb o u>=0$

lo cual dice que la fuerza es ortogonal a la cuadrivelocidad en el espacio de Minkowski:

$0=<f,u>$, pero como $f= (f\sp o, \gamma \vec f)$ y $u= (\gamma, \gamma \vec v)$ entonces

$0 = <(f\sp o, \gamma \vec f),(\gamma, \gamma \vec v)> = -c\sp 2 f\sp o \gamma + \gamma \vec f \cdot \gamma \vec v$

Despejando tenemos: $f\sp o = (\gamma/c\sp 2) \vec f\cdot\vec v$ lo cual dice: la componente temporal en el espacio -tiempo es proporcional a la potencia en 3:d. Al substituir nos da que el cuadrivector fuerza tiene los siguientes componentes:

$$f=\gamma(c\sp{-2}\vec f\cdot \vec v, \vec f) \eqno{(\theecu)}$$

donde podemos ver la traza de la teoría clásica.

\

Ahora podemos  traducir la ley de Lorentz en 3:d, $\vec f= q[\vec E + (\vec v/c) \times \vec B]$ al espacio de Minkowski.

\

Para eso, primero debemos notar que $\vec v$ es ortogonal en 3:d a $(\vec v/c)\times \vec B$, por lo tanto, el cuadrivector fuerza de Lorentz queda:

$f= \gamma(c\sp{-2}\vec f\cdot \vec v, \vec f)= \gamma(c\sp{-2}(q[\vec E + (\vec v/c) \times \vec B])\cdot \vec v, (q[\vec E + (\vec v/c) \times \vec B])) $

y simplificando

\addtocounter{ecu}{1}

$$f= \gamma q(c\sp{-2}\vec E\cdot \vec v, \vec E + (\vec v/c) \times \vec B)\eqno{(\theecu)} $$

la cual podemos reescribir en el lenguaje de las formas usando nuestra experiencia en 3:d. Al cuadrivector fuerza le corresponde una 1-forma, pues la fuerza representa la reacción de la naturaleza. La traducción requiere dos detalles: primero, cada término debe ser una 1-forma. Por eso, en la parte temporal hay que agregar $dt$, pues de lo contrario el primer término sería un escalar. En segundo lugar, dado que en el producto interior a la primera coordenada se le multiplica por $-c\sp 2$, al pasar de vectores a formas hay que dividir por esa misma cantidad. Eso se debe a la dualidad existente entre formas y vectores. Nos queda:

\addtocounter{ecu}{1}

$$f\sp 1 = -(\gamma  q/c\sp 2)[i\sb{\vec v}\epsilon \sp 1]dt + \gamma q[\epsilon \sp 1 - i\sb{\vec v/c}\beta \sp 2]\eqno{(\theecu)}$$

donde la velocidad de la luz $c$ puede medirse en cualquier tipo de unidades. Si se usa la convención $c=1$,  la fuerza queda:

$f\sp 1 = -(\gamma  q)[i\sb{\vec v}\epsilon \sp 1]dt + \gamma q[\epsilon \sp 1 - i\sb{\vec v}\beta \sp 2]$

\bigskip

 Tenemos ahora un teorema impresionante:

\addtocounter{ecu}{1}

\textit{Teorema (\theecu): Si definimos el tensor de campo que es la 2-forma fuerza del campo electromagnético como}

$F\sp 2 =  \epsilon \sp 1 \wedge dt + c\sp{-1} \beta\sp 2$

$=E\sb 1 dx \wedge dt + E\sb 2 dy \wedge dt + E\sb 3 dz\wedge dt
 + B\sb 1 dy \wedge dz +B\sb 2 dz \wedge dx + B\sb 3 dx \wedge dy $

\textit{ entonces la fuerza de Lorentz se escribe como la 1-forma en el espacio-tiempo}

 $f\sp 1 = -qi\sb{u} F\sp 2$

\bigskip

Lo que  parece impresionante de este teorema es que  de en medio de toda la sofisticación física y matemática que hemos visto, la cual fue construida para ser generalizada automáticamente  a variedades, hace que emerja una formulación intuitivamente simple y directa: visto infinitesimalmente, el producto interior de un vector con una 2-forma reemplaza lo que en álgebra lineal corresponde a la multiplicación de una matriz por un vector. Por lo tanto, lo que este teorema dice es: el experimentador estimula la naturaleza con un cuadrivector velocidad, $u$,  la naturaleza lo detecta por medio del campo $F\sp 2$, y su reacción es producir una resistencia, una fuerza, $f\sp 1$ que es simplemente una transformación lineal del estímulo por medio del campo. Esto es lo que siempre hemos buscado. Y lo estamos encontrando.

Demostración: $-qi\sb{u} F\sp 2= -q i\sb u ( \epsilon \sp 1 \wedge dt + c\sp{-1} \beta\sp 2)$

$=-q i\sb u ( \epsilon \sp 1 \wedge dt)  -(q/c)  i\sb u \beta\sp 2$

$=-q ( i\sb u \epsilon \sp 1 \wedge dt-\epsilon \sp 1 \wedge i\sb u (dt))  -(q/c)  i\sb u \beta\sp 2$

$=-q ( \epsilon \sp 1 (u) \wedge dt-\epsilon \sp 1 dt(u))  -(q/c)  i\sb u \beta\sp 2$

tengamos en cuenta que

$u=\gamma \partial/\partial t + \gamma\vec v = \gamma \partial/\partial t + \gamma[v\sb 1\partial/\partial x +v\sb 2\partial/\partial y +v\sb 3\partial/\partial z]$

por consiguiente $dt(u) = \gamma$, pues $dt(\partial \partial t)=1$ y las demás proyecciones dan cero. Por tanto:

$-qi\sb{u} F\sp 2=-q ( \epsilon \sp 1 (u) )\wedge dt+\gamma q \epsilon \sp 1   -(q/c)  i\sb u \beta\sp 2$

Averigüemos qué es $\epsilon \sp 1 (u)$ donde $\epsilon \sp 1 = E\sb 1 dx + E\sb 2 dy + E\sb 3 dz$,

$\epsilon \sp 1 (u)= E\sb 1 dx(u) + E\sb 2 dy (u) + E\sb 3 dz(u) = E\sb 1 \gamma v\sb 1 + E\sb 2 \gamma v\sb 2 + E\sb 3  \gamma v\sb 3 = \gamma \vec E \cdot \vec v $

$= \gamma i\sb{\vec v}\epsilon \sp 1 $ el cual es un escalar. Reemplazando:

$-qi\sb{u} F\sp 2=-\gamma q (i\sb{\vec v}\epsilon \sp 1 ) dt+\gamma q \epsilon \sp 1   -(q/c)  i\sb u \beta\sp 2$

Especifiquemos qué es $i\sb u \beta\sp 2$ donde $\beta \sp 2 = B\sb 1 dy \wedge dz +B\sb 2 dz \wedge dx + B\sb 3 dx \wedge dy $. Tenemos:

$i\sb u \beta\sp 2= i\sb{\gamma \partial /\partial t}  \beta\sp 2 + i\sb{\gamma \vec v} \beta\sp 2 =i\sb{\gamma \vec v} \beta\sp 2 =\gamma i\sb{ \vec v} \beta\sp 2  $, puesto que el campo magnético no tiene componente temporal y además el producto interno es lineal en todas sus entradas. Substituyendo:

$-qi\sb{u} F\sp 2=-\gamma q (i\sb{\vec v}\epsilon \sp 1 ) dt+\gamma q \epsilon \sp 1   -(q/c)  \gamma i\sb{ \vec v} \beta\sp 2 =-\gamma q (i\sb{\vec v}\epsilon \sp 1 ) dt+\gamma q \epsilon \sp 1   -q  \gamma i\sb{ \vec v/c} \beta\sp 2 $

En conclusión, hemos probado la igualdad requerida:

$-qi\sb{u} F\sp 2=-\gamma q (i\sb{\vec v}\epsilon \sp 1 ) dt+\gamma q [\epsilon \sp 1
                  - i\sb{ \vec v/c} \beta\sp 2] = f\sp 1$

\bigskip

De ahora en adelante, $c$ se será la unidad de velocidad : $c=1$.

Observemos que el tensor de campo

 $F\sp 2=E\sb 1 dx \wedge dt + E\sb 2 dy \wedge dt + E\sb 3 dz\wedge dt
 + B\sb 1 dy \wedge dz +B\sb 2 dz \wedge dx + B\sb 3 dx \wedge dy $

$=-E\sb 1 dt \wedge dx - E\sb 2 dt \wedge dy - E\sb 3 dt\wedge dz
 + B\sb 1 dy \wedge dz -B\sb 2 dx \wedge dz + B\sb 3 dx \wedge dy$

puede escribirse en forma de arreglo pseudo-matricial como sigue:

\addtocounter{ecu}{1}

$$F\sp 2= \bordermatrix{&dt &dx&dy&dz\cr
dt& 0&-E\sb 1&  -E\sb 2&  -E\sb 3\cr
dx& &0& B\sb 3& -B\sb 2\cr
dy &&&0&  B\sb 1 \cr
dz&&&&0\cr} \eqno{(\theecu)}
$$

\bigskip

Ahora bien, esa no es una matriz. Para asociarle una matriz hay que primero fijar una base.   Luego se  evalúa $F^2$ sobre los 2-plps generados por los vectores de la base. Nosotros siempre hemos trabajado en la base canónica. En nuestro caso tenemos, por ejemplo, que el término $B\sb 3 dx \wedge dy$ evaluado sobre el 2-plp $(\partial /\partial x, \partial /\partial y)$ da $B\sb 3$ y sobre todos los demás 2-plps da cero. Eso significa que a la componente $B\sb 3 dx \wedge dy$ de $F^2$ le corresponde la entrada $B\sb 3$ en la matriz dada por la base canónica. Por lo tanto, la matriz asociada , la de evaluaciones es:

\addtocounter{ecu}{1}

$$[F]=\bordermatrix{&\partial /\partial t &\partial /\partial x&\partial /\partial y&\partial /\partial z\cr
\partial /\partial t& 0&-E\sb 1&  -E\sb 2&  -E\sb 3\cr
\partial /\partial x&E\sb 1 & 0& B\sb 3& -B\sb 2\cr
\partial /\partial y  &E\sb 2 &-B\sb 3&0&  B\sb 1 \cr
\partial /\partial z&E\sb 3 &B\sb 2&-B\sb 1&0\cr} \eqno(\theecu)
$$

A este arreglo se le denomina \index{tensor electromagnético}  \textbf{'tensor electromagnético'} o bien 'tensor de fuerza del campo electromagnético'. Se dice que es una 'fuerza' puesto que puede cambiar la velocidad de una partícula cargada, sea en su módulo o sea en su dirección.

Obsérvese la distribución de trabajo adjudicada por las dos-formas al campo electromagnético:  la parte que tiene coordenada temporal, el campo eléctrico, es una fuerza que puede cambiar el módulo de la velocidad. En tanto que las componentes espaciales forman un campo, el campo magnético, el cual no puede  cambiar el módulo de la velocidad sino solamente su dirección.

Las coordenadas de la fuerza de Lorentz pueden entonces formularse, usando esta matriz, como:

$f\sb i = q\sum \sb j F\sb{ij} u\sb j$

que es una expresión muy directa. Lo es a tal punto que uno debería preguntarse cuál es el sentido de todo el tiempo que hemos invertido en una elaboración elemental de las formas diferenciales.  La respuesta está en la anunciada reescritura de las leyes de Maxwell, para lo cual necesitamos extender el operador derivada (exterior) de 3:d al espacio -tiempo o 3+1:d.

\bigskip

\addtocounter{ecu}{1}

\textit{Definición (\theecu): La derivada en el espacio-tiempo de una forma diferencial se define como}

$d(\cdot)= d\sb 3(\cdot) + \partial (\cdot)/\partial t \wedge dt$.

\textit{donde $d_3$ es la derivada en el espacio ordinario 3:d. La dimensión temporal se toma como parámetro.}

\bigskip

Nuestra definición se motiva diciendo que tiene que extender una derivación normal, por eso aparece la derivada partial en t, y además tiene que dar una 1-forma sobre funciones, por lo cual aparece la $dt$. Por supuesto, nuestra definición cumple todas las propiedades de la derivada.

Operando sobre el tensor de campo, $F\sp 2 = \epsilon \sp 1 \wedge dt + \beta \sp 2$, tenemos:

$dF\sp 2 = d\sb 3 \epsilon \sp 1 \wedge dt + d\sb 3 \beta \sp 2 +  \partial \epsilon \sp 1/\partial t\wedge dt \wedge dt + \partial \beta \sp 2/\partial t\wedge dt$

$= (d\sb 3 \epsilon \sp 1 + \partial \beta \sp 2/\partial t) \wedge dt + d\sb 3 \beta\sp 2$

puesto que el tercer término desaparece por repetición de multiplicandos.

Ahora, consideremos la ecuación $dF\sp 2 = 0$ . Como la expresión calculada para $dF\sp 2$ está separada en sus partes independientes de tiempo y espacio, tenemos por  la parte temporal: $d\sb 3 \epsilon \sp 1 + \partial \beta \sp 2/\partial t$ o sea $d\sb 3 \epsilon \sp 1 =- \partial \beta \sp 2/\partial t$. Mientras que en su parte espacial queda: $d\sb 3 \beta\sp 2=0$, y obtuvimos dos leyes de Maxwell. Las otras dos las obtendremos usando las representaciones complementarias de los campos dados.

\addtocounter{ecu}{1}

  \textit{Teorema} (\theecu):

$d\sb 3 *\epsilon \sp 2 -(d\sb 3 *\beta \sp 1 - \partial *\epsilon\sp 2/\partial t)\wedge dt =  4\pi (\sigma \sp 3 - j\sp 2 \wedge dt)$
$= 4 \pi \jmath\sp 3 $ 

\textit{donde} $\jmath\sp 3 =\sigma \sp 3 - j\sp 2 \wedge dt$

\bigskip

Demostración:

La versión 3:d de las otras dos leyes de Maxwell, las cuales describen el efecto de las corrientes y cargas, son:

$d\sb 3 *\beta \sp 1 = 4 \pi j\sp 2 + \partial *\epsilon\sp 2/\partial t $

$d\sb 3 *\epsilon \sp 2 = 4 \pi \sigma \sp 3$

Estas ecuaciones se pueden reescribir como:

$d\sb 3 *\epsilon \sp 2 = 4 \pi \sigma \sp 3$

$d\sb 3 *\beta \sp 1 - \partial *\epsilon\sp 2/\partial t = 4 \pi j\sp 2 $

Deseamos restar estas dos ecuaciones, pero hay un problema: la primera ecuación es entre 3-formas, mientras que la segunda es entre 2-formas.  Para balancearlas, podemos multiplicar la segunda ecuación por una 1-forma, pero como dicha forma no puede tener ninguna coordenada espacial, las cuales ya están en esa ecuación, entonces nos queda la opción de insertarlas en 3+1:d y multiplicar la segunda ecuación por una 1-forma temporal. Las dos ecuaciones quedan:

$d\sb 3 *\epsilon \sp 2 = 4 \pi \sigma \sp 3$

$(d\sb 3 *\beta \sp 1 - \partial *\epsilon\sp 2/\partial t)\wedge dt = 4 \pi j\sp 2\wedge dt $

y al restarlas, nos da:

$d\sb 3 *\epsilon \sp 2 -(d\sb 3 *\beta \sp 1 - \partial *\epsilon\sp 2/\partial t)\wedge dt = 4 \pi \sigma \sp 3 - 4 \pi j\sp 2 \wedge dt = 4\pi (\sigma \sp 3 - j\sp 2 \wedge dt)$

$= 4 \pi \jmath\sp 3 $

donde hemos definido la  3-forma $\jmath\sp 3 =\sigma \sp 3 - j\sp 2 \wedge dt$ se le denomina la 3-forma corriente.

\bigskip

\addtocounter{ecu}{1} 

\textit{Definición (\theecu).  Definimos ahora una representación equivalente del tensor de campo en 3+1:d con ayuda de las representaciones equivalentes de los campos eléctrico y magnético en 3:d:}

$*F\sp 2 = -*\beta \sp 1 \wedge dt + *\epsilon \sp 2$

\textit{tengamos en cuenta que en esta definición la * en el lado izquierdo es un operador en 3+1:d, mientras que en el lado derecho es en 3:d.}

\bigskip

\addtocounter{ecu}{1}

\textit{Teorema} (\theecu): $d*F\sp 2 = 4 \pi \sigma \sp 3 -4 \pi j\sp 2\wedge dt = 4 \pi \jmath \sp 3$

 Puesto que

 $*F\sp 2 = -*\beta \sp 1 \wedge dt + *\epsilon \sp 2$

 derivando obtenemos una 3-forma en 3+1:d:

$d*F\sp 2 = d( -*\beta \sp 1 \wedge dt + *\epsilon \sp 2) =[ d\sb 3(\cdot) + \partial (\cdot)/\partial t \wedge dt](-*\beta \sp 1 \wedge dt + *\epsilon \sp 2)$.

$=d\sb 3(-*\beta \sp 1 \wedge dt) + d\sb 3( *\epsilon \sp 2) + \partial (-*\beta \sp 1 \wedge dt)/\partial t \wedge dt + \partial (*\epsilon \sp 2)/\partial t \wedge dt$

$=-d\sb 3(*\beta \sp 1 \wedge dt) + d\sb 3( *\epsilon \sp 2) +  \partial (*\epsilon \sp 2)/\partial t \wedge dt$

$= d\sb 3( *\epsilon \sp 2)-d\sb 3(*\beta \sp 1 \wedge dt) +  \partial (*\epsilon \sp 2)/\partial t \wedge dt$

$= d\sb 3( *\epsilon \sp 2)-[d\sb 3*\beta \sp 1 -  \partial (*\epsilon \sp 2)/\partial t] \wedge dt$

Usando el teorema anterior, podemos concluir que:

$d*F\sp 2 = 4 \pi \sigma \sp 3 -4 \pi j\sp 2\wedge dt = 4 \pi \jmath \sp 3$

\

Tomando otra derivada más:

$d\sp 2*F\sp 2 = 0 = 4 \pi d \jmath \sp 3$

que dice que la carga ni se crea ni se destruye, por lo tanto, lo que escapa de un recinto cerrado debe registrarse como un flujo hacia el exterior.

En resumen, las 4 leyes de Maxwell pueden reescribirse en 3+1:d como

\addtocounter{ecu}{1}

$$dF\sp 2 = 0 \eqno{(\theecu)}$$

\addtocounter{ecu}{1}

$$d*F\sp 2 = 4 \pi \jmath \sp 3 \eqno{(\theecu)}$$

Las 4 ecuaciones se han convertido en 2. Y por qué no en 1? Debería poderse pues después de todo, el campo electromagnético es una entidad que funciona como un todo.

\bigskip

\section{EL VECTOR POTENCIAL}

Todo lo anterior puede reinventarse, por poco, a partir de la nada, al menos localmente, si tan sólo uno se permite explorar las consecuencias de admitir que la 2-forma tensor de campo sea la derivada de una 1-forma, a la que notaremos $A\sp 1$ y a la que podemos llamar la 1-forma correspondiente al \index{vector potencial} \textbf{vector potencial}. (Esto está garantizado siempre que no haya problemas de singularidades).

Una 1-forma $A$ en el espacio-tiempo tiene 4 coordenadas, $A= A\sb o dt + A\sb 1 dx + A\sb 2 dy + A\sb 3 dz$ las cuales podemos separar en su parte temporal $A\sb o dt $ y espacial, $\Lambda \sp 1= A\sb 1 dx + A\sb 2 dy + A\sb 3 dz$ que sigue siendo una 1-forma:

$ A= A\sb o dt+ \Lambda \sp 1$

Este potencial 1-forma  permite la libertad gauge y por lo tanto no es observable. Lo que se observa es el efecto del campo electromagnético, el cual obtendremos de $A$ por la única operación que conocemos, la diferenciación: $F=dA$. Miremos la forma general de $dA$:

$dA=F=() dt \wedge dx+ () dt \wedge dy +() dt \wedge dz +() dx \wedge dy+() dx \wedge dz+() dy \wedge dz  $

esta forma escrita de lleno, y cambiando ligeramente la notación,  es:

$$\bordermatrix{&dx\sb o &dx\sb 1&dx\sb 2&dx\sb 3\cr
dx\sb o& 0&\partial A\sb 1 /\partial x\sb o- \partial A\sb o /\partial x\sb 1
&  \partial A\sb 2 /\partial x\sb o- \partial A\sb o /\partial x\sb 2 &
  \partial A\sb 3 /\partial x\sb o- \partial A\sb o /\partial x\sb 3\cr
dx\sb 1& &0& \partial A\sb 2 /\partial x\sb 1- \partial A\sb 1 /\partial x\sb 2& \partial A\sb 3 /\partial x\sb 1- \partial A\sb 1 /\partial x\sb 3\cr
dx\sb 2 &&&0&   \partial A\sb 3 /\partial x\sb 2- \partial A\sb 2 /\partial x\sb 3\cr
    dx\sb 3 & &&&0\cr}
$$

comparemos este arreglo con aquel obtenido anteriormente para $F\sp 2$:

$$\bordermatrix{&dt &dx&dy&dz\cr
dt& 0&-E\sb 1&  -E\sb 2&  -E\sb 3\cr
dx& &0& B\sb 3& -B\sb 2\cr
dy &&&0&  B\sb 1 \cr
dz&&&&0\cr}
$$

Comparando el primer renglón vemos que:

$-E\sb 1 = \partial A\sb 1 /\partial x\sb o- \partial A\sb o /\partial x\sb 1 $

$  -E\sb 2=  \partial A\sb 2 /\partial x\sb o- \partial A\sb o /\partial x\sb 2  $

$ -E\sb 3 =  \partial A\sb 3 /\partial x\sb o- \partial A\sb o /\partial x\sb 3$

o bien que

$E\sb 1 =  \partial A\sb o /\partial x\sb 1 -\partial A\sb 1 /\partial x\sb o$

$  E\sb 2=   \partial A\sb o /\partial x\sb 2-\partial A\sb 2 /\partial x\sb o  $

$ E\sb 3 =   \partial A\sb o /\partial x\sb 3-\partial A\sb 3 /\partial x\sb o$

Tengamos en cuenta que

$d\sb 3 A\sb o dt = (\partial A\sb o /\partial x\sb 1) dx_1 \wedge dt +
(\partial A\sb o /\partial x\sb 2) dx_2 \wedge dt +
(\partial A\sb o /\partial x\sb 3) dx_3 \wedge dt$

el cual tiene la forma de un gradiente multiplicado por $dt$, de lo cual se concluye que $A\sb o$ debe ser el potencial eléctrico, escalar, $\phi$.

\

Por otro lado,

$\partial \Lambda\sp 1/\partial t = (\partial/\partial t) ( A\sb 1 dx + A\sb 2 dy + A\sb 3 dz)$

cuyos coeficiente  pueden escribirse como

$\partial A\sb 1 /\partial x\sb o$

$\partial A\sb 2 /\partial x\sb o$

$\partial A\sb 3 /\partial x\sb o$

Podemos concluir que el campo eléctrico, una 1-forma en 3:d, se relaciona con $A$ como sigue:

$\epsilon \sp 1 = d\sb 3 A\sb o dt- \partial \Lambda\sp 1/\partial t $

 La anterior ecuación diferencial es entonces equivalente a la ecuación vectorial

$\vec E =\nabla \phi - \partial \Lambda^1/\partial t$

Por otro lado, en las restantes componente de $dA\sp 1$ uno encuentra el rotacional del vector correspondiente a $\Lambda \sp 1$. Por tanto, esa 1-forma es la representación diferencial del potencial vector, $\vec A$ y, recordando la equivalencia entre rotacional de un vector y la derivada de su correspondiente 1-forma, tenemos la relación entre campo magnético y vector potencial:

$\beta \sp 2 = d\sb 3 \Lambda \sp 1$

lo cual es equivalente a : $\vec B = \nabla \times \vec A$. Con todas esas equivalencias, $dA $ puede escribirse como el tensor de campo $F$, que ya sabemos que es una 2-forma:

$A^1 = A_o dt + A_1dx + A_2 dy + A_3 dz$

$dA^1 = F^2$ donde

$F^2=E\sb x dx \wedge dt+E\sb y dy \wedge dt + E\sb z dz \wedge dt +B\sb x dy \wedge dz +B\sb y dz \wedge dx +B\sb z dx \wedge dy  $

Esto implica que, necesariamente, $dF=d\sp 2 A = 0$, lo cual tiene la forma de una ley de conservación: las variaciones permitidas son aquellas que, en cierto sentido, no cambian a $F$. De esto se obtienen las dos primeras leyes de Maxwell. Aunque eso ya lo sabíamos, veámoslo desde un punto de vista muy expandido. El tensor de campo:

$F^2=E\sb x dx \wedge dt+E\sb y dy \wedge dt + E\sb z dz \wedge dt +B\sb x dy \wedge dz +B\sb y dz \wedge dx +B\sb z dx \wedge dy  $

al derivarlo produce:

$dF^2= -[(\partial E\sb x/\partial t)dt +(\partial E\sb x/\partial x)dx+(\partial E\sb x/\partial y)dy+(\partial E\sb x/\partial z)dz ]\wedge dt \wedge dx$

$-[(\partial E\sb y/\partial t)dt +(\partial E\sb y/\partial x)dx+(\partial E\sb y/\partial y)dy+(\partial E\sb y/\partial z)dz ] \wedge dt \wedge dy $

$-[(\partial E\sb z/\partial t)dt +(\partial E\sb z/\partial x)dx+(\partial E\sb z/\partial y)dy+(\partial E\sb z/\partial z)dz ]\wedge dt \wedge dz $

$+[(\partial B\sb x/\partial t)dt +(\partial B\sb x/\partial x)dx+(\partial B\sb x/\partial y)dy+(\partial B\sb x/\partial z)dz ] \wedge dy \wedge dz $

$+[(\partial B\sb y/\partial t)dt +(\partial B\sb y/\partial x)dx+(\partial B\sb y/\partial y)dy+(\partial B\sb y/\partial z)dz ]\wedge dz \wedge dx $

$+[(\partial B\sb z/\partial t)dt +(\partial B\sb z/\partial x)dx+(\partial B\sb z/\partial y)dy+(\partial B\sb z/\partial z)dz ]\wedge dx \wedge dy  $

Teniendo en cuenta que las repeticiones en un producto wedge se aniquilan, que al cambiar de orden una vez, se cambia de signo, y que al hacer dos permutaciones no pasa nada,  podemos simplificar:

\

$dF^2  = - (\partial E\sb x/\partial y)dy \wedge dt \wedge dx$
$  - (\partial E\sb x/\partial z)dz \wedge dt \wedge dx$

$  - (\partial E\sb y/\partial x)dx \wedge dt \wedge dy$
$  - (\partial E\sb y/\partial z)dz \wedge dt \wedge dy$

$  - (\partial E\sb z/\partial x)dx \wedge dt \wedge dz$
$  - (\partial E\sb z/\partial y)dy \wedge dt \wedge dz$

$  + (\partial B\sb x/\partial t)dt \wedge dy \wedge dz$
$ + (\partial B\sb x/\partial x)dx \wedge dy \wedge dz$

$  + (\partial B\sb y/\partial t)dt \wedge dz \wedge dx$
$  + (\partial B\sb y/\partial y)dy \wedge dz \wedge dx$

$  + (\partial B\sb z/\partial t)dt \wedge dx \wedge dy$
$ +  (\partial B\sb z/\partial z)dz \wedge dx \wedge dy$

Numerando todos esos términos de la forma natural, observamos que los términos 8, 10 y 12  tienen la misma coordenada diferencial, por lo tanto se pueden agrupar. De igual forma, se  puede asociar los términos 4, 6 y 7, al igual que 2, 5  y 9, lo mismo que los 1, 3, 11  :

$dF= (\partial B\sb x/\partial x + \partial B\sb y/\partial y +\partial B\sb z/\partial z )dx \wedge dy \wedge dz$

$+ (\partial B\sb x/\partial t + \partial E\sb z/\partial y -\partial E\sb y/\partial z )dt \wedge dy \wedge dz$

$+(\partial B\sb y/\partial t + \partial E\sb x/\partial z -\partial E\sb z/\partial x )dt \wedge dz \wedge dx$

$+(\partial B\sb z/\partial t + \partial E\sb y/\partial x -\partial E\sb x/\partial y )dt \wedge dx \wedge dy$

Si $dF=0$, de la primera coordenada $dx\wedge dy \wedge dz $ sacamos que $\nabla \cdot B = 0$.

En las coordenadas temporales obtenemos 3 ecuaciones escalares que se pueden resumir en una ecuación vectorial si hacemos una correspondencia: a la 3-forma $dt\wedge dy \wedge dz $ le corresponde un vector en la dirección $i$, a la 3-forma $dt\wedge dz \wedge dx $ le corresponde un vector en la dirección $j$, a la 3-forma

$dt\wedge dx \wedge dy $ le corresponde un vector en la dirección $k$. Todo se resume entonces en:

$$-\partial B/\partial t =det \pmatrix{i&j&k\cr
         \partial/\partial x&\partial/\partial y&\partial/\partial z\cr
         E\sb x&E\sb y&E\sb z\cr}=\nabla \times E
$$

Tenemos pues dos leyes. Nos falta empaquetar en el nuevo lenguaje las otras dos leyes, las cuales  predicen el efecto del campo electromagnético sobre cargas eléctricas en movimiento. Esas ecuaciones salen de una construcción que a cada  tensor $F$ le asocia su dual $\sp * F$, de tal manera que para el caso del tensor de campo electromagnético tenemos:

$F=-E\sb x dt \wedge dx-E\sb y dt \wedge dy -E\sb z dt \wedge dz +B\sb x dy \wedge dz +B\sb y dz \wedge dx +B\sb z dx \wedge dy  $

$\sp * F=E\sb x dy \wedge dz+E\sb y dz \wedge dx +E\sb z dx \wedge dy +B\sb x dt \wedge dx +B\sb y dt \wedge dy +B\sb z dt \wedge dz  $

Nosotros hemos venido definiendo  el operador $*$ de una manera ad hoc, artificial, pero hay una definición general, basada en una métrica. Podemos ver, sin embargo, que este operador calcula la forma de completar las diferenciales hasta completar un múltiplo del elemento de volumen, y por eso a una k-forma le hace corresponder una (n-k)-forma.

\bigskip

\addtocounter{ecu}{1}

\textit{Teorema (\theecu) La ecuación tensorial}

$d\sp * F = 4\pi \sp * J$

\textit{que también se escribe como el sistema de ecuaciones}

$\partial_\nu F^{\mu \nu} = 4 \pi J^\mu$

\textit{es equivalente a  las leyes de Maxwell dadas por}

$\nabla \cdot E = 4 \pi \rho $

\textit{y}

$\partial E / \partial t = \nabla \times B - 4 \pi J$

\bigskip

Demostración:

Recordemos que la \index{4-corriente} \textbf{4-corriente} es $J=[ \rho, J\sb x, J\sb y, J\sb z]$, en donde $\rho$ es la densidad de carga y las otras componentes denotan corrientes. A $J$ lo modelamos como una 1-forma:

$J= \rho dt + J\sb x dx+ J\sb y dy+ J\sb z dz$

Su forma diferencial dual es:

$\sp * J= \rho dx\wedge dy \wedge  dz - J\sb x dt\wedge dy \wedge  dz - J\sb y dt\wedge dz \wedge  dx - J\sb z dt\wedge dx \wedge  dy $

Ahora tenemos:

$d\sp * F = d (E\sb x dy \wedge dz+E\sb y dz \wedge dx +E\sb z dx \wedge dy +B\sb x dt \wedge dx +B\sb y dt \wedge dy +B\sb z dt \wedge dz)$

$=[(\partial E \sb x/\partial t)dt +(\partial E \sb x/\partial x)dx +(\partial E \sb x/\partial y)dy +(\partial E \sb x/\partial z  )dz]\wedge dy \wedge dz $
$+[(\partial E \sb y/\partial t)dt +(\partial E \sb y/\partial x)dx +(\partial E \sb y/\partial y)dy +(\partial E \sb y/\partial z  )dz]\wedge dz \wedge dx + .. .. $

$=(\partial E \sb x/\partial x +\partial E \sb y/\partial y +\partial E \sb z/\partial z  )dx \wedge dy \wedge dz$

$+(\partial E \sb x/\partial t -\partial B \sb z/\partial y +\partial E \sb y/\partial z )dt \wedge dy \wedge dz $

$+.. ..$

$=4\pi (\rho dx \wedge dy \wedge dz)$

$-4 \pi J\sb x dt \wedge dy \wedge dz$

$-4 \pi J\sb y dt \wedge dz \wedge dx$

$-4 \pi J\sb z dt \wedge dx \wedge dy$

Estas ecuaciones también se empaquetan en la forma

$\partial_\nu F^{\mu \nu} = 4 \pi J^\mu$

Que más oficialmente se escribe como

$ d\sp * F= 4 \pi \sp * J$

En la primera coordenada de esta ecuación tenemos la tercera ley de Maxwell:

$\nabla \cdot E = 4 \pi \rho $

mientras que en las otras coordenadas tenemos:

$\partial E / \partial t = \nabla \times B - 4 \pi J$

\

Para reinventar la ley de Lorentz, tenemos simplemente que recordar nuestro esquema de interpretación de toda la maquinaria de las formas diferenciales:

$$\pmatrix{Resultado:\cr
    cambio \: en \: la \:energ\acute {\i}a\cr
    y \: momento\cr}=
\pmatrix{Reacci\acute on:\cr
    fuerza  \cr
    en \: ejercicio\cr}
\pmatrix{Est\acute {\i} mulo:\cr
    desplazamiento \: de \cr
    una \:carga\cr}
$$

lo cual se traduce en:

$$\pmatrix{Resultado:\cr
    cambio \: en \: la \:energ\acute {\i}a\cr
    y \: momento\cr}=
\pmatrix{Reacci\acute on:\cr
    tensor \: de \:campo \cr
    2-forma \cr}
\pmatrix{Est\acute {\i} mulo:\cr
    4-velocidad \: de \cr
    una \:carga\cr}
$$

lo cual, pasando a coordenadas aplicando el tensor sobre 2-plps formados por elementos de la base natural, puede notarse como

$d\vec p /d\tau = e F \vec u $, donde $e$ es la carga, que es la constante de acople entre el campo y la masa, $\vec u$ es el cuadrivector velocidad, $\vec p$ es el cuadrivector energía-momento, y $\tau$ es el tiempo propio, o sea el tiempo marcado por un reloj que va pegado a la partícula. Esa misma ecuación también se escribe como:

$dp\sp \mu/d\tau = (e/mc)\sum \sb \nu F\sp{\mu \nu}p\sb \nu $

Nuestra relectura de la ley de Lorentz presenta una  inconsistencia: $F$ es una 2-forma, es decir es un operador que a un 2-plp infinitesimal le asocia un número real. Si dicha 2-forma opera sobre un vector, aún sobrará espacio para otro vector. Por lo tanto, una 2-forma sobre vector produce una 1-forma. Por consiguiente, en el lado derecho de nuestra ecuación hay una 1-forma. Pero, por otra parte, en el lado izquierdo hay un vector.

Para alegría nuestra, esa aparente inconsistencia se resuelve diciendo que el vector ha sido codificado como una 1-forma, lo cual es una opción que siempre tenemos.

\section{EL PRINCIPIO VARIACIONAL}

A partir del potencial vector, hemos tomado la diferencial para obtener el campo, $F$, y la segunda diferencial para obtener $d^2F=0$. Esta ecuación es equivalente a dos de las leyes de Maxwell. Decimos que dichas leyes son de naturaleza geométrica. Vimos también que las otras dos leyes pueden considerarse como necesarias dada la ambivalencia en la representación de un vector en el espacio-tiempo, sea como una 1-forma, sea como una 2-forma.  Ahora vamos a probar que estas dos   leyes también son equivalentes a un principio variacional. Con esto ilustramos la forma como se \index{inventar principio variacional} \textbf{inventa un  principio variacional} para deducir unas ecuaciones dadas.

Para poder proseguir necesitamos hacer una definición. Como ya sabemos, el tensor electromagnético, el cual es una dos forma, también genera  la siguiente matriz que lo representa:

$$F^{\mu \nu}= \pmatrix{0&-E_x&-E_y&-E_z\cr
          E_x&0&-B_z&B_y\cr
          E_y&B_z&0&-B_x\cr
          E_z&-B_y&B_x &0}$$

Ahora definimos un nuevo tensor bajando los subíndices como sigue:

$$F_{\mu \nu}= \pmatrix{0&E_x&E_y&E_z\cr
          -E_x&0&-B_z&B_y\cr
          -E_y&B_z&0&-B_x\cr
          -E_z&-B_y&B_x &0}$$

Observemos que al campo eléctrico le hemos cambiado de signo mientras que al magnético no le hemos hecho nada. Esta definición de $F_{\mu \nu}$ entraña otra para las derivadas del  potencial vector con subíndices cambiados de nivel que no vamos a detallar.

Ahora prosigamos en busca de nuestro principio variacional. Si utilizamos la nomenclatura $\partial_\nu $ para designar al operador  $\partial /\partial x^\nu$ y la extendemos también a sus generalizaciones naturales, las  dos leyes que se pueden escribir como $d\sp * F = 4\pi \sp * J$ con $J = 0$ son equivalentes a

$\partial_\nu F^{\mu \nu} =0$

Multipliquemos $\partial_\nu F^{\mu \nu} =0$  por $\delta A_\mu(x,t)$ con la condición de que en los tiempos $s$ y $r$ dichas perturbaciones valgan cero: $\delta A_\mu(x,s) = \delta A_\mu(x,r) =0$

$\partial_\nu F^{\mu \nu} \delta A_\mu(x,t) =0$

Ahora integramos entre $r$ y $s$ en la coordenada temporal, y al mismo tiempo sobre todo el espacio $\Re^3$ en las demás coordenadas espaciales:

$0=\int^s_r \int_{\Re^3} \partial_\nu F^{\mu \nu} \delta A_\mu(x,t) dV dt$

Si tenemos en cuenta que

$\delta(\partial_\nu A_\mu) = \partial_\nu (A_\mu + \delta A_\mu) -\partial_\nu (A_\mu) = \partial_\nu(A_\mu + \delta A_\mu - A_\mu )= \partial_\nu(\delta A_\mu)$

podemos hacer integración por partes, aplicar las evaluaciones nulas en los bordes temporales,  y esto se convierte en

$0=-\int^t_s \int_{\Re^3}  F^{\mu \nu} \delta (\partial_\nu A_\mu(x,t)) dV dt$

pero,  como se demuestra más abajo,

$ F^{\mu \nu} \delta (\partial_\nu A_\mu(x,t))  = (1/2)F^{\mu \nu} \delta F_{\mu \nu}$
$=(1/4)\delta F^{\mu \nu}  F_{\mu \nu}$

entonces reemplazando tenemos

$0=-\int^s_r \int_{\Re^3} (1/2)F^{\mu \nu} \delta F_{\mu \nu} dVdt
= -(1/4) \int^t_s \int_{\Re^3}  \delta F^{\mu \nu} F_{\mu \nu} dVdt$

$0= -(1/4)\delta \int^t_s \int_{\Re^3}  F^{\mu \nu} F_{\mu \nu} dVdt $

con lo cual hemos demostrado que de las dos leyes de Maxwell no geométricas se deduce un punto crítico de un funcional con  lagrangiano $-(1/4)F^{\mu \nu} F_{\mu \nu}$. El camino es reversible si tan sólo se tiene en cuenta que la perturbación en $A_\mu(x,t)$ debe hacerse independientemente coordenada por coordenada.

\

Probemos ahora que

$ F^{\mu \nu} \delta (\partial_\nu A_\mu(x,t))  = (1/2)F^{\mu \nu} \delta F_{\mu \nu}$
$=(1/4)\delta F^{\mu \nu}  F_{\mu \nu}$.

Demostración: Comencemos probando que $ 2 F^{\mu \nu} \delta (\partial_\nu A_\mu(x,t))  = F^{\mu \nu} \delta F_{\mu \nu}$. En efecto:

$F^{\mu \nu} \delta F_{\mu \nu}= F^{\mu \nu} \delta (\partial_\nu A_\mu -\partial_\mu A_\nu)$

$=F^{\mu \nu} \delta (\partial_\nu A_\mu) - F^{\mu \nu} \delta (\partial_\mu A_\nu)$

$=F^{\mu \nu} \delta (\partial_\nu A_\mu) - (\partial^\nu A^\mu -\partial^\mu A^\nu) \delta (\partial_\mu A_\nu)$,

$=F^{\mu \nu} \delta (\partial_\nu A_\mu) + (\partial^\mu A^\nu - \partial^\nu A^\mu ) \delta (\partial_\mu A_\nu)$

$=F^{\mu \nu} \delta (\partial_\nu A_\mu) + F^{\nu \mu} \delta (\partial_\mu A_\nu)$

pero hay que ver que en el segundo término, los subíndices denotan variables mudas y por lo tanto se les puede intercambiar el nombre. Continuando:

$=F^{\mu \nu} \delta (\partial_\nu A_\mu) + F^{\mu \nu} \delta (\partial_\nu A_\mu)$

$=2F^{\mu \nu} \delta (\partial_\nu A_\mu)$

o lo que es lo mismo:

$ F^{\mu \nu} \delta (\partial_\nu A_\mu(x,t))  = (1/2)F^{\mu \nu} \delta F_{\mu \nu}$

\

Por otro lado,

$\delta (F^{\mu \nu}  F_{\mu \nu}) = \delta (F^{\mu \nu})  F_{\mu \nu} +F^{\mu \nu}   \delta (F_{\mu \nu}) =  \delta (F^{\mu \nu})  F_{\mu \nu} +F_{\mu \nu}   \delta (F^{\mu \nu}) =  2  F^{\mu \nu} \delta F_{\mu \nu}$

es decir

$\delta (F^{\mu \nu}  F_{\mu \nu}) = 2  F^{\mu \nu} \delta F_{\mu \nu}$

por lo que

$F^{\mu \nu} \delta F_{\mu \nu} = (1/2) \delta (F^{\mu \nu}  F_{\mu \nu}) $

como ya sabíamos que

$ F^{\mu \nu} \delta (\partial_\nu A_\mu(x,t))  = (1/2)F^{\mu \nu} \delta F_{\mu \nu}$

entonces

$ F^{\mu \nu} \delta (\partial_\nu A_\mu(x,t))  = (1/2)F^{\mu \nu} \delta F_{\mu \nu}$
$=(1/4)\delta F^{\mu \nu}  F_{\mu \nu}$

donde hemos utilizado el siguiente hecho: para bajar o subir los subíndices se multiplica bien por uno o bien por menos uno. Por lo tanto, si se bajan en un lado y se suben en otro, se multiplica siempre por mas uno.

Ahora expliquemos eso de la integración por partes en el espacio tetradimensional .

Tengamos presente que en la expresión

$\int^s_r \int_{\Re^3} \partial_\nu F^{\mu \nu} \delta A_\mu(x,t) dV dt$

se encuentran 16 términos y que la variable subindicada con $\nu$ en algunos casos es  espacial y en otros puede ser el tiempo.

Supongamos que se trata de la variable temporal. El término correspondiente se escribe como

$\int^s_r \int_{\Re^3} \partial_{(t)} F^{\mu (t)} \delta A_\mu dV dt$

donde hemos encerrado a $t$ entre paréntesis para indicar que no hay suma alguna por estar un subíndice una vez arriba y otra abajo. Intercambiando el orden de integración

$\int_{\Re^3} \int^s_r \partial_{(t)} F^{\mu {(t)}} \delta A_\mu dt dV$

podemos tomar el término interior  $ \int^s_r \partial_{(t)}F^{\mu {(t)}} \delta A_\mu dt$  e integrarlo por partes tomando $u$ como $\delta A_\mu$ y $dv$ como $\partial_{(t)} F^{\mu {(t)}}dt$ y nos queda

$\int_{\Re^3} \int^s_r \partial_{(t)} F^{\mu {(t)}} \delta A_\mu dt dV = (\delta A_\mu)(F^{\mu {(t)}})|^s_r-$
$\int^s_r  F^{\mu {(t)}} \delta (\partial_{(t)} A_\mu )dt$

pero como $\delta A_\mu(x,r)=\delta A_\mu(x,s)=0$ el término sin integral desaparece y obtenemos:

$\int_{\Re^3} \int^s_r \partial_{(t)} F^{\mu {(t)}} \delta A_\mu dt dV = $
$-\int^s_r  F^{\mu {(t)}} \delta (\partial_{(t)} A_\mu )dt$

Reinsertando este término de donde fue sacado nos queda finalmente que

$\int^s_r \int_{\Re^3} \partial_{(t)} F^{\mu (t)} \delta A_\mu dV dt
=-\int^s_r \int_{\Re^3}  F^{\mu {(t)}} \delta (\partial_{(t)} A_\mu) dV dt$

Pero si  el subíndice $\nu$ es una variable espacial, digamos $x$, razonemos como sigue:

$\int^s_r \int_{\Re^3} \partial_\nu F^{\mu \nu} \delta A_\mu dV dt = $
$\int^s_r \int_{\Re^2} \int_{-\infty}^{\infty} \partial_{(x)} F^{\mu {(x)}} \delta A_\mu dxdA dt$

y entonces podemos proceder como en el caso temporal siempre y cuando pongamos la condición adicional de que en el infinito $\delta A_\mu$ también es cero. Eso no me parece pedir demasiado. (O puede recurrirse al lenguaje de las funciones generalizadas e integrar en el sentido de las distribuciones).

Es relajante darse cuenta que detrás de todo el complique exhibido  hay  algo que suena familiar:

 $F^{\mu \nu} F_{\mu \nu}=-2(E^2-B^2)$

Antes de probar  esa identidad, debemos anotar que no se trata de una multiplicación matricial. Podemos entenderla como el producto punto entre dos vectores $F^kF_k$, donde $k = \mu \nu$. De esa forma:

$F^{\mu \nu} F_{\mu \nu}=
0-E^2_x- E^2_y-E^2_z
          -E^2_x-0+B^2_z+B^2_y
          -E^2_y+B^2_z+0+B^2_x
          -E^2_z+B^2_y+B^2_x +0$

$ = -2(E^2-B^2)$

 Para ser totalmente sinceros, esa igualdad es más una definición que cualquier otra cosa y para llegar a ella fue que se definió $F_{\mu \nu}$ a partir de $F^{\mu \nu}$ tal como fue hecho. Por qué se hace así? Porque andábamos  buscando un principio variacional, es decir un funcional cuyo punto crítico sea equivalente a las leyes que queremos explicar.  Sin embargo,  por un punto en el plano pueden pasar miles de parábolas que tengan su mínimo en dicho punto y de igual manera las leyes de Maxwell podrían ser puntos críticos de miles de funcionales. Por lo tanto, no se trata de buscar uno cualquiera. Se requiere poner condiciones. La más natural es que el Lagrangiano encontrado sea un invariante relativista. Eso quiere decir que desde todos los sistemas inerciales el Lagrangiano sea exactamente el mismo, o mejor dicho, que sea invariante  ante transformaciones de Lorentz. Pues bien, a partir de los campo magnético y eléctrico sólo hay dos invariantes. El primero es $E^2 -B^2$  y el segundo es $E\cdot B$, entendiéndose como producto punto de vectores. De esas dos posibilidades, la que funciona en el proyecto variacional es la primera.

\section{EL GAUGE DE FEYNMAN}

Por razones de la teoría de campos, se modifica el lagrangiano

 $L = -(1/4)F^{\mu \nu}F_{\mu \nu}$ 
 
 introduciendo  un $\alpha$ artificial quedando el siguiente lagrangiano:

$$L = -(1/4)F^{\mu \nu}F_{\mu \nu} - (\alpha/2)(\partial \cdot \vec A)^2$$

este $\alpha$ artificial  dañará las predicciones de este lagrangiano a menos que desaparezca de algún modo:  después de tomar las ecuaciones de Euler Lagrange, se pone $\alpha =1$ y de esa forma hay paz tanto con la teoría como con el experimento. A esta selección se le llama el \index{gauge de Feynman } \textbf{gauge de Feynman}.

\section{COHOMOLOGIA}

Existe una estructura matemática, llamada \index{cohomología} \textbf{cohomología}, que juega un papel importante en las discusiones avanzadas. De lo que hemos visto, ya podemos saber de qué se trata (las definiciones rigurosas pueden mirarse en un libro de geometría).

La teoría electromagnética expresa una interacción que depende de una 1-forma, la cual no es la diferencial exterior de una función escalar, de otra manera, $d\sp 2 = 0$ predeciría un tensor de campo nulo.

En cohomología, se dice que una forma $\alpha$ es exacta cuando su diferencial es cero: $d\alpha = 0$, que es lo mismo que decir que una forma es exacta cuando pertenece al núcleo del   diferencial. Cuando decimos que la diferencial de la 1-forma que describe el potencial vector  crea el campo electromagnético, $dA = F$, se dice en cohomología  que $A$ no pertenece al kernel de $d$, y que es un elemento no nulo del grupo de cohomología $H^1$.

\bigskip

Algo parecido ocurre con la gravitación.

\section{CONCLUSION}

La formulación propuesta por Grassmann y refinada por Cartan nos ha permitido reescribir las leyes de Maxwell en un formalismo con gran poder de concretar los resultados. Todo se deduce de asumir que las leyes del magnetismo se derivan de la existencia de una única entidad: el  potencial $A$, el cual es una 1-forma, pues describe el campo que afecta nuestro experimento, digamos a una partícula que hemos puesto a moverse.

El tensor de campo $F$ es $F=dA$

Como $d\sp 2=0$, entonces tenemos que para un sistema aislado, en el cual no haya corrientes, debe cumplirse que $dF=0$ de lo cual se deduce las dos primeras leyes de Maxwell. Definiendo el tensor dual $\sp * F$, las otras dos leyes de Maxwell quedan en forma compacta:

$ d\sp * F= 4 \pi \sp * J$

La ley de Lorentz quedó muy natural: una carga eléctrica en movimiento en presencia de un campo electromagnético cambia su energía-momento 'proporcionalmente' a $F$, lo cual, usando coordenadas, se escribe:

$dp\sp \mu/d\tau = (e/mc)\sum \sb \nu p\sb \nu F\sp{\mu \nu}$

La libertad Gauge queda inmediata a partir de $d\sp 2=0$:

$d(A+ d\phi)=dA + d\sp 2 \phi = dA$, donde $A$ es la 1-forma potencial y $\phi$ es una función escalar cualesquiera. La invariancia gauge dice entonces: las leyes del electromagnetismo en la versión relativista son invariantes ante un suma cualesquiera de una diferencial exacta de una función escalar $d\phi$, (la cual es una 1-forma).

\section{REFERENCIAS}

Las dos primeras referencias están orientadas hacia la física, la tercera es más matemática gozando al mismo tiempo de mucho sentido común.

1. Arnold, V.I. \textit{Mathematical Methods of Classical Mechanics}, Springer, 1978.

2. Misner, C., Thorne, K., and Wheeler, J. \textit{Gravitation, Freeman}, 1970.

3. Frankel T. \textit{The Geometry of Physics, An Introduction}. Cambridge University Press, 2001.

\chapter{GEOMETRIA Y GRAVITACION }        

\date{}

\Large

\centerline{RESUMEN}

\bigskip

El objetivo de la física moderna es describir las interacciones fundamentales de manera irreducible. Teniendo cuatro interacciones fundamentales y cuatro (familias de) teorías para describirlas, es imposible saber si hay conexiones entre ellas a menos que todas ellas se describan en un lenguaje común. Por eso,   en la presente sección  formularemos, parcialmente, una descripción geométrica  de la interacción gravitatoria, en parte para conocerla, en parte porque se cree que la \index{geometría} \textbf{geometría} es la lingua franca que permitirá la unificación de todas las interacciones.

Una descripción geométrica es aquella que, en suma, describe una interacción elemental como un elemento creador de fuerza, es decir, de cambio de curvatura en las trayectorias, lo cual se debe a una influencia de la curvatura del espacio-tiempo  sobre la partícula en estudio.

Necesitaremos el concepto de \index{geodésica} \textbf{geodésica}, que es el camino más corto a lo largo de un espacio curvo entre dos puntos muy cercanos, y también el de conexión o derivada covariante, la cual es un tipo de derivada que sólo contabiliza aquello que es de importancia geodésica. A dicha derivada se le asociará una segunda derivada o curvatura   teniendo en mente que a más curvatura del espacio esperamos  una mayor curvatura de las trayectorias, o como quien dice, una mayor fuerza gravitatoria. La relación con el  electromagnetismo  será estudiada en el próximo capítulo.

\
\
\bigskip
\bigskip
\bigskip
\normalsize
\color{black}
\section{INTRODUCCION}

El gran sue\~{n}o de la física-matemática es el de demostrar que las 4 interacciones que se registran en el laboratorio son casos particularizados  de una interacción general producida por un único campo. A tal sue\~{n}o se le denomina la teoría del \index{campo unificado} \textbf{campo unificado}. Un prerrequisito para lograr tal unificación es demostrar que todas las interacciones se pueden describir por el mismo formalismo matemático. Ha resultado muy inspirador el hecho de que la versión geométrica de  la gravitación, o relatividad general, haya resultado muy exitosa y que además otras interacciones también se hayan podido reformular en el mismo lenguaje. Por eso, es imposible entender la formulación geométrica del electromagnetismo sin entender las ideas básicas subyacentes a la \index{geometrización de la gravedad}

\textbf{geometrización de la gravedad}.

En esta sección conoceremos los principios básicos que permiten entender cómo se entreteje la red geométrica que permitió cazar a la gravitación y cuyo mejoramiento permitirá atrapar, en la próxima sección, también al electromagnetismo. No estaremos en capacidad de ser rigurosos ni de demostrar los grandes teoremas, lo cual exige bastante más que una introducción, pero si podremos entender algunos detalles importantes del funcionamiento de la maquinaria geométrica.

Existen variadas versiones, no siempre equivalentes, de la gravitación. La nuestra es  estándar, ideológicamente muy simple, y puede edificarse sobre cinco pilares como sigue:

1)El experimento de Galileo a raíz de la discusión sobre si un cuerpo más pesado caería más rápido que uno liviano. Todo el mundo decía que sí, pero Galileo dijo que no y lo demostró: llevó a sus amigos a la torre de Pisa y dejaron caer dos piedras, una grande y otra pequeña. Cayeron al tiempo. Para que esto suceda se requiere que la masa gravitatoria sea proporcional a la masa inercial, es decir que la resistencia que una masa ofrece a ser acelerada por un empujón sea proporcional a la que ofrece a ser acelerada por un campo gravitatorio. Como la tal constante de proporcionalidad no depende de nada, se acostumbra a tomar como uno. A esa igualdad se le llama el \index{principio de equivalencia } \textbf{principio de equivalencia }(entre la masa gravitatoria y la inercial).

2)La precesión de la órbita de \index{Mercurio} \textbf{Mercurio} que consiste en que la órbita de este planeta alrededor del sol no es una elipse como lo predecía las leyes de Kepler sino que era casi una elipse que no alcanzaba a cerrarse del todo de tal forma que la órbita de Mercurio parece describir una roseta. Ante este fracaso de la teoría aceptada, la cual era una teoría de acción a distancia, era necesario formular otros puntos de vista.

3)La alternativa de Einstein: puesto que desde un mismo punto todos los cuerpos caen cerca de la tierra con la misma aceleración, la fuerza gravitatoria puede reinterpretarse como una propiedad de la interacción del espacio mismo con el cuerpo que cae y no del cuerpo con la tierra. Acá es donde entra la geometría con interacciones locales en vez de interacciones a distancia como en la teoría clásica. Por supuesto, el espacio mismo no es el 3:d de la teoría clásica sino algo más complicado que en el caso límite en ausencia de masas se convierte en el espacio -tiempo. Einstein introdujo una convención que usaremos en esta sección porque ahorra mucho trabajo de tipografía: siempre que haya un índice repetido, una vez arriba y otra abajo, se entenderá una suma sobre todos los valores permitidos del índice. Por ejemplo, un vector en $\Re \sp n$ se escribirá como $\vec x = x\sp i \vec e \sb i$

4)El \index{vuelo del pintor}  \textbf{vuelo del pintor}. Einstein trabajaba en una oficina de patentes, donde le quedaba mucho tiempo libre para pensar. Según la leyenda, cierto día llegó alguien a pintar la fachada de una casa al frente de su oficina. Y de pronto Einstein vio que el pintor se caía  desde su  andamio. Toda caída produce heridas y traumas, reflexionó Einstein. Pero eso no sucede antes del choque con la tierra. Lo importante, por tanto, es reconocer que antes del golpe, el pintor volaba y mientras que volaba no sentía ninguna fuerza, aparte de una minúscula resistencia del aire. Su bitácora de vuelo estaba perfectamente determinada y la denominamos geodésica. Podríamos decir que la gravitación es una fuerza ficticia y tan irreal como la fuerza centrífuga: no hay un experimento local que una sola partícula aislada pueda hacer para saber si está cayendo hacia un hueco negro o si está en movimiento rectilíneo uniforme o si está quieta. La soledad y la puntitud de la partícula son ambas importantes: un campo gravitacional intenso podría descuartizar un experimentador humano, debido a que la cabeza sería muchísimo menos atraída que los pies. De manera similar, dos partículas podrían saber si hay un campo gravitatorio externo midiendo en diversas circunstancias la forma como varía la distancia entre ellas.

5)La medición de Gauss: debido al estudio de la curvatura, al que Gauss era tan aficionado, el se preguntó hacia 1850,  cómo se podría averiguar si el espacio en que vivimos era curvo. El problema es reminiscente del de medir la curvatura de la tierra. Las estrellas han sido utilizadas desde la antigüedad, según cuenta Platón,  como referentes exteriores que denuncian la \index{curvatura} \textbf{curvatura} de la tierra en dirección  transversal al ecuador: desde Egipto no se veían las mismas estrellas que desde Macedonia. Eso es una descripción no intrínseca. Pero nosotros podemos demostrar que la tierra es curva desde la misma tierra: sus montañas  se ocultan al uno ir  abandonando la costa y adentrándose en el mar. A eso se llama una descripción intrínseca.  En referencia al espacio mismo, lo que Gauss hizo fue medir los ángulos internos de un triángulo cuyos vértices eran tres picos montañosos. El resultado que obtuvo fue que la suma le dió 180 grados aparte de errores inherentes al proceso de medición. Por eso, concluyó que nuestro mundo era plano.  El análisis de esta medición ya contiene lo que hoy día llamamos descripción intrínseca (sin referentes exteriores), transporte paralelo y seguimiento de geodésicas, conceptos que con los de curvatura y derivada covariante permiten una formulación cuantitativa de la versión geométrica de la gravitación.

6) La teoría de gravitación de Einstein es un hito entre las creaciones científicas, pero con todo, vista desde la perspectiva presente, la teoría de la gravitación de Einstein es similar a una receta de cocina casera hecha por un gran chef. Pero sucede que después de los chefs siempre vienen los expertos en dietética enseñando cómo se diseñan  recetas confeccionadas científicamente. En referencia con la gravitación, baste decir, por ahora, que para que hoy en día algo se considere bien hecho debe estar basado ab initio en la teoría de grupos y no sólo en ideas físicas que se amalgaman de forma inteligente. 

\section{CURVATURA EN VEZ DE FUERZA}

La órbita de la luna alrededor de la tierra es curva. Podría ser recta, escapándose eternamente lejos de nosotros. Pero es curva por el efecto gravitacional de la tierra: la curvatura de trayectorias es la marca de fábrica de la gravitación. Entre más fuerza, más curvatura de las trayectorias. Menos fuerza, mas suavidad de las trayectorias. Cero fuerza, cero curvatura.

En una órbita elíptica de un planeta uno ve la curvatura. Y en una piedra que cae? También, sólo que hay que involucrar al tiempo:

Cuando uno gráfica la posición de un cuerpo que cae libremente, en línea recta, poniendo en un eje el tiempo y en el otro la posición, la gráfica es una parábola y su curvatura está dada por la segunda derivada: mientras que la primera derivada mide la pendiente, la segunda mide la forma como la primera derivada cambia, o sea cómo cambia la pendiente, es decir, la segunda derivada mide la curvatura. Para el caso de un cuerpo que cae cerca del nivel del mar, la curvatura (segunda derivada del espacio recorrido contra el tiempo) es siempre la misma, proporcional a la gravedad, $9.8 m/s\sp 2$.

La idea de la geometrización es  la de quitar el sentido fundamental dado a las fuerzas y concedérselo a la curvatura.

Veremos que desde el punto de vista cuantitativo la \index{curvatura vs métrica} \textbf{curvatura depende de la métrica} así: una métrica es una forma de medir distancias. Primero se define la métrica para poder definir la curva que traza  el camino más corto entre dos puntos infinitesimalmente cercanos, a la cual se le llama geodésica. A dichas curvas, lo mismo que a cualquier otra, se le puede medir la curvatura y, si dicha curvatura no es nula, entonces uno puede decir, por definición,  que  hay una fuerza  en acción.

Por supuesto que el espacio en que operamos no es el 3:d sino el 3+1:d con la métrica del espacio-tiempo, en la cual  se encuentra impresa la inseparabilidad de las 4 coordenadas y por ende la posibilidad de que una de ellas influya en cualquiera de las demás. Pero para implementar esa unidad hay que darse cuenta de que eso no sería posible si medimos el tiempo en segundos y el espacio en metros. Necesitamos la misma unidad para medir tanto el tiempo como el espacio. Uno podría elegir como unidad el dólar: 3 km equivale a 1 dólar pues es lo que me cuesta recorrerlo en taxi. Igualmente, media hora podría valer 24 dólares, pues es lo que alguien cobraría por trabajar media hora. Pues bien, la unidad que más sencillez produce es la que se define por la ecuación:

\addtocounter{ecu}{1}

$$ds\sp 2 = -c\sp 2 dt\sp 2 + dx\sp 2 + dy \sp 2 + dz\sp 2 \eqno{(\theecu)}$$

en esta ecuación, $ds$ significa 'el intervalo' que separa dos eventos en el espacio tiempo, pues evento es el nombre que se le da a dos puntos en el espacio de Minkowski. Observamos que el tiempo se mide en unidades de distancia, pues va acompañado de la velocidad de la luz, y velocidad por tiempo da espacio.

\section{DETALLE TECNICO }

Amable lector:  gracias al trabajo de varias generaciones de matemáticos y físicos podemos gozar de una formulación de la geometría diferencial que es poderosa, consistente, fructífera. Y aunque no lo parezca, muy intuitiva. Lleva muchos ejercicios familiarizarse con los conceptos básicos  y por eso es mejor remitirse a cualquier libro de geometría y a un curso formal. Con todo, repasemos lo más básico del formalismo.

El concepto de \index{variedad} \textbf{variedad} es básico y su definición, que ya vimos, describe  simplemente el proceso de hacer un mapa de la tierra: se hace el mapa de una pequeña región, después de otra adjunta y se verifica que en la intersección haya perfecta compatibilidad. Se sigue el proceso hasta acabar. Si pensamos en una esfera cada vez que se nos habla de variedad, podremos intuir de qué se habla, al menos en lo que se refiere a las presentes notas. Será necesario también mantener presente que al caer una piedra en línea recta en el espacio 3:d, en el espacio-tiempo ( que podemos entender como una gráfica de espacio contra tiempo ) describe una trayectoria curva.

Deseamos introducir el concepto de \index{vector tangente} \textbf{vector tangente} sobre una variedad, el cual ha de generalizar la idea de vector tangente sobre una esfera. La definición se basa en la \index{dualidad} \textbf{dualidad} que existe en  $\Re \sp n$ entre vectores y operadores de derivación direccional. Todos los vectores tangentes a una variedad en un punto dado forman un espacio vectorial, que generaliza la definición de plano. A dicho espacio se le llama espacio tangente en el punto dado. La reunión de todos esos espacios se llama el espacio tangente. Sobre la esfera, uno se la imagina  con todos los planos tangentes posibles, los cuales están formados de vectores tangentes.

Si deformamos la esfera, la cual nos la imaginamos plástica, se crea una superficie quizá compleja. Decimos que la deformación es \index{diferenciable} \textbf{diferenciable} cuando los planos tangentes de la esfera al ser deformados con ella se convierten en los planos tangentes a la nueva superficie. De forma equivalente, la deformación es diferenciable si la deformación inducida de cada vector tangente se convierte en un vector tangente a la nueva superficie. A la deformación inducida entre vectores  se le llama la \index{diferencial} \textbf{diferencial} y si la deformación es $F$ a su diferencial se le  nota $dF$ o también $F\sb *$. Todo se extiende naturalmente a variedades. La diferencial de una función escalar es una 1-forma.

Mientras que la diferencial fue definida con el objetivo de trazar la suerte de los vectores tangentes ante deformaciones, se necesita una maquinaria que haga lo mismo pero para las formas.  De eso se encarga el \index{pull-back} \textbf{pull-back }(puede traducirse como 'halado hacia atrás', pero en realidad no se acostumbra a traducir). La diferencial opera hacia adelante (toma un vector tangente y lo transforma en otro), mientras que el pull-back opera hacia atrás: toma una forma $\phi$  y calcula aquella forma que al sufrir la deformación inducida por $f$ se convierte en $\phi$. Toda la maquinaria funciona con un aceite que se llama 'dualidad'. Veamos como lucen algunos detalles de las definiciones. Estas definiciones se usan  a lo largo del texto sin previo aviso, e incluso ya se han usado anteriormente, pero muy intuitivamente.

\addtocounter{ecu}{1}

\textit{Definición (\theecu): Un conjunto $M$ se llama \index{variedad} \textbf{variedad} si puede cubrirse con retazos de $\Re\sp n$, y cada retazo lleva un mapa de un pedazo de $M$, y si hay compatibilidad entre los mapas. Concretamente: $M=U \sb 1 \cup  U\sb 2\cup .. ..$ donde cada $U\sb i$ está en una correspondencia 1:1 $\phi \sb U: U\rightarrow \Re\sp n$ con un abierto (sin frontera) $\phi(U)$ de $\Re\sp n$. Para cumplir la compatibilidad se requiere que en cada intersección $U \sb i \cap  U\sb j$ se cumpla : $\phi(U \sb i \cap  U\sb j)$ debe ser un abierto y además que sobre cada intersección la función $\phi\sb{U \sb i} \circ \phi \sb{  U\sb j}\sp{-1}$ sea diferenciable. A las funciones $\phi$, las llamamos mapas.}

\bigskip

\addtocounter{ecu}{1}

\textit{Definición (\theecu): Cada \index{vector tangente} \textbf{vector tangente} lo confundimos con un operador de derivación parcial que opera sobre funciones escalares definidas sobre la variedad, del tipo $f: M \rightarrow \Re $. Cada vector $\vec v $ es de la forma  $\vec v = v\sp i \partial / \partial x\sp i$ y al operar sobre $f$ produce $v\sp i \partial f(\phi(\vec x)) / \partial x\sp i$, donde $\vec x\in \Re\sp n$ tiene por coordenadas naturales $\vec x= (x\sp 1,.. ..,x\sp n)$. Por lo tanto, el operador realmente va de $\Re\sp n$ a $\Re\sp n$, pero utiliza la variedad como lugar de tránsito.}

\bigskip

\addtocounter{ecu}{1}

\textit{Definición (\theecu): El \index{espacio tangente} \textbf{espacio tangente} sobre un punto $P$ de la variedad $M$ es el conjunto de todos los vectores tangentes a ese punto. Tales vectores se pueden sumar, restar, multiplicar por una constante y siguen siendo vectores. Por eso, son en realidad un espacio vectorial. A la reunión de todos esos espacios se le llama el espacio tangente.}

\bigskip

\addtocounter{ecu}{1}

\textit{Definición (\theecu): Tomemos dos variedades $M, N$ cada una con sus mapas, $\phi$ para $M$ y $\rho$ para $N$.  Lo que nosotros llamamos deformación es una función entre variedades $F : M\rightarrow N$. Decimos que la función es \index{diferenciable}  \textbf{diferenciable} si las funciones compuestas del tipo: $ \phi \circ F \circ \rho \sp{-1}$ que van de subconjuntos de $\Re\sp n$ a $\Re\sp n$ son diferenciables. Sean $M\sb P$ el espacio tangente de $M$ en $P$ y $N\sb{F(P)}$  el espacio tangente de $N$ en $F(P)$. Los vectores tangentes $\vec v$ de $M \sb P$ operan sobre funciones escalares $f$ definidas sobre $M$, mientras que los vectores $\vec w$ de $N\sb{F(P)}$ operan sobre funciones $g$ definidas sobre $N$. Toda la maquinaria gira en torno al hecho que $g(F())$ es una función escalar sobre la cual  puede operar un vector de $M\sb P$. Si $F$ es diferenciable, llamamos diferencial de $F$ en $P \in M$ a la transformación lineal $dF = F\sb * : M\sb P \rightarrow  N\sb{F(P)}$ que opera sobre vectores tangentes así:}

\textit{$dF(\vec v)=\vec w$ tal que $\vec w(g)= \vec v(g\circ F)$ donde $g: N\sb{F(P)} \rightarrow \Re$ }

\

Ahora elaboraremos el concepto de pull-back de una forma diferencial.

 Acuñar la definición de \index{pull-back} \textbf{pull back }de formas diferenciales por medio de $F : M\rightarrow N$ requiere tener en mente los siguientes hechos: 0) Obsérvese que si la curva tiene entrecruzamientos, entonces a dos vectores en el espacio tangente del dominio, $\Re$, pueden corresponderle dos vectores diferentes en el mismo punto del espacio tangente del codominio. Eso quiere decir que en general es conflictivo hablar del pull-back de un vector. Pero se puede hablar del pull-back de una forma diferencial: 1) una n-forma opera sobre n-plps, los cuales los hemos tomado como conjuntos ordenados de vectores del espacio tangente en un punto dado. 2) El pull back toma una n-forma $\sigma$ que opera sobre plps en el espacio tangente  en el codominio y fabrica una n-forma  $\omega$ que pueda transformarse en ella por medio de $F$. 3)  $dF$ transforma, hacia adelante, vectores en vectores y por consiguiente, también transforma plps en plps, también hacia adelante. Procediendo tenemos:

\addtocounter{ecu}{1}
\textit{
Definición (\theecu): Sea la deformación $F : M\rightarrow N$, $dF$ su diferencial, $[\vec v\sb 1,.. .., \vec v\sb n]$ un n-plp basada en algún punto del dominio, $M\sb P$. El \index{pull-back} \textbf{pull-back} de $F$ se nota $F\sp *$ y su trabajo es retro-convertir formas diferenciales en formas diferenciales: si $\sigma$ es una forma diferencial basada en el codominio, entonces  $F\sp *(\sigma) = \omega$ está basada en el domino, tal que $\omega ([\vec v\sb 1,.. .., \vec v\sb n])= \sigma([dF(\vec v\sb 1),.. .., dF(\vec v\sb n)])$.}

\bigskip

Todas estas definiciones generalizan los conceptos ordinarios de \index{cambio de variable} \textbf{cambio de variable}. Como ejemplo de un cálculo, consideremos una trayectoria parametrizada $c(t)$ definida por una función $c: \Re \rightarrow \Re\sp 3$. Primero expliquemos el espacio tangente junto con las formas (que forman el espacio cotangente) y luego el cambio de variable.

En el cálculo vectorial, al elemento de longitud en $\Re$ se nota $dt$ y se toma como un infinitesimal. El vector velocidad se puede notar  $dc/dt$, y el elemento de desplazamiento  sobre la curva, un vector, es $(dc/dt)dt$.  En el lenguaje que usamos todo eso se relee como sigue: El vector $\partial/\partial t = d/dt$ denota un operador de derivación. La expresión 'dt' denota una 1-forma sobre el dominio, la 1-forma que operada sobre $d/dt$ da uno: dt(d/dt) = 1. Por eso se dice que $dt$ es la forma dual de $d/dt$. La función $c$ tiene una diferencial, $dc$, la cual toma vectores en el espacio tangente del dominio y los transforma en vectores en el  espacio tangente del codominio. Tenemos que si $c(t)= (x\sp 1, x\sp 2, x\sp 3)$ entonces, $dc(d/dt)$ tiene que ser un vector del espacio tangente de la curva, la cual está en $\Re \sp 3$. Pero dicho espacio tangente es 1-dimensional. Veamos: $dc(d/dt)$ es tal que $dc(d/dt)$ como operador de derivación aplicado sobre funciones $f: \Re \sp 3 \rightarrow \Re$ produce

$dc(d/dt)(f(\vec x)) = (d/dt)(f(c(t))$

$= (\partial f/\partial x\sp 1) (dx\sp 1/dt)  + (\partial f/\partial x\sp 2) (dx\sp 2/dt)  + (\partial f/\partial x\sp 3) (dx\sp 3/dt)$

$= (dx\sp 1/dt)(\partial f/\partial x\sp 1)   + (dx\sp 2/dt) (\partial f/\partial x\sp 2)  + (dx\sp 3/dt)(\partial f/\partial x\sp 3) $

$(dx\sp 1/dt)(\partial /\partial x\sp 1)   + (dx\sp 2/dt) (\partial /\partial x\sp 2)  + (dx\sp 3/dt)(\partial /\partial x\sp 3)(f)$

donde las evaluaciones se ejecutan en el punto $\vec x$.

Por lo tanto, concluimos que  $dc(d/dt) $ es el vector cuya representación por operador  de derivación direccional es:

$dc(d/dt) =  (dx\sp 1/dt)(\partial /\partial x\sp 1)   + (dx\sp 2/dt) (\partial /\partial x\sp 2)  + (dx\sp 3/dt)(\partial /\partial x\sp 3)$

\

Podemos decir que

$dc(d/dt) = c'(t)$

donde $c'(t)$ es el vector velocidad a la trayectoria dada por $c$.

\

El espacio tangente a la curva en el punto $c(t)$ es el conjunto de vectores de la forma

$dc(\alpha d/dt) =$
$\alpha (dx\sp 1/dt)(\partial /\partial x\sp 1)   + \alpha (dx\sp 2/dt) (\partial /\partial x\sp 2)  + \alpha (dx\sp 3/dt)(\partial /\partial x\sp 3)$

que es realmente unidimensional. Se ruega al lector no pensar en utilizar aquí el 'vector gradiente' sin enterarse primero de que el gradiente realmente no es un vector sino una 1-forma (esto tiene que ver con reglas de cambios de base).

\

Para definir una forma, primero uno fija un punto, sube al espacio tangente (la aproximación plana a la variedad en dicho punto), y sobre el espacio tangente forma un plp. Una forma le hace corresponder al plp un número real. Sobre las aristas de los plps, las formas son multilineales alternadas, como los determinantes, de donde ellas salen.

Mientras que una forma en el dominio de $c: \Re \rightarrow \Re^3$,  es una 1-forma, cuya expresión general es $\lambda dt$, en el codominio,  todo $\Re^3$, una 1-forma tiene la expresión $\alpha dx + \beta dy + \gamma dz$. Como se ve, tiene 3 grados de libertad sobre cada punto. Pero podemos restringirnos a la curva, cuyo espacio cotangente opera sobre vectores tangentes, los cuales forman en cada punto un espacio vectorial de dimensión uno. Por ello, las formas que operan sobre el espacio tangente a la curva deben formar, sobre cada punto, un espacio de dimensión uno. Es decir, deben ser un subespacio del formado por $\alpha dx + \beta dy + \gamma dz$.

 Tratemos de entender ahora el pull-back de  una 1-forma $f$ sobre el espacio tangente de la curva en un punto dado. Tal 1-forma opera sobre vectores tangentes sobre la curva. La 1-forma $f$ puede descomponerse en la base de las 1-formas de $\Re \sp 3$

 $f= f \sb 1 dx\sp 1 + f\sb 2 dx\sp 2 + f\sb 3 dx\sp 3$

$f$ opera sobre un vector como sigue:

$ f  (\alpha (dx\sp 1/dt)(\partial /\partial x\sp 1)   + \alpha (dx\sp 2/dt) (\partial /\partial x\sp 2)  + \alpha (dx\sp 3/dt)(\partial /\partial x\sp 3))  $

 $=f\sb 1 \alpha (dx\sp 1/dt)   + f\sb 2 \alpha (dx\sp 2/dt)  + f\sb 3 \alpha (dx\sp 3/dt)$

Ahora bien, el pull-back de $f$, $c\sp * (f)$, es una 1-forma  sobre el espacio tangente del dominio, es decir sobre $\Re$, tal que al operar sobre la base del espacio tangente del dominio   $   d/dt$ debe dar un número, el cual, por definición es igual a:

$  [c\sp * f](d/dt ) = f (dc(d/dt)) =  f (c'(t)))$

Explícitamente,

$   [c\sp * (f)](d/dt )=   f (
  (dx\sp 1/dt)(\partial /\partial x\sp 1)   + (dx\sp 2/dt) (\partial /\partial x\sp 2)  + (dx\sp 3/dt)(\partial /\partial x\sp 3))$

$=  (f \sb 1 dx\sp 1 + f\sb 2 dx\sp 2 + f\sb 3 dx\sp 3) (
  (dx\sp 1/dt)(\partial /\partial x\sp 1)   + (dx\sp 2/dt) (\partial /\partial x\sp 2)  + (dx\sp 3/dt)(\partial /\partial x\sp 3))$

 $=[f\sb 1  (  dx\sp 1/dt)   + f\sb 2  (\alpha dx\sp 2/dt)  + f\sb 3  (\alpha dx\sp 3/dt)]$

$= \vec f \cdot \vec dc/dt $

donde hemos asociado el vector $\vec f = (f\sb 1 ,f\sb 2 ,f\sb 3 )$ a la 1-forma $f$.

Toda esta maquinaria se ha hecho popular porque unifica los grandes teoremas del cálculo vectorial y además permite sustanciales generalizaciones al teorema de Stokes. Y es muy elegante: veamos cómo se calcular el trabajo hecho en contra de un campo de fuerzas $\vec F$ a lo largo de la trayectoria dada por $c$. En el lenguaje antiguo:

Trabajo = $\int_\gamma \vec F.\vec{ds}$ = $\int^{t_1}_{t_o} (\vec F.\vec{dc/dt})dt$

Para reescribir esto en el nuevo lenguaje, representemos la fuerza $\vec F$ por una 1-forma: $f^1 = f\sb 1 dx\sp 1 + f\sb 2 dx\sp 2 + f\sb 3 dx\sp 3$ y lo que antes se notaba como

Trabajo = $\int_\gamma \vec F.\vec{ds}$

ahora se nota como

Trabajo = $\int_\gamma f $ = $\int_c f$ pues el camino $\gamma$ está parametrizado por $c$, mientras que

  $\int^{t_1}_{t_o} (F.\vec{dc/dt})dt$

corresponde a

 $\int c\sp * f $

La conclusión es que el trabajo hecho en contra de $\vec F$ a lo largo del camino $\gamma$ se expresa  como

$\int_c f =\int c\sp * f $

donde se sobreentiende que se calcula entre los límites adecuados.

No hemos hecho una teoría de integración de formas, pero debe ser claro que todo se define con el objetivo de que la nueva notación exprese los mismos conceptos que la vieja. Retornemos ahora al estudio de la gravitación.

\section{GEODESICAS}

\textit{Una \index{geodésica}  \textbf{geodésica }es el camino más corto entre dos lugares 'cercanos'}, en la jerga se dice: localmente. Para lugares alejados, hay que considerar el funcional que da la longitud de arco, encontrar un camino que de un valor crítico de dicho funcional y después seleccionar el mínimo, si existe. Intuitivamente, podemos ver que sobre el espacio ordinario, las geodésicas están dadas por  rectas. Sobre la tierra, o una esfera, las geodésicas corresponden a círculos máximos, cuyo centro coincide con el de la tierra. Así, si un avión tiene que ir de una ciudad a otra y ambas quedan sobre el paralelo 45, entonces el avión tiene que tomar primero hacia el norte y después hacia el sur, a medida que el círculo máximo lo lleve.

\subsection{Geodésicas del plano}

En la métrica ordinaria, las geodésicas en el plano son líneas rectas. Es necesario demostrar eso con nuestra maquinaria para poder estar seguros de que vamos bien.  

La distancia $ds$ en el plano está dada por Pitágoras:

$ds = \sqrt{dx^2 + dy^2}$

si además $y$ es localmente función diferenciable de $x$, tenemos

$ds = \sqrt{1+(dy/dx)^2} dx$

lo cual se acostumbra a escribir como

$S = \int^{x_2}_{x_1} \sqrt{1+(f'(x))^2} dx$

Nuestro problema es minimizar este funcional. Para ello, primero debemos encontrar una ecuación que determine un extremo y después hay que probar que da un mínimo.

Nuestro problema  ya ha sido considerado en forma general: debemos  minimizar el funcional dado por

$S =  \int^{x_2}_{x_1} L dx$

por lo que podemos aplicar las ecuaciones de Euler Lagrange para hallar un camino que de un valor crítico del funcional:

$\frac{\partial}{\partial y} L - \frac{d}{d x} \frac{\partial}{\partial y'} L = 0$

Sucede que en este caso $L$ no depende de $y$. Cuando eso sucede el problema se resuelve usando la identidad de Beltrami, que entramos a enunciar y a probar.

\bigskip

\addtocounter{ecu}{1}

Teorema (\theecu): Identidad de Beltrami.
Sea $y(x)$ una función $\Re \rightarrow \Re $. Notamos la derivada de $y$ con respecto a $x$ como $\dot{y}$ y sea $L = L(y,\dot{y},x)$. Para hallar un extremo de $L$ se cumple la ecuación de Euler Langrange

\bigskip

$\frac{\partial}{\partial y} L - \frac{d}{d x} \frac{\partial}{\partial \dot{y}} L = 0$

\bigskip

y si además $L$ no depende de $x$, se cumple que

\bigskip

$ L - \dot{y} \frac{\partial}{\partial \dot{y}} L = C$

\bigskip

Demostración:

\bigskip

Podemos multiplicar la ecuación de Euler-Lagrange

\bigskip

$\frac{\partial}{\partial y} L - \frac{d}{d x} \frac{\partial}{\partial \dot{y}} L = 0$

\bigskip

por $\dot{y}$:

\bigskip

$\dot{y}\frac{\partial}{\partial y} L - \dot{y}\frac{d}{d x} \frac{\partial}{\partial \dot{y}} L = 0$

\bigskip

pero el término $ \dot{y}\frac{d}{d x} \frac{\partial}{\partial \dot{y}} L$ puede ser substitu\'ido por el equivalente que viene de un producto:

\bigskip

$\frac{d}{dt}(\dot{y} \frac{\partial}{\partial \dot{y}}L) = $
$ \ddot{y}   \frac{\partial}{\partial \dot{y}} L + \dot{y}\frac{d}{dt} \frac{\partial}{\partial \dot{y}}L$

\bigskip

por lo que

\bigskip

$ \dot{y}\frac{d}{dx} \frac{\partial}{\partial \dot{y}}L = \frac{d}{dx}(\dot{y} \frac{\partial}{\partial \dot{y}}L) -  \ddot{y}   \frac{\partial}{\partial \dot{y}} L$

\bigskip

Reemplazando esto en

\bigskip

$\dot{y}\frac{\partial}{\partial y} L - \dot{y}\frac{d}{d x} \frac{\partial}{\partial \dot{y}} L = 0$

\bigskip

tenemos:

\bigskip

$\dot{y}\frac{\partial}{\partial y} L -  \frac{d}{dx}(\dot{y} \frac{\partial}{\partial \dot{y}}L) +  \ddot{y}   \frac{\partial}{\partial \dot{y}} L = 0$

\bigskip

Tengamos ahora en cuenta que para una función $L(y,\dot{y},x)$ la derivada total, que se haya por la regla de la cadena, es

\bigskip

$\frac{d}{dx}L = \dot{y}\frac{\partial}{\partial y}L + \ddot{y}\frac{\partial}{\partial \dot{y}}L + \frac{\partial}{\partial t} L$

\bigskip

pero cuando $L$ no depende explícitamente de $x$, eso se convierte en

\bigskip

$\frac{d}{dx}L = \dot{y}\frac{\partial}{\partial y}L + \ddot{y}\frac{\partial}{\partial \dot{y}}L  $

\bigskip

Reemplazando esto en

\bigskip

$\dot{y}\frac{\partial}{\partial y} L -  \frac{d}{dx}(\dot{y} \frac{\partial}{\partial \dot{y}}L) +  \ddot{y}   \frac{\partial}{\partial \dot{y}} L = 0$

\bigskip

obtenemos

\bigskip

$\frac{d}{dx}L -  \frac{d}{dx}(\dot{y} \frac{\partial}{\partial \dot{y}}L) =0$

\bigskip

lo cual puede integrarse con respecto a $x$:

\bigskip

$ L -   \dot{y} \frac{\partial}{\partial \dot{y}}L =C$

\bigskip
lo cual termina la demostración.

\bigskip

Utilicemos ahora esta identidad para terminar la tarea de demostrar que el camino más corto entre dos puntos del plano con la métrica Euclidiana es una  línea recta.

\bigskip
En este caso, el funcional a minimizar es
\bigskip

$S =  \int^{x_2}_{x_1} L dx =   \int^{x_2}_{x_1} \sqrt{1+(f'(x))^2} dx$

\bigskip

Reemplazando en la identidad de Beltrami

\bigskip

$ L -   \dot{y} \frac{\partial}{\partial \dot{y}}L =C$

\bigskip

obtenemos:

\bigskip

$\sqrt{1+\dot{y}^2} - \dot{y} \frac{\partial}{\partial \dot{y}}\sqrt{1+\dot{y}^2} =C$

\bigskip

$\sqrt{1+\dot{y}^2} - \dot{y} \frac{2\dot{y}}{2\sqrt{1+\dot{y}^2}}   =C$

\bigskip

$\sqrt{1+\dot{y}^2} - \dot{y} \frac{ \dot{y}}{ \sqrt{1+\dot{y}^2}}   =C$

$\frac{1+\dot{y}^2 -\dot{y}^2 }{\sqrt{1+\dot{y}^2}} = C$

$\frac{1  }{\sqrt{1+\dot{y}^2}} = C$

$1 =  C\sqrt{1+\dot{y}^2} $

$\frac{1}{C^2}= 1+  \dot{y}^2$

$ \dot{y}^2 = \frac{1}{C^2} -1$

$ \dot{y} = \sqrt{\frac{1}{C^2} -1}$

\bigskip

Con esto hemos demostrado que el camino más corto en el plano con la métrica euclidiana es una curva de pendiente constante, es decir, es una recta.

\bigskip

Anotamos que hay infinidad de métricas que pueden variar de punto a punto y que por lo tanto pueden hacer que las geodésicas sean familias de curvas que pueden ser rarísimas. Hay una manera cómoda de imaginar infinidades de métricas: Sobre el plano imaginemos una montaña cualesquiera, la cual define geodésicas en el espacio tridimensional pero restringidas a la montaña. Las geodésicas serían los caminos por los cuales uno iría sobre la  montaña de un lugar a otro. Podemos proyectar dichas geodésicas sobre el plano y tendríamos geodésicas sobre el plano. Pero con cual métrica? Con la métrica determinada por la asociación:

$(x,y) \rightarrow (x,y,f(x,y)) \rightarrow $ métrica en el espacio

\bigskip
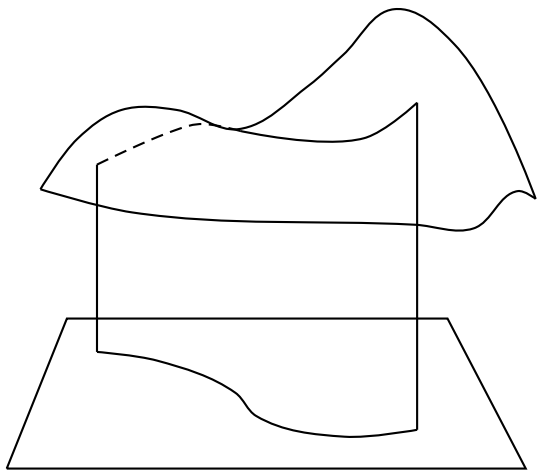
\psset{unit=0.4cm}
\begin{center}
 \begin{pspicture} (0,0)(17,20)

 \psline(1, 0.25)(18.25,0.25)(15.65,5.24)(3,5.24) (1,0.25)
 
 \psline(4, 4.13)(4,10.36)
 \psline(14.64,1.54)(14.64,12.42)
 \pscurve[showpoints=false](4, 4.13) (6,3.84)(8.60,2.77)(9.29,2)(12.23,1.31)(14.64,1.54)

 \pscurve[showpoints=false](2.12,9.54) (5.17,8.77)(14.58,8.36)(16.52,8.24)(17.99,9.48)(18.58,9.22)

 \pscurve[showpoints=false](2.12,9.54) (3.35,11.22)(4.82,12.18)(6.64,12.18)(8.70,11.54)(11.05,13)(12.17,14)(13.88,15.53)(15.88,14.36)(18.58,9.22)
 
 \pscurve[showpoints=false,linestyle=dashed](4,10.36)(6.53,11.48)(7.41,11.71)(8.76,11.48) 
 
 \pscurve[showpoints=false] (8.76,11.48)(12.88,11.24)(14.64,12.42)

 \end{pspicture}
 
 \textit{Con la métrica euclidiana del espacio, cada montaña tiene sus propias geodésicas e induce una métrica específica y geodésicas correspondientes sobre el piso }
 
\end{center}

Concretamente, la distancia entre dos puntos sería:

\bigskip

$d_{Plano}((x_1,y_1), (x_2,y_2)) = $

$d_{Espacio}((x_1,y_1, f((x_1,y_1)), (x_2,y_2,f((x_2,y_2)))$

\bigskip
Otras métricas vienen de deformar el plano para hacer superficies cualesquiera o, lo que es lo mismo,  de cambios de coordenadas. Veamos como luce la métrica sobre  la esfera unidad, utilizando, por supuesto,   coordenadas esféricas.

\bigskip

La esfera es una variedad no trivial, es decir, uno no puede tomar un subconjunto abierto del plano y deformarlo mediante una función continua, biyectiva y diferenciable para obtener una esfera  pues siempre faltará o sobrará algo. Estudiemos este problema desde el punto de vista de  las coordenadas esféricas con   radio $\rho$,   el ángulo de barrido horizontal $\theta$ y el ángulo de barrido vertical $\phi$. Se tiene:

$x = \rho sen \phi cos \theta $

$y = \rho sen \phi sen \theta $

$z = \rho cos \phi $

si $\rho = 1$, entonces el conjunto de desigualdades

$0 < \theta < 2\pi $,

$0 < \phi < \pi $

denota un abierto del plano. Las coordenadas esféricas transforman este abierto en    una esfera (la sola cáscara), a la cual le falta un meridiano,  el meridiano cero. 

Uno no se preocupa en la vida real de este problema, por dos razones: 1) La esfera sin el meridiano  forma un conjunto abierto y denso en la esfera completa con la topología inducida del espacio tridimensional. 2) Las coordenadas esféricas definen naturalmente otro mapa   para cubrir otra parte de la esfera que contenga el meridiano que faltaba, basta tomar

$-\pi < \theta < \pi $,

$0 < \phi < \pi $

y uno hace el ajuste correspondiente instintivamente.

\bigskip

Consideremos dos puntos sobre la esfera, muy cercanos el uno del otro. Queremos hallar la distancia entre ellos sobre la esfera. Como los puntos están a distancia infinitesimal, podemos aproximar la esfera por el plano tangente. Sean los puntos $P$ $Q$, separados por un ángulo de barrido horizontal $d\theta$ y por un ángulo de barrido vertical $d\phi$. Ahora bien, existe un punto $R$ que queda en el mismo meridiano que $P$ y en el mismo paralelo que $Q$. Los tres puntos 'están' todos    sobre el plano tangente a $R$ y determinan un triángulo rectángulo, cuya hipotenusa da la distancia entre los dos puntos.

Ahora bien, el cateto $RP$ está expandido por el ángulo de barrido horizontal $d\theta$, en tanto que el cateto $RQ$ lo está por el ángulo $d\phi$ de barrido vertical. El cateto $RQ$ es la longitud de arco en el círculo que da el meridiano correspondiente, por tanto mide $\rho d\phi$, pero como $\rho = 1$, entonces $RP$ mide $ d\phi$. Similarmente, el cateto $RP$ es la longitud de arco en el círculo que da el paralelo correspondiente, cuyo radio es $\rho sen\phi$. Por tanto el cateto $RP$ mide  $\rho sen\phi d\theta$. De donde concluimos que la hipotenusa mide

\bigskip

$d(P,Q) = \sqrt{(d\phi)^2 +   sen^2 \phi (d\theta)^2}$

\bigskip

y esa es la métrica. Es una métrica sobre el plano y es una métrica sobre la esfera unitaria. Como ya hemos dicho, lo hemos hecho para un manto que casi cubre la esfera, pero puede tomarse un manto complementario sobre el cual da exactamente lo mismo.

La métrica dada está completamente determinada por la forma cuadrática cuya matriz es

$$
\left(
  \begin{array}{cc}
    1 & 0 \\
    0 & sen^2\phi \\
  \end{array}
\right)
$$

Debemos deducir ahora que en la esfera las geodésicas están descritas por los círculos máximos. Pero en vez de estudiar este caso particular, vamos a estudiar el caso general, válido para cualquier variedad. Y  como test para ver que hemos hecho las cosas bien, aplicaremos los resultados a la esfera a ver si recobramos los círculos máximos. 

\section{EL TENSOR METRICO}

Generalizando, una métrica es una función que a cada punto de la variedad le hace corresponder una matriz simétrica de coeficientes reales, la cual da origen a una forma cuadrática. Si además, al diagonalizar dicha matriz se obtienen valores propios positivos, entonces en todas las bases, dicha matriz produce una forma cuadrática positiva y tenemos una métrica con todas las de la ley. Pero si hay algún valor propio negativo, estamos en el caso de una pseudométrica. 

Con la métrica puede medirse distancias sobre el plano tangente al punto dado. Como la métrica mide distancias sobre el plano tangente, dicha métrica puede usarse para medir distancias  a escala infinitesimal sobre la variedad.

 Cómo se hace esto?

 Para entenderlo,  consideremos un camino parametrizado que conduce de $P$ a $Q$ sobre la variedad. Sea la  parametrización del camino dada por

\bigskip

$\gamma : \Re \rightarrow M $

\bigskip

donde $M$ es la variedad.

La manera que vamos a utilizar para medir distancias se basa en la siguiente idea: cuando uno va en un carro, uno puede calcular la distancia recorrida mirando el velocímetro e integrando después. En cálculo, eso se escribe como es bien sabido:

\bigskip

$ds  = vdt$

\bigskip

donde todas son cantidades escalares. Pero sucede que el velocímetro no existe in abstracto sino que va sobre un carro, el cual lleva una cierta forma de viajar. Supongamos pues que la trayectoria del carro se describe por una parametrización:  

\bigskip

$\gamma : \Re \rightarrow M $

$\gamma(t) = ((x^1(t), x^2(t),....,x^n(t)$

\bigskip

Tenemos que  $ds  = vdt$   se convierte en

\bigskip

$ds = ||d\vec{s}||  = ||\vec{v}||dt = \sqrt{||\vec{v}||^2  }dt=  \sqrt{\vec{v}.\vec{v} } dt$

\bigskip

pero

$$\vec{v}.\vec{v} =  \vec v^T \pmatrix{1 & 0 \cr 0 & 1} \vec v$$

En donde vemos la matriz de la métrica euclidiana. Cambiando esa matriz por una general simétrica, $g$,  que puede depender del punto,  llegamos a la expresión:

\bigskip

$ds  = \sqrt{g(\dot{\gamma}, \dot{\gamma})}dt $

$ds^2 = g(\dot{\gamma}, \dot{\gamma})dt^2$

\bigskip

por lo que el espacio total recorrido a lo largo del camino $\gamma$ desde el tiempo $a$ hasta el tiempo $b$ es

\bigskip

$S =   \int^{b}_a \sqrt{g(\dot{\gamma}(t), \dot{\gamma}(t))}dt$

\bigskip

Maticemos un poco la naturaleza de la métrica $g$. La expresión $ g(\dot{\gamma}, \dot{\gamma})$ es una declaración intrínseca, libre de coordenadas, de lo que en una base dada puede escribirse como $\vec v^T [g] \vec v$, donde $\vec v$ es el vector velocidad, que pertenece al tangente, y $[g]$ es la matriz de $g$. Eso significa que la métrica toma dos vectores del tangente y le asocia un número real, de tal forma que la correspondencia genera una forma bilineal, simétrica. Es decir, $g$ es una 2-forma simétrica.

Podemos ahora entender la naturalidad de las siguientes definiciones:  

\bigskip

\addtocounter{ecu}{1}

\textit{Definición (\theecu): Un tensor métrico $g$ es una función que a cada punto de la variedad le hace corresponder una 2-forma simétrica. Por antonomasia, una métrica tiene, en cada punto, todos sus valores propios positivos. Pero es posible que se use la palabra métrica para denotar también pseudo-métricas y eso podrá hacerse cuando no importe la confusión. }

\addtocounter{ecu}{1}

\textit{Definición (\theecu): Una \index{superficie de Riemann} \textbf{superficie de Riemann} es una variedad 2-dimensional con un tensor métrico $M$.}

\bigskip

 Cuando asumimos que $M$ tiene una métrica, $M$ se llama Riemanniana. En relatividad, $M$ debe ser una generalización de la métrica de Minkowski, por lo que al espacio tiempo correspondiente se le llama variedad pseudo-Riemanniana.

\section{EL PROBLEMA VARIACIONAL}

El funcional que mide la distancia recorrida por un carro que lleva la trayectoria $\gamma$ es

$S =   \int^{b}_a \sqrt{g(\dot{\gamma}(t), \dot{\gamma}(t))}dt$

Este funcional no cambia de valor si cambiamos la forma de viajar: la distancia recorrida es invariante ante la reparametrización. Estaremos atentos para ver si en algún momento nos convenga  elegir una parametrización que facilite los cálculos. 

\bigskip

\addtocounter{ecu}{1}
 
\textit{ Definición (\theecu): una \index{geodésica} \textbf{geodésica} es  el camino trazado por cualquier trayectoria $\gamma$ que anule, a primer orden, la variación del funcional }

$S =   \int^{b}_a \sqrt{g(\dot{\gamma}(t), \dot{\gamma}(t))}dt$

\textit{Es decir, para hallar una geodésica, lo que tenemos que hacer es estudiar la variación a primer orden de $S$ e igualarlo a cero para encontrar un sistema de ecuaciones cuya solución sea una  trayectoria que de la geodésica.}

\addtocounter{ecu}{1}
 
\textit{ Definición (\theecu): Decimos que un \index{espacio curvo} \textbf{espacio es curvo} cuando sus geodésicas no son líneas rectas.} 

\bigskip

Como ya hemos visto, el plano puede ser muy curvo, pues se le puede adjudicar toda suerte de métricas muy distintas de la euclidiana.

\bigskip

El problema general de aniquilar a primer orden la variación de un funcional ya lo resolvimos y dan las ecuaciones de Euler-Lagrange:

\bigskip

$\frac{\partial}{\partial x^i} L - \frac{d}{d t} \frac{\partial}{\partial \dot{x^i}} L = 0$

\bigskip

 En nuestro caso $L = \sqrt{g(\dot{\gamma}(t), \dot{\gamma}(t))}$. 
 
 \bigskip

Tenemos un problema con las pseudométricas, para así poder aplicar los resultados a la relatividad (compare el artículo sobre la geodésica (geodesic) de Wikipedia), pues $ g(\dot{\gamma}(t), \dot{\gamma}(t)) $ puede ser positivo o negativo. Tomamos entonces el valor absoluto de $g$. En realidad, el signo positivo o negativo que pone el valor absoluto se esfuma en el cálculo a seguir, lo cual podemos probarlo si tomamos el signo menos:

\bigskip

$\frac{\partial}{\partial x^i} \sqrt{-g(\dot{\gamma}(t), \dot{\gamma}(t))} - \frac{d}{d t} \frac{\partial}{\partial \dot{x^i}} \sqrt{-g(\dot{\gamma}(t), \dot{\gamma}(t))} = 0$

\bigskip

Podemos especificar $ g(\dot{\gamma}(t), \dot{\gamma}(t))$ en el sistema de coordenadas adoptado por la parametrización. En ese caso, $ g(\dot{\gamma}(t), \dot{\gamma}(t))$ se evalúa en cada punto por medio de la matriz de $g$ en la base originada en ese punto por el sistema de coordenadas. Tenemos:

\bigskip

$ g(\dot{\gamma}(t), \dot{\gamma}(t)) \rightarrow  g_{\mu \nu } \dot{x}^\mu \dot{x}^\nu$

\bigskip

por lo que las ecuaciones de Euler Lagrange se releen como

\bigskip

$\frac{\partial}{\partial x^i} \sqrt{-g_{\mu \nu } \dot{x}^\mu \dot{x}^\nu} - \frac{d}{d t} \frac{\partial}{\partial \dot{x^i}} \sqrt{-g_{\mu \nu } \dot{x}^\mu \dot{x}^\nu} = 0$

\bigskip

Calculando las derivadas obtenemos

\bigskip

$ \frac{ - g_{\mu \nu,\iota   }\dot{x}^\mu  \dot{x}^\nu}{2\sqrt{-g_{\mu \nu } \dot{x}^\mu \dot{x}^\nu}} - \frac{d}{d t} ( \frac{ - g_{\mu \nu } \frac{\partial \dot{x}^\mu}{ \partial\dot{x}^i }  \dot{x}^\nu - g_{\mu \nu } \dot{x}^\mu \frac{\partial \dot{x}^\nu}{ \partial\dot{x}^i }  }{2\sqrt{-g_{\mu \nu } \dot{x}^\mu \dot{x}^\nu}}) = 0$

\bigskip

donde $ g_{\mu \nu,\iota   }$ significa la derivada de $g_{\mu \nu}$ con respecto a $x^\iota$.

Podemos multiplicar por menos uno  en todo lado para obtener exactamente lo mismo que hubiésemos obtenido si el valor absoluto hubiese sido positivo, como en una métrica:

\bigskip

$ \frac{   g_{\mu \nu,\iota   }\dot{x}^\mu  \dot{x}^\nu}{2\sqrt{-g_{\mu \nu } \dot{x}^\mu \dot{x}^\nu}} - \frac{d}{d t} ( \frac{   g_{\mu \nu } \frac{\partial \dot{x}^\mu}{ \partial\dot{x}^\iota }  \dot{x}^\nu + g_{\mu \nu } \dot{x}^\mu \frac{\partial \dot{x}^\nu}{ \partial\dot{x}^\iota }  }{2\sqrt{-g_{\mu \nu } \dot{x}^\mu \dot{x}^\nu}}) = 0$

\bigskip

Pero ahora bien,  $\frac{\partial \dot{x}^\mu}{ \partial\dot{x}^\iota} = \delta^\mu_\iota$ y $\frac{\partial \dot{x}^\nu}{ \partial\dot{x}^\iota} = \delta^\nu_\iota$ por lo que la ecuación anterior se convierte en

\bigskip

$ \frac{   g_{\mu \nu,\iota   }\dot{x}^\mu  \dot{x}^\nu}{2\sqrt{-g_{\mu \nu } \dot{x}^\mu \dot{x}^\nu}} - \frac{d}{d t} ( \frac{   g_{\mu \nu }  \delta^\mu_i    \dot{x}^\nu + g_{\mu \nu } \dot{x}^\mu \delta^\nu_i }{2\sqrt{-g_{\mu \nu } \dot{x}^\mu \dot{x}^\nu}}) = 0$

\bigskip

$ \frac{   g_{\mu \nu,\iota   }\dot{x}^\mu  \dot{x}^\nu}{\sqrt{-g_{\mu \nu } \dot{x}^\mu \dot{x}^\nu}} - \frac{d}{d t} ( \frac{   g_{\iota\nu }       \dot{x}^\nu + g_{\mu \iota } \dot{x}^\mu }{\sqrt{-g_{\mu \nu } \dot{x}^\mu \dot{x}^\nu}}) = 0$

\bigskip

donde hemos aplicado que $g_{\mu \nu }  \delta^\mu_i$ implica suma sobre $\mu$, la cual se reduce tan sólo al término  en el cual $\mu = \iota$. Y algo semejante con el otro término. Ahora podemos calcular la derivada con respecto al tiempo, en el segundo término, el cual es  un cociente:

\bigskip

$ \frac{   g_{\mu \nu,\iota   }\dot{x}^\mu  \dot{x}^\nu}{\sqrt{-g_{\mu \nu } \dot{x}^\mu \dot{x}^\nu}} $
$-\frac{\sqrt{-g_{\mu \nu } \dot{x}^\mu \dot{x}^\nu} \frac{d}{dt} ( g_{\iota\nu }       \dot{x}^\nu + g_{\mu \iota } \dot{x}^\mu) - ( g_{\iota\nu }       \dot{x}^\nu + g_{\mu \iota } \dot{x}^\mu) \frac{d}{dt} \sqrt{-g_{\mu \nu } \dot{x}^\mu \dot{x}^\nu  }}{-g_{\mu \nu } \dot{x}^\mu \dot{x}^\nu}=0 $

\bigskip

Al derivar con respecto al tiempo, hay que usar la regla de la cadena sobre $g$ pues se cambia de punto con el tiempo y la métrica cambia de punto a punto: 

\bigskip

$ \frac{   g_{\mu \nu,\iota   }\dot{x}^\mu  \dot{x}^\nu}{\sqrt{-g_{\mu \nu } \dot{x}^\mu \dot{x}^\nu}} $
$=\frac{\sqrt{-g_{\mu \nu } \dot{x}^\mu \dot{x}^\nu}   
( g_{\iota\nu,\mu } \dot{x}^\nu \dot{x}^\mu + g_{\iota\nu }       \ddot{x}^\nu + g_{\mu \iota, \nu } \dot{x}^\mu\dot{x}^\nu +  
 g_{\mu \iota } \ddot{x}^\mu) 
- ( g_{\iota\nu }       \dot{x}^\nu + g_{\mu \iota } \dot{x}^\mu)
 \frac{\frac{d}{dt}(-g_{\mu \nu } \dot{x}^\mu \dot{x}^\nu) }{2 \sqrt{-g_{\mu \nu } \dot{x}^\mu \dot{x}^\nu  }}                   }
 {-g_{\mu \nu } \dot{x}^\mu \dot{x}^\nu}$

\bigskip

$ \frac{   g_{\mu \nu,\iota   }\dot{x}^\mu  \dot{x}^\nu}{\sqrt{-g_{\mu \nu } \dot{x}^\mu \dot{x}^\nu}} $
$=\frac{ -2g_{\mu \nu } \dot{x}^\mu \dot{x}^\nu   
( g_{\iota\nu,\mu } \dot{x}^\nu \dot{x}^\mu + g_{\iota\nu }       \ddot{x}^\nu + g_{\mu \iota, \nu } \dot{x}^\mu\dot{x}^\nu +  
g_{\mu \iota } \ddot{x}^\mu) 
- ( g_{\iota\nu }       \dot{x}^\nu + g_{\mu \iota } \dot{x}^\mu)
 \frac{d}{dt}(-g_{\mu \nu } \dot{x}^\mu \dot{x}^\nu)                  }
 {-g_{\mu \nu } \dot{x}^\mu \dot{x}^\nu
 (2 \sqrt{-g_{\mu \nu } \dot{x}^\mu \dot{x}^\nu  })}$

\bigskip

$     g_{\mu \nu,\iota   }\dot{x}^\mu  \dot{x}^\nu   $
$=\frac{ -2g_{\mu \nu } \dot{x}^\mu \dot{x}^\nu   
( g_{\iota\nu,\mu } \dot{x}^\nu \dot{x}^\mu + g_{\iota\nu }       \ddot{x}^\nu + g_{\mu \iota, \nu } \dot{x}^\mu\dot{x}^\nu +  
 g_{\mu \iota } \ddot{x}^\mu) 
- ( g_{\iota\nu }       \dot{x}^\nu + g_{\mu \iota } \dot{x}^\mu)
 \frac{d}{dt}(-g_{\mu \nu } \dot{x}^\mu \dot{x}^\nu)                  }
 {-2g_{\mu \nu } \dot{x}^\mu \dot{x}^\nu}$

\bigskip

$     g_{\mu \nu,\iota   }\dot{x}^\mu  \dot{x}^\nu   $

$=  ( g_{\iota\nu,\mu } \dot{x}^\nu \dot{x}^\mu + g_{\iota\nu }       \ddot{x}^\nu +g_{\mu \iota, \nu } \dot{x}^\mu\dot{x}^\nu +  
 g_{\mu \iota } \ddot{x}^\mu) 
+ \frac{1}{2} ( g_{\iota\nu }       \dot{x}^\nu + g_{\mu \iota } \dot{x}^\mu)
 \frac{\frac{d}{dt}(-g_{\mu \nu } \dot{x}^\mu \dot{x}^\nu)}                  
 {g_{\mu \nu } \dot{x}^\mu \dot{x}^\nu}$

 \bigskip
 
Ahora recordamos que con la métrica se suben y se bajan los índices:

 \bigskip
 
$     g_{\mu \nu,\iota   }\dot{x}^\mu  \dot{x}^\nu   $
$=  ( g_{\iota\nu,\mu } \dot{x}^\nu \dot{x}^\mu + g_{\iota\nu }       \ddot{x}^\nu +g_{\mu \iota, \nu } \dot{x}^\mu\dot{x}^\nu +  
 g_{\mu \iota } \ddot{x}^\mu) 
+ \frac{1}{2} (         \dot{x}_\iota +   \dot{x}_\iota)
 \frac{\frac{d}{dt}(-g_{\mu \nu } \dot{x}^\mu \dot{x}^\nu)}                  
 {g_{\mu \nu } \dot{x}^\mu \dot{x}^\nu}$
 
Como la métrica está dada por un tensor simétrico, se tiene que $g_{\iota \mu} = g_{\mu \iota }   $ y cuando $\mu$ es un índice mudo puede reemplazarse por $\nu$. Tenemos por tanto que 
 
 \bigskip

 $g_{\mu \iota } \ddot{x}^\mu = g_{\iota\mu  } \ddot{x}^\mu =g_{\iota\nu  } \ddot{x}^\nu  $. 
 
  \bigskip
 
 Reemplazando obtenemos:

 \bigskip
 
$     g_{\mu \nu,\iota   }\dot{x}^\mu  \dot{x}^\nu   $
$=  ( g_{\iota\nu,\mu } \dot{x}^\nu \dot{x}^\mu + 2g_{\iota\nu }       \ddot{x}^\nu +g_{\mu \iota, \nu } \dot{x}^\mu\dot{x}^\nu  ) 
+ \frac{1}{2} (         \dot{x}_\iota +   \dot{x}_\iota)
 \frac{\frac{d}{dt}(-g_{\mu \nu } \dot{x}^\mu \dot{x}^\nu)}                  
 {g_{\mu \nu } \dot{x}^\mu \dot{x}^\nu}$

\bigskip

$     g_{\mu \nu,\iota   }\dot{x}^\mu  \dot{x}^\nu   $
$-   g_{\iota\nu,\mu } \dot{x}^\nu \dot{x}^\mu -2g_{\iota\nu }       \ddot{x}^\nu -g_{\mu \iota, \nu } \dot{x}^\mu\dot{x}^\nu   
=        -  \dot{x}_\iota  
 \frac{\frac{d}{dt}(g_{\mu \nu } \dot{x}^\mu \dot{x}^\nu)}                  
 {g_{\mu \nu } \dot{x}^\mu \dot{x}^\nu}$

\bigskip

Ordenamos, factorizamos  y bajamos algunos índices:

\bigskip

$   -(   g_{\iota\nu,\mu }+ g_{\mu \iota, \nu }- g_{\mu \nu,\iota   }     ) \dot{x}^\mu  \dot{x}^\nu   $
$ -2       \ddot{x}_\iota  
=         - \dot{x}_\iota  
 \frac{\frac{d}{dt}(  \dot{x}_\nu \dot{x}^\nu)}                  
 {  \dot{x}_\nu \dot{x}^\nu}$

\bigskip

Multiplicando por menos uno, recordando la simetría de $g$, y  recuperando en el lado izquierdo la derivada de un logaritmo obtenemos:

\bigskip

$   (   g_{\iota\nu,\mu }+ g_{\mu \iota, \nu }- g_{\mu \nu,\iota   }     ) \dot{x}^\mu  \dot{x}^\nu   $
$ +2       \ddot{x}_\iota  
=          \dot{x}_\iota  
 \frac{d}{dt}ln| \dot{x}_\nu \dot{x}^\nu| $

\bigskip

$  \frac{1}{2} (   g_{\iota\nu,\mu }+ g_{\mu \iota, \nu }- g_{\nu \mu,\iota   }     ) \dot{x}^\mu  \dot{x}^\nu   $
$      +   \ddot{x}_\iota  
=     \frac{1}{2}       \dot{x}_\iota  
 \frac{d}{dt}ln| \dot{x}_\nu \dot{x}^\nu| $
 
  \bigskip
  
  Observando la ecuación anterior, uno nota que en el lado derecho aparece el logaritmo del valor absoluto de  la norma cuadrado del vector velocidad: $|\dot{x}_\nu \dot{x}^\nu| =|g_{\mu \nu } \dot{x}^\mu \dot{x}^\nu|$.  Ahora bien, si el vector velocidad tuviese norma constante, se aniquilaría el lado derecho de la ecuación de las geodésicas y quedaría una ecuación mucho más simple. ¿Podrá hacerse eso?
  
La respuesta es afirmativa: debido a que la longitud de arco es invariante ante la parametrización, lo cual se deduce de la invariancia de la integral $S$ ante el cambio de variable,  podemos elegir una parametrización que vaya a velocidad constante. Naturalmente que dicha parametrización es precisamente la definida por la longitud de arco y recorre una unidad de espacio por unidad de tiempo. Quizá convenga explicitar esto. La longitud de arco se define como:

$s(t) =    \int^{t}_a \sqrt{g(\dot{\gamma}(u), \dot{\gamma}(u))}du$

o, derivando y utilizando el teorema fundamental del cálculo y el concepto de norma:

$\dot s(t) = ||\dot{\gamma}(t)|| = \sqrt{g(\dot{\gamma}(u), \dot{\gamma}(u))} $

pero

$||\dot{\gamma}(t)|| = || \frac{d\gamma (t)}{dt} || = || \frac{d\gamma (t)}{ds} ||  \frac{d s(t)}{dt} $

pues se toma la raíz positiva para definir $s(t)$. Simplificando $ || \frac{d\gamma (t)}{dt} ||$ con $ \frac{d s(t)}{dt}$, obtenemos:

$|| \frac{d\gamma (t)}{ds} || = 1$, 

lo que dice que cualquier trayectoria reparametrizada en longitud de arco tiene vector velocidad de magnitud constante uno. 

\bigskip

Retomando la ecuación de las geodésicas, obtenemos:

$$  \frac{1}{2} (   g_{\iota\nu,\mu }+ g_{\mu \iota, \nu }- g_{\mu \nu,\iota   }     ) \dot{x}^\mu  \dot{x}^\nu    
       +   \ddot{x}_\iota  
=   0 $$

o bien

$$  \frac{1}{2} (  \partial_\mu g_{\iota\nu }+ \partial_\nu g_{\mu \iota }- \partial_\iota g_{\mu \nu  }     ) \dot{x}^\mu  \dot{x}^\nu    
       +   \ddot{x}_\iota  
=   0 $$

donde $g_{\alpha\beta,\gamma } = \partial g_{\alpha\beta } /\partial x^\gamma = \partial_\gamma g_{\alpha\beta}$. 

Estas condiciones también suelen reescribirse teniendo en cuenta que $ \ddot{x}_\iota  = g_{\iota \gamma} \ddot{x}^\gamma$, que al reemplazar produce: 

$$  \frac{1}{2} (  \partial_\mu g_{\iota\nu }+ \partial_\nu g_{\mu \iota }- \partial_\iota g_{\mu \nu  }     ) \dot{x}^\mu  \dot{x}^\nu    
       +    g_{\iota \gamma} \ddot{x}^\gamma  
=   0 $$

y multiplicando en todo lado por la inversa de  $g_{\iota \gamma}$ que se nota $ g^{ \gamma\iota}$ obtenemos:

$$  \frac{1}{2}g^{ \gamma \iota} (  \partial_\mu g_{\iota\nu }+ \partial_\nu g_{\mu \iota }- \partial_\iota g_{\mu \nu  }     ) \dot{x}^\mu  \dot{x}^\nu    
       +     \ddot{x}^\gamma  
=   0 $$

Resumamos todo el trabajo anterior en el siguiente teorema:
 
\bigskip

\addtocounter{ecu}{1}
\textit{
\textit{Teorema (\theecu):  Sea $M$ una variedad con una métrica   dada $g$. Una \index{geodésica} \textbf{geodésica} sobre $M$ es una curva o trayectoria parametrizada por longitud de arco  que da un extremo del funcional distancia recorrida y se  caracteriza por que cumple con el sistema de ecuaciones:}}
 
 $$  \frac{1}{2}g^{ \gamma \iota} (  \partial_\mu g_{\iota\nu }+ \partial_\nu g_{\mu \iota }- \partial_\iota g_{\mu \nu  }     ) \dot{x}^\mu  \dot{x}^\nu    
       +     \ddot{x}^\gamma  
=   0 $$

\bigskip
  
¿ Y en donde quedaron los círculos máximos como geodésicas de la esfera?

La salida rápida a esa pregunta es como sigue: siempre existe un sistema coordenado que convierte un círculo máximo en el ecuador de la esfera unidad. En coordenadas esféricas, el ecuador tiene como ecuación $\phi = \pi / 2$, cuyo coseno es cero y cuyo seno es uno. Así que verifiquemos que estos valores satisfacen la ecuación geodésica.

La métrica sobre la esfera está dada por la matriz

$$g = \pmatrix{1 & 0 \cr 0 & sen^2\phi}$$

la derivada de dicha matriz con respecto a $\theta $ es cero y con respecto a $\phi$ es

$$\partial_\phi g = \pmatrix{0 & 0 \cr 0 & 2 sen \phi cos \phi}$$

la cual es idénticamente cero sobre el ecuador, es decir, cuando $\phi = \pi/2$ . Al reemplazar todos esos  ceros en el sistema de ecuaciones de la geodésica, obtenemos simplemente:

$\ddot{x}^\gamma  = 0$

lo cual da $x^\gamma = as+b$ y que para todos los efectos se confunde con $x^\gamma = s$. Como $\phi = \pi/2$, la ecuación $x^\gamma = s$ corresponde a $\theta = s$. Obtenemos así la parametrización del ecuador que es entonces una geodésica.

\bigskip

Es supremamente interesante que este sistema de ecuaciones de las geodésicas se preste a una interpretación geométrica que siendo más bien artificiosa ha resultado de inmensa utilidad, tanto   que ha  llegado a ser fundamental no sólo en el entendimiento que tenemos de la relatividad general y sus actualizaciones sino también de las teoría gauge. Veámosla.

\section{LOS SIMBOLOS DE CHRISTOFFEL}

Podemos entender una superficie $S$  como una termodeformación de  un retazo    $U$ tomado de un  plano hecho de material plástico. Modelamos la termodeformación mediante una función:

\addtocounter{ecu}{1}

$$\vec x : U \rightarrow \Re \sp 3\eqno{(\theecu)}$$

$U$ tiene como coordenadas $u=u\sp 1$ y $v= v\sp 2$. Si $v$ se mantiene fijo y se varía $u$, se obtendrá una curva parametrizada por $u$. Por tanto, a $u$  se le puede entender como 'tiempo horizontal' y por tanto la derivada

\addtocounter{ecu}{1}

$$\vec x\sb u = \partial \vec x/\partial u \eqno{(\theecu)}$$

es  un vector velocidad de la curva, el cual es tangente a la superficie $S$. Algo similar ocurre con $v$, quien nos daría el 'tiempo vertical'.

Destacamos entonces el doble rol de la función $\vec x$: conlleva una deformación y al mismo tiempo un sistema de coordenadas.

Si la termodeformación transforma 2-plps en 2-plps sin reducirlos a una punto o línea, entonces los dos vectores velocidad expandirán un plano, llamado el espacio tangente, del cual ellos formarán una base. Ahora bien, si tomamos las segundas derivadas y requerimos su continuidad, entonces ellas serán iguales:

\addtocounter{ecu}{1}

$$ \vec x\sb{\alpha \beta}= \partial \sp 2 \vec x /\partial u\sp \beta \partial u\sp \alpha =\vec x\sb{\beta \alpha}\eqno{(\theecu)}$$

Estas segundas derivadas están relacionadas con la aceleración, que en general no es tangente a la superficie. Para verlo, consideremos el movimiento sobre una espiral sumergida sobre una escalera de caracol: al moverse sobre la espiral siempre se está cambiando de plano tangente, por lo que debe haber una componente de la aceleración no contenida en dicho plano. Descompongamos la aceleración en su parte tangencial a la superficie y en su parte normal a ella. La parte tangencial será una combinación lineal de los vectores de la base del plano tangente, de la forma $\vec x\sb \gamma \Gamma\sp{\gamma}$, donde $\vec x\sb \gamma$ indica un vector tangente a la superficie, dado por la velocidad de una trayectoria generada por el sistema de coordenadas.

\bigskip

Por otro lado, la componente normal será un alargamiento del vector normal unitario,

$\vec N = \frac{\vec x \sb u \times \vec x\sb v}{||\vec x \sb u \times \vec x\sb v ||}$.

El escalar que da la componente normal se calcula por proyecciones, es decir, usando el producto interior, $ <\vec x \sb {\alpha \beta},\vec N>$. Tenemos:

\addtocounter{ecu}{1}

$$\vec x\sb{\alpha \beta} = \partial \sb  \beta \partial \sb \alpha\vec x = \vec x\sb \gamma \Gamma\sp{\gamma} \sb{\beta \alpha} + <\vec x \sb {\alpha \beta},\vec N>\vec N \eqno{(\theecu)}$$

Definiendo

\addtocounter{ecu}{1}

$$b\sb{\alpha \beta} = <\vec x \sb {\alpha \beta},\vec N>\eqno{(\theecu)}$$

obtenemos

\addtocounter{ecu}{1}

$$\vec x\sb{\alpha \beta} =  \vec x\sb \gamma \Gamma\sp{\gamma} \sb{\beta \alpha} + b\sb{\alpha \beta}\vec N \eqno{(\theecu)}$$

Ahora determinaremos los coeficientes $\Gamma\sp{\gamma} \sb{\beta \alpha}=\Gamma\sp{\gamma} \sb  {\alpha \beta }$, cuya simetría en los subíndices inferiores está garantizada por la igualdad de las derivas cruzadas. Utilizando la perpendicularidad entre un vector tangente y el vector normal:

$<\vec x\sb{\alpha \beta} , \vec x \sb \mu >=<\vec x\sb \gamma \Gamma\sp{\gamma} \sb{\beta \alpha} + b\sb{\alpha \beta}\vec N ,\vec x \sb \mu> = <\vec x\sb \gamma \Gamma\sp{\gamma} \sb{\beta \alpha} ,\vec x \sb \mu> = <\vec x\sb \gamma  ,\vec x \sb \mu>\Gamma\sp{\gamma} \sb{\beta \alpha}$

Usamos la notación :

\addtocounter{ecu}{1}

$$<\vec x\sb \gamma  ,\vec x \sb \mu> = g\sb{\gamma \mu}\eqno{(\theecu)}$$

teniendo en mente que los vectores velocidad que definieron el plano tangente forman una base, la cual no tiene porque ser ortogonal. Todos esos valores forman una matriz que representan una 2-forma a la cual se le llama el tensor métrico del plano tangente a la superficie en el punto dado. También notamos:

$<\vec x\sb \gamma \Gamma\sp{\gamma} \sb{\beta \alpha} ,\vec x \sb \mu> = <\vec x\sb \gamma  ,\vec x \sb \mu>\Gamma\sp{\gamma} \sb{\beta \alpha}= g\sb{\gamma \mu}\Gamma\sp{\gamma} \sb{\beta \alpha}= \Gamma  \sb{\beta \alpha , \mu} $

Por tanto tenemos:

$<\vec x\sb{\alpha \beta} , \vec x \sb \mu >=\Gamma  \sb{\beta \alpha , \mu} $

Estos coeficientes tienen un nombre oficial: símbolos o \index{coeficientes de Christoffel} \textbf{coeficientes de Christoffel} (de segundo tipo).
Puesto que las derivadas mixtas son iguales, se tiene que estos coeficientes son simétricos en los dos primeros sub-índices:

$\Gamma  \sb{\beta \alpha , \mu} = \Gamma  \sb{\alpha \beta  , \mu} $

Por otro lado y teniendo en cuenta que $\vec x\sb{\alpha }$ representa una derivada:

$\partial \sb \beta g\sb{\alpha \mu} = \partial \sb \beta  <\vec x\sb{\alpha },\vec x\sb \mu> =  <\partial \sb \beta\vec x\sb{\alpha },\vec x\sb \mu> +   <\vec x\sb{\alpha },\partial \sb \beta\vec x\sb \mu> $

$= <\vec x\sb{\alpha \beta},\vec x\sb{\mu}>  + <\vec x\sb{\alpha },\vec x\sb{\mu \beta}> = $
$\Gamma \sb{\beta \alpha , \mu} + \Gamma \sb{\beta \mu, \alpha}$

Esta ecuación es la primera del siguiente sistema que resulta de rotar los   subíndices  del término del lado izquierdo:  

$\partial \sb \beta g\sb{\alpha \mu} =\Gamma \sb{\beta \alpha , \mu} + \Gamma \sb{\beta \mu, \alpha}$

$\partial \sb \mu g\sb{\beta \alpha} =\Gamma \sb{\mu \beta , \alpha} + \Gamma \sb{\mu \alpha, \beta  }$

$\partial \sb \alpha g\sb{\mu \beta} =\Gamma \sb{\alpha \mu , \beta} + \Gamma \sb{\alpha \beta, \mu }$

Combinemos teniendo en cuenta la simetría de los dos primeros subíndices:

$\partial \sb \beta g\sb{\alpha \mu} + \partial \sb \mu g\sb{\beta \alpha } - \partial \sb \alpha g\sb{ \mu \beta} = 2  \Gamma \sb{\mu \beta , \alpha} = 2 g\sb{\gamma \alpha} \Gamma \sp \gamma \sb{\mu \beta }$

despejamos teniendo en cuenta que la matriz inversa de $g\sb{\gamma \alpha} $ es $g\sp{ \alpha\gamma}$

\addtocounter{ecu}{1}

$$\Gamma \sp \gamma \sb{\mu \beta }= (1/2)g\sp{\alpha \gamma}(\partial \sb \beta g\sb{\alpha \mu} + \partial \sb \mu g\sb{\beta \alpha } - \partial \sb \alpha g\sb{ \mu \beta})\eqno{(\theecu)}$$

De esa forma podemos conocer los coeficientes de Christoffel a partir de la métrica. Recordemos que la métrica opera sobre el espacio tangente, por lo tanto da las relaciones de ortogonalidad de  los vectores velocidad dados por el sistema de coordenadas.

\section{LA DERIVADA COVARIANTE}

A partir de la derivada de un campo vectorial definido sobre una superficie, quitaremos la componente normal de dicha derivada para obtener la derivada covariante.

Un campo vectorial es una asignación que a cada punto de un subconjunto del espacio le hace corresponder un vector. Por ejemplo, el campo gravitatorio es un campo vectorial que a cada punto del espacio le hace corresponder el vector fuerza que corresponde al peso de una masa de un 1kg. Consideremos una curva $C$ parametrizada por $t$, dentro de una superficie $S$. Consideremos que hay un campo vectorial $\vec X$ definido a lo largo de la curva y tangente a ella.  Atención:  El campo vectorial es tangente a  la curva,  pero su derivada puede tener cualquier dirección. Por ejemplo, uno ata un piedra de un hilo y la hace girar: el campo de velocidades es tangente a la curva, pero la derivada de la velocidad da la aceleración, que bien puede ir hacia la mano que la hace girar.

\addtocounter{ecu}{1}

\textit{ Definición (\theecu): Sea $\vec X$ un campo vectorial definido sobre una curva, la cual está contenida en una superficie con normal $\vec N$.  La \index{derivada covariante} \textbf{derivada covariante} es  la derivada ordinaria menos su componente normal:}

$$\nabla \vec X/dt= d\vec X/dt - <d\vec X/dt,\vec N>\vec N$$

\addtocounter{ecu}{1}

\textit{Definición (\theecu) Una \index{geodésica} \textbf{geodésica} es una curva $\vec x =\vec x(s)$ parametrizada por longitud de arco, tal que su vector tangente unitario $\vec T$ no tiene derivada covariante:}

$$\nabla \vec T/ds =0$$

Observemos que geodésica es una definición local y no global. Por otro lado, es muy bueno tener presente que una cosa es ser normal a la curva y otra es ser normal a la superficie que contiene la curva: la derivada covariante resta la componente normal a la superficie, pero no a la curva. Esto nos ayuda a ver cómo queremos cuadrar las cosas para que nuestra definición local quede ligada a una  minimización:

El vector tangente unitario es paralelo a la velocidad, de tal manera que su derivada indica aceleración, es decir una fuerza ejercida desde el exterior sobre la partícula. Esa fuerza tiene poder pues   hace curvar la trayectoria. Pero su poder es vectorial, y hay una parte que puede ser tangencial a la superficie. Si imaginamos la curva  como una cuerda elástica que se sostiene sobre la superficie por efecto de la tensión, esa componente tangencial tiende a hacer que la curva se reacomode para experimentar menos tensión.  La tensión llega a ser mínima cuando  no haya ninguna fuerza tangencial creada por el campo externo.

Por tanto, una derivada covariante diferente de cero indica una fuerza que no ha sido realizada y que busca reacomodar la curva sobre la superficie hasta que ya no experimente ninguna fuerza transversal. Pero hay que ver de qué forma sale la fuerza a partir de la métrica. Para verlo, lo mejor es hallar  la ecuación diferencial que deben cumplir las coordenadas de una curva  parametrizada con longitud de arco $s$ sobre una superficie para que sea geodésica.

La superficie  tiene parámetros $u\sp \beta$, de tal forma que $\partial \vec x/\partial u\sp \beta$ representa la velocidad con respecto a la parametrización generada por el parámetro $\beta$ de la superficie.

Para una curva $\vec x = \vec x(u(s))$ sobre la superficie

$\vec T = d\vec x/ds = (\partial \vec x/\partial u\sp \beta)(d u\sp \beta/ds) = \vec x\sb \beta (d u\sp \beta/ds)$

$ d\vec T/ds = d\sp 2 \vec x/ds\sp 2 = \vec x\sb {\beta \alpha} (d u\sp \alpha /ds)(d u\sp \beta/ds) + \vec x\sb \beta (d \sp 2  u\sp \beta/ds \sp 2 )$

$= (\vec x\sb \gamma \Gamma\sp{\gamma} \sb{\beta \alpha} + b\sb{\alpha \beta}\vec N ) (d u\sp \alpha /ds)(d u\sp \beta/ds) + \vec x\sb \gamma (d \sp 2  u\sp \gamma/ds \sp 2 )$

quitando la componente normal obtenemos la derivada covariante:

$\nabla \vec T/ds= (\vec x\sb \gamma \Gamma\sp{\gamma} \sb{\beta \alpha}  ) (d u\sp \alpha /ds)(d u\sp \beta/ds) + \vec x\sb \gamma (d\sp 2  u\sp \gamma/ds \sp 2 )$

$\nabla \vec T/ds= \vec x\sb \gamma[d \sp 2  u\sp \gamma/ds \sp 2 + \Gamma\sp{\gamma} \sb{\beta \alpha}(d u\sp \alpha /ds)(d u\sp \beta/ds)]$

Vemos que la derivada covariante puede escribirse como una combinación lineal de las velocidades de las coordenadas, todas tangentes a la superficie. Por ello, lo que tenemos es realmente tangente a la superficie.

$ $
Hemos demostrado el siguiente

\bigskip

\addtocounter{ecu}{1}

\textit{Teorema (\theecu): Una curva $\vec x = \vec x(u(s))$ parametrizada por longitud de arco es \index{geodésica} \textbf{geodésica} en una superficie cuando su derivada covariante es cero, es decir, cuando:}

$d \sp 2  u\sp \gamma/ds \sp 2 + \Gamma\sp{\gamma} \sb{\beta \alpha}(d u\sp \alpha /ds)(d u\sp \beta/ds)=0$

\textit{Este sistema de ecuaciones de segundo orden  puede reescribirse como un sistema de primer orden:}

$d u\sp \gamma/ds=X\sp \gamma$

$dX\sp \gamma/ds= - \Gamma\sp{\gamma} \sb{\beta \alpha} X\sp \alpha X\sp \beta$

\bigskip

Esta sistema de  ecuaciones define un problema que tiene solución única para condiciones iniciales dadas, un punto inicial y el vector tangente unitario inicial, al menos durante un intervalo abierto. Eso significa que a partir de un punto y en una dirección, siempre sale una geodésica. Dicha geodésica minimiza absolutamente la distancia dada por el producto interior en una pequeña vecindad. A gran escala, no hay garantía de nada, pero la curva resultante es, al menos, un valor crítico para el funcional de longitud de arco.

Por ejemplo: se puede ir de un punto a otro en una esfera por el círculo máximo que los une, pero yendo por un lado se obtendrá un mínimo absoluto, mientras que por el otro lado no hay  un mínimo local, sino el equivalente a un punto de inflexión. Eso se debe a que uno puede perturbar dicha trayectoria para aumentar el recorrido y también para disminuirlo.

\section{TRANSPORTE PARALELO}

Uno entiende instintivamente lo que significa transportar paralelamente un vector en un plano, como quien dice, transportar un lápiz paralelamente a sí   mismo mientras que permanece sobre una mesa. Pero  es muchísimo más difícil imaginarse lo que significa transportar un lápiz sobre una naranja, de tal forma que sus propiedades de tangencia y orientación permanezcan invariantes. La solución a este problema con un propósito geodésico fue propuesta por Levi-Civita. Veámosla.

En $\Re \sp 2$ el deslizamiento o \index{transporte paralelo} \textbf{transporte paralelo} de un vector a lo largo de una línea recta es exactamente lo que uno se imagina: el vector se desliza a lo largo de la línea sin cambiar su magnitud  ni ángulo con la línea. Tal operación se puede generalizar a cualquier curva suave,  con vector tangente bien definido, si uno se preocupa de no cambiar la magnitud del vector deslizado ni tampoco el ángulo del vector a transportar con el vector tangente a la curva. Lo que es válido para el plano también es verdadero para, digamos, una esfera:  podemos transportar el vector tangente unitario de un círculo máximo sobre una esfera a todo lo largo y al irlo deslizando siempre coincidirá con el vector tangente local. Después de una vuelta, la posición  final y la inicial coincidirán.

En una superficie cualquiera, para transportar paralelamente un vector a lo largo de una geodésica, se traslada el vector de tal forma que siempre permanezca en el espacio tangente, que su norma permanezca constante y lo mismo el ángulo con el vector tangente a la geodésica. Válido en el espacio y también en cualquier superficie de Riemann. Formalmente,

\bigskip

\addtocounter{ecu}{1}

\textit{Definición (\theecu):  Sea una superficie parametrizada según $\vec x$. Sea un campo vectorial $\vec X$ tangente a la superficie:}

  $\vec X = X\sp \alpha (t)\vec x\sb \alpha (u(t))$

\textit{  Sea una curva parametrizada por $t$.}

  \textit{Definimos la \index{derivada covariante} \textbf{derivada covariante} de $\vec X$  a lo largo de la curva parametrizada por $t$ como la proyección sobre el espacio tangente de la derivada ordinaria y la notamos $\nabla X\sp \gamma/dt$ .}

  Para hallar la derivada covariante se halla la derivada ordinaria y se le resta la componente normal.

 Siguiendo los mismos pasos que en la sección pasada obtenemos que en el sistema de coordenadas del plano tangente generado por la parametrización se tiene que las componentes de la derivada covariante se escriben como:

$\nabla X\sp \gamma/dt = dX\sp \gamma /dt +  \Gamma\sp \gamma \sb{\beta \alpha}X\sp \alpha X\sp \beta$

\bigskip

\addtocounter{ecu}{1}

\textit{Definición (\theecu):
  Definimos el \index{transporte paralelo} \textbf{transporte paralelo} de $\vec X\sb 0$ a lo largo de una curva, como aquella  extensión de $\vec X\sb 0$ sobre la curva cuya  derivada covariante sea cero. Concretamente, el transporte paralelo es la única solución al problema de valor inicial:}

$ d \vec X\sp \gamma /dt +  \Gamma\sp \gamma \sb{\beta \alpha}X\sp \alpha X\sp \beta=0$

$\vec X(0) = \vec X\sb o$

\textit{donde $\vec X\sb o$ es un vector tangente a la superficie, no a la curva.}

\

La solución que vamos a obtener  nos dará un vector en cada punto que será el resultado del transporte paralelo a lo largo de la curva dada. Puesto que esto también es válido para un vector tangente unitario a una geodésica, esto nos demuestra que el vector tangente unitario a una curva parametrizada con longitud de arco puede transportarse paralelamente dando los vectores tangentes unitarios en cada punto. En general, el transporte paralelo conserva la norma y los cosenos o ángulos. El transporte paralelo es fuertemente dependiente del camino tomado.

Como un ejercicio para  afinar la intuición, tomemos una esfera y simulemos el transporte paralelo en varios casos:

1. Transportemos alrededor del ecuador un vector que originalmente mira verticalmente hacia abajo y que, por supuesto, es tangente a la esfera. A medida que los transportamos conservamos tanto la norma como el ángulo con el vector tangente local. El resultado final es que al llegar al punto de partida, el vector final coincide con el inicial.

2. Transportemos a lo largo de un círculo máximo un vector que originalmente está en el polo norte formando un ángulo de 90 grados con el vector tangente al meridiano. Lo vamos corriendo conservando la norma y los 90 grados del ángulo con el vector tangente local. Al dar una vuelta completa, se llega tal como se salió.

3. Tomemos un vector que originalmente está en el polo norte formando un ángulo de 90 grados con el vector tangente al meridiano por donde comenzamos a transportarlo. Lo llevamos hasta el ecuador, por donde comenzamos a transportarlo hasta completar media vuelta. Después subimos por el meridiano local hasta el polo norte, de donde partimos. Al llegar, encontramos que el vector final y el inicial están en dirección contraria.

\section{GENERALIZACIONES}

La caracterización de geodésicas, curvas de distancia mínima local, como curvas cuya derivada covariante sea cero necesita de la existencia de una superficie que contenga la curva, de su vector normal y de un producto interno sobre el espacio tangente que permita especificar la componente normal de un vector dado. (El producto interno sobre el tangente también permite calcular la norma de desplazamientos infinitesimales que al ser integrados a lo largo de un camino nos da el funcional de distancia entre dos puntos a lo largo de un camino).

La primera generalización que podemos hacer es una consecuencia directa del hecho de que  en la evaluación de la derivada covariante tan sólo se uso la bilinealidad y simetría del producto interno. Jamás se usó la positividad del producto interno de $\Re \sp 3$. Eso significa que podemos hablar de derivada covariante en el espacio -tiempo, en el cual hay pares de vectores cuyo cuadrado interno es negativo. También podremos usarlo en deformaciones de dicho espacio, tales como los espacios resultantes de incluir perturbaciones gravitatorias en el -tiempo.

En segundo lugar nuestra definición de derivada covariante es local y considera solamente el comportamiento de los campos en la inmediata cercanía de un punto dado. Por eso podemos generalizar nuestra definición a cualquier entidad que permita hacer lo mismo localmente, en particular en superficies de Riemann, que son aquellas que tienen estructura de variedad y para cada punto existe un abierto que lo contiene en el cual para el espacio tangente en cada punto hay un producto interno, el cual puede ser positivo, métricas de Riemann, o de cualquier signo, métrica pseudo-riemannianas. Las ecuaciones diferenciales que deben cumplir las geodésicas permanecen iguales.

En tercer lugar tenemos que liberar el concepto de geodésica de estar ligado a curvas sobre superficies específicas.  El remedio a dicho problema comienza con la siguiente definición en la cual nos liberamos  de una curva en particular.

\bigskip

\addtocounter{ecu}{1}

\textit{Definición (\theecu): La   derivada covariante de una campo vectorial $ \vec v$ a lo largo de un vector $\vec X$ es la derivada covariante del mismo campo a lo largo de una curva cualquiera, parametrizada por $t$, cuya velocidad, o vector tangente, sea el vector $\vec X$. Definida de esa forma, el operador derivada covariante $\nabla$ tiene dos operandos: el campo que se deriva y el campo que da la dirección de derivación. La notamos  $\nabla \sb{\vec X} \vec v$.}

\bigskip

\addtocounter{ecu}{1}

\textit{Teorema (\theecu): La derivada covariante definida sobre $\Re \sp 3$, $\nabla$, es un operador que en un punto dado toma vectores del espacio tangente y les asocia vectores en el espacio tangente (pues le quita todo lo que no es tangente) y como operador es lineal en ambos operandos y cumple la ley de derivación de un producto:}

$\nabla \sb{\vec X}(a\vec v + b\vec w) = a \nabla \sb{\vec X}\vec v  + b \nabla \sb{\vec X}\vec w $

$\nabla \sb{a \vec X + b \vec Y}\vec v= a \nabla \sb{\vec X}\vec v  + b \nabla \sb{\vec Y}\vec v $

$\nabla \sb{\vec X}(f \vec v) = \vec X(f)\vec v   + f \nabla \sb{\vec X}\vec v $

\textit{donde hemos usado la dualidad entre vectores y operadores de derivación.}

\bigskip

\addtocounter{ecu}{1}

\textit{Definición (\theecu): La \index{derivada covariante}  \textbf{derivada covariante} $\nabla \sb{\vec X} \vec v$ de una campo vectorial $ \vec v$ a lo largo de un vector $\vec X$  es un operador que en sus dos operandos $\vec X$, $\vec v$ cumple las propiedades del teorema anterior.}

\bigskip
\bigskip

Ahora hemos reunido suficiente material para entender de qué manera la teoría geométrica de la gravitación, o relatividad general, predice la curvatura de la trayectoria de la luz cuando pasa por un campo gravitatorio fuerte. Tal efecto se verifica al registrar la luz de estrellas que están detrás del sol en el momento de un eclipse, o con instrumental apropiado, a cualquier hora: la luz le da la vuelta al sol para llegar a nosotros.  Todo lo que necesitamos es demostrar que la forma de medir geodésicas en nuestro espacio tridimensional es diferente de la traza sobre nuestro espacio de la forma de medir distancias en el espacio-tiempo con gravitación. Por consiguiente las trayectorias vistas en nuestro espacio trazadas por la luz, que sigue geodésicas en el espacio tiempo, serán diferentes de una línea recta, que son las geodésicas en $\Re \sp 3$.

\section{INTERLUDIO}

El día señalado para la práctica de campo,  fui a dar una vuelta, un tour cerrado. Tenía la intención de medir la curvatura verdadera de todo lo que existe.  Por entre las montañas subí al cielo, pero cuando busqué la bajada, me encontré perdido. Siguiendo las delicadas instrucciones de varios ángeles pude encontrar una salida. Y  donde mi camino traspasaba la frontera entre el cielo y la tierra encontré una carretera. Y justo a la salida de mi camino estaba sentada una quimera que era mitad ángel  y mitad mujer, la cual estaba ocupada anotando algo en un registro con formato.

Al acercarme, ella me fue diciendo: yo anoto aquí las obras de los hombres y de las mujeres, pues la libertad es una de las fuerzas que curvan la trayectoria del universo.  Y al ver yo que ella era receptiva, me animé a preguntarle: ¿Por esta carretera a cuál  pueblo llego? Sin dejar su trabajo, me dijo:  todos los caminos llevan a Moniquirá. Y tratando de sobreponerme al asombro que me causó su  respuesta, de nuevo la interrogué: ¿y cuánto me demoro? Un poco más analítica, mirándome me dijo:  una hora bien andada.

Después de agradecerle, tomé mi camino  y no me atreví a volver la vista atrás.

\section{CURVATURA}

A enormes distancias nosotros, como observadores externos, podemos visualizar curvaturas estudiando la trayectoria de una partícula. Contrario a observaciones externas, el punto de vista intrínseco se preocupa por un estudio tan local y tan pobre en equipo como se pueda. Ahora bien, si una partícula viaja en caída libre describe una geodésica y no siente absolutamente nada: Cómo puede entonces una partícula saber que está bajo la acción de un campo gravitatorio? Simplemente no puede. Pero si se acompaña de otra partícula entonces puede determinar exactamente la efectividad del campo si tan sólo mide la distancia entre geodésicas. Esto significa estudiar la distancia entre dos partículas gemelas que van en caída libre y que parten de posiciones y velocidades iniciales dadas.

Un poco de atención, por favor: las líneas son geodésicas en el plano con la métrica euclidiana. Y si no son paralelas, ellas se separarán a una cierta velocidad. Pero dicha velocidad es constante y su aceleración es cero. Por tanto, la curvatura estará ligada a una aceleración no nula de la distancia entre geodésicas. Verifiquemos esto sobre una esfera con la métrica heredada de $\Re \sp 3$, en la cual las geodésicas son círculos máximos, como el ecuador.

Dos geodésicas cualesquiera sobre la esfera se cortarán en un punto. A ese punto lo llamamos el polo norte y enseguida ponemos las coordenadas esféricas que constan de la magnitud del radio $\rho$, el ángulo de barrido vertical $\phi$ que empieza en el polo norte, y el ángulo de barrido horizontal $\theta $. Aproximamos la distancia entre geodésicas $J$ por la magnitud del arco subtendido por un paralelo, el cual es perpendicular a las geodésicas que se han vuelto meridianos en la esfera. El ángulo entre los planos que generan los meridianos es  fijo, sea $\theta$ y cercano a cero.  Tenemos

 $J= (\rho sen \phi)sen  \theta  = (\rho sen \phi) \theta $.

 Al mismo tiempo, la longitud de arco, $s$, medida a lo largo de los meridianos, es $s=\rho \phi$, por lo tanto, $\phi =  s/\rho$, entonces $J=\rho \theta sen(s/\rho)$. Derivando con respecto a $s$ tenemos $dJ/ds = \theta cos (s/\rho)$, si derivamos otra vez, $d\sp 2 J/ds \sp 2 = (-1/\rho) \theta sen (s/\rho)$ y combinando $J$ con su segunda derivada obtenemos:

$d\sp 2 J/ds \sp 2 + (1/\rho\sp 2)J=0$

La interpretación del coeficiente $(1/\rho\sp 2)$ es esta: un círculo es mas curvo entre más pequeño sea. Por lo tanto, la medida más inmediata de curvatura de un círculo de radio $\rho$ debe ser $1/\rho$. Entonces, si queremos una medida de la curvatura sobre la esfera en un punto dado, trazamos, sobre el punto, una cruz, la cual genera dos círculos máximos, ambos de igual curvatura, $1/\rho$, y el producto de esas dos curvaturas $1/\rho\sp 2$ es el coeficiente de la ecuación encontrada y se conoce como curvatura de Gauss, $K$. La ecuación queda de tal forma que la aceleración a la que se separan las geodésicas  depende de un sólo parámetro: la curvatura.

\addtocounter{ecu}{1}

$$d\sp 2 J/ds \sp 2 + KJ=0 \eqno{(\theecu)}$$

Trataremos enseguida de encontrar el equivalente de $K$ en el caso general. Es algo asombroso que esta ecuación, válida para una esfera, es válida en condiciones extremadamente generales. En vista de ello, esa ecuación ha recibido su nombre propio: se denomina la \index{ecuación de Jacobi} \textbf{ecuación de Jacobi}, quien la pudo demostrar para superficies en el espacio. Lo probaremos en su generalidad aunque es  un trabajo dispendioso pero de mucho provecho. Por ahora, probaremos un detalle técnico que mas luego necesitaremos.

Las geodésicas son curvas que cumplen una cierta propiedad consignada en una ecuación diferencial. Lo que vamos a probar es válido para dos soluciones cualesquiera de una ecuación diferencial, en particular para las geodésicas. Consideremos la ecuación diferencial de las curvas integrales de un campo vectorial: $dx/dt = \vec X\sb x$, la cual tiene sentido en una vecindad de x. Una solución dada empieza en $x\sb o $ y otra empieza en $x\sb 1$, y sus soluciones respectivas se notarán con los correspondientes subíndices.

La distancia entre esas dos condiciones iniciales es un $\delta x\sb o$. Cada curva solución satisface la ecuación diferencial. Como son curvas integrales, su solución está dada en términos de un parámetro, $t$.  Si imaginamos que el campo vectorial describe las velocidades de un fluido, entonces, $\delta=\delta(x\sb 1,x\sb o) = \delta(t)= x\sb 1(t)-x\sb o(t)$ es la distancia entre dos bolitas que son arrastradas por la corriente. Acá hay dos operadores posibles: primera, restar valores de funciones entre curvas diferentes, $\delta$, y segundo, comparar valores de funciones en la misma curva tomando la derivada respecto al tiempo $d/dt$. Probaremos que estas dos operaciones conmutan entre ellas.

\addtocounter{ecu}{1}

\textit{Teorema (\theecu): si $\delta$ es una comparación entre curvas mediante resta de valores y $d/dt$ es la comparación dentro de una curva mediante la derivada, se tiene que  $\delta \circ d/dt = (d/dt)\circ \delta$}

\bigskip

Demostración:

$\delta (dx/dt) = d(x\sb 1)/dt -dx/dt =d(x+\delta x)/dt -dx/dt $

$=dx/dt+d(\delta x)/dt -dx/dt = d(\delta x)/dt$

\section{EL CONMUTADOR}

Estamos tratando de caracterizar el efecto de la curvatura sobre el distanciamiento de geodésicas. Pero no de cualquier forma, sino haciendo los ajustes necesarios para generalizar la ecuación de Jacobi. El gran progreso obtenido se desprendió de preguntas muy sencillas, todas elaboraciones de la siguiente: si una ciudad se construye en una planicie, se puede diseñar para que sus cuadras queden cuadradas y la ciudad quede cuadriculada. Pero eso no puede hacerse si se construye en terreno montañoso con valles y cerros. Curiosamente, puede hacerse si la ciudad se construye sobre una ladera, bien recta, o bien en forma de caneca, con estructura cilíndrica. Pero en general, la ciudad queda pseudo-cuadriculada, sea porque las calles sean rectas y queden rombos, o bien porque las calles sean curvas y queden cuadroides.

Pues bien, una superficie siempre se presta para ser cuadriculada por geodésicas. Ahora definamos las direcciones norte, sur, occidente, oriente. Demos un paseo generando un cuadroide: un paso hacia el norte, otro hacia el occidente, pero de exactamente igual magnitud, el siguiente hacia el sur y el último hacia el oriente: terminaremos donde empezamos? Digamos lo mismo pero  de otra manera: Usted y yo salimos del mismo punto y nos proponemos generar un cuadroide dando un paseo, cada uno por dos aristas complementarias. Tomamos todas las precauciones para que los pasos sean todos iguales y las geodésicas opuestas sean 'paralelas'. La pregunta es: nos encontraremos en la esquina opuesta?

La elaboración de esta temática, cuyo estudio general se llama holonomía, con el objetivo de generalizar la ecuación de Jacobi ha dado por lo menos dos resultados que son equivalentes pero que utilizan perspectivas diferentes. El uno es el transporte paralelo  y el otro está basado directamente sobre el conmutador de la derivada covariante. Lo de conmutador se refiere a esto: en la pregunta de encontrarse en la esquina del cuadroide, las dos personas ejecutan las mismas dos acciones pero en orden inverso. Una persona hace su paseo primero al norte, luego al occidente. La otra primero al occidente, luego al norte. En general, nos conviene tener la

\addtocounter{ecu}{1}

\textit{Definición (\theecu): Si dos operaciones de cualquier naturaleza, A y B, pueden componerse en cualquier orden, entonces se denomina \index{conmutador} \textbf{conmutador} de A y B a la operación [A,B] = AB -BA.}

\bigskip

Un detalle técnico recién  demostrado dice:

\addtocounter{ecu}{1}

\textit{Teorema (\theecu): $[\delta , d/dt]=0$}

 \bigskip

\textit{Teorema  importante para demostrar como ejercicio (\theecu): El conmutador de dos vectores del espacio tangente, considerados como operadores de derivación es también un elemento del espacio tangente, es decir un vector tangente.}

\bigskip

\addtocounter{ecu}{1}

 Este teorema es sorpresivo porque un conmutador es un operador de segundo orden mientras que un vector tangente es un operador de primer orden.
 
 \bigskip

Nosotros vamos a estudiar el conmutador de derivadas covariantes a lo largo de cuadroides y de eso sacaremos un forma de 'medir' la curvatura que generaliza la ecuación de Jacobi. Aunque siempre hemos tratado de evitar sobrecargar el material, es el momento de introducir  unas pocas definiciones, lo cual  puede redundar en claridad.

\section{REPASO}

La derivada covariante que definimos anteriormente sobre superficies en el espacio, con ayuda de curvas en una superficie y vectores normales, pudo liberarse de la curva y redefinirse en la dirección de un vector. Eso lo hicimos diciendo que la nueva derivada covariante es la derivada antigua a lo largo de cualquier curva que al ser parametrizada tenga un vector velocidad igual al vector dado. Definida as\'i, la derivada covariante cumple las siguientes propiedades:

\addtocounter{ecu}{1}

\textit{Teorema (\theecu): La derivada covariante  es un operador sobre el espacio tangente que es lineal en ambos operandos y cumple la ley de derivación de un producto:}

$\nabla \sb{\vec X}(a\vec v + b\vec w) = a \nabla \sb{\vec X}\vec v  + b \nabla \sb{\vec X}\vec w $

$\nabla \sb{a \vec X + b \vec Y}\vec v= a \nabla \sb{\vec X}\vec v  + b \nabla \sb{\vec Y}\vec v $

$\nabla \sb{\vec X}(f \vec v) = \vec X(f)\vec v   + f \nabla \sb{\vec X}\vec v $

\textit{donde hemos usado la dualidad entre vectores y operadores de derivación.}

\bigskip

Observemos que estas propiedades no dependen de ningún producto interno, de ninguna superficie, de ninguna normal. Propiedades tan liberales son al mismo tiempo muy restrictivas, lo suficiente para caracterizar la derivada covariante. Definimos entonces como derivada covariante cualquier operador que cumpla esas mismas propiedades y no cree problemas de continuidad ni derivación. Es decir, si el campo que da el sentido a la derivación es diferenciable, entonces la derivada covariante respecto a él, también. Por supuesto que la utilizaremos en superficies, variedades riemannianas o pseudo-riemannianas, con o sin producto interno, aplicada sobre campos vectoriales o sobre cualquier cosa que tenga sentido.

\addtocounter{ecu}{1}

\textit{Definición (\theecu): Un \index{marco} \textbf{marco} de campos vectoriales en una región abierta $U$ consiste de $n$ campos vectoriales diferenciables y linealmente independientes que en cada punto forman una base $\vec e = (\vec  e\sb 1, .. ..,\vec  e\sb n)$ del espacio tangente. Se dice que el marco es un marco coordenado si dicho marco resulta ser el campo de velocidades generado por un sistema coordenado, es decir, si $ \vec e\sb i = \partial /\partial x\sp i$.}

\bigskip

Un sistema coordenado es una manera de demarcar calles y carreras sobre U. En general es una deformación $\phi$ invertible de un abierto $V$ de $\Re\sp n$ que se inmersa en $U$:

$\phi : V \rightarrow U $

Las calles y carreras son las imágenes de las calles y carreras naturales de $\Re \sp n$. En ese caso: $\vec e\sb i = \partial \phi/\partial x\sp i $ tal que sólo la ordenada $i$ puede variar. Como siempre reconocemos la dualidad entre vectores y operadores de derivación, $\vec e\sb i =\partial /\partial x\sp i$ que operan sobre funciones densidad sobre $U$, es decir, sobre funciones del tipo $f: U\rightarrow \Re$. La forma de operar es a través de $\phi$ de la siguiente forma:  $\vec e\sb i (f) = \partial f/\partial x\sp i= \partial f( \phi(x\sb 1,.. ..,x\sb n))/\partial x\sp i $. No hacemos mayor diferencia entre $(x\sb 1,.. ..,x\sb n)$ en $\Re \sp n$ y su imagen $\phi(x\sb 1,.. ..,x\sb n)$, de hecho, uno tiende a pensar que ambos puntos tienen el mismo nombre y eso explica la notación tan abusiva.

\bigskip

Una condición necesaria y suficiente para que un marco sea un marco coordenado es que todos los conmutadores se aniquilen: $[\vec e\sb i, \vec e\sb j]=0$. Esta condición es una puerta abierta a opciones poderosas en la teoría de ecuaciones diferenciales parciales.

\addtocounter{ecu}{1}

\textit{Teorema (\theecu): Una derivada covariante o \index{conexión} \textbf{conexión} es una transformación afín, es decir, es de la forma $f(u) = ku + c$. Por eso, a lo que haga las veces de $c$ se le llamará coeficientes de la conexión afín.}

\bigskip

Demostración: Sea $\vec e = (\vec e\sb 1, .. ..,\vec e\sb n)$ un marco de campos vectoriales sobre $U$. Como forman una base en cualquier punto, cualquier campo vectorial $\vec X$ puede descomponerse en el marco: $\vec X =  X\sp j \vec e\sb j $ , índices repetidos indican suma. Apliquemos la derivada covariante a un campo vectorial $\vec v =  v\sp k \vec e\sb k$ en la dirección $\vec X$. Al utilizar las propiedades  de la derivada covariante tengamos presente que $\vec v$ cambia de punto a punto y  por tanto todas las coordenadas también, en realidad son funciones, y por lo tanto debe aplicarse la siguiente  propiedad   de la derivada covariante:

$\nabla \sb{\vec X}(f \vec v) = \vec X(f)\vec v   + f \nabla \sb{\vec X}\vec v $

Por tanto:

$\nabla \sb{\vec X}(\vec v)= \nabla \sb{\vec X}(v\sp k \vec e \sb k ) $
$= \vec X(v\sp k)\vec e\sb k   + v\sp k \nabla \sb{\vec X}(\vec e \sb k) $
$= (X\sp j \vec e\sb j )(v\sp k)\vec e\sb k   + v\sp k \nabla \sb{X\sp j \vec e\sb j }(\vec e \sb k) $
$=X\sp j \vec e\sb j (v\sp k)\vec e\sb k   + v\sp k X\sp j \nabla \sb{ \vec e\sb j }(\vec e \sb k) $, donde aplicamos la linealidad

Ahora bien, la derivada covariante de vectores en el espacio tangente produce vectores en el espacio tangente, por tanto, puede descomponerse en el marco dado:

 $\nabla \sb{ \vec e\sb j }(\vec e \sb k) =\omega \sp i \sb{jk}\vec e\sb i$

lo cual se lee : $\omega \sp i \sb{jk}\vec e\sb i$ es la componente $k$ de la variación covariante de $\vec e\sb i$ en la dirección de $\vec e\sb j$. Se denomina a esta familia de coeficientes los coeficientes de la conexión afín.

En un marco coordenado se tiene que $\omega \sp i \sb{j,k} = \Gamma \sp i \sb{j,k}$, los símbolos de Christoffel.

Con esta notación queda:

$\nabla \sb{\vec X}(\vec v) = X\sp j \vec e\sb j (v\sp k)\vec e\sb k   + v\sp k X\sp j \omega \sp i \sb{jk}\vec e\sb i $

$=X\sp j \vec e\sb j (v\sp i)\vec e\sb i   + v\sp k X\sp j \omega \sp i \sb{jk}\vec e\sb i $,  cambiando la $k$ por $i$ en el primer término,

$=[ \vec e\sb j (v\sp i)   + v\sp k \omega \sp i \sb{jk}]X\sp j\vec e\sb i $

Entonces

$$\nabla \sb{\vec X}(\vec v) =[\partial v\sp i /\partial x\sp j + v\sp k \omega \sp i \sb{jk}] X\sp j\vec e\sb i$$

Debido a que por dualidad se define $dg(\vec x) = \vec x(g)$ para funciones escalares $g$, entonces $\partial v\sp i /\partial x\sp j =dv\sp i(\partial /\partial x\sp j)=  dv\sp i(\vec e\sb j)$ y que la 1-forma dual de $e\sb j$, notada $\sigma \sp j$ operada sobre $\vec X$ toma la j-ésima coordenada, $\sigma \sp j(\vec X) = X\sp j $, por tanto podemos reescribir

$\nabla \sb{\vec X}(\vec v) = [dv\sp i (e\sb j) +v\sp k \omega \sp i \sb{jk}] X\sp j\vec e\sb i$
$=[dv\sp i (X\sp j e\sb j) + v\sp k \omega \sp i \sb{jk}X\sp j ] \vec e\sb i$

$=[dv\sp i (\vec X) + v\sp k \omega \sp i \sb{jk}\sigma \sp j (\vec X) ] \vec e\sb i$

En conclusión

$$\nabla \sb{\vec X}(\vec v) =[dv\sp i + v\sp k \omega \sp i \sb{jk}\sigma \sp j ] (\vec X) \vec e\sb i$$

Observemos que lo que está entre corchetes es una 1-forma, que al operar sobre $\vec X$ da un real que es la coordenada i-ésima de la derivada covariante de $\vec v$ en la dirección de $\vec X$. Cada componente está dada por una 1-forma, pero hay $n$ componentes, por lo que la derivada covariante debe poderse reemplazar por una entidad complicada que tiene 1-formas en todo lado.

\addtocounter{ecu}{1}

\textit{Definición (\theecu): Definimos la torsión de una conexión $\tau$ como:}

$\tau (\vec X, \vec Y) = \nabla \sb{\vec X} \vec Y - \nabla \sb{\vec  Y }\vec X-[\vec X,\vec Y] $

\bigskip

Cuando una conexión no tiene torsión se dice que es simétrica y en ese caso $\omega \sp i \sb{jk} = \omega \sp i \sb{kj}$. La torsión de una curva se interpreta como la tendencia de la curva a salirse del plano osculador, es decir, del plano que mejor se ajusta a la curva en un determinado punto. Por eso, en algunos casos la torsión de una conexión se interpreta como la no cerradura de las geodésicas sobre ellas mismas.

Cuando la conexión viene de un marco coordenado se llama de Levi-Civita y es simétrica. Otra conexión muy importante: la de Riemann. Esa conexión  se define en una variedad de Riemann, es decir, una que tiene un producto interno en cada espacio tangente, y cumple la importante igualdad:

\addtocounter{ecu}{1}

$$d/dt<\vec X , \vec Y>= <\nabla\sb t  \vec X,\vec Y> + <\vec X,\nabla \sb t \vec Y> \eqno{(\theecu)}$$

donde $\vec X$ y $\vec Y$ son dos campos vectoriales definidos a lo largo de una curva parametrizada por $t$.

\section{LA ECUACION DE JACOBI}

Definiremos lo que se entiende por curvatura y veremos  emerger la conocida ecuación de Jacobi pero en una forma muy general.

Nuestro escenario será una superficie parametrizada en una variedad. La podemos entender como una operación semejante a la que ejecuta una persona que toma un sábana y la bate en el aire, pero en un determinado instante queda congelada. La sábana se llama $U$ y el aire es la variedad. A cada punto de la sábana (un conjunto abierto del plano) le corresponde un punto $x$ de la variedad. $U$ tiene sus propias coordenadas $u, v$.

Consideramos un campo vectorial sobre la superficie, es decir, una asignación que a cada par de valores $(u,v)$ le asigna un elemento del espacio tangente en el $x$ correspondiente. Dos ejemplos de tales campos $\partial x/\partial u, \partial x/\partial v$.

\addtocounter{ecu}{1}

\textit{Teorema (\theecu): Sea $S$ una superficie dentro de una variedad provista de una conexión simétrica. Entonces}

$\frac{\rm \nabla}{\rm \partial u } (\frac{\rm \partial x}{\rm \partial v}) =\frac{\rm \nabla}{\rm \partial v} (\frac{\rm \partial x}{\rm \partial u})  $

\textit{donde $\frac{\rm \nabla}{\rm \partial u }$ y  $\frac{\rm \nabla}{\rm \partial v }$
denotan las derivadas covariantes a lo largo de las trayectorias generadas por las coordenadas $u$ y $v$  que parametrizan la superficie.}

\bigskip

Este teorema especifica una forma de conmutar derivadas, conmutando los denominadores.  El teorema es tanto más importante al notar que tanta libertad no es posible a cada paso. Precisamente

\addtocounter{ecu}{1}

\textit{Definición(\theecu): Definimos como curvatura de una conexión al operador $\Omega $ dado por}

$\Omega(\vec X, \vec Y) = [\nabla_{ \vec X}, \nabla_{ \vec Y}] - \nabla \sb{[\vec X,  \vec Y] }$

\bigskip

Sobre   campos vectoriales es lineal y además sobre uno específico  $\vec v$ este operador cumple:

$\Omega (\vec X, \vec Y)\vec v = \nabla \sb{\vec X} (\nabla \sb{\vec Y}\vec v) -\nabla \sb{\vec Y}(\nabla \sb{\vec X} \vec v )- \nabla\sb{[ \vec X,  \vec Y]}\vec v$

Puesto que $\Omega $ se basa enteramente sobre la derivada covariante, toma vectores del espacio tangente en un punto dado y  siempre produce elementos del mismo espacio.

\addtocounter{ecu}{1}

\textit{Teorema (\theecu): Si se toma un marco coordenado, con elementos $\partial \sb i $, para expandir la curvatura, se obtiene:}

$\Omega(\vec X, \vec Y)\vec v = \Omega \sp i \sb{jkl} X\sp k Y\sp l v\sp j \partial \sb i$

\textit{donde }$\Omega \sp i \sb{jkl} = \partial \sb k \omega \sp i \sb{lj}- \partial \sb l \omega \sp i \sb{kj}+\omega \sp i \sb{kr}\omega \sp r \sb{lj}-\omega \sp i \sb{lr} \omega \sp r \sb{kl}$

\addtocounter{ecu}{1}

\textit{Teorema (\theecu): Cuando $[ \vec X,  \vec Y]=0$, como en el caso de que estos campos definan un marco coordenados, entonces}

$\Omega (\vec X, \vec Y) = [\nabla_{ \vec X}, \nabla_{ \vec Y}]$, es decir

 $\Omega(\vec X, \vec Y)\vec v = \nabla \sb{\vec X} (\nabla \sb{\vec Y}\vec v) -\nabla \sb{\vec Y}(\nabla \sb{\vec X} \vec v )$
 
 \bigskip

Si estamos en la superficie parametrizada por $(u,v)$, y la parametrización es invertible, entonces con toda seguridad los vectores velocidad coordenados $\vec X = \partial x/\partial u$ y
$\vec Y = \partial x/\partial v$ conmutan entre ellos y se cumple entonces que:

$\Omega(\vec X, \vec Y) \vec w=\Omega (\partial x/\partial u, \partial x/\partial v) \vec w =
\frac{\rm \nabla}{\rm \partial u } (\frac{\rm \nabla \vec w}{\rm \partial v}) -\frac{\rm \nabla}{\rm \partial v} (\frac{\rm \nabla \vec w}{\rm \partial u})  $

\addtocounter{ecu}{1}

\textit{Definición (\theecu): Decimos que una curva $\vec r$ parametrizada por $t$ es una geodésica si su vector velocidad $d\vec r/dt$ tiene derivada covariante con respecto a ella misma igual a cero. }

\bigskip

Obsérvese la notación:

$\frac{\rm \nabla}{\rm dt} (\frac{\rm d\vec r}{\rm dt})=0$

\bigskip

Queda sobreentendido que siempre parametrizamos con longitud de arco.

\bigskip

La curvatura que hemos definido  no es una 2-forma, sino que dados dos vectores que se utilizan para definirla,  es una transformación lineal que toma un vector del espacio tangente y le hace corresponder otro vector del mismo espacio tangente. Por lo tanto, nuestra curvatura funciona como una 3-forma con valores en el espacio tangente. En otros contextos, al vector imagen se le multiplica en producto punto o interior por otro vector para que quede una 4-forma.

\section{CAMPOS DE JACOBI}

Estamos buscando una ecuación que describa la distancia  entre geodésicas ( es decir, entre dos partículas que se muevan por una geodésica) en términos de curvatura. Puesto que la curvatura está relacionada con   la derivada covariante, y nuestras geodésicas vienen de una métrica, se presume la existencia de una superficie y de un vector normal a ella. Pero, de dónde sale una superficie?

Pues, bien, dadas dos geodésicas fabricamos una superficie con ellas: el segmento que une las dos partículas que generan las geodésicas describe una superficie, que además es suave y todo lo que se necesite. Ese es el comienzo de todo. Después, parametrizamos dicha superficie. Nos resulta el ambiente en el que hemos trabajado para definir geodésicas. La parametrización tiene la forma $\phi : U \rightarrow V$, donde $U$ es un abierto del plano que contiene el cuadrado unitario con vértice inicial  en (0,0), y $V$ es una abierto de alguna variedad y $\phi (u,v) = x$. Dicha parametrización debe cumplir: $\phi(u,0) = x\sb 1(u) $, la primera geodésica, y que para algún $\alpha $ se tenga que $\phi(u, \alpha) = x\sb 2$, la segunda geodésica. Además, cada curva $C\sb{\alpha} = \phi(u, \alpha)$ también es una geodésica. Y por supuesto, $u$ parametriza todas las curvas en longitud de arco.

Es justo preocuparse por la realización de la anterior definición, pues quizá estemos exigiendo demasiado. Para empezar, podemos imaginar que por cada punto que queda en el segmento que conecta nuestras dos partículas existe una partícula y que a todo ese bunch de partículas se les sincronizar para salgan al tiempo de una posición inicial.

La separación $\delta(u) = x\sb 1(u) -x\sb 2 (u)$ en general puede ser una función complicada. Restrinjámonos  a la parte lineal de tal función, por lo que los resultados que obtengamos han de ser válidos, por lo menos, para tiempos y separaciones cortos -lo cual habrá de ser más que suficiente.  Con ese fin, substituimos tal separación, $\delta(u)$,  por el vector en dirección $\vec J = \partial \phi /\partial v $  pero evaluando tal derivada en $v=0$, es decir, sobre  la curva $x\sb 1$ y estudiamos la influencia de $\vec J$ sobre la derivada covariante de una geodésica. Al vector tangente, que es unitario, lo notamos $\vec T$.

Puesto que cada $C\sb \alpha$ es una geodésica, tenemos que $\nabla \vec  T/d u=0$ para todo $\alpha $, lo cual significa que $\vec T$ permanece invariante ante el transporte paralelo a lo largo de la geodésica. Por tanto, teniendo en cuenta que $\vec J$ es el vector transversal cuyo parámetro es $v$, y $\vec T$ es el vector tangente, cuyo parámetro es $u$, tomamos la segunda derivada y obtenemos que la derivada de cero también es cero:

$ \frac{\rm \nabla}{\rm dv} (\frac{\rm \nabla \vec T}{\rm du})=0$

Pero como la curvatura $\Omega$ es tal que para cualquier vector $\vec w$ se cumple que:

$\Omega(\vec X, \vec Y) \vec w=\Omega (\partial x/\partial u, \partial x/\partial v) \vec w =
\frac{\rm \nabla}{\rm \partial u } (\frac{\rm \nabla \vec w}{\rm \partial v}) -\frac{\rm \nabla}{\rm \partial v} (\frac{\rm \nabla \vec w}{\rm \partial u})  $

entonces,

$\Omega(\vec J, \vec T) \vec T=\Omega (\partial x/\partial v, \partial x/\partial u) \vec T =
\frac{\rm \nabla}{\rm \partial v} (\frac{\rm \nabla \vec T}{\rm \partial u}) -\frac{\rm \nabla}{\rm \partial u } (\frac{\rm \nabla \vec T}{\rm \partial v})  $

Por tanto,

$0 = \frac{\rm \nabla}{\rm dv} (\frac{\rm \nabla \vec T}{\rm du})= \Omega(\vec J, \vec T) \vec T + \frac{\rm \nabla}{\rm \partial u } (\frac{\rm \nabla \vec T}{\rm \partial v}) $

$= \Omega(\vec J, \vec T) \vec T + \frac{\rm \nabla}{\rm \partial u } \frac{\rm \nabla }{\rm \partial v}(\frac{\rm \partial x}{\rm \partial u}) $

$= \Omega(\vec J, \vec T) \vec T + \frac{\rm \nabla}{\rm \partial u } \frac{\rm \nabla }{\rm \partial u}(\frac{\rm \partial x}{\rm \partial v}) $

$=\Omega(\vec J, \vec T) \vec T + \frac{\rm \nabla}{\rm \partial u } (\frac{\rm \nabla \vec J}{\rm \partial u})$

que escrito ordenadamente queda:

\addtocounter{ecu}{1}

$$(\frac{\rm \nabla \sp 2 \vec J}{\rm \partial u\sp 2}) + \Omega(\vec J, \vec T) \vec T =0 \eqno{(\theecu)}$$

Esta \index{ecuación de Jacobi}  \textbf{ecuación de Jacobi} realmente se reduce a la que obtuvimos con una esfera. Para mayor generalidad, consideremos, pues, una superficie cualquiera en el espacio, parametrizada como sabemos. El producto interno es el normal, el cual genera la conexión Riemanniana, aquella que conserva los ángulos. Sea $C$ una geodésica con vector tangente unitario $\vec T$ y sea $\vec T \sp{\perp}$ aquel vector del plano tangente que es perpendicular a $\vec T$.

Puesto que $C$ es una geodésica, $\nabla \vec T/ds =\nabla\sb s \vec T =0$. Probemos que con una conexión Riemanniana, entonces se cumple que $\nabla \sb s \vec T\sp{\perp}=0$. En efecto, sea

$d/ds<\vec X , \vec Y>= <\nabla\sb s  \vec X,\vec Y> + <\vec X,\nabla \sb s \vec Y> $

y en nuestro caso, con $<\vec T , \vec T\sp{\perp}>=0$, tenemos:

$0= d/ds<\vec T , \vec T\sp{\perp}>= <\nabla\sb s  \vec T,\vec T\sp{\perp}> + <\vec T,\nabla \sb s \vec T\sp{\perp}> $

y reemplazando

$0=d/ds<\vec T , \vec T\sp{\perp}>= <0,\vec T\sp{\perp}> + <\vec T,\nabla \sb s \vec T\sp{\perp}> $

por lo que

$<\vec T,\nabla \sb s \vec T\sp{\perp}> =0$

por tanto se deduce que $\nabla \sb s \vec T\sp{\perp}$   es cero, porque no puede ser perpendicular a $\vec T$,  pues debería tener al menos una componente paralela a $\vec T$.

Ahora mostremos que hemos generalizado la ecuación de Jacobi para una esfera.
Sea $\vec J$ el vector que linealmente aproxima la distancia entre $C$ y una geodésica vecina. Como $\vec T$ con su complementario perpendicular forman una base para el espacio tangente, entonces cualquier vector de dicho espacio se puede descomponer en dicha base. Sea $s$ el parámetro de $C$ y que también es parámetro de $\vec T$, de su perpendicular y de $\vec J$. Tenemos:

$\vec J (s) = x(s)\vec T + y(s)\vec T\sp{\perp}$

podemos  calcular la segunda derivada covariante de $\vec J$, teniendo en cuenta que hay que derivar como un producto y que tanto  $ \nabla\sb s \vec T$ como $\nabla \sb s \vec T\sp{\perp}$ son cero :

$(\frac{\rm \nabla  \vec J}{\rm \partial s }) = x' (s)\vec T + y' (s)\vec T\sp{\perp}$

$(\frac{\rm \nabla \sp 2 \vec J}{\rm \partial s\sp 2}) = x''(s)\vec T + y''(s)\vec T\sp{\perp}$
$=-\Omega(\vec J, \vec T) \vec T = $
$=-\Omega(x(s)\vec T + y(s)\vec T\sp{\perp}, \vec T) \vec T $

o sea

$x''(s)\vec T + y''(s)\vec T\sp{\perp}=-\Omega(x(s)\vec T + y(s)\vec T\sp{\perp}, \vec T) \vec T $

$=-x(s)\Omega(\vec T , \vec T) \vec T $
$               -y(s)\Omega( \vec T\sp{\perp}, \vec T) \vec T $

donde utilizamos la linealidad de $R$. Pero el primer término del lado derecho es cero, pues en todos lados queda derivada covariante del vector tangente con respecto a el mismo, lo que es cero.

$x''(s)\vec T + y''(s)\vec T\sp{\perp}=-y(s)\Omega( \vec T\sp{\perp}, \vec T) \vec T $

Multiplicando en ambos lados en producto interior por $\vec T\sp{\perp}$, y simplificando tenemos:

$y''(s)=-y(s)<\Omega( \vec T\sp{\perp}, \vec T) \vec T,\vec T\sp{\perp}> $

Dando a $<\Omega( \vec T\sp{\perp}, \vec T) \vec T,\vec T\sp{\perp}>$ el nombre de $K$ la ecuación se reescribe como esperábamos: $y''(s)=-K y(s)$, o sea,

\addtocounter{ecu}{1}

$$y''(s)+K(s)y(s)=0 \eqno{(\theecu)}$$

$K(s)$ puede ser constante como en el caso de un plano, un cilindro o una esfera.

\section{TRANSPORTE  Y CURVATURA}

En una superficie cualquiera, para transportar paralelamente un vector a lo largo de una geodésica, se traslada el vector de tal forma que siempre permanezca en el espacio tangente, que su norma permanezca constante y lo mismo el ángulo con el vector tangente a la geodésica. Válido en el espacio y también en cualquier superficie de Riemann.

En la tierra, considerándola esférica, los círculos máximos son las geodésicas. A lo largo de un viaje cerrado con transporte paralelo, un vector cualesquiera queda invariante en éste caso. Pero si tomamos en la tierra un viaje cerrado empezando con el vector tangente a un meridiano desde el polo norte, tomando a lo largo del mismo meridiano hasta el ecuador y luego por el ecuador hasta un cuarto de vuelta y después subimos por el meridiano correspondiente, entonces el vector original y el final estarán a un ángulo de $\pi /2$.

En efecto, bajando del polo norte hacia el ecuador, el vector tangente termina mirando hacia el polo sur. Al ser deslizado horizontalmente, seguirá mirando hacia el polo sur. Al tomar hacia el polo norte, seguirá mirando al polo sur. Al llegar de vuelta al polo norte, estará horizontal, su sombra corresponderá exactamente con el meridiano por donde subió, pero formará un ángulo recto con el vector tangente inicial.  Este desfase es claramente un resultado de la curvatura y el propósito de la presente sección es demostrar que este tipo de mediciones es de importancia universal y que se puede utilizar para estimar la curvatura de una superficie dada.

El transporte paralelo está por todos lados. Por ejemplo, en la medición de Gauss sobre los ángulos internos de un triángulo formado por tres picos montañosos. En efecto, uno siempre supone que la luz le llega al observador en línea recta, es decir definida por el  vector tangente a la geodésica por donde viaja la luz y que dicho vector al ser deslizado paralelamente a la curva llega a ser el vector tangente local. Decimos que la derivada covariante del vector tangente a lo largo de la curva es cero. Todo eso es fácil por ser una geodésica.

Cuando haya que transportar en paralelo a un vector a lo largo de una curva que no sea una geodésica, se aproxima la curva por geodésicas, se transporta el vector por dicha aproximación y se toma el límite hasta que la aproximación por geodésicas tienda a la curva dada (por favor, investigar en qué métrica y qué pasa en las esquinas).

El transporte paralelo  es un concepto que se puede tomar como fundamental, en el sentido que permite una redefinición de la derivada covariante y curvatura. Veamos.

\section{REDEFINICION DE DERIVADA COVARIANTE}

Un objetivo general es estudiar cómo cambia un campo vectorial a medida que es arrastrado por otro campo. El problema básico es que para comparar dos vectores se requiere que pertenezcan al mismo espacio vectorial. En una variedad, a puntos diferentes le corresponde espacios tangentes que son diferentes así sean isomorfos. La derivada covariante resuelve esta preocupación basándose en el transporte paralelo, mientras que la derivada de Lie lo hace con el pull-back.

\subsection{Por transporte paralelo}

Consideremos los campos $\vec u$ y  $\vec v$. Queremos averiguar cómo cambia $\vec v$ a medida que es arrastrado por $\vec u$. Esa terminología significa lo siguiente: podemos imaginar que $\vec u$ es el campo de velocidades de un fluido que no tiene ni propiedades eléctricas ni magnéticas. Quizá nos interese saber cómo sería la dinámica de una suspensión de partículas sensibles al campo electromagnético. No es tan fácil predecirlo. Pero seguramente nos ayude entender cómo varía el campo electromagnético por las líneas de flujo del fluido. Por eso decimos que queremos saber  cómo cambia $\vec v$ a medida que es arrastrado por $\vec u$.

Sobre una línea de flujo, parametrizada por $\phi(t)$, comparamos dos posiciones cercanas, una denotada por $\phi (t + \delta t)$ y la otra por $\phi(t)$. El campo electromagnético $\vec v$ toma los valores $\vec v(\phi (t + \delta t))$ y  $\vec v(\phi (t ))$ respectivamente. Para comparar esos dos valores  de  $\vec v$, tomamos un punto de operación, el dado por $\phi(t)$. Luego transportamos en transporte paralelo inverso al vector  $\vec v(\phi (t + \delta t))$ hasta nuestro punto de operación, a través de la línea de flujo de  $\vec u$. Al resultado de ese transporte antiparalelo lo notamos  $T\sb{-\delta t}( \vec v(\phi (t + \delta t)))$ y definimos la derivada covariante $\nabla \sb{\vec u}\vec v$ del campo $\vec v$ respecto a $\vec u$ como sigue:

\addtocounter{ecu}{1}

\textit{Definición} (\theecu):

\textit{Definimos la}   \textbf{derivada covariante} \textit{del campo $\vec v$ a lo largo del campo $\vec v$ como:}

 $\nabla \sb{\vec u}\vec v = lim \sb{\delta t\rightarrow 0}\frac{\rm T\sb{-\delta t}( \vec v(\phi (t + \delta t)))- \vec v(t)}{\rm \delta t}$

\textit{donde $\phi $ es la curva integral, o línea de flujo, de $\vec u$ que pasa por P.}

 Hay que estar atentos, porque nuestra definición puede aplicarse, sin previo aviso,  a cosas mucho más generales que la derivada covariante del campo electromagnético a lo largo de las líneas de flujo de un fluido.

\subsection{Por derivada de Lie}

Esta derivada ataca la comparación de un campo vectorial a medida que es arrastrado por otro, mediante el pull back. La líneas de flujo de un campo generan un homomorfismo $\phi$ de la variedad en sí misma. Tal homomorfismo es diferenciable, si el campo es continuo y su diferencial se nota $d\phi$. La diferencial de un homomorfismo $\phi$ es una transformación que a escala infinitesimal funciona como $\phi$ pero que transforma plps infinitesimales en plps infinitesimales, al contrario de $\phi$ que en general no es lineal ni afín. La diferencial está definida sobre los espacios tangentes correspondientes al punto de partida y de llegada del homomorfismo. Si tomamos el pull-back de $d\phi$ podemos comparar dos vectores que pertenecen a espacios tangentes distintos:

\addtocounter{ecu}{1}

\textit{Definición} (\theecu):

\textit{Definimos la} \index{derivada de Lie} \textbf{derivada de Lie} $L \sb{\vec u}\vec v$ \textit{del campo $\vec v$ respecto a $\vec u$ como sigue:}

 $L \sb{\vec u}\vec v = lim \sb{\delta t\rightarrow 0}\frac{\rm (d\phi)\sp *( \vec v(\phi (t + \delta t)))- \vec v(t)}{\rm \delta t}$

\textit{donde $\phi $ es la curva integral de $\vec u$ que pasa por P.}

\

\addtocounter{ecu}{1}

\textit{Teorema }(\theecu):

\textit{1. La deriva covariante es bilineal, es decir, es lineal tanto en el campo que es arrastrado y comparado como en el que arrastra.}

2. $\nabla \sb{\vec u}f\vec v = (L \sb{\vec u}f)\vec v + f(x)\nabla \sb{\vec u}\vec v$ donde $f$ es una función escalar.

3. $L \sb{\vec u}<\vec v, \vec w>= <\nabla \sb{\vec u}\vec v, \vec w(x)>
+ <\vec v (x), \nabla \sb{\vec u}\vec w >$

4. $\nabla \sb{\vec u}\vec v- \nabla \sb{\vec v}\vec u = [\vec v, \vec u]$

\section{REDEFINICION DE CURVATURA}

Mientras que Gauss trabajó con ángulos internos, también se puede trabajar con áreas. Las siguientes observaciones argumentan informalmente un resultado relacionado, muy interesante, enunciado como una redefinición de curvatura.

1. Lo que nos interesa es estudiar el desfase o ángulo que sufre un vector después de un tour cerrado en transporte paralelo. En primer término, tenemos que dicho desfase es función aditiva de regiones en el siguiente sentido: tomamos dos regiones  que compartan  el punto de partida, que es el mismo de llegada. Al rodear la primera región, se creará un desfase dado, $\phi\sb 1$ y al rodear la segunda, se creará  $\phi\sb 2$, entonces al rodear las dos regiones se creará un desfase total igual a $\phi\sb 1 + \phi\sb 2$. Obsérvese que esto es válido, aún si las regiones comparten parte de la frontera, pues en este caso, el desfase sobre la primera región se  puede descomponer en dos partes $\phi \sb 1=\phi \sb{11} + \phi\sb{12} $ correspondientes a la parte que sólo pertenece a la primera región, $\phi \sb{11}$, y a la que comparten las dos, $\phi\sb{12}$. Similarmente,   $\phi \sb 2=\phi \sb{22} + \phi\sb{21} $ Por lo tanto, el desfase total es, como se esperaba, $\phi \sb{22} + \phi\sb{11} $ pues el desfase compartido se cancela, una vez andado y otra vez desandado.

2. Para poder integrar desfases sobre superficies a partir de áreas infinitesimales, debemos creer que el desfase creado por una región tiende a cero si el diámetro de dicha región  es infinitesimal, comparada con la región a integrar. Aquí se usa fuertemente el hecho de que la variedad es suave y que por lo tanto, a caminos cercanos le corresponde desfases cercanos.

3. En un plano, el desfase a lo largo de cualquier camino cerrado es cero.

\bigskip

Todo esto nos motiva la siguiente redefinición constructiva: la curvatura es una 2-forma que en un punto dado toma un 2-plp infinitesimal del espacio tangente, y mide el desfase creado al recorrerlo en transporte paralelo en una dirección dada. Formalmente, tenemos mediciones sobre plps infinitesimales:

\addtocounter{ecu}{1}

\textit{Re-definición de curvatura} (\theecu):

\textit{La} \index{curvatura}  \textbf{curvatura }$\Omega $ \textit{es una 2-forma, que toma un 2-plp $(\xi, \eta)$, lo miniaturiza produciendo $(\epsilon \xi, \epsilon \eta)$, lo orienta, lo recorre midiendo el desfase creado en transporte paralelo, $\phi (\epsilon)$ y calcula el límite:}

$\Omega (\xi, \eta)= lim \sb{\epsilon \rightarrow 0}\frac{\rm \phi(\epsilon)} {\rm \epsilon}$

\bigskip

\addtocounter{ecu}{1}
\textit{Teorema }(\theecu):

\textit{Si sobre un mismo plp, $dS$ es la 2-forma que mide el área del plp y $\Omega$ es la 2-forma curvatura, entonces se tiene que $\Omega = \kappa dS$ donde la constante $\kappa$  es la curvatura escalar o de Riemann.}

La curvatura escalar y la 2-forma curvatura se relacionan por

$\kappa = <\Omega(\xi, \eta)(\xi),\eta> $

La manera de implantar el transporte paralelo  y el ángulo de desfase correspondiente en electromagnetismo se deduce del análisis del experimento de Bohm-Aharonov, donde el ángulo que se contabiliza es la fase de la función de onda.

\subsection{ Curvatura a la Gauss}

Ahora vamos a  atacar  la relación entre deslizamiento paralelo y curvatura. Nuestro objetivo estará orientado por la medición de Gauss, relacionando el desfase total de un vector a lo largo de un camino cerrado y la curvatura. Siempre trabajaremos sobre una superficie parametrizada biunívocamente con un solo retazo. Pero es evidente que todo puede generalizarse.

El siguiente resultado utiliza la bilinealidad del producto interno:

\addtocounter{ecu}{1}

\textit{Lema (\theecu): Sea $\vec X(t)$ y $\vec Y(t)$ dos campos vectoriales definidos a lo largo de una curva parametrizada por $u = u(t)$. Entonces:}

$(d/dt) <\vec X(t), \vec Y(t)> = <d \vec X(t)/dt, \vec Y(t)> + <\vec X(t), d\vec Y(t)/dt> $

\bigskip

\addtocounter{ecu}{1}

\textit{Corolario (\theecu):  Si $\vec X(t)$ y $\vec Y(t)$ dos campos vectoriales definidos a lo largo de una curva parametrizada por $u = u(t)$ y tangentes a una superficie, entonces:}

$(d/dt) <\vec X(t), \vec Y(t)> = <\nabla \vec X(t)/dt, \vec Y(t)> + <\vec X(t), \nabla \vec Y(t)/dt> $

\textit{donde $\nabla \vec X(t)/dt $ denota la derivada covariante de  $  \vec X(t)$ a lo largo de la curva $u = u(t)$, parametrizada por longitud de arco.}

\bigskip

Demostración: Teniendo en cuenta que la derivada covariante resta de la derivada ordinaria el componente normal al espacio tangente, podemos escribir

$\nabla \vec X(t)/dt = d \vec X(t)/dt - \alpha \vec N$, por  tanto
$d \vec X(t)/dt = \nabla \vec X(t)/dt + \alpha \vec N$

$\nabla \vec Y(t)/dt = d \vec Y(t)/dt - \beta \vec N$,  por  tanto
$d \vec Y(t)/dt = \nabla \vec Y(t)/dt + \beta \vec N$

Reemplazando en la regla del producto,

$(d/dt) <\vec X(t), \vec Y(t)> $

$= < \nabla \vec X(t)/dt + \alpha \vec N, \vec Y(t)> + <\vec X(t), \nabla \vec Y(t)/dt + \beta \vec N> =$

$< \nabla \vec X(t)/dt, \vec Y(t)>+  <\alpha \vec N, \vec Y(t)> + <\vec X(t), \nabla \vec Y(t)/dt> $

$+ <\vec X(t), \beta \vec N> $

Utilizando la perpendicularidad del vector normal y el espacio tangente, donde están $\vec X$ y $\vec Y$, nos queda:

$(d/dt) <\vec X(t), \vec Y(t)> = <\nabla \vec X(t)/dt, \vec Y(t)> + <\vec X(t), \nabla \vec Y(t)/dt> $

\

\addtocounter{ecu}{1}

\textit{Teorema (\theecu):  Si $\vec X(t)$, $\vec Y(t)$  son dos campos vectoriales tangentes a una superficie, entonces, utilizando la dualidad entre vectores y operadores de derivación tenemos para un campo $T(t)$:}

$\vec T<\vec X(t), \vec Y(t)> = <\nabla \sb{\vec T}  \vec X(t), \vec Y(t)> + <\vec X(t), \nabla  \sb{\vec T}  \vec Y(t)>$

\bigskip

Demostración:
Sea $\vec T = \alpha \sp i \vec e\sb i = \alpha \sp i \partial /\partial x \sp i$

Reemplazando en el lado derecho del teorema a demostrar:

 $= <\nabla \sb{\alpha \sp i \vec e\sb i}  \vec X(t), \vec Y(t)> + <\vec X(t), \nabla  \sb{\alpha \sp i \vec e\sb i}  \vec Y(t)>$

utilizando la linealidad de la derivada covariante, queda:

$ =\alpha \sp i<\nabla \sb{ \vec e\sb i}  \vec X(t), \vec Y(t)> + \alpha \sp i <\vec X(t), \nabla  \sb{\vec e\sb i}  \vec Y(t)>$

Consideremos la curva integral de $\vec e\sb i= \partial /\partial x \sb i$, es decir, la imagen del  eje coordenado $i$ en el retazo del cual se hizo la superficie. Tenemos que $x\sp i$
es un parámetro y sin pérdida de generalidad, puede considerarse como longitud de arco. En estas circunstancias podemos aplicar

$\nabla \sb{ \vec e\sb i}  \vec X(t) = \nabla   \vec X(t)/\partial x\sp i$

$\nabla \sb{ \vec e\sb i}  \vec Y(t) = \nabla   \vec Y(t)/\partial x\sp i$

Substituyendo en la expresión que estamos elaborando, nos queda

$ =\alpha \sp i<\nabla   \vec X(t)/\partial x\sp i, \vec Y(t)> + \alpha \sp i <\vec X(t),\nabla   \vec Y(t)/\partial x\sp i>$

$ =\alpha \sp i[<\nabla   \vec X(t)/\partial x\sp i, \vec Y(t)> +  <\vec X(t),\nabla   \vec Y(t)/\partial x\sp i>]$

y aplicando el corolario anterior

$=\alpha \sp i \partial /\partial x \sp i [<\vec X(t),\vec Y(t)>]$

lo cual se puede simplificar recuperando a $\vec T$, terminando la demostración:

$ =\vec T<\vec X(t), \vec Y(t)> $

\

\addtocounter{ecu}{1}

\textit{Corolario} (\theecu):

$d<\vec e\sb i, \vec e\sb j> (\vec e\sb k)= <\nabla \sb{\vec e\sb k}  \vec e\sb i, \vec e\sb j> + <\vec e\sb i, \nabla  \sb{\vec e\sb k} \vec e\sb j>$

Demostración:

Por el teorema anterior y utilizando la dualidad $\vec e\sb k= \partial /\partial x\sp k$ tenemos:

$\partial /\partial x\sp k <\vec e\sb i, \vec e\sb j> = <\nabla \sb{\vec e\sb k}  \vec e\sb i, \vec e\sb j> + <\vec e\sb i, \nabla  \sb{\vec e\sb k} \vec e\sb j> $

pero por dualidad para una 1-forma $df$ y un vector $\vec v $ considerado como operador $\vec v f = df(\vec v)$, por tanto si $f = <\vec e\sb i, \vec e\sb j> $, $df = d<\vec e\sb i, \vec e\sb j>$, y $\vec v =\vec e\sb k = \partial /\partial x\sp k $ obtenemos:

$\partial /\partial x\sp k <\vec e\sb i, \vec e\sb j> = (d<\vec e\sb i, \vec e\sb j>)(\partial /\partial x\sp k) = (d<\vec e\sb i, \vec e\sb j>)(\vec e\sb k) $

y por transitividad

$(d<\vec e\sb i, \vec e\sb j>)(\vec e\sb k) = <\nabla \sb{\vec e\sb k}  \vec e\sb i, \vec e\sb j> + <\vec e\sb i, \nabla  \sb{\vec e\sb k} \vec e\sb j>$

tal como se había solicitado.

Es conveniente darse cuenta que la expresión $<\vec e\sb i, \vec e\sb j>$ ha de entenderse como una función de la variedad en los reales. Es decir, no es función ni de $\vec e\sb i$ ni de $ \vec e\sb j$, sino del punto donde está anclado el espacio tangente. Por lo que en realidad tenemos una familia de funciones, subindicadas por $i$ y por $j$.

\

\addtocounter{ecu}{1}

\textit{Teorema (\theecu):  Dada una conexión sobre una superficie de Riemann, es decir una derivada covariante, se le puede asociar una familia de 1-formas $ \omega \sp i \sb{j}$ cuya matriz  es $ [\omega  \sb{ij}]$, y eligiendo una base ortonormal para el espacio tangente, dicha matriz  es antisimétrica, $\omega \sb{ij}=\omega \sb{ji}$. Como consecuencia, toda conexión en una superficie de Riemann está completamente determinada por el elemento $w\sb{12}$.}

\bigskip

Demostración:

Recordando que una derivada covariante  toma vectores en el espacio tangente y produce vectores en el mismo espacio y que por lo tanto se pueden descomponer en una base de dicho espacio tenemos:

$\nabla \sb{\vec e\sb k}  \vec e\sb i = \omega \sp r \sb{ki} \vec e\sb r$

$\nabla  \sb{\vec e\sb k} \vec e\sb j = \omega \sp r \sb{kj} \vec e\sb r$

además el tensor métrico se define como $<\vec e\sb i, \vec e\sb j>= g\sb{ij}$ y al substituir en el teorema anterior:

$(d<\vec e\sb i, \vec e\sb j>)(\vec e\sb k) = <\nabla \sb{\vec e\sb k}  \vec e\sb i, \vec e\sb j> + <\vec e\sb i, \nabla  \sb{\vec e\sb k} \vec e\sb j>$

$(dg\sb{ij})(\vec e\sb k) = <\omega \sp r \sb{ki} \vec e\sb r, \vec e\sb j> + <\vec e\sb i, \omega \sp r \sb{kj} \vec e\sb r>= \omega \sp r \sb{ki} g\sb{rj} + \omega \sp r \sb{kj} g\sb{ir}$

Definimos la familia de 1-formas $\omega \sp r \sb i $ mediante la expresión: $\omega \sp r \sb i (\vec e\sb k) = \omega \sp r \sb{ki}$ y  $\omega \sp r \sb j (\vec e\sb k) = \omega \sp r \sb{kj}$. Por tanto, la expresión anterior se reescribe:

$dg\sb{ij} = \omega \sp r \sb{i} g\sb{rj} + \omega \sp r \sb{j} g\sb{ir}$

Ahora definimos: $\omega \sb{ij} =\omega \sp r \sb{ j} g\sb{ir} $ y tenemos:

$dg\sb{ij} = \omega \sb{ij} + \omega \sb{ji}$

Si además la base es ortonormal tenemos $g\sb{ij} = \delta \sb{ij}$. Eso quiere decir que sobre la variedad, cada función $g\sb{ij}$ o bien vale cero o bien vale uno. En ambos casos, la diferencial es cero y por tanto

 $0=\omega \sb{ij} + \omega \sb{ji}$

 de  lo cual se deduce $\omega \sb{ij} =- \omega \sb{ji}$. En particular, la matriz tiene diagonal cero. Y como es antisimétrica, en una matriz 2x2, sólo hace falta saber $w\sb{12}$ para saber toda la matriz.

\addtocounter{ecu}{1}

\textit{Teorema (\theecu):  La \index{curvatura} \textbf{curvatura} sobre una superficie está determinada por una 2-forma $\theta \sb{12}$ la cual coincide con la diferencial de aquella 1-forma que determina la derivada covariante: $\theta \sb{12}=d\omega \sb{12} =-d\omega \sb{21} = -\theta \sb{21}= K dS$, donde $K$ es la curvatura escalar de la ecuación de Jacobi y $dS$ es la 2-forma que mide el elemento de área de la superficie.}

\bigskip

Argumentación: si la curvatura $\Omega$ se extrae de la derivada covariante y toda la derivada covariante está resumida en su componente $\omega \sb{12}$, la cual es una 1-forma, entonces toda la curvatura debe poderse derivar de dicho término.  La única forma posible debe ser entonces $\theta \sb{12}=d\omega \sb{12}$. Por la antisimetría de $\omega$ se deduce que
$\theta \sb{12}=d\omega \sb{12} =-d\omega \sb{21} = -\theta \sb{21}$.

Por otro lado, una 2-forma funciona sobre un 2-plp, el cual tenemos que asumirlo tangente a la superficie. Y en ese caso, la 2-forma funciona como un determinante, es decir como un área. Por lo tanto,  la 2-forma área, $dS$, y la dos forma $\theta \sb{12}$ deben ser múltiplos la una de la otra: $\theta \sb{12}= kdS$.

\bigskip

Para cambiar esta argumentación plausible por un procedimiento riguroso hay que proceder según la siguiente tónica. Primero se calcula $d\omega \sb{12}$. Como es una 2-forma, y las 2-formas son un espacio vectorial cuya base está formada por formas elementales $dx\sp r \wedge dx\sp s$ entonces $d\omega \sb{12}$ debe poderse descomponer en tal base:

$d\omega \sb{12}= (1/2)R\sp 1 \sb{2rs} dx\sp r \wedge dx\sp s$

Al evaluar estas dos 2-formas sobre plps generados por $\partial/\partial x\sp i$ se destila el valor de los coeficientes $R\sp 1 \sb{2rs}$. Se encuentra que coinciden exactamente con los de la curvatura definida como conmutador de derivada covariante.

Toda los resultados anteriores los vamos a utilizar para demostrar el siguiente teorema terminal:

\addtocounter{ecu}{1}

\textit{Teorema (\theecu): Sea $U $ una región compacta en una superficie de Riemann con frontera $\partial U$ formada por una curva suave. $U$ está contenida en un conjunto parametrizado por un solo retazo en el plano, con parámetros $x\sp i$. Por lo tanto, existe un marco coordenado $ [\vec e\sb 1, \vec e\sb 2]$ para $U$ el cual define una orientación para $U$: positivo del primer elemento de la base al segundo. Dicha orientación induce otra sobre  $\partial U$, positivo contrario a las manecillas del reloj. Dicho marco siempre se puede ortogonalizar restando proyecciones, por lo que lo consideramos ortonormal y sobre el ponemos coordenadas polares. Sea $\vec v$ un vector unitario tangente a la superficie, puesto en algún lugar de la frontera. Se transporta paralelamente a $\vec v$ alrededor de $U$ a lo largo de $\partial U$, empezando con $\vec v\sb o$ y terminando con $\vec v\sb f$. Sea  $\Delta \alpha = \vec v\sb f-\vec v\sb i$. Entonces se tiene que}

$\Delta \alpha = \int \int \sb U KdS = \int \int \sb U K dx\sp r \wedge dx\sp s$

\textit{donde $K$ es la curvatura escalar de Gauss, la misma de la ecuación escalar de Jacobi.}

Demostración:

La curva cerrada $\partial U$ puede parametrizarse en longitud de arco. Sea $\vec T$ el vector tangente unitario, sea $\alpha $ el ángulo entre $\vec e\sb 1$ y  $\vec v$, pero orientado desde el primer vector hasta el segundo. A dicho ángulo se le puede asociar la 1-forma $d\alpha = (d\alpha /ds)ds$. De esta manera $\Delta \alpha = \oint  \sb{\partial U} d\alpha $ . Sobre el marco ortogonal descomponemos en coordenadas polares a $\vec v$:

$\vec v = \vec e \sb 1 v\sp 1 + \vec e \sb 2 v\sp 2=  \vec e \sb 1 cos \alpha + \vec e \sb 2 sen \alpha $

Recordemos una fórmula que hallamos anteriormente para la derivada covariante

$$\nabla \sb{\vec X}(\vec v) =[dv\sp i + v\sp k \omega \sp i \sb{jk}\sigma \sp j ] (\vec X) \vec e\sb i$$

y aplicándola a nuestro caso se obtiene:

$\nabla \sb{\vec T} (\vec v)  = (dv\sp 1 + \omega \sb{12} v\sp 2)(\vec T)\vec e \sb 1  +  (dv\sp 2 + \omega \sb{21} v\sp 1)(\vec T) \vec e \sb 2 $

$=  (-sen \alpha d\alpha  + \omega \sb{12} sen \alpha) (\vec T) \vec e \sb 1 + (cos \alpha d  \alpha + \omega \sb{21} cos \alpha)(\vec T) \vec e \sb 2 $

$= (-\vec e \sb 1 sen \alpha + \vec e \sb 2 cos \alpha )(d\alpha - \omega \sb{12})(\vec T)$

en donde utilizamos la antisimetría de $\omega$.

Decir que $\vec v $ es deslizado paralelamente alrededor de $\partial U$ es decir que $\nabla \sb{\vec T} (\vec v) = 0$ . Igualando a cero la última ecuación, y puesto que  el seno y el coseno nunca son cero al mismo tiempo, entonces se concluye que $d\alpha - \omega \sb{12}= 0$. Lo cual significa que $d\alpha = \omega \sb{12}$, es decir, $d\alpha (\vec T)= \omega \sb{12} (\vec T)$:

$\Delta \alpha = \oint  \sb{\partial U} d\alpha = \oint  \sb{\partial U} \omega \sb{12}  $

necesitamos ahora utilizar el teorema de Stokes en 3:d, queda:

$\oint  \sb{\partial U} \omega \sb{12}  = \int \int \sb U d\omega \sb{12}$

pero   de $\omega \sb{12}$ se deduce la curvatura

 $d\omega \sb{12} = \theta \sb{12} = K dx\sp 1 \wedge x \sp 2 $

 tenemos:

$\Delta \alpha = \int \int \sb U d\omega \sb{12} = \int \int \sb U K dx\sp 1 \wedge x \sp 2 $

Esto termina la demostración propuesta. Hagamos una verificación: en un viaje desde el polo norte hasta el ecuador y vuelta después de un octante, el ángulo generado es $\pi /2$. El área de una esfera es $4\pi r\sp 2$. Por lo tanto el área de un octante es $\pi r\sp 2/2$ y la curvatura es  $1/r\sp 2$. Reemplazando:

$\pi /2 =\Delta \alpha =  \int \int \sb U K dx\sp 1 \wedge x \sp 2 $

$ = K \int \int \sb U  dx\sp 1 \wedge x \sp 2 = K$(Area de un octante)$ = (1/r \sp 2) (\pi r\sp 2/2) = \pi /2  $.

\section{CONCLUSION}

Conservando la manía de expresar las leyes de la naturaleza en la forma de ecuaciones diferenciales que nos dice cuál es el próximo paso que ella dará, hemos preparado el terreno para poder hablar de relatividad general. Para ello hemos definido la derivada covariante y hallado la ecuación seguida por una geodésica, definida en términos de transporte paralelo dada una métrica. Definimos   la curvatura  de dos formas, una basados en la derivada covariante y otra basados en la idea de Jacobi de que en un espacio curvo la tasa de separación de las geodésicas está directamente relacionada con la curvatura.

\section{REFERENCIAS}

Son las mismas que las que aparecen al final del capítulo siguiente.

\chapter{RELATIVIDAD Y GEOMETRIA}

El propósito del presente capítulo es el de mostrar que la relatividad general es un lenguaje construido con términos que  ya conocemos:  geodésicas, derivada covariante, y curvatura, densidad de energía -momento transportada por la masa o los campos, densidad del equivalente energético de la masa,  presión y todo lo que cree energía potencial.

La \index{relatividad general} \textbf{relatividad general} ha resultado relativamente muy exitosa, considerando que es difícil tener éxito.  Con todo, en la aproximación geométrica adoptada se tiene una gran esperanza, pues ha sido posible reformular la teoría del campo electromagnético, de la interacción débil y aún de la fuerte en los mismos términos,  aunque por supuesto,  en un ambiente un poco distinto al del espacio-tiempo.

Veremos la idea central de la relatividad general, y dos detalles que nos permitirán pasar de conceptos generales a una visión relativista de la gravitación de Newton.

\section{EL ESQUEMA GENERAL}

El fundamento de todo es que las trayectorias de las partículas en un campo gravitatorio se describen matemáticamente como \index{geodésica} \textbf{geodésicas}, las cuales son trayectorias que obedecen la ecuación:

\addtocounter{ecu}{1}

$$\partial \sp 2  u\sp \gamma/ds \sp 2 + \Gamma\sp{\gamma} \sb{\beta \alpha}(\partial u\sp \alpha /ds)(\partial u\sp \beta/ds)=0 \eqno{(\theecu)}$$

Los \index{coeficientes de Christoffel} \textbf{coeficientes de Christoffel}, $\Gamma\sp{\gamma} \sb{\beta \alpha}$, dependen de la \index{métrica} \textbf{métrica} de la siguiente manera:

\addtocounter{ecu}{1}

$$\Gamma \sp \gamma \sb{\mu \beta }= (1/2)g\sp{\alpha \gamma}(\partial \sb \beta g\sb{\alpha \mu} + \partial \sb \mu g\sb{\beta \alpha } - \partial \sb \alpha g\sb{ \mu \beta})\eqno{(\theecu)}$$

Toda métrica tiene una curvatura descrita por el \index{tensor de Riemann} \textbf{tensor de Riemann} $\Omega$ que se calcula como

\addtocounter{ecu}{1}

$$\Omega (\vec X, \vec Y)\vec v = \nabla \sb{\vec X} (\nabla \sb{\vec Y}\vec v) -\nabla \sb{\vec Y}(\nabla \sb{\vec X} \vec v )- \nabla\sb{[ \vec X,  \vec Y]}\vec v \eqno{(\theecu)}$$

La curvatura afecta la forma de las geodésicas pues en presencia de curvatura las geodésicas se distancias o se juntan: en un plano con la métrica normal, las geodésicas son líneas rectas, que si al comienzo son paralelas, seguirán siendo paralelas. Pero en una esfera con la métrica usual, si tomamos dos geodésicas, o meridianos,  cerca del polo norte, ellas comenzarán a distanciarse muy rápidamente. En general, si modelamos la distancia entre geodésicas por medio de $\vec J$, entonces $J$ evoluciona de acuerdo a:

\addtocounter{ecu}{1}

$$(\frac{\rm \nabla \sp 2 \vec J}{\rm \partial u\sp 2}) + \Omega(\vec J, \vec T) \vec T =0 \eqno{(\theecu)}$$

donde $\vec T$ es el vector tangente unitario a la geodésica de referencia.

Esta ecuación tiene valor científico, pues nos dice cómo hacer un experimento en un laboratorio reducido, estudiando la trayectoria de pares de partículas, viendo qué tan rápido se acercan o se alejan.

\bigskip

Para conectar estas  ideas con el mundo físico hay que fabricar un canal que  permita enlazar la \index{materia} \textbf{materia} con la curvatura.

\bigskip

Dicho canal se relaciona con el tensor de Riemann, el de la curvatura, como sigue:

\addtocounter{ecu}{1}

$$R_{\mu \nu} = R^\alpha_{\mu \alpha \nu} \eqno{(\theecu)}$$

\addtocounter{ecu}{1}

$$R = g^{\mu \nu} R_{\mu \nu}\eqno{(\theecu)}$$

\addtocounter{ecu}{1}

$$G_{\mu \nu } = R_{\mu \nu} - \frac{1}{2}g^{\mu \nu}R\eqno{(\theecu)}$$

Las componentes $G_{\mu \nu }$ corresponden a un tensor $G$ que es el que percibe el efecto de la materia. De qué forma? Lo que de la materia importa para crear curvatura está codificado en una 2-forma, el \index{tensor de tensión-energía} \textbf{tensor de tensión-energía}, $T$. Si uno sabe el tensor de tensión -energía, uno puede saber el efecto sobre la \index{curvatura} \textbf{curvatura} mediante la siguiente relación:

\addtocounter{ecu}{1}

$$G = \kappa T\eqno{(\theecu)}$$

Esa ecuación dice que toda fuente de curvatura sale de la materia y que toda materia crea curvatura.

Las predicciones de la teoría  a escala planetaria son muy buenas, pero no a escalas de clusters de galaxias o mayores. Para tratar de cuadrar, infructuosamente hasta hoy, las observaciones astrofísicas con la teoría, la anterior ecuación se ha modificado en

$G + \Lambda g = \kappa T$

donde $\Lambda $ es una constante, llamada \index{constante cosmológica} \textbf{constante cosmológica} y $g$ es el tensor métrico.

\bigskip

Cómo se manufactura el tensor de tensión-energía?

Se manufactura como una 2-forma de tal manera que contenga información sobre lo viejo, para generalizar la teoría de Newton, y sobre lo nuevo, para que explique los misterios de la vieja teoría.

Lo viejo se refiere a la masa, que en la teoría antigua era la única fuente de gravitación. Por supuesto, después de haber demostrado que toda masa contiene energía, ya no se hablará de masa sino de su componente energético. Una predicción crítica de la teoría consiste precisamente en alegar que la energía de   la luz, que no tiene masa, es influenciada por la gravedad.

Para poner a prueba dicha predicción, se han ideado varios tests. Uno de ellos se hace en el laboratorio y consiste en estudiar un \index{laser} \textbf{laser} cualquiera. Se mide su frecuencia al salir del tubo laser y se mide a un metro de distancia pero de dos maneras distintas. La primera es horizontalmente y la segunda es verticalmente, sea hacia arriba o sea hacia abajo. La predicción relativista alega: la luz pesa, por lo tanto, será atraída por la tierra y si va hacia arriba, pierde energía y su longitud de onda aumentará. Si es hacia abajo, ganará energía y su  longitud de onda disminuirá. Si es horizontalmente, a un metro de distancia no habrá cambios percibibles.

Esta predicción, revestida de los detalles cuantitativos, se considera ampliamente verificada.

\bigskip

Otro arreglo experimental que   conjuga la teoría cuántica de partículas con una visión clásica de la gravitación consiste en enviar un haz de \index{neutrones}  \textbf{neutrones},    por las aristas de una mesa inclinada. Salen de un esquina y se hacen interferir en la esquina opuesta después de recorrer dos caminos en L distintos, una parte por 'arriba' y otro por 'abajo'. Puesto que en la mecánica cuántica, hay un operador de evolución que es la concatenación de evoluciones instantáneas, que también son descritas por operadores y ellos no son conmutativos, entonces podemos predecir que los dos caminos producen efectos distintos, lo cual se evidenciará por el patrón de interferencia. Se sabe y sabe explicar que el patrón de interferencia depende del ángulo de inclinación la mesa. Utilizando la jerga moderna decimos: el espacio ordinario se convierte en un \index{espacio no conmutativo} \textbf{espacio no conmutativo} en presencia de la gravitación.

\bigskip

Bueno, pero también hay que explicar los misterios de la gravitación de Newton. Por supuesto que Mercurio es el misterio indicado.

Qué le pasa a \index{Mercurio} \textbf{Mercurio}? Resulta que su órbita precesa. Qué quiere decir eso? Imaginemos un trompo que baila: al mismo tiempo que el trompo gira sobre sí mismo, su eje se balancea describiendo un círculo. Mercurio es similar a una manchita sobre la superficie del trompo: si proyectamos la órbita de la manchita sobre el piso, veremos que ella describe una  trayectoria que es como una elipse que no se cierra sino que genera una roseta.  Eso es precesión.

La precesión de Mercurio es un misterio en la teoría clásica pues la predicción dice que las órbitas de los planetas   se cierran formando una elipse. Como podemos centrar el alegato con respecto a si la elipse se cierra o no, veamos por qué la teoría clásica espera una órbita cerrada.

La idea de Newton era que no hay diferencia entre el cielo y la tierra, es decir, que las leyes que rigen la caída libre de un cuerpo son aquellas que rigen el movimiento de los planetas.

Para verlo mejor, pensemos en una pelota que uno deja caer contra un piso bien plano y muy firme. A veces uno encuentra una pelota tan elástica que al dejarla caer y rebotar libremente, ella casi llega hasta su punto de partida. Idealmente, en el límite de elasticidad perfecta, ella debe llegar después de rebotar al mismo lugar de donde partió. Ahora bien, un planeta es una pelota cósmica que cae contra el sol y que rebota: esperamos que llegue al mismo lugar de donde partió, o sea que esperamos una órbita cerrada.

Pero por qué el planeta no cae al sol? Porque lleva una velocidad que no apunta hacia el sol, sino que siempre pasa a una distancia prudencial de él. En un caso así, uno puede descomponer la velocidad en una componente que cae hacia el sol y en otra componente perpendicular a la segunda. En ambas componentes tenemos caída libre. El resultado es que la caída y el rebote se convierten en una movimiento orbital que se cierra sobre sí mismo.

Vemos que cerrarse sobre sí misma, es la característica de las órbitas en el movimiento Kepleriano explicadas por una dinámica regida por una fuerza que viene de un potencial.  Además, la pelota que se deja caer, vuela más rápido cerca del suelo. En el caso del planeta, éste viajará más rápido cerca del sol que alejado de él. Es decir, el momento del planeta es mayor cerca del sol que lejos de él.

 Cómo se explica el misterio de la roseta de Mercurio en relatividad? Pues manufacturando un mecanismo que la explique. Cuál podrá ser? Miren que fácil:

Ya hemos visto que la energía pesa, es decir puede ser influenciada por el campo gravitatorio. Pero sabemos que la energía es un concepto sin mucho peso. Lo que si importa e importa mucho es el cuadrivector energía momento que ve la energía y el momento como una unidad, a la cual se le adjudica peso. ¿Explicará eso la precesión de Mercurio? Si: como los planetas van más rápido cerca del sol, ellos pesan más en la cercanía del sol que en su lejanía. Y por lo tanto, el momento adicional hace que la atracción sea más fuerte que si fuera debida sólo a la masa. Eso implica que el planeta Mercurio tiende a permanecer cerca del sol más tiempo de lo esperado. El resultado es que la elipse predicha por Kepler se frustra y el planeta describe una curva que a lo mejor ni se cierra.

Detalle técnico: si las órbitas no se cierran en relatividad, es porque no existe el equivalente de potencial, es decir, el tensor energía-momento no ha de poderse expresar como la derivada exterior de una forma de rango menor. En la jerga se dice: el tensor energía-momento no es una forma exacta.

Ahora lo que hay que tener en cuenta es que la teoría de la relatividad general no es un libro de hermosas historietas, como la que acabamos de ver, sino un cuerpo teórico que produce predicciones cuantitativas. Por eso, hay que pasar necesariamente a la  difícil labor manufacturera del tensor  energía-momento, para después poder calcular dinámicas. Es necesario que en este texto nos contentemos  con dos detalles. El primero nos permite entender cuantitativamente por qué la luz se curva al pasar cerca de una estrella. El segundo nos permite relacionar la teoría de gravitación de Newton con el nuevo lenguaje de las curvaturas.

\section{LA GRAVEDAD CURVA LA LUZ}

El punto de partida de la visión geométrica de la gravitación es el axioma de que una partícula que vuela en caída libre describe una geodésica. El gran trabajo que sigue es especificar la forma de medir distancias para que las geodésicas coincidan con las trayectorias observadas en la naturaleza. Ese trabajo es una labor de manufacturación delicada. Por supuesto que la teoría deberá coincidir en el caso de campos gravitacionales débiles con la teoría clásica o de Newton, la cual supone que el campo gravitatorio es el gradiente de un potencial, función escalar, y cuya divergencia es cero donde no haya masa.

Al mismo tiempo que la teoría debe extender la de Newton, también debe divergir de ella para poder explicar la curvatura de las trayectorias de la luz y la precesión de Mercurio. Y todo eso debe salir naturalmente y de argumentos muy sencillos. La estrategia elegida es meter todo en la estructura del espacio-tiempo mismo, especificando cómo se miden distancias infinitesimales. Comencemos.

Mientras que en un plano y con un sistema de coordenadas rectangulares medimos distancias por medio del teorema de Pitágoras

$ds\sp 2 = dx\sp 2 + dy \sp 2$

en el espacio-tiempo la forma más general posible es demasiado general. Simplifiquemos diciendo que no hay términos mixtos del tiempo con el espacio, lo cual básicamente significa que el tiempo está intrínsecamente de frente contra nosotros. Por otro lado, sabemos que el universo es una entidad dinámica: la constelación de la osa mayor dibujada en un vaso griego puede reconocerse inmediatamente pero la posición relativa de una de las estrellas es diferente de la actual. Con todo, la diferencia no es traumática. Eso quiere decir que a pequeña escala de tiempo, podemos aproximar la realidad variable por una ficción estática.  La forma de medir distancias es entonces:

\addtocounter{ecu}{1}

$$ds\sp 2 = g\sb{oo}(\vec x)dt\sp 2 + g\sb{\alpha \beta}(\vec x)dx\sp \alpha dy \sp \beta \eqno{(\theecu)}$$

donde los índices repetidos se suman, en este caso de 1 a 3, pues el tiempo se distingue aparte con el subíndice cero. Como en ausencia de campos gravitatorios, el espacio-tiempo con gravitación se convierte en nuestro viejo amigo, el espacio -tiempo de la relatividad especial, $g\sb{oo}(\vec x)$ será negativo y por eso estamos en un caso de métrica pseudo-riemanniana.

En el siguiente párrafo está la clave de todo lo que sigue.

Queremos simplificar todo lo posible, pero no hasta lo imposible. Lo imposible es medir gravitación con una sola partícula: puesto que va en caída libre, ella no percibe ninguna fuerza: ella no podrá decirnos que pasa en su entorno. Lo posible es evidenciar la gravitación con, por ejemplo,     dos manzanas que  vuelan en caída libre. Además de caer, ellas se van acercando la una a la otra. Es decir, la distancia entre dos geodésicas varía en presencia de un campo gravitatorio central. Si dicha variación tiene una aceleración no nula, indica una curvatura. Eso implica que una teoría de la gravitación debe contener no sólo geodésicas sino pares de geodésicas. Por tanto, el instrumento de medida no puede ser una 1-forma, pues en ella no hay cabida para efectos colaterales entre geodésicas.

Necesitamos una 2-forma que permita estudiar el resultado sobre un 2-plp (con un vector tomado de de cada una de las dos geodésicas) para que diga que efecto cruzado debe haber. Por  lo tanto, la 2-forma no ha de ser diagonal, aunque los valores diagonales deberán ser predominantes para ser aproximable por la teoría de Newton. El nombre oficial para esa 2-forma es el de tensor métrico, pues ella es la que nos va a decir la forma de medir las distancias entre eventos en presencia de campos gravitatorios y cuya matriz en la base natural está formada por los coeficientes $g\sb{\alpha\beta}$. El tensor métrico opera sobre vectores del espacio tangente en cada punto en particular.

Para poder utilizar la convergencia de la teoría con la de Newton es necesario considerar campos débiles y velocidades muy bajas (comparadas con las de la luz). Eso quiere decir que el espacio-tiempo es casi el producto cartesiano del espacio por el tiempo y que, ante una mirada desatenta, todo luce común y corriente, en particular, deberíamos encontrar que no hay diferencia notable entre tiempo propio y tiempo del observador. Para poder parametrizar en longitud de arco, el tiempo propio, marcado según un reloj que la partícula lleva, lo medimos como:

$d\tau \sp 2 = -ds\sp 2 = -g\sb{oo}(\vec x)dt\sp 2 - g\sb{\alpha \beta}(\vec x)dx\sp \alpha dy \sp \beta$

Dividiendo esta ecuación por $dt\sp 2$ :

$(d\tau/dt) \sp 2 = -g\sb{oo} - g\sb{\alpha \beta}(dx\sp \alpha/dt) (dy \sp \beta/dt)$

Como la velocidad es muy baja comparada con la de la luz, nosotros caricaturizamos la situación diciendo que $\vec v = d\vec x/dt \approx \vec 0$. Sin embargo, a bajas velocidades puede haber aceleraciones casi infinitas, que es  lo que nos interesa, y por eso nuestras aproximaciones pueden producir algo interesante.  En fin, reemplazando la velocidad por cero, tenemos:

$(d\tau/dt) \sp 2 = -g\sb{oo}$

$d\tau/dt = (-g\sb{oo})\sp{1/2}$

por tanto, el cuadri-vector velocidad es:

$u=dx/d\tau = (dt/d\tau)[1,d\vec x/dt] \approx (dt/d\tau)[1,0] =  (-g\sb{oo})\sp{1/2}[1,\vec 0] $

Recordando que estamos prácticamente en el espacio -tiempo vacío, y poniendo la velocidad de la luz como 1, tenemos que $g\sb{oo}\approx -1$. Por tanto, $u\approx [1,\vec 0]$ quedando demostrado que necesariamente $d\tau \approx dt \approx ds $. Con todas esas simplificaciones, ya nos queda más fácil estudiar las geodésicas, las cuales deben cumplir:

$\partial \sp 2  u\sp \gamma/ds \sp 2 + \Gamma\sp{\gamma} \sb{\beta \alpha}(\partial u\sp \alpha /ds)(\partial u\sp \beta/ds)=0$

que cambiando a nuestras aproximaciones queda:

$\partial \sp 2  x\sp i/dt \sp 2 \approx - \Gamma\sp{i} \sb{jk}(\partial x\sp j /dt)(\partial x\sp k/dt) \approx -\Gamma \sp i \sb{oo}$

(el cero cero se debe a que ambas derivadas son con respecto al tiempo, que es el parámetro número cero).

Recordemos la definición de los $\Gamma$ :

$\Gamma \sp \tau \sb{\mu \beta }= (1/2)g\sp{\alpha \tau}(\partial \sb \beta g\sb{\alpha \mu} + \partial \sb \mu g\sb{\beta \alpha } - \partial \sb \alpha g\sb{ \mu \beta})$

Por tanto

$\Gamma \sp i \sb{oo }= (1/2)g\sp{\alpha i}(\partial \sb o g\sb{\alpha o} + \partial \sb o g\sb{o \alpha } - \partial \sb \alpha g\sb{oo})$

pero como nuestro tensor es casi diagonal, lo que no sea diagonal es casi cero:

$\Gamma \sp i \sb{oo }= (1/2)g\sp{\alpha i}( - \partial \sb \alpha g\sb{oo})$

y substituyendo

$\partial \sp 2  x\sp i/dt \sp 2 \approx -(1/2)g\sp{\alpha i}( - \partial \sb \alpha g\sb{oo})=(1/2)g\sp{\alpha i}(  \partial \sb \alpha g\sb{oo})$

Teniendo en cuenta que $g\sb{\alpha \beta}$ es una 2-forma, ¿cómo es posible que al lado izquierdo tengamos vectores? Eso es posible porque  $g\sp{\alpha i}$ se encarga de la conversión. En términos de vectores, tenemos  que, para cada componente, el lado derecho es una suma de derivadas parciales con respecto a las coordenadas espaciales, o sea, el gradiente de $g\sb{oo}$:

$\partial \sp 2  x\sp i/dt \sp 2 \approx [(1/2)\nabla  g\sb{oo}]\sp i$

o sea que al tener en cuenta la aproximación dada por el potencial gravitatorio newtoniano:

$\partial \sp 2  \vec x /dt \sp 2 \approx (1/2)\nabla  g\sb{oo}] = - \nabla \phi$

de donde deducimos que $g\sb{oo}/2 \approx \phi + c $

la constante $c$ es constante con respecto a las coordenadas espaciales, por lo que podría ser una función del tiempo. Esa posibilidad la obviamos en nuestra primera aproximación. Observemos que toda la teoría antigua cabe en el coeficiente cero cero y que por ende, al transcribir dicha teoría a relatividad, los demás coeficientes son cero.

Para averiguar $c$, asumimos que en infinito no hay gravitación y recobramos el espacio -tiempo vacío: $\phi(\infty)=0$ y  $g\sb{oo}(\infty)=-1 $  y por tanto $-1/2 = 0 +c$, es decir $g\sb{oo} \approx 2(\phi -1/2)$. En conclusión:

$g\sb{oo}=2\phi -1 $

Habíamos supuesto que $g\sb{oo}\approx -1$ en el vacío,  pero ahora encontramos que en presencia de un campo gravitatorio débil, $g\sb{oo}\approx 2\phi -1 $. Ahora bien, la métrica en el espacio -tiempo vacío multiplica la coordenada del tiempo por $c\sp 2$, pero en presencia de una débil gravedad, nosotros  la multiplicamos por $g\sb{oo}$. Sabiendo que la velocidad de la luz depende del medio donde esté, nosotros no vamos a inventar nada nuevo, sino simplemente a decir que la gravitación modifica el medio de tal manera que la luz deja de viajar de velocidad $c=1$ a velocidad

$c=\sqrt{-g\sb{oo}}\approx \sqrt{1-2\phi} \approx (1-\phi)$

Lo cual releemos: el potencial gravitatorio funciona como un \index{índice de refracción} \textbf{índice de refracción} que disminuye la velocidad de la luz. Quizá convenga aclarar que una geodésica es la trayectoria seguida por una partícula que no cambia la estructura del espacio tiempo a su alrededor. Por supuesto que tal cosa no existe. Pero  existen elementos que causan muy poco desorden. La luz puede ser uno de ellos. Por eso, no es grave error imaginar que un rayo de luz débil siga una geodésica. Ahora bien, cuando decimos que la gravedad disminuye la velocidad de la luz, estamos diciendo que la  luz sigue una trayectoria curva, pues para ir de un lugar a otro, a la luz le rinde más alejándose de la fuente de gravedad así le toque dar la vuelta. Este giro en el lenguaje se debe a que una geodésica también es un extremo de un funcional que mide el tiempo propio a lo largo de las trayectorias. Pasemos ahora a una reformulación de la gravitación de Newton desde el punto de vista de Jacobi.

\section{NEWTON A LA JACOBI}

El espacio-tiempo que estamos considerando está perturbado por un débil potencial gravitatorio y todas las velocidades, comparadas con la de la luz, son casi nulas. Eso significa que con respecto a primeras derivadas no hay diferencia con el espacio -tiempo. La diferencia será visible en las segundas derivadas y la consiguiente curvatura.

En un débil potencial, la ley de gravitación clásica es aceptable. A eso debe sacársele alguna información. Consideremos, pues, dos partículas que caen libremente en el espacio. Ellas describen geodésicas en el espacio-tiempo, pero en el espacio 3:d pueden describir parábolas o elipses, o incluso espirales complicadas, por ejemplo en el caso de un asteroide girando alrededor de otro más grande mientras que éste gira alrededor del sol.

Ya habíamos definido dos formas de comparar curvas que son solución de ecuaciones diferenciales: comparaciones dentro de una misma solución, medida por $d/dt$ y comparaciones entre geodésicas, dada por $\delta $. Además, estas dos operaciones conmutan. Consideremos
pues las dos trayectorias descritas por las dos partículas, las cuales cada una obedece la ley de Newton: $d\sp 2 \vec x/dt\sp2= -\nabla V = \nabla \phi$. Si ahora comparamos dos trayectorias por medio de $\delta$ obtenemos:

$\delta (d\sp 2 \vec x/dt\sp2) = \delta (\nabla \phi)$

tomando coordenada por coordenada, en el espacio,

$\delta (d\sp 2  x \sp{\alpha} /dt\sp2) = \delta (\partial \phi/\partial x\sp{\alpha})$

o bien

$\delta ((d\sp 2   /dt\sp2)(x \sp{\alpha})) = \delta (\partial \phi/\partial x\sp{\alpha})$

conmutando en el lado izquierdo

$(d\sp 2/dt\sp2)(\delta   x \sp{\alpha} ) = \delta (\nabla \phi) = \delta (\partial \phi/\partial x \sp{\alpha})$

elaborando el lado derecho con la regla de la cadena:

$(d\sp 2/dt\sp2)(\delta   x \sp{\alpha} ) = (\partial \sp 2 \phi/\partial x \sp{\alpha} \partial x \sp{\beta})\delta x\sp{\beta}$

esta es una ecuación que nos dice como cambia la distancia en el espacio entre dos partículas que caen libremente. Ahora, veamos la versión del mismo fenómeno pero en el espacio tiempo. Cada partícula vuela describiendo una geodésica, la cual parametrizamos con el tiempo propio $\tau$ y cuyo vector tangente es la misma cuadri-velocidad unitaria que ya averiguamos anteriormente $u=[1, \vec 0]$ pues $\gamma =1$. Como las velocidades consideradas son muy bajas, el tiempo propio puede substituirse por un tiempo absoluto $t$, medido por el observador. Por lo tanto, el vector $\vec J$ que traza la componente lineal de la distancia entre geodésicas no tiene componente temporal, pues todas las partículas tienen el mismo tiempo absoluto. Pero la ecuación de Jacobi sigue siendo válida: $\nabla \sp 2 \vec J/dt\sp 2 = -R(\vec J,\vec T) \vec T. $

Como el espacio-tiempo considerado se diferencia levemente del espacio-tiempo plano, tenemos:

$\nabla J\sp{\alpha}/dt = d J\sp{\alpha}/d \tau =d J\sp{\alpha}/dt $

$\nabla \sp 2 J\sp{\alpha}/dt\sp 2 = d \sp 2 J\sp{\alpha}/d \tau  \sp 2 =d \sp 2 J\sp{\alpha}/dt \sp 2 $

donde  sólo las coordenadas espaciales pueden ser no nulas, pues $J\sp o = 0$ .

Ahora reescribamos la ecuación de Jacobi en coordenadas. Para ello, a lo largo de la geodésica ponemos una base móvil, de elementos mutuamente perpendiculares formando una base en el espacio-tiempo de tal forma que el vector velocidad espacial es siempre uno de  los vectores de la base espacio-temporal. Esta base móvil es una base para el espacio tangente en cada punto de la geodésica, por lo tanto todos los vectores de dicho espacio pueden descomponerse en ella, en particular el vector tangente $\vec T  $,  y el vector de Jacobi,  $\vec J$. Tal base se numera $(\vec e\sb{\alpha})$. Tomemos la ecuación de Jacobi $\nabla \sp 2 \vec J/dt\sp 2 =d\sp 2 \vec J/dt\sp 2 = -R(\vec J,\vec T) \vec T$ y multipliquemos por $\vec  e\sb{\alpha}$:

$<d \sp 2 \vec J/dt\sp 2, \vec  e\sb{\alpha}> = (d \sp 2/dt\sp 2)< \vec J, \vec  e\sb{\alpha}> = -<R(\vec J,\vec T) \vec T, \vec  e\sb{\alpha}>$

$(d \sp 2/dt\sp 2)< \vec J,\vec  e\sb{\alpha}> = -<R(\vec J,\vec T) \vec T,\vec  e\sb{\alpha}> $

Atención: el operador segunda derivada  parece que saliera por linealidad o algo así. En realidad, eso es un espejismo: lo que pasa es que la derivada de un producto tiene dos términos, pero el faltante es cero, pues la derivada covariante de un marco ortogonal  con uno de sus componentes a lo largo de una geodésica siempre es cero, pues   conserva la norma y los ángulos con la geodésica.

Puesto que $< \vec J,\vec  e\sb{\alpha}> = J^\alpha$ y expandiendo en la base:

$d \sp 2 J \sp{\alpha}/dt\sp 2  = - <R( J\sp{\beta} \vec  e\sb{\beta}, T\sp{\gamma} \vec e\sb{\gamma})  T\sp{\gamma} \vec e\sb{\gamma},e\sb{\alpha}>. $

Utilizando la linealidad de $R$ en todas sus entradas:

$d \sp 2 J \sp{\alpha}/dt\sp 2  = - J\sp{\beta}T\sp{\gamma} T\sp{\gamma}<R(  \vec  e\sb{\beta},  \vec e\sb{\gamma})  \vec e\sb{\gamma},e\sb{\alpha}>. $

Teniendo en cuenta que sólo $T\sp 0 =1$ y que todos los demás son cero:

$d \sp 2 J \sp{\alpha}/dt\sp 2  = - J\sp{\beta}<R(  \vec  e\sb{\beta},  \vec e\sb{0})  \vec e\sb{0},e\sb{\alpha}>. $

$d \sp 2 J \sp{\alpha}/dt\sp 2  = - J\sp{\beta}R \sb{0 \beta 0} \sp \alpha $

 Puesto que consideramos potenciales muy débiles, $J\sp{\alpha} = \delta x \sp{\alpha}$ y $J\sp{\beta} = \delta x \sp{\beta}$. Substituyendo:

$d \sp 2 \delta x\sp{\alpha}/dt\sp 2  = - R \sb{0 \beta 0} \sp \alpha\delta x \sp{\beta} $

En la primera coordenada se tiene $0=0$ pues $J\sp 0=0$ como resultado de una sincronización del experimento. Comparando las otras tres coordenadas de esta ecuación con las de  la ecuación en 3:d que encontramos anteriormente

$(d\sp 2/dt\sp2)(\delta   x \sp{\alpha} ) = (\partial \sp 2 \phi/\partial x \sp{\alpha} \partial x \sp{\beta})\delta x\sp{\beta}$

concluimos que

$- R \sb{0 \beta 0} \sp \alpha = \partial \sp 2 \phi/\partial x \sp{\alpha} \partial x \sp{\beta}$

 pero puesto que $\phi$ es una función escalar, cuyo gradiente en 3:d es el campo de fuerzas y cuya divergencia es cero a menos que haya una fuente de campo, o sea, una masa:

$\nabla \sp 2 \phi = \sum \partial \sp 2 \phi/\partial x \sp{\alpha} \partial x \sp{\alpha}= -4\pi \kappa \rho = - R \sb{0 \alpha 0} \sp \alpha $,

índices repetidos, uno arriba y otro abajo, es una suma.

Y, por ahora, es todo lo que podemos decir sobre la reescritura de la gravitación de Newton en el nuevo lenguaje de la relatividad general. No ha sido mucho, pero nos ha dado una idea de cómo se manufactura una teoría que hasta el día de hoy causa la admiración de todo el que la estudia.

\section{CONCLUSION}

Hemos visto de qué manera se geometriza la gravitación: hemos enriquecido el espacio-tiempo con una métrica cuyas implicaciones físicas deben coincidir con las del espacio-tiempo cuando no hay masa y que debe extender la teoría clásica sobre la gravitación cuando la función potencial es débil.  Logramos probar que en presencia de materia, el tensor de curvatura no es nulo y que la velocidad de la luz disminuye, por lo que la materia funciona como una entidad refractaria. Eso implica que las geodésicas, que son curvas de distancia mínima local, se acercan unas a otras en presencia de materia. Todo esto nos permite comenzar a entender por qué la luz de las estrellas se curva pasando cerca del sol. Avanzamos lo suficiente como para entender que existen dos planteamientos equivalentes para geometrizar la gravitación: el de la derivada covariante o conexión, que es la derivada de interés geodésico, y el transporte paralelo, que está implícito en la medición de los ángulos internos de un triángulo para determinar la curvatura del espacio donde el triángulo está inmerso.

\section{Remodelación}

Ningún genio es indispesable para la ciencia, pues a la postre todo ha de ser necesariamente descubierto por un proceso  natural de maduración  que consiste simplemente en hacer que los muy diversos enfoques que se tienen formen un todo coherente. Sin embargo, a Einstein se le considera uno de los grandes genios de la humanidad por habernos enseñado, algunos años antes de su tiempo normal de maduración, algo tan extrañamente alocado y hermoso como es que el espacio-tiempo y la gravedad forman un todo inseparable. Y todo eso, el lo hizo artesanalmente, es decir a pura intuición y aplicando una cierta generalización del análisis dimensional combinado con tratamientos de casos especiales y asintóticos. 

Las nuevas generaciones de matemáticos y físicos han tratado de superar el elemento artesanal de la relatividad general y de imponer una manera elegante y limpia de reinventarla e incluso de proponer nuevas alternativas. Lo que eso quiere decir, hoy en día,  es que una teoría sobre las interacciones fundamentales debe nacer naturalmente de la teoría de grupos.  

\bigskip

Los siguientes son ejemplos de grupos que han sido relacionados con la gravitación:

\begin{itemize}
	\item El grupo GL(n, R) de matrices invertibles $n\times n$ con entradas reales   (Mansouri and Chang, 1976).
	
	\item El grupo de las transformaciones de Lorentz extraído de la relatividad especial, el cual corresponde al conjunto de todas las matrices que conservan el intervalo relativista (Utiyama, 1956; Chamseddine, 2005).

 \item El grupo de Poincaré que es el grupo de Lorentz mas las traslaciones (A. López.Pinto, A. Tiemblo  and R. Tresguerres, 1996).

\item Es intrigante una severa crítica en contra de la gravitación  a la Einstein escrita por Pommaret (1987), quien además propone su propia teoría utilizando una estructura muy poderosa: los pseudo-grupos que son grupos de transformaciones definidos sólo localmente.

\end{itemize}
   
La gravitación ha resultado ser la más desafiante de todas las interacciones. Por eso, no hay mucha prisa para entender los trabajos que la estudian, sino que primero hay que entender lo más fácil. Nuestro compromiso es que la interacción más sencilla, la electromagnética, quede bien entendida a un nivel que le permita a uno aventurarse con las otras interacciones y con la gravitación.  Nosotros ligaremos el electromagnetismo con el grupo $U(1)$, el grupo de los números complejos con módulo UNO. Esto nos abre la puerta al estudio de la interacción débil, con grupo $SU(2)$, el grupo de las matrices con entradas complejas $2\times 2$ con determinante +1 y unitarias (una matriz es unitaria cuando ella por su transpuesta conjugada da la identidad). La siguiente generalización es la interacción fuerte con grupo $SU(3)$. Siguiendo por esta línea se llega directamente a una teoría de la gravitación basada en $U(2) \times U(2)$, donde el grupo $U(2)$ consta de  matrices unitarias pero sin restricción en cuanto al signo del determinante.

\section{REFERENCIAS}

\begin{enumerate}

\item N.N \textit{An introduction to differential geometry}

\item Ali H. Chamseddine (2005) Applications of the Gauge Principle to Gravitational Interactions. arXiv:hep-th/0511074v1 7 Nov 2005

\item Misner, C., Thorne, K., and Wheeler, J. \textit{Gravitation, Freeman}, 1970.

\item  Frankel T. \textit{The Geometry of Physics, An Introduction}. Cambridge University Press, 2001.

\item Kibble T. Lorentz invariance and the gravitational field. \textit{Journal of Mathematical physics} 2, 212 (1961).

\item A. López.Pinto, A. Tiemblo(.) and R. (1996) Hamiltonian Poincaré gauge theory of gravitation. arXiv: gr-qc/9603023v1.
 
\item Mansouri F, Chang L.N (1976) Gravitation as a gauge theory. \textit{Phys. Rev. D} 13, 3192 - 3200 (1976)
 
\item Pommaret J.F. (1987) \textit{Lie Pseudogroups and Mechanics}. Gordon and Breach. NY, London, Paris. 

\item Utiyama R. (1956), Invariant theoretical interpretation of interaction. \textit{Physical Review }101, 1597.

\end{enumerate}

\chapter{LA CONEXION ELECTROMAGNETICA}

\Large

\centerline{RESUMEN}

\bigskip

Se discute el problema de la geometrización del campo electromagnético. Es decir, se busca demostrar que una partícula cargada en un campo electromagnético  traza una geodésica, o lo que es equivalente, que su dinámica está descrita por una ecuación diferencial que involucra una   \textbf{derivada covariante}. El espacio de operación no es 3:d ni 3+1:d, sino una amalgama entre el espacio-tiempo y el espacio de la fase de la función de onda.   

\normalsize

\section{INTRODUCCION}

La teoría clásica de la gravitación predecía que las orbitas de los planetas debían ser elipses, algo a lo que la órbita de Mercurio no se ajustaba. Esa era una razón, entre muchas, para fabricar una nueva teoría. En respuesta a ese desafío, la teoría geométrica de la gravitación, o \index{relatividad general}  \textbf{relatividad general}, ha llegado a ser el ejemplo de una teoría exitosa: detrás de los sofisticados métodos matemáticos  se esconde una teoría ideológicamente sencilla, coincide con la teoría clásica cuando ésta es correcta, explica los misterios de la vieja teoría y predice efectos que no cabe imaginar  en el antiguo contexto, como que la luz se curva al pasar por una estrella. Y además, tiene sus propios misterios, algo que la hace más encantadora, como el significado y valor de la constante cosmológica.

La teoría geométrica tiene varios elementos. Algunos son:

1) Todo lo que contenga masa, energía o momento interactúa por medio de la gravitación. Al incluir la energía y el momento, además de la masa, podemos  predecir que la luz de las estrellas se curva al pasar cerca del sol. Pero como hay infinitas maneras de curvarse, tenemos como desafío predecir la manera correcta como eso sucede.

2) La interacción gravitatoria no es a distancia a través del vacío, sino que está mediada por el espacio y es local. Es decir, la dinámica inmediata de una partícula test  que no alcanza a perturbar las predicciones y que está bajo la acción de un campo gravitatorio debe poder predecirse correctamente a partir de la  vecindad ocupada por la partícula test. Por supuesto, aprovechando la experiencia de la relatividad especial, el espacio no es el 3:d habitual, sino el espacio-tiempo que forma una unidad indivisible evidenciada en la métrica del espacio-tiempo, que ahora es una pseudo-métrica, la cual será perturbada por las entidades gravitacionales. La dinámica se deduce de asumir que la partícula test sigue una geodésica, es decir, sigue una trayectoria que es un mínimo local (mínimo absoluto en una pequeña vecindad) de la distancia dada por  la  pseudo-métrica.

3) Puesto que la acción es local, la dinámica debe poderse describir por una ecuación diferencial, cuyo único objetivo es decir precisamente que (localmente) se minimiza la distancia a medida que se avanza. Es decir, la solución a tal ecuación diferencial, con condiciones iniciales dadas, da una geodésica. Pero puesto que todo está en la métrica, y se debe poder reproducir la teoría clásica, entonces en las geodésicas debe estar escrito de algún modo   que en presencia de campo gravitatorio las trayectorias se curvan, es decir se aceleran.

4) Pero la curvatura  debe poderse evidenciar localmente, por lo que debe existir un mecanismo natural que ligue curvatura con geodésicas. Precisamente, se inventó un lenguaje que permite decir: una aceleración no nula de la distancia entre geodésicas (entre dos partículas que caen libremente)  evidencia curvatura. Y vice-versa. El nuevo lenguaje describe la distancia entre geodésicas  por una ecuación de segundo orden, la ecuación de Jacobi, pero las derivadas no son derivadas ordinarias, sino que son derivadas  covariantes.

 \bigskip 

Gracias a una estupenda maquinaria matemática, los anteriores conceptos permiten demostrar que en presencia de campo gravitatorio, las geodésicas se acercan con una aceleración covariante no nula . Un completo desarrollo de la teoría geométrica de la gravitación la convierte en una impresionante escultura en honor del buen gusto de la ciencia moderna y de la sencillez y gracia de la naturaleza.
\bigskip
Con justa razón debemos preguntarnos si es posible extender la teoría geométrica de la gravitación, o relatividad general, para que incluya las otras interacciones.  Sin embargo, y aprendiendo la lección de la historia, esto debe hacerse naturalmente a partir de la teoría de grupos. Comenzamos tal discusión con la interacción electromagnética.

\textbf{Geometrizar} \index{geometrizar} la interacción electromagnética no puede significar armar rancho aparte. Debe extender la geometrización de la interacción gravitatoria. Pero cómo?

Puede ser de ayuda el tratar de explicitar el cambio de perspectiva que estamos reclamando. Cuando uno se imagina que hay partículas que interactúan a distancia, uno podría pensar que las partículas son como naves espaciales que saben y deciden a donde van. En contraste, cuando uno piensa en geometrizar, uno se imagina un  río que se hace cargo de todo lo que arrastra. En relatividad, el río es el espacio tiempo, el cual se recarga de gravitación en presencia de materia.
\bigskip

Qué haremos ahora para \index{geometrizar el   electromagnetismo} \textbf{geometrizar el campo electromagnético}?

\bigskip
 Puesto  que la gravedad y el electromagnetismo son diferentes (las constantes de acople y las cargas son distintas), podemos decir que cada una de ellas vive en su propio espacio. El resultado es que tenemos que  reconocer otros grados de libertad además de los del espacio-tiempo, es decir, hay que extender el espacio-tiempo  y formular una métrica para el nuevo espacio-extendido  que extienda la métrica gravitatoria y que prediga la dinámica de una partícula masiva y cargada en presencia de un campo gravitatorio y otro electromagnético.

Este programa suena complicado pero hay una esperanza debida a la sencillez de la interacción electromagnética.  Veamos lo que eso significa.

La 1-forma   del potencial electromagnético $A$ produce, al derivarla exteriormente, la 2-forma tensor de campo, $F=dA$, la cual rige la dinámica de una partícula cargada según aparece en la ecuación que describe la fuerza de Lorentz: una partícula cargada en un campo electromagnético modifica su momento, es decir su velocidad o su dirección. Ambas cosas, cambio de velocidad escalar o de su dirección, se llaman curvatura en geometría. Todo eso quiere decir sólo una cosa: la 2-forma $F$ debe poder predecir la forma asintótica de la  curvatura de la métrica extendida  al ir apagando el  campo gravitatorio.

Tendríamos dos campos con dos métricas sobre un espacio-tiempo extendido y perturbado que describen dos conductas asintóticas, y hay que buscar una sola métrica que en casos especiales reproduzca las dos métricas dadas. Toda la experiencia indica que, a las escalas de energía al alcance humano, uno puede considerar  los dos campos por separado y sumar los efectos parciales para obtener el efecto total. Veamos esto un poco mejor.

 La forma como sabemos que hay gravitación es por el estudio de las geodésicas. Pero ellas generan su    propio sistema de coordenadas: ¿Cómo se ve el campo electromagnético desde una geodésica? Respuesta: se ve como campo electromagnético. Por lo tanto, el electromagnetismo tiene vida propia y parece ser distinto del campo gravitatorio.

  Pero esto es un resultado clásico. No hay por qué pensar que eso siga siendo válido para grandes energías, cuando se requiera de descripciones cuánticas. Precisamente, el ideal del campo unificado solicita que a grandes energías los dos campos se mezclen.

 Sin embargo, si consideramos que la geometría unificará todo, estamos prediciendo que la  gravitación es primordial, y que crea el espacio tiempo en el cual las demás interacciones pueden expresarse a su gusto. Esto equivale a decir que las tales métricas funcionan como operadores que actúan unos sobre los efectos de otros.

\section{IMPLEMENTACION}

Puesto que hemos convenido en que los distintos campos funcionan como operadores, unos modificando el resultado de otros, eso quiere decir que el mundo donde han de vivir dichas interacciones en su ambiente geométrico ha de ser aquel donde el lenguaje de operadores sea el natural. Ese es precisamente el mundo de la mecánica cuántica. Nuestro programa ahora consiste en asociar una derivada covariante a la interacción electromagnética, primero apagando la interacción gravitatoria y luego teniéndola en cuenta. De esa forma estaremos diciendo que la dinámica se desarrolla siguiendo geodésicas, las cuales presuponen una métrica y un principio de minimización de distancias estrictamente local.

\subsection{Sin gravitación}

Al geometrizar la gravedad construimos un lenguaje cuyo corazón era la \textbf{derivada covariante}. Recordemos la fórmula que asociamos a dicha derivada:

\addtocounter{ecu}{1}

$$\nabla \sb{\vec X}(\vec v) =[dv\sp i + \omega \sp i \sb{jk}\sigma \sp j v\sp k] (\vec X) \vec e\sb i \eqno{(\theecu)}$$

donde la  1-forma $\sigma \sp j$  denota $dx\sp j$ y $\omega \sp i \sb{jk}$ es la expresión coordenada de la derivada covariante.  Veamos ahora a dónde queremos llegar con el electromagnetismo en relación con derivadas covariantes.

Observemos que la definición de derivada covariante  también se le puede dar sentido para campos que toman valores complejos. Por tanto, si en vez del vector $\vec v$ ponemos una función de onda $\psi$ que a cada punto del espacio-tiempo le asocia un escalar complejo,  entonces $dv$ se convierte en  $d\psi= (\partial \psi/ \partial x\sp j) dx\sp j$. La 1-forma $dx\sp j$ notada también  $\sigma \sp j$  proyecta sobre la coordenada $j$. Para un función de onda general,    la derivada covariante que buscamos se lee

$$\nabla \sb{\vec X}(  \psi) =[(\partial \psi^i/ \partial x\sp j) dx\sp j + \omega \sp i \sb{jk}\sigma \sp j \psi\sp k] (\vec X) \vec e\sb i \eqno{(\theecu)}$$

o bien

$$\nabla \sb{\vec X}(  \psi^i) =[(\partial \psi^i/ \partial x\sp j) dx\sp j + \omega \sp i \sb{jk}\sigma \sp j \psi\sp k] (\vec X)   \eqno{(\theecu)}$$

Cuando la función de onda tiene una única componente, la derivada covariante se lee:

$$\nabla \sb{\vec X}(\psi) =[(\partial \psi/ \partial x\sp j) dx\sp j+ \omega  \sb{j} \sigma \sp j \psi] (\vec X)  $$

Pero $\sigma_j $ y $dx_j$ son lo mismo y su efecto sobre $\vec X$ es el de tomar la proyección correspondiente. Por lo que también podemos decir:

$$\nabla \sb{\vec X}(\psi) =[(\partial \psi/ \partial x\sp j)  + \omega  \sb{j} \psi] x^j  $$

Vemos que  la expresión general $dv\sp i + \omega \sp i \sb{jk}\sigma \sp j v\sp k$ se especializa en

\addtocounter{ecu}{1}

$$\partial \psi /\partial x\sp j + \omega \sb j \psi \eqno{(\theecu)}
$$

Acá es a donde queremos llegar y estemos atentos a reconocer dicha  expresión en nuestro siguiente tratamiento del electromagnetismo.

\

La ecuación genérica de la mecánica cuántica es de la forma

\addtocounter{ecu}{1}

$$i\hbar \frac{\rm \partial \psi }{\rm \partial t} = H\psi \eqno{(\theecu)}$$

donde $\psi$ es la función de onda, que toma valores complejos, y $H$ es el operador de energía, el mismo que hace evolucionar el mundo, es decir, a $\psi$. Por lo tanto, $H$ tendría una parte correspondiente a la energía cinética y otra a la potencial.

Una partícula que tiene momento y se desplaza en un potencial $V$ está descrita por la ecuación:

$i\hbar \frac{\rm \partial \psi }{\rm \partial t} = -(\frac{\rm \hbar \sp 2}{\rm 2m})\sum \frac{\rm \partial \sp 2 \psi}{\rm (\partial x\sp \alpha)\sp 2} + V\psi $

si la partícula está acoplada al campo electromagnético a través de una carga $e$ y el campo electromagnético tiene su 1-forma $A$ en 3:d y su potencial escalar $\phi$, la ecuación que   describe la dinámica está dada por:

$i\hbar \frac{\rm \partial \psi }{\rm \partial t} = (\frac{\rm 1}{\rm 2m})\sum [-i\hbar \frac{\rm \partial }{\rm \partial x\sp \alpha} -eA\sb{\alpha}]\sp 2 \psi  + V\psi - e\phi\psi $

Esta ecuación puede maquillarse como

$i\hbar \frac{\rm \partial \psi }{\rm \partial t} = (\frac{\rm 1}{\rm 2m})\sum [-(i\hbar \frac{\rm \partial }{\rm \partial x\sp \alpha} +eA\sb{\alpha})]\sp 2 \psi  + V\psi +( i\sp 2)(\frac{\rm \hbar}{\rm \hbar}) e\phi\psi $

$i\hbar \frac{\rm \partial \psi }{\rm \partial t} = (\frac{\rm 1}{\rm 2m})\sum [i\hbar \frac{\rm \partial }{\rm \partial x\sp \alpha} +eA\sb{\alpha}]\sp 2 \psi  + V\psi +( i\sp 2)(\frac{\rm \hbar}{\rm \hbar}) e\phi\psi $

$i\hbar \frac{\rm \partial \psi }{\rm \partial t} = (\frac{\rm 1}{\rm 2m})\sum [i\hbar \frac{\rm \partial }{\rm \partial x\sp \alpha} -(i\sp 2)e(\frac{\rm \hbar}{\rm \hbar})A\sb{\alpha}]\sp 2 \psi  + V\psi +( i\sp 2)(\frac{\rm \hbar}{\rm \hbar}) e\phi\psi $

$i\hbar \frac{\rm \partial \psi }{\rm \partial t} -( i\sp 2)(\frac{\rm \hbar}{\rm \hbar}) e\phi\psi = (\frac{\rm 1}{\rm 2m})\sum [i\hbar \frac{\rm \partial }{\rm \partial x\sp \alpha} -(i\sp 2)e(\frac{\rm \hbar}{\rm \hbar})A\sb{\alpha}]\sp 2 \psi  + V\psi$

$i\hbar [\frac{\rm \partial }{\rm \partial t} -(\frac{\rm ie}{\rm \hbar}) \phi]\psi = -(\frac{\rm \hbar \sp 2}{\rm 2m})\sum [ \frac{\rm \partial }{\rm \partial x\sp \alpha} -(\frac{\rm ie}{\rm \hbar})A\sb{\alpha}]\sp 2 \psi  + V\psi $

Puesto que estamos trabajando con operadores, podemos reescribir explícitamente.

$i\hbar [\frac{\rm \partial }{\rm \partial t} -(\frac{\rm ie}{\rm \hbar}) \phi]\psi = -(\frac{\rm \hbar \sp 2}{\rm 2m})\sum [ \frac{\rm \partial }{\rm \partial x\sp \alpha} -(\frac{\rm ie}{\rm \hbar})A\sb{\alpha}][ \frac{\rm \partial }{\rm \partial x\sp \alpha} -(\frac{\rm ie}{\rm \hbar})A\sb{\alpha}] \psi  + V\psi $

Comparemos los operadores entre corchetes, $[\frac{\rm \partial }{\rm \partial t} -(\frac{\rm ie}{\rm \hbar}) \phi]$ y
$[ \frac{\rm \partial }{\rm \partial x\sp \alpha} -(\frac{\rm ie}{\rm \hbar})A\sb{\alpha}]$. Tienen la misma forma geométrica, $t$ es la coordenada cero, y los demás $x$ son las coordenadas espaciales. Se unifican como  $[ \frac{\rm \partial }{\rm \partial x\sp j} +\omega \sb j]$, donde para las coordenadas espaciales
$\omega \sb j =- (\frac{\rm ie}{\rm \hbar})A\sb{\alpha}$
y para la coordenada temporal $\omega \sb o = -(\frac{\rm ie}{\rm \hbar}) \phi$

Reconociendo la forma de derivada covariante en estos operadores, podemos con toda justicia adjudicarles el nombre correspondiente:

$\nabla \sb o = \frac{\rm \partial }{\rm \partial t} -(\frac{\rm ie}{\rm \hbar}) \phi$

$\nabla \sb \alpha =  \frac{\rm \partial }{\rm \partial x\sp \alpha} -(\frac{\rm ie}{\rm \hbar})A\sb{\alpha}$

 Teniendo en cuenta que $A$ es una 1-forma definida por : $A\sp  1 = \phi dt + A\sb \alpha dx\sp \alpha$,  se puede simplificar más. Resumamos ese resultado en el siguiente teorema.

\bigskip

\addtocounter{ecu}{1}

\textit{Teorema (\theecu). El electromagnetismo tiene una \index{derivada covariante} \textbf{derivada covariante} asociada dada por:}

$\nabla \sb j =  \frac{\rm \partial }{\rm \partial x\sp j} +\omega \sb j$

\textit{donde}

$\omega \sb j = -(\frac{\rm ie}{\rm \hbar})A\sb{j}$

\textit{y los subíndices corren por todos las vierbein (nombre en alemán, que también se usa para denominar las coordenadas del espacio-tiempo).}

\bigskip

Podemos observar que la conexión electromagnética no opera solamente sobre el espacio-tiempo sino que incluye una variable compleja, completamente imaginaria. El significado de dicha variable es el de la fase de la función de onda. Podemos decir que es $U(1)$.

Hemos deducido que el campo electromagnético, capturado en su potencial 1-forma, es exactamente la conexión que define la derivada covariante que geometriza la interacción electromagnética. Por eso, se acostumbra a decir que la \index{conexión electromagnética}  \textbf{conexión electromagnética} es
$\omega = -(\frac{\rm ie}{\rm \hbar})A$  o bien,  $\omega \sb j = - (\frac{\rm ie}{\rm \hbar})A\sb{j}$

Con la nueva equivalencia, la ecuación para una partícula cargada se escribe:

\addtocounter{ecu}{1}

$$i\hbar \nabla \sb o \psi = -(\frac{\rm \hbar \sp 2}{\rm 2m})\sum \nabla \sb \alpha \nabla \sb \alpha \psi  + V\psi \eqno{(\theecu)}$$

en donde la suma del lado derecho se corre sobre las coordenadas espaciales.

Como habíamos previsto, el campo electromagnético está totalmente descrito por su 1-forma $A\sp  1 = \phi dt + A\sb \alpha dx\sp \alpha$,  en el sentido que todo su efecto sobre una partícula cargada se deriva de la ecuación anterior definida en su totalidad por $A\sp 1$. La curvatura es entonces la diferencial de la conexión y es la representación del tensor electromagnético:

$\theta = dw = -(\frac{\rm ie}{\rm \hbar})dA\sp 1 = -(\frac{\rm ie}{\rm \hbar})F\sp 2 $

$=-(\frac{\rm ie}{\rm \hbar})(\epsilon\sp 1 \wedge dt + \beta\sp 2)$

\subsection{Con pseudo-gravitación}

Habiendo considerado que los campos funcionan como operadores nos cambiamos a la mecánica cuántica que es el mundo natural de los operadores. Ya implementamos una  geometrización del campo electromagnético y logramos asociarle una derivada covariante o conexión $\omega$. Quisiéramos ahora considerar el efecto de la gravitación. Nosotros ya sabemos que  la gravitación es creada por una métrica pseudo-Riemanniana, con posibles valores negativos, y que su efecto se hace sentir sobre las geodésicas a través de su conexión expresada en los símbolos de Christoffel $\Gamma \sp r \sb{k j}$. Lo más natural sería estudiar el efecto de los dos campos asumiendo que el campo total tiene una conexión que es la suma de las dos.

Pues no tan rápido: la ecuación de la mecánica cuántica que utilizamos es la de Schr\"{o}dinger, y cuando derivamos dicha ecuación asumimos que el espacio y el tiempo eran tan inmiscible como agua y aceite, tal ecuación no involucra el espacio-tiempo. La ecuación que logra eso se llama la ecuación de Dirac y su compatibilidad con gravitación no parece tener esperanza. Eso se debe a que las partículas descritas por dicha ecuación tienen spin, o si uno se permite un lenguaje clásico, giran en torno de sí mismas, y como están cargadas crean un campo magnético, lo cual nos permite ofrecer una explicación al magnetismo del hierro. Lo malo es que el tal spin y la gravitación no han encontrado lenguaje común.

Como consuelo, consideremos una historia en la cual la gravitación no tuviese efecto temporal y todo su poder se redujera a curvar el espacio 3:d. Tenemos entonces la conexión electromagnética, $\omega$, que es un múltiplo de la 1-forma potencial, y por otro lado, la conexión métrica, dada por los $\Gamma \sp r \sb{k j}$. Ahora, superponemos los dos campos, pero hay que hacerlo de forma que tenga sentido, el cual resulta de considerar la siguiente pregunta: si sobre una variedad, el espacio, tenemos dos estructuras, el espacio tangente que permite hablar de velocidades y la función de onda, dada por un número complejo en cada punto, entonces cómo debe considerarse sobre toda la estructura la derivada covariante que extienda de manera simple las derivadas covariantes en cada estructura? A falta de no querer estudiar una respuesta juiciosa a esta pregunta, presentamos una receta equivalente.

Teníamos la derivada covariante del electromagnetismo:

$\nabla \sb j =  \frac{\rm \partial }{\rm \partial x\sp j} +\omega \sb j$

$\nabla \sb j \psi  =  \frac{\rm \partial \psi }{\rm \partial x\sp j} +\omega \sb j \psi $

ahora comenzamos a abrirle campo al efecto de la métrica:

$\nabla \sb j \psi  =  \frac{\rm \partial \psi }{\rm \partial x\sp j} + \omega \sb j \psi - \Gamma \sp r \sb{k j} \psi $
$=\partial \sb j \psi  + \omega \sb j \psi - \Gamma \sp r \sb{k j} \psi$

\
Con esta familia de términos podemos armar un operador de segundo grado sin subíndices para fabricar la ecuación de Schr\"{o}dinger en \index{espacio curvo} \textbf{espacio curvo}, es decir, bajo el efecto de la gravitación:

$$i\hbar \nabla \sb o \psi = -(\frac{\rm \hbar \sp 2}{\rm 2m})\nabla \sp 2 \psi  + V\psi $$

donde

$$\nabla \sp 2 \psi = g\sp{jk} \nabla \sb j \nabla \sb k \psi$$

y

$$ \nabla \sb j \psi  =  \frac{\rm \partial \psi }{\rm \partial x\sp j} + \omega \sb j \psi - \Gamma \sp r \sb{k j} \psi $$
$$ \  \  \  \ =\partial \sb j \psi  + \omega \sb j \psi - \Gamma \sp r \sb{k j} \psi$$

Al operador $\nabla \sp 2  $  se le llama el Laplaciano:

$\nabla \sp 2 \psi = g\sp{jk} \nabla \sb j \nabla \sb k \psi$

\

El  mérito de la nueva ecuación es el de mostrar el efecto total como una suma del efecto en solitario de la gravitación, mas el efecto en solitario del campo electromagnético, mas un efecto de interacción de los dos campos. Es decir, en mecánica cuántica, los dos campos se entremezclan. Uno podría tener la esperanza de que con estos términos se prediga el efecto de un campo gravitatorio muy débil sobre el electromagnético, pero eso no tiene interés experimental directo: el campo eléctrico es unas  $10\sp{40}$ veces más fuerte que el gravitatorio. Eso es lo que permite que existan estructuras eléctrica y localmente (casi) neutras pero que sin embargo  desafíen la gravitación, como un árbol de 70 metros de altura.

\section{LA LIBERTAD GAUGE}

Todas las propiedades del campo electromagnético las hemos podido deducir, tanto en mecánica clásica, como en relatividad, como en mecánica cuántica, a partir del \index{potencial vector } \textbf{potencial vector}, $A$. En mecánica clásica, $A$ era un vector con tres componentes cuyo rotacional determinaba el campo magnético. Pero como $\nabla \times \nabla \phi =0$ para cualquier función escalar $\phi$, entonces, $\nabla \times (A + \nabla \phi) = \nabla \times A  $, por lo tanto, el vector potencial no se puede especificar completamente, sino que dos maneras diferentes de determinarlo difieren por el gradiente de una función escalar.

Gracias al lenguaje de las formas, la libertad gauge pudo escribirse en el espacio tiempo simplemente como $d\sp 2 =0$, donde $d$ es la derivada exterior. En efecto, si $A\sp 1$ es la 1-forma potencial, entonces $dA$ da el tensor electromagnético, que es una 2-forma. Pero entonces $dA\sp 1 = dA\sp 1 + d\sp 2 \phi$ donde $\phi$ es cualquier función escalar.

En mecánica cuántica,  la libertad gauge asociada al campo electromagnético (sin gravitación) consiste en la capacidad de  elegir, en cada punto del espacio de manera independiente, la fase cero del valor de cada función de onda. La única precaución es que todas las funciones de onda tengan el mismo modo de calibración (la palabra gauge significa medir, calibrar). Aunque la libertad gauge es inherente a cada punto, se acostumbra a acomodarla dentro de la teoría de los haces fibrados, los cuales son complicaciones asociadas a las variedades. Lo que ahora importa resaltar al respecto es esto: las variedades están definidas por retazos con intersecciones no vacías y la unidad resulta de la compatibilidad de las descripciones en
dichas intersecciones. Esto debe ser especialmente válido para  la libertad gauge, la cual permite cambiar de fase de un retazo a otro. A esta construcción se le llama fibrado. En el próximo capítulo veremos el \index{fibrado electromagnético}  \textbf{fibrado electromagnético}.

\bigskip

Finalmente, debemos decir que carecemos de una teoría geométrica que unifique electromagnetismo y gravitación, pero fabricamos una ecuación que presumiblemente hace esto  para campos gravitatorios débiles. Queda por demostrar que tal ecuación conserva las probabilidades. Respondamos olímpicamente a tal requerimiento: el efecto de la gravitación en coordenadas rectilíneas es equivalente a carecer de ese efecto pero en coordenadas curvas. Y un cambio de coordenadas ni crea ni destruye partículas.

\section{CONCLUSION}

Con la magia de un gran partido de fútbol hemos geometrizado el campo electromagnético y hasta produjimos una teoría de juguete que unifica la gravitación y electromagnetismo.

\section{REFERENCIAS}

\begin{enumerate}

\item  Frankel T. \textit{The Geometry of Physics, An Introduction}. Cambridge University Press, 2001.

\end{enumerate}

\chapter{EL FIBRADO ELECTROMAGNETICO}        

\Large

\centerline{RESUMEN}

\bigskip
Demostramos que la geometrización del electromagnetismo viene naturalmente de un grupo. Para \\ ello, se introduce la noción de fibrado principal como la maquinaria matemática encargada de implementar la directiva geométrica para  la formulación de interacciones:  si las arbitrariedades matemáticas de un formalismo que describe una interacción forman un grupo, entonces dicho grupo es el responsable de la interacción. Se desarrolla  el caso del electromagnetismo con grupo $U(1)$.

\normalsize

\section{INTRODUCCION}

Tenemos el gusto de presentar en este capítulo  una introducción  a uno de los temas más sofisticados y más productivos de la geometría diferencial y de la física matemática: el concepto de fibrado principal (principal bundle) y su aplicación directa a las teorías gauge,  específicamente al electromagnetismo.

Lo que haremos es implementar un cambio de perspectiva. Hasta ahora teníamos: el grupo gauge es una característica de las arbitrariedades admitidas por el formalismo matemático que describe la interacción. Pero de ahora en adelante, el grupo gauge será para nosotros el corazón mismo de la interacción. Advertimos que el viejo punto de vista no es incompatible con el nuevo y que los dos coexisten, y que ambos se utilizan a cualquier momento. Nuestro cambio de perspectiva ha de originar una reestructuración de lo que ya sabemos y algo tendrá que  añadir. ¿Qué teníamos?

Habíamos  representado al potencial vector como una 1-forma diferencial $A_\mu (x)dx^\mu$, la cual para cada punto del espacio-tiempo es una transformación lineal que toma vectores de $\Re^4$, el espacio tangente, y los asocia con un escalar real. Codificado de esa forma, $A_\mu (x)$ contiene dos de las leyes de Maxwell, mientras que las otras dos se relacionan con el dual del tensor de campo, la 2-forma $F=dA$. Todo esto se hizo dentro de la mecánica relativista  no cuántica. Pero por otro lado, la mecánica cuántica está ligada a la fase y en particular al grupo $U(1)$, cuyos elementos son de la forma $e^{i\theta}$.

Nuestro nuevo enfoque  mezcla todo eso, inescrupulosamente, en una sola entidad llamada \index{fibrado} \textbf{fibrado principal}.  Para la interacción electromagnética la situación es  sencilla, perfectamente entendible, tal como entramos a explicar. Queda preparado el terreno para el estudio de las interacciones  débil y fuerte. El juego  consiste en que  en vez de tomar elementos $e^{i\theta}$ en $U(1)$ se toma la correspondiente generalización en $SU(2)$ o en $SU(3)$, lo cual significa que en vez de tomar   $ \theta$ real, se toma una  matriz, por lo que los estados del campo de la otras interacciones ya no podrán representarse por una sola coordenada, sino que se requerirá de dos o tal vez tres.

\section{FIBRADOS PRINCIPALES}

Nuestra nueva maquinaria puede motivarse reflexionando sobre el concepto de fuerza y cambiando el espacio externo por el espacio interno que es el grupo gauge. El concepto de fuerza lo vemos como sigue:

La derivada de una parametrización $M$ es la velocidad,  la cual es un vector. Para un punto dado, dichos vectores forman un espacio vectorial llamado el espacio tangente a la curva en el punto dado, y la unión de todos esos espacios tangentes se denomina el espacio o \index{fibrado tangente} \textbf{fibrado tangente}, notado $TM$. La aceleración a su vez es la derivada de la velocidad, por lo tanto está en el espacio tangente del espacio tangente, el cual es $T(TM)$. La aceleración es un múltiplo escalar de la fuerza, el motor de todo. Por lo tanto, la física fundamental clásica que describe un sistema pasivo en el espacio externo bien puede estar contenida en  el siguiente diagrama:

$$ M\rightarrow TM \rightarrow T(TM)$$

Nuestro punto de quiebre dice así: mientras que la fuerza es  una instrucción que le dice a la masa cómo cambiar su aceleración en el espacio externo, un campo de una interacción es una instrucción que dice a la partícula cómo cambiar el espacio interno y cómo resultado cambiar la aceleración en el espacio externo.

Cuando queremos pasar de la descripción de una partícula movida por una fuerza   en un espacio  externo a la descripción de una partícula en un campo   entonces se usa la potente maquinaria del \index{fibrado principal} \textbf{fibrado principal} $P$, el cual está asociado con el siguiente diagrama, que es una modificación gramatical del diagrama asociado a la fuerza:

$$M\rightarrow  P \rightarrow T(P)$$

donde $P$ es un grupo, el \index{grupo gauge} \textbf{grupo gauge}, y $T(P)$ es el espacio tangente de $P$. El objetivo es asegurar que la física  sea una consecuencia de todo el diagrama. En particular, debemos encontrar una expresión para la fuerza, el cual vendrá por supuesto, en términos de curvatura. Presentaremos dos visiones: una que involucra  una forma diferencial y otro que da un conmutador.

 Puesto que la curvatura se relacionaba, en el primer diagrama,   con aceleraciones o segundas derivadas, la curvatura debía ser una forma de por lo menos grado 2. Nuestra implementación dio una forma de grado 3 con valores en el espacio tangente (  con opción de convertirse en una de grado 4, con valores en los reales). Pero en el segundo diagrama, el espacio tangente es reemplazado por un grupo de Lie. Por lo tanto, lo que antes era una segunda derivada ahora deberá ser una primera derivada. Por consiguiente, las nuevas curvaturas podrán ser formas de orden menor que las del diagrama de las fuerzas, el diagrama uno.

 \bigskip

 Explicitemos todo este programa para la interacción electromagnética.

 \bigskip

El fibrado principal $P$ para la interacción electromagnética  es extremadamente simple y nada sorpresivo: se trata del producto cartesiano del espacio -tiempo por el grupo de la fase $U(1)$, notado como
 $P=\Re^4 \times U(1)$.

\bigskip
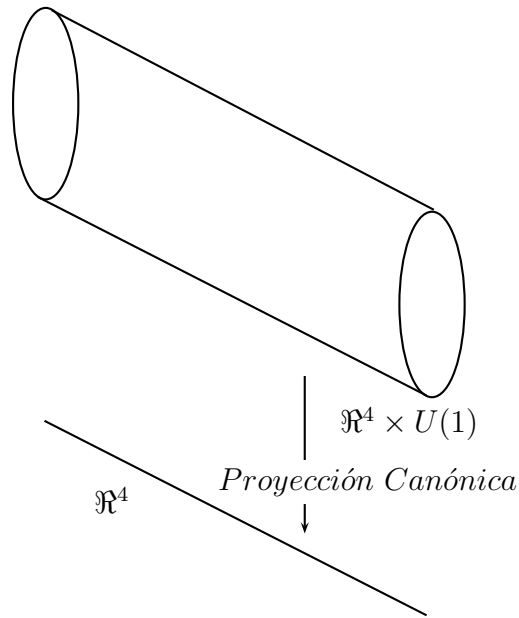
\psset{unit=0.4cm}
\begin{center}
 \begin{pspicture} (1,1)(17,20)  
 \psline(2.76,7.75)(15.46,1.28)
 \psline(2.53,15.2)(15.46,8.59)
 \psline(3,21.41)(15.7,14.77)
 \psellipse(2.8,18.3)(1.12,3.22)
 \psellipse(15.64,11.62)(1.12,3.12)
 \psline{->}(11.41,9.24)(11.41,4)
 \rput*(13.5,5.71){$Proyecci\acute{o}n \ Can\acute{o}nica$}
  
 \rput(5,5.1){$\Re^4$}
 \rput(14.88,7.65){$\Re^4 \times U(1)$}
 
 \end{pspicture}
 
 \textit{El haz electromagnético  $\Re^4 \times U(1)$, consta de un aro en cada punto del espacio-tiempo.}
 
\end{center}

Es usual imaginar la estructura compuesta $\Re^4 \times U(1)$ por medio de una caricatura sencilla de manejar: al espacio-tiempo se lo imagina uno de una sola dimensión, como un hilo, y al grupo $U(1) $ como un  anillo. El producto cartesiano de esos dos elementos da un cilindro el cual es el  fibrado principal.

 \bigskip

 ¿De qué manera se extrae la física de esta entidad geométrica? Veamos.

\bigskip

Consideremos una trayectoria sobre $U(1)$, es decir una función $\theta (t)$ que a un $t$ le asocia un elemento $\theta (t)$ $=e^{i\alpha(t)}$, donde $\alpha(t)$ es una función de $\Re$ en $\Re$. Su velocidad es $i\dot{\alpha} e^{i\alpha}$ . Cada uno de estos vectores es tangente a $U(1)$ y variando $\alpha$ se tienen todos los posibles vectores tangentes que forman el espacio tangente en un punto dado. Tomemos una curva $\theta (t)$ que en $t=0$ pasa por $1= (1,0)$, es decir $\alpha (0)=0$: la velocidad en ese momento es  $i\dot{\alpha}(0) e^{i\alpha(0)} = i\dot{\alpha}(0)e^0 = i\dot{\alpha}(0)$.

Variando la curva se obtienen diferentes realizaciones de $i\dot{\alpha}(0)$ y con todas las trayectorias se obtienen todos los vectores tangentes a $U(1)$ en 1, que es su elemento unidad. Esos vectores forman un espacio vectorial de dimensión 1, isomorfo a $\Re$, y que notamos $i\Re$, al cual se le denomina el \index{álgebra de Lie} \textbf{álgebra de Lie} de $U(1)$. En general, el \textbf{álgebra de Lie} de un grupo es el espacio tangente al grupo en la unidad del grupo. 

Observemos también que el espacio tangente a $U(1)$ en cualquier otro punto es también isomorfo a $i\Re$. Uno podría representar un vector tangente a $U(1)$ en coordenadas cartesianas, pero, como hemos visto, sale muy fácil  representar un tal vector en coordenadas polares, de la forma $iK\partial /\partial \theta$, donde $K$ es la norma del vector tangente. Notemos también que en un punto dado, el espacio tangente a $U(1)$ es isomorfo al álgebra de Lie que es el espacio tangente sobre el punto 1 = (1,0).  En realidad, es la misma álgebra de Lie  pero trasladado en  transporte paralelo sobre $U(1)$, que es un círculo, desde $(1,0)$ hasta el punto de operación.

\section{EL ESPACIO TANGENTE}

Hemos hablado de dinámica pero hemos visto el espacio tangente tan sólo de $U(1)$, pero nada más. Ahora bien, el \index{espacio tangente} \textbf{espacio tangente} al fibrado principal  $\Re^4 \times U(1)$, que es como un cilindro, es  $T(\Re^4) \times T(U(1))$, el cual  en realidad tiene 10 dimensiones, 5 donde se ancla, en el espacio base, el cilindro, y otras 5 creadas por los vectores tangentes. Pero el espacio tangente en un punto específico del cilindro tiene 5 dimensiones y  lo visualizamos como un 'plano'. Sobre dicho 'plano' trazamos 'dos ejes:  El primer eje es el formado por el espacio tangente a $U(1)$ que es isomorfo a $i\Re$ y al cual se le llama \index{espacio vertical} \textbf{espacio vertical}, $V(P)$. El segundo 'eje'  da coordenadas al espacio tangente a $\Re^4$ que es otro $\Re^4$. Así que el segundo 'eje'  es en realidad un manojo de 4 ejes.

\section{LA FORMA GENERADORA}

Ahora viene lo artístico. En primer término, fijamos un punto del fibrado $P$, nuestro cilindro. Pero debido  a que todos los espacios tangentes a $U(1)$ son isomorfos a $i\Re$, y a que desde el punto de vista rotacional $U(1)$ no tiene ningún punto preferido, podemos proceder con tanta sencillez como si todo sucediese en 1, el elemento neutro  de $U(1)$. 

En  $T(P)$,  el espacio vectorial tangente a $P$ en un punto dado que no se especifica, tenemos el espacio vertical, isomorfo al álgebra de Lie, y le buscamos un complemento  a dicho espacio. A dicho complemento lo  llamamos \index{espacio horizontal} \textbf{espacio horizontal} y lo notamos $H(P)$. Tendremos $T(P)=V(P)+H(P)$ en cada punto del fibrado principal $P$, que estamos visualizando como un cilindro. Pero atención, el complemento $H(P)$ no ha de ser $\Re^4$. Precisamente, el problema es que dicho complemento no es único. Es lo mismo que pasa en el plano cartesiano: el único complemento al eje $\vec Y$ no es el eje $\vec X$, sino que además de éste hay muchos más. ¿Cómo elegiremos el apropiado?

\bigskip
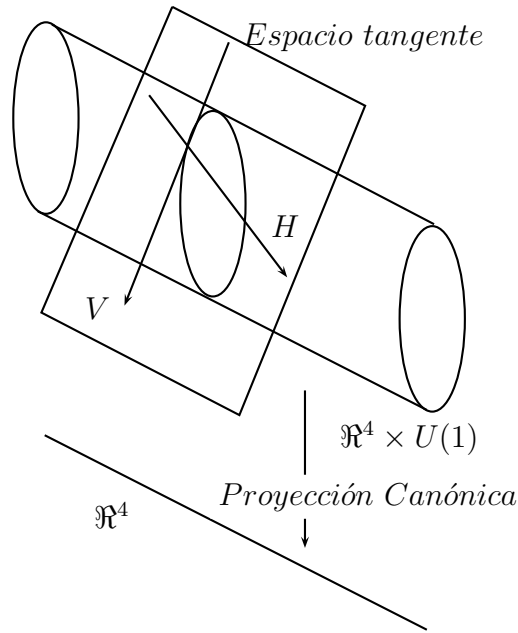
\psset{unit=0.4cm}
\begin{center}
 \begin{pspicture} (1,1)(17,22)

 \psline(2.76,7.75)(15.46,1.28)
 \psline(2.53,15.2)(15.46,8.59)
 \psline(3,21.41)(15.7,14.77)
 \psline(7,22)(13.41,18.71)(9.23,8.42)(2.65,11.83) (7,22)
 \psline{->}(8.88,20.82)(5.41,12.06)
 \psline{->}(6.23,19.06)(10.82,13,83)
 \psellipse(2.8,18.3)(1.12,3.22)
 \psellipse(8.35,15.45)(1.12,3.12)
 \psellipse(15.64,11.62)(1.12,3.12)
 \psline{->}(11.41,9.24)(11.41,4)
 \rput*(13.5,5.71){$Proyecci\acute{o}n \ Can\acute{o}nica$}
  
 \rput(5,5.1){$\Re^4$}
 \rput(14.88,7.65){$\Re^4 \times U(1)$}
 \rput(4.53,12,18){$V$}
 \rput(10.76,14.71){$H$}
 \rput(13.52,21){$Espacio \ tangente$ } 

 \end{pspicture}
 
 \textit{El haz electromagnético  $\Re^4 \times U(1)$, con el espacio tangente a un punto, generado por su espacio vertical $V$ y el espacio horizontal $H$.}
 
\end{center}

El espacio horizontal que nos conviene es aquel que nos permita  recobrar la física fundamental del campo electromagnético.  Pero lo haremos indirectamente como sigue. De igual manera que en el plano cartesiano  un eje $W$ que sea complemento   al eje $\vec Y$ está determinado por una proyección, la proyección a lo largo de dicho eje $W$ y que proyecta dicho espacio sobre cero, así mismo nosotros fabricaremos una \index{proyección} \textbf{proyección} que proyecte una sombra sobre el eje vertical $V(P)$ que aniquile cierto espacio, el espacio $H(P)$. Hagámoslo paso a paso.

Una curva en el espacio-tiempo se denota como $\gamma (t) = (x_\mu (t))$, en tanto que una curva en $P$, el cilindro,  se denota como
$\Gamma (t) = (x_\mu (t), \theta (t))$.

Recordemos que nosotros utilizamos un isomorfismo o equivalencia entre vectores y operadores de derivada direccional. Por ejemplo, al vector en el plano $(3,5)$ le asociamos el operador de derivada direccional $3\partial/\partial x + 5 \partial /\partial y$. O bien, si una curva se describe por $x^\mu $ y su velocidad por  $\dot{x} ^\mu$  entonces el operador de derivada direccional asociado es $\dot{x} ^\mu \partial /\partial x^\mu $. Similarmente, al vector tangente a $U(1)$ dado por $\dot{\theta}$ le asociamos el operador $\dot{\theta}\partial /\partial \theta .$ Por tanto, los vectores tangentes a $P$ son de la forma (velocidad en el $\Re^4$, velocidad en el espacio $U(1)$). Concretamente, un vector tangente a $P$ es de la forma: $\dot{x} ^\mu \partial /\partial x^\mu + \dot{\theta} \partial /\partial \theta $.

Pasemos ahora a determinar el espacio horizontal $H(P)$. Como hemos dicho, eso es lo mismo que determinar una proyección que aniquile sobre cero a dicho espacio. La proyección está dada por una forma diferencial especial $\omega $ que al operar sobre cualquier vector horizontal produzca el elemento cero. Pero atención, la sombra proyectada cae sobre el espacio vertical. Por eso, ese cero es el cero del álgebra de Lie $i\Re$. En la jerga oficial, a la proyección  se le llama forma diferencial con valores en el álgebra de Lie, que no es lo mismo que las formas diferenciales ordinarias que toman  vectores, los procesa y producen un número. Todo eso suena mucho más extraño que la siguiente definición de dicha forma, la \index{forma electromagnética} \textbf{forma electromagnética}:

$\omega =  -ieA_\mu dx^\mu +d\theta $

Para comprobar que esta forma toma vectores en $T(P)$ y produce vectores en el álgebra de Lie $i\Re$ simplemente hagamos un cálculo directo.

Las formas diferenciales operan sobre los vectores tangentes de $P$ obedeciendo  las siguientes reglas de ortonormalidad local:

1) $d\theta(\partial /\partial \theta)=1$

2) $dx^\mu(\partial /\partial x^\mu)=1$

3) Todas las demás combinaciones dan cero y la operación se extiende por linealidad tanto a todo $T(P)$ como a las combinaciones lineales de formas diferenciales.

Ahora bien, un vector en $T(P)$ tiene la forma $\dot{x} ^\mu \partial /\partial x^\mu + \dot{\theta} \partial /\partial \theta $ que reescribimos como  $\dot{x} ^\mu \partial /\partial x^\mu + iK \partial /\partial \theta $. Operemos a $\omega$ sobre dicho vector recordando las reglas de ortonormalidad local. La localidad significa que aunque los coeficientes de los vectores de $T(P)$ varíen de punto a punto del fibrado, realmente sólo importa la descomposición en el $T(P)$ sobre el punto fijo elegido:

$\omega (\dot{x} ^\mu \partial /\partial x^\mu + \dot{\theta} \partial /\partial \theta )= $
$(-ieA_\mu dx^\mu +d\theta)(\dot{x} ^\mu \partial /\partial x^\mu + iK \partial /\partial \theta ) $

$= (-ieA_\mu dx^\mu)(\dot{x} ^\mu \partial /\partial x^\mu + iK \partial /\partial \theta) +  (d\theta ) (\dot{x} ^\mu \partial /\partial x^\mu + iK \partial /\partial \theta )$

$=i(-eA_\mu dx^\mu)(\dot{x} ^\mu \partial /\partial x^\mu ) +  (d\theta) (  iK  \partial /\partial \theta )$

$=i[(-eA_\mu dx^\mu)(\dot{x} ^\mu \partial /\partial x^\mu ) +  (d\theta) (  K  \partial /\partial \theta )] $

$ = i[(-eA_\mu \dot{x} ^\mu (dx^\mu)( \partial /\partial x^\mu ) + K (d\theta) (    \partial /\partial \theta )]  = i[(-eA_\mu \dot{x} ^\mu  + K ]$

lo cual demuestra que el resultado de operar la forma $\omega$ sobre un elemento de $T(P)$ da $i$ veces un real, es decir,  el resultado total está en $i\Re$, el álgebra de Lie.

Ahora bien, lo que tenemos en mente es que $A=A_\mu dx^\mu$ sea el  potencial vector electromagnético y que de $\omega $  saquemos las leyes que rigen la conducta de dicho  campo.  Eso se logra a partir de $\omega = -ieA_\mu dx^\mu +d\theta$ tomando la deriva exterior:

$d\omega = -ied(A_\mu dx^\mu) +d^2\theta = -ied(A_\mu dx^\mu) = -iedA$

Como ya sabemos, el tensor de campo $F$, una 2-forma, y el potencial vector $A$, una 1-forma, se relacionan por

$F=dA = (\partial_\nu A_\mu -\partial_ \mu A_\nu)dx^\mu \wedge dx^\nu$.

Por tanto

$d\omega = -ieF$.

Hemos demostrado así que \emph{la forma diferencial con valores en el álgebra de Lie de $U(1)$ contiene la física del campo electromagnético en su parte geométrica}: de aquí se sacan dos leyes de Maxwell (pero quedan faltando dos que salen de análisis variacional). 

Hasta ahora no hemos visto cuál es el papel de la parte angular del fibrado, además de complicar todo. Su utilidad se verá al estudiar el espacio que es aniquilado por $\omega$.

\section{EL ESPACIO HORIZONTAL}

El \index{espacio horizontal}  \textbf{espacio horizontal} es el conjunto de todos los vectores $h$ que son aniquilados por $\omega$:

$H(P)=\{ h \in T(P) : \omega(\vec h) =0$ \}

Puesto que queremos que $T(P) = V(P) \times H(P)$ entonces cada vector en  $T(P)$ podrá descomponerse en una base determinada por dicha descomposición. Tal descomposición es única debido a que lo único en común entre el espacio horizontal y el vertical es el vector cero.  Hallemos pues una base apropiada. Un vector en $V(P)$ es un vector tangente a $U(1)$, un elemento de la forma $iK \partial /\partial \theta$. Un elemento del espacio horizontal $h$ también es un elemento  del espacio tangente de $P$, el cual escrito en la base natural es una combinación lineal de $ \{ \partial /\partial x^\mu,\partial /\partial \theta\}$. En general un elemento cualquiera $\vec X$ de $T(P)$ tiene la forma:

$\vec X = \dot{x} ^\mu \partial /\partial x^\mu + iK \partial /\partial \theta $
 $= \epsilon^\mu  \partial /\partial x^\mu+ \lambda  \partial /\partial \theta$.

quedando claro  que los coeficientes de la combinación puede cambiar de punto a punto.

Lo anterior dice también que una base para $T(P)$ tiene 4 vectores de la forma  $\partial /\partial x^\mu$ y otro de la forma $\partial /\partial \theta$.  Demostremos ahora que 4 vectores de la forma 

$ \partial /\partial x^\mu+\kappa_\mu \partial /\partial \theta$ 

conforman una base para el espacio horizontal.  Para demostrarlo, comencemos notando que dichos vectores son linealmente independientes pues van en direcciones distintas del espacio-tiempo. Por lo tanto, esta base genera un espacio de dimensión 4.  Para demostrar que son una base para el espacio horizontal nos resta demostrar que  cualquier combinación lineal con estos vectores   es  aniquilada por $\omega$. Veamos:
De acuerdo con nuestro supuesto, un vector cualquiera del espacio horizontal puede escribirse como:
$h=  \beta^\mu ( \partial /\partial x^\mu+  \kappa_\mu \partial /\partial \theta )$

Por otra parte $\omega (h) =  <\omega, h > =0$, es decir:

$0= <\omega, h > = <-ieA_\mu dx^\mu +d\theta , \beta^\mu ( \partial /\partial x^\mu+  \kappa_\mu \partial /\partial \theta )>$

Utilizando la bilinealidad de la acción de las formas diferenciales sobre los vectores, eso se reduce a:

$0=    -ieA_\mu  \beta^\mu + \beta^\mu \kappa_\mu $

como eso debe ser cierto para todo $\beta^\mu $, es decir, para cualquier elemento del espacio horizontal,  tenemos que:

$  -ieA_\mu+\kappa_\mu     =0$

o sea $\kappa_\mu = ieA_\mu  $

Este resultado da la condición para que la combinación lineal propuesta esté dentro del espacio horizontal.

\bigskip

Por lo tanto, el espacio horizontal está generado por 4 elementos de la forma $ \partial /\partial x^\mu+  \kappa_\mu \partial /\partial \theta =  \partial /\partial x^\mu  + ieA_\mu  \partial /\partial \theta$.

\bigskip

Para qué nos tomamos tanto esfuerzo en caracterizar el espacio horizontal? Para definir la Derivada covariante $D$ como la proyección de la derivada ordinaria sobre el espacio horizontal. Con más exactitud tenemos.

\section{LA DERIVADA COVARIANTE}

La derivada ordinaria es aquel operador que de una curva en $P$  saca un vector tangente. La expresión de dicho operador lo sacamos de la identidad siguiente, que se deriva de la regla de la cadena:

$d/dt =\dot{x} ^\mu \partial /\partial x^\mu + \dot{\theta} \partial /\partial \theta $

Este operador produce vectores en $T(P)$, cuya descomposición en $V(P) + H(P)$ se escribe como

$\vec X= \alpha \partial /\partial \theta + \beta^\mu ( \partial /\partial x^\mu  + ieA_\mu  \partial /\partial \theta)  $

La primera parte está sobre el espacio vertical y la segunda sobre el horizontal. Definimos la   \textbf{derivada covariante $D$} como la proyección sobre el espacio horizontal de la derivada ordinaria. Es decir, olvidamos la primera parte. Su expresión componente por componente es:

$D_\mu =  \partial /\partial x^\mu  + ieA_\mu  \partial /\partial \theta  $

Obsérvese que la expresión para la derivada covariante se dedujo a partir de la expresión para la forma $\omega$ sobre $T(P)$ que contiene el potencial vector.

\section{EL CONMUTADOR}

De qué  manera  está escondida la física fundamental que regula el campo electromagnético en el nuevo formalismo del fibrado principal? Lo que tenemos que lograr es extraer la ley $dA=F$. Eso se puede sacar de dos maneras, la primera ya la vimos y es tomando la deriva exterior de $\omega$. La segunda forma es tomando el conmutador de la derivada covariante,  el cual está ligado con la curvatura, la cual   da la fuerza del campo. Hagamos el cálculo del conmutador  usando la convención opcional de notar $iA_\mu \partial /\partial \theta$ con el símbolo  $iA_\mu$, el cual es compatible con nuestro isomorfismo básico de identificar vectores con operadores. Eso equivale a notar a

$D_\mu =  \partial /\partial x^\mu  + ieA_\mu  \partial /\partial \theta  $

como

$D_\mu =  \partial /\partial x^\mu  + ieA_\mu = \partial_\mu  + ieA_\mu$

Podemos entonces calcular el \index{conmutador} \textbf{conmutador}

$[D_\mu,D_\nu]= D_\mu D_\nu - D_\nu D_\mu$

$=(\partial_\mu  + ieA_\mu) (\partial_\nu  + ieA_\nu)- (\partial_\nu  + ieA_\nu) (\partial_\mu  + ieA_\mu)$

$=\partial_\mu \partial_\nu + ie \partial_\mu A_\nu + ie A_\mu \partial_\nu - e^2 A_\mu  A_\nu$

 $-(\partial_\nu   \partial_\mu + ie \partial_\nu A_\mu+ ieA_\nu \partial_\mu-e^2 A_\nu A_\mu)$

$=ie( \partial_\mu A_\nu -  \partial_\nu A_\mu ) + ie(A_\mu \partial_\nu -A_\nu \partial_\mu )$

para poder interpretar correctamente la expresión anterior es suficiente tener en cuenta que todos los términos son operadores que operan sobre una función escalar, digamos $\psi$. Calculemos entonces usando la regla del producto:

$[D_\mu,D_\nu]\psi = ie( \partial_\mu A_\nu \psi -  \partial_\nu A_\mu \psi ) + ie(A_\mu \partial_\nu \psi -A_\nu \partial_\mu \psi )$

$ = ie( \partial_\mu (A_\nu \psi) -  \partial_\nu (A_\mu \psi) ) + ie(A_\mu \partial_\nu \psi -A_\nu \partial_\mu \psi )$

$=ie(\partial_\mu A_\nu) \psi+ie A_\nu \partial_\mu \psi) -  ie((\partial_\nu A_\mu) \psi+A_\mu \partial_\nu \psi)  + ie(A_\mu \partial_\nu \psi -A_\nu \partial_\mu \psi )$

$=ie( (\partial_\mu A_\nu) \psi- (\partial_\nu A_\mu)\psi) = ie( \partial_\mu A_\nu -\partial_\nu A_\mu)(\psi) $

por lo que podemos concluir que:

$[D_\mu,D_\nu]\psi=-ieF_{\mu \nu} (\psi)$

es decir:

$[D_\mu,D_\nu] = -ieF_{\mu \nu}$.

Así hemos demostrado que la física del electromagnetismo está en el conmutador de la derivada covariante.

\section{CAMBIO DE GAUGE}

La teoría que hemos expuesto es totalmente solipsista: ha tenido en cuenta tan sólo lo que un experimentador solitario podría describir. Pero qué pasaría si otro experimentador desea poner a prueba las ideas aquí expuestas? En primer término, el debería sentirse con el derecho de cambiar de escala, de coordenadas, de fijar toda arbitrariedad a su antojo. Se necesita pues unas reglas que permitan hacer comparaciones ante esta situación.

Por otro lado, no habría en principio ninguna objeción a que un experimentador solitario investigue el universo. Lo que hay que tener en cuenta es que todo el universo, el externo y el interno, se pueda o no describir por un solo sistema de coordenadas. Veamos. Por un lado, el fibrado principal asociado al campo electromagnético es trivial en el sentido que dicho fibrado es el producto de dos espacios.  El espacio -tiempo  puede describirse con un sólo abierto que genera un único sistema de coordenadas. No pasa así con  $U(1)$ ni mucho menos con las generalizaciones apropiadas que se requieren para describir la interacción fuerte o débil. En todos esos casos, un experimentador solitario estaría obligado a usar varios  sistemas de coordenadas locales para cubrir todo su universo. Y entonces el estudio de cambio de escala, de coordenadas o de gauge sería imprescindible. Por fortuna en el caso del electromagnetismo esto es elemental.

Supongamos entonces que hay un experimentador estudiando el electromagnetismo.  En cuanto al cambio de coordenadas espacio-temporales no nos ocuparemos en ente momento (parte del problema ya lo resolvimos en el estudio de la relatividad).  Para describir al grupo $U(1)$  es suficiente particionarlo en dos abiertos,  $U$ y $V$, cada uno que cubra más de la mitad de la circunferencia.   El problema es que en $U\cap V$ se tiene  que lograr un acuerdo sobre la física observada.

Supongamos pues que   sobre el primer abierto utilizamos la forma diferencial $\omega =  -ieA_\mu dx^\mu +d\theta $, mientras que en el otro usamos la  forma diferencial $\omega '=  -ieA'_{\mu} dx^\mu +d\theta ' $. Pero la física está contenida no en las formas diferenciales de orden uno, $A$  o $ A'$, sino en su derivada. Por tanto, para que la física sea la misma, 

 $d\omega =  -iedA_\mu dx^\mu + d^2 \theta = -iedA_\mu dx^\mu $

 debe ser igual a

 $d\omega '=  -iedA'_{\mu} dx^\mu + d^2 \theta '  =  -iedA'_{\mu} dx^\mu $
 
 Es decir:
 
$ -iedA_\mu dx^\mu  = -iedA'_{\mu} dx^\mu $

O bien 

$dA = dA' $

 por lo que ambas expresiones deben diferir por un cero. En términos de derivación, cero se escribe como $0 = d^2 $. Por lo tanto, concluimos que

$dA = dA'  + 0 = dA'  + d^2 \phi = d(A'  + d\phi )$

 $A=A ' + d\phi$

donde $\phi$ es una función escalar, para que $d\phi$ sea una 1-forma.  

\bigskip

Vemos  que un cambio en la   arbitrariedad de la fase, pues $A$  va sobre la parte imaginaria,  no produce ningún cambio en  el  tensor de campo. Además, la \index{libertad gauge} \textbf{libertad gauge} que ya conocíamos desde hace muchas páginas sigue vigente:  al potencial vector se le puede aumentar una forma exacta  y la física no cambia.

\section{PODER CREATIVO}

Es importante darse cuenta que ahora tenemos ideas claras sobre cómo  proseguir para tratar de geometrizar las otras interacciones. Veamos.

Para geometrizar la interacción gravitatoria hemos curvado el espacio-tiempo. Para geometrizar la interacción electromagnética, hemos visto al espacio-tiempo como al adolescente que se pone un zarcillo en cada oreja, con la aclaración de que el espacio-tiempo tiene una oreja en cada punto. Todos los zarcillos son del mismo tipo $\textbf{U(1)}$.

Sabemos que existen otras interacciones, \index{interacción débil}  \index{interacción fuerte} la \textbf{débil} y la \textbf{fuerte}. ¿ Qué podemos hacer para geometrizarlas? Nuestra obligación natural ha de ser la de preguntarnos si existen estructuras que generalicen a $\textbf{U(1)}$. La respuesta es afirmativa: las generalizaciones de $\textbf{U(1)}$ que combinan con la mecánica cuántica son los grupos de la forma \index{$\textbf{SU(n)}$} $\textbf{SU(n)}$. Si $\textbf{U(1)}$ representaba el grupo de arbitrariedades del cambio de fase para funciones de onda con una sola coordenada, $\textbf{SU(n)}$ representa el espacio de cambios de fase para funciones de onda con $n$ coordenadas. Por tanto, los elementos del álgebra de Lie de dichos grupos serán representados por matrices.

Al tratar de generalizar lo visto, desde $\textbf{U(1)}$ a grupos superiores, lo que hay que tener en cuenta fundamentalmente es que $\textbf{U(1)}$ es conmutativo mientras que $\textbf{SU(n)}$ no, pues el producto de matrices no es conmutativo.

El teorema fundamental de la teoría gauge que enseña el papel de la  conmutatividad es el siguiente: cuando el \index{grupo conmutativo} \textbf{grupo es conmutativo}, los mediadores de la interacción no interactúan directamente entre ellos. En el caso del electromagnetismo el mediador es el fotón: para que dos fotones interactúen se requiere que se materialicen en electrones y positrones y que luego ellos se desintegren en otros dos fotones. Pero cuando el  \index{grupo no  conmutativo} \textbf{grupo es no conmutativo}, los mediadores de la interacción sí interactúan entre ellos directamente. Para la interacción fuerte los mediadores son los gluones y para la interacción débil son los bosones $W$ y $Z$. Ellos no necesitan transformarse en quarks para interactuar.

Tenemos la certeza que con lo visto el estudio de esos temas  será mucho más sencillo.

\section{  REFERENCIAS}

\begin{enumerate}

\item Nash C., S. Sen \textit{Topology and the  Geometry of physics}, AP, 1983

\item Scott S. 'Some notes on geometry and quantization'  en \textit{Primer Encuentro de Geometría Diferencial en Física }, Editado por F. Torres Ardila del Tambor de Feynman, UniAndes, Bogotá 1994.

\end{enumerate}

\chapter{CONCLUSION}        

\date{}

\large

Hemos estudiado al electromagnetismo como prototipo de las teorías gauge. Gozamos ahora de bastante claridad conceptual sobre la estructura, objetivos y métodos de una teoría gauge, la cual es aquella que, por un lado, presenta indeterminaciones inherentes a la matemática que hemos inventado para estudiarlas o quizás al proceso de observación y abstracción, y por el otro, puede reescribirse en el lenguaje de la geometría copiado, adaptado o extendido  de la relatividad general.

Gracias a este nuevo enfoque pudimos ver que el \index{grupo gauge} \textbf{grupo  gauge}, introducido como el grupo de arbitrariedades de una teoría, también es el grupo que describe la naturaleza íntima de la interacción.

En medio de tanto tecnicismo, también se deja entrever una cierta magia, algo artístico en todo esto. Es mejor tenerlo presente.

\normalsize

\printindex

\end{document}